\documentclass [11pt,makeidx]{amsbook}

\usepackage{comment}
\usepackage{amssymb} \usepackage{amsfonts} \usepackage{amsmath}
\usepackage{amsthm} \usepackage{epsfig} 
\usepackage[cmtip,arrow]{xy}
\usepackage{amsmath,amssymb,enumerate,pb-diagram,pb-xy}

\usepackage[applemac]{inputenc}
\input xy
\xyoption{all}

\usepackage{caption}

\newtheorem{theorem}{Theorem}[chapter]
\newtheorem{lemma}[theorem]{Lemma}
\newtheorem{question}[theorem]{Question}
\newtheorem{corollary}[theorem]{Corollary}
\newtheorem{definition}[theorem]{Definition}
\newtheorem{example}[theorem]{Example}

\newtheorem{proposition}[theorem]{Proposition}

\newtheorem{remarks}[theorem]{Remarks}
\newtheorem{remark}[theorem]{Remark}
\numberwithin{section}{chapter}
\numberwithin{equation}{chapter}

\newcommand{\Aa}{{\mathcal A}}
\newcommand{\Bb}{{\mathcal B}}
\newcommand{\Cc}{{\mathcal C}}
\newcommand{\Dd}{{\mathcal D}}
\newcommand{\Ee}{{\mathcal E}}
\newcommand{\Ff}{{\mathcal F}}
\newcommand{\Gg}{{\mathcal G}}
\newcommand{\Hh}{{\mathcal H}}
\newcommand{\Ii}{{\mathcal I}}

\newcommand{\Kk}{{\mathcal K}}
\newcommand{\Ll}{{\mathcal L}}
\newcommand{\Mm}{{\mathcal M}}
\newcommand{\Nn}{{\mathcal N}}
\newcommand{\Oo}{{\mathcal O}}
\newcommand{\Pp}{{\mathcal P}}
\newcommand{\Qq}{{\mathcal Q}}
\newcommand{\Rr}{{\mathcal R}}
\newcommand{\Ss}{{\mathcal S}}
\newcommand{\Tt}{{\mathcal T}}
\newcommand{\Uu}{{\mathcal U}}
\newcommand{\Ww}{{\mathcal W}}
\newcommand{\Vv}{{\mathcal V}}

\newcommand{\Xx}{{\mathcal X}}

\newcommand{\CC}{{\bf C}}
\newcommand{\DD}{{\bf D}}
\newcommand{\BB}{{\bf B}}
\newcommand{\TT}{{\bf T}}

\newcommand{\WW}{{\bf W}}
\newcommand{\EE}{{\bf E}}
\newcommand{\VV}{{\bf V}}
\newcommand{\HH}{{\bf H}}

\newcommand{\PP}{{\bf P}}

\newcommand{\ZZ}{{\bf Z}}
\newcommand{\KK}{{\bf K}}

\newcommand{\bb}{{\bf b}} 
\newcommand{\UU}{{\bf U}}

 \newcommand{\pp}{{\bf p}}

\newcommand{\cG}{{\mathfrak c}}
\newcommand{\HG}{{\mathfrak H}}

\newcommand{\gG}{{\mathfrak g}}
\newcommand{\sG}{\mathfrak s}
\newcommand{\lG}{\mathfrak l}
\newcommand{\pG}{\mathfrak p}
\newcommand{\mG}{\mathfrak m}
\newcommand{\eG}{\mathfrak e}
\newcommand{\rG}{\mathfrak r}

\newcommand{\GG}{{\mathfrak G}}
\newcommand{\vG}{\mathfrak v}
\newcommand{\wG}{\mathfrak w}

\newcommand{\PG}{\mathfrak P}

\newcommand{\qG}{\mathfrak q}
\newcommand{\iG}{\mathfrak i}
\newcommand{\MG}{\mathfrak M}

\newcommand{\VG}{\mathfrak V}

\newcommand{\tG}{\mathfrak t}
\newcommand{\nG}{\mathfrak n}

\newcommand{\hG}{\mathfrak h}
\newcommand{\fG}{{\mathfrak f}}

\newcommand{\C}{{\mathbb C}}
\newcommand{\Q}{{\mathbb Q}}
\newcommand{\Z}{{\mathbb Z}}
\newcommand{\R}{{\mathbb R}}

\newcommand{\N}{{\mathbb N}}

\newcommand{\connectedsum} {{\#}}

\def\cvd{\hfill$\Box$}

\def\Dim{\emph{Proof : }}
\def\cvd{\nopagebreak\par\rightline{$_\blacksquare$}}

 \newcommand{\x}{\times}
    
 \def \cs{\#}

\author{Riccardo Benedetti}

\address{Dipartimento di Matematica \\
Universit\`a di Pisa \\
Largo B.~Pontecorvo 5 \\
56127 Pisa, Italy}

\email{riccardo.benedetti@unipi.it}

\title{Lectures on Differential Topology}

\subjclass[2000]{}


\thanks{}

\begin{document}

\frontmatter

\maketitle

\chapter*{Preface}
Along the years  I have teached several times courses of differential topology for the  master 
curriculum in Mathematics at the University of Pisa. Typically the class was attended by students
who had accomplished (or were accomplishing) a first three years curriculum in mathematics, together with a few peer
physicists and a few beginner PhD students. With the constraints imposed by their presumable
knowledges, time after time a certain body of topics, in different combinations, as well as a certain 
way to present them stabilized. This textbook summarizes such  teaching
experiences, so it keeps a character of ``lecture notes" rather than of a comprehensive and systematic treatise.
It happens in a class to prefer a shortcut towards some interesting application,
giving up the largest generality. Similarly in this text, for example, we will mainly focus on {\it compact} manifolds
(especially when we consider the sources of  smooth maps). This allows  simplifications
in dealing for instance with function spaces or with certain ``globalization procedures'' of maps. There is
already a plenty of interesting facts concerning compact manifolds, so we will do it without remorse.

There is a lot of classical wellknown references (like \cite{M1}, \cite{GP}, 
\cite{H},  \cite{M2}, \cite{M3}, \cite{Mu}, $\dots$)  which I used  in  preparing  the courses 
and have strongly influenced these pages.  
So, why a further texbookt on differential topology?   An important motivation came to 
me by looking at the personal 
polished notes of a few good students glimpsing the lines of a reasonable text,  
together with the remark made by someone of them  that 
``they had not been able to find anywhere some of the topics treated in the course".  
It would be very hard to claim any `originality' in dealing with such a classical matter.  However,
the last sentence has perhaps a  grain of truth, at least referring to textbooks mainly 
addressed to undergraduate readers. Let us indicate a main example. A theme of this text (alike others) 
is the synergy between {\it bordism} and {\it transversality}.  One of the beforehand mentioned constraints is that
we cannot assume any familiarity with algebraic topology  or homological algebra (besides perhaps the very basic 
facts about homotopy groups); on another hand, it is very useful and meaningful to dispose of a 
(co)-homology theory suited to embody several differential topology constructions.  We will 
show that (oriented or non oriented) bordism provides instances of covariant so called
``generalized'' homology theories for arbitrary pairs $(X,A)$ of topological spaces, constructed via geometric means.
Then, by restricting to compact smooth manifolds $X$, and after a reidexing of the bordisms modules by the {\it codimension}
(so that they are now called cobordism modules),
transversality allows to incorporate the bordism modules
into  a {\it contravariant cobordism functor} with as target the category of {\it graded rings};  
also the product on cobordism modules is defined  by direct geometric means. This multiplicative structure is a substantial enrichement 
and will lead to several important  and often very classical applications. For example it is  the natural contest for unavoidable
topics like the degree theory or the Poincar\'e-Hopf index theorem. Once the cobordism product has been
well-defined, we are exempted from reproving it case by case for each specialization/application; moreover, the emphasis
is on the ``invariance up to bordism"  rather than on the ``invariance up to homotopy" as it happens for most
established references. Not assuming any familiarity with algebraic 
topology,  this presentation could also be useful as an intuitive, geometrically based introduction to 
some themes of that discipline. 

Overall this text, besides the very foundation topics, is a collection
of themes - whose choice certainly is also matter of personal preference - 
in some case advanced and of historical importance, with the common feature that they
can be treated with ``bare hands".  This means by just combining  specific differential-topogical  
cut-and-paste procedures and applications of transversality, mainly through the cobordism 
multiplicative structure.  This geometric constructive character provides the `tone' of this text, 
would be accessible to motivated master undergraduate students, to PhD students and also useful to a 
more expert reader in order to recognize 
very basic reasons for some facts already known to her/him as resulting from more advanced  
theories and/or technologies. 

\bigskip
 
 {\it Dedico tutto questo ai miei nipoti Pietro(lino) e Martin(in)a}   

\bigskip

\bigskip
 
Riccardo Benedetti              
\bigskip

Sassetta, July, 2019 

\tableofcontents


\mainmatter

\chapter*{Introduction}
These lecture notes have been conceived having in mind
a typical class of rather good and motivated students who have accomplished (or are 
accomplishing) a first three years curriculum in mathematics, then  being aware that their mathematical
background is likely limited. For example, besides very basic facts about homotopy groups, 
it is not assumed any familiarity with algebraic topology and 
homological algebra; even with respect to general topology
we assume  some knowledge about compactness in Hausdorff  second-countable topological spaces but not
about paracompactness. 

In a sense, the most natural reading of this text
is the linear one from the beginning to the end. However, different subpaths
and different combinations of the matter are possible and meaningful.
Referring to the teaching experiences which have originated these pages, 
never a single course has covered the whole  content of the book as well as the main part of every 
chapter has been sometimes experimented in a course implementation. 

The text (alike the lectures it derives from) intends to be uniformly accurated for what
concerns definitions, statements and description of the main constructions; it also aims
to develop an articulated and coherent discourse. On the other hand
it is intentionally not uniform for what concerns details in proofs. At some points (especially in the
final few chapters) the attive participation of the reader is  required in order to 
complete some argument or to check some claim. Overall it is assumed a quite collaborative and motivated
reader. For this reason we have found not necessary to add a list of exercises chapter by chapter.

We have made a  contained use of figures; basically only pictures containing substantial,
not only allusive information have been introduced. Drawing pictures is often an useful
support in order to follow  a geometric argument, but also this is left to
reader's initiative.

The bibliography is very far from being exhaustive; besides a few classical refrences which have 
certainly influenced these pages  (like \cite{M1}, \cite{GP}, \cite{H},  \cite{M2}, \cite{M3}, \cite{Mu})
we just list the texts which have been actually cited.

We will make some (very moderate indeed) use of the language of categories.
We collect in an appendix the few necessary notions.
Differential Topology concerns  the category of {\it smooth manifolds} 
and  {\it smooth maps}; this includes the study of smooth
manifolds considered up to {\it diffeomorphism} that is the equivalence in that category.
A first unpreventable task is to define these objects and morphisms. We do it from scratch 
in Chapters \ref{TD-LOCAL}, \ref {TD-DIFF} and \ref {TD-SMOOTH-MAN} by extending progressively
the category, from the one of open sets in euclidean spaces and smooth maps, passing through
the category of  {\it embedded} smooth manifolds in some euclidean space and ending with
the category of ``abstract" smooth manifolds defined by abstraction of some properties of embedded
ones. Along these successive generalizations  we develop the accessory notions of {\it submanifold},
manifold with {\it boundary},  (orientable) {\it oriented} manifold with {\it oriented boundary}. 
Basic notions like the one of immersion, submmersion, embedding,  
smooth homotopy, isotopy or diffeotopy between smooth maps are also introduced.

In most applications we will  focus on {\it compact} manifolds, especially when we consider the sources of  smooth maps. 
Then intentionally we will not present the most general version of many results; 
there is already plenty of interesting facts concerning compact manifolds and the compacteness assumption
simplifies a lot many arguments for example in dealing with function space topology or with cut-and-paste
constructions where one can use only {\it finite} partitions of unity, avoiding any reference to paracompactness. 
Moreover, exploiting the fact that every abstract compact manifold can be eventually embedded in some euclidean space, 
several important issues like a tubular neighbourhood theory are developed first for embedded
manifolds and then extended to the whole range of compact manifolds.
These claims will be substantiate below.

In Chapter \ref{TD-LOCAL} we assume the knowledge of basic several variables
differential calculus and  we collect some  facts
concerning smooth maps between open sets of euclidean spaces.
Some of these facts (such as the {\it inverse map} theorem and its 
geometric applications to the local normal form of immersions and summersions) 
should be familiar to the reader. Others are presumably less
familiar such as Morse's lemma, the isotopic linearization of diffeomorphisms
of $\R^n$, bump functions, the smooth homogeneity of connected open sets (which later extends
to arbitrary connected manifolds).  Two characteristic features of differential
topology already manifest themselves. On one side there is a sort
of ``local rigidity'':  up to local change of smooth coordinates,
linear algebra provides the actual local
models in many `generic and stable' smooth situations. 
On another side, smooth maps are very ``flexible'', the existence of
bumb functions being a typical instance of it. This will be the key for globalization
procedures and cut-and-paste constructions. Flexibility is a quality one expects 
from a topological theory, but this is moderated by that sort of local rigidity 
which (at least in suitable generic and stable situations) allows to have a good geometric control;
this moderate flexibility eliminates too ``wild'' phenomena occurring in 
general topology, even dealing with merely  topological manifolds, or allows
simple proofs of facts (like the invariance of dimension up to diffeomorphism)
whose topological counterpart holds as well but is much more delicate.
Moreover, especially the homogeneity property indicates that the veritable questions
in differential topology concern the {\it global} structure of manifolds.

In Chapter \ref{TD-DIFF}, we straighforwardly extend the
notions of smooth map and diffeomorphism to {\it arbitrary}
topological subspaces of some euclidean spaces; then an embedded smooth 
manifold $M$ of dimension $m$ is defined as a topological subspace of 
some $\R^n$ which is locally diffeomorphic to open subsets of $\R^m$.  
Although not so demanding, this extension leads to a plenty of embedded manifolds beyond 
the open sets, including very familiar objects like the graphs of smooth maps between open sets
which ultimately are  the local model for every embedded smooth manifolds.
Embedded smooth manifolds are naturally endowed with a maximal {\it atlas} of smooth {\it charts}
(with corresponding smooth {\it local coordinates}) and smooth maps
between embedded manifolds have natural {\it representations in local coordinates}. These notions
are the key for the final abstraction made in Chapter \ref{TD-SMOOTH-MAN}.

After having stressed in Chapter \ref{TD-LOCAL}
the fuctorial character of the elementary {\it chain rule}, we follows
along the successive generalizations the costruction of the fundamental covariant {\it tangent functor}
which associates to every  manifold its {\it tangent bundle} and to every smooth map its {\it tangent map}.
The tangent functor is a main source of invariants of smooth manifolds. In the embedded category 
this is a direct and transparent construction, being completely internal to the category. 
Modeled on them, in Chapter  \ref {TD-DIFF}  we
state the general notions of {\it embedded smooth fibre bundle} (in particular {\it vector bundle}) 
and fibred maps and we elaborate on different notions of fibred bundle equivalence.
In a sense tangent bundles and maps for embedded manifolds 
impose by themselves, starting from the basic ones for open sets, based on the chain rule. 
In the abstract case they must be somehow ``invented'', with the constraint to agree  with what we have 
already done in the embedded category. Probably this is the most demanding extension 
passing from the embedded to the abstract setting. Eventually this leads us in Chapter \ref{TD-SMOOTH-MAN}
to the general notions of {\it principal bundles} with a given {\it structural group} $G$ and its {\it associated
bundles}, governed by suitably defined $G$-valued {\it cocycles}.

Our typical student is probably already aware of the topology of the uniform convergence on compact sets
of  continuous maps between open sets of euclidean spaces. This directly extends to 
$\Cc^r$ maps, $r\geq 0$, in terms of the uniform convergence on compact  
sets of the maps and  their partial derivatives up to the order $r$. This restricts to the set of smooth 
maps and we can also consider on it the union topology over $r\in \N$. 
By using the representation in local coordinates, the definition of these function spaces
extends to smooth maps $f: M \to N$ between smooth manifolds (embedded as well as abstract) 
giving us the spaces $\Ee^r(M,N)$ endowed with the so called $\Cc^r$ {\it weak topology} and the space
$\Ee(M,N)$ endowed with the union topology.  The adjective ``weak" alludes
to further function space topologies, the so called {\it strong topologies}. These coincide with the weak ones if the source
manifold is compact while they are much finer otherwise and are aimed  to have a control `at infinity'. 
However, we will not treat the strong topology because in the relevant applications considered in this text
the source manifold $M$ will be {\it compact}. For example we show that if $M$ is compact, $f:M\to N$
is an embedding if and only if it is an injective immersion and that  immersions, summersions
and embeddings respectively form a (possibly empty) {\it open set} in  $\Ee(M,N)$.

Chapters  \ref{TD-DIFF}  to   \ref{TD-COMP-EMB} deal with the embedded category.
Chapter \ref{TD-STIEF-GRASS} is devoted to a detailed presentation of two distinguished
families of manifolds, that is {\it Stiefel and Grassmann  manifolds}. Stiefel manifolds are naturally
embedded and we provide also embedded models for the Grassmann ones. 
Besides the fact that
they are non trivial examples of (embedded) smooth manifolds, they will be crucial in the study
of vector (and frame) bundles on arbitrary manifolds.

In Chapter  \ref{TD-EMB-VB} we introduce the  fundamental  {\it
  pull-back} construction of embedded smooth fibred bundles. Then we apply it 
  to the so called {\it tautological (vector or frame) bundles} over the Grassmann
manifolds. This can be considered as a powerful machine that produces
embedded vector bundles (and the associated frame bundles)
over embedded smooth manifolds and naturally incorporates the
tangent bundles and their tensorial relatives.  All the bundles constructed
in this way, partitioned by the rank, are considered up to `strict equivalence'.  
After having constructed via a suitable 
limit procedure the {\it infinite grassmannian}  $\GG_{\infty,k}$ of $k$ planes in $\R^\infty$
with its limit tautological bundles, a main  result of the chapter is the classification of 
such rank-$k$ vector bundles  up to  strict equivalence  over a {\it compact} manifold $M$;
this establishes a bijection with $[M,\GG_{\infty,k}]$ the set of {\it homotopy classes} of smooth
maps from $M$ to $\GG_{\infty, k}$. The compactness assumption
simplifies the discussion and we will not touch any generalization of this result. A typical way to
get algebraic/topological invariants is to construct functors from some subcategory of the  topological spaces
to some category of algebraic structures (groups, rings, vector spaces, $\dots$). At the end of  Chapter  \ref{TD-EMB-VB} 
we present a non trivial  implementation of this idea based on  such vector bundles. 
By augmenting the strict equivalence to a suitable  {\it stable equivalence}, we realize that the quotient set $\KK_0(M)$
of the whole collection of the vector bundles considered above (all ranks confused) carries a natural {\it ring structure};  
together with the pull-back construction, this eventually builds a {\it contravariant functor} from the (sub) category of compact manifolds
to the category of abelian rings which verifies the {\it homotopy invariance property}.

In Chapter \ref{TD-COMP-EMB} we focus on embedded {\it compact} manifolds. We develop in such a framework
a theory of {\it tubular neighbourhoods} of submanifolds and of {\it collars} for the boundary of a manifold with
boundary. This is ultimately based on the pull-back construction considered above and on embedded {\it normal bundles}
orthogonal to the tangent bundle with respect of a given riemannian metric (for instance the standard one) on the ambient
euclidean space. Then we give some applications of this technology. For simplicity let us limit here to boundaryless
manifolds. If $M$ is compact and $N$ is also compact or more generally is embedded in some $\R^m$ being
furthermore a closed subset, then we prove that smooth maps are {\it dense} in $\Cc^r(M,N)$ for every $r\geq 0$.
Primary topological invariants as the fundamental group or more generally higher homotopy groups are
defined in general in terms of homotopy classes of {\it continuous} maps defined on spheres. 
As an application of the above density theorem we see that they can be equivalently defined in terms
of smooth homotopy between smooth maps $f:S^n\to N$. Another important application is for every $r\geq 1$, the approximation
of compact $\Cc^r$-manifolds $M\subset \R^h$ by smooth embedded manifolds and the existence and uniqueness
up to diffeomorphism of a smooth structure on every such a manifold $M$. Exploiting the fact that our embedded
grassmannian are not only embedded smooth manifolds but actually {\it regular real algebraic sets}, and that the tautological bundles are
also real algebraic, we outline Nash's celebrated result that every embedded smooth manifold $M\subset \R^h$ can be approximated
by a regular sheet  of a real algebraic set of $\R^h$ (shortly   by a {\it Nash manifold}) and that every 
compact embedded smooth manifold admits a structure of Nash manifold unique up to Nash diffeomorphism.
We state the smooth {\it Sard-Brown} theorem which is the base of {\it transversality} which will be more systematically
developed in Chapter \ref {TD-TRANSVERSE}; here we anticipate some manifestation.
We discuss also a version of Sard-Brown in the category of embedded Nash manifolds.
In the general settings the result is expressed in measure-theoretic terms while in the Nash case it is  purely a geometric
statement, as well as the proof. Thanks to the above approximation theorem, in many situations faced in differential topology, smooth and Nash 
Sard-Brown theorem  can be used indifferently. By using the restriction to $M\subset \R^h$
of `generic' linear projections of $\R^h$ onto lines, we show that {\it Morse functions} form an {\it open and dense} subset of $\Ee(M,\R)$.
We study also some instance of generic linear projections onto hyperplanes and prove {\it `easy' Whitney's immersion/embedding theorem}:
if $m=\dim M$, $h\geq 2m$, then $M$ can be immersed into $\R^{2m}$, if $h\geq 2m +1$, then $M$ can be embedded into
$\R^{2m +1}$. 

At the end of Chapter  \ref{TD-SMOOTH-MAN} we show that every abstract compact smooth manifold can be embedded in some
$\R^n$ as well as abstract vector bundles on it are strictly equivalent to embedded ones. Hence all the results of 
Chapters  \ref{TD-DIFF}  to   \ref{TD-COMP-EMB} hold as well for abstract compact smooth manifolds, provided that they are
considered up to diffeomorphism. In particular we have the remarkable fact  that all compact manifolds of a given dimension can be 
immersed/embedded into a  same euclidean space.
As we are mainly concerned with compact manifolds, the abstraction made in Chapter 
\ref{TD-SMOOTH-MAN} would appear a bit superfluous. However there are natural constructions 
to build new compact manifolds, starting from given ones; it would be artificial to force them to deal from 
the beginning in the embedded setting. It is more convenient to use the embedding result {\it a posteriori}, in order to exploit 
the facts already established for compact embedded manifolds.

Chapters \ref{TD-LOCAL} to  \ref{TD-SMOOTH-MAN} (with perhaps the exception of the end of Chapter \ref{TD-COMP-EMB} about the rings $\KK_0(*)$)
form the strictly fundational part of this text. The subsequent chapters articulate somehow a more advanced discourse. 
 
In Chapter \ref{TD-CUT-PASTE} we collect several constructions that produce a new compact manifold by modifying
a given one. At first we prove a so called {\it Thom's lemma} about the extension of any isotopy defined on a compact source
manifold to an ambient diffeotopy; this is the main tool in order to prove that every implementation of such constructions is uniquely
well defined up to diffeomorphism. Among basic cut-and-paste procedures we recall {\it gluing along diffeomorphic
boundary components}, {\it connected sum} with a discussion about the related notion of {\it twisted spheres}, 
{\it attaching a  $p$-handle} i.e. a standard handle $D^p\times D^{m-p}$ of index $p$  to an $m$-manifold $M$ along an embedding in
$\partial M$ of the {\it attaching tube} $S^{p-1}\times D^{m-p}$ (in particular, by attaching a $0$-handle one creates a new
component $D^m$ with boundary $S^{m-1}$, by attaching a $m$-handle we cap a spherical component of $\partial M$ (if any)
with an $m$-disk). 
In many cases the immediate result 
is rather a {\it smooth manifold with corners}. Corners also arise by taking the product of two manifolds with non empty boundary.
Then we discuss a standard procedure of {\it smoothing the corners} that produces ordinary smooth manifolds well defined
up to diffeomorphism. We discuss also the {\it strong Whitney embedding/immersion  theorem} of any $m$-dimensional compact
manifold $M$ in $\R^{2m}$ and in $\R^{2m-1}$ respectively. The main difference with respect to the above `easy' Whitney theorems 
is that the 
strong ones are not  enterely based on `generic position arguments' (i.e. transversality); in fact they are achieved by performing
a robust modification of the `generic' immersions provided by the easy immersion theorem, or of determined `generic' maps of
$M$ in $\R^{2m-1}$, respectively. The proof of the strong embedding  introduces the so called {\it Whitney's trick} in order to eliminate 
pairs of selfintersection points in the image of a generic immersion in $\R^{2m}$;
this will be
reconsidered in Chapter \ref {TD-HIGH} and in Chapter \ref{TD-4}. By elaborating on the strong immersion theorem, we present
{\it Rohlin's embedding theorem in $\R^{2m-1}$ up to surgery}; in particular this shows that for every $M$ as above there is $M'$
such that the disjoint union $M\amalg M'$ is the boundary of a compact $(m+1)$-manifold $W$, and $M'$ can be {\it embedded} 
into $\R^{2m-1}$. In the last section of the chapter we describe the modification obtained by {\it blowing up a manifold $M$ along
a smooth centre $X\subset M$}; this replaces $X$ with its {\it projectivized normal bundle} in $M$.

In Chapter  \ref {TD-TRANSVERSE} we develop the {\it transversality} concept in a more systematic way. As usual the source manifold $M$
is compact possibly with non empty boundary and for simplicity we assume here that the target manifold $N$ is also compact 
and boundaryless; $Z$ is a boundaryless compact submanifold of $N$. There are two kinds of {\it basic tranversality theorems}.
The first kind concerns a certain geometric tameness under transversality hypothesis:
if $f:M \to N$ is transverse to $Z$, then $(Y,\partial Y) = (f^{-1}(Z), (\partial f)^{-1}(Z))$ is a
nice `proper submanifold' of $(M,\partial M)$ of the same {\it codimension} of $Z$ in $N$. 
There is also a specialization within the category of oriented manifolds.
The second kind states that transverse maps are generic and stable that is they
form an {\it open and dense} set in $\Ee(M,N)$ (moreover there is a relative version concerning maps which concide
on $\partial M$, provided that this restriction is already transverse to $Z$ by itself). In other words, by means of an arbitrarily small 
perturbation, every map $f:M\to N$ becomes transverse (hence with a nice geometric behaviour) 
and transversality is a {\it stable property} with respect to small perturbations. The bridge between the two kinds of theorems
is represented by the so called {\it parametric transversality} whose proof is substantially based on the Sard-Brown theorem. 
These so called basic transversality theorems suffice for most later applications in this text. However, transversality i.e. `general position'
reasoning  is a profound, potent and pervasive paradigm beyond such basic results. 
Without any pretention of completeness in the second part of the chapter we collect a few instances of further applications
(including the notion of `generic immersion', already employed while discussing Whitney's strong embedding theorem). 

In Chapter \ref {TD-HANDLE} we formalize the notion of a {\it smooth triad} $(M,V_0,V_1)$ where $M$ 
is a compact smooth $m$-manifold (possibly with empty boundary, so that $(M,\emptyset,\emptyset)$ is allowed) and $V_0$ and $V_1$ 
are union of connected components of 
$\partial M$ in such a way  that the boundary is the disjoint union $\partial M = V_0    \amalg  V_1$; we define generic Morse functions $f: M\to [0,1]$
on a triad ($f^{-1}(j)=V_j$, $j=0,1$, $f$ has only non degenerate critical points placed outside a neighbourhood of $\partial M$)
whose density and stability is assured by the basic results of Chapter \ref{TD-TRANSVERSE}. A main achievement of  Chapter  
\ref {TD-HANDLE} is that every such a Morse function carries a {\it handle decomposition of the triad} that is a way to recostruct
the triad (up to diffeomorphism) starting from a collar of $V_0$ in $M$ and by attaching succesively a handle of index $p$
for every non degenerate critical point of index $p$ of $f$. In a sense Morse functions on any triad are used as a tool to prove the {\it existence}
of such handle decompositions. Then handle decompositions are managed by themselves as we do not face the issue  whether
every decomposition is carried by some Morse function. 
Associated to every decomposition of a triad $(M,V_0,V_1)$ there is a {\it dual decomposition} of the triad $(M,V_1,V_0)$ 
where every $p$-handle is converted into a $(m-p)$-handles and these are attached backward starting from a collar of $V_1$
in $M$. If the decomposition is associated to a Morse function $f$, then the dual  is associated to $1-f$.
We point out two {\it basic moves} which modify any given decomposition
without changing the triad (up to diffeomorphism): the so called {\it sliding handles} which is nothing else than the possibility
of modifying any attaching map up to isotopy already treated in Chapter  \ref{TD-CUT-PASTE}; the {\it elimination/insertion of
pairs of complementary handles}. We show some elementary instances of specialization (`reordering') or simplification
(`elimination of $0$- and $m$-handles') of handle decompositions by means of applications of the basic moves. 
As a simple but important application we get the classification up to diffeomorphism of compact $1$-dimensional 
manifolds, confirming the intuition: a connected compact $1$-manifold either is diffeomorphic to $S^1$
or to the $1$-disk $[-1,1]$.

In Chapter \ref{TD-BORDISM} we develop {\it bordism}. There is an unorieted version and an oriented one.
Two (oriented) compact boundaryless $m$-manifold $M_0$ and $M_1$ are (oriented) bordant manifolds
if  ($M_0 \amalg -M_1$) $M_0\amalg M_1$ is the (oriented) boundary of a compact (oriented) manifold
$W$. The quotient set of the relation generated by `being bordant' and (oriented) diffeomorphisms 
is denoted by ($\Omega_m$) $\eta_m$ and is a ($\Z$-module) $\Z/2\Z$-vector space,
the operation being induced by the {\it disjoint union}. If $X$ is any topological space, a continuous map
$f:M \to X$ is called a {\it singular} smooth $m$-manifold in $X$ and we can extend the definition of bordism
to such singular manifolds, hence of the above modules
to $\Omega_m(X)$ or  $\eta_m(X)$, also denoted by $\Bb_m(X;R)$, $R=\Z, \Z/2\Z$.
When $X= {\rm pt}$ we recover the initial modules because the maps $f$ are immaterial in this case. 
Moreover, we can define also relative versions  $\Bb(X,A;R)$ for topological pairs $(X,A)$ ($X$ being as usual
identified with $(X,\emptyset)$).
 We prove that in this way we define a {\it covariant functor} from the category of topological pairs
 to the category of $R$-modules  which turns out to be a {\it generalized homology theory}: this means that
 all {\it Eilenberg-Steenrod axioms} are satisfied with the possible exception of `dimension'; its failure depends
 on the non triviality of  $\Bb_m({\rm pt};R)$, $m\geq 1$, an issue that will be somehow discussed along the rest of the text.   
We discuss some relationships between bordism and homotopy group functors.

In Chapter \ref{TD-COBORDISM} we specialize bordism assuming that $X$ is a compact boundaryless smooth
manifold. First, alike the homotopy groups,  thanks to the approximation theorems of Chapter
\ref {TD-COMP-EMB} it is not restrictive to deal only with smooth maps $f:M\to X$. The bordism modules
$\Bb_m(X;\Z/2\Z)$ are indexed over $\Z$ by postulating that they are the trivial module $0$ if $m<0$.
We formally reindex them by the {\it codimension}, by setting $\Bb^r(X;\Z/2\Z)=\Bb_m(X;\Z/2\Z)$, $r= \dim X - m$,
so that they are trivial if $r> \dim X$ and are now called {\it cobordism modules}. The key point is that by combining a 
slight extension of the basic transversality
theorems of Chapter  \ref {TD-TRANSVERSE} with variations on the {\it pull-back construction} already used
in the framework of Chapter  \ref{TD-EMB-VB}, we incorporate $X \Rightarrow \oplus_r \Bb^r(X;\Z/2\Z)$
into a {\it contravariant functor} from the (sub)category of compact boundaryless smooth manifolds
to the category of {\it graded rings}; this means that    $ \oplus_r \Bb^r(X;\Z/2\Z)$ is endowed
with a  multiplicative structure which distributes itself into a family of  $\Z/2\Z$-bilinear maps
$ \sqcup: \Bb^r(X;\Z/2\Z)\times \Bb^s(X;\Z/2\Z)\to \Bb^{r+s}(X;\Z/2\Z)$
defined geometrically via transversality and an implementation of the pull-back construction.
If $X$ is oriented we can perform all the construction within the oriented category, that is 
in terms of the $\Z$-modules $\Bb^r(X;\Z)$. If $\alpha = [M_1]$ and $\beta =[M_2]$
are represented by submanifolds of $X$, then $\alpha \sqcup \beta$ is represented by
any transverse intersection $M'_1 \cap M'_2$ where $M'_j$ is a suitable small perturbations of $M_j$,
$j=1,2$. Over $R=\Z/2\Z, \Z$, case by case, the product
verifies the non commutativity relation $\alpha \sqcup \beta = (-1)^{rs} \beta \sqcup \alpha$
which can be also cecked in geometric way.
If $X= {\rm pt}$, then the product reduces to $[M]\sqcup [N]= [M\times N]$. 
 If $r+s=\dim X=n$ and $X$ is connected (possibly oriented), then $\Bb^n(X;R)=R$ and the product $\sqcup$
induces a linear map $\phi^r: \Bb^r(X;R) \to {\rm Hom}(\Bb_r(X;R), R)$; in many situations it
is convenient to consider the quotient module $\Hh^r(X;\R)= \Bb^r(X;R)/\ker (\phi^r)$ with the induced
linear injection $\hat \phi^r: \Hh^r(X;R)\to  {\rm Hom}(\Hh_r(X;R), R)$. In particular if $X$ is oriented
$\Hh^r(X;\Z)$ is torsion free. If $X$ is again connected (possibly oriented) and $\dim X= 2m$, 
then we have the {\it intersection form} $\sqcup : \Hh^m(X;R)\times \Hh^m(X;R)\to R$ which is symmetric
if either $R=\Z/2\Z$ or $R=\Z$ and $m$ is even, it is antisymmetric otherwise. Sometimes
it is expressed as $ \bullet : \Hh_m(X;R)\times \Hh_m(X;R) \to R$.

The cobordism multiplicative structure is a substantial enrichement.
In Chapter \ref{TD-CB-APPL} we collect a few classical applications: the {\it fundamental class} $[X]\in \Hh^0(X;R)$
when $X$ is connected an possibly oriented; {\it Brouwer's fixed point theorem} for continuous maps $f:D^n \to D^n$ , $n\geq 1$;
a {\it separation theorem} for hypersurfaces in $S^n$, $n>1$; {\it intersection and linking numbers};  the $R$-{\it degree}, 
case by case, of continuous maps $f:M \to N$ between (possibly oriented) compact, connected,  boundaryless smooth manifolds;
a proof of the {\it fundamental theorem of algebra}; {\it Borsuk-Ulam theorem}. We define also 
the {\it Euler class} $\omega(\xi) \in \Bb^k(X;R)$ of a rank-$k$ vector bundle $\xi$ over $X$, defined by the transverse 
self-intersection of the zero section of $\xi$ in its total space. A non zero Euler class is a primary obstruction
to the existence of a nowhere vanishing section of $\xi$.

In Chapter \ref{TD-LINE-BUND} we focus on line bundles (i.e. rank-$1$) on $X$, on oriented
rank-$2$ vector bundles provided that also $X$ is oriented, and on their Euler classes in $\Bb^1(X;\Z/2\Z)$,
$\Bb^1(X;\Z)$, $\Bb^2(X;\Z)$. A key point here is that $\PP^\infty(\R)$ is a $\KK(1,\Z/2\Z)$ space,
$S^1$ is a $\KK(1,\Z)$, $\PP^\infty(\C)$ is a $\KK(2,\Z)$. This eventually gives precise information, case by case,
about $\Hh^1(X;R)$ and $\Hh^2(X;\Z)$. For example: every class in $\Hh^1(X;R)$ is the Euler class of a 
(possibly oriented) line bundle over $X$; it can be represented by an embedded (possibly oriented) hypersurface $S$
of $X$; $[S_0]=[S_1]$ if and only if the associated bundles are strictly equivalent, if and only if the (oriented) bordism between $S_0$ and $S_1$
is realized by means of a (oriented) triad $(W,S_0,S_1)$  properly embedded into $X\times [0,1]$.
Similarly for $ \Bb^2(X;\Z)$.

In Chapter \ref{TD-EP} we focus at first on the Euler class in $\Bb^m(M;\Z) = \Z$ of the tangent bundle
of a compact oriented connected boundaryless smooth $m$-manifold $M$. This integer is denoted
by $\chi(M)$ and called the {\it Euler-Poincar\'e characteristic} of $M$. Essentially by definition, it can be computed by means of
any section of $T(M)$ transvese to the zero section, that is by means of any tangent vector field on $M$
with only non degenerate zeros. This can be extended to any tangent vector field on $M$ with only isolated (not necessarily
non degenerate) zeros. This is the content of the  {\it Index Theorem}; in fact the key point is the reformulation of the sign
of a non degenerate zero in terms of the $\Z$-degree of a suitably map $f: S^{m-1}\to S^{m-1}$ defined locally
at the zero by means of the vector field; this reformulation by the degree makes sense also for any isolated zero
and well defines its {\it index}. Then $\chi(M)$ is eventually equal to the sum of such indices. Invariance of the degree
up to bordism does play a crucial role in this achievement. The characteristic is multiplicative
with respect to the product of compact boundaryless manifolds.
The value of $\chi(X)$ does not depend on the choice of the orientation of $X$; eventually $\chi(M):= \frac{1}{2}\chi (\tilde M)$
is well defined also if $M$ is not orientable, $\tilde M \to M$ being the orientation $2$-to-$1$ covering map.
We extend the index formula to well define the {\it relative characteristic} $\chi(M,V_0)$ of a triad $(M,V_0,V_1)$
by using suitable tangent vector fields on $M$, transverse to the boundary and with only isolated zeros.
The characteristic has certain homotopy invariance properties so that for example if $B$ is the total space of a disk bundle
over a boundaryless $M$, then $\chi(M)= \chi(B,\emptyset, \partial B)$. In the special case when $M$ is embedded into
$\R^h$ and $B$ is a tubular neighbourhood of $M$ in $\R^h$, this leads to the classical fact that $\chi(M)$
coincides with the degree of the {\it Gauss map} $M\to S^{h-1}$.
The extended characteristic has also remarkable additivity properties with respect to the composition of triads.
Moreover, $\chi(M,V_0)$ can be computed by means of any gradient vector field of any Morse function $f:M \to [0,1]$ on the triad.
By combining these facts we obtain for example that if $M$ is boundaryless and {\it odd} dimensional,
then $\chi(M)=0$ (use both $f$ and $1-f$ to compute $\chi(M)$ in two ways); if $V$ is even dimensional and is the 
boundary of some $M$, then   $\chi(V)\equiv 0$ mod $(2)$. It follows for example that for every even $m$, $\eta_m$
is non trivial because $\chi(\PP^m(\R))= 1$. At the end of the chapter we shortly discuss other ways (combinatorial or
algebraic/topological) to define the EP characteristic.

In Chapter \ref{TD-SURFACE} we apply several tools developed in the previous Chapters in order to classify
the compact surfaces (i.e. smooth $2$-manifolds) up to diffeomorphism and also to determine both bordisms
$\eta_2$ and $\Omega_2$.  If $M$ is a connected boundaryless compact surface, we show that $\eta_1(M)$ is a finite
dimensional $\Z/2\Z$-vector space and that the symmetric intersection form $\bullet : \eta_1(M)\times \eta_1(M)\to \Z/2\Z$
is non degenerate. We focus on its isometry class as the main invariant up to diffeomorphism. After having
established the abstract algebraic classification up to isometry of non degenerate symmetric bilinear forms on
finite dimensional $\Z/2\Z$-spaces, we show that step by step there is a perfect 2D topological counterpart: finally
every isometry class can be realized as the intersection form of some surface $M$ and two surfaces are diffeomorphic
if and only if they have isometric intersection forms. In particular we can derive from that isometry class
whether $M$ is orientable or not and the value of $\chi(M)$. At the end of the day, if $M$ is orientable then 
$M$ is the connected sum of $S^2$ with $g$ copies of $S^1\times S^1$, where  $\chi(M)=2-2g$; if $M$ is non
orientable, then $M$ is a connected sum of copies of $\PP^2(\R)$ whose number is determined by $\chi(M)$.
$\Omega_2=0$, while $\eta_2=\Z/2\Z$ generated by $[\PP^2(\R)]$. We also discuss some aspect
of the {\it stable equivalence} generated by diffeomorphisms and the elementary stabilization
consisting in performing the connected sum with $\PP^2(\R)$; in particular we refer the relationship
with so called {\it Nash rationality question} in dimension $2$.
In the subsequent chapters it will emerge the theme of the {\it quadratic  enhancements} of the intersection
form associated to the immersion of a surface into a higer dimensional manifold.
At the end of Chapter \ref{TD-SURFACE} we develop a bit the abstract theory of such quadratic
enhancements of non degenerate symmetric bilinear forms on finite dimensional $\Z/2\Z$-spaces, including the introduction
of the {\it Arf} and the {\it Arf-Brown} invariants. 

The Euler-Poincar\'e characteristic mod$ (2)$  is a first
$\eta_m$-characteristic number for every $m\geq 0$, that is it defines   homomorphisms  $\chi_{(2)}: \eta_m \to \Z/2\Z$,
surjective for even $m$. Pontryagin remarked that there is a natural way to construct
plenty of so called {\it stable} $\eta_m$-characteristic number as follows. Let $\alpha \in \eta^m(\GG_{n, m+1})$ ($n$ big enough), 
then define
$c_\alpha: \eta_m \to \Z/2\Z$, $c_\alpha(\beta)= [M]\sqcup s_M^*(\alpha) \in \Z/2\Z$,
where $M$ is any connected representative of $\beta$, $s_M: M \to \GG_{n, m+1}$
is a classifying map (uniquely defined up to homotopy if $n$ is big enough) of the {\it stable tangent bundle}
$T(M)\oplus \epsilon^1$ (the last being the product line bundle over $M$), $s_M^*(\alpha)\in \eta^m(M)$
is the pull-back of $\alpha$. In \cite{T}, Thom computed the ring $\oplus_m \eta_m$; by means of the
Pontryagin-Thom construction that is  treated in Chapter \ref{TD-PT}, this is reduced to the computation of
the homotopy groups of certain `Thom's spaces', an this can be eventually achieved by means of the powerful tools
introduced in homotopy theory since Serre's thesis \cite{Se}. A byproduct of Thom's work is the 
{\it completeness} of stable $\eta$-characteristic numbers. In other words, $\beta \in \eta_m$ is equal to zero if and only
if for every stable $\eta_m$-characteristic number $c_\alpha$ as above, $c_\alpha(\beta)=0$. In Chapter \ref{TD-ETA-CHAR} we
propose an enterely geometric proof due to \cite {BH} only based on transversality of this remarkable completeness of $\eta$-characteristic 
numbers. An analogous result for $\Omega_m$ holds as well but is more complicated
and even its formulation is beyond the limit of the present text. By similar geometric means we limit to deal with a special case,
that is we prove that if $M$ is {\it parallelizable} (i.e. its tangent bundle is strictly equivalent to the product bundle - hence $M$ is orientable) 
then $[M]=0 \in \Omega_m$, for every choice of the orientation of $M$.

Chapter \ref{TD-PT} is devoted to the {\it Pontryagin-Thom} contruction. 
The original Pontryagin construction was inventend to rephrase the
study of the homotopy groups of spheres $\pi_{n+k} (S^n)$, $k\geq 0$, $n >1$, in terms of a certain more
geometric (hence presumably more accessible at that time, about 1938)
codimension $n$ embedded oriented bordism theory with target $S^{n+k}$. 
This so called {\it framed bordism}
makes sense for arbitrary compact target $M$ in both an oriented and a non oriented version
and is related to $[M,S^n]$.
Viceversa, later Thom's extension of Pontryagin
construction was mainly intended as a way to rephrase the study of
the cobordism rings in terms of the homotopy groups
(becomed more accessible at that time, about 1954, after Serre's
Thesis) of certain so called Thom's spaces which in a sense generalize
the spheres.  So the P-T construction is a powerfull bridge between
two different ways to approach a same mathematical reality.
Concerning the determination of  $\pi_{n+k} (S^n)$,
Pontryagin succeeded for $k\leq 2$ and in Chapter \ref{TD-PT} we outline these results.
For $k=0$  one realizes that the $\Z$-degree establishes an isomorphism
between $\pi_n(S^n)$ and $\Z$. As a corollary we show that a compact connected boundaryless 
manifold $M$ is {\it combable} (i.e. it admits a nowhere vanishing tangent vector field)
if and only if $\chi(M)=0$; in particular every odd dimensional $M$ is combable.
Difficulty increases with $k$. For $k=2$ a key ingredient is the 
 Arf invariant of the quadratic enhancement of the intersection form
of every framed orientable surface in $S^6$. The hardest application of this
geometric  way is for $k=3$ and is due to Rohlin. We limit to state the result.
This is of major importance for its consequences in the theory of $4$-manifolds
and will be reconsidered in Chapter \ref{TD-4}.

In differential topology there is a precise distinction between `high' (i.e. greater or equal $6$)
dimensions and low dimensions less or equal $4$, $5$ rather being on the border between the two
regimes. The main reason is that for $d\geq 6$, Smale's (simply connected) {\it $h$-cobordism theorem}
holds and moreover there is a ``stable proof'' that is working uniformly  for all high dimensions.
This is a main application of handle decomposition theory.
The same proof does not apply to low dimensions, in some case the theorem fails, in some case
it is still an open question. In Chapter \ref{TD-HIGH} we briefly discuss this issue.
We do not give a  proof of the stable $h$-cobordism theorem; rather we focus on a main step
where the high dimension assumption is crucial. This is related to the possibility of applying
the Whitney's trick (early introduced for the strong embedding theorem) in order to eliminate pairs
of intersection points of opposite sign between transverse submanifolds of complementary dimension into
a simply conneted ambient manifold of dimension greater or equal to $5$,
and to certain  `unlinkig of spheres into spheres'.  

For `very low' dimensions $0\leq d \leq 2$ we have achieved a complete classification
of compact manifolds up to diffeomorphism. This is essentially `hopeless' for $d>2$,
even for $d=3,4$. In Chapters \ref{TD-3} and \ref{TD-4} we face some aspects of these
low dimensional theories. We stress that in both cases we do not touch the mainstream
themes of the last decades (the {\it geometrization conjecture} (now a theorem) of $3$-manifolds
or the use of powerful {\it gauge theories} applied to the study of $4$-manifolds).
We limit to develop a few classical differential topological results by applying several tools
established in the previous chapters.

In Chapter \ref{TD-3} we give a few elementary and selfcontained proofs of the primary
fact that compact orientable boundaryless $3$-manifolds are parallelizable and study
combing and framing.
An important amount of the chapter is devoted to several proofs
of ``$\Omega_3=0$" and of the equivalent {\it Lickorish-Wallace theorem}
about $3$-manifolds up to `londitudinal' Dehn surgery equivalence respectively. Every
proof will illuminate different facets of the matter. We determine the bordism 
semigroup (which turns out to be a group) of immersions of surfaces into
a given compact connected boundaryless $3$-manifold $M$. If $M$ is orientable,
a key ingredient will be the Arf-Brown invariant of the quadratic enhancement
of the intersection form associated to every immersion of a surface in $M$
(endowed with an auxiliary framing).  We also classify
compact boundaryless $3$-manifolds up to certain equivalence relations generated by
diffeomorphisms and blow-up-down along smooth centres (a notion introduced
in Chapter \ref{TD-CUT-PASTE}). The subtler so called `tear' equivalence, in the
non orientable case also involves instances of quadratic enhancement of the intersection form 
of  charcteristic surface, i.e. representing  the Euler class of the determinant bundle of 
the ambient $3$-manifold. We also discuss an application to a solution of the so called
{\it Nash rationality question} in dimension $3$.

In Chapter \ref{TD-4}, by analogy to the case of  surfaces, we focus on (the isometry class of) the intersection form
$\sqcup_M: \Hh^2(M;\Z)\times \Hh^2(M;\Z) \to \Z$ as the main invariant of every
compact connected oriented boundaryless $4$-manifold $M$. It turns out that this is a symmetric
unimodular $\Z$-bilinear form on the finite rank free $\Z$-module $\Hh^2(M;\Z)$.
We proof Rohlin's theorem that the signature $\sigma$ of the intersection form
determines an isomorphism $\sigma: \Omega_4 \to \Z$, so that $\Omega_4$ is generated
by $[\PP^2(\C)]$. In fact we follows his original
geometric proof.  Trying to pursuing the analogy with surfaces, first we face the problem of the abstract
arithmetic classification of such symmetric unimodular forms. 
A first main difference is that it is complete only in the {\it indefinite} case. Then we try
to develop as much as possible a parallel 4D counterpart at least in the indefinite case,
by restricting in fact to {\it simply connected} $4$-manifolds. We establish a classification
up to {\it odd stabilizations}, the elementary ones being the connected sum with $\pm \PP^2(\C)$.
We just outline a more subtle classification up to {\it even stabilization} i.e. up to connected
sum with $S^2\times S^2$. The arithmetic tells us that there are {\it characteristic elements} 
$\beta \in \Hh^2(M;\Z)$ such that for every $\alpha \in \Hh^2(M;\Z)$, 
$\alpha \sqcup \alpha = \beta \sqcup \alpha$ mod$(2)$, and that $\sigma = \beta\sqcup \beta$ mod$(8)$.
Every $\beta$ can be represented by an oriented surface $F$ embedded in $M$, called a {\it characteristic surface}.
We prove the congruence 
$ \sigma - \beta \sqcup \beta = 8\alpha (F)$ mod $(16)$, where $\alpha (F) \in \Z/2\Z$
is the Arf invariant of a  quadratic enhancement of the intersection form of $F$ which
represents an obstruction to surgery $F$ within $M$ to an embedded $2$-sphere.
If the intersection form is even, we can take $F=\emptyset$, so that we recover the
original celebrated Rohlin congruence $\sigma = 0$ mod$(16)$ (early obtained as
a corollary of the fact that $\pi_{n+3}(S^n)=\Z/24\Z$ for $n$ big enough).
This implies in particular
that there are unimodular symmetric forms which cannot be realized as the intersection form
of any simply connected $4$-manifold. We propose an elementary proof due to \cite{Mat}
and based on the classification up to odd stabilizations. We end the chapter with an 
informative and discorsive section about more recent achievements in the realm of
$4$-manifolds.   

Overall this text is a collection of themes, in some case advanced and of historical importance, 
with the common feature that they can be treated with ``bare hands''. 
This means by just combining specific differential-topogical cut-and- paste procedures and applications of 
transversality, mainly through the cobordism multiplicative structure. Tools widely developed along these
pages. Of course the choice of the themes is also matter of personal preference. 
It is aimed to be accessible and useful to motivated and collaborative
master undergraduate students, to PhD students and also to a more expert reader in order to 
recognize very basic reasons for some facts already known to her/him as resulting from more advanced 
theories and/or technologies.

\chapter{The smooth category of open subsets of euclidean spaces}\label{TD-LOCAL}
We will be concerned with manifolds. Roughly, a manifold is a 
topological space locally modeled
on some euclidean space $\R^n$, $n\in \N$. So let us recall a few facts about
our favourite local models. Many of them  should be familiar to the readers,
so sometimes we will omit the proofs or just sketch them.

\section{Basic structures on $\R^n$}\label{basic-structure}
Every space $\R^n$, $n\in \N$, is endowed with a variety of structures
that case by case will be involved in the discussion.

$\R^n$ is the {\it vector space}
of colums vectors (with $n$ rows). {\it We stipulate that 
if $x\in \R^n$ occurs as a vector
in any linear algebra formula then it is considered as a column}.  

The space $\Ll(\R^n,\R^m)$ of linear maps $L: \R^n \to \R^m$ coincides
with the space of matrices $m \times n$, $M(m,n,\R)$, so that for
every $x \in \R^n$, $x\to Lx$ via the usual ``lines by column''
product. By using the lexicographic order on the entries of any matrix
$L=(l_{i,j})_{i=1,\dots, m; j=1,\dots, n}$, we fix also the
identification of $M(m,n,\R)$ with $\R^{mn}$. As every vector space,
$\R^n$ has a canonical {\it affine space} structure determined by the
map that associates to every couple of {\it points} $(x,y)\in \R^n
\times \R^n$ the {\it vector} $\overrightarrow {xy}:= y-x$. Every {\it
  affine map} $f: \R^n \to \R^m$ is of the form $f(x)= w+ Lx$ where
$w\in \R^m$ and $L \in \Ll(\R^n,\R^m)$.
 
$\R^n$ is a {\it complete metric space} endowed with the
{\it euclidean distance}
$d=d_n$ defined by
$$d(x,y)= \sqrt{\sum_{j=1}^n (x_j-y_j)^2} \ . $$

The {\it standard positive definite scalar product}
$(*,*)=(*,*)_n$ on $\R^n$ is defined by
$$ (x,y):= \sum_{j=1}^n x_jy_j= x^tIy$$
with the associated {\it norm} $||x||= \sqrt{(x,x)}$. We note that
$$d^2(x,y)=(x-y,x-y)$$
and that the familiar formula
$$ (x,y) = ||x||\cdot||y||\cos\theta$$
allows to recover the measure of the angle formed by the ordered and oriented lines
spanned by two non zero vectors $x$, $y$; in particular they are othogonal
iff $(x,y)=0$. Hence many basic objects of elementary geometry
can be expressed analytically by means of the standard scalar
product.

$\R^n$ is a {\it topological space} endowed with the topology
$\tau=\tau_n$ induced by the distance $d_n$. As for any metrizable
topological space, a subset $U$ of $\R^n$ is {\it open} if and only if
for every $x\in U$, there is $r>0$ such that the {\it ``open''
  $n$-ball of center $x$ and radius $r$}
$$ B^n(x,r):=\{y\in \R^n; \ d(x,y)<r \}$$
is contained in $U$. We will denote by 
$$D^n=\overline B^n(0,1)$$ 
the {\it closed unitary $n$-ball}
also called the unitary {\it $n$-disk},
and by 
$$S^{n-1}=\partial D^n = \{x\in \R^n; \ d(0,x)=||x||=1\}$$ 
the {\it unitary
sphere}. One verifies that the ``open'' balls are indeed open sets and the open
balls with center in $\Q^n\subset \R^n$ 
and rational radius form a {\it countable basis of open sets of $\tau$}
(every open set is union of such balls). Any other
scalar product $$(x,y)_A:= x^tAy$$ 
defined by a positive definite
symmetric matrix $A=A^t$, determines (by the same formulas as above)
a norm $||.||_A$, a distance $d_A$ and an associated topology $\tau_A$.
In fact all these distances are {\it topologically equivalent}, that is every $\tau_A = \tau$.
This can be proved by means of the version of elementary {\it spectral theorem}
stating that there exists a basis of $\R^n$ which is simultaneously orthonormal
for $(*,*)$ and orthogonal for $(*,*)_A$.  Another
topologically equivalent distance on $\R^n$ is defined by
$$\delta(x,y):= \max \{ |x_j-y_j|; \   j=1,\dots,n \} \ . $$

Accordingly to general topological definitions, for every $X\subset \R^n$, 
$$\tau \cap X=\{U\cap X; \ U\in \tau\} $$
is the  topology on $X$ that makes it a {\it topological subspace} of $(\R^n,\tau)$;
given subspaces $X\subset \R^n$, $Y\subset \R^m$, a map
$f:X\to Y$ is {\it continuous} if for every open set $U\subset Y$, the {\it inverse
image} $f^{-1}(U)=\{x\in X; \ f(x)\in U \}$ is an open set of $X$.
A continuous map $f:X\to Y$ is a {\it homeomorphism} if it is bijective and 
also the {\it inverse map} $f^{-1}: Y\to X$ is continuous.
 
Every subspace $X\subset \R^n$ is metrizable (hence in particular {\it
  Hausdorff}) by the restriction to $X$ of the distance $d$ (or of any distance
topologically equivalent to $d$); the restriction of any
(countable) basis of open sets of $\tau$ is a (countable) basis of
$\tau \cap X$.

As for every Hausdorff space with a countable basis, a subspace $X$ of
$\R^n$ is {\it compact} (i.e. every open covering of $X$ admits a {\it
  finite} sub-covering) if and only if it is {\it sequentially compact}
(i.e. every sequence $a_n$ of points of $X$ admits a sub-sequence $a_{j_n}$
converging to some point $x$ of $X$).  A subspace is compact if and
only if it is {\it closed} (i.e. the complementary is open) and {\it
  bounded} (i.e. it is contained in some ball $B^n(0,r)$). $\R^n$ is {\it
  locally compact} (for every $x\in \R^n$ the family of {\it closed
  balls} $\overline B^n(x,r)=\{y\in \R^n; \ d(x,y)\leq r\}$, when
$r>0$ varies, is a basis of compact neighbourhoos of $x$). The same
holds for every subspace $X$ which is a closed subset of $\R^n$.

We have

\begin{proposition}\label{connected} A non empty open subset $U\subset \R^n$ is 
{\rm connected} (i.e. $U$ is the only open-and-closed non empty subset of $U$) 
if and only if it is {\rm path connected} (i.e. for every two points $x_0$, $x_1$ of $U$, there
is a continuous path $\alpha: [0,1] \to U$ such that $\alpha(0)=x_0$,
$\alpha(1)=x_1$).  
\end{proposition}
\Dim The ``if'' implication holds in general for arbitrary topological spaces and is due to the
basic fact that intervals in the real line are connected; for ``only
if'', note that ``being connected by a continuous path'' defines
an equivalence relation on $U$. The equivalence classes are
called the {\it path connected components} of $U$. As every open ball 
$B^n(x,r)\subset U$ is contained in the path
connected component of $U$ which contains $x\in U$, then every path
connected component of $U$ is open, hence there is only one if $U$ is
connected.

\cvd

\medskip
\section{Differential calculus}\label{diff-calc}
Another fundamental structure carried by the spaces $\R^n$
is the {\it differential calculus}. 
Let $U\subset \R^n$, $W\subset \R^m$ be open sets. A map 
$$f=(f_1,\dots, f_m): U \to W$$ is said to be $\Cc^0$ if it is
continuos.  The map is {\it differentiable} at $x\in U$ if there is a
(necessarily unique) linear map $d_xf \in \Ll(\R^n,\R^m)$ that
``well'' approximates $g(h)=f(x+h)-f(x)$ in a neighbourhood of
$h=0$. Precisely, for every $\epsilon >0$, there is $\delta>0$ such
that for every $h$ such that $||h||<\delta$, $x+h \in U$ and
$$||g(h)-d_xf(h)|| \leq \epsilon||h|| \ .$$ The linear map $d_xf$ is
called the {\it differential of $f$ at $x$}. The map $f$ is (globally)
{\it differentiable} if it is differentiable at every point $x\in
U$. In such a case it is defined the {\it differential map}
$$df: U\to M(m,n,\R), \  df(x):=d_xf \ . $$
We say that $f$ is $\Cc^1$ if it is differentiable and $df$ is continuous 
 (being $\Ll(\R^n,\R^m)=M(m,n,\R)$ confused with $\R^{mn}$
as above). Every $\Cc^1$ map is $\Cc^0$. By induction, for every $r\geq 1$,
we say that $f$ is $\Cc^r$ if $df$ is $\Cc^{r-1}$.
In practice, $f$ is $\Cc^r$, $r\geq 1$, if and only
if it is $\Cc^0$  and for every multi-index $J=j_1\dots j_n$ of {\it order}
$|J|:= j_1+\dots + j_m \leq r$, for every $i=1,\dots, m$,
it is defined and is continuous
the {\it partial derivative function}
$$ \frac{\partial ^J f_i}{\partial^{j_1}x_1\dots
  \partial^{j_n}x_n}:U\to \R \ . $$
Then for every $x\in U$, the partial derivatives of the first order can be organized
in a $m\times n$ matrix so that
$$ d_xf := \left (\frac {\partial f_i}{\partial x_j}(x)\right )_{i=1,\dots, m; \ j=1,\dots , n}\in M(m,n,\R) \ . $$
This is a consequence of the ``chain rule'' (see below).

A map $f$ is $\Cc^\infty$ or, equivalently, {\it smooth} if it is $\Cc^r$ for every $r\geq 0$.
If $f$ is smooth, then also $df$ is smooth. So we can define inductively for every $r\geq 1$,
$d^rf= d(d^{r-1}f)$.

If $f$ is (at least) $\Cc^1$ we have the following {\it uniform} version of the above property
that define the differentials $d_xf$: for every $x\in U$ there exists a neighbourhood $W$ of $x$
in $U$ (we can take as $W$ a compact closed ball $\bar B^n(x,\rho)\subset U$), such that
for every $\epsilon >0$, there is $\delta>0$ such that for every $y\in W$ and for every $h$, $||h||<\delta$ we have 
$y+h\in U$ and
$$ || f(y+h)-f(y) - d_yf(h))|| \leq \epsilon \delta \ ; $$
in other words
$$\lim_{h\to 0} \frac{g(y,h)-d_yf(h)}{||h||}=0 \ . $$
uniformly with respect to $y\in W$. 

\medskip

{\it From now on we will be mainly concerned with smooth maps.}

\medskip

{\bf (Taylor polynomials.)}  
A {\it homogeneus polynomial maps of degree $k\geq 1$} $$\pG: \R^n \to \R^m$$
is by definition of the form $\pG(x)=\phi(x,\dots,x)$, where $\phi: (\R^n)^k\to \R^m$ is a (necessarily unique)
{\it symmetric} $k$-linear map ($\phi$ is called the ``polarization'' of $\pG$).  
It follows that the set $\Pp_k(n,m)$ of these homogeneus polynomial maps 
has a natural structure
of finite dimensional real vector space. 
A {\it polynomial map of degree} $d$, $p:\R^n \to \R^m$, is of the form 
$$p= p_0+p_1+\dots + p_d$$
where $p_0 \in \R^m$ and for $j\geq 1$, $p_j$ is homogeneous polynomial of degree $j$ and $p_d$ is not zero. 

Let $f: \R^n \to \R^m$ be a smooth map. 
Then for every $k\geq 1$ there is a smooth map 
$$T_k(f): U\to \Pp_k(n,m)$$
such that for every $k\geq 1$, for every $x\in U$, there is a neighbourhood $W$ of $x$ in $U$ such that
for every $\epsilon >0$, there is $\delta >0$ such for every $y\in W$ and every $h$, $||h||<\delta$,
we have $y+h \in U$ and
$$ ||f(y+h)-(f(y)+T_1(f)(y)(h)+\dots + T_k(f)(y)(h))|| \leq \epsilon ||h||^k \ . $$
The maps $T_k(f)$ are uniquely determined by these conditions.  Clearly 
$$T_1(f)(x)=d_xf \ . $$
More generally, every $\pG=(\pG_1,\dots , \pG_m) \in \Pp_k(n,m)$ is of the form
$$ \pG_i(h)= \sum_{|J|=k} a_i^J h^{j_1}_1\cdots h^{j_n}_n$$
where the coefficients $a^J_i \in \R$.  Then one verifies that $T_k(f)(x)$ is uniquely determined
by the formulas
$$ a^J_i =  \frac{1}{k!}  \frac{ \partial^J f_i}{\partial ^{j_1}x_1\dots \partial^{j_n}x_n}(x) \ . $$
In other words, $T_k(f)(x)$ is determined by means of $\frac{1}{k!}d^k_x(f)$.
$T_k(f)(x)$ is the {\it homogeneous degree-$k$ Taylor polynomial of $f$ at $x$}. 
Setting $f(x)=T_0(f)(x)$, the polynomial map (of the variable $h$) 
$$\Tt_k(f)(x):= T_0(f)(x) + T_1(f)(x) + \dots + T_k(f)(x)$$
is called {\it Taylor polynomial} of $f$ at $x$ of degree $\leq k$.

\medskip

\section{The  category of open subsets of euclidean spaces and smooth maps}\label{open-set-cat}
Let $f:U\to W$, $g:U' \to W'$ be smooth maps between open subsets of some (possibly variable)
euclidean spaces. The composition $g \circ f$ is defined when $W\subset U'$.
The fundamental well known {\it chain rule} for the composition of differentiable maps 
states that for every $x\in U$, $y=f(x)$, $g\circ f$ is differentiable at $x$ and
$$d_x(g\circ f)= d_yg \circ d_x f \ .$$
It follows immediately that if $f$ and $g$ are smooth then also $g\circ f$ is smooth.
Then we can consider the category whose {\it objects} are the open subsets of euclidean
spaces and for every couple $(U,W)$ of objects, the {\it ``arrows''} (that is the {\it morphisms})
are the smooth
maps $\Cc^\infty(U,W)$.  

For every object $U\subset \R^n$, the unit map $1_U$ is the {\it identity}
$$ {\rm id}_U: U\to U, \ {\rm id}_U(x)=x \ $$
which is obviously smooth. For every $x\in U$, 
$$d_x{\rm id}_U=  {\rm id}_{ \R^n} = I_n \in {\rm End}(\R^n)=M(n,\R) \ . $$

If $U'\subset U$ then the inclusion $i: U' \to U$ is smooth and for
every $f\in \Cc^\infty(U,W)$, the {\it restriction} $f|_{U'}= f\circ i$ is smooth.

The {\it equivalences} in this category are the {\it diffeomorphisms}. 
Let $U\subset \R^n$ and $W\subset \R^m$
be open sets. Then $f\in C^\infty(U,W)$ is a diffeomorphism 
if it is a homeomorphism and 
also the inverse map $f^{-1}: W\to U$ is smooth. In such a case, 
by applying again the chain rule, we have that 
for every $x\in U$, $y=f(x)$, 
$d_y f^{-1}\circ d_xf = I_n$, $d_xf \circ d_y f^{-1} = I_m$,
then by elementary linear algebra 
both inequalities $n\leq m$ and $m\leq n$ hold, so that $m=n$; finally
$d_xf \in {\rm GL}(n,\R)$ is invertible and
$$ d_y f^{-1} = (d_x f)^{-1} \ . $$ 
Hence we have proved the {\it invariance of the dimension up to diffeomorphism}
and this is reduced to the basic invariance of dimension up to linear isomorphism.

Another consequence of these considerations (based on the chain rule):
\smallskip

{\it If a smooth homeomorphism $f:U\to W$ has differentiable inverse
map $f^{-1}$ then it is a diffeomorphism (i.e. $f^{-1}$ is smooth indeed).}
\medskip

\section{The chain rule and the tangent functor}\label{tang-fuct}
The chain rule can be rephrased in the language of {\it functors}
between categories. A way is to consider the category of {\it pointed} open subsets
of some euclidean spaces and pointed smooth maps. Then by setting
$$ (U,x),\  U\subset \R^n \Longrightarrow \R^n$$
$$ f\in \Cc^\infty ((U,x),(W,y)), \ U\subset \R^n, \ W\subset \R^m \Longrightarrow d_xf \in \Ll(\R^n,\R^m) $$
we define a {\it covariant} functor from the smooth pointed category to the category
of finite dimensional real vector spaces and linear maps.

Avoiding to deal with the pointed category, another way is by defining the so called 
{\it tangent functor} which is a covariant
functor from our favourite category to itself. Set
$$ U\subset \R^n  \Longrightarrow  T(U):=U\times \R^n\subset \R^n \times \R^n$$
$$ f\in \Cc^\infty(U,W) \Longrightarrow  Tf\in \Cc^\infty(T(U),T(W)),  \ Tf(x,v):= (f(x), d_xf(v)) \ . $$
The chain rule can be rewritten as
$$ T(g \circ f)=Tg \circ Tf$$
$$ T{\rm id}_U = {\rm id}_{T(U)}$$
if $f\in \Cc^\infty(U,W)$ is a diffeomorphism, then also $Tf$ is a diffeomorphism.

There is a natural projection $$\pi_U:T(U)\to U,  \ \pi_U(x,v)=x \ . $$
$(T(U),\pi_U)$ is called the {\it tangent bundle of $U$}. For
every $x\in U$, the {\it fibre} $$T_xU:= \pi_U^{-1}(x)$$ is naturally identified with the
{\it vector space} $\R^n$ and is called the {\it tangent spaces to $U$ at $x$}.
Every $v\in T_xU$ is a tangent vector at $x$. This notion of $T_xU$ is essentially the
one we get by considering $U$ as an open set in the {\it affine} space $\R^n$. 

The map $Tf$ is called the {\it tangent map} of $f$. Clearly, 
$$\pi_W\circ Tf= f \circ \pi_U$$
that is $Tf$ sends every fibre $T_xU$ to the fibre $T_{f(x)}W$ by means of the linear map $d_xf$
which varies smoothly when $x$ varies in $U$.

\section{Tangent vector fields, riemannian metrics, gradient fields}\label{fields} 
A {\it tangent vector fields} on $U$ (often we will omit to say ``tangent'') is a smooth map of the form
$$ X: U\to T(U),\ X(x)=(x,v_X(x)) $$
so that $\pi_U\circ X = {\rm id}_U$. Such a map is also called a (smooth) {\it section} of the tangent
bundle. Hence $X$ selects a family (a ``field'') of vectors $\{v_X(x)\in T_xU\}_{x\in U}$ which vary smoothly
with the point $x\in U$. In practice $X$ could be confused with the smooth map $v_X: U\to \R^n$; however,
if $\phi:U \to W$ is a diffeomorphim, as a map $v=v_X$ is transported on $W$ by the composition $v \circ \phi^{-1}$,
while  the vector field $X$ is transported to $\phi_*X$ on $W$ by the composition $T\phi\circ X$, that is for every $y=\phi(x)\in W$
$$ \phi_*X(y) = (y, d_x\phi(v_X(x)) )\ . $$

Denote by $\Gamma(T(U))$ the set of  vector fields on $U$. For every $X,Y \in \Gamma(T(U))$,
every $f\in \Cc^\infty(U,\R)$, and every $x\in U$, set
$$X+Y(x)=(x, v_X(x)+v_Y(x)), \ fX(x)=(x,f(x)v_X(x)) \ . $$
This defines on $\Gamma(T(U))$ a natural structure of {\it module} over the commutative
ring $\Cc^\infty(U,\R)$ which induces (by restriction to the constant functions) a structure of
$\R$-vector space. Let us denote by $\eG_i(x)=(x,e_i)$, $i=1,\dots,n$, 
the {\it constant vector field} on $U$ such that
$e_i=(0,0,\dots, 1,\dots , 0)^t$ is the $i$th-vector of the canonical basis of $\R^n$. Sometimes
$\eG_i$ is also denoted by $\frac{\partial}{\partial x_i}$. Then for every $X\in \Gamma(T(U))$,
$$X= \sum_i v_{X,i} \eG_i$$ 
that is the fields $\eG_i$ form a basis of such a module.

A riemannian metric on $U\subset \R^n$, is a smooth map
$$g: U \to M(n,\R)$$
such that for every $x\in U$, the matrix $g(x)$ is symmetric and positive definite.
Then $\{g(x)\}_{x\in U}$ is a smooth {\it fields of positive definite scalar products} 
$(*,*)_{g(x)}$
defined on each tangent space $T_xU$. Denote by $S(n,\R)$
the space of symmetric $n\times n$ matrices ($S(n,\R)$ can be identified
with $\R^{\frac{n(n+1)}{2}}$). By setting 
$$U\to U\times S(n,\R), \ x\to (x,g(x))$$ 
then the riemannian metric can be interpreted as a section 
of the ``product bundle'' $U\times S(n,\R) \to U$.

If $g$ is a riemannian metric on $U$
and $X,Y\in \Gamma(T(U))$, then 
$$x\to (v_X(x),v_Y(x))_{g(x)}$$
defines a smooth function  on $U$. 

If $g_0$ and $g_1$ are riemannian metrics on $U$, then $g_t= (1-t)g_0 + tg_1$, $t\in [0,1]$,
is a path of riemannian metrics.

An {\it isometry} $\phi: (U,g)\to (W,h)$ ($g,h$ being riemannian metrics) 
is by definition a diffeomorphism such that for every $x\in U$, every $v,w \in T_xU$,
$$ (v,w)_{g(x)} = (d_x\phi(v),d_x\phi(w))_{h(\phi(x))} \ . $$

Given $(U,g)$ and a diffeomorphism $\phi : U\to W$, this transports $g$ 
to the riemannian metric
$\phi_* g$ on $W$ such that $\phi$ is {\it tautologically} an isometry.
If $y\in W$, set $P(y)= d_y\phi^{-1}$, then
$$ \phi_*g: W \to M(n,\R), \ y \to  P(y)^tg(\phi^{-1}(y))P(y) \ . $$

If $f\in \Cc^\infty(U,\R)$, its differential function $df: U \to M(1,n)$
can be considered as a smooth {\it field of linear functionals} 
$\{d_xf: T_xU \to \R\}_{x\in U}$, $d_xf$ belonging to the dual 
space $T^*_xU$;
in other words, it is identified with
the section $x\to (x,d_xf)$ of the {\it cotangent bundle} 
$$\pi^*: T^*(U)=U\times M(1,n)\to U \ . $$
Every such a section $\Omega(x)=(x,\omega(x))$ is called a 
{\it differential form} on $U$.
For every form $\Omega$  and every vector field $X$ on $U$,
$$x\to \omega(x)(v(x)) $$ defines a smooth function on $U$. 
If $\phi: U\to W$ is a diffeomorphism, $\Omega$ a differential
form on $U$, then $\phi$ transports $\Omega$ to the form
$\phi_*\Omega$ on $W$ such that for every $y\in W$, every
$w\in T_yW$, then 
$$\phi_*\Omega(y)=(y,\omega(d_y\phi^{-1}(w))) \ . $$
Denote by $\eG^j$, $j=1,\dots, n$ the field constantly equal to the
functional $e^j$ such that
$$ (e^i(e_j))_{i,j} = I_n \in M(n.\R) \ . $$
Then every $\Omega \in \Gamma(T^*(U))$ is a unique linear combination
$$ \Omega = \sum_j a_j \eG^j$$
the $a_j$ being smooth functions on $U$.
Sometimes one writes $\partial x_j$ instead of $\eG^j$.

If $g$ is a riemannian metric on $U$, then by setting
for every $v, w \in T_xU$, 
$$\psi_v(w):=g(x)(v,w)\in \R$$
one defines a smooth {\it field of linear isomorphisms}
$\Psi_g:=\{ \Psi_{g,x}: T_x(U)\to T^*_x(U)\}_{x\in U}$. This transforms 
vector fields into differential forms. For every $f\in \Cc^\infty(U,\R)$,
let $\nabla_g f$ be the unique vector field on $U$ such that $\Psi_g(\nabla_gf) = df$,
so that for every $x\in U$, $v\in T_x(U)$, then 
$$d_xf(v)=(\nabla_g(x),v)_{g(x)} \ . $$
The field $\nabla_gf$ is called the {\it gradient of $f$ with respect to the metric $g$}.
Clearly for every $x\in U$, $d_xf(\nabla_g(x))=(\nabla_g(x),\nabla_g(x))_{g(x)} \geq 0$
and is strictly positive if and only if $d_xf\neq 0$.

Obviously every $U$ admits riemannian metrics, for example any constant one $g_A(x)=A$
where $A$ is a symmetric positive definite matrix. In particular $g_0:=g_{I_n}$ is called the {\it standard
riemannian metric}. We have that for every smooth function $f$ on $U$, 
$$\nabla_{g_0}f(x) = d_xf^t \ . $$
  
\section{Inverse function theorem and applications}\label{inverse-funct}
Let $L\in \Ll(\R^n,\R^m)$ be a linear map of {\it maximal rank} $r$. There are a few
possibilities and by elementary linear algebra, for every case there is a {\it normal form}
up to pre or post composition with linear isomorphisms.
\begin{itemize}
\item If $r=n=m$, then $L\in GL(n,\R)$ is invertible and the normal form is $I_n$ abtained
as $$I_n= L\circ L^{-1} = L^{-1} \circ L \ . $$

\item If $n<m$, then the rank $r$ is equal to $n$ and $L$ is {\it injective}.
Let us fix a direct sum decomposition
$$ \R^m= L(\R^n)\oplus V$$
 and a basis $\Bb = \Bb' \oplus \Bb"$ of $\R^m$ adapted to the decomposition.
This determines a linear isomorphism $\phi_\Bb: \R^m\to \R^n\times \R^{m-n}$
such that for every $x=(x_1,\dots, x_n)^t \in \R^n$, we have
$$\phi_\Bb \circ L(x)= (x_1,\dots, x_n, 0,\dots, 0)^t$$
that is the standard inclusion $j=j_{n,m}: \R^n \to \R^m=\R^n\times \R^{m-n}$. This is the normal
form in this case.

\item If $n>m$, then the rank $r$ is equal to $m$ and $L$ is {\it surjective}.
Fix a direct sum decomposition
$$ \R^n = V \oplus \ker(L) \ $$
and an adapted basis $\Bb= \Bb' \oplus \Bb"$ of $\R^n$.
This determines a linear isomorphism (in fact the inverse of the above defined $\phi_\Bb$)
$\psi_\Bb: \R^m\times \R^{n-m} \to \R^n$ such that for every $x=(x_1,\dots, x_n)^t \in \R^m\times \R^{n-m}$,
we have
$$ L\circ \psi_\Bb(x)=(x_1,\dots, x_m)$$
that is the natural projection $\pi_{n,m}:\R^m \times \R^{n-m} \to \R^m$.
This is the normal form in this case.
\end{itemize} 

Let us consider now a morphism $f\in \Cc^\infty(U,W)$ in our favourite category,
$U\subset \R^n$, $W\subset \R^m$, $p\in U$. Assume that $d_pf$ has maximal rank
$r$. The following fundamental theorems state that locally in a neighbourhood of $p$
in $U$, the map $f$ takes the same normal form of the linear map $d_pf$, up to pre or
post composition with smooth diffeomorphisms. As a first step, let us remark that
the {\it punctual} hypothesis has in fact a {\it local} valence. By a well known criterion
$d_pf$ has maximal rank $r$ if and only if there is a $r\times r$ submatrix $A(p)$
of $d_pf$ such that $\det A(p) \neq 0$. By taking the same submatrix $A(x)$ of
$d_xf$ for every $x\in U$, we define the smooth function
$$ \det A: U\to \R, \  x \to \det A(x) \ . $$
Then by the ``sign permanence'', there exists an open neighbourhood $U'\subset U$ of $p$
in $U$, such that for every $x\in U'$, $d_xf$ has maximal rank $r$.

A map $f\in \Cc^\infty(U,W)$ such that for every $x\in U$, $d_xf$ is injective is called an
{\it immersion}. If $d_xf$ is surjective for every $x\in U$, then $f$ is called a {\it summersion}.
If $n=m$ the two notions coincide.
We can state now the theorem mentioned in the title.

\begin{theorem}\label{inverse-function} {\rm (Inverse function theorem)} 
Let $f\in \Cc^\infty(U,W)$, $U, W \subset \R^n$, such that
for every $p\in U$, the differential $d_pf$ is invertible.
Then $f$ is a {\rm local diffeomorphism}, that is for every
$p\in U$ there is a open neighbourhood $U'$ of $p$ in $U$
such that $W'=f(U')$ is an open subset of $W$ and the restriction
$f|_{U'} \in \Cc^\infty(U',W')$ is a diffeomorphism.
\end{theorem}

\begin{corollary}\label{immersion}{(\rm (Local immersion theorem)}
Let $f\in \Cc^\infty(U,W)$, $U\subset \R^n$,
$W\subset \R^m$, $n<m$, be an immersion. Then for every $p\in U$
there exist 
\begin{itemize}
\item An open neighbourhood $U'$ of $0$ in $\R^n$; 
\item an open neighbourhood
$W'$ of $q=f(p)$ in $W$; 
\item an open neighbourhood $W"$ of $0$ in $\R^m$ and a diffeomorphism 
$$\phi: (W',q)\to (W",0)$$ 
\end{itemize}

\noindent such that  for every $x\in U'$, $x+p \in U$, 
$f(x+p)\in W'$ and
$$ \phi \circ f (x+p)=j_{n,m}(x) \ . $$
\end{corollary}

\begin{corollary}\label{summersion}{\rm (Local summersion theorem)}
 Let $f\in \Cc^\infty(U,W)$, $U\subset \R^n$,
$W\subset \R^m$, $n>m$, be a summersion. Then for every $p\in U$ there exist
\begin{itemize}
\item An open neighbourhood $U'$ of $p$ in $U$; 
\item an open neighbourhood $U"$ of 
$0$ in $\R^n$ and a diffeomorphism $$\psi: (U",0)\to (U',p)$$ 
\end{itemize}
\noindent such that $f(U')\subset W$
and $$f\circ \psi (x)-f(p) = \pi_{n,m}(x) \ . $$
\end{corollary}

\medskip

{\it Proof of the Corollaries.}  In both cases it is not restrictive
to assume that $p=0$ and $f(0)=0$.  We will use the notations
introduced at the beginning of the section, by replacing $L$ with
$d_0f$.

{\it (Immersions)} Given a direct sum decomposition of $\R^m=
d_0f(\R^n)\oplus V$, with adapted basis $\Bb=\Bb'\oplus \Bb"$ and
associated linear isomorphism
$$\psi_\Bb: \R^n\times \R^{m-n} \to \R^m$$ consider the smooth map
$$ g: U\times \R^{m-n}\to \R^m, \ g(x,h)=f(x)+\psi_\Bb(0,h) \ . $$
It is easy to verify that $d_{(0,0)}g$ is invertible and we can apply the inverse function theorem to $g$ on a neighbourhood
$U'\times A$ of $(0,0)$. By construction, for every $x\in U'$,  $f(x)= g\circ j_{n,m}(x)$, so that
$g^{-1}\circ f (x)=j_{n,m}(x)$.

\medskip

{\it (Summersions)} Given a direct sum decomposition $\R^n= V\oplus
\ker (d_0f)$ with adapted basis $\Bb=\Bb'\oplus \Bb"$ and associated
linear isomorphism
$$\phi_\Bb: \R^n \to \R^m \times \R^{n-m}$$  set $p:\R^m \times \R^{n-m} \to \R^{n-m}$
the natural projection.
Define
$$g:U \to \R^m \times \R^{n-m}, \  g(x)=(f(x), p(\phi_\Bb (x)) \ . $$
One verifies that $d_0g$ is invertible, so we can apply the inverse function theorem to $g$ on a 
neighbourhood $U'$ of $0$. By construction, for every $x\in U'$,  $f(x)=\pi_{n,m}\circ g(x)$,
and we conclude similarly to the case of immersions above.

\cvd

\medskip

Corollaries \ref{immersion} and 
\ref{summersion} are  instances of the following general {\it constant rank theorem}. 
The proof, based again on the Inverse Function Theorem, is a simple variation  and is left
as an exercise.

\begin{theorem}\label{constant-rank} {\rm (Constant rank theorem)}
Let $f:U\to W$ be a smooth map, $U\subset \R^n$, $W\subset \R^m$ be 
open sets. Assume that $d_xf$ is of constant rank $k\leq \min \{n,m\}$.
Then for every $p\in U$, $q=f(p)$ up to pre and post composition with
local diffeomorphisms $\psi: U'\to U$, $\psi(0)=p$, $\phi: W'\to W$, $\phi(0)=p$
we have that 
$$\rho:=\phi^{-1} \circ f \circ \psi: U' \to W',\ \rho(u_1,\dots,u_n)=(u_1,\dots, u_k)  \ . $$
\end{theorem}

\cvd







Strictly related to the local summersion theorem there is another
corollary known as the {\it implicit function theorem}.
A consequence of the proof of Corollary \ref{summersion} is that 
there is a diffeomorphism $\rho: A\times B \to U'$, where
$A\times B\subset \R^{n-m}\times \R^m$ is an open neighbourhood of 
$(x_0,y_0)=(0,0)$ and $U'$ is an open neigbourhood of $0$ in $U\subset \R^n$,
such that restriction of $g=f\circ \rho$ to $A\times B$ verifies:
\begin{enumerate}
\item $g(x_0,y_0)=0$;
\item The restriction $\tilde g$ of $g$ to $\R^m=\{x_0\}\times \R^m$ 
has invertible differential $d_{y_0}\tilde g$ at $y_0$
\end{enumerate}

We take such a situation as the hypotheses of the implicit function theorem.
 
\begin{corollary}\label{implicit-funct}{\rm (Implicit function theorem)}
Let $A\times B \subset \R^k\times \R^m$ be an open set.
Let $g: A\times B\to \R^m$ be a smooth map and $(x_0,y_0)\in A\times B$
such that $g(x_0,y_0)=0$. Let $\tilde g$ be the restriction of $g$ to $\R^m = \{x_0\}\times \R^m$.
Assume that $d_{y_0} \tilde g$ is invertible. Then there exist an open neighbourhood $A'\times B'$
of $(x_0,y_0)$ in $A\times B$, and a smooth maps $h: A'\to B'$ such that 
$${\rm Graph}(h)= f^{-1}(0)\cap A'\times B' \ . $$
\end{corollary}
\medskip

It follows that $f(x,h(x))=0$ for every $x\in U'$ and $h$ is said to be (locally)
{\it  implicitly defined by the equation} $f(x,y)=0$.

\medskip

{\it Sketch of proof.}
We use similar arguments as in the proofs of the previous
corollaries.
Consider the smooth map 
$$ G: A\times B \to \R^k \times \R^m, \ G(x,y)=(x,g(x,y)) \ . $$
It is immediate to check that $d_{(x_0,y_0)}G$ is invertible,
so we can apply the inverse function theorem to $G$ in
a neighbourhood $A_1\times B'$ of $(x_0,y_0)$, and
the inverse  map is necessarily of the form
$$G^{-1}(x,y)=(x,l(x,y)) $$
for a suitable smooth map $$l: G(A_1\times B')\to B' \ . $$
Take $$A'=\{x\in U; \ (x,0)\in G(U_1\times W'\}$$ and define
$h(x)=\l(x,0)$. The reader can complete by exercise the verification that
$A'\subset A_1$ and this eventually achieves the proof.

\cvd
 
\medskip

A proof of the inverse function theorem should be known to the reader.
A current conceptual proof is based on {\it Banach's principle} for contractions 
on complete metric spaces. This is suited for generalizations to
infinite dimensional Banach spaces. However we just sketch one
 in our finite dimensional situation, based on elementary properties
of continuos functions on compact spaces.
\medskip

{\it Sketch of a proof of the inverse function theorem.}
We can assume for simplicity, and it is not restrictive, that $p=0$ and $f(p)=0$.
Possibly by composing $f$ with $(d_0f)^{-1}$ we can also assume that $d_0f=I_n$. 
\medskip

The proof is achieved  by following the next sequence of claims.
\medskip
 
{\bf Claim 1.} There is a sufficiently small closed ball $\overline B = \overline B^n(0,\epsilon)\subset U$
such that
\begin{enumerate}
\item For every $x\in \overline B$, $d_xf$ is invertible;
\item  For every $x\in \overline B$, $x\neq 0$, then $f(x)\neq 0$;
\item For every $x,z \in \overline B$, $2|| f(x)-f(z)||\geq ||x-z||$.
\end{enumerate}

Assuming these facts, by the continuity of the function and the compactness of $\partial \overline B$,
there is $\delta >0$ such that for every $x\in \partial \overline B$, $||f(x) ||\geq \delta$.
Consider the open ball $B'= B^n(0,\delta/2)$.

\medskip

{\bf Claim 2.} Set $A=  B\cap f^{-1}(B')$. Then the restriction $\phi: A \to B'$ of $f$ to the open set $A$ is bijective.

\medskip

{\bf Claim 3.} $\phi$ is a homeomorphism. 

\medskip

{\bf Claim 4.} $\phi$ is a diffeomorphism

\medskip

{\it Proof of Claim 1.}  The first point is evident. Assuming that the second point fails, there would
be a sequence $x_n$ in $U$, converging to $0$, such that $f(x_n)=0$ for every $n$.
Hence  $||\frac{f(x_n) - x_n}{x_n}||=1 $ against the fact that $d_0f=I_n$. As for the third point, 
consider the function $g(x)=f(x)-x$, so that 
$$||x-z||-||f(x)-f(z)|| \leq ||g(x)-g(z)|| \ . $$
As 
$$ \frac{\partial g_i}{\partial x_j}(x)= \frac{\partial f_i}{\partial x_j}(x)- \frac{\partial f_i}{\partial x_j}(0) $$
we can take $\epsilon$  in order to make $|\frac{\partial g_i}{\partial x_j}(x)|<\frac{1}{2n^2}$
uniformely on $\overline B$. Then the conclusive inequality
$$||g(x)-g(z)||\leq \frac{1}{2} ||x-z||  $$
is obtained by 
applying several times the Main Value Theorem for functions of one variable,

\medskip

{\it Proof of Claim 2.} It is enough to prove that for every $y\in B'$ there is a unique
$x\in U$ such that $f(x)=y$.  The smooth function $h(x)=||y-f(x)||^2$ has a minimum point $p$
on the compact set $\overline B$ and by construction $p$ belongs necessarily to the open ball $B$.
A simple computation then shows that $d_pf(y-f(p))=0$, hence $y-f(p)=0$ because $d_pf$ is invertible.
As for the uniqueness, this follows by the inequality $||p_1 -p_2||\leq 2||f(p_1)-f(p_2)||$, so that $p_1=p_2$
if $f(p_1)=f(p_2)=y$.

\medskip

{\it Proof of Claim 3.} The same inequality implies that $||\phi^{-1}(y_1)-\phi^{-1}(y_2|| \leq 2||y_1-y_2||$
and the continuity of $\phi^{-1}$ follows.

\medskip

{\it Proof of Claim 4.} As we know, it is enough to show that $\phi^{-1}$ is differentiable.
In fact by using directly the definition of the differential one can prove that
$d_y \phi^{-1} = (d_x\phi)^{-1}$, where $y=\phi(x)$. The details are left to the reader.

\cvd

\section {Topologies on spaces of smooth maps}\label{funct-space}
Let $U\subset \R^n$, $W\subset \R^m$ be  open sets.
For every $k\geq 0$ we define a topology $\delta_k$ on $\Cc^k(U,W)$;
we will denote by $\Ee^k(U,W)$  the set $\Cc^\infty(U,W)\subset \Cc^k(U,W)$
endowed with the subspace topology.
We determine $\delta_k$ by giving for every $f=(f_1,\dots, f_m) \in \Cc^k(U,W)$
a basis of open neighbourhoods $\Uu_k(f, K, \epsilon)$ where the varying
arguments are a compact subset $K\subset U$ and a real $\epsilon >0$.
Then, by definition, $g\in \Cc^k(U,W)$ belongs to $\Uu_k(f,K,\epsilon)$ if
and only if 
\begin{enumerate}
\item For every $x\in K$, $||g(x)- f(x)|| <\epsilon$;
\item For every multi-index $J$ such that $|J|=r \leq k$, for every $i=1,\dots, m$,
for every $x\in K$, we have
$$|| \frac {\partial^J(g_i-f_i)} {\partial x_1^{j_1}\dots \partial x_n^{j_n}} (x)|| < \epsilon \ . $$
\end{enumerate}

We omit the proof that this actually defines bases of neighbourhoods of some topologies.

We denote by $\Ee(U,W)$ the set $\Cc^\infty(U,W)$ endowed with the {\it union topology}
$\delta = \cup_k \delta_k$. 

All these are called {\it weak topologies}. This understands
the existence of other {\it strong} topologies, say $\sigma_k$, on the same sets. 
By considering for example
$\Ee(\R^n,\R)$, we can control the difference of two functions, up to an arbitrarily prescribed
order on an arbitrarily given compact set $K$, but we have not any control ``at infinity''.
The strong topologies $\sigma_k$, which contain $\delta_k$ being {\it heavily finer}, allow instead
such a control at infinity. On another hand, the weak topologies $\delta_k$ have nice
properties, for example one can prove that they are metrizable, hence every $f$ has a countable
basis of open neighbourhoods. On the contrary this is not the case for the strong topologies;
for example if a sequence $g_n \to f$ in  $\Cc^k (\R^n,\R)$ with the strong topology, then there exists a compact
set $K$ in $\R^n$ such that $g_n$ definetly equal $f$ on the complement of $K$.
However, {\it we do not define the strong topologies}. To our aims, the control at compact sets
will suffice.

\medskip

Let us recall also (a particular case of) the classical {\it Stone-Weierstrass theorem} (see \cite{Stone}).

\begin{theorem}\label{stone-weier} For every $f\in \Cc^k(U,\R^m)$, for every $k\geq 0$, for every neighbourhood $\Uu=\Uu_k(f,K,\epsilon)$,
there exists a polynomial map $p: \R^n \to \R^m$ such that the restriction of $p$ to $U$ belongs to $\Uu$. In other words, the polynomial
maps are {\rm dense} in  $\Cc^k(U,\R^m)$ for every $k\geq 0$ and in $\Ee(U,W)$.
\end{theorem}

\section{ Stability of summersions and immersions at a compact set}\label{stab-summ}
Let $f\in \Cc^k(U,W)$ be as above, $k\geq 1$, $K\subset U$ a compact set. We say that $f$ is
{\it a summersion (resp. an immersion) at} $K$ if for every $x\in K$, $d_xf$ is surjective
(resp. injective). This is equivalent to the fact that there exists an open neighbourhood 
$K\subset U' \subset U$ such that the restriction of $f$ to $U'$ is a 
summersion (immersion). Here is the stability results. 

\begin{proposition}\label{stab}
If $f$ is either (1) a summersion, (2) an immersion or (3) an injective immersion  at $K$, 
then there is a neighbourhood $\Uu=\Uu_1(f,K,\epsilon)$ such that every $g\in \Uu$ 
shares the same properties of $f$ respectively. 
\end{proposition}

\Dim If $n\geq m$ (resp. $n<m$), then every $m\times n$ matrix $A$   has $\binom{n}{m}$ $m\times m$
(resp. $\binom{m}{n}$) submatrices say $A_j$; in any case define 
$$\delta(A)= \sum_j (\det A_j)^2 \ . $$
In both cases (1) and (2) the hypothesis is equivalent to $d(x):=\delta(d_xf) >0$ for every $x\in K$.
As  $d$ is continuous and $K$ is compact, then 
$$\sup_{x\in K} \{d(x)\} = \max_{x\in K} \{d(x)\} =d_0 >0 \ . $$
Then it is clear that if $g$ is close enough to $f$ at $K$ in $\Cc^1(U,W)$, then $\delta(d_xg)>0$
for every $x\in K$. As for (3), assume that the thesis fails. Then there would exist a seguence $g_n\in \Cc^1(U,W)$,
sequences of points $x_n$, $y_n$ in the compact set $K$ such that: 
\begin{enumerate}
\item Every $g_n$ is an immersion at $K$ (by (2));
\item $g_n \to f$ and  $dg_n \to df$ uniformly
on $K$; 
\item $x_n\to x$, $y_n \to y$ in $K$, $x_n \neq y_n$ and $g_n(x_n)=g_n(y_n)$ for every $n$.
\item  $v_n:=\frac{x_n -y_n}{||x_n -y_n||}\to v \in S^{n-1} $ (by the compactness of the unitary sphere $S^{n-1}$).
\end{enumerate}

Then: $g_n(x_n)\to f(x)$, $g_n(y_n)\to f(y)$, hence $x=y$ because $f$ is injective. 
Hence
$$ || g_n(x_n) - g_n(y_n) - d_{y_n}g_n(x_n-y_n))||/||x_n-y_n||\to 0$$
uniformly, so that 
$$||d_{y_n}g_n(v_n)|| \to ||d_xf(v)||=0 \ . $$ 
This is absurd because $d_xf$ is injective.

\cvd

\section {An elementary division theorem}\label{division}
By definition a {\it convex} subset $C$ of $\R^n$ has the property that for every $x_0,x_1\in C$,
the (parametrized) segment $\gamma: [0,1]\to \R^n$, $\gamma(t)=(1-t)x_0 + tx_1$ is enterely
contained in $C$. We have
\begin{theorem}\label{elem-div} {\rm (Elementary division theorem)} Let 
$f=(f_1,\dots, f_m) \in \Cc^\infty(U,\R^m)$ where $U\subset \R^n$ is a convex
open subset. Assume that $0\in U$ and $f(0)=0$. Then there are smooth
maps $g_j =(g_{j1},\dots, g_{jm}): U\to \R^m$, $j=1,\dots, n$, such that for every $x\in U$,
$ f(x)=\sum_j  x_jg_j(x) $, and (necessarily) $g_{ji}(0)= \frac{\partial f_i}{\partial x_j}(0)$. 
\end{theorem}

\Dim   It is a basic property of elementary integration that
for every smooth function $h: U \to \R$, the function $\tilde h: U\to \R$
defined by $\tilde h (x) = \int_0^1 h(tx)dt$ is smooth. 
By the funtamental theorem of elementary integration for continuous functions, we have
that 
$$f(x)= \int_0^1 \frac {d f(tx)}{dt}dt = (\int_0^1 \frac {d f_1(tx)}{dt}dt, \dots, \int_0^1 \frac {d f_m(tx)}{dt}dt) \ . $$
By the chain rule, for every $i=1,\dots, m$,
$$\int_0^1 \frac {d f_i(tx)}{dt}dt = \int_0^1 (\sum_j  x_j\frac{\partial f_i}{\partial x_j}(tx))dt = 
\sum_j x_j \int_0^1  \frac{\partial f_i}{\partial x_j}(tx))dt \ . $$ 
We achieve the proof by setting 
$$g_{ji}(x):= \int_0^1  \frac{\partial f_i}{\partial x_j}(tx))dt $$

\cvd

\medskip

\begin{remarks}\label {division2} {\rm (1) The same arguments in the above proof
work as well by assuming only that $U$ is 
{\it starred} with center at $0$. 
\smallskip

(2) In the setting of the division theorem, if $n=m=1$, we have that $f(x)=xg(x)$,
that is the coordinate function $x$ divides $f$. Assume now $m=1$,  $f$ is defined
on an open set of the form $U=A\times (-1,1) \subset \R^{n-1}\times \R$ and that
$\{f=0\}= U\cap \{x_n=0\}$. Then by applying fibre by fibre the same construction
of the above proof, we get that $f(x)=x_{n}g(x)$. Moreover, if $f$ is a summersion,
then $g(x)$ is nowhere vanishing.}
\end{remarks}

\medskip

We will see several applications of the division theorem.

\section{A differential interpretation of the tangent spaces: derivations}\label{derivation}
Above we have introduced the tangent spaces $T_xU$,  mainly by considering
$U$ as an open set of the {\it affine} space $\R^n$. Here we give a genuine differential
interpretation, compatible with the already defined tangent functor.

Let $p\in U$. Consider the set of smooth functions $f:U'\to \R$
defined on some open neighbourhood $U'$ of $p$ in $U$. On this set put the
equivalence relation such that $(U_1,f_1)\sim (U_2,f_2)$ if and only if there is
$(U_3, f_3)$ such that $U_3\subset U_1\cap U_2$ and for every $y\in U_3$
$f_3(y)=f_1(y)=f_2(y)$. Denote by $\Ee_p$ the quotient set. Note that
$U$ is immaterial for this purely local definition, as we would get the same $\Ee_p$
by taking for instance the whole of $\R^n$ instead of $U$. Similarly
also $T_p=T_pU$ essentially does not depend on the choice of the open set containing $p$. 
Denote by $[f]=[f]_p$ an equivalence
class. The usual sum and product defined on every $\Cc^\infty(U',\R)$
induce well defined sum and product on $\Ee_p$ which make it a {\it commutative ring}
as well a real {\it vector space} with compatible operations. The translation
$x\to x-p$ determines
a canonical isomorphism between $\Ee_p$  and $\Ee_0$, then the considerations we are going to
do for $\Ee_0$ can be straighforwardly transported to $\Ee_p$ by this translation. 
Let $v=(v_1,\dots, v_n)^t \in T_0\sim \R^n$.
By  means of the usual   
{\it derivative at $0$ in the direction $v$}, we define the function 

$$\delta_v: \Ee_0 \to \R, \ \delta_v([f])= \sum_j \frac {\partial f}{\partial x_j}(0)v_j \ . $$
One verifies easily that
$\delta_v$ is well defined (it does not depend on the choice
of the representative $f$), is $\R$-linear, and moreover verifies
the {\it Leibniz identity}:
$$ \delta_v ([f][g]) =  f(0)\delta_v([g]) + g(0)\delta_v([f]) \ . $$
Let us call a {\it derivation on $\Ee_0$} any map $\delta: \Ee_0 \to \R$
that verifies the same properties. Set ${\rm Der}(\Ee_0)$ the set of these derivation.
It has a natural structure of real vector space, so that the map
$$ L: T_0 \to {\rm  Der}(\Ee_0), \ L(v)=\delta_v$$
is $\R$-linear. Let us prove that {\it $L$ is a linear isomorphism}. For every
derivation $\delta$ we will find a unique $v\in T_0$ such that $\delta=\delta_v$.
It follows immediately from the derivation properties that for every constant germ
$[f]$ (i.e. with a constant representative), $\delta([f])=0$. For every $[f]$
we can take a representative $f$ defined on a small open ball $B^n(0,\epsilon)$
(which is convex). By the division theorem, for every $x$ in such a ball,
$$ f(x)-f(0)= \sum_j g_j(x)x_j$$
for some smooth functions $g_j$. Then, by using again the derivation properties, we have
$$\delta([f])= \sum_j \frac{\partial f}{\partial x_j}(0)\delta([x_j])$$
hence we conclude by taking $v= (\delta([x_1]),\dots, \delta([x_n]))$.

\cvd

The above discussion can be {\it globalized} by replacing $T_p$
with the set $\Gamma(T(U))$ of (tangent) vector fields on $U$, and
$\Ee_p$ with the commutative ring $\Cc^\infty(U,\R)$ with the induced
compatible structure of $\R$-vector space. $\Gamma(T(U))$ is also in a natural way
a vector space and this extends to a natural structure of $\Cc^\infty(U,\R)$-{\it module}.
For every vector field $X\in \Gamma(T(U))$, define
$$\delta _X: \Cc^\infty(U,\R)\to \Cc^\infty(U,\R), \ \delta_X(f)(x)= \delta_{X(x)}([f]_x) \ . $$
It is $\R$-linear and verifies the Leibniz rule
$$ \delta_X(fg)= f\delta_X(g) + \delta_X(f)g$$
hence, by definition, it is a {\it derivation} on $\Cc^\infty(U,\R)$.
Finally the map
$$L: \Gamma(T(U))\to {\rm Der}(\Cc^\infty(U,\R)), \ L(X)=\delta_X$$
establishes an isomorphism of $\Cc^\infty$-modules.

Note that if $\delta, \delta' \in {\rm Der}(\Ee_0)$ (resp. $\in {\rm Der}(\Cc^\infty(U,\R))$)
then $\delta\delta'$ is not in general a derivation, while $\delta\delta' - \delta'\delta$ is a
derivation. In particular for every couple $X,Y\in \Gamma(T(U))$ there is a unique
vector fields $[X,Y]$ such that
$$ L([X,Y]) = L(X)L(Y) - L(Y)L(X) \ . $$

\section{Morse lemma}\label{morse}
Let $f\in \Cc^\infty(U,\R)$, $U$ open set of $\R^n$. 
A point $p\in U$ is {\it regular} for $f$ if $d_pf\neq 0$ (that is $f$ is a summersion near $p$);
otherwise we say
that $p$ is {\it critical} (or also {\it singular}). We are interested to
the local behaviour of $f$ at $p$ (actually to the germ $[f]_p$).
Up to pre or post composition with a translation we can normalize the situation
so that $p=0$ and $f(0)=0$. Moreover we can assume that $U=B^n(0,\epsilon_0)$
for some $\epsilon_0>0$ and, case by case, we can restrict $f$ to any
$U_\epsilon= B^n(0,\epsilon)$, $0< \epsilon \leq \epsilon_0$.
For every $\epsilon$, the commutative ring $\Cc^\infty(U_\epsilon,\R)$ has a canonical ideal
$$ \mG_\epsilon = \{g\in \Cc^\infty(U_\epsilon,\R); \ g(0)=0 \}$$
so that we are assuming that $f\in \mG_\epsilon$.  It is an immediate
corollary of the division theorem that $\mG_\epsilon$ is generated
by the coordinate functions $x_j$, $j=1,\dots, n$; that is every
$g\in \mG_\epsilon$ is a $\Cc^\infty(U_\epsilon,\R)$-linear combination
of the coordinate functions. Hence we have that for $x\in U$,
$$f(x)=\sum_j g_j(x)x_j, \  d_0f=(g_1(0),\dots, g_n(0))\ .$$

If $0$ is a regular point for $f$, the particular case of theorem \ref{summersion} 
can be rephrased
by saying that, up to pre composition with a local diffeomorphism,  $f$ locally
coincides with $d_0f$ that is its first Taylor polynomial $T_1(f)(0)$.

If $0$ is critical, then all the smooth functions $g_j$ vanish at $0$, and we can apply again 
to each of them the division theorem and eventually we get that on $U$
$$ f(x)= \sum_{|J|=2} g_{J}(x) x^J, \ x^J:= x_1^{j_1}\dots x_n^{j_n}$$
that is it has the form of a homogeneus polynomial of degree $2$ whose coefficients
are  smooth functions. Moreover
$$ T_2(f)(0)=\sum_{|J|=2} g_{J}(0)x_1^{j_1}\dots x_n^{j_n} \ . $$
In fact we can express $ T_2(f)(0)$ in the form
$$T_2(f)(0)= \frac{1}{2}x^tH_0(f)x:= Q_0(f)(x)$$
where $H_0(f)$ is the {\it symmetric} (by Schwartz Lemma) {\it Hessian matrix} of $f$  at $0$
$$ H_0(f)=  \left (\frac{\partial^2 f}{\partial x_i \partial x_j}(0)\right )_{i,j=1, \dots, n} $$
while $Q_0(f)$ is the associated {\it quadratic form}.
We can organize the above functions $g_{i,j}$ to rewrite $f$ on $U$ as
$$f(x)= x^t G(x) x$$
where 
$$G: U\to M(n,\R)$$ 
is a smooth map such that $G(x)=G(x)^t$ is symmetric
for every $x\in U$, and $G(0)=T_2(f)(0)$.

We say that the critical point $x=0$ is {\it non degenerate} if
$$\det H_0(f)\neq 0 \ . $$
We have the following characterization of non degenerate critical points.
For every $U_\epsilon$, denote by $J(f,\epsilon)$ the {\it Jacobian ideal} of $\Cc^\infty(U_\epsilon,\R)$
generated by the partial derivative functions  $\frac{\partial f}{\partial x_j}$, that is the ideal
of the $\Cc^\infty(U_\epsilon, \R)$-linear combinations $\sum_j h_j\frac{\partial f}{\partial x_j}$,
$h_j\in \Cc^\infty(U_\epsilon, \R)$. If $0$ is a critical point, then
$J(f,\epsilon) \subset \mG_\epsilon$ . Then we have

\begin{lemma}\label{J=m} $0$ is a non degenerate critical point of $f\in \Cc^\infty(U,\R)$, $f(0)=0$,
if and only if there exists $0<\epsilon \leq \epsilon_0$ such that $J(f,\epsilon)=  \mG_\epsilon$.
\end{lemma}
\Dim It is enough to prove the inclusion ``$\supseteq$'', then it is enough to show that the generating
coordinate functions $x_j$ belong to $J(f,\epsilon)$ for some $\epsilon$. As $0$ is non degenerate, the smooth map 
$$x\to (\frac{\partial f}{\partial x_1}(x), \dots,  \frac{\partial f}{\partial x_n}(x))$$
has invertible differential at $0$, then we can apply the inverse map theorem locally on a neighbourhood
$U_\epsilon$ of $0$, so that there are smooth functions $F_j$ such that for every $j=1,\dots,n$, $F_j(0)=0$ and 
$$x_j=F_j(\frac{\partial f}{\partial x_1}, \dots,  \frac{\partial f}{\partial x_n}) \ . $$ 
Again by the division theorem we finally get
$$ x_j=\sum_i G_{j,i}(x) \frac{\partial f}{\partial x_i}(x)$$
and the Lemma is proved.

\cvd

\medskip

Assume that $0$ is a non degenerate critical point for $f$.
We are going to prove  that up to pre composition with local diffeomorphisms at $0$,  $f$ locally
coincides with $T_2(f)(0)$. More precisely, the Hessian matrix $H_0(f)$  has a certain
{\it index of negativity} $0\leq\lambda\leq n$ (i.e. the maximal dimension of the linear
subspaces of $\R^n$ on which the restriction of the quadratic form $Q_0(f)$ is negative).
By definition $\lambda$ {\it is the index of the non degenerate critical point $0$}. 
This notion is stable under local diffeomorphism.

\begin{lemma} If $\ 0 \ $ is a non degenerate critical point of index $\lambda$ of $f\in \Cc^\infty(U_\epsilon,\R)$,
      $f(0)=0$, and $\phi:W\to U_\epsilon$ is a diffeomorphism, $\psi(0)=0$, then $0$ is a non degenerate
      critical point of $f':=f\circ \phi$ of index $\lambda$.
\end{lemma}
\Dim By direct computation, using the chain rule an the fact that $d_0f=0$, we have
$$H_0(f')= d_0\phi^t H_0(f) d_0\phi$$
hence the symmetric matrices $H_0(f')$ e $H_0(f)$ are congruent so they are both
non singular and have the same signature.

\cvd

\medskip

Let $0$ be a non degenerate critical point of $f$ of index $\lambda$.
Up to composition with a linear isomorphism $x=Pu$, we have that
$$ Q_0(f)(Pu)= -(\sum_{j=1}^\lambda u_j^2)+(\sum_{j=\lambda+1}^n u_j^2) = u^tI_{n,\lambda}u $$
where $I_{n,\lambda}$ is the suitable diagonal matrix with $\pm 1$ entries.
Finally we   can state
\begin{theorem}\label{Morse-lemma} {\rm (Morse Lemma)} Let $0$ be a non degenerate critical point of index
  $0\leq \lambda \leq n$ of $f\in \Cc^\infty(U,\R)$, $f(0)=0$.
  Then there is a local diffeomorphism $x=\phi(u)$, $0=\phi(0)$, such that $\psi:=\phi^{-1}$
is defined on some $U_{\epsilon}$ and
$$ f(\phi(u))= u^t I_{n,\lambda}u \ . $$
\end{theorem}
\Dim It is not restrictive to assume that
$$ H_0(f)=I_{n, \lambda} \ . $$
Let us take as above on $U$ an expression
$$f(x)= x^t G(x) x \ . $$
If $\epsilon>0$ is small enough, we have that on $U_\epsilon$ every symmetric matrix
$G(x)$ has negativity index $\lambda$, and by applying to the canonical basis of
$\R^n$ the usual algorithm producing a normalized othogonal basis with respect to
the scalar product $(*,*)_{G(x)}$, we eventually get a smooth map
$$P: U_\epsilon \to {\rm GL}(n,\R)$$
such that:
\begin{enumerate}
\item $P(0)=I_n$;
\item For every $x\in U_\epsilon$, the linear isomorphism $x=P(x)u$ is such that
$$ P(x)^tG(x)P(x)=I_{n,\lambda} \ . $$
\end{enumerate}

Then consider the smooth map

$$\psi: U_\epsilon \to \R^n, \ \psi(x)=P(x)^{-1}x$$
one verifies that $\psi$ has invertible differential at $0$, so by the inverse map theorem,
possibly by shrinking $\epsilon$, $u=\psi(x)$ is a diffeomorphism onto its open image
and finally
$$f(x)=x^tG(x)x= u^tI_{n,\lambda}u \ $$
as desidered.

\cvd
  
Let us state, without proof, an interesting generalization of Morse's Lemma.
With the usual notation, for every $k\geq 1$, define
$\mG_\epsilon^k$ as the ideal of $\Cc^\infty(U_\epsilon,\R)$ generated by the monomials
$x^J=x^{j_1}_1\dots x^{j_n}_n$, $J$ be an arbitrary multi-index with $|J|=k$.
Clearly $\mG_\epsilon = \mG_\epsilon^1\subset \mG_\epsilon^2 \subset \dots$.
We have
\begin{proposition}\label{finite-det} Let $f\in \Cc^\infty(U,\R)$, $f(0)=0$, be such that
  $0$ is a critical point, and there is $k\geq 1$ such that $\mG_\epsilon^k\subset J(f,\epsilon)$.
  Then up to pre composition with a local diffeomorphism at $0$, $f$ locally coincides with
  the Taylor polynomial $T_k(f)(0)$.
\end{proposition}

\section{Bump functions and partitions of unity}\label{bump}
Consider the function $\alpha: \R \to \R$ defined by $\alpha(x)=0$ if
$x\leq 0$, $\alpha(x)=e^{-\frac{1}{x}}$ if $x>0$.  One verifies that
$\alpha$ is smooth and that for every $k\geq 1$,
$\frac{d^kf}{dx^k}(0)=0$. Then we say that $\alpha$ is {\it flat} at
$0$, although $\alpha$ is not locally constant at $0$. This phenomenon
is an important feature of the ``flexibility'' of smooth functions
that makes them suited for topological applications. On the contrary,
for example, analytic functions are much more rigid: an analytic
function on $\R$ which is flat at some points is constant.

Let us fix two real numbers $0< a < b$. Define $\beta=\beta_{a,b}: \R
\to \R$, $$\beta(x)=\alpha(x-a)\alpha(b-x) \ . $$ Hence $\beta$ is
smooth, $\beta(x)=0$ on $\{x\leq a\}\cup \{x\geq b\}$, is strictly
positive on $\{a<x<b\}$ with a unique maximum; $\beta$ is flat at $a$
and $b$.

Define $\gamma= \gamma_{a,b}: \R \to \R$ by
$$ \gamma(x)= \frac{\int_{|x|}^b \beta(t)dt}{\int_a^b \beta(t)dt}
\ . $$ Then $\gamma$ is smooth, $\gamma(x)=1$ if $|x|\leq a$,
$\gamma(x)=0$ if $|x|\geq b$, $0 \leq \gamma(x) \leq 1$ and is
monotone on each connected interval of $\{a<|x|<b\}$; $\gamma$ is flat
at $\pm a$ and $\pm b$.  For every $n\geq 1$ we can define $\gamma_n:
\R^n \to \R, \ \gamma_n=\gamma_{n,a,b}(x)=\gamma_{a,b}(||x||)$,
however we will omit the index $n$ whenever the dimension is clear by
the contest. Such a function $\gamma_{a,b}: \R^n \to \R$ is called a
{\it bumb function} on $\R^n$ with center $0$ and rays $a,b$. If
$\tau_p(x)= x-p$, then
$$\gamma_{p,a,b}=\gamma_{a,b}\circ \tau_p$$ is a bump function with
center $p$; when the center is clear from the context we will omit
also to indicate it.

Recall that the {\it support} of a function is the closure of the set
where it is not zero. Hence $\overline B^n(p,b)$ is the support of
$\gamma_{p,a,b}$.

We introduce also bump functions ``at infinity'' as follows. Let $\R^n
\subset \R^{n+1}$ as the hyperplane with equation $x_{n+1}=0$. Denote
by $\pi^+: S^n\setminus \{e_{n+1}\} \to \R^n$
($e_{n+1}=(0,\dots,0,1)$) the {\it stereographic projection} defined
geometrically by
$$\pi^+(x)= r(x,e_{n+1})\cap \R^n$$ where $r(x,e_{n+1})$ is the
straight line passing through the two points. Similarly define the
projection $\pi^-:S^n\setminus \{-e_{n+1}\} \to \R^n$. One easily
verifies by direct computation that
$$\rho := \pi^- \circ (\pi^+)^{-1}: \R^n \setminus \{0\} \to \R^n \setminus \{0\} $$
is a diffeomorphism. Then a {\it bump function at infinity} is by definition of the form
$ \gamma_\infty(x) = \gamma \circ \rho (x)$ if $x\in \R^n \setminus \{0\}$,  $\gamma_\infty(0)=0$
which clearly is smooth. 

We extend now the definition to bump functions at an arbitrary compact
subset of $K\subset \R^n$, as follows. Let $K\subset \R^n$ be a
compact set, $U$ an open neighbourhood of $K$. Then we can find
$W_0:=U_{\infty,a_\infty}:= \R^n \setminus \overline B^n(0,a_\infty)$,
some $W_j:= B^n(p_j,b_j)$, $j=1,\dots, k$, and some $0<a_j<b_j$,
$a_\infty < b_\infty$ such that:
\begin{enumerate}
\item   $\overline W_0 \cap K = \emptyset$;
\item The open balls $U_j:=B^n(p_j,a_j)$ together with
  $U_0:=U_{\infty,b_\infty}$ make a finite open covering $\Uu$ of
  $\R^n$;
\item The union of the above open balls that intersect $K$ is an open
  neighbourhood $U'\subset U$ of $K$.
\end{enumerate}

Denote by $\gamma_0$ the bump function at infinity with support equal
to $\overline W_0$ and constantly equal to $1$ on $U_0$; by $\gamma_j$
the bump function at $p_j$ with rays $a_j,b_j$.  For every $j=0,\dots,
k$, define the smooth function
 $$\lambda_j:= \frac {\gamma_j}{\sum_j \gamma_j} \ . $$ By the
properties of the covering $\Uu$ and of the bump functions, the
denominator is strictly positive everywhere.  Clearly, for every $x\in
\R^n$,
 $$ \sum_j \lambda_j(x)=1 \ . $$ Such a family of function
$\{\lambda_j\}$ is called a {\it partition of unity subordinate to the
  (finite) covering} $\Uu$.  Now we define ``local'' constant
functions $c_j: W_j\to \R$, such that $c_j=1$ if $U_j\cap K$ is non
empty, $c_j=0$ otherwise. Finally set
$$\gamma_K = \sum_j \lambda_j c_j \ . $$ By construction it is smooth,
it is constantly equal to $1$ on $U'$ and has compact support
contained in $U$.  Any $\gamma_K$ constructed in this way is called a
{\it bump function at $K$}.

Bump functions are an important device. A basic use is the following:
let $U$ be an open neighbourhood of a compact set $K$ as above,
$f:U\to \R$ be a smooth function {\it locally} defined at $K$. In
certain cases it is useful to find a {\it globally} defined smooth
function $\hat f: \R^n \to \R$ with compact support and which locally
agrees with $f$ at $K$, that is there is an open neighbourhood
$K\subset U' \subset U$ such that $f(x)=\hat f(x)$ for $x\in U'$.
Take any bump function $\gamma=\gamma_K$ at $K$ constructed as above;
then $\hat f$ defined by $\hat f(x)= \gamma(x)f(x)$ if $x\in U$, $\hat
f(x)=0$ if $x\in \R^n \setminus U$, does the job.

These partitions of unity provide also a very flexible way to
construct riemannian metrics on $\R^n$.  Let $\{\lambda_j\}$ be as
above. Fix on every $U_j$ an arbitrary riemannian metric $g_j$ (for
instance a constant one varying with $j$). Then
$$g= \sum_j \lambda_jg_j$$
is a well defined riemannian metric on the whole of $\R^n$.

\medskip

In the next sections we will see a few other concrete applications.

\begin{remark}\label{part-unity}
  {\rm As $\R^n$ is metrizable, locally compact and with a countable
    basis of open sets, one can prove that for every open set
    $U\subset \R^n$, for every open covering $\Aa$ of $U$ there exist
    a countable family of open balls $\Bb=\{B_j=B^n(p_j,b_j)\}_{j\in
      \N}$, and for every $j\in \N$, $0<a_j<b_j$ such that
\begin{enumerate}
\item $\Bb$ is a {\it refinement} of $\Aa$, that is it is a open covering of $U$ 
and every $B_j$ is contained in some $A\in \Aa$;
\item $\Bb$ is {\it locally finite}, that is for every $p\in U$, there is a ball $B=B^n(p,r)$ 
which intersects only finitely many $B_j$'s;
\item Also $\Uu=\{ B^n(p_j,a_j)\}_{j\in \N}$ is an open covering of $U$.
\end{enumerate}

Take the corresponding family of bump functions $\{\gamma_j= \gamma_{p_j,a_j,b_j}\}$.
Set, for every $j\in \N$
$$\lambda_j =\frac{\gamma_j}{\sum_{j=1}^\infty \gamma_j} \ . $$ This
is well defined and smooth because, by the local finiteness, the
denominator reduces at every point $p$ to a strictly positive sum of a
finite number of terms.  Clearly for every $p$, $$\sum_j \lambda_j(p)
= 1 \ . $$ The family $\{\lambda_j\}$ is called a partition of unity
subordinate to the covering $\Bb$ (which refines the given $\Aa$).
For example, if $K \subset U$ is a compact set as above we could apply
the construction to the open covering of $\R^n$, $\Aa=\{U, \R^n
\setminus K\}$, and use the resulting partition function over $\Bb$ to
construct as well a bump function $\gamma_K$ at $K$.  These more
general partitions of unity rely on a topological property called {\it
  paracompacteness}; however, {\it we will not really need them}.}
\end{remark} 

\section{Homotopy, isotopy, diffeotopy}\label{H-I-D}
Here we fix a few notions and terminology that shall be widely employed 
and developed. $U$ and $V$ are open sets in Euclidean spaces. A map
$$F: U\times [0,1] \to V$$ is {\it smooth} if it is the restriction of a smooth map
defined on the open set $U\times J$, $J$ being an open interval and $[0,1] \subset J$. 
For every $t\in [0,1]$, set $f_t$ the restriction of $F$ to $U\times \{t\}$.
Then $F$ is called a (smooth)  homotopy between $f_0$ and $f_1$.
It can be considered as a continuous path in $\Ee(U,V)$ 
joining $f_0$ and $f_1$.

We say that $f:U\to V$ is an {\it embedding} if $f$ is an injective immersion and is a homeomorphism onto
its image. If $f_t$ is an embedding  for every $t\in [0,1]$,
then $F$ is called an {\it isotopy} between $f_0$ and $f_1$.

If $U=V$ and $f_t$ is a diffeomorphism for every $t\in [0,1]$,
then $F$ is called a {\it diffeotopy}.
In this case 
$F$ can be reconsidered as follows: consider 
the map 
$$H:U\times [0,1]\to U\times [0,1], \ H(p,t)=(f_0(p),t) \ . $$
Then $G:=F\circ H^{-1}$ is a diffeotopy between ${\rm id}_U$ and $f_1\circ f_0^{-1}$,
and $F= G\circ H$. This formal manipulation suggests nevertheless the following
specialization of homotopy. If $G: V\times [0,1] \to V$
is a diffeotopy between $g_0={\rm id}_V$ and $g_1$, and 
$\phi:U\times [0,1]\to V\times [0,1]$ is of the form $\phi(p,t)=(f(p),t)$ for some
$f: U\to V$, then
$G\circ F$ is called a {\it diffeotopy} between $f_0:=f$ and $f_1:= g_1\circ f$;
sometimes one also says  that $f_0$ and $f_1$ are homotopic through an
{\it ambient isotopy}.

Let $f:U \to U$ be a diffeomorfism. The {\it support of $f$} is
the closure of the subset of $U$ on which $f(x)\neq x$. 
If $F$ is a diffeotopy between $f$ and ${\rm id}_U$ the support of
$F$ is the closure of the union of supports of the $f_t$'s. 

Homotopy and its relatives define equivalence relations on the pertinent
space of maps. Clearly $F(p,t)=(f(p),t)$ is a homotopy between $f$ and itself.
If $F$ is a homotopy between $f_0$ and $f_1$, then $\hat F(p,t):= F(p,1-t)$
is a homotopy between $f_1$ and $f_0$. As for the transitivity: by using the
$1$-dimensional bump functions, we see that there exist 
a smooth function
$s:[0,1]\to [0,1]$ and $1/3> \epsilon >0$ such that $s(t)=0$ on $[0,\epsilon]$,
$s(t)=1$ on $(1-\epsilon, 1]$, and $s$ is a diffeomorphim on $[\epsilon, 1-\epsilon]$.
If $F$ is any homotopy between $f_0$ and $f_1$, then replace it by
$\tilde F(p,t)=F(p,s(t))$. If a homotopy $\tilde F'$ connects $f_0$ and $f_1$,
while $\tilde F"$ connects $f_1$ and $f_2$, then set
$$\tilde F(p,t)=\tilde F'(p,2t), \ t\in [0,1/2]$$
$$\tilde F(p,t)=\tilde F"(p,2t-1), \ t\in [1/2,1] \ . $$
It is a {\it smooth} homotopy between $f_0$ and $f_2$.
For isotopies and diffeotopies we argue similarly.


\section{Linearization of diffeomorphisms of $\R^n$ up to isotopy}\label{linearization}

We have

\begin{proposition}\label{lin1} Every diffeomorphism $f:\R^n \to \R^n$, $f(0)=0$,
  is diffeotopic to the differential $d_0f\in {\rm GL}(n,\R)$, through diffeomorphisms $f_t$ such that
  $f_t(0)=0$ for every $t\in \R$.
\end{proposition}
\Dim Define $F:\R^n\times \R \to \R^n$, by $F(x,t)= f(tx)/t$ if $t\neq 0$,
$F(x,0)= d_0f$. It follows from the very definition of the differential that
$F$ is continuous; clearly it is smooth where $t\neq 0$. To check that it is fully
smooth we note that by the division theorem $F(x,t)= \sum_j g_j(y)x_j, \ y=tx$
the $g_j$ being smooth maps of $y$.

\cvd

We can strenghten the above Proposition. Let us set ${\rm GL}^\pm={\rm
  GL}^{\pm}(n,\R)$ the open subsets of ${\rm GL}(n,\R)$ formed by the
matrices $A$ such that either $\det A >0$ or $\det A <0$.
Take the identity $I_n$ and the matrix $I_{n,1}$ (the
notation has been introduced in the proof of Morse's Lemma) as base
points of the two sets respectively. We have

\begin{theorem}\label{lin2}
  Every diffeomorphism $f:\R^n \to \R^n$, $f(0)=0$, such that $d_0f\in
  {\rm GL}^+$ (resp. $d_0f\in {\rm GL}^-$) is diffeotopic to the linear
  isomorphism $I_n$ (resp. $I_{n,1}$), through diffeomorphisms $f_t$
  such that $f_t(0)=0$ for every $t\in \R$.
\end{theorem}

\Dim  If $U$ is a connected open set of some $\R^n$, then
it follows easily from the proof of Proposition \ref{connected}
that  any two points of $U$ can be connected by
 a piecewise smooth path in $U$. In fact 
 it is not hard to see that one can take a globally smooth path
 (use bump functions in order to get a smoothing).
 By using this remark, it is enough to prove that both open sets
 ${\rm GL}^\pm$ are connected. In fact
it is enough to show that ${\rm GL}^+$ is connected. For if
$A\in {\rm GL}^-$, then $I_{n,1}A$ is in ${\rm GL}^+$; if $A_t$ is a
path connecting $I_{n,1}A$ with $I_n$ in ${\rm GL}^+$, then
$I_{n,1}A_t$ is a path connecting $A$ and $I_{n,1}$ in ${\rm
  GL}^-$. 
  
 Let us show first that there is a path $B_t$ in ${\rm GL}^+$
connecting any given $A=B_0$ with some $B=B_1$ which belongs to
$${\rm SO}(n):=\{P\in {\rm GL}(n,\R); \ P^{-1}=P^t, \ \det P = 1\} \ . $$
Let $<*,*>$ be the
positive definite scalar product on $\R^n$ determined by imposing that
the ordered columns of $A$ form an orthnormal basis $\Bb$ of $\R^n$
with respect to such a scalar product. Set $$(*,*)_t= (1-t)<*,*>+
t(*,*)$$ where $(*,*)$ is the standard euclidean scalar product, $t\in
[0,1]$. Then $(*,*)_t$ is a path of positive definite scalar
products. For every $t\in [0,1]$, apply the usual Gram-Schmidt
othogonalization
algorithm to the basis $\Bb$ that produces an othonormal basis $\Bb_t$
for $(*,*)_t$; by considering the ordered vectors of $\Bb_t$ as
columns of a matrix $B_t$, we eventually get a path of matrices such
that $B_0=A$ and $B_1\in {\rm SO}(n)$. It remains to show that every
$B\in {\rm SO}(n)$ can be connected to $I_n$ by a path in ${\rm
  SO}(n)$.  Let us consider $B: \R^n \to \R^n$ as a linear isometry
with respect to $(*,*)$. By linear algebra we know that $\R^n$ can be
decomposed as the orthogonal direct sum of $B$-invariant linear
subspaces $V_i$ of dimension either $1$ or $2$. In the first case the
restriction of $B$ to $V_i$ is the identity; in the second case $B$
acts on $V_i$ as a rotation. Then we are reduced to prove that a
rotation on $\R^2$ can be connected to $I_2$ by a path of rotations,
and this is immediate.

\cvd

\section{Homogeneity}\label{Homogeneity}
We have
\begin{proposition}\label{Homog} Let $p,q\in \R^n$ such that
  $||p-q||=d>0$.  Then for every $\epsilon >0$ there is a
  diffeomorphism $f: \R^n \to \R^n$  such that
\begin{enumerate}
\item $f(p)=q$
\item $f$ is diffeotopic to the identity of $\R^n$ by a diffeotopy
of compact support contained
in $B^n(p,d+\epsilon)$.
\end{enumerate}
\end{proposition}
\Dim  In this proof we use some tools that will be developed in Chapter \ref{TD-CUT-PASTE}.
Without the requirement about the supports the proof is
immediate: set $v= q-p$, then $f_t(x):= x+ tv$, $t\in \R$, $f=f_1$
verify the thesis. Note that for every $x\in \R^n$, $f_t(x)$ is the integral
line defined on the whole
real line of the vector field on $\R^n$ constantly equal to $v$. 
Now we use a bump function to modify this vector field making it with
compact support. Let $d+\epsilon/3 < a <b <d+\epsilon/2$, and
consider the bump function $\gamma=\gamma_{p,a,b}$. Take the smooth vector field
on $\R^n$ defined by $\gamma(x)v$. For every $x \in \R^n$ there is a unique
maximal parametrized integral curve denoted again $f_t(x)$ such that $f_0(x)=x$;
as the field has  compact support also in this case every $f_t(x)$ is defined on the whole
real line. The $f_t(x)$ for $t\in [0,1]$ realizes the required isotopy.

\cvd

The above proposition is a sort of local case of the following
more general result

\begin{theorem}\label{Homog2} Let $U\subset \R^n$ be a
connected open set. Then for every $p\neq q \in U$
there is a diffeotopy $F$  of $U$ 
between $f_0={\rm id}_U$ and $f=f_1$ such that $f(p)=q$,
and $F$ has compact support.
\end{theorem}
\Dim The proof is qualitatively similar to the one of
Proposition \ref{connected}. Being `connected' via a diffeotopy
with compact support as in the statement of the theorem
defines an equivalence relation on $U$. By applying  Proposition \ref{Homog}
on a chart diffeomorphic to $\R^n$ at every $p\in U$
we realize that every equivalence class is an open set, hence there is
only one because $U$ is connected.

\cvd

\chapter{The category of embedded smooth manifolds}\label{TD-DIFF}
Let us begin by widely extending the notions of smooth map and diffeomorphism
 to {\it arbitrary} topological subspaces of some $\R^n$, $n\in \N$.
 
 \medskip
 
 Let $X\subset \R^n$, $Y\subset \R^m$ be arbitrary subspaces. Then
 $f: X \to Y$ is $\Cc^k$, $k\geq 0$, if for every $x\in X$ there exist an open
 neighbourhood $U$ of $x$ in  $\R^n$ and a map $g_U \in \Cc^k(U,\R^m)$
 such that for every $y \in U$, $f(y)=g_U(y)$. Such a map $g_U$
 is called a {\it local $\Cc^k$ extension of $f$ at $x\in X$}.
 
 $f$ is $\Cc^\infty$ (i.e. {\it smooth}) if for every $x\in X$
 there are smooth local extensions of $f$ at $x$.  
 
 A map $f: X \to Y$ is a {\it diffeomorphism} if it is a homeomorphism
 and both $f$ and $f^{-1}$ are smooth maps. 
 
 It is easy to verify  by using the results of Chapter \ref{TD-LOCAL}  
 that  $\Cc^k$ maps, smooth maps and diffeomorphisms are stable
  under composition of maps. By using this very general notion
  of diffeomorphism we can readly define embedded smooth manifolds.
  
  \begin{definition}\label{embedded}{\rm 
  For every $0 \leq k \leq n$, a topological subspace $M\subset \R^n$
  is an {\it embedded smooth $k$-manifold} ($k$ is called the {\it
    dimension} of $M$) if for every $p\in M$, there exist an open
  neighbourhood $W$ of $p$ in $M$, an open set $U$ of $\R^k$ and a
  diffeomorphism $\phi: W \to U$.
  
  Every such a $(W,\phi)$ is called a  {\it chart} of $M$; set $\psi = \phi^{-1}$, 
  then $(U,\psi)$ is called a {\it local parametrization} of $M$.
  The family of all charts
  is called {\it the atlas} $\Aa= \Aa_M$ of $M$. Hence by definition $\Aa$ incorporates 
  an open covering of $M$. {\it An atlas} $\ \Uu \subset \Aa$ of $M$ is any family
  of charts that incorporates an open covering of $M$. 
  
  The {\it category of smooth embedded manifolds} has as {\it objects}
  the embedded smooth manifolds in some $\R^n$, $n\in \N$; the {\it
    morphisms} are the smooth maps between embedded smooth manifolds;
  the diffeomorphisms are the {\it equivalences} in the category.}
  \end{definition}
  
 \section {Basic properties and examples} \label{basic-exa}
  {\rm We are going to list a few basic examples or properties that follow immediately
  from the definitions or are consequence of  results of Chapter \ref{TD-LOCAL}.
  
  $\bullet$ A $0$-manifold in $\R^n$ is a subset of isolated points. It is compact if and
  only if it is finite; otherwise it is countable.
  
  $\bullet$ In order to show that $M\subset \R^m$ is a smooth manifold
  (sometimes we will omit to say ``embedded'') it is enough to exhibit
  an atlas $\Uu$. The whole atlas $\Aa$ is implicitly determined by
  $\Uu$.  For example for every $(W,\phi)\in \Uu$, for every open
  subset $U'\subset U$, the restriction $(U',\phi':= \phi_{|U'})$
  belongs to $\Aa$.
  
 $\bullet$ Every open set $U\subset \R^n$ is a $n$-manifold: the inclusion $j: U\to \R^n$
  forms an atlas of $U$ with only one chart. Hence the category discussed in Chapter
  \ref{TD-LOCAL} is a subcategory of the present category. More generally
  an open set in a $k$-manifold $M$ is also a $k$-manifold.
  
  $\bullet$ Let $U$ be an open set in $\R^n$, $f:U \to \R^m$ a smooth map. Then its
  {\it graph}
  $$G(f):=\{(x,y)\in U \times \R^m; \ y=f(x) \}$$
  is a $n$-smooth manifold embedded in $\R^{n+m}$.
  In fact  $W=G(f)\cap (U\times \R^m)= G(f)$,
  $\phi : W\to U$, $\phi(x, f(x))=x$ form an atlas of $G(f)$ with only one chart;
  the inverse
  parametrization is $\psi: U \to W$, $\psi(x)=(x,f(x))$.
  
  $\bullet$ Let $V$ be a linear (or affine) $k$-subspace of $\R^n$. It is a $k$-manifold,
  in fact the atlas $\Aa$ contains any  linear (affine) isomorphism $L: V \to \R^k$. 
   
$\bullet$ Let $M\subset \R^m$, $N\subset \R^n$ be embedded smooth manifolds.
 Then the product $M\times N$ is a smooth manifold embedded into $\R^{n+m}$,
 and $$\dim (M\times N) = \dim M + \dim N \ . $$ In fact if $(W,\phi)$ is a chart
 of $M$ at $p$, $(W',\phi')$ of $N$ at $q$, then $(W\times W', \phi \times \phi')$
 is a chart of $M\times N$ at $(p,q)$.
 
 $\bullet$ If $(W,\phi), (W',\phi')\in \Aa$ are charts of a $k$-manifold $M$, 
 and $W\cap W'\neq \emptyset$, then 
  $$ \beta_{W,W'}:=\phi' \circ \psi : \tilde U  \to \tilde U'$$
  is a {\it diffeomorphism between open sets of $\R^k$}
(that is $\tilde U =\phi(W\cap W')\subset U$ and $\tilde U' = \phi'(W\cap W')\subset U'$).  
 It is  called indifferently 
  {\it change of charts} or {\it of local parametrizations}
  or also {\it of local coordinates}. 
  
  $\bullet$ If $f:M \to N$ is a smooth map between embedded smooth manifolds,
 $(W,\phi)$ is a chart of $M$, $(W',\phi')$ of $N$ such that $f(W)\subset W'$,
  then
  $$f_{U,U'}:= \phi'\circ  f \circ \psi: U \to U'$$ is a smooth map between open sets
 of euclidean spaces  called {\it a representation
 of $f$ in local coordinates} or shortly {\it a local representation of $f$}.
 
 $\bullet$ {\it The dimension of embedded smooth manifolds is invariant up to
 diffeomorphism.} This follows immediately from the above items and the
 ``invariance of dimension'' already discussed in Chapter \ref{TD-LOCAL}.

 \begin{lemma}\label{path-component}
   (1) An embedded smooth $k$-manifold $M\subset \R^n$ is connected
   if and only if it is path connected.

   (2) Every path connected component
   of $M$ is a $k$-manifold. $M$ is the disjoint union of its path
   connected (equivalently connected) components.
 \end{lemma}

 \Dim It is a general topological fact that a path connected space is connected.
 For the other implication we can repeat the argument already used for
 the open sets in $\R^k$. In fact by using a chart around any point $p\in M$
 we can argue that the path connected component
 of $p$ is open in $M$, hence there is only one if $M$ is connected.
 This proves (1) and also (2) indeed.

 \cvd

\medskip  
 
The definition of embedded smooth manifold $M\subset \R^n$
implies some strong {\it local} constraint on the relative
configuration of the pair $(\R^n,M)$. We have

\begin{lemma}\label{relative} Let $M\subset \R^n$ be an embedded 
 smooth $k$-manifold; $p\in M$. Then there exist a chart $(\Omega, \beta)$
 of $\R^n$, $p\in \Omega$, such that $(\Omega\cap M, \beta_|)$
 is chart of $M$ and moreover 
 $$\beta(\Omega, \Omega \cap M,p)= (B^n(0,1), B^n(0,1)\cap \R^k,0)$$
 (where $\R^k \subset \R^k \times \R^{n-k}=\R^n$ as usual).
 Such a  $\beta$ is called a {\rm relative normal chart} of the pair
   $(\R^n,M)$.
 \end{lemma}

\Dim It follows immediately from the definition of embedded manifold
 that there exist an open neighbourhood $\Omega$ of $p$ in $\R^n$, an
 open set $U$ of $\R^k$, and an {\it injective immersion} $\psi: U \to
 \Omega$, such that $\psi(U)=\Omega \cap M := W$.  By Theorem
 \ref{immersion} on local normal form of immersions, possibly by
 shrinking $\Omega$, there is chart $(\Omega,\beta)$ of $\R^n$ that
 verifies the statement of the Lemma.
 
 \cvd 
 
 The above argument can be somehow reversed. 
 
 \begin{lemma}\label{inj-imm} Let $U$ be an open set of $\R^k$, $\psi: U \to \R^n$
be an injective immersion such that $\psi:U \to \psi(U)$ is a {\it homeomorphism}. Then
$M=\psi (U)$ is a smooth manifold embedded in $\R^n$, and $\psi: U\to M$
is a (global) smooth parametrization of $M$.
 \end{lemma} 
 \Dim By using again Theorem \ref{immersion} and the fact that $f$ is
a homeomorphism onto its image, we readly see that at every $p\in M$
one can  find relative normal charts of $(\R^n,M)$, and eventually
$\psi$ is a diffeomorphism onto $M$.

\cvd

The condition that $\psi$ is a homeomorphism onto its image is {\it necessary} 
as it is shown by the following example:

\begin{example}\label{no-emb}{\rm Consider the smooth map
$$E: \R^2 \to \R^2 \times \R^2, \ E(x,y)= (\cos(2\pi x),\sin(2\pi x),
    \cos(2\pi y),\sin(2\pi y)) \ . $$ For every $a\in \R$, $a\neq 0$,
    consider the map
    $$f: \R \to \R^2\times \R^2, \ f(x)= E(x,ax) \ . $$
    This is an injective
    immersion but if $a$ is {\it not a rational number}, then
    it is not a homeomorphism onto its image in
    $S^1\times S^1$.  In fact one can verify that $f(\R)$
    is {\it dense} in $S^1\times S^1$ (every non empty open set of
    $S^1\times S^1$ intersects $f(\R)$), hence $f(\R)$ is not an
    embedded manifold in $\R^2 \times \R^2$.}
\end{example}

\medskip

{\bf Submanifolds.} If $Y\subset M$ are embedded smooth manifolds in $\R^n$ , we say that
$Y$ is a {\it submanifold} of $M$.  In particular both $Y$ and $M$ are
submanifolds of $\R^n$. By extending the argument of Lemma
\ref{relative} (details are left as an exercise) we can prove

\begin{lemma}\label{relative2}
  Let $Y$ be a submanifold of $M\subset \R^n$, of dimension $k$ and
  $m$ respectively. Let $p\in Y$. Then there exist relative normal
  charts (for triples)
$$\beta: (\Omega,\Omega\cap M, \Omega \cap  Y,p) \to
(B^n(0,1), B^n(0,1)\cap \R^m, B^n(0,1)\cap \R^k,0) $$
where as usual we consider $\R^m= \R^k \times \R^{m-k}$, $\R^n = \R^m \times \R^{n-m}$.
\end{lemma}

\medskip

By using the immersions, we have indicated above a way to get embedded
manifolds (endowed with global smooth parametrizations).  Now we show
how embedded manifolds can be defined {\it implicitly}.

\begin{lemma}\label{surj-summ}
  If $f: U \to W$ is a surjective smooth summersion between open sets
  of euclidean spaces, $U\subset \R^n$, $W\subset \R^m$.  Then for
  every $q\in W$, $M=f^{-1}(q)$ is an embedded smooth manifold in
  $\R^n$ and $\dim M= n-m$.
\end{lemma}
\Dim Being an embedded manifold is a local property. Hence the lemma
is an immediate consequence of Theorem \ref{summersion} on
local normal form of summersions or (equivalently)
of the implicit function Theorem
\ref{implicit-funct}.
  
 \cvd
  }

\begin{remark}\label{strange-sub}{\rm In spite of the existence of relative normal
charts at every point of a submanifold, the relative position of two
submanifolds of some $\R^n$ can look stranger than one could expect.
This is mainly due to the fact that submanifolds are not necessarily closed
subsets. Consider for example the map $f: (0, +\infty)\to \C\sim \R^2$
$$f(x)=\frac{x}{1+x}e^{ix} \ . $$
This is an immersion and a homeomorphism onto its image say  $N$.
Then the unitary circle $S^1$ and $N$ are disjoint $1$-submanifolds of $\R^2$.
Nevertheless, two points $p\in S^1$ and $q\in N$ cannot be separated
by normal charts of $S^1$ and $N$ at $p$ and $q$ respectively. In other words
$N\cup S^1$ is not an embedded submanifold.}
\end{remark}

 \begin{example}(Spheres) \label{sfere}{\rm 
 Let us show, in several ways, that the unitary sphere $S^n\subset
 \R^{n+1}$, $n\in \N$, is an embedded smooth $n$-manifold. Let
 $\R^{n+1}=\R^n \times \R$.  Let $W^+= S^n \setminus \{e_{n+1}\}$,
 $\phi_+: W^+\to \R^n$ be the {\it stereographic projection} with
 center $e_{n+1}$. It is defined geometrically by $\phi_+(x)=
 r(x,e_{n+1})\cap \R^n$ where $r(x,e_{n+1})$ is the straight line
 passing through the two points. Analytically we have
 $$\phi_{+}(x)= \frac{1}{1-x_n}(x_1,\dots,x_{n-1}) \ . $$ 
 This is a diffeomorphism onto $\R^n$ with inverse given by 
 $$\psi_{+}(y)= \left(\frac{2y}{1+||y||^2},\frac{||y||^2-1}{||y||^2+1}\right) \ . $$
Then $(W^+,\phi_+)$ is a chart of $S^n$ at every points different from $e_{n+1}$.
By using the similar projection with center $-e_{n+1}$, we get a chart $(W^-,\phi_-)$
which misses only $-e_{n+1}$. Hence $\{(W^\pm,\phi_\pm\}$ is an atlas of $S^n$
(formed by two charts).

For every $p\in S^n$, let $p^\perp$ the subspace of $\R^{n+1}$
orthogonal to $p$. Then by using the projection of $S^n \setminus
\{p\}$ onto $p^\perp$
with center $p$ (followed by any linear chart of
$p^\perp$ onto $\R^n$) then we obtain other charts of the atlas
$\Aa_{S^n}$.

Further charts are obtained as graphs of functions defined on the
unitary open disk of $p^\perp$ with center $p$.  The basic example for
$p=e_{n+1}$ is the function $h:B^n\to \R$,
$$h(x)=\sqrt{1-\sum_{i=1}^{n-1} x_i^2} \ . $$

$S^n= f^{-1}(1)$, where $f:\R^{n+1}\setminus \{0\} \to \R$, $f(x)=||x||^2$.
As $df_x=(2x_1,\dots, 2x_{n+1})$ then $f$ is a summersion and this implies
again (implicitly) that $S^{n}$ is a $n$-manifold by Lemma \ref{surj-summ}.

$S^n$ is a {\it compact}
manifold.  In fact it is closed because $S^n= f^{-1}(1)$ as above;
obviously it is bounded.

$S^n$ is path connected: given $x\neq y \in S^n$, let $P$ be the $2$-plane
spanned by these two vectors. Then $P\cap S^n$ is a maximal circle,
$x,y$ separate it into two arcs both connecting $x$ and $y$.}
  \end{example}

  \medskip
  
  Important examples of embedded
  smooth manifolds  (widely generalizing the spheres) 
  are discussed in Chapter \ref{TD-STIEF-GRASS}.

\section{The embedded tangent functor}\label{emb-tang}
Let us fix a setting we will refer to all along the rest of this Chapter.

\begin{itemize}
\item $M\subset \R^h$ is an embedded smooth manifold
of dimension $m$, $p\in M$; $N\subset \R^k$ is an embedded
smooth manifold of dimension $n$, $q\in N$; 
\item $f:M \to N$ is a smooth map, $f(p)=q$.
\item $\phi:W\to U\subset \R^h$ is a chart of $M$ at $p$, $\phi(p)=a$,
with inverse local parametrization $\psi:U\to W \subset M$.
\item $  f_{U,U'}: U\to U'$ is a representation of $f$ in local coordinates at $p$; 
recall  that this is obtained as follow:
we take a local chart of $M$ at $p$ for semplicity still denoted $(W,\phi)$,
and a local chart $(W',\phi')$ of $N$ at $q$, $\phi'(q)=b$, such that 
$f(W)\subset W'$; then 
$$  f_{U,U'}= \phi' \circ f \circ \psi: U\to U'$$
($U$ and $U'$ being open set of $\R^h$ and $\R^k$ respectively).
\item Possibly by shrinking $W$ we can also assume that there are an open
neighbourhood $\Omega$ of $p$ in $\R^h$ such that  $\Omega \cap M=W$,
a local smooth extension $\Phi: \Omega \to \R^m$ of $\phi$ and a local 
smooth extension $F: \Omega \to \R^k$ of $f$. 
\end{itemize}

\medskip

The facts collected in the following Lemma
are easy consequences of 
the very definitions and of the results
of Chapter \ref{TD-LOCAL}. The reader would like to make the useful
exercise to fill the details.

\begin{lemma}\label{emb-diff}

(1)  The differential $d_a \psi$ is injective so it is a linear isomorphism onto
its image $d_a \psi (\R^m)$, $\R^m = T_a U$,  which is a $m$-linear subspace of $\R^h=T_p\R^h$ . 
This image  does not depend
on the choice of the local parametrization $\psi$ of $M$ at $p$.
Hence 
$$T_pM= d_a \psi (\R^m)$$ 
is well defined and is called the
{\rm tangent space to $M$ at the point $p$}.

\medskip

(2) The restriction of the differential $d_p \Phi$ to $T_pM$
is the inverse isomorphism $(d_a \psi)^{-1}$.
Hence it does not depend on the choice of the
local extension $\Phi$ of $\phi$, and
$$d_p\phi := d_p \Phi_{|T_pM}$$
is a well defined linear isomorphism
$$d_p \phi: T_pM \to T_aU \ . $$

\medskip

(3) The restriction of $d_p F$ to $T_pM$ does not depend on the choice
of the local extension of $f$ and is valued in $T_qN$.
Hence it is well defined
$$ d_pf := d_pF_{|T_pM}$$
it is a linear map
$$d_pf: T_pM \to T_qN$$
and is called the {\rm differential of $f$ at $p$}.
We have
$$ d_a  f_{U,U'} = d_q\phi' \circ d_pf \circ d_a\psi: T_aU \to T_bU'$$
and this is the {\rm representation in local coordinates} of $d_pf$.
In particular this applies when  $M=W$, and $f=\phi': W\to U' \subset \R^m$
is another chart of $M$ at $p$.

\medskip

(4) If $g\circ f$ is a compostion of smooth maps between embedded smooth
manifolds, $f(p)=q$, then
$$ d_p(g\circ f)= d_qg \circ d_p f \ . $$
If $f$ is a diffeomorphism, then $d_pf$ is a linear isomorphism and
$d_q f^{-1}= (d_p f)^{-1}$.  If $f={\rm id}$, then $d_p f = {\rm id}_{T_pM}$. 

\medskip

(5) If $M=G(g)$ is the graph of a smooth map $g:U\to \R^s$ defined
on an open set $U\subset \R^m$, then 
$$ T_{(x,g(x))}M=G(d_x g) \ . $$
\medskip

(6) If $M= g^{-1}(q)$, where $g: \Omega \to \R^s$ is a summersion, $p\in M$,
then
$$ T_pM = \ker d_p g \ . $$
\end{lemma}

\cvd

\medskip

Set
$$ T(M)=\{(x,v)\in \R^h \times \R^h; \ x\in M, \ v \in T_xM\} \ . $$
The restriction of the projection of $\R^h\times \R^h$ onto the first
factor $\R^h$ defines a smooth projection
$$\pi_M: T(M)\to M \ . $$

\begin{example} \label{TSn}{\rm 
$$T(S^n)= \{(x,v)\in \R^{n+1} \times \R^{n+1}; \ x\in S^{n}, \ v \in x^\perp\}$$
Check it!}
\end{example}

\medskip

As a set $T(M)= \cup_{x\in M} T_xM$. 
Note that for every open set $W\subset M$, $T(W)$ coincides with
$\pi_M^{-1}(W)$, it is naturally included in
$T(M)$ as an open set, and $\pi_W= (\pi_M)_{|T(W)}$.

We are going to show that
\medskip

{\it $T(M)$ is an embedded smooth manifold of dimension $2m$, of a special
  nature indeed}.
\medskip

Every chart $\phi: W\to U\subset \R^m$ of $M$ can be enhanced
to a chart 
$$T\phi: T(W)\to T(U)=U\times \R^m, \ T\phi (x,v):= (\phi(x), d_x\phi (v)) \ . $$
The inverse parametrization is
$$ T\psi: U\times \R^m \to T(W), \ T\psi(y,w)=(\psi(y), d_y \psi (w)) \ . $$
If $\pi_U$ is the natural projection onto $U$, it is immediate that
the following diagram denoted $[\psi,T\psi]$ commutes

$$ \begin{array}[c]{ccc}
U\times \R^m&\stackrel{T\psi}{\rightarrow}&T(W)\\
\downarrow\scriptstyle{\pi_U}&&\downarrow\scriptstyle{\pi_W}\\
U&\stackrel{\psi}{\rightarrow}&W\end{array}$$

We say that $\pi_M: T(M)\to M$ is {\it locally trivial (a product)} over $W$
and that the above diagram is a {\it local trivialization}. 
By varying the chart $(W,\phi)$ in the atlas $\Aa$ of $M$ we get an atlas
$$T\Aa = \{ (T(W), T\phi)\}$$ of $T(M)$. The local coordinates
for $T\Aa$ changes in a special way as they are of the form

$$ T\beta:= T\phi' \circ T\psi: \tilde U\times \R^m \to \tilde U'\times \R^m$$

$$ T\beta(x,v)= (\phi'\circ \psi(x), d_x(\phi'\circ \psi)(v))
=(\beta(x), d_x\beta(v)) \ . $$
Hence, for every $x$ varying in $M$, it is a linear isomorphism on the second argument
which ``varies smoothly'' with the point $x$. This means that
\medskip

{\it The intrinsic linear structure of every fibre $T_xM=\pi_M^{-1}(x)$
of the projection $\pi_M$ is respected by the changes of coordinates
for the atlas $T\Aa$.} 
\medskip

We can encode the same information by lifting $T\Aa$ at the level
of the open covering $\{W\}$ of $M$; that is we have the locally trivializing
commutative diagrams

$$ \begin{array}[c]{ccc}
W\times \R^m&\stackrel{\tilde T\psi}{\rightarrow}&T(W)\\
\downarrow\scriptstyle{\pi_W}&&\downarrow\scriptstyle{\pi_W}\\
W&\stackrel{{\rm id}_W}{\rightarrow}&W\end{array}$$
where
$$ \tilde T\psi = T\psi \circ (\phi, {\rm id}_{\R^m}) \  . $$
Any change of local trivialization for $\tilde T\Aa$
is of the form
$$ \tilde T\beta: (W\cap W')\times \R^m \to (W\cap W')\times \R^m, \ (x,v)\to (x, d_x\beta(v)) \ . $$
We summarize all these facts by saying that
$$\pi_M: T(M)\to M$$ 
is the {\it tangent vector bundle} of the embedded 
smooth manifold $M$ and that $T\Aa$ (actually and equivalently 
$\tilde T\Aa$) is its {\it vector bundle atlas}.

\medskip

In section \ref{emb-fib-bundle}
we will formalize these notions in a more general setting.

\medskip

Now we extend the definition of the {\it tangent map}
already considered in Chapter \ref{TD-LOCAL} in the case of open sets in some
$\R^n$. Let $f:M\to N$ be our smooth map between embedded smooth manifolds, then set:

$$ Tf: T(M)\to T(N), \ Tf(x,v)= (f(x), d_xf(v)) \ . $$

Note that the defining inclusion $T(M)\subset \R^h\times \R^h = T(\R^h)$
is nothing else than $Tj$, $j:M\to \R^h$ being the inclusion.
Clearly the following diagram, denoted $[f,Tf]$, commutes

$$ \begin{array}[c]{ccc}
T(M)&\stackrel{Tf}{\rightarrow}&T(N)\\
\downarrow\scriptstyle{\pi_M}&&\downarrow\scriptstyle{\pi_N}\\
M&\stackrel{f}{\rightarrow}&N\end{array}$$
that is $Tf$ sends every fibre $T_xM$ linearly to the fibre $T_{f(x)}N$,
by af `smooth field' of linear maps. 

If $g\circ f$ is a composition of smooth maps between embedded
smooth manifolds, then 
$$T(g\circ f)= Tg \circ Tf$$ 
$$T{\rm id_M}={\rm id_{T(M)}}$$
if $f$ is a diffeomorphism, then $Tf$ is a diffeomorphism and 
$$Tf^{-1}=(Tf)^{-1} \ . $$  
All  verifications are local and follows immediately from
Lemma \ref{emb-diff} and the properties of the tangent map
in the category of open sets in euclidean spaces.

We can summarize these considerations a follows:

\medskip

{\it The {\rm tangent category} of the category of embedded
smooth manifolds has as {\rm objects} the tangent vector bundles
of embedded smooth manifolds and as {\rm morphisms} the tangent
maps of smooth maps between embedded smooth manifolds. 
Then 
$$ M \  \Rightarrow \ \pi_M: T(M)\to M$$
$$ f:M\to N \ \Rightarrow \ [f,Tf] $$
define a {\rm covariant functor} from the category of embedded
smooth manifolds to its tangent category.}

\section{Immersions, summersions, embeddings, Monge charts}\label{imm-summ-emb}
The notions of immersion and summersion extend immediately to map between
embedded smooth manifolds: $f:M\to N$ is an {\it immersion} (resp.
{\it summersion})
if for every $x\in M$, $d_xf$ is {\it injective} ({\it surjective}).
We say that $f:M\to N$ is an {\it embedding} if $f$ is a diffeomorphism
onto its image (in particular the inclusion $M\subset \R^h$ is an embedding).
The proof of the following proposition is of local nature and follows
from Lemmas \ref{inj-imm} and \ref{surj-summ}.

\begin{proposition}\label{global-imm-summ}
(1) Let $f: M\to N$ be a surjective summersion; then for every $q\in N$,
  $Y=f^{-1}(q)$ is a submanifold of $M$, $\dim Y = \dim M - \dim N$.

  (2) If $f: M\to N$ is an embedding then $f(M)$ is a submanifold
  of $N$.

  (3) $f:M \to N$ is an embedding if and only if $f$ is an immersion
  and  a homeomorphism onto its image.

  (4) If $f:M\to N$ is both an immersion and a summersion, then it is a local
  diffeomorphism.

\end{proposition}

\cvd

\medskip

We have seen in example \ref{sfere} a distinguished local graph
chart of $S^n$. Here we show that such a kind of charts exists
for every embedded smooth $m$-manifold $M\subset \R^h$ at every point.
For every multi-index $J=(j_1,\dots,j_m)$, $|J|=m$,
let $J'$, $|J'|=h-m$ be its complementary multi-index. Denote by
$\R_J$ the subspace of $\R^h$ generated by $(e_{j_1},\dots,e_{j_m})$;
hence we have the orthogonal direct sum decomposition
$\R^h=\R_J\oplus \R_{J'}$ and the orthogonal projection onto $\R_J$,
$\pi_J(x)=(x_{j_1},\dots,x_{j_m})$. For every $p\in M$, denote by
$\pi_{J,p}: \R^h \to \R_J$ the composition of the translation
$x\to x-p$, followed by $\pi_J$. Denote by $\phi_{J,p}$ the restriction
of $\pi_{J,p}$ to (any suitable subset of) $M$. We have

\begin{proposition}{\rm (Monge charts)}
  For every embedded smooth $m$-manifold $M\subset \R^h$,
  for every $p\in M$, there exist $J$, $|J|=m$, and an open neighbourhood
  $W$ of $p$ in $M$ such that $(W,\phi_{J,p})$ is a chart of $M$ at $p$.
  The inverse local parametrization is of the form $\psi_{J,p}:U \to W$,
  $U\subset \R_J$, $\psi_{J,p}(y)=(y, f_{J,p}(y))$ (by using the above decomposition
  $\R^h= \R_J \oplus \R_{J'}$). Hence at every point $p$, $M$ is locally a graph
  of a smooth function defined on some $\R_J$.
  \end{proposition}
  
\Dim By elementary linear algebra, there exist $J$ such that the
restriction of $\pi_J$ to $T_pM$ is a linear isomorphism onto $\R_J$.
As $d_p \phi_{J,p}$ coincides with such a restriction, then $\phi_{J,p}$
is a local diffeomorphism.

\cvd

\section{ Topologies on spaces of smooth maps}\label{weak-top}
Let $M\subset \R^h$ $N\subset \R^k$ be smooth manifolds as usual.  We
define the {\it weak topology} on every set $\Cc^r(M,N)$, $r\geq 0$, the
topological spaces $\Ee^r(M,N)$ (subspaces of $\Cc^r(M,N)$ formed by
the smooth maps) and the space $\Ee(M,N)$ that is $\Cc^\infty(M,N)$
equipped with the union of the $\Ee^r$ topologies.  This extends the
case of open sets treated in Chapter \ref{TD-LOCAL} which is actually
used in order to do it. There are two equivalent ways; both
determine a basis of open neighbourhoods of every element in the
pertinent map space. We leave to the reader the verification that the two
topologies defined in these ways actually are the same one.
\medskip

(1)  For every $f\in \Cc^r(M,N)$ we consider neighbourhoods of the
following form
$$\Uu_r(f, \hat f, K, \epsilon)$$
where
\begin{itemize}
\item $\hat f: \Omega \to \R^k$ is a local $\Cc^r$ extension of
$f_{|W}:W\to N$, $W=\Omega\cap M$, $\Omega \subset \R^h$
being open;
\item $K\subset W$ is a compact set;
\item $\epsilon >0$.
\end{itemize}

Then $g\in \Cc^r(M,N)$ belongs to  $\Uu_r(f, \hat f, K, \epsilon)$
if and only if there exists a $\Cc^r$ extension $\hat g: \Omega \to \R^k$
of $g_{|W}$ such that $\hat g \in \Uu_r(\hat f, K, \epsilon) \subset \Cc^r(\Omega,\R^k)$.

\medskip

(2) For every $f\in \Cc^r(M,N)$ we consider neighbourhoods of the
following form
$$\Uu_r(f, f_{U,U'}, K, \epsilon)$$
where
\begin{itemize}
\item $f_{U,U'}: U' \to U$ is a  (necessarily $\Cc^r$)
 representation of $f$ in local coordinates ($U\subset \R^m$, 
$U'\subset \R^n$ being open sets);
\item $K \subset U$ is a compact set;
\item $\epsilon >0$.
\end{itemize}

Then $g\in \Cc^r(M,N)$ belongs to  $$\Uu_r(f, f_{U,U'}, K, \epsilon)$$
if and only if it admits a local representation (over the same open sets
$U,U'$) $g_{U,U'}:U\to U'$ such that $g_{U,U'}\in \Uu_r(f_{U,U'}, K,\epsilon)\subset
\Cc^k(U,U')$. 

\section{Homotopy, isotopy, diffeotopy, homogeneity}\label{H-I-D}
These notions already introduced in Chapter \ref{TD-LOCAL} within the smooth category
of open sets, extend {\it verbatim} to embedded smooth manifolds.
They correspond to continuous paths in appropriate map spaces and bring equivalence
relations along.

The proof of the homogeneity Theorem \ref{Homog2} is essentially of local nature.
and extends straightforwardly. 

\begin{theorem}\label{Homog3} Let $N$ be a connected embedded smooth manifold.
Let $p,q\in N$. Then there is a  diffeotopy with compact support between $f_0={\rm id}_N$
and $f=f_1$ such that $f(p)=q$.  
\end{theorem}

\cvd

\section {Embedded fibre bundles}\label{emb-fib-bundle}
The tangent vector bundle is a first fundamental  example
of the general notion of {\it fibre bundle}. We will encounter several instances of all along this text.
Chapter \ref{TD-EMB-VB} will develop this topic. Here we state the basic facts. 

An {\it embedded smooth  fibre bundle} with {\it base space} $X$, {\it total space} $E$ 
and {\it fibre} $F$, is a surjective summersion $f: E \to X$ between embedded smooth manifolds
such that every fibre
$f^{-1}(q)$, $q\in X$, is a submanifold of $E$ diffeomorphic to a given manifold $F$, and which is {\it locally
trivial(izable)} at every point $q$ of $X$. This means that for every $q\in X$, there is a open
neighbourhood $\Omega$ in $X$ and trivializing commutative diagram 
of the form

$$ \begin{array}[c]{ccc}
\Omega \times F &\stackrel{ \Phi }{\rightarrow}& \tilde \Omega\\
\downarrow\scriptstyle{\pi_\Omega}&&\downarrow\scriptstyle { f_|}\\
\Omega&\stackrel{{\rm id}_\Omega}{\rightarrow}&\Omega \end{array}$$ 
where $\tilde \Omega:= f^{-1}(\Omega)$, $\Phi$ is a diffeomorphism (with inverse say $\Psi)$.
If $E=X\times F$ and $f=\pi_X$ is the natural projection then it is a
{\it trivial} (also called `product') fibre bundle.
The family of all {\it local trivializations} as above form the
{\it maximal fibred  atlas} $\Ff$ of the fibre bundle. A fibred
atlas is a subfamily of $\Ff$ such that the $\Omega$'s form an open
covering of $X$, hence the $\tilde \Omega$'s of $E$.
Every fibred atlas is contained in a unique maximal one, so
it is enough to give a fibred atlas in order to determine a fibre
bundle structure. Every change of local trivialization is of the
form
$$ \Phi' \circ \Psi: (\Omega \cap \Omega') \times F \to     (\Omega \cap \Omega') \times F$$
$$(p, y)\to (p, \rho(p)(y))$$
where $\rho(p)$ belongs to the group ${\rm Aut}(F)$ of the smooth automorphisms of the fibre $F$.
\medskip

$\bullet$ In many cases the fibre $F$ has an additional structure 
which is preserved by a subgroup $G$ of  ${\rm Aut}(F)$
 (for example $F$ is a linear subspace of some $\R^n$ and $G= {\rm GL}(F)$);
if the $\rho(p)$'s as above belong to $G$ then we have a $G$-{\it fibre bundle}
({\it vector bundle}, $\dots$). 

\medskip

$\bullet$ A particular case is when $\dim F = 0$. In such a case a fibration
$f: E \to X$ is also called a {\it covering map} (of {\it degree} $d$ if $F$
is compact hence finite, and $d=|F|$). For every local trivialization,
the restriction of $f$ to every connected component of $\tilde \Omega$
is a diffeomorphism onto $\Omega$, provided that $\Omega$ is connected.
\medskip


$\bullet$ A {\it fibred map} between fibre bundles is a commutative diagram of smooth maps
$[g,\tilde g]$ of the form
$$ \begin{array}[c]{ccc}
E &\stackrel{ \tilde g }{\rightarrow}& E'\\
\downarrow\scriptstyle{f}&&\downarrow\scriptstyle { f'}\\
X&\stackrel{g}{\rightarrow}&X' \end{array}$$
so that every fibre $E_x\sim F$ is mapped to the fibre $E'_{g(x)}\sim F'$.
It is a {\it fibred diffeomorphism} if both $g$ and $\tilde g$ are diffeomorphisms.
In such a case $F=F'$.
The diagrams $[f,Tf]$ of the tangent functor are basic examples of fibred maps.\medskip

\medskip

{\bf Fibred equivalences.} Consider the set $\Ff(X,F)$ of fibred bundles over a given base space
$X$, with given fibre $F$.  There are two natural equivalence relations on $\Ff(X,F)$:

(1) The {\it full equivalence}: it is generated by the fibred diffeomeorphisms $[g,\tilde g]$
such that $g$ belongs to the group ${\rm Aut}(X)$ of smooth automorphisms
of $X$. 

(2) The {\it strict equivalence} (often we will omit to say ``strict''): it is generated by 
the fibred diffeomorphism of the form $[{\rm id}_X,\tilde g]$.

\medskip

This specializes directly to the case of $G$-fibred bundles.

\section{Tensor functors}\label{other-functor}
Let us recall some elementary facts of finite dimensional {\it multi-linear algebra}.
Every finite dimensional
real vector space $V$ has an infinite family of associated {\it tensor spaces}
$T^p_q(V)$, $p,q\in \N$ - also denoted $(V)^{\otimes^p}\otimes (V^*)^{\otimes^q}$ -
formed by the multilinear forms
$$ \alpha: \prod_{i=1}^p V^* \times \prod_{j=1}^q V\to \R \ . $$
Hence the {\it dual space} $V^*=T^0_1(V)$, while $V$  is ``equal'' to 
$T^1_0(V)$ via the canonical identification of $V$ with its {\it bidual space} $(V^*)^*$.
If $\dim V= m$, then
$$\dim T^p_q(V)= m^{pq} \ . $$
Moreover, to every
basis $\Bb$ of $V$, we can associate in a canonical way a basis
$\Bb^p_q$ of $T^p_q(V)$, we can say that the basis $\Bb$ ``propagates'' to every tensor space. 
The linear group GL$(V)$ acts on 
$T^p_q(V)$ by
$$ (g,\alpha)\to g(\alpha)$$
$$g(\alpha)(w^1,\dots, w^p, v_1,\dots, v_q)=\alpha((g^t)^{-1}(w^1),\dots,
(g^t)^{-1}(w^p), g(v_1),\dots, g(v_q)) \ . $$
By applying this to $V=\R^m$ (endowed with the canonical basis $\Cc$) and
to $T^p_q(\R^m)$ (with the canonical basis
$\Cc^p_q$) we get a homomorphism of group (that is a {\it representation})
$$ \rho_{p,q}:{\rm GL}(m,\R)\to {\rm GL}(T^p_q(\R^m))\sim {\rm GL}(m^{pq}, \R)$$
which is an explicit {\it regular rational} map. 
The basic example is 
$$\rho_{0,1}(A)=(A^t)^{-1} \ . $$
As another example: $T^0_2(\R^m)$ can be identified with $M(m,\R)$
by associating to every matrix $B$ the form 
$$ (v,w)\to v^tBw \ . $$
Then 
$$ \rho_{0,2}(P)(B)=P^tBP \ . $$ 

In some case it is interesting to consider suitable subspaces $W$ of
$T^p_q(V)$, $\dim W=w$ say, which are invariant for the action of GL$(V)$
and are endowed as well with a basis $\Bb_W$  canonically associated to $\Bb$.
By appling this to $V=\R^m$, this gives rise to other representations
$$ \rho_W: GL(m,\R)\to {\rm GL}(W)\sim {\rm GL}(w,\R) \ . $$
For example consider the subspace $W=S^2_0(V)\subset T^0_2(V)$ of
{\it symmetric bilinear form on $V\times V$} (i.e. the space of {\it scalar
products} on $V$). In this case the representation $\rho_W$ is just
the ``restriction'' of $\rho_{0,2}$. Another example  
is the subspace  $\Lambda^0_q(V)\subset T^0_q(V)$
of {\it alternating multilinear forms}. 
As a particular case $\Lambda^0_m(\R^m)$ is $1$-dimensional
with canonical basis 
$$\det: M(m,\R)\to \R, \ X \to \det(X)$$
considered as
$m$-linear function of the columns of $X$.  
This gives rise to the representation 
$$ \delta_m: {\rm GL}(m,\R)\to {\rm GL}(1,\R), \ \delta_m(P)= \det P \  . $$

\medskip

We are going to show that for every
embedded smooth $m$-manifold $M\subset \R^h$, the tangent vector bundle 
$$ \pi=\pi_{M}: T(M)\to M$$
has naturally associated  a family of further embedded  vector bundles over $M$
$$\pi_{p,q}=\pi_{p,q,M}: T^p_q(M)\to M $$
such that for every $x\in M$, $\pi_{p,q}^{-1}(x)= T^p_q(T_xM)$, and
clearly $T(M)=T^1_0(M)$.

Let us start with the {\it cotangent bundle}
$$T^*(M):=T^0_1(M) \ . $$
Recall that $(\R^h)^*= M(h,1,\R)\sim \R^h$. For every $x\in M$,
denote by $V_x$
the orthogonal complement of $T_xM$ in $\R^h$,
so that we have the orthogonal direct sum decomposition
$\R^h= T_xM \oplus V_x$. For every functional $\gamma\in T^*_xM$,
extend it to a functional on the whole of $\R^h$  by imposing
that $\gamma(v+w)=\gamma(v)$ for every $w\in V_x$. In this way we have
identified $T^*_xM$ as a linear subspace of $(\R^h)^*$.
For every open subset $U\subset \R^m$, the 
cotangent bundle is the product bundle $U\times (\R^m)^*\to U$.
By copying the definition of the tangent bundle, set
$$ T^*(M)=\{ (x,\gamma) \in \R^h\times (\R^h)^*; \  x\in M, \  \gamma \in T^*_xM\}$$
endowed with the natural projection 
$$\pi^*_M: T^*(M)\to M \ . $$ For every open set $W\in M$,
$T^*(W)=(\pi^*_M)^{-1}(W)$, it is an open set of $T^*(M)$ and
$\pi^*_W$ is the restriction of $\pi^*_M$.  We define the vector
bundle atlas $T^*\Aa$ of $T^*(M)$; for every chart $(W,\phi)$ of $M$
with inverse local parametrization $\psi$, set $(T^*(W), T^*\phi)$,
$$T^*\phi: T^*(W)\to U\times (\R^m)^*,
\ (x, \gamma)\to (\phi(x), \gamma \circ d_{\psi(x)}) \ . $$
The changes of local (fibred) coordinates for $T^*\Aa$
are of the form
$$ T^*\beta: \tilde U \times (\R^m)^* \to \tilde U'\times (\R^m)^*$$
$$T^*\beta(x,\gamma)=(\beta(x), \rho_{0,1}(d_x \beta)(\gamma)) \ . $$

If $f:M\to N$ is a {\it diffeomorphism}, we define
$$T^*f: T^*(N)\to T^*(M), \ (y,\gamma) \to (f^{-1}(y), \gamma \circ d_{f^{-1}(y)})$$  
Then 
$$ M  \ \Rightarrow \ \pi^*_M: T^*(M)\to M$$
$$ f:M\to N \ \Rightarrow \ [f,T^*f]$$
define the {\it contravariant cotangent functor} from the {\it restricted category of embedded smooth manifolds}
to its {\it cotangent category} (`restricted' means that only the diffeomorphisms are allowed as morphisms).
To get a covariant version of the same functor it is enough to replace $T^*f$ with $T^*(f^{-1})$.
\medskip

{\it (The $T^0_2$ functor)} For every $x\in M$, identify $T^0_2(T_xM)$ as a subspace of $T^0_2(\R^h)$  
by  extending every bilinear form $\alpha$ over $T_xM$  to a bilinear form over the whole of $\R^h$
by imposing that for every $v+w,u+z \in T_xM \oplus V_x$, $\alpha(v+w,u+z)=\alpha(v,u)$. 
By the usual scheme, set
$$ T^0_2(M)=\{(x,\alpha)\in \R^h\times T^0_2(\R^m) ; \ x\in M, \ \alpha \in T^0_2 (T_xM)  \}  $$
$$ \pi_{0,2, M}: T^0_2(M)\to M $$
the natural projection.
We have the vector bundle atlas $T^0_2\Aa$  obtained by associating to every
chart $(W,\phi)$ of $M$, with inverse local parametrization $\psi$, the chart $(T^0_2(W), T^0_2\phi)$
$$ T^0_2\phi(x,\alpha)=(\phi(x), \alpha \circ (d_x \psi \times d_x \psi)) \ . $$
The changes of coordinates for $T^0_2\Aa$ are of the form
$$T^*\beta(x,\alpha)=(\beta(x), \rho_{0,2}(d_x \beta)(\alpha)) \ . $$
If $f: M\to N$ is a diffeomorphism, we can define 
$$T^0_2 f: T^0_2(N) \to T^0_2(M)$$
$$ T^0_2f (y,\alpha)= (f^{-1}(y), \alpha \circ (d_{f^{-1}(y)}f \times d_{f^{-1}(y)}f) \ . $$
This leads  to the contravariant functor defined on the restricted category of 
embedded smooth manifolds:
$$ M  \ \Rightarrow \ \pi_{0,2,M}: T^0_2(M)\to M$$
$$ f:M \to N \ \Rightarrow \ [f, T^0_2f] \ . $$
As above we can obtain a covariant version by replacing $T^0_2f$ with $T^0_2(f^{-1})$.

\medskip

{\it (The $T^p_q$ functors) } The construction of $T := T^1_0$ , $T^*:= T^0_1$, $T^0_2$ functors 
can be generalized straightforwardly (with the same formal features)
to every $(p,q)$, getting the {\it tensorial functors}

$$ M \ \Rightarrow \pi_{p,q,M} : T^p_q(M) \to M$$
$$ f:M\to N \ \Rightarrow \ [f T^p_q f]$$
where we can stipulate to take always the covariant version (and we refer to
the restricted smooth category when necessary).

\medskip

{\it (The determinant bundle)} By using the spaces $\Lambda^0_m(T_xM)$
we get the  determinant bundle of $M$ (with $1$-dimensional fibre)
$$ \delta_M: \det T(M) \to M$$
with changes of $\det T\Aa$ coordinates
$$\det T\beta(x,r)=(\beta(x), (\det d_x \beta ) r) \ . $$

\section{Tensor fields, unitary tensor bundles}\label{tensor-field}
We can extend and generalize the content of section \ref{fields} of Chapter \ref{TD-LOCAL}
to embedded smooth manifolds.
\medskip

Let $\pi: E(M)\to M$ be any  tensor vector bundle as above, with fibre $E_xM$ over $x\in M$ of dimension say $r$. 
A {\it section} of this bundle
is a smooth map
$$\sigma: M\to E(M)$$
such that for every $x\in M$, $\pi(\sigma(x))=x$. 
In other words $\sigma$ determines a smooth {\it field of tensors} of a certain type
on $M$.  
Denote by 
$$\Gamma (E(M))$$
the set of these sections.
As for every vector bundle, every $\Gamma(E(M))$ has a canonical {\it zero section}
$$\sigma_0(x)=(x,0), \ x\in M \ . $$
In this way $M$ is canonically included into $E(M)$.
Every $\Gamma(E(M))$ is a module over the commutative ring $\Cc^\infty(M,\R)$, hence a real
vector space.

\medskip
 
$\bullet$ An element of $\Gamma(T(M))$ is called a {\it vector field} on $M$.
Generalizing verbatim section \ref{derivation}, $\Gamma(T(M))$ is isomorphic
to the vector space of {\it derivations} on $\Cc^\infty(M,\R)$, ${\rm Der}(\Cc^\infty(M,\R))$.

\medskip

$\bullet$ An element in $\Gamma(T^*(M))$ is called a $1$-{\it differential form} on
$M$.  If $f:M\to \R$ is a smooth function, then $df \in \Gamma(T^*(M))$. 

\medskip

$\bullet$ A section $g\in \Gamma(S^0_2(M))$ such that $g(x)$ is positive definite
for every $x\in M$ is called a {\it riemannian metric} on $M$. Every $M$ admits riemannian
metrics: for every riemannian metric $\hat g$ on $\R^h$ (for instance the standard $g_0$),
then the restriction of $\hat g_x$ to $T_xM$ for every $x\in M$ defines a riemannian
metric $g$ on $M$. 

$f: (M,g)\to (N,g')$ is an {\it isometry} if it is a diffeomorphism and for every $x\in M$, $v,w \in T_xM$, then
$g_x(v,w)=g'_{f(x)}(d_xf(v),d_xf(w))$.

If $(W,\phi)$ is a chart of $(M,g)$, with inverse parametrization $\psi: U \to W$,
then by imposing that $\psi$ is tautologically an isometry we get a representation
$g_U$ of $g$ in local coordinates; $g_U$ is an instance of riemannian metric 
on the open set $U\subset \R^m$ as defined in Chapter \ref{TD-LOCAL}.

\medskip

$\bullet$ Given a riemannian metric $g$ on $M$, 
for every smooth  function $f:M\to \R$ there is a unique vector field $\nabla_g f$
(called the {\it gradient of $f$ with respect to $g$}) such that for every $x\in M$, every $v\in T_xM$,
$$d_xf(v)=g_x(\nabla_gf(x),v) \ . $$

\medskip

$\bullet$ {\it (Other functors)} By setting

$$ M \ \Rightarrow \ \Gamma(T^*(M))$$
$$ f: M \to N \ \Rightarrow \  f^*: \Gamma(T^*(N)) \to \Gamma(T^*(M))$$
where
$$f^*(\omega)(x)(v)= \omega(f(x))(d_xf(v))$$ 
ones defines a contravariant functor from the category of embedded smooth manifolds to the category
of real vector spaces.

\medskip

By allowing only the diffeomorphisms as morphisms, then by setting

$$ M \Rightarrow \Gamma(T(M))$$
$$ f: M \to N \Rightarrow  f_*: \Gamma(T(M)) \to \Gamma(T(N))$$
where
$$ f_*(X)(y) :=  d_{f^{-1}(y)}(X(f^{-1}(y))   $$  
one defines a covariant functor from the `restricted' category of embedded smooth manifolds to the category
of real vector spaces.

$\bullet$ Let $(W,\phi)$ and $(U,\psi)$ be chart/parametrization of $M$ as above, then for every 
$X\in \Gamma(T(M)$, every $\omega \in \Gamma(T^*(M))$,
by using either $\phi_*$ or $\psi^*$ we get local representantions in the coordinates of $U$ of 
the type described in section \ref{fields}.  Representations in local coordinates can be straightforwardly 
developed for every field of tensors of arbitrary type on $M$.

\subsection {Unitary tensor bundles}\label {riem-bundle}
Let $(M,g)$ be endowed with the riemannian metric restriction of the standard metric $g_0$ on $\R^h$.
Set
$$UT(M)=\{ (x,v)\in T(U); ||v||_{g_x}=1 \}$$
with the restriction
$$u\pi_M : UT(M)\to M$$
of $\pi_M: T(M)\to M$.
Then $UT(M)$ is a submanifold of $T(M)$ of dimension $m(m-1)$,
and $u\pi_M$ is a surjective summersion with every fibre diffeomorphic (isometric indeed) to
the unitary sphere  $S^{m-1}$.  More precisely, the local trivializations of $T(M)$,
$$ \begin{array}[c]{ccc}
U\times \R^m&\stackrel{T\phi}{\rightarrow}&T(W)\\
\downarrow\scriptstyle{\pi_U}&&\downarrow\scriptstyle{\pi_W}\\
U&\stackrel{\phi}{\rightarrow}&W\end{array}$$

restrict to ``unitary'' local trivializations

$$ \begin{array}[c]{ccc}
U\times S^{m-1} &\stackrel{UT\phi}{\rightarrow}&UT(W)\\
\downarrow\scriptstyle{\pi_U}&&\downarrow\scriptstyle{u\pi_W}\\
U&\stackrel{\phi}{\rightarrow}&W\end{array}$$

Then $u\pi_M : UT(M)\to M$ is called the {\it unitaty tangent bundle of $M$}.

Let $\pi: E(M)\to M$ be as before any of our tensor bundles. For every $x\in M$, the positive scalar product
$g_x$ on every $T_xM$ canonically propagates to  a positive definite scalar product
$g^E_x$ on the fibre $E_xM$; this is defined as follows: given one $g_x$-{\it othonormal} basis $\Bb_x$ of $T_xM$,
$g^E_x$ is determined by imposing that the basis $\Bb^E_x$ of $E_xM$ canonically associated to $\Bb_x$
is $g^E_x$-othonormal (one verifies that this does not depend on the choice of the basis $\Bb_x$).
Then by the very same procedure we get the {\it unitary tensor bundle}
$$ u\pi: UE(M)\to M$$ 
with fibre isometric to the unitary sphere $S^{r-1}$.

\begin{remark} \label{unitary-isot}{\rm 
We have defined the unitary tangent bundle (and its relatives) by
using the restriction of the standard riemannian metric on the ambient
euclidean space.  However, if $f:M\to M$ is a diffeomorphism then in
general the unitary tangent bundle is {\it not} preserved; moreover
the costruction of {\it a} unitary tangent bundle works as well if $M$
is endowed with an arbitrary riemannian metric; {\it from a differential
topological view point, there is not a privileged riemannian
metric}. So we dispose indeed of an infinite family of unitary bundles. The total
spaces of two unitary bundle defined with respect to two metrics $g_0$
and $g_1$ are canonically diffeomophic via radial diffeomorphisms fibre
by fibre, centred at the origine of each $T_xM$.  Moreover by using
the path of riemannian metrics $g_t = (1-t)g_0 + tg_1$ this
diffeomorphism is connected to the identity by a smooth path (an {\it
  isotopy}) through diffeomorphisms of unitary bundles of the same
type.  This considerations ``propagate'' to all tensor bundles. Every unitary tensor bundle is
well defined up to isotopy.}
\end{remark}

\section{Parallelizable, combable and orientable manifolds}\label{orientable}
An embedded smooth manifold $M\subset \R^h$ of dimension $m\geq 1$ is
said {\it parallelizable} if there are $m$ sections
$\Sigma=(\sigma_1,\dots , \sigma_m)\in \Gamma(T(M))^m$ such that for
every $x\in M$, $\Sigma(x)$ is a basis of $T_xM$. This property
``propagates'' to every of our favourite tensor bundles say $\pi:
E(M)\to M$ with fibres $E_xM$ of dimension say $r$.  In fact for every
$(p,q)$, the canonical correspondence $\Sigma(x) \to \Sigma(x)^p_q$
determines
$$\Sigma^p_q \in \Gamma(T^p_q(M))^{m^{pq}}$$ such that for every $x\in
M$, $\Sigma(x)^p_q$ is a basis of $T^p_q(T_xM)$; similarly we have a
{\it nowhere vanishing} section $\det \Sigma$ of the determinant
bundle $\delta_M: \det(T(M))\to M$. In generical notations, denote
$\Sigma \in \Gamma(E(M))^r$ such a distinguished field of bases. We
can define
$$ t_\Sigma : M\times \R^r \to E(M), \ t_\Sigma(x,v)= (x, \sum_j
v_j\sigma_j(x)) $$ clearly this is a diffeomorphism and also a vector
bundle map in the sense that for every $x\in M$, it induces a {\it
  linear isomorphism} $\{x\}\times \R^r \to E_xM$.  Moreover the
following diagram obviously commutes
$$ \begin{array}[c]{ccc}
M\times \R^r&\stackrel{t_\Sigma}{\rightarrow}& E(M)\\
\downarrow\scriptstyle{p_M}&&\downarrow\scriptstyle{\pi}\\
M&\stackrel{{\rm id}_M}{\rightarrow}&M\end{array}$$
 Then $t_\Sigma$ is called a {\it global trivialization of the bundle $E(M)$}.
 
\medskip

So $M$ is parallelizable if and only if its tangent bundle is strictly
equivalent to a product bundle, and a {\it necessary} condition in
order that $M$ is parallelizable is that the determinat bundle of $M$
has a nowhere vanishing section. Let us say that $M$ is {\it orientable} if it verifies such 
a necessary condition. Obviously, if $M$
is parallelizable, then it is ``{\it combable}", that is it carries a nowhere
vanishing tangent vector field. Every open set of $\R^n$ is parallelizable, hence
orientable and combable. The same facts hold {\it locally} on every manifold
$M$.  So we have here a bunch of crucial genuine global questions concerning 
the structure of a generic smooth manifold $M$ in terms of the existence of 
suitable patterns of sections of natural fibre bundles over $M$. 

\medskip

Let us explicate now the definition of orientability. It is clear
that $M$ is orientable if and only
if every connected component of $M$ is orientable; so let us assume that $M$
is connected. Consider the {\it unitary} determinant bundle. 
The fibre is $S^0=\{\pm 1\}$, so we can write it as 
$$\pG: \tilde M \to M$$
where $\tilde M$ is a $m$-manifold, $\pG$ is a covering map of degree $2$
called the {\it orientation covering of $M$}. 
The fibre over every $x\in M$ is  $\{ (x, \pm 1) \}$. 
There are two possibilities: either $\tilde M$ is connected
or it has two connected components 
$$\tilde M = \tilde M_+ \cup \tilde M_- $$
where 
$$\tilde M_\pm= \{ (x, \pm 1); x\in M \} \ . $$
Obviously the restriction of $\pG$ to $\tilde M_\pm$ is a diffeomorphism
(basically it is the identity). If $x \to (x,\sigma(x)) $ is a nowhere vanishing section
of the determinant bundle,
as $M$ is connected the  sign $\frac{\sigma(x)}{|| \sigma(x)||_{g(x)}}$ is constant.
So we have proved 
\begin{proposition}\label{orient-1}
  $M$ is orientable if and only if $\tilde M = \tilde M_+ \cup \tilde M_-$
is not connected.
\end{proposition}

\begin{example}{\rm Referring to section \ref{stif-gram-fib}, examples of {\it connected}
$\pG: \tilde M \to M$ are the natural covering maps 
$$ S^n \to \PP^n(\R)$$
when $n$ is {\it even}. Then such projective spaces are not orientable.}
\end{example}

The alternative ``$M$ orientable/non-orientable'' can be reformulated
as follows: a {\it signature} $\sG$ on an atlas $\Uu$ of $M$
assign to every chart a sign
$\sG(W,\phi)=\pm 1$. Given such an $\sG$, modify $\Uu$ to $\Uu_\sG$ by post composing
every chart with negative sign with a linear reflection of $\R^m$ (which has the determinant equal to $-1$).
An atlas $\Uu$ is {\it oriented} if all changes of coordinates for $\Uu$ have the determinant sign
constantly equal to $1$. Then we have

\begin{proposition} \label{orient-2}
  The following facts are equivalent to each other:
\begin{enumerate}
 \item $M$ is orientable;
 \item There exists an oriented atlas $\Uu$ of $M$;
 \item For every atlas $\Uu$ of $M$ there exists a signature $\sG$ such that
   $\Uu_\sG$ is oriented.
\end{enumerate}

\end{proposition}  

 We leave the proof to the reader as an useful exercise on this
 complex of definitions. The condition of point (2) in the Proposition
 is often given as the {\it very definition of orientability}. A
 reader can do it without effecting the rest of our discussions. Here
 is some further remarks on these notions.

 $\bullet$ If $M$ is connected and orientable, then every oriented atlas $\Uu$
 is contained in an unique maximal oriented atlas. There are exactly two
 maximal oriented atlas say $\Aa^\pm$. Any signature $\sG$ on $\Aa$ such that
 $\Aa_\sG$ is oriented produces one among $\Aa^\pm$; $\sG$ produces $\Aa^+$
 if and only if the opposite signature $-\sG$ produces $\Aa^-$.
 By definition $\Aa^\pm$ define two opposite {\it orientations} of $M$ and make
 it (in two ways) an {\it oriented} manifold. If $M$ is oriented, $-M$
 denotes $M$ endowed with the opposite orientation. The two components
 of $\tilde M$ are naturally oriented and correspond to the two orientations
 of $M$.

\medskip

$\bullet$ The definition via oriented atlas allows us to recover
the elementary notion of orientation of $\R^m$ as a vector space.
By definition two bases $\Bb$ and $\Dd$ of $\R^m$ are {\it
  co-oriented} if the determinant of the change of linear coordinates passing
fro $\Bb$ to $\Dd$ is positive.  By the multiplicative properties of
the determinant, this defines an {\it equivalence relation} on GL$(m,\R)$
(considered as the space of bases of $\R^m$); then an {\it orientation} on $\R^m$
is an equivalence class of bases. Let us call {\it standard orientation} 
the class $[\Cc]$ of the canonical basis $\Cc$. If $U$ is a (connected) open set
of $\R^m$ we get the {\it standard field of orientations} by giving
each $T_xU=\R^m$ the standard orientation. $U$ is obviously an
orientable manifold and we can take the maximal oriented atlas say $\Aa^+$ of $U$ which
contains the chart ${\rm id}:U \to U$. Let $\psi: U' \to U" \subset U$
the local parametrization associated to a chart of $\Aa^+$. By taking
the standard field of orientations on $U'$, $d\psi$ transforms it to
the field of orientations $\{ [d_y\phi (\Cc)]\}_{x=\psi(y)}$ on
$U"$. The fact that $\psi$ belongs to $\Aa^+$ just means that this
last field coincides with the standard one on $U''$. Extenting this
considerations to an arbitrary manifolds $M$, an orientation on $M$,
if any, can be considered as a ``locally coherent'' field of orientations
on each $T_xM$.

\medskip

$\bullet$ Let $f: M \to N$ be a diffeomorphism. If $\Uu=
\{(W,\phi)\}$ is an atlas of $M$, then
$$f(\Uu):= \{(f(W), \phi \circ f^{-1}) \}$$ is an atlas of $N$.
The proof of the following Lemma follows immediately from the definitions.

\begin{lemma}\label{oriented-diff}
  Let $f:M\to N$ be a diffeomorphism between connected oriented manifolds with
  maximal oriented atlas say $\Aa^+_M$ and $\Aa^+_N$ respectively.
  The following facts are equivalent to each other.
\begin{enumerate}
\item $f(\Aa^+_M) = \Aa^+_N$ .
\item There exist an oriented atlas $\Uu \subset \Aa^+_M$ such that $f(\Uu) \subset \Aa^+_N$.
\item For every representation in local coordinates $f_{U,U'}: U \to U'$ of $f$ relative to charts
in $\Aa^+_M$ and $\Aa^+_N$ and for every $x\in U$, then $\det d_xf_{U,U'} >0$.
\end{enumerate}

If one (hence all) of the above conditions is verified, then we say
that $f$ is an {\rm oriented diffeomorphism}.
\end{lemma}

$\bullet$  By specializing the {\it objects} to oriented manifolds we get a sub-category 
of our favourite one.

\begin{remark}\label{0-manifolds} (Oriented $0$-Manifolds) {\rm A $0$-manifold is a discrete
set of points, hence  just one point if connected. We stipulate that
it is orientable and is {\it oriented} by giving it a sign $\pm 1$.}
\end{remark}

\section{Manifolds with boundary, oriented boundary, proper submanifolds}\label{boundary}
By definition an embedded smooth $m$-manifold $M\subset \R^n$
is locally diffeomorphic to open sets of the basic model $\R^m$.
Let us change this last by taking instead the {\it half-space}
$$ \HH^m = \{x\in \R^m; \ x_m\geq 0 \}$$
with the {\it boundary}
$$\partial \HH^m= \{x\in \HH^m; \ x_m=0\} \ . $$

\begin{definition}\label{b-embedded}{\rm 
  For every $0\leq m \leq n$, a topological subspace $M\subset \R^n$
  is an {\it embedded smooth $m$-manifold with boundary} if for every $p\in M$,
  there exist an open
  neighbourhood $W$ of $p$ in $M$, an open set $U$ of $\HH^m$ and a
  diffeomorphism $\phi: W \to U$. The notions of ``chart'', ``local parametrization'',
  ``atlas'' extend straightforwardly. By definition, the {\it boundary} $\partial M$ is the set of points $p\in M$
  such that there exists a chart $(W,\phi)$ at $p$ such that $\phi(p)\in \partial \HH^m$.
  }
\end{definition}

The following Lemma provides a basic way to produce manifolds with boundary.

\begin{lemma}\label{>=} Let $X$ be a $m$-manifold with empty boundary, $f:X\to J$
  a surjective summersion, where $J$ is an open interval of $\R$, and $0\in J$.
  Then $M=\{x\in X; \ f(x)\geq 0\}$ is a $m$-manifold with boundary $\partial M = \{f(x)=0\}$.
\end{lemma}
\Dim The question being of local nature one can reduce to summersions in normal form
for which the result is evident.

\cvd

The following Lemma contains by the way an extension of Lemma \ref{relative}
and similarly is an application of the inverse map theorem (and its corollaries).
  
\begin{lemma}\label{b-relative} Let $M\subset \R^n$ be an $m$-manifold with
  boundary. Then

  (1) If $p\in \partial M$, then for every chart $(W,\phi)$ of $M$ at $p$,
  $\phi(p)\in \partial \HH^m$.

  (2) ${\rm Int}(M):= M \setminus \partial M$ is an open set in $M$ and
  a manifold with empty boundary (called the {\rm interior} of $M$). For every $p\in{\rm Int}(M)$ there are
  normal relative charts of $(\R^n,{\rm Int}(M))$ at $p$ that do not intersect
  $\partial M$.

  (3) For every $p\in \partial M$, there are normal relative charts
  of $(\R^n,M,\partial M)$ at $p$:
  $$\beta: (\Omega,\Omega \cap  M, \Omega \cap \partial M,p)\to
  (B^n(0,1), B^n(0,1)\cap \HH^m, B^n(0,1)\cap \partial \HH^m ,0)$$

  (4) If $\partial M \neq \emptyset$, then it is $(m-1)$-manifold with empty boundary.
  
 \end{lemma}

\cvd

The definition of ``embedded smooth manifold with boundary'' does not
exclude that $\partial M = \emptyset$. We have early considered such a {\it boundaryless}
case. It is formally convenient to stipulate that the empty set $\emptyset$ is a
$k$-boundaryless manifold for every $k\in \N$. In such a way for example point (4)
of the last Lemma holds even if $\partial M = \emptyset$. By setting $M=(M,\emptyset)$
for every boundaryless manifold, the early category of {\it embedded smooth manifolds}
extends to the category of  {\it embedded smooth manifolds with boundary}.
Let us briefly retrace within such an extension the main facts developed so far .

\medskip

$\bullet$ The tangent functor and its relatives extend verbatim.
If $\partial M$ is non empty, the inclusion $j:\partial M \to M$ leads to
a vector bundle embedding $[j,Tj]$ of $\pi_{\partial M}:T(\partial M)\to \partial M$
into $\pi_M: T(M)\to M$. The total space  $T(M)$ is a manifold with boundary
equal to the restriction over $\partial M$  of the tangent bundle of $M$
(with the notions that we will introduce in Chapert  \ref{TD-EMB-VB} it is the pull-back
$j^*T(M)$ over $\partial M$).
Similarly for the other tensors bundles.

\medskip

$\bullet$ Also ``orientability/orientation'' estends directly. The boundary $\partial M$
of an {\it oriented} $M$ is {\it orientable} and we can fix the following procedure
in order to make it the {\it oriented boundary} of $M$:

\medskip

{\it (``First the outgoing normal'')}
Take an oriented atlas $\Uu$ of $M$ made by normal charts. Post compone every chart along
the boundary $\partial M$ with a trasformation  $r\in SO(m)$ such that
$r(e_1,\dots,e_m)=(-e_m,r(e_1,\dots,e_{m-1})$. The so obtained atlas, say $r\Uu$ is again
an oriented atlas of $M$ and its restriction to $\partial M$ is an oriented atlas which
carries a determined orientation of the boundary.

\medskip

By the usual convention $M=(M,\emptyset)$, the category of {\it oriented}  boundaryless manifolds extends
to the category of {\it oriented} manifolds with oriented boundary.

\medskip

$\bullet$ {\it (Submanifolds)} Alike the boundaryless case, let us stipulate that if $Y, M \subset \R^n$
are embedded smooth manifolds with boundary and $Y\subset M$, then $Y$ is a {\it submanifold}
of $M$. By extending the Remark \ref{strange-sub}, because of  the presence of the boundary
there are several qualitatively different ways of being a submanifold; let us list a few
examples:

\begin{enumerate}
\item $(Y\subset M) = (\overline B^n(0,1) \subset B^n(0,2))$: $\partial Y\neq \emptyset$ and
  $Y$ is contained in the interior of $M$.

\item $(Y\subset M) =  ({\rm Int}(M)\subset M)$; if $\partial M \neq \emptyset$, then $Y$
is not closed in $M$.

\item  $(Y\subset M) = (N \subset B^n(0,1))$, where $N$ is defined in Remark \ref{strange-sub}:
  $Y$ is boundaryless, is contained in the interior of $M$, and every point of $\partial M$
  is in the closure of $Y$; again $Y$ is not closed in $M$.

\item  $(Y\subset M)$ where $Y= \overline B^n(0,1)$, $M=\{ x_n\geq -1\}$. Then $\partial Y$ is tangent
  to $\partial M$, while the interior of $Y$ is contained in the interior of $M$.

\item Let $\gamma:=\gamma_{1,2}:\R \to \R$ the bump function defined in Chapter \ref{TD-LOCAL}.
  $(Y\subset M) = (N \subset \HH^2)$, where $N=\{(x,y)\in \HH^2; \ y\geq \gamma(x)\}$. 
Then $\partial Y$ is partially contained in the interior of $M$, partially into $\partial M$.
\item $\dots$

\end{enumerate}

Among this wide typology there is a particularly clean type which deserves to be pointed out
by a definition.
\begin{definition}\label{proper-sub}
{\rm Let $Y\subset M \subset \R^n$ smooth manifolds with
boundary. Then $Y$ is a {\it proper submanifold} of $M$ if
\begin{enumerate}
\item $Y$ is closed in $M$;
\item $\partial Y = Y\cap \partial M$;
\item $Y$ is {\it transverse} to $\partial M$.  This means that for every $p\in Y\cap \partial M$
$$ T_pM = T_pY+T_p\partial M \ . $$
\end{enumerate}
}
\end{definition}

All the above examples are not proper.  Every $M$ is a proper submanifold of itself. 
The properness implies for instance that 
every boundaryless component of $Y$ is contained in the interior
of $M$; if $\partial M = \emptyset$, then also $\partial Y=\emptyset$;
if $\dim Y = \dim M$ then $Y$ is union of connected components of $M$.

The following Proposition extends (1) of Proposition \ref{global-imm-summ}
in two ways, to manifolds with boundary and to oriented manifolds.

\begin{proposition}\label{b-summ}
Let M be a manifold with boundary and $N$ a boundaryless
one. Let $f:M \to N$ be a surjective {\rm relative summersion}
(that is both $f$ and $\partial f:= f_{|\partial M}$ are
summersions). Then:

(1) For every $q\in N$,  $Y=f^{-1}(q)$ is a proper 
submanifold of $M$, $\dim Y = \dim M - \dim N$.

(2) If both $M$ and $N$ are oriented, then $Y$
is orientable, and we can fix a procedure
to orient it, in such a way that the 
orientation  of $\partial Y$ as oriented boundary of $Y$
coincides with the orientation obtained by applying the
procedure to $\partial f$, provided that $\partial M$
is the oriented boundary of $M$.
\end{proposition}
\Dim  
Assume that $M\subset \R^h$, $\dim M= m$, $\dim N = n$.
If $q$ does not belong to the image of $\partial f$, then we apply directly
Proposition \ref{global-imm-summ} so that $Y$ is a closed boundaryless
submanifold of the interior of $M$. Assume now that $q$ belongs to
the image of $\partial f$. The question being of local nature,
we reduce to analyze a representation (called $f$ as well)  of $f$ in local coordinates
which are normal for $(M,\partial M)$:
$$ f:  (B^m(0,1)\cap \HH^m, B^m(0,1)\cap \partial \HH^m)\to U\subset \R^n$$
$q=0\in U$. Moreover we can assume that $f$ is the restriction
of a smooth map  $g: B^m(0,1)\to U$ defined on the whole of $B^m(0,1)$, 
which a surjective summersion.
By applying again  Proposition \ref{global-imm-summ} to $g$,
we have that $\tilde Y= g^{-1}(0)$ is a boundaryless submanifold
of $B^m(0,1)$ of the correct dimension, such that $Y=f^{-1}(0)$
is $Y=\tilde Y \cap \HH^m $. As $f$ is a relative summersion, 
one readly checks that $\tilde Y$ is transverse to $\partial \HH^m$
and that the restriction say $\pi$ to $\tilde Y$ of the projection
onto the $x_m$ coordinate is a summersion onto its image and that
$Y=\{y\in \tilde Y; \ \pi(y)\geq 0\}$. We conclude by applying Lemma \ref{>=}.

Let us come to the orientation. First consider the case $f={\rm id}_M$.
Then $Y=\{p\}$ is just a point of $M$. Let us orient it by giving it the sign
$+1$.  By applying the rule to $\partial f$ we get the same sign.
In the general case. For every $p\in Y$
let 
$$\nu(p)= (T_pY)^\perp \cap T_pM$$
clearly
$$ T_p(M)= T_pY\oplus \nu(p)$$
and $\nu(p)$ varies ``smoothly'' when $p$ varies along $Y$
(by using the contents of next Chapter \ref {TD-EMB-VB}
this means precisely that $\nu: Y \to \GG_{k,n}$ is a smooth
map). In our hypotheses, for every $p\in Y$, the restriction
of $d_pf$ to $\nu(p)$ is a linear isomorphism onto $T_{f(p)}N$.
Let us consider the orientation on $N$ as a field of orientations
on the $T_yN$, $y\in N$, (i.e. a field of equivalence classes of
bases of $T_yN$) which is locally coherent). Take an orienting (say ``positive'')
basis $\Bb_q$ of $T_qN$. For every $p\in Y$, lift it to a basis
$\Bb_p$ of $\nu(p)$  by means of the restriction of the differential of
$f$. This determines a field of ``transverse orientations'' $[\Bb_p]$
along $Y$. At every $p$, take a basis $\Dd_p$ of $T_pY$
such that the basis $\Dd_p \oplus \Bb_p$ of $T_pM$ (compatible
with the above direct sum decomposition of $T_pM$)
is positive with respect to the given orientation of $M$.
This determines a field $[\Dd_p]$ of orientations on the $T_pY$,
eventually the desidered orientation of $Y$. 
This procedure could be finalized in terms of the construction
of a suitable oriented atlas for $Y$; we leave it to the reader. 
One can check that
the restriction of this procedure to $\partial f$ is compatible
in the sense of the last statement of the proposition.

\cvd  

\medskip

$\bullet$ Also the {\it topologies} of spaces of smooth maps between
manifolds with boundary extend word by word.

\section{Product, manifolds with corners, smoothing}\label{corner}
We know that the product of two boundaryless manifolds is a boundaryless
manifold. The situation is more complicated if we consider non empty boundaries.
The following Lemma is immediate.

\begin{lemma}\label{easy-prod} Let $M$ be a boundaryless (embedded smooth) $m$-manifold,
  $N$ be a $n$-manifold with $\partial N \neq \emptyset$. Then $M\times N$ is a
  $(m+n)$-manifold with $\partial(M\times N)= M\times \partial N$
\end{lemma}

\cvd

However, if both $\partial M$ and $\partial N$ are non empty, then $M\times N$ {\it is no longer}
an embedded  smooth manifold with boundary.

\begin{example}\label{square} {\rm
As a basic example, consider the square 
$$Q=D_1\times D_2:=[-1,1]\times [-1,1] \subset \R^2 \ . $$
Its topological frontier is
$$\partial Q = (\partial D_1 \times D_2) \cup (D_1\times \partial D_2) \ ;$$
its interior 
$$ Q\setminus \partial Q = {\rm Int}(D_1)\times {\rm Int}(D_2) $$
is an open set of $\R^2$ hence a $2$-manifold with empty boundary;
$$ Q\setminus (\partial D_1 \times \partial D_2)$$
is a $2$-manifold with boundary equal to
$$\partial Q \setminus (\partial D_1 \times \partial D_2)\ ; $$
$ \partial D_1 \times \partial D_2$
is a $0$-manifold. The points where $Q$ fails to be a manifold with boundary are
the ``corner'' points which form $\partial D_1 \times \partial D_2$.
  }
\end{example}

The behaviour of such a simplest example is qualitatively the general one:

\begin{proposition}\label{prod-corner} Let $(M,\partial M)\subset \R^h$ and $(N,\partial N)\subset \R^k$
  be an $m$-manifold and an $n$-manifold with boundary respectively. Then $M\times N \subset \R^h\times \R^k$
  verifies the following properties:
  \begin{itemize}
  \item  Set
    $$\partial (M\times N):= (\partial M \times N)\cup (M\times \partial N)\ . $$
    Then
    $$ (M\times N) \setminus \partial (M\times N)$$
    is a boundaryless $(m+n)$-manifold;
  \item $$(M\times N) \setminus (\partial M \times \partial N)$$
    is a $(m+n)$-manifold with boundary equal to
    $$\partial (M\times N) \setminus (\partial M \times \partial N) \ ;$$
  \item $\partial M \times \partial N$ is a boundaryless $(m+n-2)$-manifold.
  \end{itemize}
\end{proposition}

\cvd

Hence $M\times N$ fails to be a manifold with boundary at the ``corner locus''
$\partial M \times \partial N$. This means that the category of embedded smooth manifolds
with boundary is {\it not closed} with respect to the product. This is somehow unpleasant.
A way to fix this fact is to
enlarge our category by  extending the sets of basic models, incorporating
the corners. We do it in the minimal way suited to incorporate
such product manifolds.

\begin{definition}\label{basic-corner}{\rm The {\it basic $m$-corner models} is

    $$\CC^m= \{ x\in \R^m; \ x_m\geq 0, \ x_{m-1}\geq 0 \}$$
    that is the intersection between $\HH^m$ with another halfspace.
    Its {\it boundary} (in fact its topological frontier) is
    $$ \partial \CC^m = \{x\in \CC^m; \ x_m=0 \} \cup   \{x\in \CC^m; \ x_{m-1}=0 \} \ . $$ 
    $\CC^m \setminus \{x_m=0, x_{m-1}=0\}$ is a manifold with boundary and the last set
    is its {\it corner locus}.
    
  }
\end{definition}

\begin{definition}\label{man-with-corner}{\rm 
For every $0\leq m  \leq n$, a topological subspace $M\subset \R^n$
  is an {\it embedded smooth $m$-manifold with corners} if for every $p\in M$,
  there exist an open
  neighbourhood $W$ of $p$ in $M$, an open set $U$ of $\CC^m$ 
  and a
  diffeomorphism $\phi: W \to U$. The notions of ``chart'', ``local parametrization'',
  ``atlas'' extend straightforwardly. The {\it boundary} $\partial M$ is the set of points $p\in M$
  such that there exists a chart $(W,\phi)$ at $p$ such that $\phi(p)\in \partial \CC^m$.
  The {\it corner locus} is where $M$ is not locally a smooth manifold with boundary.    }
\end{definition}

\medskip

The following properties clearly hold for the basic models and descend easily
to every manifold with corners.

\medskip

(i) Every manifold with corners is naturally {\it stratified} by means
  of the disjoint locally finite  union of boundaryless connected smooth manifolds
  (of varying dimension $m-2\leq d \leq m$)
  called the {\it strata}; the top dimensional strata are the components of the
  boundaryless smooth $m$-manifold   $M\setminus \partial M$; the $(m-1)$-strata
  are the componets of $\partial M$ from which we have removed the corner locus;
  the $(m-2)$-strata are the components of the corner locus which  is a boundaryless
  manifold of dimension $m-2$ contained in the boundary of $M$.
  The closure of every stratum is union of strata, as well as the maximal smooth manifold with boundary
  contained in the closure of every stratum.

(ii) The product of two smooth manifolds with boundary is a manifold with corners.
\medskip

However, manifolds with ``codimension $2$'' corners are not closed under the product (take for instance
the cube $[-1,1]^3$). So we have only shifted the difficulty and we should extend furthermore our category of
manifolds. This would bring us a bit far away from our original objects of interests.
  Fortunately there is another way that leads back manifolds with corners (according with the above
  restrictive definition) to ordinary manifolds with
  boundary, even though {\it up to diffeomorphism}. To introduce such a  {\it ``smoothing the corner''}
  procedure, let us consider again our simplest square example. The function
  $$ f: \R^2\to, \ f(x)=(x_1-1)(x_2-1)(x_1+1)(x_2+1)$$
  has the property that $Q$ is the closure of a connected component of
  $$ \R^2 \setminus f^{-1}(0) $$
  and for every $x\in {\rm int}(Q)$, $f(x)>0$. For every $\epsilon >0$, sufficiently small,
  there is a connected component $Q_\epsilon$ of $f(x)\geq \epsilon$ contained in the interior
  of $Q$, and which is a smooth manifolds with boundary {\it homeomorphic} to $Q$.
  Moreover, we can construct a ``piece-wise smooth'' radial 
  homeomorphism (centred at $0$) $s: Q_\epsilon \to Q$ such that the natutal stratification of $Q$ lifts to
  a stratification by smooth submanifolds of $Q_\epsilon$ and the restriction to
  the maximal manifold with boundary contained in the closure of every stratum is a diffeomorphism
  onto its analogous image in $Q$. Finally, up to diffeomorphism, the result of such a smoothing
  does not depend on the specific implementation (in particular on the choice of the small $\epsilon$).

  This basic idea can be generalized. By applying it to $\CC^m$, by using
  $M_\epsilon = \{ x_mx_{m-1} \geq \epsilon\}\cap \CC^m$, $\epsilon>0$ small enough,
  we get nice local smoothing homeomorphism $s: M_\epsilon \to \CC^m$ with the same qualitative
  properties as above. 
  Then one should have to prove that such local smoothings
  can be patched to give a global smooth atlas. This could be a bit technically demanding (with simplifications
  if the manifolds are compact) and we do not further push in that direction.  In Section \ref{corner-smooth} we will reconsider and
  properly establish such a smoothing procedure in a more flexible ``abstract'' setting.  
  Anyway we already state the following
  
  \begin{proposition}\label{smoothing} For every $m$-manifold with corner $M\subset \R^h$, then
    
    (1) by implementing a determined {\rm ``smoothing the corner''}
    procedure, we get a smooth manifold with boundary $\tilde M
    \subset \R^h$ and a piece-wise smooth homeomorphism
    $$\sG: (\tilde M, \partial \tilde M)\to (M,\partial M)$$
    such that the natural stratification of $M$ lifts to a stratification of $\tilde
    M$ by boundaryless smooth submanifolds, and the restriction of
    $\sG$ to the maximal smooth manifold with boundary contained in
    the closure of every stratum of $\tilde M$ is a diffeomorphism
    onto its analogous image in $M$.

    (2) $\tilde M$ is uniquely determined up to diffeomorphism (i.e. it does not
    depend on the actual implementation of the procedure).

  \end{proposition}

  \cvd

  \medskip

  Coming back to our motivating problem, the product of two smooth manifolds with boundary
  as a smooth manifold with boundary is well defined {\it up to diffeomorphism}.

\chapter{Stiefel and Grassmann manifolds}\label{TD-STIEF-GRASS}
The tensorial vector bundles contructed in Chapter \ref{TD-DIFF}
belong to a wide category of ``embedded vector bundles'' that we will consider 
in Chapter  \ref{TD-EMB-VB};  the core of that discussion 
will  consist in remarkable families of embedded smooth
manifolds and smooth maps between them that we are going to study by
themselves.

\section{Stiefel manifolds}\label{stiefel}
 We introduce first the  {\it Stiefel manifolds}. There are two
 versions that we call {\it linear} and {\it orthogonal}
 respectively. For every $n\in \N$ and every $0\leq k \leq n$, the
 {\it linear Stiefel manifold} $L_{n,k}$, as a {\it set}, is the set
 of ordered $k$-uple $(v_1,\dots,v_k)$ of linearly independent vectors
 in $\R^n$.  By arranging each of them in a $n\times k$ matrix $A$ (so
 that $v_j$ is the $j$-column of $A$), $L_{n,k} \subset M(n,k,\R)$. In
 fact it is {\it an open subset}: consider the smooth function
 $\delta: M(n,k,\R)\to \R$ defined in the proof of Proposition
 \ref{stab}, then $L_{n,k}= M(n,k,\R)\setminus \delta^{-1}(0)$.  This
 specifies the embedded smooth manifold nature of $L_{n,k}$. As a
 particular case we have ${\rm GL}(n,\R)= L_{n,n}$.
 For every $P\in {\rm GL}(n,\R)$, $A\to PA$ defines a  diffeomorphism
 (restriction of a linear map)
 $L_{n,k}\to L_{n,k}$, and it is well known that this action is transitive;
 in particular for every $A\in L_{n,k}$, there exists $P\in {\rm GL}(n,\R)$
 such that $PI_{n,k}=A$
 where $I_{n,k}$ is the matrix whose columns are $e_1,\dots,e_k$, the
 first $k$ vectors of the canonical basis of $\R^n$. 
 
 \medskip
 
 Now, let $S_{n,k} \subset L_{n,k}$ be the closed subset defined as
 $f^{-1}(I_k)$ where
 $$f: L_{n,k}\to S(k,\R)$$ is the smooth map $f(A)=A^tA$ with values
 in the space $S(k,\R)$ of $k\times k$ {\it symmetric} matrices which
 can be identified with $\R^{\frac{k(k+1)}{2}}$. In other words, we
 require that the columns of any $A\in S_{n,k}$ form an orthonormal
 system. As particular cases we have $S_{n,1}=S^{n-1}$, $S_{n,n}=O(n)$
 the classical (real) {\it orthogonal groups}. As
 $M(n,k,\R)=(\R^n)^k$, we see immediately that $S_{n,k}\subset
 (S^{n-1})^k$, hence $S_{n,k}$ is compact. 
 The above action of GL$(n,\R)$ on $L_{n,k}$ restricts to a transitive
 action of $O(n)$ on $S_{n,k}$: for every $A\in S_{n,k}$,
  there exists $P\in O(n)$ such that $PA=I_{n,k}$.
 It follows that in order to prove that $S_{n,k}$ is an embedded smooth manifold 
 in $(\R^n)^k$,
 it is enough to prove that there is a chart $(W,\phi)$ of $S_{n,k}$
 at $J:= I_{n,k}$. Hence it is enough to prove that $d_Jf$ is
 surjective and conclude by applying again Theorem
 \ref{summersion}. Let us compute $d_Jf$ by the very definition of the
 differential. Then
 $$ df_J(B)= \lim_{t\to 0} \frac{ (J+tB)^t(J+tB)-I_k}{t} =$$
 $$ \lim_{t\to 0} ( J^tB + B^tJ + tB^tB)= J^tB + B^tJ \ . $$
 We have to prove that for every symmetric matrix $C\in S(k,\R)$
 there exists $B\in M(n,k,\R)$ such that
 $J^tB + B^tJ=C$.   
 Set $B= \frac{1}{2} JC$. Then 
 $$J^tB + B^tJ=
 \frac{1}{2}J^tJC + \frac{1}{2}C^tJ^tJ= \frac{1}{2}C+\frac{1}{2}C^t = C $$
 because $C=C^t$. Summarizing, $S_{n,k}$ is a compact embedded
 smooth manifold in $L_{n,k}\subset M(n,k,\R)=(\R^n)^k$,
 of dimension
 $$ \dim S_{n,k} = nk - \frac{k(k+1)}{2} \ . $$
 $S_{n,k}$ is called a {\it orthogonal Stiefel manifold}.
 In particular the orthogonal group $O(n)$  is a compact 
 embedded smooth submanifold of $(S^{n-1})^n$ of dimension 
 $$\dim O(n) = n^2 -  \frac{n(n+1)}{2} \ . $$
 
\begin{remark} {\rm The operation $(A,B)\to AB$, 
and $A\to A^{-1}$ that define the group structure of ${\rm GL}(n,\R)$
are smooth (for $A^{-1}$ recall the determinantal formula based on
{\it Cramer's rule}).  These restrict to smooth operations giving the
group structure of the manifold $O(n)$. Hence ${\rm GL}(n,\R)$ and
$O(n)$ are basic examples of {\it Lie group}. 
$O(n)$ is a Lie subgroup of ${\rm GL}(n,\R)$, in the sense
that the first is a submanifold of the second and the smooth
operations are compatible.}
\end{remark}

\medskip
 
 The {\it Gram-Schmidt orthonormalization algorithm} applied to the
 ordered columns of every $A\in L_{n,k}$ defines a smooth map
 $$ \rG_{n,k}: L_{n,k} \to S_{n,k}$$ 
 which is onto and such that
 $\rG_{n,k}(A)=A$ for every $A\in S_{n,k}$. The map $\rG_{n,k}$ is the
 {\it canonical retraction} of $L_{n,k}$ onto $S_{n,k}$.
 
 \section {Fibrations of Stiefel manifolds by Stiefel manifolds}
 \label{Stiefel-fib}
 For every $0\leq h < k \leq n$, $L_{n,k}$ is a submanifold (an open set)
 in the product $L_{n,h}\times L_{n,k-h}$ and denote by
 $$l_{k,h}: L_{n,k} \to L_{n,h}$$
 the restriction of the natural projection onto the first factor.
 This map is {\it equivariant} for the above actions of {GL}$(n,\R)$
 on both Stiefel manifolds (i.e. $l_{k,h}(PA)= Pl_{k,h}(A)$), hence
 in order to study local properties such as the smoothness of the map,
 it is enough to study the restriction of $l_{k,h}$ on $l_{k,h}^{-1}(\Omega)$
 where $\Omega$ is a neighbourhood of $I_{n,h}$. 
 Clealy $l_{k,h}(I_{n,k})=I_{n,h}$. The fibre $F_{k,h}:=l_{k,h}^{-1}(I_{n,h})$
 over $I_{n,h}$ is made by the  $2\times 2$ block matrices of the form
$$Y(S,D):=
\begin{pmatrix}
I_h&S\\
0&D
\end{pmatrix}
$$ where $(S,D)\in M(h,k-h,\R) \times L_{n-h,k-h}$.  If $P\in {\rm
  GL}(n,\R)$ is such that $PI_{n,h}=A$, then
$P(l_{k,h}^{-1}(I_{n,h}))= l_{k,h}^{-1}(A)$, all fibres are
diffeomorphic to each other.  Let $\Omega$ be the open neighbourhood
of $I_{n,h}$ made by matrices of the form
$$ X=
\begin{pmatrix}
B\\
R
\end{pmatrix}
$$
where $B\in {\rm GL}(h,\R)$. 
We define the smooth  map $X\to P(X)\in {\rm GL}(n,\R)$
$$ P(X)=
\begin{pmatrix}
B&0\\
R&I_{n-h}
\end{pmatrix}
$$
such that $P(X)I_{n,h}=X$.
Finally we have the following commutative diagram of smooth maps 
$$ \begin{array}[c]{ccc}
\Omega \times F_{k,h} &\stackrel{ \Psi }{\rightarrow}& l_{k,h}^{-1}(\Omega)\\
\downarrow\scriptstyle{\pi_\Omega}&&\downarrow\scriptstyle { l_{k,h}}\\
\Omega&\stackrel{{\rm id_\Omega}}{\rightarrow}&\Omega \end{array}$$ 
such that the first row is the diffeomorphism defined by
$$ (X,S,D)\to P(X)Y(S,D) \ . $$
The costant section of the product on the left, $X\to (X,0,I_{n-h,k-h})$
is transformed into the section of $l_{k,h}$ over $\Omega$:
$$ s(X)=
\begin{pmatrix}
B&0\\
R&I_{k-h}
\end{pmatrix}
$$

A similar construction can be performend for the orthogonal
Stiefel manifolds.
 For every $0\leq h < k \leq n$, $S_{n,k}$ is a submanifold 
 in the product $S_{n,h}\times S_{n,k-h}$ and denote by
 $$h_{k,h}: S_{n,k} \to S_{n,h}$$
 the restriction of the natural projection onto the first factor.
 This map is {\it equivariant} for the above actions of $O(n)$
 on both Stiefel manifolds. 
 Clealy $h_{k,h}(I_{n,k})=I_{n,h}$. The fibre $h_{k,h}^{-1}(I_{n,h})$
 over $I_{n,h}$ is made by the  $2\times 2$ block matrices of the form
$$Y(D):=
\begin{pmatrix}
I_h&0\\
0&D
\end{pmatrix}
$$
where $D \in S_{n-h,k-h}$.
If $P\in O(n)$ is such that $PI_{n,h}=A$, then $P(h_{k,h}^{-1}(I_{n,h}))= h_{k,h}^{-1}(A)$,
all fibres are diffeomorphic to each other.
Let $\Omega$ be the open neighbourhood of $I_{n,h}$ in $S_{n,h}$
made by matrices of the form
$$ X=
\begin{pmatrix}
B\\
R
\end{pmatrix}
$$
where $B\in O(h)$. 
Recall the ``Gram-Schmidt''  retractions $\rG_{n,k}$ defined above. Then
we  define the smooth  map $X\to P(X)\in O(n)$
$$ P(X)=  \rG_{n,n} (
\begin{pmatrix}
B&0\\
R&I_{n-h}
\end{pmatrix})
$$
such that $P(X)I_{n,h}=X$.
Finally we have the following commutative diagram of smooth maps 
$$ \begin{array}[c]{ccc}
\Omega \times S_{n-h,k-h} &\stackrel{ \Psi }{\rightarrow}& h_{k,h}^{-1}(\Omega)\\
\downarrow\scriptstyle{\pi_\Omega}&&\downarrow\scriptstyle { h_{k,h}}\\
\Omega&\stackrel{{\rm id}_\Omega}{\rightarrow}&\Omega \end{array}$$ 
such that the first row is the diffeomorphism defined by
$$ (X,D)\to P(X)Y(D) \ . $$
The costant section of the product on the left, $X\to (X,0,I_{n-h,k-h})$
is transformed into the section of $h_{k,h}$ over $\Omega$
$$ s(X)= \rG_{n,k}(
\begin{pmatrix}
B&0\\
R&I_{k-h}
\end{pmatrix})
$$
 
 Summing up:
 
 \medskip
 
 {\it All these restriction of natural projections  onto Stiefel manifolds are
 locally trivial(izable) fibrations with a transitive action of either the group
 GL$(n,\R)$ or $O(n)$ respectively, which sends fibres into fibres.
 In the case of othogonal Stiefel manifolds, the fibre is a Stiefel manifold
 itself.}
 
 \medskip

 $\bullet$ A case of particular interest is when $n=k$. In the linear case we have a fibration
 of the linear group GL$(n,\R)$ over $L_{n,h}$ with fibre the {\it subgroup}
 of GL$(n,\R)$ made by the matrices of the form
 $$Y(S,D):=
\begin{pmatrix}
I_h&S\\
0&D
\end{pmatrix}
$$
where $(S,D)\in M(h,n-h,\R) \times {\rm GL}(n-h,\R)$.
 
 In the orthogonal case  we have a fibration of the othogonal group $O(n)$
 over $S_{n,h}$ with fiber  the  orhogonal group $O(n-h)$. Sometimes
 this is summarized by writing
 $$ S_{n,h}= O(n)/O(n-h) \ . $$
 \medskip
 
 $\bullet$ Another useful fibration is $h_{k,1} : S_{n,k} \to S^{n-1}$
 with fibre $S_{n-1,k-1}$.
 \medskip
 
 $\bullet$ Recall that $O(n)$ has two connected components and that
 the component containing $I_n$ is the special othogonal group
 $SO(n)$. If $h<n$, also the action of $SO(n)$ on $S_{n,h}$ is transitive,
 hence we can specialize all the discussion obtaining a fibration
 $$sh_{n,h}: SO(n)\to S_{n,h}$$
 with fibre  $SO(n-h)$, so that
 $$S_{n,h}=SO(n)/SO(n-h)$$
 in particular this implies that
 \medskip
 
 {\it For $h<n$, the Stiefel manifold
 $S_{n,h}$ is connected.}
 
\section{Grassmann manifolds}\label{Grass}
 
 For every $(n,k)$ as above, we are going to define now the
 {\it Grassmann manifold} $\GG_{n,k}$.
 
 Denote by $G_{n,k}$ the {\it set} of linear subspaces of $\R^n$ of
 dimension $k$.  Let $\GG_{n,k}$ be the closed subset of
 $S(n,\R)=\R^{\frac{n(n+1)}{2}}$ defined by the polynomial matrix equations
 $$  A^2-A=0, \ {\rm trace}(A)=k \ . $$
 If $A\in S(n,\R)$ verifies $A^2-A$ then its spectrum of eigenvalues is $\{0,1\}$, and
 by the spectral theorem for real symmetric matrices,
 the respective eigenspaces provide an orthogonal direct sum decomposition of $\R^n$;
 the last condition on the trace
 is equivalent to the fact that the eigenspace for the eigenvalue $\lambda=1$ has dimension
 equal to $k$, and also to the fact that $A$ has rank equal to $k$.
 
 We fix a bijection $V\to A_V$ from $G_{n,k}$ onto $\GG_{n,k}$ as
 follows.  For every $V\in G_{n,k}$ we have the orthogonal direct sum
 decomposition $\R^n= V \oplus V^\perp$, ($V^\perp$ being the
 orthogonal space to $V$ with respect to the standard euclidean scalar
 product) and the linear map $A_V\in \Ll(\R^n,\R^n)=M(n,\R)$ such that
 $A_V(v+v')=v$.  One readly verifies that $A_V\in \GG_{n,k}$. The
 inverse map $A\to V_A$ is defined by setting $V_A$ equal to the
 eigenspace of $A$ relative to the eigenvalue $\lambda=1$.
 
 Next we prove that $\GG_{n,k}$ is an embedded smooth manifold in
 $S(n,\R)$, of dimension $k(n-k)$. Note that the action by smooth
 diffeomorphisms of $O(n)$ on $S(n,\R)$ given by $(P,A)\to P^tAP$,
 restricts to an action on $\GG_{n,k}$: for every $A\in \GG_{n,k}$
 $(P^tAP)^2 -P^tAP=P^t(A^2-A)P=0$; as $P^t=P^{-1}$, then ${\rm
   trace}(PAP^{-1})=k$ because the trace is invariant up to
 conjugation.  This action corresponds via the above bijection $V\to
 A_V$ to the action of $O(n)$ on the set $G_{n,k}$ defined by
 $(P,V)\to PV$.  These actions are transtive, hence for every $A\in
 \GG_{n,k}$ there exists $P\in O(n)$ such that $P^tAP=H$
where $H$ is the $2 \times 2$ block diagonal matrix

$$
H= \begin{pmatrix}
I_k&0\\
0&0
\end{pmatrix}
$$

So it is enough to find a chart of
 $\GG_{n,k}$ at $H$.  First note that the space of symmetric matrices
 of rank $k$ ( denote it by $S(n|k,\R)$) is a submanifold of $S(n,\R)$
 of dimension $\frac{k(k+1)}{2}+k(n-k)$.  A local parametrization of
 $S(n|k,\R)$ at $H$ is given by
 $$ (S(k,\R)\cap {\rm GL}(k,\R)) \times M(k,n-k,\R)\to
 \Ww \subset S(n|k,\R), \ (D,B)\to Z(D,B) $$
where $Z(D,B)$ is the $2\times 2$ block symmetric matrix
$$
Z(D,B)= \begin{pmatrix}
D&B\\
B^t&B^tD^{-1}B
\end{pmatrix}
$$

To see that 
$Z(D,B)$ is of rank $k$, consider the non singular matrix
$$ X(D,B)=
\begin{pmatrix}
I_k&0\\
-B^tD^{-1}&I_{n-k}
\end{pmatrix} $$
then
 $$X(D,B)Z(D,B)=
 \begin{pmatrix}
D&B\\
0&0
\end{pmatrix} $$

\noindent This last matrix has the same rank of $Z(D,B)$ and this is equal to ${\rm rank}(D)=k$.
The same argument shows that if one changes the second block along the diagonal of $Z(D,B)$
by any one different from $B^tD^{-1}B$, then the resulting matrix would have rank $>k$.   
Clearly $Z(I_k,0)=H$. Hence $\Ww\cap \GG_{n,k}$ is given
by restriction to $\Ww$ of the matrix equation $A^2-A=0$. The matrix equation carried by the first
$k\times k$ block along the diagonal reads:
$$BB^t+D^2 -D=0$$
and by replacing $BB^t=D-D^2$ into the equations carried by the other blocks,
a direct computation shows that they are automatically satisfied.
We are reduced to study the map
$$h:( S(k,\R)\cap {\rm GL}(k,\R)) \times M(k,n-k,\R)\to S(k,\R), \ (D,B)\to BB^t+D^2-D$$
which is a summersion at $(I_k,0)$; hence, possibly shrinking $\Ww$, we conclude
that $Z(h^{-1}(0))=\Ww\cap \GG_{n,k}$ is an embedded smooth manifold of dimension
$k(n-k)$.  

An alternative way to get the same conslusion is to provide a local
parametrization of $\GG_{n,k}$ at $H$.  Let $\tilde U$ be the subset
of $G_{n,k}$ formed by the $k$-linear subspaces $V$ of $\R^n=\R^k
\times \R^{n-k}$ such that $V\cap \R^{n-k}=\{0\}$. Every $V\in \tilde
U$ is the graph of a uniquely determined linear map $L_V: \R^k \to
\R^{n-k}$. In fact the restriction to $V$ of the projection onto
$\R^k$ is a linear isomorphism; hence the inverse isomorphism is of
the form $x \to (x,L_V(x))$. Then $\tilde U$ can be identified with
$M(n-k,k,\R)$.  The restriction to $M(n-k,k,\R)$ of the above map
$V\to A_V $ can be explicitely computed as follows. For every $L\in
M(n-k,k,\R)$, let $V=V_L$ be the graph of $L$.  Consider the ordered
basis of $\R^n$
$$\Bb_L=\{(e_1,L(e_1)),\dots , (e_k,L(e_k)),e_{k+1},\dots, e_n) \} $$
such that the first $k$-vectors form a basis of $V$. Apply to $\Bb_L$
the Gram-Schmidt orthogonalization algorithm which produces an
orthonormal basis $\Dd_L$ of $\R^n$, whose first $k$ vectors are a
orthonormal basis of $V$ and the last $n-k$ of $V^\perp$.  By
organizing as usual $\Dd_L$ in a $n\times n$ matrix, we get $P_L\in
O(n)$. Finally $A_L=A_V = P_L^tHP_L$. The map $L \to A_L$ is clearly
smooth; by a bit of direct computation we see that it is indeed an
immersion. Finally, if $\Omega$ is a sufficiently small neighbourhood
of $H$ in $S(n,\R )$, and $W= \Omega \cap \GG_{n,k}$, then for every
$A\in W$, $V_A$ belongs $\tilde U$; the restriction to $W$ of $ A \to
L_V$ is a chart of $\GG_{n,k}$ with values in a open neighbourhood $U$
of $0\in M(n-k,k)$.  We have eventually proved that $\GG_{n,k}$ is an
embedded smooth manifold of dimension $k(n-k)$ in $S(n,\R)$.
 
\section{Stiefel manifolds as fibre bundles over Grassmann manifolds}
\label{stif-gram-fib}
 There are natural surjective maps
 $$ \lG_{n,k}: L_{n,k}\to \GG_{n,k}$$
 $$ s_{n,k}: S_{n,k}\to \GG_{n,k}$$ defined in both cases by $B \to
 A_{[B]}$ where $[B]$ denotes the linear $k$-subspace of $\R^n$
 generated by the columns of $B$.

 Let us concentrate on the map
 $s_{n,k}$.  Note that $[B]=[C]$ if and only if there exists $Q\in
 O(k)$ such that $C=BQ$, and that $A_{[B]}=H$ if and only if it is of
 the form
 $$ B=
\begin{pmatrix}
Q\\
0
\end{pmatrix}, \ Q\in O(k)  \ . $$
It follows that every fibre of $s_{n,k}$ is diffeomorphic to $O(k)$
and there is a transitive action (on the {\it right}) of $O(k)$ itself
on every fibre.
  
  The map $s_{n,k}$ is {\it equivariant} with respect to the actions
  of $O(n)$: $(P,B)\to PB$ on $S_{n,k}$, $(P,A)\to P^tAP$ on
  $\GG_{n,k}$, respectively.  Recall that `equivariant' means that for
  every $(P,B)$, $A_{[PB]}= P^tA_{[B]}P$. Then it is enough to analyse
  the behaviour of the restriction of the map to the inverse image
  $\tilde \Omega:= s_{n,k}^{-1}(\Omega)$ (which is a open
  neighbourhhod of $J$ in $S_{n,h}$) of some open neighbourhood
  $\Omega$ of $H$ in $\GG_{n,h}$.  For every $B\in S_{n,k}$, if $P$ is
  the top $k\times k$ submatrix of $B$, let us express this by writing
  $B=(P|D)$. Let $\tilde \Omega$ be the open neighbourhood of $J$ in
  $S_{n,k}$ formed by the matrices $B=(P|D)$ such that $P$ is non
  singular. If $B\in \tilde \Omega$ then $[B]\cap \R^{n-k} =\{0\}$,
  hence its image say $\Omega$ in $\GG_{n,k}$ is an open
  set. Moreover, If $[(P|D)]=[(R|S)]$, then there is $Q\in O(k)$ such
  that $(P|D)=(RQ|SQ)$. If $P$ is non singular, then also $R$ is
  necessarily non singular. This means that $\tilde \Omega=
  s_{n,k}^{-1}(\Omega)$ that is a {\it satured} open set of $S_{n,k}$
  with respect to the surjective map $s_{n,k}$.  We can make explicit
  $s_{n,k}(B)$ on $\tilde \Omega$ by applying to every $[B]$ and its
  orthonormal basis given by $B$ itself the construction already used
  above in order to construct a local parametrization of $\GG_{n,k}$
  at $H$.  This shows that $s_{n,k}$ is smooth. Moreover, define
  $\phi: \tilde \Omega \to M(k,n-k,\R)$ by $\phi((P|D))=DP^{-1}$. If
  $(P|D)=(RQ|SQ)$ as above, then $SQQ^{-1}R^{-1}=SR^{-1}$. Then there
  is an induced smooth map $\Omega \to M(k,n-k,\R)$ whose inverse map
  is
$$\psi: M(k,n-k) \to \Omega, \ \psi(Z)= A_{[\rG_{n,k}(I_k|Z)]}$$
  providing once again a local parametrization of $\GG_{n,k}$ at $H$.
  We can summarize this discussion by saying that there is a locally
  trivializing commutative diagram at $H$
$$ \begin{array}[c]{ccc} \Omega \times O(k) &\stackrel{ \Psi
    }{\rightarrow}& \tilde
    \Omega\\ \downarrow\scriptstyle{\pi_\Omega}&&\downarrow\scriptstyle
             { s_{n,k}}\\ \Omega&\stackrel{{\rm
                 id}_\Omega}{\rightarrow}&\Omega \end{array}$$ where
  $\Psi(A,Q)=\psi(Z)Q$, $A=\psi(Z)$.  Its orbit by the action of
  $O(n)$ provides a fibred atlas for the summersion $s_{n,k}$.
  Summing up we have proved:
 \begin{proposition}\label{s-g-fib} The map $s_{n,k}: S_{n,k}\to \GG_{n,k}$
   is a fiber bundle with fibre $O(k)$.  Every change of
   trivialization
   $$ \Phi' \circ \Psi(\Omega \cap \Omega')\times 0(k) \to
   (\Omega \cap \Omega')\times O(k)$$
 is of the form
 $$ (p,P)\to (p, PQ(p))$$
 where $p\to Q(p)$ defines a smooth map $\Omega \cap \Omega' \to O(k)$.
 \end{proposition}
 
 We have also the following topological corollaries
 
 \begin{corollary} Every $\GG_{n,k}$ is a compact and connected embedded smooth manifold. 
 As a topological space it has the quotient space topology
 $S_{n,k}/s_{n,k}$.
 \end{corollary}  
  
  \medskip

$\bullet$ {\it Real projective spaces.}  A particular case of the
  above discussion is when $k=1$. In such a case $\GG_{n,1}$ is also
  denoted by $\PP^{n-1}(\R)$ and called the (real) {\it
    $(n-1)$-projective space}.  $S_{n,1}=S^{n-1}$, and the map
  $s=s_{n,1}:S^{n-1}\to \PP^{n-1}(\R)$ is a {\it smooth covering map
    of degree $2$}.

\medskip

$\bullet$ {\it Complex Stiefel and Grassmann manifolds.}  As a smooth
manifold $\C^n= \R^{2n}$, hence $M(n,\C)$ is a submanifold of
$M(2n,\R)$ etc. All along the above discussion let us replace:
\begin{itemize}
\item $\R^n$ with $\C^n$.
  The real linear subspaces of $\R^n$ with the {\it complex} linear
subspaces of $\C^n$.
\item The standard positive definite scalar product on $\R^n$ with the standard positive definite
{\it Hermitian product} on $\C^n$, $<v,w>= v^t\bar w$. 
\item The (real) orthogonal groups $O(n)$ with the {\it unitary groups} 
$$U(n):=\{A\in {\rm GL}(n,\C); \ A^{-1}=A^*:= \bar A^t\} \ . $$
\item The spaces of real symmetric matrices $S(n,\R)$  with the spaces of {\it Hermitian matrices}
$$H(n, \C)=\{A\in M(n,\C); \ A = A^* \} \ . $$
\item The spectral theorem for real symmetric matrices with the {\it spectral theorem for
complex hermitian matrices}.
\end{itemize}

\noindent Then by repeating verbatim the above constructions, for every $(n,k)$ as above, we realize the (unitary) 
{\it complex Stiefel manifold} $S_{n,k}(\C)$ as a compact embedded smooth
manifold in $M(n,k,\C)$, the complex Grassmannian manifold $\GG_{n,k}(\C)$ as a compact embedded
smooth manifold in $H(n,\C)$ (defined by the usual equations $A^2-A=0, \ {\rm trace}(A)=k$),
the complex projective spaces $\PP^{n-1}(\C)=\GG_{n,1}(\C)$, 
and so on. Although we are dealing with spaces based  on the {\it complex numbers},
we stress that in this way we have actually realized them 
as {\it real} embedded smooth manifolds.

\medskip

{\it We understand that also all next considerations about Stiefel and Grassmann
manifolds would have a counterpart for the complex version.}

\section{ A cellular decomposition of the Grassmann manifolds}
\label{cell-decomp}
 We describe a natural partition of $\GG_{n,k}$ by a finite number of
 subsets each one diffeomorphic to some $\R^h$, $0\leq h \leq \dim
 \GG_{n,k}$, (i.e. an {\it open $h$-cell}) and such that its closure
 in $\GG_{n,k}$ is union of cells of lower dimension. Let $L\in
 \GG_{n,k}$, that is $L$ is a $k$-dimensional linear subspace of
 $\R^n$ (here we confuse $G_{n,k}$ and $\GG_{n,k}$).  For every
 $i=0,\dots, n$, denote by
 $$p_i: \R^n \to \R^{n-i}$$ the projection onto the first $n-i$
 coordinates, $p_i((x_1,\dots,x_n)^t)=(x_1,\dots, x_{n-i})^t$.  The
 dimensions of $p_i(L)\subset \R^{n-i}$ decrease from $k$ to $0$ in
 exactly $k$ steps; that is there are integers
 $$1\leq \sigma_1<\sigma_2 < \dots < \sigma_k\leq n$$
 such that for $j$ that decreases from $k$ to $1$,
 $$\dim p_{\sigma_j+1}(L) -\dim p_{\sigma_{j}}(L)=1 \ . $$
Then
 $$\sigma(L):=(\sigma_1,\dots,\sigma_k)$$ is called the {\it Schubert
  symbol} of $L$. There is a concrete elementary way to determine
$\sigma(L)$:
 \smallskip
 
 $\bullet$ Fix any rank $k$, $n\times k$ matrix $A\in L_{n,k}$ which
 projects to $L\in \GG_{n,k}$.
\smallskip
 
 $\bullet$ Apply to $A$ the Gauss algorithm via elementary operations
on the columns and get a matrix $$\hat A \in L_{n,k}$$ in {\it column
  echelon form} which also projects to $L$. So for every $j=1,\dots, k$,  
  the $(\sigma_j,j)$ entry of $\hat A$ is equal to $1$ and is a so
called `pivot' of $\hat A$; the (transposed of the) $\sigma_j$th row of
$\hat A$ is the $\sigma_j$th vector of the standard basis of $\R^k$;
beyond the pivots, for every $1\leq j \leq k$, an $(s,j)$ entry of $\hat A$ is possibly non zero
only if $\sigma_j < s \leq n$ and $s$ is not the row index of any
pivot row. The computation of
$\sigma(L)$ by means of $\hat A$ is
immediate from the very definition. This means in particular that the
initial choice of the matrix $A$ is immaterial to this computation of
$\sigma(L)$; $\sigma(\hat A):=\sigma(L)$ is the symbol of the matrix
$\hat A$ and two matrices in column echelon form have the same index
if and only if they share the pivot positions. 

We claim furthermore that the whole matrix $\hat A$ does not depend on
the choice of $A$ as it is completely determined by $L$.  For, given
$\sigma = \sigma (L)$, denote by $p_\sigma$ the projection of $\R^n$
onto the $k$ coordinates $(x_{\sigma_1},\dots, x_{\sigma_k})$; then
the restriction of $p_\sigma$ to $L$ is a linear isomorphism and the
columns of $\hat A$ are characterized as the vectors of $L$ which are
mapped in the order by $p_\sigma$ to the vectors $e_{\sigma_1},\dots,
e_{\sigma_k}$ of the standard basis of $\R^k$.
 
Summarizing, there are $\binom{n}{k}$ Schubert symbols. For every such
a symbol $\sigma$, the subset $C_\sigma$ of $\GG_{n,k}$ formed by the
$k$-subspaces of $\R^n$ which share the symbol $\sigma$ is in
bijection with the subset $\hat C_\sigma$ of $L_{n,k}$ formed by the
matrices in columns echelon forms which also share the symbol $\sigma$. 
$\hat C_\sigma$ has a natural base point, that is the matrix $J_\sigma$
whose entries different from the pivots are zero; then
$$ \hat C_\sigma = J_\sigma + \VV_\sigma$$
and it is easy to check that $\VV_\sigma$ is a linear subspace of
$M(n,k,\R)$ formed by the matrices with a given pattern of zero
entries determined by the symbol $\sigma$. The other entries
contain free parameters. By counting the free parameters 
column by column, we readly verify that 
$$d_\sigma := \dim \VV_\sigma = \sum_{j=1}^k \left (n-\sigma_j -(k -j)\right) \ . $$
It follows that $C_\sigma \subset \GG_{n,k}$ admits a smooth parametrization
$$ \psi_\sigma: \R^{d_\sigma}\to C_\sigma \ . $$
By varying the symbols we have obtained
a partition of $\GG_{n,k}$ by open cells. We claim that:
\smallskip

{\it The closure of every $C_\sigma$ in $\GG_{n,k}$
is formed by the $C_{\sigma'}$'s such that for every $j$, $\sigma'_j
\geq \sigma_j$}.
\smallskip

This claim is not obvious. We omit the proof, however next item 4)
should help the reader to reconstruct such a proof.

{\bf Remarks and examples.}

1) There is one top dimensional (i.e. of dimension $k(n-k)$) 
cell of $\GG_{n,k}$ corresponding to
the symbol $(1,2,3,\dots, k)$. This covers a chart around 
the image of $I_{n,k}$ in $\GG_{n,k}$. In general every cell
$C_\sigma$ has a natural base point, that is the image in $\GG_{n,k}$ 
of the  the matrix $J_\sigma \in \hat C_\sigma$. There is one $0$-cell corresponding to
the symbol $(n-k+1, n-k,\dots, n)$.

2) In the case of projective spaces $\PP^n(\R)=\GG_{n+1,1}$, there are $n+1$ cells, one cell for
every dimension $n, \dots, 0$ corresponding to the symbols $(1)$,$(2),\dots, (n+1)$.
The closure of every cell of dimension $d$ say is a copy of $\PP^d(\R)$ linearly
embedded into $\PP^n(\R)$.
\smallskip

3) For example  $\GG_{4,2}$ has six cells corresponding to the Schubert symbols 
$(1,2)$, $(1, 3)$, $(1, 4)$, $(2, 3)$, $(2, 4)$, $(3, 4)$, and these cells have dimensions 
$4, 3, 2, 2,1,0$  respectively. 
\smallskip

4) The cells of $\GG_{n,k}$ can be described also in terms of the orthogonal
Stiefel manifold $S_{n,k}$. A matrix $\tilde A\in S_{n,k}$ is in {\it orthogonal} 
column echelon form of symbol $\sigma$ if its standard column echelon form
$\hat A$ is of symbols $\sigma$ and $\tilde A$ may differ from $\hat A$ only by:
1) the pivot entries of $\tilde A$ are non zero not necessarily equal to $1$;
2) the entries of a pivot row of $\tilde A$ on the left of the pivot are not necessarily
equal to $0$; 3) the last non zero entry of every column is positive.
One can verify that for every $L\in \GG_{n,k}$ there is only one $\tilde A \in S_{n,k}$
which projects to $L$; in fact if $\hat A$ is the unique matrix in standard column echelon
form which projects to $L$, then we can  obtain $\tilde A$  by applying the Gram-Schmidt
algorithm to the columns of $\hat A$  considered in the backward order (normalized
to achive also the condition 3) above). The subset $\tilde C_\sigma$ of $S_{n,k}$
formed by the matrix in echelon form of symbol $\sigma$ is 
diffeomorphic to $\hat C_\sigma \subset L_{n,k}$ and 
 maps diffeomorphically onto $C_\sigma \subset \GG_{n,k}$. 
 One can prove that the closure of $\tilde C_\sigma$
in $S_{n,k}$ is diffeomorphic to a {\it closed} disk of dimension $d_\sigma$
which maps onto the closure of $C_\sigma$ in $\GG_{n,k}$.
\smallskip  

5) Referring to Section \ref{LIMITS}, the cell decompositions respect the
inclusions $$j_n:\GG_{n,k}\to \GG_{n+1,k}$$ in the sense that the cells
of $\GG_{n,k}$ are also cells of $\GG_{n+1,k}$; hence we have also a cell
decomposition of the limit infinite Grassmannian $\GG_{\infty,k}$.

\section{Stiefel and Grassmannian manifolds as regular real algebraic sets}
\label{algebraic-set}
For the notions and basic results of (real) algebraic geometry mentioned  
in this section we can refer for example to \cite{BCR} or to \cite{BR}.
\medskip
 
By definition a {\it real algebraic set} $Z\subset \R^m$, for some
$m\in \N$, is of the form $Z=F^{-1}(0)$ for some {\it polynomial} map
$F: \R^m \to \R^h$. Hence the Stiefel and Grassmannian manifolds (even
in the complex version) are also examples of real algebraic sets. We
are going to outline a way to recover that they are embedded smooth
manifolds by the means of algebraic geometry, obtaining indeed a
stronger result.

For every algebraic set $Z$ as above,
$$I(Z):= \{p(X)\in \R[X_1,\dots,X_m]; \ p(x)=0 \ {\rm for \ every \ }
x\in Z\}$$ is called the (defining) {\it ideal of $Z$}. By a theorem
of Hilbert, $I(Z)$ is {\it finitely generated}, that is there exist
some polynomials $p_1(X),\dots , p_k(X) \in I(Z)$ such that $I(Z)$
coincides with the set of linear combinations of the $p_j(X)$'s with
polynomials coefficients in $\R[X_1,\dots,X_m]$.  Consider the
polynomial map
 $$P: \R^m\to \R^k, \ P(x)=(p_1(x),\dots ,p_k(x)) \ . $$
 For every $p\in Z$, set 
 $$r(p)= {\rm rank}\  d_pP \ . $$
 It is not too hard to show that $r(p)$ {\it does not depend on the choice
 of the generators $p_1,\dots,p_k$}.
 So it is well defined
 $$ r(Z)= \max \{r(p); \ p \in Z \} \ . $$
 
 Assume for simplicity that $Z$ is {\it irreducible} that is it cannot
 be expressed as $Z=Z_1\cup Z_2$ where $Z_1$ and $Z_2$ are algebraic
 sets both different from $Z$ (one can prove that the connected
 Stiefel and Grassmannian algebraic sets are irreducible - that is all
 with the exception of the othogonal groups $O(n)$). Then $p\in Z$ is
 a {\it regular point} if $r(p)=r(Z)$. Note that by the definition,
 the set $R(Z)$ of regular points of $Z$ is {\it non empty}. A {\it
   Zariski open set} in $\R^m$ is of the form $\R^m \setminus Y$ where
 $Y$ is an algebraic set in $\R^m$.  The following is a non trivial
 result.
 \begin{theorem}
Let $Z\subset \R^m$ be an irreducible algebraic set of rank $r=r(Z)$.
Then for every $p\in R(Z)$ there exist a Zariski open set $U$ of
$\R^m$ and a polinomial map $F=(F_1,\dots, F_r ): \R^m \to \R^r$ such
that:
\begin{enumerate}
\item $p\in U$.
\item $F_j\in I(Z)$, $j=1,\dots, r$.
\item $Z\cap U = U\cap F^{-1}(0)$.
\item For every $x\in U\cap Z$,
$$ {\rm rank} \ d_xF = r \ . $$
\end{enumerate}

In particular $R(Z)$ is an embedded smooth manifold in $\R^m$ of
dimension $m-r$.
\end{theorem} 

Assuming this fundamental theorem, we can prove

\begin{corollary} Let $Z\subset \R^m$ be one of our favourite
  (Stiefel or Grassmannian) algebraic sets.
 Then $Z=R(Z)$. In particular $Z$ is an embedded smooth manifold of
 dimension $m-r(Z)$.
 \end{corollary}
 \Dim We know that $R(Z)\neq \emptyset$. Let $p\in R(Z)$. By using the
 suitable transitive action on $Z$ of orthogonal (unitary) groups, we
 realize that for every $q\in Z$ there is a particularly simple linear
 diffeomorphism $\phi: \R^m\to \R^m$ such that $\phi (Z)=Z$ and
 $\phi(p)=q$. Although this is a particular case of a general result
 on the invariance of $R(Z)$ up to ``algebraic isomorphism'', these
 diffeomorphisms are so simple that one can check directly that since
 $p$ is regular then also $q$ is regular. Then $Z=R(Z)$.

 \cvd
 \medskip
 
 Note that the linear Stiefel manifolds are in fact Zariski open sets
 of the pertinent matrix space.
 
 \medskip
 
 \begin{remark}\label{algebraic-map}
   {\rm  We stress that the notion of regular point is rather a
     delicate one. For example it can happen that for some irreducible
     algebraic set $X\subset \R^m$ which is an embedded smooth
     manifold, nevertheless $R(X)\neq X$.}
 \end{remark}

\chapter{Tautological bundles and pull-back}\label{TD-EMB-VB}
The basic notions about fibred bundles have been already introduced
in Section \ref{emb-fib-bundle}, and we will use them.
The tensorial vector bundles and their relatives, defined in Chapter \ref{TD-DIFF}
belong to a wide category of ``embedded fibred bundles''
constructed via the {\it pull-back} of
{\it tautological bundles} over
Grassmann manifolds.  
We are going to state these matters.

\section{Tautological bundles}\label{tautological}
We are going to construct so called {\it tautological fibre bundles}
over the grassmannian $\GG_{n,k}$.
\medskip

$\bullet$ {\it (The tautological vector bundle)} 
Define
$$ \Vv(\GG_{n,k})= \{(A,v)\in \GG_{n,k}\times \R^n; \ v\in V_A\}$$
i.e. $v$ belongs to the $k$-linear subspace $V$ of $\R^n$
such that $A=A_V$, via the usual bijection $G_{n,k}\cong \GG_{n,k}$.
The restriction of the projection onto the first factor defines the smooth
surjective map
$$ \tau_{n,k}: \Vv(\GG_{n,k})\to \GG_{n,k} \ . $$
It is clear that for every $A\in \GG_{n,k}$, the inverse
image $\tau_{n,k}^{-1}(A)=V_A$. We have
\begin{proposition}\label{taut-vect}  $ \tau_{n,k}: \Vv(\GG_{n,k})\to \GG_{n,k} $
is an embedded smooth vector bundle with fibre $\R^k$. It is called the 
{\rm tautological vector bundle} over $\GG_{n,k}$.
\end{proposition}

\medskip

$\bullet$ {\it (The tautological linear frame  bundle)}
Define
$$ \Ll(\GG_{n,k})= \{(A,X)\in \GG_{n,k}\times L_{n,k}; \ \lG_{n,k}(X) = A \}$$
i.e. $X$ spans the $k$-linear subspace $V$ of $\R^n$
such that $A=A_V$. 
The restriction of the projection onto the first factor defines the smooth
surjective map
$$ l \tau_{n,k}: \Ll(\GG_{n,k})\to \GG_{n,k} \ . $$
It is clear that for every $A\in \GG_{n,k}$, the inverse
image $l\tau _{n,k}^{-1}(A)$ consists of all {\it linear frames}
of $V_A$. We have
\begin{proposition}\label{taut-lin-frame}  $ l \tau _{n,k}: \Ll(\GG_{n,k})\to \GG_{n,k} $
is an embedded smooth fibre bundle with fibre GL$(k,\R)$. It is called the 
{\rm tautological linear frame bundle} over $\GG_{n,k}$.
\end{proposition}

$\bullet$ {\it (The tautological orthogonal frame  bundle)}
Define
$$ \Ss(\GG_{n,k})= \{(A,X)\in \GG_{n,k}\times S_{n,k}; \ s_{n,k}(X) = A \}$$
i.e. $X$ spans the $k$-linear subspace $V$ of $\R^n$
such that $A=A_V$. 
The restriction of the projection onto the first factor defines the smooth
surjective map
$$ s\tau _{n,k}: \Ss(\GG_{n,k})\to \GG_{n,k} \ . $$
It is clear that for every $A\in \GG_{n,k}$, the inverse
image $s\tau_{n,k}^{-1}(A)$ consists of all {\it orthonormal frames}
of $V_A$. We have
\begin{proposition}\label{taut-ort-frame}  $ s\tau _{n,k}: \Ss(\GG_{n,k})\to \GG_{n,k} $
is an embedded smooth fibre bundle with fibre $O(k)$. It is called the 
{\rm tautological orthogonal frame bundle} over $\GG_{n,k}$.
\end{proposition}





{\it Proofs:} Let us prove Proposition \ref {taut-vect}. Recall that $\GG_{n,k}$
is endowed with an atlas $\{(\Omega_V, \phi_V)\}_{V\in G_{n,k}}$
where  
$$\Omega_V=\{ A\in \GG_{n,k}; \ V_A\cap V^\perp = \{0\} \}$$
equivalently, $V_A$ is the graph of a uniquely determined
linear map $L_A: V\to V^\perp$. Set as usual
$\tilde \Omega_V= \tau_{n,k}^{-1}(\Omega_V)$.
Then a vector bundle atlas of $\tau_{n,k}$ is given by
the locally trivializing commutative diagrams ($V$ varying in $G_{n,k}$, 
$\Bb=\{v_1,\dots,v_k\}$ varying in
the linear frames of $V$)

$$ \begin{array}[c]{ccc}
\Omega_V \times \R^k &\stackrel{ \Psi_\Bb }{\rightarrow}&  \tilde \Omega_V\\
\downarrow\scriptstyle{\pi_{\Omega_V}}&&\downarrow\scriptstyle { \tau_{n,k}}\\
\Omega_V&\stackrel{{\rm id}_{\Omega_V}}{\rightarrow}&\Omega_V \end{array}$$ 
where 
$$\Psi_\Bb(A,x)= (A, \sum_{i=1}^k x_iv_i + \sum_{i=1}^k x_iL_A(v_i)) \ . $$
It is immediate that for every couple $(V,\Bb)$, $(V', \Bb')$ there is a smooth
map
$$ \lambda_{\Bb,\Bb'}: \Omega_V \cap \Omega_{V'}\to {\rm GL}(k,\R)$$
such that the corresponding change of local trivialization is of the form
$$ (\Omega_V\cap \Omega_{V'})\times \R^k \to    (\Omega_V\cap \Omega_{V'})\times \R^k$$
$$(A, v)\to (A,\lambda_{\Bb,\Bb'}(A)v) \ . $$
\medskip

{\bf Remark.} {\it By restricting to {\it orthogonal} frames $\Bb$ of the 
$V$'s, we get a sub-fibred atlas such that the change of local trivializations 
are governed by smooth maps 
$$ \lambda_{\Bb,\Bb'}: \Omega_V \cap \Omega_{V'}\to O(k) \ . $$}
\medskip

The proof of the other two propositions is similar and left to the reader.
Note that the change of local trivializations for the frame bundles are governed by {\it the same}
smooth maps $\lambda_{\Bb,\Bb'}$ as above, with values in GL$(k,\R)$
for $l\tau_{n,k}$, or in $O(k)$ for $s\tau_{n,k}$ respectively;
the groups ${\rm GL}(k,\R)$ or $O(k)$ act on themselves by left multiplication. 

\cvd   

\section{Pull-back}\label{pull-back}
We introduce a fundamental construction on embedded smooth fibred bundles. 
We state it in wide generality; later we apply it to the tautological bundles
of Section \ref{tautological}. 

Let us give an embedded smooth fibre bundle 
$$\xi:= f:E \to X$$
with fibres $E_x$ diffeomorphic to the manifold $F$ (recall section \ref{emb-fib-bundle}).

Let $g\in \Ee(M,X)$. Then set
$$ g^*E=\{ (p,y)\in M\times E; \ g(p)= f(y)\} \ $$
$$ g^*: g^*E\to E, \ g^*(p,y)=y $$
$$ g^*f:  g^*E \to M, \  g^*f(p,y)=p \ . $$
Obviously we have the commutative diagram  of smooth maps, denoted by $[g,g^*]$

$$ \begin{array}[c]{ccc}
g^*E &\stackrel{ g^* }{\rightarrow}&  E\\
\downarrow\scriptstyle{g^*f}&&\downarrow\scriptstyle { f}\\
M&\stackrel{g}{\rightarrow}& X \end{array}$$
Moreover, for every $p\in M$, $x=g(p)$, then $g^*E_p:= (g^*f)^{-1}(p)$ is equal to
the fibre $E_x$. Hence, also every $g^*E_p$ is diffeomorphic to $F$. In fact we have

\begin{proposition} (1) For every fibre bundle  $\xi:= f:E \to X$ with fibre $F$,
for every $g\in \Ee(M,X)$,
$$g^*\xi:=  g^*f: g^*E\to M$$ 
is an embedded smooth fibre bundle with fibre $F$. It is called the
{\rm pull-back of $\xi$ via $g$}. Moreover, $[g,g^*]$ is a fibred
map between fibred bundles.

(2)  For every $h\in \Ee(N,M)$, every $g\in \Ee(M,X)$, then
$$ (g\circ h)^*\xi = h^*(g^*\xi) \ . $$

(3) $$ (g\circ h)^* = g^* \circ h^* \ . $$ 

\end{proposition}
\Dim  The second and third points follow from the very definitions.
As for the first; consider a fibre bundle atlas of $\xi$.
This is formed as usual by locally trivializing diagrams
$$ \begin{array}[c]{ccc}
\Omega\times F &\stackrel{ \Psi}{\rightarrow}&  \tilde \Omega\\
\downarrow\scriptstyle{\pi_\Omega}&&\downarrow\scriptstyle { f}\\
\Omega&\stackrel{{\rm id}_\Omega}{\rightarrow}& \Omega \end{array}$$
and any change of local trivializations is of the form 
$$(\Omega\cap \Omega')\times F \to  (\Omega\cap \Omega')\times F$$
$$ (x,y)\to (x, \rho(x)(y)$$
$$x\to \rho(x)\in {\rm Aut}(F) \ . $$
The $\Omega$'s form an open covering of $X$. Fix an open covering
$\{W\}$ of $M$ such that $g(W)$ is contained in some $\Omega$.
The for every $W$ we have the locally trivializing commutative diagram
$$ \begin{array}[c]{ccc}
W\times F &\stackrel{ \Psi \circ (g,{\rm id}_F)}{\rightarrow}&  \tilde W\\
\downarrow\scriptstyle{\pi_W}&&\downarrow\scriptstyle { g^*f}\\
W&\stackrel{{\rm id}_W}{\rightarrow}& W \end{array}$$
The chage of local trivialization is of the form
$$(W\cap W')\times F \to  (W\cap W')\times F$$
$$ (w,y)\to (w, \rho(g(w))(y))$$
$$w\to \rho(g(w))\in {\rm Aut}(F) \ . $$ 

\cvd

\begin{remark}\label{G-pull-back}{\rm If $F$ has an additional structure
preserved by a subgroup $G \subset {\rm Aut}(F)$, and  $x\to \rho(x)$
as above is a smooth map with values in $G$ (i.e. $\xi$ is a ``$G$-bundle'')
then also the pull-back $g^*\xi$ has the same property. For example il $\xi$
is a vector bundle (with fibre $\R^k$) then also $g^*\xi$ is so.}
\end{remark} 


\section{Categories of vector bundles }\label{pull-back-vb}
Let $M$ be an embedded smooth manifold (possibly with boundary).
Let $$f:M \to \GG_{n,k}$$ be a smooth map. Then we can consider the
pull-back vector bundle $f^*\tau_{n,k}$, that is
$$ \begin{array}[c]{ccc}
f^*\Vv(\GG_{n,k}) &\stackrel{ f^* }{\rightarrow}&  \Vv(\GG_{n,k}) \\
\downarrow\scriptstyle{f^*\tau_{n,k}}&&\downarrow\scriptstyle { \tau_{n,k}}\\
M&\stackrel{f}{\rightarrow}& \GG_{n,k} \end{array} \ . $$
By the strict definition,  the total space of  ${\rm id}_{\GG_{n,k}}^*\tau_{n,k}$
is a submanifold of $\GG_{n,k}\times (\GG_{n,k}\times \R^n)$; however,  the projection onto the product of the first and third factors 
gives a canonical fibred diffeomorphim onto the total space of $\tau_{n,k}$. Modulo this {\it normalized embedding}, we can stipulate that
$${\rm id}_{\GG_{n,k}}^*\tau_{n,k} = \tau_{n,k}  \ . $$
Similarly, for every $f:M\to \GG_{n,k}$ as above, the total space of $f^*\tau_{n,k}$ has a canonical embedding
into $M\times \R^n$; modulo this normalization we can state that
$${\rm id}_M^*(f^*\tau_{n,k})=f^*\tau_{n,k} \ . $$
We stipulate that such a normalization is performed by default.
Note also that the composition of $f^*$ with the natural projection of $\Vv(\GG_{n,k})$ to $\R^n$
gives a map which is linear and injective at every fibre of $f^*(\Vv(\GG_{n,k}))$, from which we can reconstruct
tautologically the map $f$.

\medskip

Denote $\Nn=\{(n,k)\in \N \times \N; \ 0\leq k \leq n \}$.
For  every
$(n,k)\in \Nn$ set
$$ \Vv_{n,k}(M) := \{ f^*\tau_{n,k}; \ f\in \Ee(M,\GG_{n,k}) \} \ ; $$ 
and
$$ \Vv(M)= \cup_{(n,k)\in \Nn} \ \Vv_{n,k}(M) \ . $$

Then we see immediatly that

$$ M \ \Rightarrow \ \Vv(M)$$

$$g:N\to M \ \Rightarrow  \ g^{\bullet}: \Vv(M)\to \Vv(N), \ g^\bullet(f^*\tau_{n,k})= (f\circ g)^*\tau_{n,k}$$
so that 
$$(g\circ h)^{\bullet}= h^\bullet \circ g^\bullet$$
define
a contravariant functor from the category of embedded smooth manifolds (with boundary)
to this category of embedded smooth vector bundles.  Moreover, for every $f$ and every
$g$ as above there is the natural vector bundle map 
$$[g,g^*]: g^\bullet (f^*\tau_{n,k})\to f^*\tau_{n,k} \ . $$
If $g: N\to M$ is a diffeomorphism, then $g^\bullet : \Vv(M)\to \Vv(N)$ is a bijection
(with inverse $(g^{-1})^\bullet$), and for every $f$, $[g,g^*]$ is a vector bundle
isomorphism between $g^\bullet (f^*\tau_{n,k})$ and $f^*\tau_{n,k}$.

The tangent bundle of a manifold $M\subset \R^n$ as well all its tensorial relatives belong to $\Vv(M)$. 
For example $\pi_M: T(M)\to M$
is the pull-back of the (tautological) map 
$$ t_M: M \to \GG_{n,m}, \ t_M(p)=T_pM \ . $$
More generally we have

\begin{lemma}\label{emb-to-emb} If $\xi:= f:E\to M$ is a smooth vector bundle with fibre $\R^k$  such that the total space $E$
is a submanifold of some $\R^n$, and every fibre $E_x$ is a linear
$k$-subspace of $\R^n$, then the bundle $\xi$ belongs to $\Vv(M)$.
\end{lemma}

\Dim In fact $\xi$ is the pull-back of the (tautological) map 
$$e_M: M \to \GG_{n,k}, \ e_M(x)=E_x \ . $$

\cvd

\subsection{ Bundle equivalences}\label{B-equivalence}
We are going to refine the above constructions by introducing suitable quotient sets of $\Vv(M)$.

For every $f: M \to \GG_{n,k}$, and every inclusion $j_{n}: \GG_{n,k}\to \GG_{n+1,k}$ (see Section \ref {LIMITS}),  
 the total space of $(j_{n}\circ f)^*\tau_{n+1,k}$ 
 is  embedded in $M\times \R^n$ and coincides with the total space of $f^*\tau_{n,k}$.
 This gives us a {\it canonical identification} between these formally different points of $\Vv(M)$.
 A first mild quotient of $\Vv(M)$ is obtained by means of such canonical
 identifications. Let us keep for it the name $\Vv(M)$. For every equivalence class,
 there is one representative $f^*\tau_{n,k}$ with {\it minimum} $n$.
 
 More substantially we can restrict to $\Vv(M)$  the {\it full equivalence} between vector bundles
defined  in Section \ref{emb-fib-bundle}, generated by arbitrary vector bundle
isomorphisms of the form $[g,\tilde g]$. Denote by $\VV(M)$ the quotient set.

\begin{example} {\rm For example, if $g\in {\rm Aut}(M)$, then for every $f: M\to \GG_{n,k}$, the corresponding $[g,g^*]$ 
realizes a full equivalence between $f^*\tau_{n,k}$ and $g^\bullet(f^*\tau_{n,k})$.
By the way this establishes an action of ${\rm Aut}(M)$ on $\Vv(M)$, so that  $\VV(M)$
is a quotient set of $\Vv(M)/{\rm Aut}(M)$.}
\end{example}  

\medskip
 
We can restrict to $\Vv(M)$ the {\it strict equivalence} between vector bundles
defined  in Section \ref{emb-fib-bundle}, generated by isomorphisms of the form
$[{\rm id}_M,\tilde g]$. Denote by $\VV_0(M)$ the quotient set.
Clearly $\VV(M)$ is a quotient of $\VV_0(M)$.
 
\begin{example}\label{equiv-ex}{\rm   (i) If $f,g: M\to \GG_{n,k}$ are two different constant maps,
then $f^*\tau_{n,k}$ and $g^*\tau_{n,k}$ are different points of $\Vv(M)$ which obviously are strictly equivalent. 

(ii) Let  $g: M \to N$ be a diffeomorphism; then $[g^{-1},Tg^{-1}] \circ [g,g^*] $
is a strict equivalence between $T(M)$ and $g^*T(N)$.

(iii) By generalizing the above item, let  $[g,\tilde g]$ realize a full equivalence between bundles in $\Vv(M)$;
then also $[g,g^*]$ as in the above example realizes such an equivalence. Moreover, $[g^{-1}, \tilde g^{-1}]\circ [g,g^*]$  
realizes instead a strict equivalence.
}
\end{example}

\medskip

$\bullet$ {\it By associating to every $f^*\tau_{n,k}$
its class in the preferred quotient set of $\Vv(M)$, we get variants of the basic
pull-back functor defined above}.

\medskip

We will concentrate on $\VV_0(M)$. In particular we pose the following natural
question: set
$$ \Ee(M,\GG):=  \cup_{(n,k)\in \Nn} \ \Ee(M,\GG_{n,k}) \ . $$

\begin{question}\label{V0-question} {\rm Consider the obvious surjective map
$$(.)^*: \Ee(M,\GG) \to \VV_0(M), \ f \to [f^*\tau_{n,k}] \  $$
so that tautologically 
$$ \VV_0(M)= \Ee(M,\GG)/(.)^* \ . $$
This relation on $\Ee(M,\GG)$ is only implicitly defined.
The question is to {\it make it explicit}.  An answer will be discussed later when $M$ is compact.}
\end{question}

\section{The frame bundles}\label{frame-functor}
We can repeat the above scheme by using instead the tautological frame bundles.
It is enough to replace $\Vv(M)$ either with
$$\Ll(M) = \cup_{(n,k)\in \Nn} \ \Ll_{n,k}(M)$$
$$ \Ll_{n,k}(M):= \{f^*\l\tau_{n,k}; \ f\in \Ee(M, \GG_{n,k}) \} \ . $$
or the similarly defined $\Ss(M)$ and $\Ss_{n,k}(M)$ by using
the tautological bundles $s\tau_{n,k}$.
For every $f:M \to \GG_{n,k}$, the vector bundle $f^*\tau_{n,k}$
is {\it associated} to its {\it linear frame bundle}  $f^*l\tau_{n,k}$,
provided that both are considered as GL$(k,\R)$-bundle.
By the reduction from GL$(k,\R)$ to $O(k)$, then $f^*\tau_{n,k}$
is associated to its {\it orthogonal frame bundle} 
$f^*s\tau_{n,k}$, both considered as $O(k)$-bundles.
In particular by applying this to the tangent bundle $T(M)$ of a manifolds,
we get the linear or othogonal {\it frame bundle} of $M$, say $F_l(M)$ or $F_s(M)$. 
$M$ is parallelizable if and only if $F_l(M)$ (hence $F_s(M)$) has a
section.

\section{Limit tautological bundles}\label{LIMITS}
 We will deal with a few concrete instances of the following general
 topological construction. Let $\{X_n\}_{n\in \N}$ be a countable
 family of Hausdorff topological spaces each admitting a countable
 basis of open sets. Assume that for every $n$, $X_n$ is strictly
 contained in $X_{n+1}$ as a closed subset. Then consider the 
 ``limit'' space
 $$ X_\infty = \cup_n X_n$$
 endowed with the {\it final topology} with respect to the family of
 inclusions
 $$\{i_n : X_n\to X_\infty\} \ ;$$
this means the {\it finest} topology such that every $i_n$
is continuous. In other words, $A$ is open in $X_\infty$ if
and only if for every $n$,
$A\cap X_n$ is open in $X_n$. We have

\begin{lemma}\label{final-compact}
  If $K\subset X_\infty$ is compact then there is
  $n\in \N$ such that $K\subset X_n$.
\end{lemma}
\Dim Assume that there is not, then there should be an infinite  sequence $x_n$
in $K$ such that $x_n\in X_{n+1}\setminus  X_n$.  The union of these
points of $K$ would be a closed subset of $K$ (hence compact) with
induced discrete topology (i.e. it would be a {\it compact and discrete}
space). Such a space is necessarily {\it finite} against
  our assumption.

  \cvd

  \medskip
{\it Some examples:}
\medskip

  $\bullet$ $\R^n\subset \R^{n+1}$, $(x)\to(x,0)$. Then we can define the limit space $\R^\infty$.

  \medskip

  The above inclusions induce ``equatorial" inclusions $i_n:S^{n-1}\to S^n$ of unit spheres, 
  so we can define the limit space $S^\infty$.

  \medskip

  $\bullet$ The definition of $S^\infty$ can be generalized to arbitrary Stiefel
  manifolds. The inclusions $M(n,k,\R)\to M(n+1,k,\R)$
  $$
  A \to \begin{pmatrix}
A\\
0
\end{pmatrix}
$$ 
induce inclusions of embedded smooth manifolds $i_n: S_{n,k}\to S_{n+1,k}$, and we
can define the {\it Stiefel limit space} $S_{\infty,k}$.

$\bullet$ The inclusions $S(n,\R)\to S(n+1,\R)$

  $$ A\to  \begin{pmatrix}
A&0\\
0&0
\end{pmatrix}
$$ 
induce the inclusions $j_n:=j_{n,n+1}: \GG_{n,k}\to \GG_{n+1,k}$, and we
can define the {\it limit grassmannian} $\GG_{\infty,k}$.

  \medskip

$\bullet$ Clearly we have the family of commutative diagramms of smooth maps 

$$ \begin{array}[c]{ccc}
S_{n,k}&\stackrel{i_n}{\rightarrow}& S_{n+1,k}\\
\downarrow\scriptstyle{s_{n,k}}&&\downarrow\scriptstyle{s_{n+1,k}}\\
\GG_{n,k}&\stackrel{j_n}{\rightarrow}& \GG_{n+1,k}\end{array}$$
so we can eventually define the ``limit projection'' which is continuous
  
$$ \begin{array}[c]{ccc}
S_{\infty,k}\\
\downarrow\scriptstyle{s_{\infty,k}}\\
\GG_{\infty,k}\end{array}$$

Symilarly by using the linear frames we have the limit projection
$$ \begin{array}[c]{ccc}
L_{\infty,k}\\
\downarrow\scriptstyle{l_{\infty,k}}\\
\GG_{\infty,k}\end{array}$$

\begin{example}\label{P-infinito}
  {\rm
    As a particular case we have the projection
    $$s_{\infty,1}: S^\infty \to \PP^\infty(\R) \ . $$
    We easily realizes
    that $s_\infty$ is a continuous covering map of degree $2$, alike
    every $s_{n,1}$. Thanks to lemma \ref{final-compact}, for every
    $p\in \N$, every continuous map $f: S^p \to S^\infty$ is of the
    form $i_n\circ \tilde f$, for some $\tilde f: S^p \to S^n$ such
    that the image of $\tilde f$ does not contain $e_{n+1}$. By
    considering $S^n=\R^n\cup \{\infty\}$ via the stereographic
    projection with center $\infty:= e_{n+1}$, then $\tilde f$
    factorizes through a map with values in $\R^n$ which is
    contractible. We can conclude that every such a map $f$ is {\it
      homotopically trivial}. In other words all homotopy groups
    $\pi_p(S^\infty)$ are trivial. By a theorem of Whitehead (see
    \cite{H}), it follows that $S^\infty$ is contractible, hence
    $s_\infty: S^\infty \to \PP^\infty(\R)$ is a {\it universal
      covering map}. By the theory of covering maps we eventually get that
    the fundamental group $\pi_1(\PP^\infty(\R))\sim \Z/2\Z$, while
    all other groups $\pi_p(\PP^\infty(\R))$, $p>1$, are trivial. We
    summarize these facts by saying that $\PP^\infty(\R)$ is a {\it
      $K(\Z/2\Z,1)$ spaces.}} 
\end{example}

$\bullet$ The same limit procedure applies to the tautological bundles.
We have the family of commutative diagramms of smooth maps 
$$ \begin{array}[c]{ccc}
\Vv(\GG_{n,k})&\stackrel{\tilde j_{n}}{\rightarrow}& \Vv(\GG_{n+1,k})\\
\downarrow\scriptstyle{\tau_{n,k}}&&\downarrow\scriptstyle{\tau_{n+1,k}}\\
\GG_{n,k}&\stackrel{j_n}{\rightarrow}& \GG_{n+1,k}\end{array}$$
so we eventually define the ``limit tautological vector bundle'' :
  
$$ \begin{array}[c]{ccc}
\Vv(\GG_{\infty,k})\\
\downarrow\scriptstyle{\tau_{\infty,k}}\\
\GG_{\infty,k}\end{array} \ . $$
Similarly we have the limit bundles

$$ \begin{array}[c]{ccc}
\Ll(\GG_{\infty,k})&\stackrel{ }{\ \ }& \Ss(\GG_{\infty,k})\\
\downarrow\scriptstyle{l\tau_{\infty,k}}&&\downarrow\scriptstyle{s\tau_{\infty,k}}\\
\GG_{\infty,k}&\stackrel{ }{ \ \ }& \GG_{\infty,k}\end{array}$$

\section{A classification theorem for compact manifolds}\label{class-comp}
In this section {\it we assume that $M$ is compact}.
By  Lemma \ref{final-compact},  $f\in \Cc^0(M,\GG_{\infty,k})$
if and only if there is a minimum $n$ such that it factorizes through
a continuous map $$\hat f: M \to \GG_{n,k}$$ 
followed by the inclusion
$$j_{n,\infty}: \GG_{n,k}\to \GG_{\infty,k} \ . $$ 
So it makes sense to say that such a map $f$ is {\it smooth}
if $\hat f$ is smooth in the usual sense. Moreover also the topologies on the spaces
$\Ee(M,\GG_{n,k})$ pass to the limits, giving us the topological space $\Ee(M,\GG_{\infty,k})$
of such smooth maps. 
If $f, \hat f$ are as before, we have
$$ f^*\tau_{\infty,k} = \hat f^*\tau_{n,k}$$
provided that the we have incorporated the  {\it canonical identifications}
illustrated in section \ref{B-equivalence}.
Set 
$$ \Vv_k(M):= \{f^*\tau_{\infty,k}; \ f\in \Ee(M,\GG_{\infty,k}) \} \ . $$
It is clear from the above considerations that the already defined space
$\Vv(M)$ can be described as
$$ \Vv(M)= \cup_{k=0}^\infty \ \Vv_k(M) \  $$
as well as
$$ \Ee(M,\GG)= \cup_{k=0}^\infty \ \Ee(M, \GG_{\infty,k}) \ . $$
Thus we have rephrased in terms of these limits the surjective maps
$$ (.)^*: \Ee(M,\GG) \to \Vv(M)$$
$$ [(.)^*] : \Ee(M,\GG)\to \VV_0(M) \  $$
and we stipulate that the target spaces are endowed with the
{\it quotient topology}. 

Given $f_0,f_1 \in \Ee(M,\GG)$ we say that they
are {\it smoothly homotopic} if $f_0,f_1 \in \Ee(M,\GG_{\infty,k})$ for some $k$,
and are connected by a smooth homotopy $F\in \Ee(M\times [0,1],\GG_{\infty,k})$,
provided that  $f_t: =F_{| M\times \{t\}}$. 
As usual, this defines an equivalence relation on $\Ee(M,\GG)$.
Denote by $ [M,\GG]$
the set of smoothly homotopy classes of maps of $\Ee(M,\GG)$.

\begin{proposition}\label{fromV-to-homotopy}  Let $M$ be an embedded compact smooth manifold. If
$[f_0^*\tau_{\infty,k}]=[f^*_1\tau_{\infty,k}]$
in $\VV_0(M)$, then $f_0$ and $f_1$ are homotopic. Hence it is well
defined a surjective map
$$ \vG: \VV_0(M)\to [M,\GG], \ [f^*\tau_{\infty,k}]\to [f] \ . $$
\end{proposition}
\Dim We will provide two proofs. 

{\it First proof:}  If  $[f_0^*\tau_{\infty,k}] = [f^*_1\tau_{\infty,k}] \in \VV_0(M)$, 
we can assume that they both factorize through maps
(for simplicity we keep the same names) $f_0,f_1: M \to \GG_{n,k}$, for some $n$ big enough.
Moreover, sometimes we will confuse here a point $A\in \GG_{n,k}$ with the corresponding
subspace $V_A \subset \R^n$.
For $j=0,1$, for every $p\in M$, we have the direct sum decomposition $\R^n= f_j(p)\oplus f_j(p)^\perp$.
The projections of the canonical basis $\{e_1,\dots, e_n\}$ onto $f_j(p)$, when $p$ varies, define
$n$-sections $s_{j,1},\dots, s_{j,n}$ of  $f_j^*\tau_{n,k}$ which span the fibre $f_j(p)$ over every $p\in M$.
The map $f_j$ can be reconstructed from these set of sections as follows: for every $p\in M$, the linear evaluation map 
$$\eG_{j,p}:  \R^n\to f_j(p), \  \eG_{j,p}(X)= \sum_i x_is_{j,i}(p)$$ 
is onto so that $\ker (\eG_{j,p})= f_j(p)^\perp$ and finally $f_j(p)=  \ker (\eG_{j,p})^\perp$.
A strict equivalence from $f^*_0\tau_{n,k}$  to $f^*_1\tau_{n,k}$
transports  the system of sections $s_{0,1},\dots, s_{0,n}$  to a system  $s'_{1,1},\dots, s'_{1,n}$ over
$f^*_1\tau_{n,k}$ which generate all its fibres. Denote by $\eG'_{1,p}$ the corresponding evaluation maps
and apply to it the above procedure in order to produce a map from $M$ with value in $\GG_{n,k}$; 
we realize that this recovers $f_0$. For every $p\in M$, $\ker (\eG'_{1,p})$ is a graph of a linear
map $L_p: f_1(p)^\perp \to \ f_1(p)$, while $f_1(p)^\perp$ itself is the graph of the zero map. The homotopy
$L_{p,t}=tL_p, \ t\in [0,1]$, eventually allows to define a desired homotopy between $f_0$ and $f_1$. 
\smallskip

{\it Second proof:}  We know that $f_j$ is determined by a map say $g_j$ from $f_j^*(\Vv(\GG_{\infty, k}))$
to $\R^\infty$ which is linear and injective at every fibre. Moreover, it factorizes through a map
with value in some $\R^k$ with $k$ big enough.  If  $[f_0^*\tau_{\infty,k}] = [f^*_1\tau_{\infty,k}] \in \VV_0(M)$,
we can transport the map $g_1$ to a map $g'_0$ with such a property, defined on  $f_0^*(\Vv(\GG_{\infty, k}))$
and we have to show that $g_0$ and $g'_0$ are homotopic through maps  
that are linear injections on fibers. First compose $g_0$ with the homotopy $a_t :\R^\infty\to \R^\infty$  
defined by $a_t(x_1,x_2, \dots) = (1- t)(x_1,x_2, \dots)+t(x_1,0,x_2,0, \dots)$. This moves the image of $g_0$ into the odd-numbered coordinates. 
Similarly we can move $g'_0$ into the even-numbered coordinates. By keeping the names of these maps, we eventually define
the desired homotopy  $h_t =(1- t)g_0+tg'_0$.
 
\cvd

\smallskip

Finally we can answer Question \ref{V0-question}, at least in the compact case.
A similar classification theorem holds under more general assumptions.
Compactness simplifies the proof and it will suffice to the aim of this text.

\begin{theorem}\label{VBClassification} {\rm (Classification Theorem)} 
  Let $M$ be an embedded {\rm compact} smooth manifold.
  Then the map
  $$ \vG: \VV_0(M)\to [M,\GG], \ [f^*\tau_{\infty,k}]\to [f]$$
  is bijective. That is for every $f_0,f_1\in \Ee(M,\GG)$,
$[f_0^*\tau_{\infty,k}] = [f^*_1\tau_{\infty,k}] \in \VV_0(M)$ if and only if
$f_0$, $f_1$ are smoothly homotopic. Hence
the map $[(.)^*]$ induces the inverse map of $\vG$ 
$$ \cG: [M,\GG]\to \VV_0(M), \ \cG([f])= [f^*\tau_{\infty,k}]  \ {\rm whenever}\  f\in \Ee(M,\GG_{\infty,k}) \ . $$

\end{theorem}

\Dim Thanks to Proposition \ref {fromV-to-homotopy}, it is enough to prove that
if $f_0$ and $f_1$ are homotopic,
then $f_0^*\tau_{\infty,k}$ and $f_1^*\tau_{\infty,k}$ are strictly equivalent.
We can assume that 
a homotopy factorizes through $F:M\times [0,1] \to \GG_{n,k}$, $n$ big enough. Take the pull-back $F^*\tau_{n,k}$.
The idea is to use it in order to connect $f^*_0\tau_{n,k}$ and $f^*_1\tau_{n,k}$ by a path $f^*_t\tau_{n,k}$
of bundles strictly equivalent to each other.
For every $t\in [0,1]$, $p\in M$, denote by $V_{t,p}$ the fibre of $f^*_t \tau_{n,k}$ over $p$.

{\bf Claim 1.} {\it There is $\epsilon >0$ such that for every $t\leq \epsilon$, $f_0^*\tau_{n,k}$
 is strictly equivalent to $f_t^*\tau_{n,k}$.} 
 
 To prove it, recall the elementary fact that if 
 $\R^n=V'\oplus V=V"\oplus V$ ($V$, $V'$ and $V"$ being linear subspaces), 
 then $\phi: V'\to V"$, $\phi(v')= v"$ if $v'=v"+v$, is a {\it canonical}
 linear isomorphism between $V'$ and $V"$. We have: 
 \smallskip
 
{\bf Claim 2.} {\it There is $\epsilon>0$ such that for every $0\leq t \leq \epsilon$,
for every $p\in M$,  $\R^n=V_{0,p}\oplus V_{0,p}^\perp= V_{t,p} \oplus V_{0,p}^\perp$.} 
\smallskip

\noindent Assuming {\bf Claim 2}, then, for every $t\leq \epsilon$, the ``field'' of canonical
 isomorphisms $$\phi_p: V_{t,p} \to V_{0,p}$$ 
 when $p$ varies in $M$, defines a 
 strict equivalence, as required by {\bf Claim 1}. Let us prove {\bf Claim 2.} If such an $\epsilon$ does not exist
 by compactness there would exist a converging sequence $(p_n,t_n) \to (p_0,0)$ in $M\times [0,1]$,  such that 
 for every
  $n$, $\dim V_{t_n, p_n}\cap V_{0,p_n}^\perp >0$. But this is impossible because
 $V_{0,p_0}\cap V_{0,p_0}^\perp = \{0\}$ and this is an open condition.  
 
 Set $\epsilon_0\in [0,1]$ the $\sup$ of the $\epsilon$'s verifying {\bf Caim 1}.
 We claim furthermore that $\epsilon_0$ is a {\it maximum}. In fact by applying the same
 argument, we see that there is $\epsilon>0$ such that $f^*_{\epsilon_0}\tau_{n,k}$ is strictly
 equivalent to $f^*_t\tau_{n,k}$, for $t\in (\epsilon_0-\epsilon, \epsilon_0]$. Finally we claim
 that $\epsilon_0=1$: if $\epsilon_0 <1$, we can apply again the above argument to $f_{\epsilon_0}$ and find
 $\epsilon_1=\epsilon_0+\epsilon$, for some small $\epsilon >0$, which works as well,
 against the fact ther $\epsilon_0$ is the maximum.
 
\cvd

\medskip

$\bullet$ {\it The above discussion can be repeated word by word, getting similar conclusions,
by dealing with {\it embeded frame bundles} and using the limit tautological bundles $\l\tau_{\infty,k}$
or $s\tau_{\infty,k}$, $\Ll(M)$ or $\Ss(M)$. }

\medskip

\section{The rings of stable equivalence classes of vector bundles}\label{K-theory}
The final aim of this section is to endow  a  suitable quotient space $\KK_0(M)$ of $\VV_0(M)$
with a natural {\it ring} structure, for every embedded smooth manifold $M$. This leads
to a contravariant functor from the category of embedded smooth manifolds to the category of
commutative rings. If $M$ is compact we point out more information such as the invariance up to homotopy
of the functor. 

\subsection{Grassmannian operations}\label{Grass-oper} 
The operations of the ring $\KK_0(M)$ will descend from simple `operations' defined between Grassmann
manifolds. 
\smallskip

$\bullet$ The inclusion $S(n,\R) \to S(n+m,\R)$
$$ A\to  \begin{pmatrix}
A&0\\
0&0
\end{pmatrix}
$$

induces for every $k\leq n$, a smooth inclusion
$$j_{n,n+m}: \GG_{n,k}\to \GG_{n+m,k} \ . $$

$\bullet$ The inclusion $S(n,\R)\times S(m,\R)\to S(n+m,\R)$

 $$ (A,B)\to  \begin{pmatrix}
A&0\\
0&B
\end{pmatrix}
$$

induces for every $k\leq n$, $h\leq m$ a smooth  inclusion
$$ \oplus_{n,k,m,h}: \GG_{n,k}\times \GG_{m,h}\to \GG_{n+m,k+h} \ .$$

\medskip

$\bullet$ For every $V\in G_{n,k}$ denote by $V^*$ its dual
spaces. Recall that this is considered as a subspace of $(\R^n)^*=M(n,1,\R)$ as
follows. Let $\R^n = V\oplus V^\perp$ the othogonal direct sum
decomposition, $V^\perp\in G_{n,n-k}$ being the orthogonal complement
of $V$ with respect to the standard euclidean scalar product. Then
extend every $\gamma \in V^*$ to a functional defined on the whole of
$\R^n$ by setting $\gamma(u+w)=\gamma(u)$.  $M(n,1,\R)$ is canonically
isomorphic to $\R^n$ via the transposition.

Let $(V,W)\in G_{n,k}\times G_{m,h}$.  Denote by $V\otimes W$ the
space of bilinear forms defined on $V^*\times W^*$. Its dimension is
$kh$. In fact there is the canonical bilinear map
$$\otimes: V\times W \to V\otimes W, \ v\otimes w(\gamma,\rho):=
\gamma(v)\rho(w) $$
and for every couple of bases $(\Bb, \Dd)$ of $V$
and $W$ respectively, then $\Bb\otimes \Dd =\{ v_i\otimes w_j;
\ v_i\in \Bb, \ w_j\in \Dd\}$ is a basis of $V\otimes W$.  By using
the decomposition
$$ \R^n \times \R^m = (V\oplus V^\perp)\times (W\oplus W^\perp)$$
and arguing as above we can consider $V\otimes W$ as a subspace of
$\R^n \otimes \R^m$,
hence (via canonical isomorphisms) as an element of $G_{nm,kh}$.
In this way we have defined a map (between {\it sets}):
$$ G_{n,k}\times G_{m,h} \to G_{nm,kh} \ . $$
This can be transported to a map
$$\otimes_{n,k,m,h}:  \GG_{n,k}\times \GG_{m,h}\to \GG_{nm,kh}$$
via the usual bijections $V\to A_V$, $\dots$.
One can check by direct computation that this is a {\it smooth map}
between embedded smooth manifolds.

\medskip

Similarly one can check that the set map
$$G_{n,k}\to G_{n,n-k}, \ V\to V^\perp $$
induces a {\it diffeomorphism}
$$\perp_{n,k}: \GG_{n,k}\to \GG_{n,n-k}$$
with inverse $\perp_{n,n-k}$. 

\subsection {The ring $\KK_0(M)$} \label{K-ring} 
The grassmannian operations of Section \ref{Grass-oper} induce operations
$$ \oplus: \Vv(M)\times \Vv(M) \to \Vv(M),  \ f^*\tau_{n,k}\oplus g^*\tau_{r,s}= (\oplus \circ (f,g))^* \tau_{n+r,k+s} \  $$
$$ \otimes: \Vv(M)\times \Vv(M)\to \Vv(M),  \  f^*\tau_{n,k}\otimes g^*\tau_{r,s}= (\otimes \circ (f,g))^* \tau_{nr,ks} \  $$
$$ \perp: \Vv(M) \to \Vv(M), \ \perp(f^*\tau_{n,k})= (\perp \circ f)^*\tau_{n,n-k} \ . $$

The operations $\oplus, \otimes,\perp$ descend to each quotient set
$\Vv(M)/{\rm Aut(M)}$, $\VV(M)$ and $\VV_0(M)$.
\medskip

The grassmannian operations $\oplus$ and $\otimes$ pass to the limits:
$$ \oplus : \GG_{\infty,k}\times \GG_{\infty, h} \to \GG_{\infty, k+h}$$
$$ \otimes :  \GG_{\infty,k}\times \GG_{\infty, h} \to \GG_{\infty, kh} \ $$
and are continuous in the limit topology.
The operation $\perp$ induces in fact a family
of continuous maps
$$\perp_n : \GG_{\infty, k} \to \GG_{\infty, n-k}, \ n\geq k  \ . $$

\medskip

For every embedded smooth manifold $M$, these operations define a {\it
  ring} structure on a suitable quotient of $\VV_0(M)$ that we are
going to point out. Denote by $\epsilon^k$ the class in $\VV_0(M)$ of
the trivial (product) bundle $M\times \R^k\to M$. Clearly
$$\epsilon^k\oplus\epsilon^h=\epsilon^{k+h}\ . $$

\begin{definition}\label{w-stable-equiv} {\rm 
We say that $\xi$ and $\eta$ in  $\VV_0(M)$ are
{\it weakly stably equivalent} if there exist $\epsilon^k$ and $\epsilon^h$ such that
$$ \xi \oplus \epsilon^k= \eta\oplus \epsilon^h \ . $$}
\end{definition}

This is an equivalence relation indeed. Let us just check the transitivity. If
$$ \xi \oplus \epsilon^k=\eta \oplus \epsilon^h, \ \eta \oplus \epsilon^r=\beta \oplus \epsilon^s$$
then
$$ \xi \oplus \epsilon^{k+r}=\beta \oplus \epsilon^{h+s} \ . $$

\begin{example}{\rm (1) Let $M$ be a smooth manifold with non empty boundary $\partial M$.
Let $i:\partial M \to M$ the inclusion. Then $T(\partial M)$ and
$i^*T(M)$ are weakly stably equivalent vector bundles on $\partial
M$. Fix any riemannian metric $g$ on $M$. For every $x\in \partial M$,
consider $\nu(x)=(T_x\partial M)^{\perp_{g(x)}}$; as
$T_xM=\nu(x)\oplus T_x\partial M$, this defines a vector bundle $\nu$
on $\partial M$, with $1$-dimensional fibres, such that $i^*T(M)=\nu
\oplus T(\partial M)$.  The bundle $\nu$ has a nowhere vanishing
section (for every $x\in \partial M$ take the ``outgoing'' $g$-unitary
vector in $\nu(x)$). Then $[\nu]=\epsilon^1$. In particular
$S^n=\partial B^{n+1}(0,1)$, $T(B^{n+1})$ is trivial as it is the
restriction of $T(\R^{n+1})$, hence $[T(S^n)]$ is weakly stably
trivial.}
\end{example}

Denote by $\KK_0(M)$ the quotient of $\VV_0(M)$ up to weakly stable
equivalence. It is clear that if $M=\{p\}$ is one point, then
$\KK_0(\{p\})=0$.
\medskip

\begin{proposition} The operations $\oplus$, $\otimes$ descend to $\KK_0(M)$ and make it an abelian ring.
\end{proposition}
\Dim Associativity of $\oplus$ is evident. The weakly stable
equivalence class $[\epsilon^1]$ is the zero element; for every
$[[\xi]]$, assume that $\xi \in \Vv_{n,k}(M)$, then $\xi^\perp \in
\Vv_{n,n-k}(M)$ is such that
$$[\xi\oplus \xi^\perp]=\epsilon^n$$
hence
$$[[\xi^\perp]]=-[[\xi]] \ . $$
With a bit of more work one can also check the ring structure. We leave it as an exercise.

\cvd

Summing up
\medskip

$$M \  \Rightarrow \ \KK_0(M) $$
$$g:N\to M \ \Rightarrow  \ g^{\bullet}: \KK_0(M)\to \KK_0(N),
\ g^\bullet([[f^*\tau_{\infty,k}]])= [[(f\circ g)^*\tau_{\infty,k}]]$$
define a {\it contravariant functor from the category of embedded smooth manifolds (with boundary) to the category
of abelian rings}.
\medskip

If $M$ is {\it compact}, the above construction of the ring $\KK_0(M)$ from $\VV_0(M)$ can be 
rephrased in terms of $[M,\GG]$.
So:  $[f_0]$, $f_0:M \to \GG_{\infty,s}$,
and $[f_1]$, $f_1: M \to \GG_{\infty,r}$, are weakly stably equivalent if and only if there are
constant maps $c_0: M \to \GG_{\infty,k}$, $c_1:M\to \GG_{\infty,h}$, such that
$[\oplus \circ (f_0,c_0)]=[\oplus \circ (f_1,c_1)]$ in $[M,\GG]$. Denote by $[[M,\GG]]_0$
the quotient set. We have

\begin{proposition}\label{CT-iso}  Let $M$ be compact. The operations $\oplus$ and $\otimes$ descend to $[[M,\GG]]_0$
and make it an abelian ring such that the map $\vG$ induces a  ring isomorphism
$$\tilde \vG: \KK_0(M)\to [[M,\GG]]_0 $$
with inverse
$$\tilde \cG: [[M,\GG]]_0\to \KK_0(M)$$
induced by the map $\cG$ of the Classification Theorem \ref {VBClassification}.
\end{proposition}

\cvd

\medskip

\begin{corollary}\label{Homot-invariance} {\rm (Homotopy invariance)} Let $M$, $N$ be compact smooth manifolds.
Then:

  (1) If $g_1, g_2 \in \Ee(N,M)$ are smoothly homotopic, then $g_1^\bullet=g_2^\bullet$. 

(2) If $M$ and $N$ are  smoothly homotopically equivalent, then $\KK_0(M)$ and $\KK_0(N)$ are isomorphic.
In particular if $M$ is smoothly contractible, then $\KK_0(M)\sim \KK_0(\{p\})=0$.

\end{corollary}
\Dim (1) and (2) follows from the Classification Theorem, as $[[*,\GG]]_0$ is manifestly homotopically invariant.

\cvd
 
\smallskip

We conclude this Section with a few scattered remarks.
\bigskip

\begin{remarks}\label{K-remark}{\rm 
    
    (1)  $\KK_0(*)$ is a versions in our embedded smooth framework of so called {\it reduced topological
$K$-theory} \cite{A} \cite{B}. Taking into account, for simplicity, only the additive structure, the {\it unreduced} group say $\KK(M)$
is constructed as follows. First we consider the quotient say $\tilde \VV_0(M)$ of $\VV_0(M)$ up to {\it stable equivalence};
 this is defined similarly to the above {\it weak stable equivalence} by imposing in the definition that $k=h$. 
 The operation $\oplus$ passes to the quotient, so that $(\tilde \VV_0(M),\oplus)$ is a {\it commutative monoid} with (the class
 of) $\epsilon^0$ as zero element.
 \medskip

 {\it  $(\tilde \VV_0(M),\oplus)$ verifies the ``cancellation rule''.}

 \medskip

 In fact, if $\xi\oplus \eta = \xi \oplus \alpha$, we know that there exists $\beta$ such that
 $\xi\oplus \beta=[\epsilon^n]$ (for some $n$), hence $[\epsilon^n]\oplus \eta = [\epsilon^n]\oplus \alpha$
 and finally $\eta=\alpha$.

 Then $\KK(M)$ is the {\it Grothendieck group} of this
 monoid with cancellation rule.
 It is a general construction (producing for instance $(\Z,+)$ from $(\N,+)$)
 that works as follows. Consider the product $\tilde \VV_0(M)\times \tilde \VV_0(M)$;
 often an element $(\xi,\eta)$ is written as a formal difference $\xi -\eta$.  Put on this product the equivalence
 relation such that
 $$\xi-\eta \sim \alpha - \beta$$
 if and only if
 $$\xi\oplus \beta = \alpha \oplus \eta \ .$$
 The cancellation rule is used to
 check that it is actually an equivalence relation. The addition rule on the quotient $\KK(M)$ naturally is
 $$(\xi-\eta)\oplus (\alpha - \beta)= \xi\oplus \alpha - \eta \oplus \beta \ ; $$
 the zero element is given by
 $$[\epsilon^0] - [\epsilon^0] = \xi - \xi, \ \forall \xi \in \tilde \VV_0(M) \ ; $$
The inverse of $\xi-\eta$ is $\eta-\xi$.
\medskip

{\it Every element of $\KK(M)$ can be represented by  a difference of the form $\xi - [\epsilon^n]$ (for some $n$).}

\medskip

In fact, for every $\alpha - \beta$, let $\beta\oplus \gamma = [\epsilon^n]$, then
$$\alpha - \beta = \alpha\oplus \gamma -\beta\oplus \gamma := \xi - [\epsilon^n] \ . $$

The correspondence $\xi-[\epsilon^n] \to \xi$ induces a canonical surjective homorphism $\KK(M)\to \KK_0(M)$.
It is well defined because if $\xi - [\epsilon^n] = \xi'- [\epsilon^m]$ in $\KK(M)$, then
$\xi\oplus [\epsilon^m]=\xi'\oplus [\epsilon^n]$, hence $\xi=\xi'$ in $\KK_0(M)$.
The kernel consists of the elements of the form $[\epsilon^n]-[\epsilon^m]$ which is 
isomorphic to $\Z$ so that $\KK(M)\sim \KK_0(M)\oplus \Z$ (in a non canonical way).
\smallskip
 
(2) If $M$ is compact, the construction of $\KK(M)$ from $\VV_0(M)$ can be rephrased in terms of $[M,\GG]$.
This produces a group (a ring indeed) $[[M,\GG]]$ which is isomorphic to $\KK(M)$,
via the Classification Theorem (similarly to Proposition \ref{CT-iso}).
Hence also the functor
$$M \ \Rightarrow \ \KK(M)$$
$$ \dots \ \Rightarrow \dots $$
verifies the {\it homotopy invariance} properties, similarly to Corollary \ref{Homot-invariance}.

\smallskip
  
(3) We can develop the very same constructions  by using the complex grassmannians
$\GG_{n,k}(\C)$ and the complex vector bundles; this leads to the functors 
\bigskip

$$M \ \Rightarrow \ \KK_0(M,\C), \ \KK(M,\C)$$
$$ \dots \ \Rightarrow \dots  \ . $$

\smallskip

(4) {\it Bott's periodicity theorem} \cite{B}, \cite{At} is among the fundamental results in this theory. Let us just recall a few 
related statements that we can formulate in our setting.

\begin{itemize}
\item  For every compact $M$, $\KK(M\times S^2, \C)\sim \KK(S^2, \C)\otimes \KK(M,\C)$;
\item  $\KK(S^2,\C)=\Z[X]/(X-1)^2$ where $X$ is the tautological complex bundle over $\PP^1(\C)$ (recall that $\PP^1(\C)$
is diffeomorphic to $S^2$, the ``Riemann sphere");
\item For every $m\geq 1$, $\KK_0(S^{m+8})\sim\KK_0(S^m)$, $\KK_0(S^{m+2},\C)\sim\KK_0(S^m, \C)$.
\end{itemize}
\smallskip

(5) {\bf On real algebraic vector bundles.}
We have shown that every Grassmann manifold $\GG_{n,k}$ is also a
regular real algebraic set.  Dealing with real algebraic sets, say
$X\subset \R^n$, $Y\subset \R^m$, a natural class of maps $\Rr(X,Y)$
consists of so called {\it regular rational maps} (shortly ``algebraic'')
that is restriction of rational maps $r: \R^n\to \R^m$, whose
denominators nowhere vanish on $X$.
Consider the tautological vector bundle
$$\tau_{n,k}: \Vv(\GG_{n,k})\to \GG_{n,k} \ . $$
It is immediate that also the total space $\Vv(\GG_{n,k})$
is a regular algebraic set, and that $\tau_{n,k}$ is
algebraic. Moreover, if $M$ is any regular real  algebraic sets, and
$$f:M \to \GG_{n,k}$$
is an algebraic map, then one readly checks that also the pull-back
$f^*\tau_{n,k}$ verifies the same properties. So we can consider
the family of {\it algebraic vector bundles} on $M$
$$\Vv^{{\rm alg}}(M) = \cup_{(n,k)\in \Nn} \ \Vv^{{\rm alg}}_{n,k}(M)$$
where we consider only the pull-back via algebraic maps.
The operations $\oplus$, $\otimes$, $\perp$ restrict algebraically.
We can also consider $\VV_0^{{\rm alg}}(M)$ where we impose that the strict
equivalence are realized by algebraic map. Many constructions developed so far
have a natural ``algebraic'' specialization (for instance we have $\KK_0^{{\rm alg}}(M)$,
$\KK^{{\rm alg}}(M)$). By forgetting the algebraic structure and keeping only the
one of smooth manifold, we have natural forgetting maps
$$ \Vv^{{\rm alg}}(M)\to \Vv(M), \ \VV_0^{{\rm alg}}(M)\to \VV_0(M), \dots, \ \KK^{{\rm alg}}(M)\to \KK(M) $$
and natural interesting questions (injective, surjective, $\dots$, ?)
whose answers presumably depend on the real algebraic
structure. On another hand, it is not so evident how to formulate an algebraic version of the
Classification Theorem (for example our proof that smooth homotopy defines an equivalence
relation used the bump function, and this is not very ``algebraic'' indeed).

Similar algebraic specialization holds also for the frame tautological bundles.
}
\end{remarks}

\chapter{Compact embedded smooth manifolds}\label{TD-COMP-EMB}
The hypothesis that an embedded smooth manifold $M$ is {\it compact}
usually simplifies the study of several objects associated to it.
A first example has been the  proof of the Classification Theorem
of embedded vector bundles in Chapter \ref{TD-EMB-VB}. We will develop
this theme, by considering first a few technical device.

\section{Nice atlas and finite partitions of unity}\label{nice-atlas}
Let $M$ be an embedded smooth $m$-manifold (possibly with boundary).
Recall that a normal chart $(W,\phi)$ of $M$ is either
contained in the interior of $M$ and of the form
$$ \phi: W \to B^m(0,1)$$
or it intersects $\partial M$ and is of the relative form
$$\phi (W,W\cap \partial M)\to
(B^m(0,1)\cap \HH^m, B^m(0,1)\cap \partial \HH^m) \ . $$
The bump function (recall Section \ref{bump})
$$\gamma=\gamma_{1/3,1/2}:B^m(0,1)\to \R$$
lifts to a {\it global bump function}
$$\gamma_W: M \to \R$$
with compact support
$$S_W= \phi^{-1}(\bar B^m(0,1/2) \subset W \ . $$
Denote by
$$B_W = \phi^{-1}(B^m(0,1/3))\subset S_W \ . $$
$B_W$ is a relatively compact open set in $M$.

\begin{definition}\label{nice-atlas}{\rm Let $M$ be a compact
    embedded smooth manifold.

    (1) A {\it nice atlas} of $M$ is a finite atlas
    $\Uu=\{(W_j,\phi_j)\}_{j=1,\dots, s}$ formed by normal charts, such
    that the family $\{B_j\}$ ($B_j:= B_{W_j}$) is a open covering of
    $M$.

  (2) Set $\gamma_j:= \gamma_{W_j}$,
  $$\lambda_j:= \frac{\gamma_j}{\sum_j \gamma_j} $$
  so that
  $$\sum_j \lambda_j=1 \ . $$
  Then $\{\lambda_j\}_{j=1,\dots,s}$ is the (finite) {\it partition of
    unity} subordinate to the nice atlas $\Uu$.}
  \end{definition}

It is clear that every compact $M$ admits nice atlas. In fact
we will use nice atlas {\it adapted} to a determined situation or to the
solution of a determined problem. Note for example that the finite
partitions of unity of $\R^n$ involving a bump function at infinity
used in Section \ref{bump}, are in fact restriction of partitions of
unity subordinate to a nice atlas of $S^n$, provided that 
$$\R^n \subset \R^n\cup \{ \infty \}=S^n $$ 
via a stereographic projection.

\section{Spaces of maps with compact source manifold}
\label{comp-func-space} 
We adopt the notations of Section \ref{weak-top}. The so called {\it weak topology} 
is completely adequate when the source
manifold is compact, as it allows a complete global
control over the whole of $M$. In fact, let $f\in \Ee^r(M,N)$. Let $\Uu$ be
a nice atlas of $M$ such that every $(W_j,\phi_j)$ carries a local
representation $f_j$ of $f$. Consider the neighbourhoods of $f$ of the form
$$\Uu_r(f,f_j,\bar B_j,\epsilon) \ . $$
Then every  
$$ \cap_j  \ \Uu_r(f,f_j,\bar B_j,\epsilon)$$
is an open neighbourhood of $f$, and by varying $\epsilon>0$ we get
a {\it basis of neighbourhoods} of $f$ because
$$\cup_j \bar B_j = M \ . $$

Equivalently, in a more ``embedded fashion'':  assume $M\subset \R^h$,
$N\subset \R^k$. Let $\Uu$ be a nice atlas of $M$ such that every $(W_j,\phi_j)$
supports a local smooth extension $g:\Omega_j\to \R^k$ of $f$. 
We can also assume that $\R^h\subset S^h$ as above, and that the 
$\Omega_j$ are part of a nice atlas $\tilde \Uu$ 
of $S^h$ (which restrict to the nice atlas of $M$). By using the partition
of unity subordinate to $\tilde \Uu$ we show that $f$ has a {\it global} smooth extension
$\hat f$ to the whole of $\R^h$. Then, by varying $\epsilon>0$, we 
have a basis of neighbourhoods of $f$ of the
form
$$\Uu_r(f,\hat f,M,\epsilon) \ . $$

Let us study now some remarkable subsets of $\Ee^r(M,N)$, $r\geq 1$ or $\Ee(M,N)$.
\begin{lemma}\label{inj-imm} Let M be compact. Then $f:M\to N$ is an
embedding if and only if it is an injective immersion.
\end{lemma}
\Dim One implication is evident. We know that the other is in general
false without the compactness. To prove it recall that in a compact
(Haussdorf) space a subset is compact if and only if it is closed, and that
a continuous map sends compact sets to compact sets; it follows that
since $M$ is compact, then $f$ is closed so that $f^{-1}$ is continuous
and $f$ is a homeomorphism onto its image in $N$.

\cvd

\medskip

We have

\begin{proposition}\label{open-set}
  Assume that $M$ is compact. Then the subsets of immersions,
  summersions, embeddings, diffeomorphisms are (possibly empty) open
  sets in $\Ee^r(M,N)$, $r\geq 1$ and in $\Ee(M,N)$.
\end{proposition}
\Dim An immersion or summersion $f$ is characterized by the condition
of maximum rank of $d_xf$ at every $x\in M$. If $g$ belongs to a
neighbourhood of $f$ in $\Ee^r(M,N)$, $r\geq 1$ giving a global
control on the whole of $M$ as above (with $\epsilon >0$ small enough)
then $g$ verifies the same maximum rank condition. As for embeddings,
thanks to Lemma \ref{inj-imm} it is enough to prove that if $g$ is
close enough to an injective immersion $f$ then also $g$ is so. Assume
that this thesis fails.  Then there would exist a seguence $g_n \in
\Cc^\infty(M,N)$ , sequences of points $x_n$, $y_n$ in $M$ such that:
\begin{enumerate}
\item Every $g_n$ is an immersion;
\item $g_n \to f$ and  $dg_n \to df$ uniformly
on $M$; 
\item $x_n\to x$, $y_n \to y$ in $M$, $x_n \neq y_n$ and $g_n(x_n)=g_n(y_n)$ for every $n$.
\end{enumerate}
Then: $g_n(x_n)\to f(x)$, $g_n(y_n)\to f(y)$, hence $x=y$ because $f$ is injective.
Then we can localize the situation in a chart of $M$ at $x$ and conclude (getting a contradiction)
by applying the local Proposition \ref{stab}.
Finally if $f$ is a diffeomorphism, in particular it is an embedding, hence $g$ close to $f$
is an embedding. It is enough to prove that $g$ is onto. It is not restrictive to assume that
$M$ is connected, so that also $N$ is connected. As an embedding $g$ is an open map,
its image is open in $N$; on another hand the image of $g$ is compact hence closed because
$M$ is compact. Then  the image of $g$ coincides with the whole of $N$.

\cvd

\section{Tubular neighbourhoods and collars}\label{tubular}
Let $M\subset \R^h$ be a compact boundaryless smooth $m$-manifold. Let $\R^h$ be endowed with
the standard riemannian metric $g_0$. Let us perform the following construction.

(1)  Consider the smooth map
$$\nu: M \to \GG_{h,h-m}$$
where for every $p\in M$, $\nu(p)$ is the (matrix corresponding to the)
orthogonal space $(T_pM)^\perp$  (with respect to $g_0$). 

(2) Take the pull-back 
$$\nu^*\tau_{h,h-m}:  \nu^*(\Vv(\GG_{h,h-m}))\to M  \  . $$
Every fibre $\nu(p)$ of this vector bundle is endowed with the restriction of $g_0$.
We consider $M\subset \nu^*(\Vv(\GG_{h,h-m}))$
via the canonical ``zero section''.

(3) Define the smooth map
$$f_\nu: \nu^*(\Vv(\GG_{h,h-m}))\to \R^h, \ f_\nu(p,v)=p+v \ . $$  
For every  $\epsilon >0$, set
$$ N_\epsilon(M)=\{(p,v)\in  \nu^*(\Vv(\GG_{h,h-m})); \ ||v||_{g_0}\leq \epsilon\} \ . $$ 
It is immediate to verify that
\begin{itemize}
\item $f_\nu(p)=f_\nu(p,0)=p$, for very $p\in M$;
\item there exists $\epsilon>0$ small enough such that
the restriction of $f_\nu$ to $N_\epsilon(M)$ is an immersion.
In fact, 
$ \dim \nu^*(\Vv(\GG_{h,h-m}))=\dim \R^h$, and
for every $x=(p,0)$, the image of $d_xf_\nu$ is equal  to $T_pM\oplus \nu(p)= T_{p}\R^h=\R^h$,
so that $f_\nu$ is an immersion at $M$ and the claim follows by the compactness of $M$.
\end{itemize}

(4) There exists $\epsilon>0$ small enough, such that the restriction (we keep the name)
$$ f_\nu: N_\epsilon(M) \to \R^h$$
is an embedding onto a compact $h$-submanifold of $\R^h$ with boundary,
containing $M$ in its interior. We already know that for $\epsilon >0$ small enough,
$f_\nu$ is an immersion;  it is enough to prove that it is also injective. As it is the identity on $M$, and
$M$ is compact, this follows from the same argument used above to show that the embeddings form 
an open set.

(5) Set 
$$U:= f_\nu(N_\epsilon(M))\subset \R^h$$ 
$$p: U \to M,  \ p:=  \nu^*\tau_{h,h-m}\circ (f_\nu)^{-1} \  . $$

\medskip

Let us analyze the arbitrary or inessential choices made in order to perform
this construction.

\begin{itemize}
\item Certainly $\epsilon$ is not unique.
\item The standard metric $g_0$ has nothing special from a topological differential
view point (we made a similar consideration when we discussed the unitary tangent
bundles). In fact the construction works as well by starting with an {\it arbitrary}
riemannian metric $g$ on $\R^h$.
\item What we have really used of the map $\nu$ is that it defines a 
{\it transverse distribution of $(h-m)$-planes along $M$}, that is for every
$p\in M$, 
$$\R^h=T_pM \oplus \nu(p) \ . $$ However, this is a fake generalization
because it is not hard to prove, by using as usual $\R^h\subset S^h$ and suitable
nice atlas, that for every such a transvese distribution, there is a riemannian metric
$g$ on $\R^h$ that realizes it.
\end{itemize}

Summing up, we can vary the  metric $g$ and the final choice of $\epsilon >0$.  
Let us call {\it tubular neighbourhood of $M$ in $\R^h$} any couple $(U,p)$ obtained
by any implementation of the construction. We have the following {\it uniqueness up to
isotopy} of these tubular neighbourhoods. Fix a auxiliary {\it base} tubular neighbourhood
say $(U^*, p^*)$ constructed by using the standard $g_0$ and some $\epsilon_0$.
We have
\begin{proposition}\label{unique-TN} Let $M\subset \R^h$ be a compact boundaryless $m$-manifold.
Let $(U,p)$  be a tubular neighbourhood of $M$ in $\R^h$.
Then there is a smooth map
$$H: U^*\times [0,1] \to \R^h$$
such that for every $t\in [0,1]$,
\begin{enumerate}
\item $H_t$ is an embedding of $U^*$ onto $U_t\subset \R^h$;
\item $H_t$ is equal to ${\rm id}_M$ on $M$;
\item $(U_t, p_t)$ is a tubular neighbourhood of $M$ in $\R^h$
where  $ p_t := p^*\circ H_t^{-1}$. 

Moreover
\item $H_0={\rm id}_{U^*}$;
\item $(U_1,p_1)=(U,p)$. 
\end{enumerate}

\end{proposition}
\Dim If $(U,p)$ differs from $(U^*,p^*)$ only by $\epsilon \neq \epsilon_0$,
the statement is clearly true (use a radial isotopy fibre by fibre).
Assume that $(U,p)$ has been constructed by using a metric $g$.
Take the path of riemannian metrics $g_t=(1-t)g_0 + tg$, $t\in [0,1]$.
Then there is a ``path" of tubular neighbourhoods $(U_t,p_t)$
constructed by using $g_t$ and some $\epsilon_t>0$. We can also assume that
$\epsilon_t$ is a smooth function of $t$, and that $\epsilon_1=\epsilon$.
Hence we have the family of embeddings
$$ f_{\nu_t}: N_{\epsilon_t}(M,g_t) \to \R^h \ . $$
There is also a family of strict equivalences 
$[{\rm id}_M, \rho_t]$
between $\nu_0^*\tau_{h,h-m}$ and $\nu_t^*\tau_{h,h-m}$
given for every $t\in [0,1]$ by the ``field'' of canonical
linear isomorphisms 
$$\nu_0(p) \to \nu_t(p), \ p\in M$$ 
associated to the two direct sum decompositions
$$\R^h= T_pM\oplus \nu_0(p) = T_pM \oplus \nu_t(p) \ . $$
We can assume (we are free to change $\epsilon_0$) that for every $t$, 
$$\rho_t(N_{\epsilon_0}(M,g_0))\subset N_{\epsilon_t}(M,g_t) $$
and we can define the embeddings
$$ f_{\nu_t}  \circ  \beta_t  \circ (f_{\nu_0})^{-1}: U^* \to U_t  \ . $$
This can be transformed to  $H_t$ with the required properties
by composing it with radial isotopies fibre by fibre.

\cvd

\begin{remark}\label{partial-tube}{\rm The above constructions work as well
    if $M$ is compact with non empty boundary $\partial M$. The resulting
    tubular ``neighbourhoods'' $(U,p)$ are not really neighbourhoods of $M$
    in $\R^h$. Rather they are submanifolds with corners of $\R^h$, containing
    $(M,\partial M)$ as a proper submanifold.}
  \end{remark}

\medskip

\subsection{Tubular neighbourhoods of submanifolds}\label{tub-sub}
Assume now that $Y\subset M \subset \R^h$, $\dim Y = s$, $\dim M= m$, $s<m$, $M$ and $Y$ compact.
Assume also that $M$ and $Y$ are boundaryless. Fix a riemannian metric $g$ on $\R^h$.
As above, we have the associated maps
$$\nu_M: M \to \GG_{h,h-m}$$
$$\nu_Y: Y \to \GG_{h,h-s} \ . $$
Set for every $y\in Y$, 
$$\hat \nu_Y(y):= \nu_Y(y)\cap T_yM \ . $$
This define a smooth map
$$\hat \nu_Y: Y \to \GG_{h,m-s} \ . $$
Define
$$f_{\hat \nu_Y}: \hat \nu^*_Y(\Vv(\GG_{h,m-s}))\to \R^h, \ f_{\hat \nu_Y}(y,v)=y+v \ . $$
Let $(U_M,p_M)$ be a tubular neighbourhood of $M$ constructed by means of $\nu_M$.
There is $\epsilon >0$ small enough such that the image via $f_{\hat \nu_Y}$ of
$$\hat N_{\epsilon}(Y,g)=\{(y,v)\in  \hat \nu^*_Y(\Vv(\GG_{h,m-s})); \ ||v||_g \leq \epsilon\}$$
is contained in $U_M$. Finally define  
$$f_{Y,M}: \hat N_\epsilon(Y,g) \to M, \   f_{Y,M} := p_M \circ  f_{\hat \nu_Y}    \ . $$
Arguing similarly as made above for $f_\nu$,  this verifies
\begin{itemize}
\item $f_{Y,M}(y)=f_{Y,M}(y,0)=y$, for very $y\in Y$;
\item there exists $\epsilon>0$ small enough such that
the restriction of $f_{Y,M}$ to $\hat N_\epsilon(Y,g)$ is an immersion.
\item In fact, there is $\epsilon >0$ small enough such that
the restriction of $f_{Y,M}$ to $\hat N_\epsilon(Y,g)$ is an embedding
onto a neighbourhood $U_{Y,M}$ of $Y$ in $M$.
\end{itemize}

Finally $(U_{Y,M},p_{Y,M})$, where 
$p_{Y,M}=  \hat \nu^*\tau_{h,m-s}\circ (f_{Y,M})^{-1}$
is by definition a {\it tubular neighbourhood of $Y$ in $M$}.
\medskip
  
$\bullet$ By varying $g$ and $\epsilon$,  we have again
the  {\it uniqueness of these tubular neighbourhoods of $Y$ in $M$ up to isotpy}.
We leave the details to the reader. 
 
\subsection{Collars} Consider now $M\subset \R^h$ compact with non empty boundary $\partial M$.
We would apply the above construction, by considering $\partial M$ as a ``monolateral" submanifold
of $M$. By keeping the above notations, we know that  
$$ \hat \nu^*_{\partial M}(\Vv(\GG_{h,1}))$$ is strictly equivalent to the product bundle 
$$\partial M \times \R \to \partial M$$
and a section is given by the unitary ``positive" $v$ (write ``$v>0$"), that is pointing towards the interior of $M$.
So we can define
$$\hat N^+_{\epsilon}(\partial M,g)=\{(y,v)\in  \hat \nu^*_{\partial M}(\Vv(\GG_{h,1})); \ ||v||_g \leq \epsilon, \  ``v\geq 0'' \} \ . $$
By using it, the construction can be repeated and we eventually get (by definition) a {\it collar} of $\partial M$ into $M$,
that is an embedding $C: \partial M \times [0,\epsilon] \to M$ which is the identity on $\partial M$. Again
we have the {\it unicity of  collars up to isotopy.}

\smallskip

\begin{remark}{\rm In the construction of the collars, it is not necessary that the whole $M$ is compact,
    it is enough that $\partial M$ is so.}
  \end{remark}

\begin{remark}\label{proper-tube}{\rm Assume that $Y\subset M \in \R^h$ are compact manifolds with boundary
such that $Y$ is a {\it proper} submanifold of $M$. Then we can apply again the above construction to get
tubular neighbourhoods of $Y$ in $M$ {\it relative} to the boundaries, that is which restrict to tubular neighbourhoods
of $\partial Y$ in $\partial M$.}
\end{remark}
\medskip

Tubular neighbourhoods have several interesting applications. Here is a simple
one. Assume that $M\subset \R^h$ is compact.  We
already know (by using the partitions of unity) that every $f\in
\Ee(M,N)$, $N\subset \R^k$,
extends to a smooth map $\hat f: U\to
\R^k$ defined on a neighbourhood of $M$ in $\R^h$. Let $(U,p)$ be a tubular neighbourhood of $M$.  Then
$f\circ p: U \to N$ is a smooth extension of $f$ {\it with values in
  $N$}.

\section{Proper embedding and ``double'' of manifolds with boundary}\label{partial-proper-emb}
Let $M\subset \R^h$ be a compact smooth manifold with $\partial M \neq \emptyset$.
The existence of collars suggests a variant in the definition of nice atlas.

\begin{definition}\label{atlas-collar} {\rm A {\it nice atlas with collar} of
    $(M,\partial M)$ is of the form
$$ \{(W_\partial,\phi_\partial)\}\cup \{(W_j,\phi_j)\}_{j=1,\dots, s}$$
where
\begin{enumerate}
\item $W_\partial$ is an open neighbourhood of $\partial M$ and
  $$\phi_\partial: W_\partial \to \partial M \times [0,1)$$
  is a diffeomorphism which is equal to the identity on $\partial M$.
  Denote by $B_\partial := \phi_\partial^{-1}([0,1/3))$.
\item Every $(W_j,\phi_j)$ is an normal chart contained in the interior of
  $M$, and $B_j\subset W_j$ is defined as for the usual nice atlas.
\item $\{B_\partial\}\cup\{B_j\}$ is an open covering of $M$.
The existence of nice atlas with collar is a direct consequence of the existence
of collars.
\end{enumerate}  

\smallskip

Given such a nice atlas with collar,
 every $W_j$ carries a global bump function $\gamma_j:M\to \R$ as in Definition
\ref{nice-atlas}. Define the {\it collar global bump function}
$$\gamma_\partial:M\to \R$$
such that on $W_{\partial}$ it is equal to
$\gamma \circ p_{[0,1)} \circ
  \phi_\partial$, where $p_{[0,1)}\partial M \times [0,1) \to [0,1)$
 is the projection, and $\gamma$ is the restriction to $[0,1)$
 of the $1$-dimensional bump function $\gamma_{1/3,1/2}$; on
 $M\setminus W_\partial$, $\gamma_\partial$ is constantly
 equal to $0$. Define
 $$\lambda_\partial = \frac{\gamma_\partial}{\gamma_\partial + \sum_{i=1}^s \gamma_i}$$
 $$\lambda_j = \frac{\gamma_j}{\gamma_\partial + \sum_{i=1}^s \gamma_i}\ . $$
 Then the family of functions
 $$\{\lambda_\partial\}\cup \{\lambda_j\}_{j=1,\dots,s}$$
 is the {\it partition of unity} subordinate to the given nice atlas with collar.}
\end{definition}

\smallskip

\begin{corollary}\label{partial-function} For every compact manifold $M$ with non empty boundary
  there is a smooth function $f:M\to [0,1]$ such that $\partial M = f^{-1}(0)$
  and $f$ is a summersion on a neighbourhood of $\partial M$.
\end{corollary}
\Dim Take a nice atlas with collar. Define locally the following functions
$$f_\partial: W_\partial \to \R, \ f_\partial = p_{[0,1)} \circ \phi_\partial \ ; $$
  $$f_j: W_j \to \R, \ f_j(x)=1/2, \ \forall x\in W_j \ . $$
  Finally set
  $$f= \lambda_\partial f_\partial + \sum_j \lambda_j f_j \ . $$
It is not hard to verify that it is smooth and verifies the required
  properties.

  \cvd

  \smallskip

  The following is an easy generalization

  \begin{corollary}\label{partial-function2} Let $M$ be a compact manifold with
    boundary $\partial M$ equipped with a partition $\partial M = N_0 \cup N_1$,
    where both $N_0$ and $N_1$ are union of connected components of $\partial M$.
    Then there exists a smooth function $f: M \to [0,1]$ such that
    $f^{-1}(0)=N_0$,  $f^{-1}(1)=N_1$,  and $f$ is a summersion on a neighbourhood of $\partial M$.
    \end{corollary}

  \cvd
  
  \begin{remark}\label{minimal-part}{\rm To get the above corollaries we can even use
  a simpler covering of $M$ consisting of  $(W_\partial,\phi_\partial)$ as above together 
  with an open set of the form $U=M\setminus W'$ where $W'\subset W_\partial$ is a smaller compat collar of
  $\partial  M$ contained in $B_\partial$. Hence $W'\subset W" \subset B_\partial$, where $W"$ is another
  collar of $\partial M$, so that the compact sets $B_\partial$ and $B'_\partial:= \overline {M \setminus W"}$ cover $M$.
  By playing with collar bump functions  and variants we get smooth functions $\gamma_\partial$ and $\gamma'_\partial$
  defined on $M$ where $\gamma_\partial$ is as above, while $\gamma'_\partial$ is equal to $1$ on $B'_\partial$
  and is equal to $0$ on $W'$; $\lambda_\partial$, $\lambda'_\partial$ denote the functions of the associated smooth
  partition of unity. Then to prove for instance Corollary \ref{partial-function} define $f_\partial$ as above,  $f_U$
  constantly equal to $1/2$ on $U$ and finally take   $f= \lambda_\partial f_\partial + \lambda'_\partial f_U $.} 
   \end{remark}
   
  \begin{proposition}\label{proper-emb} Let $M\subset \R^h$ be a compact smooth $m$-manifold
    with boundary $\partial M$. Then there is a diffeomorphism
    $\beta: M \to M'\subset \R^n$ (some $n$ big enough) such that $(M',\partial M')$ is a proper submanifold
    of $(\HH^n,\partial \HH^n)$.

 \end{proposition}   
  \Dim Take a nice atlas with collar.
  Define
  $$\beta=(\beta_\partial,\beta_1,\dots, \beta_s): M\to (\R^h\times \R)\times (\R^m\times \R)^s:= \R^n$$
  $$\beta_\partial = (\lambda_\partial \phi_\partial, \lambda_\partial)$$
  $$\beta_j=(\lambda_j\phi_j, \lambda_j) \ . $$
  We claim that this $\beta$ works. To show that it is an embedding it is enough to prove that it is
  an injective immersion. It is an immersion because every $x\in M$ belongs either to $B_\partial$
  or to some $B_j$. The restriction of either $\lambda_\partial \beta_\partial$ or $\lambda_j\beta_j$
  is $\phi_\partial$ or $\phi_j$. In any case it is an injective immersion, so $\beta$ is a fortiori an immersion.
  As for the injectivity, let $x\neq y$. If both belong to either $B_\partial$ or some $B_j$, then they are
  already separated by $\lambda_\partial \beta_\partial$ or $\lambda_j\beta_j$. Otherwise they are separated
  by either $\lambda_\partial$ or some $\lambda_j$. Hence $\beta$ is injective.
  Finally it follows by the construction that the image $M'$ of $\beta$ is contained in $\HH^n$
  and that $\partial \HH^n$ intersects transversely $M'$ at $\partial M= \partial M'$ (in fact $\partial M = \partial M'$
  is contained in $\partial \HH^n$, and there is a small $\epsilon>0$ such that
  $$M'\cap \{x\in \HH^n; \ x_n<\epsilon\} = \partial M \times [0,\epsilon) \ . $$

 \cvd

 \begin{remarks}\label{on-proper}{\rm (1) Corollary \ref{partial-function} is also a 
 consequence of Proposition \ref{proper-emb}.  In fact $f$ given by the composition
 of $\beta$ with the projection onto the $x_n$ coordinate has the
 required property with value in some $[0,a)$, $a>0$, and to get
   $[0,1)$ is just a simple question of reparametrization.
   
   (2) A proof of Proposition \ref{proper-emb} can be obtained by using the open covering
   with associated partition of unity of Remark \ref{minimal-part}. For one can take
   $$\beta=(\beta_\partial,\beta_U): M\to (R^h\times \R)\times (\R^h\times \R)$$
   where
  $\beta_\partial$ is as above, $\beta_U =( \lambda'_\partial j_U, \lambda'_\partial)$
  and $j_U$ is the inclusion of $M$ into $\R^h$.
 }
   \end{remarks}   
   
   \smallskip
     
 {\bf The double of $M$.} Let $M'\subset \R^n$ be obtained from $M$ as
 in the proof of Proposition \ref{proper-emb}.  Let $M''$ be the image
 of $M'$ via the reflection
     $$(x_1,\dots,x_n)\to (x_1,\dots,-x_n) \ ; $$ $\partial
 M'=\partial M'' = \partial M$. Also $M''$ is diffeomorphic to $M$ and
 is a proper submanifold of $\{x_n \leq 0\}$. Then $D(M):=M'\cup M''$
 is compact smooth {\it baundaryless} manifold, containing both $M'$
 and $M''$ as submanifolds. $\partial M$ is given by the {\it
   transverse intersection} of $D(M)$ with $\partial \HH^n$. {\it
   Considered up to diffeomorphism} $D(M)$ only depends on $M$ (also
 considered up to diffeomorphism).  In this sense it is called the
 {\it double of $M$}.

\section{A fibration theorem}\label{fibre-teo}

\begin{proposition}\label{fib}{\rm
    (Fibration Theorem)} Let $M$ be a compact boundaryless smooth manifold and
  $f:M\to N$ a surjective summersion onto the connected manifold $N$. Let $q_0\in N$,
  $F= f^{-1}(q_0)$. Then $f$ is a smooth fibre bundle with fibre $F$.
\end{proposition}
\Dim  Let $q_0\in N$ and $F= f^{-1}(q_0)$. We know that $F$ is a submanifold
of $M$. Fix a tubular neighbourhood $(U,p)$ of $F$ in $M$.
Let $D$ be a small open disk in $N$ around $q_0$ such that $f^{-1}(D)\subset U$.
Define $h: f^{-1}(D)\to F\times D$, $h(x)=(f(x),p(x))$. Clearly, $f=p_D\circ h$,
where $p_D$ is the projection onto $D$. 
Moreover, $h(x)=(x,0)$ for every $x\in F$. As $f$ is a summersion, it is easy
to verify that the differential of $h$ is invertible on $f^{-1}(D)$ (possibly shrinking $D$).
As $h$ is essentially the identity on $F$, and the fibres are compact, an
usual argument (for instance like in the costruction of the tubular neighbourhoods)
shows that if $D$ is small enough, $h$ is a diffeomorphism, hence a local trivialization 
of $f$. If $q$ is an arbitrary point of $N$, we can cover a smooth arc joining
$q_0$ and $q$ in $N$ by a ``chain'' of similar local trivializations over
a chain $D=D_0, D_1, \dots , D_k$, $D_k$ around $q$, of small disks centred at the arc, 
$D_j\cap D_{j+1}\neq \emptyset$, so that one eventually deduces that the fibre 
$F'$ over $q_1$ is diffeomorphic to $F$.
Finally we have proved that   $f$ is a smooth fibration with fibre $F$.   

\cvd

\section{Density of smooth among $\Cc^r$-maps }\label{smooth-density}
Recall that for every $r\geq 0$, $\Cc^r(M,N)$ denotes the space of
$\Cc^r$ maps endowed with the weak topology; $\Ee^r(M,N)$ is the
subspace of smooth maps. We have
\begin{proposition}\label{smooth-dense}
  Assume that $M\subset \R^h$, $N\subset \R^k$ are boundaryless
  compact smooth manifolds.  Then for every $r\geq 0$, $\Ee^r(M,N)$ is
  dense in $\Cc^r(M,N)$.
\end{proposition}
\Dim Let $(U_M,p_M)$ and $(U_N,p_N)$ be respective tubular neighbourhoods. Let $(U,p)\subset (U_M,p_M)$ be a smaller
tubular neigbourhood (it just differs by a smaller ``$\epsilon$'', so that $p$ is the restriction of $p_M$).
Let $f\in \Cc^r(M,N)$. Consider the $\Cc^r$ extension $\hat f = f\circ p_M$. Apply Stone-Weierstrass Theorem \ref{stone-weier} 
to get a {\it polynomial} map $P: U_M\to \R^k$ which uniformely approximates (in the $\Cc^r$-topology) $\hat f$ on $U$ (which is
compact); we can also require that $P(U)\subset U_N$. Finally the restriction to $M$ of $p_N\circ P$ is a {\it smooth}
map from $M$ to $N$ which approximates $f$ in the $\Cc^r$-topology.

\cvd

By a very similar argument we have also
\begin{lemma}\label{near-homotopic} Let $M\subset \R^h$, $N\subset \R^k$ be compact boundaryless manifolds.
If $f\in \Ee^r(M,N)$ is close enough to $g\in \Cc^r(M,N)$ then they are $\Cc^r$-homotopic. If they are both smooth
then they are smoothly homotopic.
\end{lemma} 
\Dim If $f$ is close enough to $g$ we can assume that for every $p\in M$, for every $t\in [0,1]$,
$(1-t)g(p)+tf(p)$ belongs to $U_N$. Then $H(p,t)=p_N((1-t)g(p) + tf(p))$ is a required homotopy.

\cvd

\begin{remark}\label{partial-tube2}{\rm By using Remark \ref{partial-tube}, Proposition \ref{smooth-dense}
    and Lemma \ref{near-homotopic} hold true if $N$ is the interior of a compact manifold with boundary
    $\bar N$. Clearly they hold also if $N$ is an open set of $\R^k$}
  \end{remark}

\section {Smooth homotopy groups - Vector bundles on spheres}\label{VB-sphere} 
The above results have the following important application. 
Fundamental topological-algebraic invariants, the {\it homotopy groups} $\pi_n(X)$ , $n\geq 1$
(considered up to isomorphism)  are defined for every path connected topological space $X$
in terms of {\it continuous homotopy} classes of {\it continuous} maps $S^n \to X$. 
If $X= N \subset \R^k$ is as in above Remark \ref{partial-tube2}, then
Proposition \ref{smooth-dense} and Lemma \ref{near-homotopic}
imply that  we can equivalently define the homotopy groups of $N$ by using {\it smooth}
maps $S^n \to N$ up to {\it smooth} homotopy. If it is necessary to deal with pointed
maps, we can do it by using the smooth homogeneity of $N$.

Let us use these facts to classify (up to strict equivalence) the
embedded vector bundles on a unit sphere $S^m\subset \R^{m+1}$, $m\geq
2$. Let $\xi = f^*\tau_{n,k}$, for some smooth map $f:S^m \to
\GG_{n,k}$. Let us fix $1>\epsilon>0$. Set $D^+= S^m\cup \{x_{m+1}\geq
-\epsilon\}$, $D^-= S^m\cup \{x_{m+1}\leq \epsilon\}$. Clearly, both
$D^\pm$ are diffeomorphic to a closed $m$-disk, $S^m = D^+\cup D^-$,
$D^+\cap D^-$ is a tubular neighbourhood of the equatorial sphere
$S^{m-1}\subset S^m$, diffeomorphic to $S^{m-1}\times
[-1,1]$. We know by the Classification Theorem that the pull-back of
$\xi$ on $D^\pm$ via the respective inclusion maps is strictly
equivalent to the product bundle $D^\pm \times \R^k \to D^\pm$. Fix
two respective trivializations. The change of trivialization on
$D^+\cap D^-$ produces a smooth map
$$\rho_\xi: D^+\cap D^- \to {\rm GL}(k,\R)$$
and we consider its restriction (we keep the name)
$$\rho_\xi: S^{m-1}\to {\rm GL}(k,\R) \ . $$ As $D^+\cap D^-$ is
connected, the image of $\rho_\xi$ is contained in one of the two
connected components of GL$(k,\R)$ and up to strict equivalence we can
assume that this is the subgroup GL$^+(k,\R)$.  The arbitrary choices
made to define $\rho_\xi$ are the positive scalar $\epsilon$, the representative $\xi$ in
its strict equivalence class, the two trivializations. It is easy to
verify (by using the Classification Theorem) that the homotopy class
$[\rho_\xi]$ does not depend on these choices so we have well defined a
map
$$ \VV_{0,k}(S^m)\to [S^{m-1},{\rm GL}(k,\R)], [\xi]\to [\rho_{[\xi]}] \ . $$
If $m-1>1$, the (smooth) $\pi_{m-1}({\rm GL}^+(k,\R))$ is abelian, the choice
of a base point is immaterial, so that $[\rho_{[\xi]}]\in \pi_{m-1}({\rm GL}^+(k,\R))$.
If $m=2$, we have to take into account the base points say $p_0=e_1$ of $S^1$
and say $x_0=I_k$ of ${\rm GL}^+(k,\R)$ and work with {\it pointed} smooth maps.
However this is a minor technical point, we can manage it by using the smooth
homogeneity of ${\rm GL}^+(k,\R)$ (we skip the details), so that we can eventually
consider again $[\rho_{[\xi]}]\in \pi_1({\rm GL}^+(k,\R))$.
Summing up, for every $m \geq 2$, for every $k\geq 1$,  we have defined a map
$$\rho: \VV_{0,k}(S^m)\to  \pi_{m-1}({\rm GL}^+(k,\R))   \ . $$
We claim that {\it this map is bijective}. In fact we can exhibit $\rho^{-1}$.
This will be a particular case of Proposition \ref{emd-abs-VB}, see 
Section \ref{more-vect-sphere}.

\begin{remark}\label{complex-VB}{\rm The same construction works as well for complex
    embedded smooth vector bundles on $S^m$, by replacing GL$^+(k,\R)$ with
    GL$(k,\C)$ (which is connected), or also for bundles with ``reduced group''
    like for instance $SO(k)$.}
  \end{remark}

\section{Smooth approximation of compact embedded $\Cc^r$-manifolds}\label{smooth-Cr}
For every $r\geq 0$ there is a natural category of embedded
$\Cc^r$-manifolds and $\Cc^r$-maps ($\Cc^r$-diffeomorphisms) between
them. When $r=0$ we have the category of (embedded) {\it topological
  manifolds} and continuos maps (homeomorphisms). This presents its
own phenomena (including ``wild'' ones) that are beyond the aims and
the possibilities of this text.  On another hand, we are going to see
that to a large extent (at least in the compact case), for $r\geq 1$,
there are not essentially new phenomena with respect to the smooth
category. Basically this depends on the density of smooth maps already
established. 

$\bullet$ For $r\geq 1$, let $M\subset \R^h$ be a boundaryless
compact $\Cc^r$-manifold. The construction of the tubular neighbourhoods
of $M$ in $\R^h$ works verbatim in the $\Cc^r$-category.
It is enough to start with a $\Cc^r$-map $\nu: M\to \GG_{h,h-m}$
defining a distribution of transverse $(h-m)$-planes along $M$.
If we use for instance the standard metric $g_0$ on $\R^h$,
we obtain only a $\Cc^{r-1}$-map. However by applying the same argument
of the proof of  Proposition \ref{smooth-dense} we can approximate it
by a $\Cc^r$-map, keeping the transversality.
Assume that we have fixed one $(U,p)$.
We can summarize this by the following commutative diagramms (where for
simplicity we have written $\tau$ instead of $\tau_{h,h-m}$):
$$ \begin{array}[c]{ccc}
U&\stackrel{ F }{\rightarrow}& \Vv(\GG_{h,h-m})\\
\downarrow\scriptstyle{p}&&\downarrow\scriptstyle { \tau}\\
M&\stackrel{\nu}{\rightarrow}&\GG_{h,h-m} \end{array}$$ 
where $F= \nu^* \circ (f_\nu)^{-1}$. $F$ is a $\Cc^r$-map
and verifies the following properties (which are easy to check):
\begin{itemize}
\item $M=F^{-1}(\GG_{h,h-m})$, where $\GG_{h,h-m}\subset \Vv(\GG_{h,h-m})$
as the zero section.
\item The image of $F$ is contained in the interior of a compact submanifold with boundary
of the form $N_\epsilon(\GG_{h,h-m})$ for some $\epsilon >0$.
\item $F$ is {\it transverse to $\GG_{h,h-m}$ }, that is
for every $p\in M$, 
$$T_{F(p)}\Vv(\GG_{h,h-m})= T_{F(p)}\GG_{h,h-m} + d_pF(T_pU) \ . $$
\end{itemize}
\smallskip

This means that $M=F^{-1}(\GG_{h,h-m})$ can be considered as a sort of
``global equation'' defining $M$, which localizes in terms of very domestic equations:
for every given triavialization $\Phi: \tau^{-1}(W) \to W\times \R^{h-m}$
of the tautological bundle, 
we can consider the restriction of $\Phi \circ F$ obtaining a map
$$ (\Phi \circ F)^{-1} (W\times \R^{h-m})\to W\times \R^{h-m} \ . $$
Let $\pi: W\times \R^{h-m} \to \R^{h-m}$
 the projection. As $F$ is transverse to $\GG_{h,h-m}$ then
 $\pi\circ \Phi \circ F$ is a summersion (possibly shrinking $U$), and 
 $$(\Phi\circ F)^{-1}(W\times \{0\})= (\pi\circ \Phi \circ  F)^{-1}(0) \ . $$ 
 By the way this confirms that $M$ is a submanifold of $U$ of the correct
 dimension thanks to  Proposition \ref{global-imm-summ}.
\smallskip

$\bullet$ By  the density Theorem \ref{smooth-dense}, see also Remark \ref{partial-tube2},
we can uniformly approximate $F$ (in the $\Cc^r$-topology) 
on a slightly smaller compact tubular neighbourhood $U'\subset U$ with a {\it smooth}
map 
$$\tilde F: U' \to \Vv(\GG_{h,h-m}) \ . $$
As the transversality is manifestly a  $\Cc^1$-open condition,
if $\tilde F$ is close enough to $F$, then it is transverse to
$\GG_{h,h-m}$ and by applying to $\tilde F$ the above construction
and again  Proposition \ref{global-imm-summ}, we conclude
that $M':= \tilde F^{-1}(\GG_{h,h-m})$ is a compact submanifold of the interior of $U'$,
$\dim M'=\dim M$. Moreover, If $\tilde F$ is close enough to $F$ 
then the restriction of $p$ to $M'$ defines a $\Cc^r$-diffeomorphism $\rho: M' \to M$.
For as $p$ is the identity on $M$ this last claim follows by the very same
argument used in the construction of the tubular neighbourhood 
to show that $f_\nu: N_\epsilon(M)\to U$ is a diffeomorphism.
Note that $M'$ can be {\it arbitrarily $\Cc^r$-close to $M$} in the sense that
the $\Cc^r$-diffeomorphism $\rho^{-1}:M\to M'$ composed with the inclusion of $M'$ in $U'$ can be arbitrarily close
to the inclusion of $M$ of in $U'$.

\medskip

Summing up:
\begin{proposition}\label{smooth-struct} {\rm (Smooth approximation theorem)} 
For every $r\geq 1$, for every embedded compact  boundaryless $\Cc^r$-manifold $M\subset \R^h$
there is a smooth manifold $M'\subset \R^h$ $\Cc^r$-diffeomorphic to $M$. Moreover
$M'$ can be chosen arbitrarily $\Cc^r$-{\rm close} to $M$ (i.e. $M'$ is a {\rm smooth approximation} of $M$ in $\R^h$).
\end{proposition}

\smallskip

These smooth
structures are {\it unique up to diffeomorphism}. Precisely

\begin{proposition}\label{unique-smooth} {\rm (Uniqueness of smooth structure)} 
If $M$, $N$ are compact boundaryless embedded smooth manifolds 
which are $\Cc^r$-diffeomorphic, for some $r\geq 1$, then they are smoothly diffeomorphic.
\end{proposition}

\smallskip

In fact, if $f: M\to N$ is a $\Cc^r$-diffeomorphism, it can be approximated by a smooth map $\tilde f$
which is an injective immersion (because $r\geq 1$), hence it is a diffeomorphism.

\cvd

\section{Nash approximation of compact embedded smooth manifolds}\label{An-Nash}
By following carefully the above construction of $\tilde F$, we have
more information about its ``degree of smoothness''. Here we use some
notions recalled il Section \ref{real-alg-bund}.  We assume also that
the reader has a few basic knowledge of real analytic maps.  For the
notions of real (semi)-algebraic geometry we refer to \cite{BCR},
\cite{BR}.

Let $X\subset \R^k$ be a compact regular real algebraic set of
dimension $r$ (as a smooth manifold).  Let us specialize the
construction of a tubular neighbourhood in this algebraic situation.
If we use the standard metric $g_0$ on $\R^k$, then the associated map
$\nu: X \to \GG_{k,k-r}$ is {\it algebraic}. The map $f_\nu:
N_\epsilon(X)\to \R^k$ is algebraic.  The pull-back bundle $\nu^*\tau$
is algebraic. Hence the tubular neighbourhood projection $p:U\to X$ is
the composition of algebraic maps and of a map obtained by {\it
  inverting an algebraic map}.  According to Remarks
\ref{partial-tube} and \ref{partial-tube2} these considerations hold
also for the tubular neighbourhoods of a compact {\it regular
  ``semilagebraic'' set with boundary}, that is obtained as in Lemma
\ref {>=}, assuming that $X$ is a regular real algebraic set and the
function $f$ is algebraic (so that also the boundary is a real
algebraic set).  Then such a projection $p$ is not any smooth map.
A basic example of function of this type is $y=\sqrt{1+x^2}$ and we note that
its graph is a branch of the hyperbole defined by the polynomial equation
$y^2-x^2-1=0$. We
would say that it belongs to the {\it smallest class of maps
  containing the algebraic maps, closed by usual algebraic operations
  and for which the inverse map theorem
  and its corollaries hold true}. As algebraic maps are {\it real
  analytic}, and the inverse map theorem holds for real analytic maps,
then {\it $p$ is at least real analytic}.  But we have more. Recall
that by definition a {\it semialgebraic set} $Y$ in some $\R^n$ is
definable as the union of a {\it finite} family of subsets of $\R^n$
each one definable as the solution of a finite system of real
polynomial inequalities.  Obviuosly this extends the notion of
algebraic set. Fixing a few technical issues, by developing these
considerations one defines the subcategory of {\it Nash manifolds and
  maps} of the category of smooth embedded manifolds.  A Nash
$m$-manifold is an embedded real analytic $m$-manifold $M\subset
\R^n$, for some $n$, which is also a semialgebraic set; in particular
this implies that $M$ is contained in a real algebraic set $X$ of the
same dimension.  A Nash map $f:M\to N$ between Nash manifolds is a
real analytic map such that its graph is a semialgebraic set.
We say that a Nash manifold $M\subset \R^n$ is {\it normal} if it is contained
in the regular part $R(X)$, $X$ being as above. A normal compact boundaryless
Nash manifold $M$ is union of connected components of $R(X)$.
Although semialgebraic and analytically smooth, in general $M$ is not normal
but it has a {\it normalization} up to Nash diffeomorphisms.
More precisely we have the following very concrete description
of Nash manifolds and maps (see \cite{AM})
\begin{proposition}\label{artin-mazur} Let $M\subset \R^n$ be a connected Nash $m$-manifold
  and $f:M\to \R^h$ be a Nash map. Then there are:
  \begin{enumerate}
  \item An irreducible $m$-dimensional real algebraic set $X\subset \R^n\times \R^k$,
    for some $k$;
  \item A polynomial map $p:X\to \R^h$;
    
  \item A Nash manifold $M'\subset M\times \R^k$, such that $M' \subset R(X)$, and it is the graph
a Nash map
$g:M\to \R^k$, so that $\sigma(x)=(x,g(x))$ is a Nash diffeomorphism;

\item $f= p \circ \sigma$.
\end{enumerate}
  
\end{proposition}

\cvd

\smallskip

If $M$ and $N$ are Nash manifolds, Nash maps form
a subspace $\Nn^r(M,N)$ of $\Ee^r(M,N)$, for $r\geq 1$ and $\Nn(M,N)$ of
$\Ee(M,N)$; thanks to the inverse map theorem which holds for Nash
maps, a compact Nash manifold $M$ has Nash tubular neighbourhoods $(U,p)$ ($U$ is a compact
Nash manifold with boundary - possibly with corners - and $p$ is a
Nash map). With the very same proof of Proposition \ref{smooth-dense}
we have the following {\it density of Nash maps}.

\begin{proposition}\label{nash-dense}{\rm (Density of Nash maps)}
  Assume that $M\subset \R^h$, $N\subset \R^k$ are 
  Nash manifolds, $M$ compact boundaryless, $N$ the interior of a compact $\bar N$.  
  Then for every $r\geq 1$, $\Nn^r(M,N)$ is dense in $\Ee^r(M,N)$, $\Nn(M,N)$ in $\Ee(M,N)$.
\end{proposition}

\cvd

Let $M\subset \R^h$ be a compact smooth boubdaryless $m$-manifold and consider again
the commutative diagramm
$$ \begin{array}[c]{ccc}
U&\stackrel{ F }{\rightarrow}& \Vv(\GG_{h,h-m})\\
\downarrow\scriptstyle{p}&&\downarrow\scriptstyle { \tau}\\
M&\stackrel{\nu}{\rightarrow}&\GG_{h,h-m} \end{array}$$

$\GG_{h,h-m}$ is a regular real algebraic set,
$N_\epsilon(\GG_{h,h-m})$ is a compact regular semialgebraic set with
boundary contained in the regular real algebraic set
$\Vv(\GG_{h,h-m})$, hence we fix for it a Nash tubular neighbourhood
say $(U_\GG, p_\GG)$. The approximating map $\tilde F$ is of the form
$$p_{\GG} \circ P$$ where $P$ is a polynomial map (by application of
Stone-Weirstrass); $\tilde F$ is eventually a Nash map close to $F$,
then $M':= \tilde F^{-1}(\GG_{h,h-m})$ is a Nash manifold
$\Cc^\infty$-close to $M$.  So by adapting the very same construction
used to give a compact $\Cc^r$-manifold a smooth structure, we have
the following celebrated result by J. Nash \cite{Na}. A first
approximation theorem in this vein is due to Seifert \cite{Seif},
concerning the case of manifolds with product tubular neighbourhood.

\begin{theorem}\label{nash-teo} (1) {\rm (Nash approximation theorem)} Let $M\subset \R^h$ be a compact 
  connected smooth boundaryless manifold. Then there is a Nash manifold $M'\subset \R^h$
  diffeomorphic to $M$ and which can be chosen arbitrarily $\Cc^\infty$-close to $M$. Up to stabilize the embedding
  $M\subset \R^h \subset \R^h\times \R^k$, for some suitable $k$, we can assume that the Nash approximation
  $M'\subset \R^{h+k}$ is normal, that is $M'$ is union of connected components of $R(X)$, $X\subset \R^{h+k}$
  being a real algebraic set of the same dimension.

(2) {\rm (Uniqueness of Nash structures)} If two compact embedded boundaryless Nash manifolds $M\subset \R^h$, $N\subset \R^k$
are smoothly diffeomorphic, then they are Nash diffeomorphic to each other.
\end{theorem} 

\cvd

\medskip

\begin{remarks}\label{partial-nash}{\rm (1) Let $M$ be compact smooth with non empty boundary $\partial M$.
We can apply the  Nash approximation to a double $D(M)$ of $M$ (realized in $\R^n$ as above)
and get a boundaryless Nash manifold $D(M)'\subset \R^n$ close to $D(M)$. Then $M':= D(M)'\cap \HH^n$ is a Nash model 
(with boundary) of $M$.

(2) In his pyoneristic paper \cite {Na}, Nash stated also a few conjectures/questions towards
potential improvements of his result. The most natural conjecture was that
$M$ can be approximated by a regular real algebraic set. We will return on it in Section \ref{nash-tognoli}.
Another question concerned the existence of {\it rational} real algebraic models, see also Sections
\ref{stable-2-Nash}, \ref{tear}). }
\end{remarks}

\medskip

\subsection{On Nash vector bundle}\label{Nash-VB} By using the
classification theorem \ref{VBClassification},
the density of Nash maps, and Lemma \ref{near-homotopic} we readily
have (details are left as an exercise) the following existence and
uniqueness of Nash structures on smooth vector bundles. This answers
in the Nash category the analogue of (more demanding) questions posed
in  Remark  \ref{K-remark} (5) about real algebraic vector bundles.

\begin{proposition} Let $M$ be a compact embedded Nash manifold. Then

  (1) Every smooth embedded vector bundle on $M$ is strictly equivalent to a Nash vector bundle.
  \smallskip
  
  (2) If two Nash vector bundles on $M$ are smoothly strictly equivalent, then they are Nash
  strictly equivalent to each other.
  \end{proposition}

\begin{remark}\label{nash-utility}
  {\rm Beside its theoretic interest, approximation by Nash manifolds
    and density of Nash maps can be also of practical
    utility. Whenever we are interested in the density of smooth maps
    verifying a certain property, and we are in condition to apply
    Nash approximation and density of Nash maps, then it will be
    enough to show that Nash maps with the given property are dense
    among Nash maps. The main advantage is that we have a much
    stronger geometric control on the {\it image} of Nash than of
    arbitrary smooth maps. We will substantiate this remark in next
    sections.}
  \end{remark}

\medskip
    
The interested reader can find a lot of information about Nash
manifolds in \cite{BCR} and mostly in \cite{Shi}.

\section{Smooth and Nash Sard-Brown theorem}\label{Sard}
Let us recall some facts of analysis.  

\smallskip

(i) Every open set $U\subset \R^n$ is endowed with the ($n$-dimensional)
{\it Lebesgue measure} and this defines the class of {\it measure zero} i.e.  
{\it negligible} subsets of $U$.

\smallskip

(ii) If $X\subset U$ is negligible and $f:U\to W$ is a $\Cc^1$-map between open sets of $\R^n$, 
then $f(X)$ is negligible in $W$.

\smallskip

(iii) If $U'\subset U$ is an open subset and $X$ is negligible in $U$ then $X\cap U'$ is negligible
in $U'$.

\smallskip

(iv)  A countable union of negligible subsets of the open set $U$ is negligible.

\smallskip

(v) If $X$ is negligible in the open set $U$, then $U\setminus X$ is dense in $U$.

\smallskip

(vi) {\it (Fubini property)} If $U\subset \R^h\times \R^k$ , $X\subset U$ and for every 
$a\in \R^h$,$X\cap \{a\} \times \R^k$ is negligible in $U\cap \{a\}\times \R^k$, then
$X$ is negligible in $U$.

\smallskip

(vii) If $M$ is a smooth embedded $m$-manifold, we say that $X\subset M$ is {\it negligible in $M$}
if for every chart $\phi: W\to U\subset \R^m$, $\phi(X\cap W)$ is negligible in $U$.
Thanks to the above properties of negligible sets it is enough to check it on the open
sets of any countable atlas of $M$ (which certainly exists). 
{\it We stress that we have not defined any measure on $M$, we have just defined the class of negligible
subsets.}

\medskip

Let $f:M\to N$ be a smooth map between embedded smooth manifolds of dimension $m$ and $n$ respectively.
By definition a point $p\in M$ is {\it critical} for $f$ if ${\rm rank} \ d_pf < n=\dim N$.
Set $C(f)\subset M$ the set of critical
points of $M$. 
$$N\setminus f(C(f))\subset N$$ 
is the set of {\it regular values} of $f$ while $q\in f(C(f))$ is said a {\it critical
value} of $f$. The set $M\setminus C(f)$ is open (possibly empty) in $M$. If $M$ is compact,
$f(C(f))$ is compact, hence closed in $N$.  
Sard's theorem is a fundamental result for differential topology; in particular it is the  base
of {\it transversality} theory that we will develop later.

\begin{theorem}\label{sard-teo} {\rm (Sard's theorem)} Let $f: M \to N$ be a smooth map between embedded 
smooth manifolds. Then $f(C(f))$ is negligible in $N$.
\end{theorem}

In fact in differential topological applications one rather uses the following corollary, also known
as Brown's theorem.

\begin{corollary} \label{brown-teo} {\rm (Brown's theorem)} Let $f: M \to N$ be a smooth map between embedded 
smooth manifolds. Then $N\setminus f(C(f))$ is dense in $N$ (open and dense if $M$ is compact).
\end{corollary}

\medskip

{\bf Easy special cases.}
A special case of Sard's theorem is when $\dim M < \dim N$. Then $C(f)=M$. In this case the proof
is easy: clearly $M$ is negligible in $M\times \R^{n-m}$ and $f(M)= f \circ p_M(M)$
$ f \circ p_M: M\times \R^{n-m}\to N$, $p_M$ being the projection onto $M$. Then we can apply
the above property (ii).
\smallskip

A special and immediate case of Brown's theorem is when $M$ is the {\it finite} union of disjoint submanifolds
of $N$ of dimensions strictly less than $\dim N$, and $f$ is the union of the inclusion maps.

\cvd

A very readable proof of Sard's theorem, which {\it fully} employes the fact that $f$ is $\Cc^\infty$,  is in \cite{M1}. 
We stress that it is a result of analytic nature and rather delicate. To appreciate better this point, 
let us recall the following Morse-Sard $\Cc^r$ generalization.  

\begin{theorem}\label{morse-sard} {\rm (Morse-Sard theorem)} Let $f:M\to N$ be a $\Cc^r$-map between embedded
 smooth manifolds. If $r> \max \{0,m-n\}$ then $f(C(f))$ is negligible in $N$.
 \end{theorem}
 
 The condition which relates the ``degree of regularity'' of $f$ and
 the dimensions of the manifolds is sharp. Whitney \cite{Whit} has
 constructed an example of a $\Cc^1$-function $f:\R^2\to \R$ such that
 $C(f)$ contains a subset $J$ homeomorphic to an open interval, and
 that $f$ is not constant on $J$. Hence $f(C(f))$ contains an open
 interval. A proof of the Morse-Sard theorem can be found in \cite{H}.
 
 \subsection{A Sard-Brown theorem in the Nash category}
 Here is a Nash version of the Sard-Brown theorem, whose statement is purely geometric.
 
 \begin{theorem} \label{nash-sard} Let $f: M\to N$ be a Nash map between embedded Nash manifolds.
 Then $f(C(f))$ is the union of a finite set of Nash submanifolds of $N$ of dimensions stricly
 less than $\dim N$. 
 \end{theorem}
 
 \begin{remark}\label{nash-suffice}
   {\rm Assume that $M$ and $N$ are embedded smooth manifolds such
     that we can apply to both the Nash approximation by means of Nash
     manifolds $M'$ and $N'$, so that $\Nn(M',N')$ is dense in
     $\Ee(M',N')$. It follows that the set of smooth maps $f: M\to N$
     which verify Brown's theorem is dense in $\Ee( M,N)$. In many
     applications this suffices}
 \end{remark}
  
 \smallskip
 
 {\it Outline of a proof of Theorem \ref{nash-sard}.} 
 Alike the statement of the theorem, 
 it is of purely geometric nature. For all details 
 one can look at \cite{BCR}.
 Let us recall the following basic facts about semialgebraic sets:
 
 (1) We know that every embedded Nash manifold is in particular a semialgebraic set.
 
 (2) Every semialgebraic set $X\subset \R^n$ is the union of a {\it finite} number
 of disjoint connected Nash embedded manifolds.
 
 (3) If $X\subset M$ is a semialgebraic subset of the embedded Nash manifold
 $M$, and $f:M\to N$ is a Nash map between Nash manifolds, then $f(X)$
 is a semialgebraic subset of $N$. This is a formulation adapted to our
 situation (and in fact a corollary) of the celebrated  
 {\it Tarski-Seidenberg theorem} that the projection in $\R^{n-1}$
 of a semialgebraic set $X$ in $\R^n$ is a semialgebraic set of $\R^{n-1}$.
 Moreover, all Nash manifolds making a partition of $f(X)$ as in (2)
 have dimension less or equal $\dim M$.
 
 \medskip
 
 Let us come to the proof of Theorem \ref{nash-sard}.  Let $f:M\to N$
 be our Nash map between embedded Nash manifolds.  As $f$ is a Nash
 map, it is not hard to check that $C(f)$ is a semialgebraic subset of
 $M$. By applying point (2), one realizes that $C(f)$ is the {\it
   finite} union of disjoint connected Nash submanifolds each one, say
 $Y$, verifying the following property: there exists $0\leq k <\dim N$
 such that for every $p\in Y$, ${\rm rank} \ d_p
 f_{|Y}=k$. $f(C(f))\subset N$ is the union of the images $f(Y)$'s
 hence it is a semialgebraic subset of $N$.  By point (2) again, it is
 the disjoint union of a finite number of disjoint connected Nash
 submanifolds of $N$. We claim that for every such a manifold, say
 $Z$, $\dim Z < \dim N$. If for example $N=\R$, then the restriction
 of $f$ on every $Y$ has vanishing differential, hence $f$ is constant
 on $Y$, so that $f(C(f))$ is a {\it finite} subset of $\R$. In
 general we can assume that $Z\subset f(Y)$ for some $Y$ as above, and
 $\dim Z = \dim N$ would be against the {\it constant rank theorem}
 \ref{constant-rank}.
 
 \cvd

 \begin{remark}\label{nash-utility2}
   {\rm We continue in the vein of Remark \ref{nash-utility}.  The
     Nash Sard-Brown theorem is an important example of application of
     the stronger geometric control on the images of Nash maps. Merely
     {\it continuous} maps (between open sets of some euclidean space)
     can have ``wild'' behaviour (i.e. anti intuitive with respect to
     an ``ordinary'' geometric intuition).  Let us recall for instance
     the so called {\it Peano's curves}, i.e.  {\it surjective}
     continuous maps $g:[0,1]\to [0,1]^2$. Wild phenomena make the
     category of topological manifolds much delicate to deal with. By
     Sard's theorem (easy case) {\it there are not smooth Peano's
       curves}.  In the Nash situation, even better the image of any
     such a Nash $g$ is a finite union of points or Nash
     $1$-manifolds. Smooth maps (and manifolds), although much more
     ``tame'' than merely $\Cc^0$ ones, are suited to topological
     considerations because they are very ``flexible''. This is
     basically due to the existence of bump functions and the {\it
       flatness} phenomenon that they incorporate. On another hand,
     this also implies for example that subsets of a smooth manifold
     defined by a finite set of smooth equations or inequalities can
     be weird: for instance one can prove that {\it every} compact
     subset of $\R^n$ can be realized as the zero set of a smooth
     function. In a sense this means that the formulation of the
     smooth Sard's theorem in {\it measure} theoretic terms, is the
     best one can say in general about the image of the critical set.
     The situation is dramatically simpler and geometrically friendly in
     the Nash case. It can be profitable to combine the flexibility
     of smooth manifolds with the Nash approximation and the density
     of Nash maps (whenever they can be applied).}
   \end{remark}

\section{Morse functions via generic linear projections to lines}\label{dense-morse}
Let $M$ be a compact boundaryless embedded smooth $m$-manifold.
\begin{definition}\label{morse}{\rm A smooth function $f:M\to \R$ is a {\it Morse function} if
    it has only non degenerate
critical points.}
\end{definition}
\smallskip

According to Chapter \ref{TD-LOCAL}, the notion of non degenerate
critical point $p$ of a determined {\it index} say $\lambda$ can be
defined on any representation in local coordinates of $f$ at $p$ (as
it does not depend on the choice of the local coordinates). By Morse
Lemma, the non degenerate critical points are isolated, hence by
compactness every Morse function on $M$ has only a finite number of
critical points. At least one of them is certainly a minimum (of index
$\lambda = 0$) at least one is a maximum (of index $\lambda = m$). A
Morse function on $M$ is {\it generic} if distinct critical points
take distinct (critical) values. In such a case we can order the
critical points $p_0,p_2,\dots, p_r$ so that $c_j:=f(p_j)<f(p_{j+1})=:
c_{j+1}$. Up to a linear reparametrization of the image, sometimes we
assume also that $f(M)=[0,1]$.

We want to prove that Morse functions exist and moreover are open and
dense in $\Ee(M,\R)$.

\begin{lemma}\label{open-morse}
  Let $M\subset \R^h$ be a compact boundaryless smooth manifold.  The
  set of Morse functions on $M$ is open in $\Ee(M,\R)$.
\end{lemma}
\Dim Let $f:M\to \R$ be a Morse function, with critical points
$p_1,\dots, p_k$. Fix a nice atlas of $M$ such that every critical
point $p_j$ is contained in a $B_j$ of some normal chart and these
$B_j$'s are pairwise disjoint. If $g$ is close enough to $f$ (in the
$\Cc^1$ topology) then it has no critical points on the compact set
$M\setminus \cup_j B_j$.  Let us analyze the local representation of
$f$, say $\hat f_j$, defined on the compact set $\bar U_j:=\phi_j(\bar
B_j)\subset \R^m$, for every $j=1,\dots, k$. On $\bar U_j$, the
positive smooth function
$$a_{\hat f_j}(x):=||d_x \hat f_j||^2+ (\det (\frac{\partial^2 \hat
  f_j}{\partial x_i \partial x_j}(x))^2$$ never vanishes, because the
first term vanishes only at $0=\phi(p_j)$, and the second term does
not vanish because the critical point is non degenerate. By
compactness, there is $d>0$ such that, for every $x\in \bar U_j$,
$a_{\hat f_j}(x)> d$.  If $g$ is close enough to $f$ in the $\Cc^2$
topology, then $a_{\hat g_j}(x) > d/2$, hence also $g$ has only non
degenerate critical points on $\bar B_j$. As there is a finite number
of critical points of $f$, we readly conclude that if $g$ is close
enough to $f$ in the $\Cc^2$ topology, then $g$ is a Morse function.

\cvd

\medskip

Let $M\subset \R^h$ be as above. For every linear function $L\in
(\R^h)^*$,
$$L(x)=a_1x_1+\dots + a_hx_h$$ corresponding to $(a_1,\dots,
a_h)\in M(1,h,\R)$) consider the restriction $L_M$ to $M$. We have

\begin{theorem}\label{generic-proj}
  Let $M\subset \R^h$ be a compact boundaryless smooth manifold. Then
  for every $f\in \Ee(M,\R)$, there is a open dense subset $\Ll_f$ of
  $(\R^h)^*$ such that for every $L\in \Ll_f$, $f+L_M$ is a Morse
  function.
\end{theorem}

\begin{corollary}\label{open-dense-morse}
  Let $M\subset \R^h$ be a compact boundaryless smooth manifold.
  Then:

  (1) There is a open dense set $\Ll$ in $(\R^h)^*$
such that for every $L\in \Ll$, $L_M$ is a Morse function.

(2) The set of generic Morse functions is a open dense set in $\Ee(M,\R)$.
\end{corollary} 

\medskip

{\it Proof of Corollary \ref{open-dense-morse}.}  (1) is a consequence
of Theorem \ref{generic-proj} applied to the costant function
$f=0$. Theorem \ref{generic-proj} together with Lemma \ref{open-morse}
implies that the set of Morse functions is open and dense in
$\Ee(M,\R)$ (if $L$ is close to zero, then $f+L_M$ is close to $f$). 
It is evident that generic Morse functions form an open
set in the set of Morse functions. Then it remains to show that
generic Morse functions are dense. Let $f:M\to \R$ be a Morse
function. Assume that there is a critical point $p$ which shares
the value with another one. It is enough to show  that arbitrarily close to $f$ there is 
a Morse function
$g$ with the same set of critical points of $f$, such that $g(p)\neq
g(p')$ for any other critical point $p'$. Then we conclude
by induction on the number of sharing  value critical points.
Let
$(W,\phi)$ be a normal chart centred at $p$, such that $W$ does not
contains other critical points of $f$. 
Let
$\gamma$ be the global bump functions on $M$ associated to this normal
chart.  For every $\epsilon \neq 0$, set $g_\epsilon= f+ \epsilon
\gamma$. Clearly, if $|\epsilon|$ is small enough, then $g_\epsilon$ is
close to $f$ (because $M$ is compact), hence it is a Morse function.
It is also clear that $g_\epsilon$ coincides with $f$ outside the
compact support of $\gamma$ (contained in $W$). 
A discrepancy between the sets of critical points
could
only occur on the support of $\gamma$.  But for every $x\in U$, $d_x
\hat g= d_x \hat f + \epsilon d_x \gamma_{1/3,1/2}$.  On $B^m(0,1/3)$
this reduces to $d_x \hat f$, hence $p$ is the only critical point of
$g_\epsilon$ on $B\subset W$ (with the usual notations about normal
charts). The function $\hat f$ has no critical points on the compact
set $\overline{B^m(0,1/2)\setminus B^m(0,1/3)}$, hence if $|\epsilon| >0$ is
  small enough the same fact holds for $g_\epsilon$. 
  Finally, by the finiteness of the critical set,  it is clear that we can take 
  $|\epsilon|$ small enough so that
  $g_\epsilon(p)$ differs from any other critical value. 
  
  \cvd

\medskip

{\it Proof of Theorem \ref{generic-proj}.} By the Nash approximation
theorem and the density of Nash functions it is not restrictive to
assume that $M\subset \R^h$ is a Nash $m$-manifold, and that $f:M\to \R$
is a Nash function. We will give a proof based on the Nash version of
Sard-Brown theorem. For a reader who would prefer a purely smooth
proof, we will indicate in parallel how to manage it by means of the
ordinary Sard-Brown theorem. Let us start with a local Lemma.
\begin{lemma}\label{local-morse} Let $f:U:=B^m(0,1)\to \R$ be a Nash
  function. Then there is a negligible subset $X$ of  $(R^m)^*\sim M(1,m,\R)$
  such that for every  $L\in (R^m)^*\setminus X$, $f+L_U$ is a Morse function.
\end{lemma}
\Dim
The differential
$$df: U \to M(1,m,\R)$$
is a Nash map. For every $L$, for every $p\in U$,
$p$ is a critical point of $f+L_U$ if and only if
$d_pf=-L$. $-L$ is regular value of $df$ if and only if
for every $p\in U$ such that $d_pf=-L$,
$$d_p(d_pf)= (\frac{\partial^2 f}{\partial x_i\partial x_j}(p))_{i,j=1,\dots m}\in M(m,\R)$$
is invertible. Hence, $-L$ is a regular value of $df$ if and only if all the
critical points of $f+L_U$ are non degenerate, that is $f+L_U$ is a Morse function.
We conclude by means of the Nash Sard-Brown theorem.

\cvd

\medskip

In the smooth case we have the same Lemma with the same proof, by using the smooth
Sard-Brown theorem. 


Let $M\subset \R^h$ be a compact Nash $m$-manifold as above. $M$
is covered by a finite set of
Nash Monge charts (this depends on the compactness of $M$ and on the inverse
function theorem which holds in the Nash category). Possibly reordering
the coordinates of $\R^h$, the corresponding Nash local Monge parametrization of
$M$
is of the form
$$ U:=B^m(0,1)\to (x,\psi(x))\in M\subset \R^m\times \R^{h-m}$$
so that the associated local representation of $f$ is the Nash function
$$\hat f(x_1,\dots x_m)= f(x_1,\dots,x_m,\psi(x_1,\dots,x_m)) \ . $$
 Let us write every $L\in M(1,h)$
in the form
$$L(x)=(a_1x_1+\dots +a_mx_m) + (a_{m+1}x_{m+1}+\dots + a_{h}x_{h}):=$$
$$\alpha(x_1,\dots,x_m)+ \beta(x_{m+1},\dots, x_h)$$
then the corresponding local representation of $f+L_M$ is
$$(\hat f(x_1,\dots,x_m)+ \beta(\psi(x_1,\dots, x_m)) + \alpha(x_1,\dots,x_m):= \hat f_{\beta}+ \alpha_U \ . $$
For every fixed $\beta\in M(1,h-m,\R)$, let us vary $\alpha\in M(1,m,\R)$ and apply Lemma \ref{local-morse}
to $\hat f_{\beta}$. Then for every $\beta$ the subset $C_\beta \subset M(1,m,\R)$ of $\alpha$'s such that
$\hat f_\beta + \alpha_U$ is not a Morse function consists of a finite number of disjoint Nash submanifolds of $M(1,m,\R)$
of dimension $<m$. Also the subset $C_f$ of $M(1,h)$ such that the restriction of $f+L_M$ to the given Monge chart
is not Morse is a semialgebraic subset, hence it is the finite union of disjoint Nash submanifolds of $M(1,h,\R)$.
It is also the union of the slices $C_\beta$, $\beta$ varying in $M(1,h-m, \R)$. As every $C_\beta$ is union
of manifolds of dimension $<m$, then $C_f$ is union of manifolds of dimension $< h$. As there is a finite number of
Monge charts, there is a finite number of such sets $C_f$ in $M(1,h,\R)$. The complement $\Ll_f$ of their union is dense
in $M(1,h,\R)$ and for every $L\in \Ll_f$, $f+L_M$ is a Morse function.

In the smooth case, the dimensional consideration about $C_f$ is replaced by the conclusion that
it is negligible, by using this information about every slices and the Fubini property (vi) recalled at the beginning
of this section.

\cvd

\subsection {Manifolds with boundary} Let $M$ be a compact smooth manifold with boundary $\partial M$, 
and let us fix a partition $\partial M = V_0\cup V_1$ as in  Corollary \ref{partial-function2}. By this Corollary
we know that the set, say $\Ee(M,V_0,V_1;\R)$, of smooth functions $f:M\to [0,1]$
such that $f^{-1}(j)=V_j$, $j=0,1$, and without critical points near $\partial M$ is non empty.
We can extend the results obtained in the boundaryless case.
\begin{proposition}\label{partial-morse-dense} The generic Morse functions belonging to 
$\Ee(M,V_0,V_1;\R)$ form an open dense set.
\end{proposition}

The only point that needs some further considerations is the existence of such relative Morse functions.
By using the notations of Remark \ref{minimal-part}, via the proper embeddings and the double of $M$, 
the results in the boundaryless case tell us that there are arbitrarily small linear projections $L$ which restrict
to Morse functions on $U$. If $f$ belongs to $\Ee(M,V_0,V_1;\R)$  and $L$ is small enough, then
$\lambda_\partial f + \lambda'_\partial L$ provides a Morse function closed to $f$; details are left
as an exercise.
\cvd

\section{Morse functions via distance functions}\label{distance-function}
The use of generic linear projections to line is a geometrically transparent way to
produce Morse functions on a compact embedded smooth manifold. Here we outline another
natural way based on distance functions. Let $M\subset \R^h$ be compact boundaryless as usual.
For every $q\in \R^h$ consider the smooth (actually polynomial) function
$$\delta_q: \R^h \to \R, \ \delta_q(x):= ||x-q||^2 \ . $$
We have
\begin{theorem}\label{d2-morse} There is an open and dense set $\Omega \subset \R^h$ such that
for every $q\in \Omega$, the restriction of $\delta_q$ to $M$ is a Morse function.
\end{theorem}

{\it Sketch of proof.}  Consider $\nu: M\to \GG_{h,h-m}$ corresponding to the distribution
of normal $(h-m)$-planes with respect to the standard metric $g_0$ on $\R^h$.
Let 
$$f_\nu: \nu^*(\Vv(\GG_{h,h-m}))\to \R^h, \ f_\nu(p,v)=p+v $$ be the map already used to construct a tubular
neighbourhood of $M$ in $\R^h$. One proves that the restriction of $\delta_q$ to $M$
has some degenerate critical point if and only if $q$ is not a regular value of $f_\nu$
(the reader can try to prove this by exercise; anyway all detalis can be found in \cite{M2} Part 1-6).
Then we conclude by applying the favourite version of Sard-Brown theorem.  

\cvd

\medskip

\subsection {Exhaustive sequences of compact submanifolds of non compact manifolds}\label{exhaustive}
The argument of Theorem \ref{d2-morse} applies also to any boundaryless non compact submanifold  $N\subset \R^h$ 
which is also a {\it closed subset} of $\R^h$. 
Then by using a generic $\delta_q$, we can find a sequence o increasing 
regular values $c_n$, $c_n\to +\infty$,
of the restriction of $\delta_q$ to $N$ such that every 
$$N_n:= \{ x\in N; \ \delta_q(x)\leq c_n \} $$
is a compact submanifold with boundary of $N$,  $N_n \subset N_{n+1}$ and $\cup_n N_n= N$.
That is we have an {\it exhaustive sequence of nested compact submanifolds with boundary of $N$.}
Every compact subset of $N$ is contained in some $N_n$. In particular,
If $f: M\to N$ is a $\Cc^r$ or a $\Ee$-map, $M$ being compact,
then there is $n$ such that
$f(M)\subset N_n$ and we can extend the density result of $\Ee(M,N)$ in $\Cc^r(M,N)$. If all involved
manifolds are Nash we have the density of $\Nn(M,N)$ in $\Ee(M,N)$ as well.  We can also extend to $N$
the notion of tubular neighbourhood. Fix a sequence a tubular neighbourhoods $\pi_n: U_n\to N_n$ constructed with respect to the standard
metric $g_0$ on $\R^h$ and a suitable decreasing sequence of $\epsilon_n >0$.  For every smooth positive function $\epsilon: N \to \R^+$ 
denote by $N_\epsilon:=  \{ x\in \R^h; \ d(x,N)<\epsilon (x)\} $ that is the $\epsilon$-{\it neighbourhood} of $N$
with respect to the euclidean distance. Then we can find such a function $\epsilon$ such that for every $x\in N_n$, $\epsilon(x)<\epsilon_n$
so that the restriction of the projections $\pi_n$ to $N_\epsilon$ match with the projection $\pi: N_\epsilon \to N$ such that $\pi(y)\in N$
is the nearest  point to $y$ on $N$. 

\section{Generic linear projections to hyperplanes}\label{generic-projection} 
Let $M\subset \R^h$ be a compact boundaryless $m$-manifold as above.
We have seen that generic linear projections of $M$ to $1$-dimensional subspaces of $\R^h$ are Morse
functions. Here we consider projections to hyperplanes, when the {\it codimension} $h-m$ is big enough. 
Precisely, let $\R^{h-1}\subset \R^{h-1}\times \R$; for every 
$v\in S^{h-1}\setminus \R^{h-1}$, let $p_v: \R^h=\R^{h-1}\oplus \ {\rm span}(v)\to \R^{h-1}$ be
the associated projection. We have

\begin{proposition}\label{2n} (1) If $h>2m$, then there is an open dense subset $I_M \subset S^{h-1}$ such that
for every $v\in I_M$, the restriction of $p_v$ to $M$ is an immersion.

(2) If $h>2m+1$, then there is  an open dense subset $E_M \subset S^{h-1}$ such that
for every $v\in E_M$, the restriction of $p_v$ to $M$ is an embedding.
\end{proposition}



\Dim  (1) Let $UT(M)\subset M\times S^{h-1}$ the total space of the
unitary tangent bundle of $M$ (constructed by using the standard metric $g_0$ on $\R^h$).
Let $t: UT(M)\to S^{h-1}$ the restriction of the projection $M\times S^{h-1}\to S^{h-1}$.
Then the restriction of $p_v$ to $M$ fails to be an immersion if and only if
$v$ belongs to the image of $t$. $\dim UT(M)=2m-1<h-1$. Hence  $S^{h-1}\setminus t(UT(M))$
is open and dense (by the easy case of Sard's theorem). This achieves point (1).

(2) The diagonal $\Delta$ is a closed subset of $M\times M$. Consider the smooth map
 defined on the complementary open set 
$$ \beta: M\times M \setminus \Delta \to S^{h-1},\  \beta(x,y)=\frac{x-y}{||x-y||} \ . $$
Then the restriction of $p_v$ to $M$ is not injective if and only if $v$ or $-v$ belongs
to the image of $\beta$. $\dim (M\times M \setminus \Delta)=2m< h-1$.
Hence $S^{h-1}\setminus {\rm Im}(\beta)$ is a dense subset. 
Its intersection with the dense open set $S^{h-1}\setminus t(UT(M))$ is also dense.
 Then we have a dense set of $v$'s  such that the restriction of $p_v$ is an injective
immersion, hence an embedding of $M$ because it is compact. Finally this
set of $v$'s is also open because the set of embeddings is open.
\cvd

The Morse projections to lines, and the above special cases of
projections to hyperplane are the simplest instances of the general
problem of understanding ``generic'' linear projections of compact
embedded smooth manifolds to lower dimensional subspaces.  An
interested reader can look at the definetly more advanced paper
\cite{Ma}.

\subsection{Truncated classifying maps}\label{truncated-VB} 
The classification theorem \ref{VBClassification}, has been formulated in terms of the limit
grassmannians $\GG_{\infty,k}$; however we know that every classifying map $f: M \to \GG_{\infty,k}$
factorizes through some $\hat f: M \to \GG_{n,k}$  (similarly for homotopies between
maps defining strictly equivalent vector bundles), but a priori $n$ might  vary with $M$.
In fact, arguing similarly to the weak immersion/embedding theorem, we show that
there is a ``uniform truncation'' depending only on the dimension.

\begin{proposition} Let $M$ be a compact embedded $m$-manifold.
(1) Then every $f:M\to \GG_{\infty,k}$ is homotopic to a map $g$ which factorizes through
a map $\hat g: M\to \GG_{m+k+1,k}$. 

(2) Two homotopic classifying maps with values in $\GG_{m+k+1,k}$ are homotopic
via a homotopy which factorizes through a map in $\GG_{m+k+2,k}$.
\end{proposition}
\Dim Start with $\hat f: M\to \GG_{n,k}$, with $n>m+k+1$. Hence the corresponding
bundle is embedded into $M\times \R^n$. Consider  linear projections
 $p_v:\R^n \to \R^{n-1}$, as above, and the maps 
$$F_v: M\times \R^n \to M\times \R^{n-1}, \   (x,v)\to (x,p_v(v)) \ . $$
For a generic $v$, $F_v$ embedds the vector bundle into $M\times \R^{n-1}$;
this corresponds to a map $M \to \GG_{n-1,k}$ homotopic to the given one by
the classification theorem. Similar considerations hold for homotopies.

\cvd

\chapter{The category of smooth manifolds}\label{TD-SMOOTH-MAN}
{\it Abstract} smooth manifolds and smooth maps between them will be introduced  by taking as 
{\it definition} some properties verified by embedded ones. We will see in Section \ref{emb-abs-comp}  
that abstract compact manifolds can be embedded in some $\R^n$.
As we are are mainly interested in compact manifolds, considered up to
diffeomorphism,  this abstraction would appear to be a bit superfluous. However there are some
good reasons to proceed. There are natural  constructions (quotients, ``cut-and-paste'', $\dots$, 
we will see them)  to build new {\it abstract} manifolds, starting from given ones (even embedded, 
even staying in the realm of compact manifolds).
It would be artificial to force them to deal from the beginning  in the 
embedded setting. It is more convenient to use the embedding result {\it ex post},
in order to exploit the facts already established for compact embedded manifolds.

\begin{definition}\label{abs-smooth-man}{\rm A topological space $M$ is a $m$-{\it smooth manifold} 
(we will omit the adjective ``abstract") if:

\begin{itemize}
\item $M$ is Hausdorff and with a countable basis of open sets.
\item $M$ admits an {\it smooth atlas} $\Uu=\{W_j,\phi_j\}_{j\in J}$ ($J$ being any set
of indices); that is 

\noindent (i) $\{W_j\}_{j\in J}$ is an open covering of $M$;

\noindent (ii) every {\it chart} $\phi_j:W_j \to U_j\subset \R^m$ is a {\it homeomorphism}
onto a open set of $\R^m$ (denote by $\psi_j:U_j \to W_j$ the inverse {\it local parametrization});

\noindent (iii) for every $i,j \in J$,
$$\phi_j\circ \psi_i : \phi_{i}(W_i\cap W_j)\to   \phi_j(W_i\cap W_j)$$
is a smooth {\it diffeomorphism}.
\end{itemize}

We summarize this item by saying that $M$ is (smoothly) {\it locally $m$-euclidean}.
}
\end{definition}

\medskip

\begin{remarks}\label{firs-rem}
{\rm (1) Every smooth atlas $\Uu$ of $M$ is contained in and implicitely determines a unique maximal smooth atlas 
$\Aa=\Aa_M$;  this is identified
with a specific {\it smooth structure on $M$}. 
In the embedded case $M\subset \R^n$, the charts of $\Aa$ were smooth by themselves,
referring to the smooth structure of the ambient euclidean space. In the abstract case every single chart is only a homeomorphism;
the smooth structure is enterely carried by the changes of local coordinates. Nevertheless, this is enough to deduce for example that
the {\it dimension} $m$ is well defined, alike the embedded case.

(2) Obvioulsy every embedded smooth manifold is a smooth manifold.
\smallskip

(3) Being locally euclidean does not imply any of the global topological requirements of the first item. For example
consider $M=\R^m \times (\R,\tau_d)$ where the second factor is endowed with the {\it discrete topology}.
Then $M$ is Hausdorff and locally   $m$-euclidean, but it has no countable basis of open sets.
On another hand, consider on $\R\times \{0,1\}$ (with the product topology) the equivalence relation
such that $(x,j)\sim (y,i)$ if and only if either $(x,j)=(y,i)$ or $x=y$ and $x>0$. Let $M$ be the quotient topological space.
Then $M$ is $1$-locally euclidean and has a countable basis of open sets, but it is not Hausdorff. In fact
the two points $[(0,0)]\neq [(0,1)] \in M$ cannot be separated by disjoint neighbourhoods. 
$M\times \R^k$ presents the same phenomenon in arbitrary dimension. 
\smallskip

(4) The fact that 
``locally euclidean" does not imply Hausdorff poses some principle question  when one uses manifolds
as model of some physical space or space-time. Local observations can support the idea that phenomena live
in a locally euclidean environment, but it is much more arbitrary to assume also the (global) separation property.
For example in some models of space-time one does not assume a priori that it is Hausdorff, and this
property is derived a posteriori as consequence of certain  global ``causality assumptions''  which look founded on some reasonable
physical (or philosophical) considerations \cite{HE}. To our aims, we do not hesitate to make these topological assumptions;
as the theory is already  rich, there are no reasons to renouce say the limit uniqueness or the
equivalence between compact and sequentially compact spaces.
 }
\end{remarks}  

\smallskip

\begin{definition}\label{abs-smooth-map}{\rm Let $f:M\to N$ be a continuous map between 
smooth manifolds of dimension $m$ and $n$ respectively. A {\it representation in local coordinates} of $f$ is of the form
 $$\hat f = \phi' \circ f \circ \psi: U \to U'$$
 where $\phi:W\to U \subset \R^m$ is a chart of $\Aa_M$, $\phi': W' \to U'\subset \R^n$ is a chart of $\Aa_N$,
 $f(W)\subset W'$. Then $f$ is {\it smooth} if for every $p\in M$ there is a local representation of $f$ such that
 $p\in W$ and $\hat f$ is a smooth map between open sets of euclidean spaces. The map $f$ is a {\it diffeomorphism}
 if it is a homeomorphism and both $f$ and $f^{-1}$ are smooth.}
 \end{definition}
 
 The following Lemma is an easy consequence of the definitions and of the basic fact that the composition of smooth 
  maps between open sets of euclidean spaces is smooth (details are left as an exercise).
  \begin{lemma} If $f: M \to N$ is a smooth maps between smooth manifolds, then every local representation
  of $f$ in local coordinates is smooth.
  \end{lemma}
  
  \cvd
  
  \smallskip
  
  Obviuosly smooth maps and diffeomorphisms between embedded  manifolds fulfill the above
  definition. So we have introduced the {\it category of smooth manifolds and smooth maps (diffeomorphisms)}
  which extends the embedded one.
\smallskip
  
  Let us describe some constructions that naturally produce (abstract) smooth manifolds.
  
  \smallskip

 (1) {\it (Quotient manifolds)} Let $\tilde M$ be a smooth manifold (even embedded). Let $G$ be
  a subgroup of the group Aut$(\tilde M)$ of smooth automorphisms of $\tilde M$.
  Assume that $G$ acts {\it freely} and {\it properly discontinuously} on $\tilde M$.
  This means that for every $p\in \tilde M$, the identity is the only element of $G$
  that fixes $p$, and that for every compact subset $K$ of $\tilde M$, the set of $g\in G$
  such that $K\cap g(K) \neq \emptyset$ is {\it finite}. Let $M:=\tilde M/G$ be
  the quotient topological space. It is known that $M$ is Hausdorff and with countable
  basis. Moreover, the projection $\pi: \tilde M \to M$ is a covering map. We can assume that
  for every $p\in M$, there is a open connected neighbourhood $W$ of $p$ such that
  the restriction of $\pi$ to every connected component $\tilde W$ of $\pi^{-1}(W)$
  is a homeomorphism, and $(\tilde W,\phi)$ belongs to $\Aa_{\tilde M}$. Then
  by varying $p$ in $M$, $\{ (W, \phi \circ \pi^{-1}) \}$ is a smooth atlas of $M$,
  such that $\pi$ becomes a smooth, locally diffeomorphic map.  
\smallskip

(2) {\it (Grassmann manifolds again)} We have already defined the projective spaces $\PP^k(\R)$
as special instances of (embedded) grassmann manifolds. There is another
classical way to obtain it. Consider $\R^{k+1}$. The multiplicative group $\R^*$
acts on $\R^{k+1}$. Consider the quotient topological space $\R^{k+1}/\R^*$.
This is not Hausdorff; the only satured open set of $\R^{k+1}$  containing $0$
is the whole of $\R^{k+1}$ and this intersects any other satured open set.
If we remove $0$, and we restrict the action of $\R^*$ things go better.
Evidently the orbits, i.e. the equivalence classes are in bijective correspondence
with the set of $1$-dimensional linear subspaces of $\R^{k+1}$.
Then one easily verifies that the quotient topological space 
$\PP^k(\R):=(\R^{k+1}\setminus \{0\})/\R^*$
is now Hausdorff and with countable basis. We see also that we get the same
quotient space if we restrict the equivalence relation to the unit sphere $S^k$,
and that the restriction of the projection onto the quotient, $\pi: S^k \to \PP^k(\R)$
is a $2:1$ local {\it homeomorphism}. In fact it is the quotient map by the action on $S^k$
of the group $G$ of order $2$ generated by the antipodal map $x\to -x$.
Then we can endow $\PP^k(\R)$ with a smooth manifold structure as a particular case
of point (1). We can do it also without resctricting to $S^k$.
A finite atlas of $\PP^k(\R)$ is formed by $\{(W_j,\phi_j)\}_{j=1,\dots,k+1}$,
where $W_j$ is the image of the satured open set $\{ x_j \neq 0\}$ of $\R^{k+1}\setminus \{ 0 \}$;
$$\phi_j([x_1,\dots, x_{k+1}])=(x_1/x_j,\dots, x_{j-1}/x_j, x_{j+1}/x_j,\dots, x_{k+1}/x_j)$$
is a homeomorphism of $W_j$ onto $\R^k$. It is immediate to check that the changes
of local coordinates are smooth (actually rational). A posteriori we can define, in a natural
way, a diffeomorphism of
this abstract model of the projective space to the embedded model already constructed.

Every grassmann manifold could be treated in a similar way. First define it
as the quotient topological space of the associated linear Stiefel manifold (which is a open
set in some euclidean space). Prove that this quotient is Hausdorff and with countable
basis and finally give it a (abstract) smooth atlas made by the image of suitable satured
open sets of the Stiefel manifold. A posteriori one can construct a diffeomorphism
onto the already constructed embedded model. 

\smallskip

\begin{example}\label{so(3)} 
{\rm Let us make a few examples.
We are going to establish that $SO(3)\sim \PP^3(\R)$.
An elegant way to see it is by using {\it quaternions}. 
Let $\HH$ be the {\it quaternion} algebra in its matrix form.
That is $\HH$ is the subalgebra of the matrix algebra
$M(2,\C)$ of the matrices of the form
$$A=  \begin{pmatrix}
a+ib& c+id\\
-c+id& a-ib
\end{pmatrix} $$
where $a,b,c,d\in \R$. Then $\HH$ is generated
by the matrix
$$A(i)=\begin{pmatrix}
i& 0\\
0& -i
\end{pmatrix}, \ A(j)= \begin{pmatrix}
0& 1\\
-1 & 0
\end{pmatrix} , \ A(k)=\begin{pmatrix}
0& i\\
i& 0
\end{pmatrix} $$
which verifies the relations
 $$A(i)^2=A(j)^2=A(k)^2=-I$$
$$A(i)A(j)=A(k)=-A(j)A(i), \ A(j)A(k)=A(i)=-A(k)A(j)$$ 
$$\ A(k)A(i)=A(j)=-A(j)A(k) \ . $$
By setting
$$A^*:= \bar A^t$$
we have 
$$(AB)^*=A^* +B^*, \ (AB)^*=A^*B^* \ $$
$$|A|^2:= AA^* = \det A \ $$
and if $A\neq 0$
$$ A^{-1}= \frac{1}{|A|^2}A^* \ . $$
Set
$$\HH_1 = \{ A\in \HH; \ |A|=1\} \ . $$
This is a group with respect to the restriction of the multiplication.
In fact $\HH_1$ is naturally identified with the  special unitary group $SU(2)$
which as a manifold is naturally identified with the
unit sphere $S^3$ in $\R^4$.
Set
$$\HH_0=\{A\in \HH; \  A^*=-A\} \ . $$
which is naturally identified with an euclidean space $\R^3$.
One verifies easily that for every $A\in \HH_1$,
$$ \alpha_A: \HH_0\to \HH_0, \ X\to AXA^{-1}$$
acts as a rotation on $\HH_0=\R^3$.
In fact this gives us a degree 2 covering map
$$ SU(2)\to SO(3), \ A\to \alpha_A$$
such that $\alpha_A=\alpha_B$ if and only if $B=\pm A$.
Hence finally 
$$SO(3)\sim SU(2)_{/\pm I} \sim \PP^3(\R) \  $$
as claimed. 

Let us consider now for every $(P,Q)\in SU(2)\times SU(2)=\HH_1 \times \HH_1$,
the map
$$\alpha_{P,Q}: \HH \to \HH, \ A\to PAQ^{-1} \ $$
by identifying $\HH \sim \R^4$, $\alpha_{P,Q}\in SO(4)$
and $\alpha_{P,Q}=\alpha_{P',Q'}$ if and only if $(P,Q)=\pm (P',Q')$.
Then similarly as above
 we get that
$$ (SU(2)\times SU(2))/\pm 1 \sim SO(4) \ . $$
}
\end{example}

(3) {\it (Grassmann manifolds of oriented spaces)} The set $\tilde G_{m,n}$
of {\it oriented} $n$-subspaces of $\R^m$ can be naturally endowed with
a smooth compact manifold structure $\tilde \GG_{m,n}$ such that the map
$$ p: \tilde G_{m,n}\to G_{m,n}$$
that forgets the orientation becomes a degree $2$ smooth covering map
$$ p: \tilde \GG_{m,n}\to  \GG_{m,n} \ . $$ 
There is a natural tautological bundle
$$\tilde \tau : \Vv(\tilde \GG_{m,n})\to \tilde \GG_{m,n}$$
which in fact equals  $p^*(\tau)$. The fibres of $\tilde \tau$ are tautologically
oriented and this is also the case for every pull-back of $\tilde \tau$.

\smallskip

(4) This example could sound a bit artificial, but it reveals nevertheless some
subtilities. let $M$ be a smooth manifolds (even embedded). Let $f:X\to M$
be any {\it homeomorphism}.  Then 
$$\Uu_f:=\{(f^{-1}(W), \phi \circ f)\}_{(W,\phi)\in \Aa_M}$$
is a {\it smooth} atlas on $X$ so that $f$ becomes tautologically a diffeomorphism.
If $X=M$ (as a topological space), the two smooth structures given by $\Uu_f$
and $\Aa_M$ are diffeomorphic to each other but they are not the same structure (in other
words ${\rm id}_M$ is {\it not} a diffeomorphism). Even if $M$ is embedded
in no natural way the structure given by $\Uu_f$ is embedded.   
 
\medskip

Let us retrace and extend a few notions already developed for embedded
manifolds. The operative principle is:

\medskip {\it Whatever has been built in terms of smooth atlas
can be done as well for abstract smooth manifolds.}
\medskip

\smallskip

{\bf Manifolds with boundary.} We extend the Definition \ref{abs-smooth-man} 
by admitting smooth atlas with charts homeomorphic to open sets
of the half space $\HH^m$. The boundary $\partial M$ is (well) defined
by the same arguments of the embedded case.

\medskip

{\bf Orientable/oriented} manifolds as well as the {\bf oriented boundary}
of an oriented manifold with boundary treated in terms of {\bf oriented atlas}
make sense verbatim also in the abstract case. Also the interpretation
in terms of the deteminant bundle will extend as soon as it shall be
defined (see below).

\medskip

Boundaryless {\bf submanifolds} of a boundaryless manifold are defined in terms of
the existence of atlas made by relatively normal charts. Relatively normal
charts are defined also at the boundary of a manifold with boundary.
As for embedded manifolds, especially if both the manifold and a submanifold have
non empty boundary, there are many possible configurations. Also in the
abstract case, one points out the notion of {\bf proper submanifold} 
(we will return on and use it diffusely later).

\medskip

{\bf Smooth fibred bundles} and related notions introduced in Section 
\ref{emb-fib-bundle} extend words by words by replacing embedded 
with arbitrary smooth manifolds and maps.

\medskip
By using the basis of neighbourhoods defined in terms of representations
of maps in local coordinates (as in (2) of Section \ref {weak-top}), then 
the definition of the {\bf function spaces} $\Ee^r(M,N)$, $\Ee(M,N)$ extends 
without any change  to the abstract case.

\medskip
{\bf Homotopy, isotopy, diffeotopy} and the {\bf homogeneity} property
(see Section \ref{H-I-D}) extend as well.

\section{The (abstract) tangent functor}\label{abs-tangent-fuct}
Probably this is the most demanding extension by dealing with abstract
smooth manifolds. In the case of embedded manifolds
tangent bundles and maps imposed by themselves, starting from the basic 
ones for open sets in some euclidean spaces. In the abstract case 
 they must be somehow ``invented'',
with the constraint to agree with already done in the embedded category.  
This also will bring us to a general notion of
{\it fibre bundle in the sense of Steenrod} \cite{Steen}.
\smallskip

{\bf Construction of the tangent bundle.}
Let $M$ be a $m$-smooth manifold with its maximal smooth atlas $\Aa=\{(W_j,\phi_j)\}_{j\in J}$.
For every $(i,j)\in J^2$, define the map
$$ \mu_{ji}: W_i\cap W_j \to {\rm GL}(m,\R), \ \mu_{ji}(x)= d_{\phi_i(x)} (\phi_j \circ \phi_i^{-1})\ . $$
This family of maps $\{\mu_{ji}\}_{(i,j)\in J^2}$ verifies the following properties:
\begin{enumerate}
\item Every $\mu_{ji}$ is smooth.
\item For every $j\in J$, for every $x\in W_j\cap W_j=W_j$,
$$\mu_{jj}(x)=I_m \ . $$ 
\item For every $(j,i)\in J^2$, for every $x\in W_i\cap W_j=W_j\cap W_i$,
$$\mu_{ji}(x)=\mu_{ij}(x)^{-1} \ . $$
\item For every $(i,j,k)\in J^3$, for every $x\in W_i\cap W_j\cap W_k$
$$\mu_{ik}(x)\mu_{kj}(x)\mu_{ji}(x)=I_m \ . $$
\end{enumerate}

\smallskip

We summarize these properties by saying that 
\smallskip

{\it $\{\mu_{j,i}\}$ is a
smooth cocycle on the open covering $\Aa$ with values in the Lie
group {\rm GL}$(m,\R)$}.

\smallskip

Note that as GL$(m,\R)$ is {\it non commutative} (if $m>1$), then the order in property 4
is not negligible.

Let us consider now the topological product $M\times \R^m \times J$, where $J$ is endowed with the discrete topology.
Let $\Tt$ be the subspace made by the triples $(x,v,j)$ such that $x\in W_j$. Hence $\Tt$ is the {\it disjoint} union
of the open sets $W_j\times \R^m \times \{j\}$, $j\in J$, each one being canonically homeomorphic to $W_j\times \R^m$.
Now let us put on $\Tt$ the relation $(x,v,j)\sim (x',v',k)$ if and only if $x=x'$ and $v'=\mu_{kj}(x)v$.
The cocycle properties 2--4 ensure that it is an equivalence relation. We set 
$$T(M):= \Tt/\sim $$ 
the topological quotient space and denote by $q:\Tt \to T(M)$ the canonical continuous
projection. We have the well defined surjective map
$$ \pi_M: T(M)\to M, \  \pi_M([x,v,j])=x $$
which is continuous. In fact for every open set $A$ of $M$, $(\pi_M\circ q)^{-1}(A)$ is the intersection of
$\Tt$ with $A\times \R^m \times J$, hence it is a satured open set, with open image in $T(M)$.
It is a topological exercise to show that $T(M)$ {\it is Hausdorff and with countable basis}, this is left to the reader.
\smallskip

{\it (Local trivializations)} For every $j\in J$, set
$$\Psi_j: W_j\times \R^m \to T(M), \ (x,v)\to q(x,v,j)=[(x,v,j)] \ . $$
One verifies that
\begin{enumerate}
\item $\Psi_j$ is continuous (because $q$ is continuous);
\item $\Psi_j$ takes values in $\pi_M^{-1}(W_j)$ and $\pi_M\circ \Psi_j=p_j$,
where $p_j:W_j\times \R^m \to W_j$ is the projection.
\item In fact $\Psi_j$ {\it is a homeomorphism onto} $\pi_M^{-1}(W_j)$.
For if $b=[x,v,k]\in \pi_M^{-1}(W_j)$ , then $b= \Psi_j(x,\mu_{jk}(x)v)$, hence
$\Psi_j$ is onto. If $[x,v,j]=[x',v',j]$, then $x=x'$ and $v=v'$ because  $\mu_{jj}=I_m$.
Hence $\Psi_j$ is injective. Finally, to show that the inverse of $\Psi_j$ is continuous, 
it is enough to show that if $A$ is open in $W_j\times \R^m$, the $q^{-1}(\Psi_j(A))$ is open in
$\Tt$. Since the $W_k\times \R^m \times  \{k\}$'s form a open covering of $\Tt$, it is enough to prove that  every
$q^{-1}(\Psi_j(A))\cap (W_k\times \R^m \times \{k\})$ is open. 
This intersection is contained in the open set $(W_j\cap W_k)\times \R^m \times \{k\}$ of $\Tt$.
On this open set $q=\Psi_j\circ r$, where $r(x,v,k)=(x, \mu_{jk}(x)v) $ which is continuous; the thesis follows.
\end{enumerate}

\smallskip

{\it (Changes of local trivializations)} These are of the form
$$ \Psi_j^{-1}\circ \Psi_i(x,v)=(x,\mu_{ji}(x)v)$$
defined on $(W_j\cap W_i)\times \R^m$ to itself. Clearly they are smooth, and
pointwise linear in the second argument.
So we have proved that
$$\pi_M: T(M)\to M$$
is a (abstract) smooth vector bundle over $M$ with fibre $\R^m$, called the {\it tangent bundle of $M$}.
For every $p\in M$, the fibre $T_pM:= \pi_M^{-1}(p)$ is {\it by definition} the {\it tangent space
of $M$ at $p$}.
It is clear that $T(M)$ is a smooth manifold because it is locally
diffeomorphic to spaces of the form $W_j\times \R^m$, $W_j$ being a open set in the smooth
manifold $M$. To be even more concrete, we can exhibit the following special smooth atlas
of $T(M)$ made of fibred maps: 
$$T\Aa=\{\pi_M^{-1}(W_j), \Phi_j)\}_{j\in J}$$ 
where $\Phi_j:= (\phi_j,{\rm id})\circ \Psi^{-1}_j$, and
$$(\phi_j,{\rm id}): W_j\times \R^m \to U_j\times \R^m\subset \R^m\times \R^m, \ (x,v)\to (\phi_j(x),v) \ . $$
The changes of local coordinates are of the form
$$\Phi_j\circ \Phi_i^{-1}(x,v)= (\phi_j\circ \phi_i^{-1}(x), \mu_{ji}(x)v)$$
which {\it ultimately is nothing else than the tangent maps of the change of coordinates on} $M$.

\smallskip

{\bf Tangent maps.} Let $f:M\to M'$ be a smooth map between smooth manifolds.
We want to define now the tangent map 
$$Tf: T(M)\to T(M')$$ in such a way that $[f,Tf]$
is a vector bundle fibred map. We have constructed the tangent
bundles by patching together the product pieces. We do similarly for $Tf$.  
Precisely, let $(\pi_M^{-1}(W),\Phi)$, $(\pi_{M'}^{-1}(W'),\Phi')$
be fibred charts of $T(M)$ and $T(M')$  which dominate charts 
$(W,\phi)$, $(W',\phi')$ of $M$ and $M'$ respectively. Assume also
that this system of charts gives us a representation in local coordinates of $f$, 
$\hat f = \phi'\circ f \circ \phi^{-1}$.
Then we {\it locally} define 
$$Tf_{W,W'}: \pi_M^{-1}(W)\to \pi_{M'}^{-1}(W'), \ Tf_{W,W'}= \Phi' \circ T\hat f \circ \Phi^{-1} \ . $$
Recalling the equivalence relation that we have used to build the tangent bundles, one readily
checks that these locally defined $Tf$'s are in fact {\it representations in local
(fibred)  coordinates  of a {\bf globally} defined fibred map } $Tf: T(M)\to T(M')$.
For every $p\in M$ the restriction say $d_pf$ of $Tf$ to $T_pM$, is a linear map
$$d_pf: T_pM\to T_{f(p)}M'$$
which  by definition is the {\it differential of $f$ at $p$}.

\medskip

{\bf Tangent functor.} The basic functorial properties of the chain rule globalize, so that we have:

\medskip{\it The tangent category of the category of smooth manifolds has as objects the 
tangent vector bundles of smooth manifolds and as morphisms the tangent maps of 
smooth maps between embedded  manifolds. This is a subcategory of smooth vector bundles 
over smooth manifolds. Then
$$M \ \Rightarrow  \ \pi_M :T(M)\to M, \  f:M\to M' \ \Rightarrow  [f,Tf]$$
define a covariant functor from the category of embedded smooth manifolds to its tangent category. 
This extends the  embedded tangent functor.}

\medskip

{\bf Immersions and embeddings.} As we dispose now of the differentials $d_pf$
for every smooth map, the notions of immersion and embedding extend
as well as  the related results of section \ref{imm-summ-emb}.

\section{Principal and associated bundles with given structural group}\label{principal-associated}
The construction of the tangent bundles is suited to a wide generalization.

Let $G$ be a Lie group (such as GL$(m,\R)$, $O(m)$, $SO(n)$, $U(n)$, $\dots$,). 
Assume that it acts on a smooth manifold
$F$. This means that there is a goup homomorphism (also called a {\it representation})
$$\rho: G \to {\rm Aut}(F) \ ; $$ 
the associate action is 
$$ G\times F \to F, \ (g,x)\to \rho(g)(x)$$
and sometimes one simple writes $gx$ instead of $\rho(g)(x)$. Sometimes one also
requires that $\rho$ is injective so that $G$ is confused with its image in
Aut$(F)$ and considered as a {\it group of transformations of $F$} (but this is
not strictly necessary).

\begin{remark}{\rm $G$ acts on itself by the injective homorphism $g\to L_g$ (i.e. by left multiplication)
$$ G\times G \to G, \ (g,h)\to L_g(h):= gh \ . $$}
\end{remark}

\medskip

Let $M$ be a smooth manifold and $\Uu= \{A_s\}_{s\in \Ii}$ be a open covering
of $M$. A {\it principal cocycle} on $\Uu$ with values in the {\it structural group}
$G$ is a family of smooth
maps
$$ \cG=\{c_{ts}: A_s\cap A_t \to G\}_{(s,t)\in \Ii^2}$$
such that

\begin{enumerate}
\item For every $s\in \Ii$, for every $x\in A_s$,
$$c_{ss}(x)= 1\in G \ . $$ 
\item For every $(s,t)\in \Ii^2$, for every $x\in A_s\cap A_t$,
$$c_{st}(x)= c_{ts}(x)^{-1} \ . $$
\item For every $(s,t,r)\in \Ii^3$, for every $x\in A_s\cap A_t\cap A_r$
$$c_{sr}(x)c_{rt}(x)c_{ts}(x)= 1 \ . $$
\end{enumerate}
\smallskip

For every representation $\rho: G \to {\rm Aut}(F)$ as above, we have an {\it associated
cocycle} with values in Aut$(F)$
$$ \{\rho_{ts}:= \rho \circ c_{ts}: A_s\cap A_t \to {\rm Aut}(F)  \}_{(s,t)\in \Ii^2}$$ 
which verifies the same properties 1-3 (by replacing $1\in G$ with $1\in {\rm Aut}(F)$).

\smallskip

Then we can repeat words by words the above construction of the tangent bundles
and get a {\it smooth fibres bundle over $M$ with structural group $G$ and fibre $F$}.
So we have a wide family of bundles which share the basic cocycle $\cG$.
When $F=G$ and $G$ acts as above by left multiplication, we get the
{\it principal bundle} of this family; all the other bundles are said 
{\it associated} to such a principal bundle.

 \subsection{Equivalent cocycles}\label{equivalent-cocycle}
The strict equivalence of fibre bundles can be rephrased in terms
of the defining cocycles. Assume that two cocycles $c$ and $c'$ with values in $G$ are 
defined on the same open covering  $\Uu= \{A_s\}_{s\in \Ii}$ of $M$.
Then they define strictly equivalent bundles if and only if there is a family
of maps
$$\{\lambda_s: A_s \to G\}_{s\in \Ii}$$
such that for every $(s,t)$, for every $x\in A_s\cap A_t $,
$$ c'_{ts}(x)= \lambda_s(x)c_{ts}(x)\lambda_t(x)^{-1} \ . $$

\subsection{Tensor bundles}\label{abs-tensor} We can apply this machinery to construct the
abstract version of the tensor bundle relatives to the tangent bundle.

In Section \ref {other-functor}, for every $(p,q)$, we have defined the representation 
$$ \rho_{p,q}:{\rm GL}(m,\R)\to {\rm GL}(T^p_q(\R^m))\sim {\rm GL}(m^{pq}, \R)$$
which is an explicit  rational regular map. By using it we get the tensor
bundle
$$\pi_{p,q}: T^p_q(M)\to M \ . $$

\medskip

The representation
$$ \det: {\rm GL}(m,\R)\to \R^*$$
leads to the {\it determinat bundle} of $M$ 
\medskip

The principal bundle of this family is the {\it frame bundle} of $M$,
once we have identified the columns of any non singular matrix 
with a basis of $\R^m$.

\medskip

{\bf Tensors fields. } The contents of Sections \ref{tensor-field} and  \ref{orientable} extend verbatim.

 \section{Embedding  abstract compact manifolds}\label{emb-abs-comp} 
Let $M$ be a compact smooth $m$-manifold possibly with boundary $\partial M$.
The notion of nice atlas makes sense in full generality. We have:

\begin{proposition}\label{emb-abs-comp} (1) Let $M$ be a compact smooth manifold.
Then there is a diffeomorphism $f: M\to M'$ onto an embedded manifold $M' \subset \R^{h}$,
for some $h$. 
\smallskip

(2) The tanget map $Tf$ establishes a vector bundle equivalence
between the respective tangent bundles of $M$ and $M'$. This equivalence propagates to
all tensor bundles and to the frame bundle.
 \end{proposition}

\Dim   (1): we argue as in the proof of Proposition   
\ref{proper-emb}, by using a nice atlas of $M$ $\{(W_j,\phi_j)\}_{j=1,\dots, s}$
including also relative normal charts along $\partial M$, instead of a nice atlas with collar. This allows to define
the embedding 
$$\beta=(\beta_1,\dots, \beta_s): M\to (\R^m\times \R)^s$$
$$\beta_j=(\lambda_j\phi_j, \lambda_j) \ . $$
The verification is the same of Proposition \ref{proper-emb}.
\smallskip

Point (2) follows from the fact that the abstract functor extends
the embedded one.

\cvd

\smallskip

By combining the last proposition with Proposition \ref{2n}
we have:

\begin{corollary}\label{Weak -Whitney} {\rm (Weak Whitney immersion/embedding theorem)}
Every $m$-dimensional compact smooth $m$-manifold
$M$ can be immersed in $\R^{2m}$ and can be embedded in $\R^{2m+1}$.
\end{corollary}
\Dim If $M$ is boundaryless it is an immediate corollary of Propositions \ref{emb-abs-comp} and \ref{2n}. If $M$ has boundary
we can reduce to the boundaryless case by using the double of $M$.

\cvd

\smallskip

So, up to diffeomorphism we can assume that every compact 
manifold $M$ is embedded. We extend now this result to
every abstract vector bundle over $M$, besides the tangent and tensor
bundles.

\begin{proposition}\label{emd-abs-VB} 
Every abstract vector bundle $\xi$ over an embedded
compact smooth manifold $M\subset \R^h$  is strictly
equivalent to an embedded vector bundle.
\end{proposition}
\Dim By compactness we can assume that the abstract bundle $p: E \to M$ is determined by
a cocycle $c_{ts}$ over a nice atlas $\Uu=\{(W_j,\phi_j)\}_{j=1,\dots,s}$ of $M$. Consider the family
of local trivializations $\Phi_j: p^{-1}{|W_j} \to W_j \times \R^n$, and let $\{\lambda_j\}$ be the partition of unity
over $\Uu$ as usual. For every $j$, denote by $q_j: W_j \times \R^n \to \R^n$ the natural projection.
Finally define
$$h: E \to M\times \R^{ns}\subset \R^{h+ns}, \ h(e)=(p(e), \lambda_1(p(e))q_1(e), \dots, \lambda_s(p(e))q_s(e)) \ .  $$
The restriction of $h$ to $M\subset E$ as the zero section is equal to the identity. Moreover every fibre of the bundle
is linearly embedded onto a $n$-subspace of $\R^{ns}$.

\cvd 

\medskip

\begin{remark}\label{frame+sphere}{\rm (1) The conclusion of Proposition \ref{emd-abs-VB}
holds as well for the frame bundle and more generally
for any abstract smooth fibre bundle over $M$ with embedded fibre.}
\end{remark}

\subsection {On vector bundles on sphere again}\label{more-vect-sphere}
Now we can complete the classification of vector bundles on the spheres stated in Section \ref {VB-sphere} . 
By combining those constructions  with the present ones, 
every 
map  $\rho_\xi : S^{m-1} \to  {\rm GL}^+(k,\R)$ extends to a cocycle
$\rho_\xi : D^+ \cap D^- \to {\rm GL}^+(k, \R)$ on the nice covering of the sphere formed by the two
smooth disks $D^+$, $D^-$. So the claimed inverse map  $\rho^{-1}$ is obtained by taking the strict equivalence
class of the embedded vector bundle over $S^{m-1}$ constructed as in Proposition \ref {emd-abs-VB} 
by using this cocycle. 

\subsection{On tubular neighbourhoods and collars again}\label{more-tube}
In Section \ref{tubular} we have constructed tubular neighbourhoods and collars
{\it unique up to isotopy} starting from an {\it embedded} compact manifold $M\subset \R^k$.
Above we have shown that every (abstract) compact manifold $M$ can be embedded
in some $\R^k$ and, {\it a priori}, that family of tubular neighbourhoods and collars, considered up to isotopy,
could depend on the embedding. However this is not the case. First every embedding $M\subset \R^k$
can be ``stabilized'' to $M\subset \R^k \subset \R^{k+h}$; moreover, by using the results of the present section
with  Proposition \ref{2n}, if $h$ is big enough, up to isotopy two  embeddings of $M$ in $\R^{k+h}$ have 
disjoint images and can be extended to an embedding of $M\times [0,1]$ so that they are isotopic to each other. 

\medskip
 
Summarizing: 
\smallskip

{\it By considering compact smooth manifolds up to diffeomorphism, we can exploit
all the results already obtained in Chapter \ref{TD-COMP-EMB} for embedded compact manifolds}.

\section{On complex  manifolds}\label{complex-man}
Another reason to introduce the abstract notion of manifold in terms of atlas with change of coordinates in
a determined class of homeomorphism (for instance smooth diffeomorphisms in our favourite setting)
is that it is suited to several interesting implementations. Abstract {\it complex $n$-manifolds} 
have as local models the open sets in $\C^n$ and change of coordinates that are complex analytic
(i.e. holomorphic) diffeomorphisms (biholomorphisms). Holomorphic maps beetween complex manifolds
are defined in terms of holomorphic local representations;  and so on, by following and specializing 
several constructions developed above (complex tangent bundle, complex submanifolds etc.). 
On the other hand, by the {\it maximum principle}, the constant functions $c: M \to \C$ are the only holomorphic functions defined on any 
{\it compact} connected complex manifold $M$. So compact complex manifolds {\it cannot be embedded
into any $\C^m$} (as complex submanifolds). This is a main difference with respect to our favourite real smooth theory. Moreover,
bumb functions do not exist in the complex setting, so the many constructions which have employed such a tool
cannot be performed on complex manifolds. Although we have introduced them as examples of embedded
smooth manifolds, {\it complex Stiefel and Grassmann manifolds} (in particular the complex projective spaces)
can be naturally endowed with an (abstract) structure of compact complex manifold.

By identifying $\C^n\sim \R^{2n}$ and considering holomorphic maps as a special kind of smooth maps,
by forgetting the complex structure, every complex $n$-manifolds $M$ can be considered as 
a smooth $2n$-manifold (as we have done for the complex Grassmannian); moreover the complex structure
induces on this $2n$-manifold a {\it natural orientation}. Especially in dimension $4$, $2$-complex manifolds
(also called {\it complex surfaces}) form an important class of oriented $4$-manifolds. 

\medskip

{\bf (The Riemann sphere)} As a basic example let us consider $\PP^1(\C)$; let us identify $\R^2 \sim \C$ and consider the two-charts
atlas of the $2$-sphere $S^2$ given by the stereographic projections from the two poles. These
can be considered as $\C$-valued charts. In order to make it a complex-manifold atlas it is enough
to compose the second projection with the complex conjugation $z\to \bar z$. Moreover it is 
immediate to identify such an atlas with the standard two-charts complex atlas of $\PP^1(\C)$.
This show in particular that $\PP^1(\C)$ is diffeomorphic to $S^2$; this last considered as a $1$-dimensional
complex manifold is called the {\it Riemann sphere}.

\chapter{Cut and paste compact manifolds}\label{TD-CUT-PASTE}
In this Chapter we  deal with compact manifolds or more generally
with possibly non compact manifolds which nevertheless can be embedded 
in some $\R^n$ being a closed subset too. Thus we can exploit the results of Chapter \ref{TD-COMP-EMB}.

\section{Extension of isotopies to diffeotopies}\label{ext-isotop}
We recall a few notions.  

Let $N$ be a smooth boundaryless $n$-manifold.
Let $M$ be a smooth $m$-manifold and 
$$F: M\times [0,1]\to N$$
a smooth map such that $f_t$ is an embedding for every $t\in [0,1]$; then $F$ is an {\it isotopy}
connecting $f_0$ and $f_1$.  

A  {\it diffeotopy} of $N$ (also called an {\it ambient isotopy})
is a smooth map
$$G: N\times [0,1] \to N$$
such that $g_t$ is a diffeomorphism for every $t\in [0,1]$. 
We will also assume that $g_0={\rm id}_N$.
Hence diffeotopies are  special  isotopies. 

\begin{definition}\label{iso-extend}{\rm We say that an isotopy $F$ as above {\it extends to an ambient
isotopy} if there is a diffeotopy $G$ of $N$ such that $f_t=g_t\circ f_0$ for every
$t\in [0,1]$. Note that $\{V_t=f_t(M)\}$ is a one parameter family of submanifolds of $N$
(each diffeomorphic to $M$), and $V_t=g_t (V_0)$, for every $t$.}
\end{definition}

\smallskip

We are going to see that under mild compactness assumptions, isotopies actually extend to
diffeotopies. This will be a key result to show that several cut-and-paste procedures below
are well defined.  To this aim, it is useful to recast diffeotopies as flow
of (suitable) vector fields. In doing it we will tacitly incorporate basic facts
about the existence, uniqueness and regular dependence on the initial data
of the solutions of {\it ordinary differential equations} (see for instance \cite{A}).
\smallskip
 
For every isotopy $F$ as above, its {\it track} is the map defined as
$$\hat F: M\times [0,1] \to N\times [0,1], \ \hat F(x,t):= (f_t(x),t) \ . $$
The {\it support}  of $F$ is the closure in $M$ of the set
$$\{x\in M| \ \exists t\in [0,1], \ f_t(x)\neq f_0(x)\} \ . $$

Given an ambient isotopy $G$ of $N$, as above,
and its track $\hat G$ (which is  a level preserving diffeomorphism),  
consider on $N\times [0,1]$ the constant ``vertical'' tangent
vector field $V$ defined by
$$ V(x,t)=(0,1)\in T_xN \times \R \ . $$
The tangent map $T\hat G$ transforms this field into another tangent vector field on
$N\times [0,1]$ of the form 
$$X_G(x,t)=(v_G(x,t),1) \ . $$ 
Then the map $\hat G$ transforms every  vertical integral
line $j_x: [0,1] \to N\times [0,1]$ of $V$ such that $j_x(0)=(x,0)$ , into the integral line 
$\hat j_x: [0,1]\to N \times [0,1]$ of the field $X_G$ such that $\hat j_x(0)=(x,0)$ . In fact, by 
construction 
$$\hat G(j_x(t))=\hat j_x(t)=(g_t(x),t) $$
that is $G$ is the {\it flow} of $X_G$, with initial data at $N\times \{0\}$. 
Hence we can reconstruct the diffeotopy $G$ by integration of
the field $X_G$. 

On the other hand, if $v(x,t)$, $t\in[0,1]$, is any {\it time depending
smooth tangent vector field on $N$}, let $X(x,t)=(v(x,t),1)$ be the corresponding
field on $N\times [0,1]$. Let us say that it has {\it complete integral lines} 
if for every initial point $(x,0)\in N\times [0,1]$, the corresponding integral line of $X$ is defined
on the whole interval $[0,1]$. If $X$ has complete integral lines then it {\it generates
a diffeotopy of $N$}, that is there is a unique diffeotopy $G=G_X$ such that $X=X_G$.
This establishes a bijection between diffeotopies and such tangent vector fields $X$ with complete
integral lines.

If $N$ is not compact, not every $X$ has complete integral lines;
by local existence and uniqueness, in general, for every $(x,0)$ there is a maximal open interval
$[0,t_x)\subset [0,1]$ on which the corresponding integral line is defined.
However, if we assume that $v(x,t)$  has {\it compact support}, 
then it is not hard to show that $X$
actually has complete integral lines, and the generated diffeotopy $G_X$ has compact
support. Recall that the support of $v(x,t)$ is defined as the closure in $N$ of
the set
$$\{x\in N| \ \exists t\in [0,1], \ v(x,t)\neq 0\} \ . $$
{\it Viceversa}, if a diffeotopy $G$ has compact support, then also
$v_G$ has compact support. This restricts the above bijection to diffeotopies and
such tangent vector fields with compact support, and
gives us a very flexible way to construct diffeotopies, under mild compactness assumptions.
Finally we can state and prove our extension theorem (sometimes known as
``Thom's lemma'').

\begin{proposition}\label{Thom_lemma} Let $F:M\times [0,1] \to N$ be
an isotopy of embeddings of the {\rm compact} boundaryless smooth $m$-manifold
$M$ into the boundaryless $n$-manifold $N$. Then $F$ extends to
an ambient isotopy of $N$ with compact support.
\end{proposition}
\Dim Consider the track $\hat F$ of the isotopy $F$. It is   
a level preserving embedding of $M\times [0,1]$ onto a compact proper submanifold
say $\hat M$ of  $N\times [0,1]$. 
Consider the constant vertical tangent vector field on $M\times [0,1]$
$$ V_M(x,t)=(0,1)\in T_xM\times \R \ . $$
The tangent map $T\hat F$ sends $V_M$ to a vector field $X_M$ of the form
$$X_M(y,t)=(v_M(y,t),1), y=f_t(x)$$
defined along $\hat M$.
Then the  idea is to extend $X_M$ to a tangent vector field $X$ of the form
$$X(y,t)=(v(x,t),1)$$
defined on the whole of $N\times [0,1]$ and such that $v(y,t)$ has compact
support. The ambient isotopy $G_X$ generated by the field $X$ will
eventually extend the isotopy $F$. Clearly this extension task only concerns the
``horizontal'' part $v_M$. 
Under the assumption made at the beginning of this
section, we know from Chapter \ref{TD-COMP-EMB} that there are a proper compact tubular neigbourhood $U$
of $\hat M$ in $N\times [0,1]$ (which restricts to a tubular neighbourhood
of $f_t(M)$ in $N\times \{t\}$ for every $t\in [0,1]$), 
and a compact submanifold with boundary $W$ of $N$ such that
$U$ is contained in ${\rm Int}(W)\times [0,1]$. 
By using the local product structure of $U$ along $\hat M$, we can cover $\hat M$ by a finite number
of smooth closed $(n+1)$ balls, each one say $B$ easily supporting a smooth extension $v_B$ of the restriction of
$v_M$ to $B\cap \hat M$, and such that their union is contained in  $U$.
Such $B$'s can be incorporated in a nice covering with collar of
$W\times [0,1]$, say $\Uu$. Locally extend $v_M$ on any open set of such a covering different from the $B$'s
by setting it constantly equal to $0$. By using a partition of unity supported by $\Uu$, we finally get
the required smooth extension of $v_M$     
to a smooth time depending field $v$ defined on the whole of $N$, constantly equal to zero on the complement
of $W$, and with compact support contained in $W$.

\cvd

\begin{remarks}\label{remark_iso_diffeo}{\rm

(1) For the sake of simplicity, we have proved Thom's lemma under the assumption that both the compact manifold $M$ and $N$ are
boundaryless.  Mild adaptations of the same construction
allow to extend the results under more general hypotheses.
Assuming that both $M$ and $N$ possibly have boundary, we can cover
for instance the following situations, getting a pertinent version of Thom's lemma (details are left to the readers):  

(a) $F$ is an isotopy of embeddings of $M$ either in $N\setminus \partial N$ or
in $\partial N$.

(b) $F$ is an isotopy of proper embeddings of
$(M,\partial M)$ in $(N,\partial N)$.

(c) Every boundary component of $\partial M$ is embedded by every $f_t$ either in $N\setminus \partial N$
or in $\partial N$, being $f_t(M)$ transverse to $\partial N$ along $f_t(M)$; for instance, this includes
the case when for every $t$, $f_t$ parametrizes a collar of a compact boundary component of $\partial N$. 

(d) For every $t\in [0,1]$, $f_t$ parametrizes a relative tubular neighbourhood of a compact proper
submanifold $(Y,\partial Y)$ of $(N,\partial N)$. 

\medskip

(2) If $M$ is not compact, in general an isotopy of embeddings of $M$ in $N$ does not
extend to any diffeotopy. For example, take $M=\R$ and $N=\R^2$, then it is easy
to construct an isotopy of embeddings connecting $f_0$ being the natural inclusion $\R_x\subset \R^2_{x,y}$
with  $f_1$ having as image the set $\{(x,y); \ x^2+(y-1)^2=1, \ (x,y)\neq (0,2)\}$. For basic topological
reasons it cannot be extended. On the other hand, what is really important to achieve the proof of
Thom's lemma is that the isotopy $F$ has compact support, even if $M$ is possibly non compact. 
}
\end{remarks} 

As a corollary, we have also the following sort of relative extension result.
\begin{corollary}\label{rel-TL}  Let $Y$ be a compact submanifold of the
manifold $M$. Let $F$ be an isotopy of embeddings of $Y$ into the 
manifold $N$ such that a version of Thom's lemma holds. Assume that $f_0$ can be extended to an 
embedding $h_0:M\to N$.
Then also $f_1$ can be extended to an embedding $h_1:M\to N$; 
moreover we can require that
$h_0$ and $h_1$ are diffeotopic to each other.
\end{corollary}
\Dim By Thom's lemma 
$F$ extends to a diffeotopy $G$ of $N$, hence $h_1:=g_1\circ h_0$
is an embedding of $M$ in $N$ which extends $f_1$ and is diffeotopic to $h_0$
by construction.

\cvd 

\section{Gluing manifolds together along boundary components}\label{partial-glueing}
Let $M_1$ and $M_2$ be $m$-compact manifolds with boundary, $V_1$ and $V_2$ be
unions of connected components of $\partial M_1$ and $\partial M_2$ respectively, and
$\rho: V_1\to V_2$ be a diffeomorphism. Consider the compact topological quotient space
$$ M_1\amalg_\rho M_2 $$
by the equivalence relation on the disjoint union $M_1\amalg M_2$
which identifies every $x\in V_1$ with $\rho(x)\in V_2$; $\rho$ is called the
{\it gluing map}. Denote by
$$q: M_1 \amalg M_2 \to M_1\amalg_\rho M_2$$ 
the projection onto the quotient space, for $s=1,2$, 
$$i_s: M_s\to M_1 \amalg M_2$$
the inclusion, and finally set 
$$j_s=q\circ i_s \ . $$ It is clear that $j_s$ is a homeomorphism onto its image. 
We have: 

\begin{proposition}\label{gluing-exists}  The quotient space  $M_1\amalg_\rho M_2$ can be endowed with a structure
of smooth $m$-manifold such that for every $s=1,2$, $j_s$ is a smooth embedding, and
$$\partial (M_1\amalg_\rho M_2) = (\partial M_1 \amalg \partial M_2)\setminus (V_0\amalg V_1)\ . $$
\end{proposition}
\Dim Fix a collar $c_1:[-1,0]\times V_1 \to M_1$ of $V_1$ in $M_1$ and a collar $c_2: V_2 \times [0,1] \to M_2$
of $V_2$ in $M_2$. Define 
$\psi_V: (-1,1)\times V_1 \to M_1\amalg_\rho  M_2$ by 
$$\psi_V(t,x)= j_1(c_1(t,x)) \ {\rm if} \ t\in (-1,0], \   \psi_V(t,x)= j_2(c_2(\rho(x),t)) \ {\rm if} \ t\in [0,1) \ . $$ 
It is clear that $\psi_V$ is a homeomorphism onto an
open neighbourhood $U$ of 
$$V:=j_1(V_1)=j_2(V_2)$$ in $M_1\amalg_\rho M_2$. By composing the charts of a
smooth atlas of $(-1,1)\times V_1$ with $\phi_V=\psi^{-1}_V$ we get a smooth atlas say $\Uu_V$ on $U$
such that $\psi_V$ becomes tautologically a diffeomorphism. Similarly, let $\Uu_s$ be a smooth atlas
on $j_s(M_s \setminus V_s)$ such that the restriction of $j_s$ to $M_s\setminus V_s$ is tautologically
a diffeomorphism. It is immediate to check that $\Uu_V \cup \Uu_1 \cup \Uu_2$ is a smooth atlas
on $M_1\amalg_\rho M_2$ that determines a smooth manifold structure with the required properties.
An equivalent way to get such a smooth structure on  $M_1\amalg_\rho M_2$ 
is as follows: take the disjoint union $(M_1 \setminus V_1)\amalg (M_2\setminus V_2)$
and identify the two open sets $c_1((-1,0)\times V_1)$ and $c_2((0,1) \times V_2)$ by identifying
$(t,x)\in (-1,0)\times V_1$ with $(1-t,\rho(x))\in (0,-1)\times V_2$.

\cvd

\medskip

Let us say that a smooth structure on $M_1\amalg_\rho M_2$ obtained so far is given by
{\it gluing $M_1$ and $M_2$ together by means of the gluing map $\rho$}.   
Such a smooth structure  depends on the choice
of  collars entering the construction. However we have the following 
{\it uniqueness up to diffeomorphism}. Precisely:
\begin{proposition}\label{gluing-unique} 
Any two smooth structures given by gluing $M_1$ and $M_2$ together
via the gluing map $\rho$ are diffeomorphic to each other, via a diffeomorphism
which is the identity at the boundary.
\end{proposition}
\Dim Assume for simplicity that two implementations of the construction
differ by the choice of two different
collars $c_2,c'_2: V_2 \times [0,1] \to M_2$. Denote by $M$ and $M'$
the respective smooth structures on $M_1\amalg_\rho M_2$.
The isotopy (relative to $V_2$) of the two collars of $V_2$ in $M_2$
extends to a diffeotopy $G$ of $M_2$. Then the map $h: M \to M'$
such that $h={\rm id}_{j_1(M_1)}$ on $j_1(M_1)$, $h=g_1\circ (j_2)^{-1}$
on $j_2(M_2)$ provides a required diffeomorphism. The general case 
is achieved by a similar argument.

\cvd 

\medskip

Hence it makes sense to denote by $M_1\amalg_\rho M_2$ such
a diffeomorphism class of smooth manifolds obtained by gluing $M_1$ and $M_2$
together. In fact {\it we will often do the abuse to confuse such a class with any
representative}.
\smallskip

In some cases we can deduce that $M_1\amalg_\rho M_2$ and 
$M_1\amalg_{\rho'} M_3$ are diffeomorphic, where 
$\rho: V_1\to V_2$, $\rho': V_1\to V_3$ are respective gluing maps. 

\begin{proposition} \label{same-gluing} (1) If the diffeomorphism 
$\rho' \circ \rho^{-1}: V_2 \to V_3$ extends to a diffeomorphism 
$h: M_2 \to M_3$. Then $M_1\amalg_\rho M_2$ and $M_1\amalg_{\rho'} M_3$
are diffeomorphic.

(2) If two gluing maps $\rho_0, \rho_1: V_1 \to V_2$ are isotopic, then the 
manifolds obtained by gluing $M_1$ and $M_2$ together by means of
$\rho_0$ and $\rho_1$ respectively are diffeomorphic to each other.
\end{proposition}

\Dim 
A collar of $V_3$ in $M_3$,
used to define a smooth structure of $M_1\amalg_{\rho'}M_3$, can be lifted by $h$
to a collar of $V_2$ in $M_2$; this can be used to define a smooth structure of
$M_1\amalg_\rho M_2$ which by construction is diffeomorphic to $M_1\amalg_{\rho'}M_3$.
This achieves (1).

 As for $(2)$,  $\rho_1\circ \rho_0^{-1}$ is diffeotopic to the identity of $V_2$ which
 obviously extends to the identity of the whole $M_2$. By Corollary \ref{rel-TL},
 then also   $\rho_1\circ \rho_0^{-1}$ extends to a diffeomorphism of $M_2$ and
 we can apply previous item (1).
 
 \cvd

\medskip

{\bf Oriented version.} Keeping the above setting, assume furthermore that $M_s$ is oriented and that $V_s$ is part
of the {\it oriented boundary} $\partial M_s$. If $\rho: V_1 \to V_2$ is a {\it orientation reversing} diffeomorphism then
$M_1\amalg_\rho M_2$ is endowed with a structure of {\it oriented} smooth $m$-manifold such that 
$j_1$ and $j_2$ are {\it orientation preserving embeddings}. Up to orientation preseving diffeomorphism it is well
defined the oriented manifold $M_1\amalg_{\rho} M_2$ which actually only depends on the isotopy class of the orientation
reversing attaching diffeomorphism $\rho$.

\section{On corner smoothing}\label{corner-smooth} Of course the notion of smooth manifold with corners (extending
Definition \ref{man-with-corner}) makes sense in the abstract setting. Making use of tubular neighbourhoods and collars 
as in the previous section, it is not hard to see that every compact smooth $m$-manifold with corner $M$ 
verifies the following properties:

\begin{itemize}
\item $M$ is  a topological $m$-manifolds and contains a  boundaryless compact smooth $(m-2)$-manifold $L$
(the corner locus) such that
$M\setminus L$ is a smooth $m$-manifold with boundary.

\item There is a open neighbourhood $U$ of $L$ in $M$ and a homeomorphism  
$$\phi: U\to L \times  [0,1)\times [0,1)$$
such that for every $x\in L$, $\phi(x)=(x,0,0)$, and the restriction of $\phi$ to $U\setminus L$ is a diffeomorphism
onto $ L \times  [0,1)\times [0,1)) \setminus L\times \{(0,0)\}$.
\end{itemize}

By using these data we can prove that

\medskip

{\it There is a natural {\rm corner smoothing} procedure that gives a smooth structure on $M$ which is compatible with the given 
smooth structures on $L$ and $M\setminus L$.}

\medskip

For let us fix a homeomorphism $\tau: [0,1)\times [0,1)\to  B^2(0,1)\cap \HH^2$ which is a diffeomorphism
outside $(0,0)$ (for instance do it by using polar coordinates). Then set
$$\tau': L \times  [0,1)\times [0,1) \to L\times  (B^2(0,1)\cap \HH^2), \ \tau'(x,y,z)=(x,\tau(y,z))$$
and take the composition $\tau' \circ \phi: U\to L\times  (B^2(0,1)\cap \HH^2)$. Take on $U$ the differential structure
such that $\tau'\circ \phi$ is tautologically a diffeomorphism. A smooth atlas of this structure 
together with a smooth atlas of $M\setminus L$ make a smooth atlas on $M$ which by construction 
is compatible with the given smooth structures.  
Note that the induced smooth structure on 
$\partial M$ coincides up to diffeomorphism with the one obtained by gluing
the closure of the components of $\partial M \setminus L$ along the common boundary.
Similarly to  Proposition   \ref{gluing-unique} {\it the corner smoothing produces a unique
smooth structure  up to diffeomorphism} (we left the details as an exercise).  

\section {Uniqueness of smooth disks up to diffeotopy}\label{unique-disk}
Let $M$ be a smooth boundaryless $m$-manifold; a smooth embedding $$\beta: D^m \to M$$ 
of the closed unitary $m$-disk is called a {\it smooth $m$-disk} in $M$. If $M$ is oriented,
two smooth $m$-disks in $M$ are {\it co-oriented} if  both preserve or reverse the orientation, provided that
$D^m$ inherits the standard orientation of $\R^m$.

We have

\begin{proposition} \label{unique-disk2} Let $M$ be a connected  smooth boundaryless $m$-manifold. Let 
$\beta_r: D^m \to D_r \subset M$, $r=0,1$ be smooth $m$-disks in $M$ . Then

(1) If $M$ is oriented and  $\beta_0$ and $\beta_1$ are co-oriented, then there is a diffeotopy of $M$ which connects $\beta_0$ and $\beta_1$.
In particular there is an oriented smooth automorphism $f$ of $M$ such that $\beta_2= f\circ \beta_1$.

(2) If $M$ is not orientable then there is a diffeotopy of $M$ which connects $\beta_0$ and $\beta_1$.
In particular there is a smooth automorphism $f$ of $M$ such that $\beta_2= f\circ \beta_1$.
\end{proposition}

\Dim
In both cases, thanks to the homogeneity of $M$, possibly by composing $\beta_1$ with a diffeotopy,
we can assume that  $x_0=\beta_0(0)=\beta_1(0)$. Possibly up to radial isotopies centred at $0$, we can
assume that both $\beta_0$ and $\beta_1$ have image contained in a chart $\phi: W\to \R^m$ of $M$ such that
$\phi(x_0)=0$. Then we are reduced to the case $M=\R^m$, $\beta_r(0)=0$. Assume that the two disks are co-oriented.
Then we can easily adapt the proof of Proposition \ref{lin2} and conclude that both $\beta_r$ are isotopic
to a same linear embedding of the disk in $\R^m$. By applying   Thom's lemma we achieve (1).

If $M$ is not orientable, a priori the two disks localized in a chart at $x_0$ as above might be not co-oriented. However, by the non-orientability
of $M$, we can find a smooth simple loop $\lambda$ based at $x_0$ such that by ``sliding'' say $\beta_1$ along $\lambda$ we return back with the opposite
orientation. Then up to isotopy we can always reduce to two co-oriented disks in $\R^m$ and conclude as before.

 \cvd  

\section{Connected sum, shelling}\label{connected-sum}
Let us describe a further cut-and-paste procedure to construct compact manifolds.

$\bullet$ Let $M_1$ and $M_2$ be boundaryless, connected, compact smooth $m$-manifolds, $m\geq 1$. 

$\bullet$ For $s=1,2$, let 
$$\delta_s: D^m\to D_s\subset M_s$$ 
be a smooth embedding. 

$\bullet$ Consider 
$\tilde M_s = M_s \setminus {\rm Int}(D_s) $. 
Then $\tilde M_s$ is a compact
smooth manifold with one boundary component $V_s$ diffeomorphic to $S^{m-1}$.

$\bullet$ Let $\rho: V_1 \to V_2$, $\rho=\rho(\delta_1,\delta_2)$
being the diffeomorphism obtained by the restriction of  $ \delta_2 \circ \delta_1^{-1}:D_1 \to D_2$.
Finally consider the compact boundaryless manifold
$$W:=\tilde M_1 \amalg_\rho \tilde M_2 \ . $$

\medskip

Here is an equivalent description of the smooth manifold
$W$. Take the disjoint union $$(M_1 \setminus \delta_1(0))\amalg  (M_2 \setminus \delta_2(0))$$
and for every $(u,t)\in S^{m-1}\times (0,1)$ identify $\delta_1(tv)$ with $\delta_2((1-t)v)$.

Every $W$ obtained by implementing this procedure is called {\it a connected sum of $M_1$ and $M_2$}.

There is a natural {\it oriented} version, where $M_1$ and $M_2$ are oriented and $\delta_2\circ \delta_1^{-1}$
is orientation reversing.
The resulting connected sum is naturally oriented in a compatible way with 
$M_1$ and $M_2$. 

Every connected sum depends on the choice of the smooth $m$-disks $\delta_j$. 
We are going to analyze to which extent it is uniquely
defined up to diffeomorphism. 

 \medskip

\begin{proposition}\label{unique-CS} Let $M_1$ and $M_2$ be boundaryless, connected, compact smooth $m$-manifolds.  Then

(1) If both $M_1$ and $M_2$ are {\rm oriented}, then the {\rm oriented connected sum} $M_1\# M_2$ is well defined up
to oriented preserving diffeomorphism (i.e. it does not depend on the choice of the embeddings $\delta_s$, provided that
$\delta_2 \circ \delta_1^{-1}$ reverses the orientation).

(2) If at least one among $M_1$ and $M_2$ is not orientable,  then the  connected sum $M_1\# M_2$ is well defined up
to diffeomorphism (i.e. it does not depend on the choice of the embeddings $\delta_s$).
\end{proposition}

\Dim If both manifolds are oriented, possibly by pre-composing the smooth disks with the reflection $(x_1,\dots,x_m)\to (-x_1,\dots, x_m)$,
we can assume that the $m$-disks in $M_1$ preserve while the $m$-disks in $M_2$ reverse the orientation;
if at least one is non orientable, say $M_1$, while $M_2$ is orientable, then we can assume that the disks in $M_2$ are
co-oriented. By  Proposition \ref{unique-disk2}, in every case the disks in $M_1$ or $M_2$ entering different implementations 
of the connected sum procedure
are diffeotopic to each other. Then the proposition follows by several applications of  Proposition \ref {same-gluing}. 

\cvd

\medskip

\begin{remarks}\label{-M}{\rm
(1) When it is well defined, strictly speaking $M_1\# M_2$ 
denotes a diffeomorphism class of smooth manifolds.
Again we will do often the abuse to confuse it with any representative.

  (2) In the oriented case, if $-M$ denotes the connected oriented
    manifold $M$ endowed with the opposite orientation, then it can
    happen that $M_1\# M_2$ is not diffeomorphic to $-M_1 \# M_2$ via an
    orientation preserving diffemorphism. They are diffeomorphic if there is
    an orientation preserving  diffeomorphism between $M_1$ and $-M_1$.
    
    (3) The discussion about the connected sum works as well for compact manifolds
    with boundary, provided that the disks are embedded in their interior.}
    
\end{remarks}

\subsection{Thick connected sum, shelling}\label{shelling}
Let us keep the above setting. Assume furthermore that $M_s$ is a
boundary component of $\partial N_s$ of the compact $(m+1)$-manifold
$N_s$.  Then we can consider the topological quotient space 
$$ N_1 \amalg_{\hat \rho} N_2$$ 
where $\hat \rho: D_1 \to D_2$ is equal to
$\delta_2 \circ \delta_1^{-1}$. Arguing similarly to Section
\ref{partial-glueing}, it is not hard to show that this quotient space
carries a natural structure of smooth $(m+1)$-manifold with corners
which, by corner smoothing, leads to a well defined smooth manifolds denoted
$$N_1 \hat \# N_2$$ 
compatible with the smooth inclusions of $N_s$; moreover
$$\partial (N_1 \hat \# N_2) =(\partial N_1 \setminus M_1)\amalg
(\partial N_2 \setminus M_2)\amalg (M_1\connectedsum M_2) \ . $$ 
Naturally everything is well
defined (only) up to diffeomorphism, possibly in the oriented category.  

\begin{definition}\label{shelling-def}{\rm  In the above setting, if
$N_2= D^{m+1}$, then we say that $N:=N_1$ and $\tilde N:= N \hat \# D^{m+1}$ are
{\it related by a shelling} (along $M:=M_1$). 
}
\end{definition}

We have
\begin{proposition}\label{shelling-diffeo} If $N$ and $\tilde N$
are related by a shelling, then they are diffeomorphic, as well as 
$M \# S^{m}$ is diffeomorphic to $M$.
\end{proposition}

The proof involves several applications of the extension of isotopies and the disk unicity
as above. We left the details as an exercise. 

\cvd

\subsection{Weak connected sum, twisted spheres}\label{twisted-sphere} 
There is a weak variant of the connected sum procedure; by keeping
the notations of the beginning of Section \ref {connected-sum}, at the end
we take
$$\tilde M_1 \amalg_\beta \tilde M_2$$ where 
$$\beta: V_1 \to V_2$$
is {\it any} diffeomorphism, that is  we do {\it not} require that it is the
restriction of the composition of $m$-disks $\delta_2 \circ \delta_1^{-1}$.  
In the oriented
situation we require also that $\beta$ reverses the orientation. The
essential difference between the original procedure is that $\beta$
does not necessarily extend to a diffeomorphism $\hat \beta: D_1 \to
D_2$ between the whole embedded smooth $m$-disks.  If we incorporate
this last requirement about $\beta$, the present weak procedure is
equivalent to the previous one.  Without such a further requirement,
it is definitely different. 

We call {\it smooth twisted $m$-sphere}
any manifold obtained by implementing the weak connected sum procedure
starting from $M_1=M_2=S^m$.  We collect below a few (non exhaustive)
important facts about this topic.

\begin{proposition}\label{on-twisted}
  (1) If $1\leq m\leq 4$, then every diffeomorphism $\beta: S^{m-1}\to S^{m-1}$ extends to a diffeomorphism
$\hat \beta: D^m\to D^m$; hence every $m$-weak (oriented) connected sum is a (oriented) $m$-connected sum.

(2) For every $m\geq 1$, every smooth twisted sphere is {\rm
    homeomorphic} to $S^m$. If $1\leq m \leq 4$ it is diffeomorphic to
  $S^m$.

(3) There are smooth twisted $7$-spheres that are not diffeomorphic to $S^7$.
\end{proposition}

We limit to a few comments about the proofs, item by item.
\medskip

(1): For every $m$, possibly by composing with a reflections along a hyperplane of $\R^{m+1}$, 
it is not restrictive to assume that $\beta$ preserves the orientation of $S^m$. 

The validity (or not) of item (1) is invariant on the isotopy class of $\beta$.  

For $m=1$, item (1)  is immediate via linear parametrizations of the interval 
$D^1$. 

For $m=2$, we prove that $\beta$ is isotopic to the identity (which obviously extends to the identity of $D^2$). 
In fact, up to isotopy it is not restrictive to 
assume that $\beta$ is the identity on an open sub-arc $J$ of $S^1$ (diffeomorphic to $(0,1)$).
Let $J'$ be another open sub-arc  of $S^1$ such that $S^1 = J\cup J'$. We get an isotopy  of $\beta$ with the identity as follows
$$ H(x,t)= x \ {\rm if} \ x\in J, \  H(x,t)= tx +(1-t) \beta(x) \ {\rm if} \ x\in J' \ . $$ 

{\it (Smale Theorem)} For $m=3$, item (1) is already non trivial and due to Smale \cite{S1}; as above it is enough to prove
that $\beta$ is isotopic to the identity. A proof can be built by using special dynamical properties of integration of {\it planar} tengent vector fields,
the so called  {\it Poincar\'e-Bendixson Theory}. Up to isotopy we can assume that $\beta$ is the identity on a hemisphere.
So, via the stereographic projection, it is enough to prove that a diffeomorphism $g: \R^2 \to \R^2$ which is the identity outside
the unitary disk $D^2$ is isotopic to the identity through diffeomorphisms sharing this property. Again up to isotopy it is not
restrictive to assume that these diffeomorphisms are equal to the identity also on a collar of $S^1=\partial D^2$ in $D^2$.  
Consider the constant unitary vertical
tangent field on $\R^2$, $\vG_0 = e_2$, and let $\vG_1$ its image by means of the differential $dg$. 
These fields can be considered as smooth maps $\vG_i: D^2\to \C^*$ (completed by a constant map outside $D^2$).
We can lift them to maps $\tilde \vG_i: D^2\to \C$ via the universal covering map $\exp: \C \to \C^*$. By taking the convex
combinations $\tilde \vG_t:= t\tilde \vG_1 + (1-t)\tilde \vG_0$, $t\in [0,1]$, and projecting them back to $\C^*$, we get a homotopy
$\vG_t$ between $\vG_0$ and $\vG_1$ through nowhere vanishing tangent vector fields which are constant outside $D^2$ 
minus a collar of $S^1$. Now one would integrate the homotopy $\vG_t$ to a diffeotopy between $g$ and the identity.
This is a rather delicate task. A key dynamical property is that in the present situation {\it no maximal integral curves of $\vG_t$
are trapped in (the compact set) $D^2$}. In particular an integral line which crosses the upper hemicircle of $S^1$
pointing inside $D^2$, after a certain time crosses the lower hemicircle pointing outside. By elaborating on this fact, one
eventually constructs a desired isotopy of diffeomorphisms (for all details se also Section 6.4. of \cite{Mart}).

For $m=4$, (1) is difficult (see \cite{Ce}).

\smallskip

(2): It is easy to extend every $\beta$ as above to a {\it homeomorphism}
$\hat \beta: D^m \to D^m$; we can  
 get such a 
$\hat \beta$ by a radial extension sending for every $x\in S^{m-1}$,
the interval $[x,0]\subset D^m$ lineraly onto the interval
$[\beta(x),0]$ (this is also known as the {\it Alexander trick}). 
By the way, this is a diffeomorphism on $D^m \setminus \{0\}$,
$0$ being in general the only non smooth point.
By using this fact it is easy to show that every
twisted $m$-sphere is homeomeorphic to $S^m$. For $1\leq m \leq 4$  
it is diffeomorphic to $S^m$ thanks to item (1).
\smallskip

(3): These are the celebrated {\it Milnor's exotic $7$-spheres} \cite{M4}.  

\begin{remark}\label{inverse con-sum}{\rm Let $M$ be a compact oriented boundaryless smooth $m$-manifold.
Let $Y\subset M$ be a submanifold diffeomorphic to $S^{m-1}$ so that $M\setminus Y= M_1 \amalg M_2$
consists of two connected non compact manifolds. The closure $\hat M_s$ of $M_s$ in $M$  is a compact manifold
$\hat M_s$ with boundary equal to $Y$. Let us glue to $\hat M_s$ a disk $D^m$ via a diffeomorphism
$\rho_s:S^{m-1} \to Y$, obtaining two oriented boundaryless manifolds $\tilde M_s$. Then 
$$M= \tilde M_1 \connectedsum \tilde M_2 \ . $$
In general this factorization of $M$ is not unique. For example the standard $S^7$ can be expressed as $S^7\connectedsum S^7$
as well as the connected sum of two exotic $7$-spheres.}
\end{remark}   

\section{Attaching handles}\label{handle}
This is a very important procedure. We will see in Chapter \ref {TD-HANDLE}
that every compact manifold admits ``handle decompositions'' that is it can be built (up to diffeomorphism)
by iterated applications of this basic attaching procedure.
\smallskip

For every $m\geq 0$, for every $0\leq q \leq m$,
$$ H^q = H^{q,m}=D^q\times D^{m-q}$$
is the {\it standard $q$-handle of dimension $m$}. If clear from the contest, we will omit to indicate the dimension; $q$ is also called the
{\it index} of the handle. Strictly speaking  such a handle $H^q$ is a manifold with corner with boundary
$$\partial H^q = (S^{q-1}\times D^{m-q})\cup (D^q \times S^{m-q-1}) \ ; $$
up to smoothing it is diffeomorphic to $D^m$ endowed with a determined decomposition by submanifolds of $\partial D^m= S^{m-1}$.

Let us fix a few terminology.

$\bullet$ $\Sigma_a:=S^{q-1}\times \{0\} \subset \Tt_a:=S^{q-1}\times D^{m-q}$ are called respectively the {\it $a$-sphere} and the 
{\it $a$-tube}  of $H^q$.

$\bullet$ $\Sigma_b:=\{0\} \times S^{m-q-1} \subset \Tt_b:=D^q \times S^{m-q-1}$ are called respectively the {\it $b$-sphere}  and the {\it $b$-tube} 
 of $H^q$.

$\bullet$  $C:=D^q\times \{0\}$ is called the {\it core} of the handle.

$\bullet$ $C^*:= \{0\} \times D^{m-q}$ is called the {\it co-core} of the handle.

\smallskip

Note that the $a$-sphere is the boundary of the core, the $b$-sphere is the boundary of the co-core; the core and the co-core intersect
transversely only at $(0,0)$. $\Tt_a$ and $\Tt_b$ intersect at the respective boundaries both equal to $S^{q-1}\times S^{m-q-1}$.

\medskip

Let $N$ be a compact smooth $m$-manifold with boundary. Given a $q$-handle $H^q$ of dimension $m$, let
$h:\Tt_a \to \partial N$ be a smooth embedding. Then $S_a:= h(\Sigma_a)$ is the {\it embedded (attaching) $a$-sphere};
$T_a:=  h(\Tt_a)$ is a tubular neighbourhood of $S_a$ in $\partial N$, endowed by means of $h$ of a {\it global trivialization}.
$T_a$ is also called the {\it embedded (attaching) $a$-tube}. Consider the topological quotient space 
$$N\amalg_{h} H^q$$
by the equivalence relation on the disjoint union $N\amalg H^q$ which identifies every $x\in \Tt_a$ with $h(x)\in T_a$.
Then
$N\amalg_{ h} H^q$ has a natural structure of manifold with corner which by smoothing leads to a smooth manifold
well defined up to diffeomorphism. Considered up to diffeomorphism, we say that  $N\amalg_{h} H^q$ is the smooth
manifold obtained by attaching a $q$-handle to $N$ via the attaching map $h$. At this point it is routine to apply 
as above the extension of isotopies to diffeotopies and get:

\begin{proposition}\label{handle-prop} Up to diffeomorphism, $N\amalg_{h} H^q $ only depends on the isotopy class of the 
attaching embedding $h$.
\end{proposition}

\cvd

\medskip

Here is a few complements about attaching handles.
\smallskip

(1) Up to diffeomorphism, the boundary of $N\amalg_h H^q$ is given by 
$$\partial (N\amalg_h H^q)= (\partial N \setminus {\rm Int}(T_a))\amalg_{h_{|\partial \Tt_b}} \Tt_b \ ; $$
sometimes we denote it by $$\sigma(\partial N,h)$$ and call it the $(m-1)$-manifold obtained by {\it surgery} on $\partial N$ with {\it surgery
data $h$}.

(2) If $N$ is oriented and $q>1$, then also $N\amalg_h H^q$ can be oriented in a compatible way. In fact as $q>1$, the $a$-tube is connected
and we can take the orientation of $H^q$ such that the gluing diffeomorphism $h:\Tt_a \to T_a$ reverses the orientation. For $q=1$,  
$\Tt_a$ is not connected and it is not always possible to make $h$ orientation reversing on both components. Attaching $1$-handles is the only
case which imposes some constraints in order to perform the construction within the oriented category.

(3) If $N$ is connected and $q>1$, then also $N\amalg_h H^q$ is connected. In fact the connected $T_a$ is contained in one connected
component of $\partial N$ and by attaching $H^q$, which is connected, connectedness is preserved. By attaching a $1$-handle we can reduce
the number of connected components by $1$. This happens if the connected components of $T_a$ belong to different components
of $\partial N$.

(4) The $a$-tube of a $0$-handle is empty; then attaching a $0$-handle to $N$ means to ``create'' a new connected component
diffeomorphic to $D^m$. The $a$-tube of a $m$-handle is the whole boundary of $D^m$. Hence by attaching a $m$-handle we fill
a spherical component of $\partial N$ (if any, otherwise we cannot attach any $m$-handle). 

(5) Up to diffeomorphism,  the thick connected sum can be rephased in terms of attaching a $1$-handle to $N_1$ and $N_2$ 
with one component of $T_a$ in $\partial N_1$ and the other in $\partial N_2$. Similarly by suitably attaching a $1$-handle to
$$(M_1\times [0,1])\amalg (M_2\times [0,1])$$ we get a manifold $W$ such that 
$$\partial W = (M_0 \amalg M_1) \amalg (M_1 \connectedsum M_2)$$
(possibly in the oriented category).

\begin{remark}\label{gen-glue}{\rm Attaching a handle is an instance of the following more general
gluing procedure: for $j=1,2$, let $Y_j$ be a $(m-1)$ sub-manifold with boundary $\partial Y_j$ of
$\partial M_j \subset M_j$. Let $\rho: Y_1 \to Y_2$ be a diffeomorphism. Then $M_1\amalg_\rho M_2$
is in a natural way a $m$-manifold with corners, hence a well defined smooth manifold up to corner
smoothing (and up to diffeomorphism).}
\end{remark}

\section{Strong embedding theorem, the Whitney Trick}\label{whitney-trick}
The aim of this section is to provide information about the following theorem,
the proof introduces the very important so called ``Whitney trick" \cite {Whit2}. 

\begin{theorem}\label{stong-emb}
Every compact boundaryless smooth $m$-manifold $M$ can be embedded into $\R^{2m}$.
\end{theorem}  

\medskip

{\it A sketch of proof.}
We limit to a rough outline of the proof, stressing anyway that it is
substantially different from the weak immersion/embedding theorem \ref {Weak -Whitney}. This last
is enterely based on so called ``general position arguments'' or,
equivalently, on {\it transversality} (concepts that we will develop
in Chapter \ref{TD-TRANSVERSE} although  we are anticipating a few applications).
By pushing the general position arguments  (see Section \ref{miscellaneaT}), we can at most
refine the weak immersion theorem and get that a ``generic
immersion'', say $\pi: M\to \R^{2m}$, of our compact boundaryless $m$-manifold in
$\R^{2m}$ has the further properties:
\smallskip

{\it The inverse image of every point in $\pi(M)\subset \R^{2m}$
  consists in at most $2$ points; if $\pi(p)=\pi(p')=q$, then
$\R^{2m}= d_p \pi(T_pM)\oplus d_{p'} \pi (T_{p'}M)$.
Then, by compactness of $M$, there is in the image of $\pi$ a finite number of such 
``simple normal crossing  points''.}  

\smallskip

We can start with such a generic immersion.  If there are normal crossing
points, they persit under any small perturbation of the
immersion. To get an embedding we must operate
a robust modification of $\pi$.  Basically there are two ``moves'':

\begin{enumerate}
\item Introduce if necessary a further crossing point. 
\item Eliminate a couple of double points by applying the so
  called {\it Whitney Trick}.
\end{enumerate}

\smallskip

As we are going to see, this scheme actually works for $m\neq 2$;
fortunately for $m=2$, the strong embedding theorem holds as a
corollary of the {\it classification of smooth compact surfaces} (see
Chapter \ref{TD-SURFACE}). So we definitively assume here that $m\neq 2$.
Moreover it is not restrictive to assume  that $M$ is connected.
\smallskip

$\bullet$ The basic local model for a single self-intersection point
is as follows:
$$\alpha: \R^m \to \R^{2m}, \ \alpha(t_1,t_2,\cdots,t_m) =
 \left(t_1 - \frac{2t_1}{u},t_2, \dots, t_m,\frac{1}{u},\frac{t_1t_2}{u},
\frac{t_1t_3}{u}, \cdots, \frac{t_1t_m}{u} \right)$$
where
$$u=(1+t_1^2)(1+t_2^2)\cdots(1+t_m^2) \ . $$
It is an embedding except
for the points $(1,0,\dots, 0), (-1,0, \dots,0)$ which are sent to
$0\in \R^{2m}$.  Moreover, when $||t||\to +\infty$, $\alpha$ tends to
the usual linear embedding $(t_1,\dots, t_m)\to (t_1,t_2,\dots ,t_m, 0,\dots 0)$
of $\R^m \subset \R^m\times \R^m=\R^{2m}$. To add such a double point to a given
immersion $\pi$, we can
do it locally in a chart at a point $q\in \pi(M)$ where at $q\sim 0$,
$\pi(M)$ looks like the image of the above linear embedding. Then by
using two suitable bump functions on $\R^{m}$ at $0$ and at infinity
respectively, and the associated partition of unity, it is not hard to
modify $\pi$ to get one with one more self-intersection point.

\begin{remark}\label{signed-dp}{\rm Give $\R^m$ and $\R^{2m}$ the standard
orientation; then the single self-intersection point has a sign. Its mirror image
has the opposite sign.}
\end{remark}

$\bullet$ The Whitney Trick applies at a {\it Whitney disk} $D$ connecting two
crossing points $q_1$, $q_2$ in $\pi(M)$.  This means that the following
pattern is realized:

\smallskip

(1) There is an embedded smooth circle $\gamma$ in $\pi(M)$ with two corners
at $q_1$ and $q_2$; these divide $\gamma$ in two arcs with closures $\gamma_1$ and $\gamma_2$
respectively; these closed arcs $\gamma_j$, $j=1,2$, are contained into smooth open
$m$-disks $U_j$ in $\pi(M)$, their union is an open neighbourhood of $\gamma$ in $\pi(M)$,
they intersect transversely each other at $\{q_1,q_2\}$, and do not
contain other crossing points of $\pi(M)$;
\smallskip

(2) There are:

- a $2$-disk $\Dd$ in $\R^2$ with boundary $\partial \Dd$ with two corners
$a_1$, $a_2$ which is contained in the union of two smooth arcs $\lambda_1$, $\lambda_2$ 
in $\R^2$ which intersect transversely at $\{a_1,a_2\}$; 

- an embedding $\psi: U \to \R^{2m}$
where $U$ is an open $2$-disk in $\R^2$ containing  $\Dd\cup (\lambda_1 \cup \lambda_2)$,
such that 
\begin{itemize}
\item $\psi(\lambda_j)\subset U_j, \  j=1,2$;  
\item $\psi(\partial \Dd,\{a_1,a_2\})=(\gamma,\{q_1,q_2\})$;
 \item for every $x\in \lambda_j$, j=1,2,  $d_x\psi(T_xU)\cap T_{\psi(x)}U_j = d_x\psi(T_x\lambda_j)$;  
\item $\psi({\rm Int}(\Dd)) \subset \R^{2m}\setminus \pi(M) $.
\end{itemize}
We summarize (1) and (2) by saying that the smooth $2$-disk with corners $D:= \psi(\Dd)$ is 
{\it properly embedded} into $(\R^{2m},\pi(M))$ and {\it connects the crossing points} $q_1$, $q_2$.

Moreover, we require:  
\smallskip

(3) We can extend the embedding $\psi$ to a  parametrization of a neigbourhood
of $D$ in $\R^{2m}$ by a {\it standard model}, that is
to an embedding
$$\Psi: U\times \R^{m-1}\times \R^{m-1}\to \R^{2m}$$
such that $\Psi(\lambda_1\times \R^{m-1}\times \{0\})=U_1$ and $\Psi(\lambda_2\times\{0\} \times \R^{m-1})=U_2$.
\medskip

 Thanks to such a standard model, it is not hard to realize that
a Whitney disk (if any) can be used as a guide to construct a $1$-parameter family of
immersions, with compact support around $D$, by ``pushing $M$ across
$D$'', eventually removing $q_1,q_2$ without modifying the
configuration of the other crossing points.

\begin{remark}\label{opposite-sign}{\rm We can fix local orientations around a Whitney disk.
The required properties implies that the two crossing points connected by the disk 
have {\it opposite signs} with respect to such orientations}
\end{remark} 

\smallskip

$\bullet$ To  conclude the proof of the embedding theorem,
we have to show that for every generic projection, possibly after having inserted a new crossing point (recall Remarks \ref{signed-dp} and
\ref {opposite-sign}) there is  
a couple of crossing points connected by a Whitney disk which can be eliminated.
For $m=1$ this follows by
somewhat subtle but elementary planar considerations.  For $m>2$, 
we will discuss this issue  within a larger range of application of the Whitney
trick in Chapter \ref{TD-HIGH} (see Remark \ref {on-WT} (2) and Proposition \ref{Conc-inductive}).

\cvd

\smallskip

\begin{remarks}\label{on-WT}{\rm 
(1) If $m=2$, the circle $\gamma$ can be constructed as well and one could construct a generically
immersed disk $D$ in $\R^{4}$, bounded by $\gamma$, but we cannot
exclude the existence of crossing points of $D$ itself  or of transverse intersection
of $D$ with $\pi(M)$ apart from $\gamma$.
\smallskip

 \begin{figure}[ht]
\begin{center}
 \includegraphics[width=5cm]{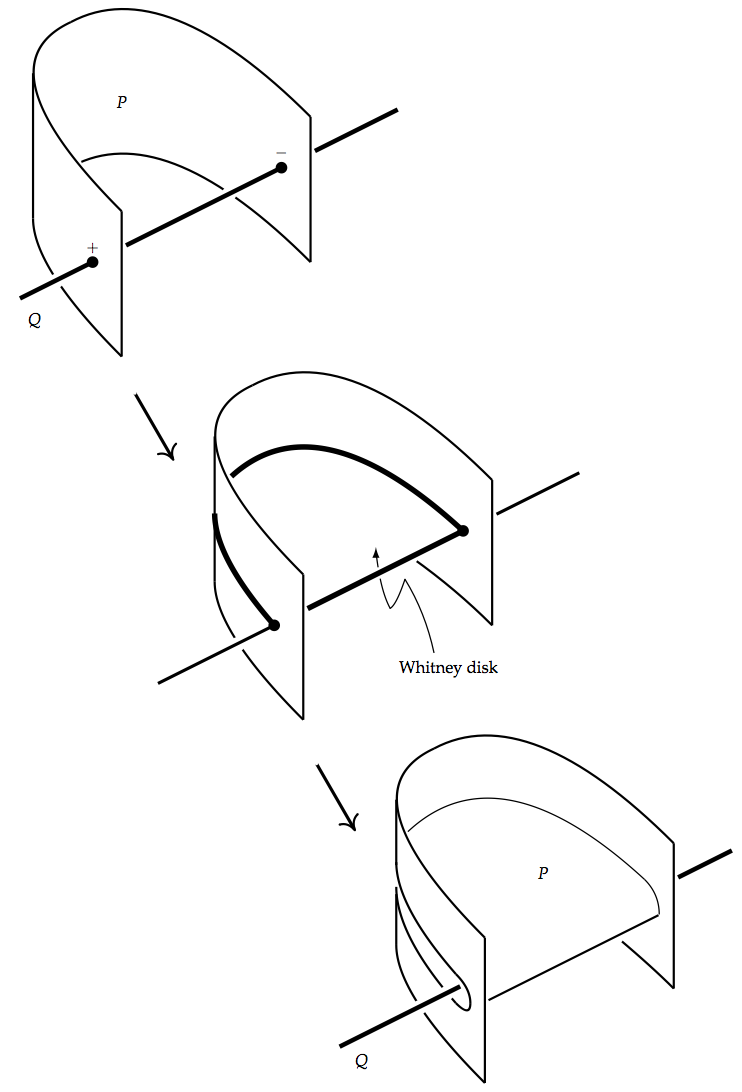}
\caption{\label{WT} Whitney's trick.} 
\end{center}
\end{figure}

(2) The notion of Whitney disk, hence the Whitney trick,
    can be extended to eliminate couple of tranverse intersections of two submanifolds
    $P, Q$ of a given manifold $M$, such that $\dim M = \dim P + \dim Q$ (the boundary
    loop $\gamma$ being formed by two arcs in $P$ and $Q$ respectively). This
    technique has been of absolute importance in the achievement of fundamental
    results for smooth manifolds of sufficiently high dimension (see Chapter \ref{TD-HIGH}).
    The fact noticed above that the scheme does not apply in the case
    $\dim M=4$ has been the ultimate reason for special and astonishing phenomena occurring
    in the realm of $4$-manifolds. We will develop these comments much later in the text (see Chapter \ref{TD-4}).}
\end{remarks}

\section{On immersions of $n$-manifolds in $\R^{2n-1} $}\label{n-in-2n-1}
The aim of this section is to provide some information about the following hard immersion theorem \cite{Whit3}.

\begin{theorem}\label{W-hard-Imm} Every compact boundaryless $n$-manifold $M$
can be immersed into $\R^{2n-1}$.
\end{theorem}

\cvd

\smallskip

It is not restrictive to assume that $M$ is connected.  Similarly as
in the discussion about the hard embedding theorem, ``hard" means that
it is not only based on general position arguments. This kind of
arguments (mostly in the spirit of ``multi-jet-transversality'' - see
Section \ref{miscellaneaT} ) allows to preliminarly determine the {\it
  generic maps} $f:M\to \R^{2n-1}$ which in general are not
immersions.  For simplicity let us give a few details for $n=2$ (the
general case is similar). The local models of such a generic map are
all realized by
$$ g: \R^2_{u,v} \to \R^3_{x,y.z}, \ x=u^2, \ y=v, \ z=uv \ . $$ The
line $\{v=0\}$ is the non injectivity locus of this map and its image
is a half lines. The image of every other line $\{v=c\}$ is the
parabole $x=(z/c)^2$ in the hyperplane $\{y=c\}$. The point $0\in
\R^2$ is the unique at which the map $g$ is not an immersion and its
image $0=g(0)$ is called the {\it Whitney point} in the model. The
transverse intersection with the image of $g$ of a small sphere around
the Whitney point is a wedge of two smooth circles.  The restriction
of $g$ to $\R^2\setminus \{0\}$ is a generic immersion, that is along
the the image of $\{v=0\}\setminus \{0\}$ there are two transverse
branches of the images of $g$.

In general, we can describe qualitatively a generic maps $f: M\to
\R^{2n-1}$ as follows. Assume first that $n\geq 3$. The image say
$\Sigma$ of the non injectivity locus is a compact $1$-dimensional
submanifold of $\R^{2n-1}$ possibly with boundary; $W=\partial \Sigma$
is formed by the so called Whitney points of $f$. The restriction of
$f$ to $\tilde W:= f^{-1}(W)$ is a bijection onto its image and $f$ is
not an immersion at every point of $\tilde W$.  $\tilde \Sigma
:=f^{-1}(\Sigma \setminus \partial \Sigma)$ is a smooth (non compact)
$1$-submanifold of $M$ and the restriction of $f$ to $\tilde \Sigma$
is a double covering map onto the interior of $\Sigma$.  The
restriction of $f$ to $M\setminus \tilde W$ is a generic immersion, so
that locally along every component of the interior of $\Sigma$, there
are two transverse branches of the image of $f$.

If $n=2$ the situation is a bit more complicated. In fact beyond the
Whitney points, $\Sigma$ has in general also a finite set of three
branches crossing points (the ``triple points" of the image) at which
the local model for the generic immersion of $M\setminus \tilde W$ is
given by three hyperplanes of $\R^3$ in general position.

These generic maps are {\it stable} in the sense that their
qualitative features are preserved up to small smooth perturbations.
Starting from a generic map $f: M\to \R^{2n-1}$, we have to perform a
robust alteration of it in order to get an immersion $\hat f: M\to
\R^{2n-1}$.  The Whitney points are partitioned by couples of points
which are connected by a smooth arc contained in $\Sigma$.  Then we
perform a kind of rather subtle ``surgery'' along each such an arc
$\gamma$. To give an idea, assume that $n=2$ and that, for simplicity,
the arc $\gamma$ connecting two Whitney points does not include triple
points.  Remove from $f(M)$ the intersection with the interior of a
small smooth ``$\epsilon$-neighbourhood" $U$ (diffeomorphic to $D^3$)
of $\gamma$ in $\R^3$ whose boundary intersects transversely $f(M)$ at
two smooth circles; then fill them by two disjoint embedded
$2$-disks. In this way we get $\Sigma'$ from which two Whitney points
have been eliminated; in fact $\Sigma'$ is the image of the non
injectivity locus of a generic map $f': M'\to \R^3$, where $M'$ is a
surface obtained from $M$ by cutting and pasting. To eventually restore a
map $f": M\to \R^3$ ones connects again the above $2$-disks by
attaching a suitably oriented $1$-handle embedded into the smooth
$3$-disk $U$. By doing it along every arcs $\gamma$ we eventually get
a desired generic immersion $\hat f: M\to \R^3$. Moreover, $\hat f$
can be obtained arbitrarily close to the given generic map $f$ in the
$\Cc^0$-topology.

\subsection{On  Smale-Hirsch immersion theory}\label{smale-H}
Whitney's hard immersions theorem has been reobtained later as a
non trivial application of Hirsch immersion theory \cite{H2}. 
Extending early Smale's results in the case when $M$ is a sphere, 
this faces the general question of the existence of immersions
$f: M \to N$, $n=\dim N > \dim M = m$, and the classification of immersions
in a given homotopy class of maps from $M$ to $N$ up to 
regular homotopy (two
immersions $f_0, f_1: M\to N$ are {\it regularly homotopic} if they
are connected by a homotopy $f_t$ such that for every $t\in [0,1]$,
$f_t$ is an immersion). Remarkably  these questions are translated 
into homotopy theoretic problems. When $N=\R^{m+k}$, $k\geq 1$,
the existence problem can be translated as follows. By the easy
Whitney immersion theorem, there are immersions $f: M \to \R^{m+m}$,
and by using the standard metric $g_0$ on $\R^{m+m}$ we have 
the induced normal map
$$\nu_f: M \to \GG_{2m,m}, \ \nu_f(x)= (d_xf(T_xM))^\perp \ . $$
Then there exists an immersion $\hat f: M \to \R^{m+k}$, $1\leq k \leq m$
if and only if there exists an immersion $f$ as above and a map
$$\hat \nu: M\to \GG_{m+k,k}$$
such that the vector bundle $\hat \nu^* (\Vv_{m+k,k})$ is weakly stably isomorphic
to $\nu_f^*(\Vv_{2m,m})$. By the classification of vector bundles on compact manifolds,
this is equivalent to establish a homotopy between classifying maps.
Moreover, given such a map $\hat \nu$, there is an immersion $\hat f$ such that
$\hat \nu = \nu_{\hat f}$.

When $N=\R^{m+k}$, all immersions are homotopic to each other; it turns out that
$f_0$ and $f_1$ are regularly homotopic if and only if the bundle maps
$[\nu_{f_0}, \nu_{f_0}^*]$ and $[\nu_{f_1}, \nu_{f_1}^*]$ are homotopic through
bundle maps over a (ordinary) homotopy connecting $f_0$ and $f_1$.  

The following corollary is immediate.

\begin{corollary} If $M$ is parallelizable then it can be immersed into $\R^{m+k}$ for every $k\geq 1$.
\end{corollary}

\cvd

For every $m\geq 0$, let $i(m)$ be the minimum
$k\geq 1$ such that every compact boundaryles $m$-manifold $M$ can be immersed into
$\R^{m+k}$. By the hard Whitney immersion theorem, we have that $i(m) \leq m-1$.
By using the above translation of the problem into (hard) homotopy theoretic ones,
we eventually know the exact value of $i(m)$, see \cite{RC}.

\begin{theorem} For every $m\geq 0$, $i(m)=m- \alpha(m)$,
where $\alpha(m)$ is the number of $1$ in the dyadic
expansion of $m$.
\end{theorem}
 
\section{Embedding $n$-manifolds in $\R^{2n-1}$ up to surgery}\label{n-in-2n-1-emb} 
By construction if we use Whitney's method or perturbing an immersion $f: M \to \R^{2n-1}$
whose existence is an application of Hirsch results, we can assume anyway to deal 
with {\it generic} immersions, $M$ being any compact connected boundaryless 
$n$-manifold. So if $n\geq 3$, adopting the above notations, $\Sigma$ is a compact
boundaryless $1$-submanifold of $\R^{2n-1}$. For every component $C$
of $\Sigma$, locally along $C$ we see two transverse branches of the
image of $f$; $\tilde C:= f^{-1}(C)$ is a compact boundaryless
$1$-submanifold of $M$ and the restriction of $f$ to $\tilde C$ is a
$2$-folds covering which a priori can be non trivial ($\tilde C$
connected) or trivial ($\tilde C$ with two connected components).  

The main aim of this section is to show that starting from a generic immersion $f: M
\to \R^{2n-1}$ as above, by attaching suitable ``round handles'' to
$M\times [0,1]$ at $M\times \{1\}$ we get a $(n+1)$-manifolds $W$ such
that $\partial W = M \amalg \hat M$ (hence $\hat M$ is obtained by a
kind of ``surgery" on $M$) and $f$ can be altered on $\hat M$ to get
an {\it embedding} $\hat f: \hat M \to \R^{2n-1}$.  This construction
is due to Rohlin (see the translations of his papers in \cite{GM}) and
will be used in Chapters \ref{TD-3} and \ref{TD-4}.

Let us analyze more closely the properties of such a generic immersion. 
$C$ has a tubular neighbourhood
$U\sim C\times D^n$ in $\R^{2n-1}$ such that $\tilde U :=
f^{-1}(f(M)\cap U)$ is a tubular neighbourhood of $\tilde C$ in $M$. A
priori there are two possibilities for $U$. Either it is identified
with the mapping cylinder of
$$h_0: D^n \times D^n \to D^n\times D^n, \ h_0(y,z)=(z,y)$$
or to the mapping cylinder of
$$h_1: D^n \times D^n \to D^n\times D^n, \ h_1(y,z)=(y,z) \ . $$
In both cases, the subset 
$$X:= (\{0\}\times D^n)\cup (D^n\times \{0\})$$ is $h_j$-invariant,
$j=0,1$, and the mapping cylinder of the restriction of $h_j$ to $X$
realizes the image $f(\tilde U)$ in $U$.  The tubular neighbourhood
$\tilde U$ can be realized respectively either as the mapping cylinder
of
$$ g_0: \{0,1\}\times D^n\to  \{0,1\}\times D^n, \ g_0(u,x)= (1-u, x)$$
or of
$$ g_1: \{0,1\}\times D^n\to  \{0,1\}\times D^n, \ g_1(u,x)= (u, x)$$
and in both cases, the restriction of $f$ to $\tilde U$ can be expressed as
$$f(u,x,t)=(ux,(1-u)x,t) \ . $$ The first case would correspond to the
non trivial covering $\tilde C \to C$; the second to the trivial
one. However, as $\R^{2n-1}$ is orientable, then also $U$ must be
orientable and one easily sees that this constraint cannot be realized
in the first case if $n$ is even. So we have proved
\begin{lemma} If $n=\dim M\geq 3$ is even only trivial coverings $\tilde C \to C$ can occur.
\end{lemma}

\cvd

\smallskip

We are going now to construct $W$, $\partial W = M \amalg \hat M$ and the embedding
$\hat f: \hat M \to \R^{2n-1}$ with the desired features.
Let $C$ be a component of $\Sigma$. Use the above models for the
neighbourhoods $U$, $\tilde U$.  Consider $\frac{1}{2}\tilde U \subset
\tilde U$ obtained as the mapping cylinder of the restriction of $g_j$
to $\{0,1\}\times \frac{1}{2}D^n$ and set
$$\tilde U' := \tilde U \setminus {\rm Int} \frac{1}{2}\tilde U \ . $$  
Define the map
$$ \hat f: \tilde U' \to U$$
by 
$$\hat f(0,x,t)=(\phi(|x|)(-x_1,x_2,\dots, x_n), x, t), \ \ \hat f(1,x,t)=(x, \phi(|x|)(-x_1,x_2,\dots, x_n), t)$$
where $x=(x_1,\dots, x_n)$ and 
$$\phi: [1/2,1]\to [0,1]$$ is a smooth strictly decreasing function
which coincides with $t\to -t +1/2$ near $t=1/2$, $\phi(1)=0$ and
$\phi$ is flat at $1$.  The image of $\hat f$ in $U$ is the mapping
cylinder of the restriction of $h_j$ to an invariant subset $\tilde X$
of $D^n \times D^n$ wich coincides with $X$ near the boundary.
$\tilde X$ is diffeomorphic to two disjoint copies of $D^n$ hence it
``desingularizes'' $X$.  The map $\hat f$ extends to the whole of $M
\setminus {\rm Int}\frac{1}{2}\tilde U$ by taking the restriction of
$f$ to $M\setminus \tilde U$.  Do it for every component of $\Sigma$
(by using pairwise disjoint tubular neighbourhoods). Thus we have
obtained a $n$-submanifold, say $\tilde M$, of $\R^{2n-1}$ which is
the image of a smooth map $\hat f: M_0 \to \R^{2n-1}$, where $M_0$ is
a submanifold with boundary of $M$ obtained by removing a system of
small open tubular neighbourhoods of the $\tilde C$'s.  It turns out
that the quotient $$\hat M:= M_0/\hat f$$ is in a natural way a
boundaryless compact manifold and the induced map (we keep the name)
$$\hat f: \hat M \to \tilde M$$
is a diffeomorphism.  For every component $C$, the identification induced by $\hat f$ at the corresponding
boundary components of $M_0$ is given by
$$ (u,x_1, x_2, \dots, x_n,t)\sim (1-u,-x_1,x_2,\dots, x_n,t) \ . $$
It remains to describe the ``handles'' attached to $M\times [0,1]$ at $M\times\{1\}$ 
producing a $(n+1)$-manifold $W$ such that $\partial W = M\amalg \hat M$.
There is one such a handle for every component $C$. If $\tilde C \to C$ is the trivial
covering, let $H$ be the mapping cylinder of the identity of $[0,1]\times D^n$.
Then attach $H$ at $M\sim M\times \{1\}$ along $\tilde U$, by means of the attaching map
which identifies $(0,x,t)$ (resp. $(1,x,t)$)  of $H$ with $(0,x,t)$ ($(1,x,t)$) of $\tilde U$.
If $\tilde C \to C$ is non trivial (recall that it happens only if $n$ is odd) then 
we do similarly by using the mapping cylinder $\tilde H$
of the map
$$k: [0,1]\times D^n \to [0,1]\times D^n, \  k(v,x_1,x_2,\dots,x_n)=(1-v, -x_1,x_2,\dots, x_n) \ . $$
This complete the construction. We stress that by the very construction: 
\smallskip

{\it If $M$ is orientable then also the $(n+1)$-manifold $W$ constructed so far 
and the manifold $\hat M$ embedded in $\R^{2n-1}$ such that $\partial W= M \amalg \hat M$ are orientable.}

\begin{remark}\label{all-W}{\rm The constructions and the considerations of this section hold by starting from
any generic immersion $f:M \to W$ from a compact (possibly orientable) boundaryless $n$-manifold into an arbitrary (possibly orientable) 
$(2n-1)$-manifold $W$.}
\end{remark}

\section{Projectivized vector bundles and blowing up}\label{blow-up}
$\R^n$ can be considered as a vector bundle over the $0$-manifold $M=\{0\}$.
The projective space $\PP^{n-1}(\R)$ can be considered
as a fibration over $M$ which ``projectivizes'' the given vector bundle.
If $$\xi:= p: E \to M$$ is  any vector bundle (for example the tangent bundle),
over a compact $m$-manifold $M$ with fibre $\R^n$, we can perform the above projectivization 
fibre by fibre and obtain a fibration 
$$\pp: \PP(E)\to M$$ 
with fibre $\PP^{n-1}$. Every local trivialization
$W\times \R^n \sim p^{-1}(W)$ of the vector bundle gives rise to a local trivialization
$W\times \PP^{n-1} \sim \pp^{-1}(W)$. If  $(E,p)$ is defined by means
of a cocycle $\{ \mu_{i,j}: W_i\cap W_j \to {\rm GL}(n,\R) \}$, then it induces
a cocycle with values in the projectivized linear group $\PP{\rm GL}(n,\R)$
that defines $(\PP(E), \pp)$. The total space $\PP(E)$ is a compact manifold of dimension
$m+n-1$. A point in $\PP(E)$ is a line $l_x$ in $E_x= p^{-1}(x)$ for some
$x\in M$.  We can pull-back $\xi$ to $\PP(E)$ via the projection $\pp$ and obtain 
the vector bundle $\pp^*(\xi)$ over $\PP(E)$. We note that the restriction of $\pp^*(\xi)$
to every fibre of $\pp$
is a product (trivial) bundle. Moreover, $\pp^*(\xi)$ has a canonical {\it tautological} sub-bundle
of rank $1$ (i.e. a line bundle) $\lambda_\xi$ : the total space is 
$$\Lambda_\xi=\{(l_x,v)\in \pp^*(\xi); \ v\in l_x\} $$
with the natural projection onto $\PP(E)$.
Its fibre over $l_x$  is the line contained in the fibre of $\pp^*(\xi)$ at $l_x$,
made by the vectors belonging  to $l_x$. By using for instance an auxiliary riemannian 
metric
on the total space of $\pp^*(\xi)$ we realize that up to strict equivalence
it canonically splits as a direct sum
$$ \pp^*(\xi)\sim \lambda_\xi \oplus \beta_\xi$$
where also the bundle $\beta_\xi$ is well defined up to strict equivalence.
By iterating this construction starting again from $\beta_\xi$, we eventually get

\begin{proposition}\label{flag} For every vector bundle $\xi: E\to M$ over a compact manifold $M$, there
is a canonical construction (via iterated projectivization of vector bundles) that 
produces a smooth compact manifold $F(\xi)$ endowed with a surjective smooth map
$$f_\xi: F(\xi) \to M$$ such that the vector bundle $f_\xi^*(\xi)$ over $F(\xi)$
splits as a direct sum of line bundles. In particular this applies to the tangent bundle
of $M$.
\end{proposition}

\cvd

\subsection{Blowing up along smooth centres}\label{BU}
Let us start with the blowing up of $\R^n$, $n\geq 1$, with centre 
the $0$-submanifold $X=\{0\}$.
Consider
$$\R^n\times \PP^{n-1}(\R)$$
where $\R^n$ is endowed with usual coordinates $x=(x_1,\dots,x_n)$, while the projective
space is endowed with {\it homogeneous} coordinates $t=(t_1,\dots, t_n)$. 
Set
$$ \BB(\R^n,0):=\{ (x,t)\in \R^n\times \PP^{n-1}(\R); \ x_it_j=x_jt_i, \ i,j=1,\dots, n \}$$
this is well defined because the  equations are homogeneous in the $t$'s.
Denote by 
$$\rho: \BB(\R^n,0)\to \R^n$$
the restriction of the projection onto $\R^n$. These objects verify several interesting
properties:
\smallskip

(1) {\it $\BB(\R^n,0)$ is a smooth $n$-manifold.}  

\noindent If $U_j$ is the standard chart
of the projective space with non-homogeneous coordinates 
$y_i=t_i/t_j$, $t_j\neq 0$, $i\neq j$, then one readily checks that
$\BB(\R^n,0)\cap (\R^n\times U_j)$ is given as the graph of the smooth function
$x_i= x_jy_i$, $i\neq j$.
\smallskip

(2) {\it The restriction 
$$\rho: \BB(\R^n,0)\setminus \rho^{-1}(0)\to \R^n\setminus \{0\}$$
is a diffeomorphism}. 

\noindent Assume that 
$((a_1,\dots,a_n), (y_1,\dots,y_n))\in \BB(\R^n,0)$ with some $a_i\neq 0$.
Then for every $j$, $y_j= (a_j/a_i)y_i$ is uniquely determined as a point of
$\PP^{n-1}(\R)$. This also shows  that 
$$(a_1,\dots, a_n)\to ((a_1,\dots, a_n),(a_1,\dots, a_n))\in \BB(\R^n,0)\setminus \rho^{-1}(0)$$
defined for $(a_1,\dots,a_n)\in \R^n \setminus \{0\}$ is the inverse diffeomorphism. 
\smallskip

(3) {\it The inverse image 
$$\rho^{-1}(0)=\{0\} \times \PP^{n-1}(\R) \sim \PP^{n-1}(\R)$$ 
and it is in natural 
bijection with the set of lines in $\R^n$ passing through $0$; hence it is the
projectivization of $\R^n$ considered as vector bundle over the $0$-dimensional
manifold $X=\{0\}$.} 

\noindent Every such a line 
$L$ has a parametric equation $x_i=a_it$, $i=1,\dots, n$. Consider $L'=\rho^{-1}(L\setminus \{0\})$.
$L'$ has parametric equations $x_i=a_it$, $t_i=a_it$, $t\neq 0$, $i=1,\dots , n$. As the $t$'s are homogeneous,
equivalently $L'$ is described by $x_i=a_it$, $y_i=a_i$, $t\neq 0$. These equations extend to define
the so called {\it strict transform} $\tilde L$ of $L$ in $\BB(\R^n,0)$, that is the closure of $L'$;
finally $\tilde L$ intersects transversely $\PP^{n-1}(\R)$ at the point $(a_1,\dots, a_n)$ and
$$L \to \tilde L \cap \PP^{n-1}(\R)$$ defines the required bijection (after all, it corresponds to the bijection
between $x\in \PP^{n-1}(\R)$ and the respective fibre in the tautological bundle over $\PP^{n-1}(\R)$).
\smallskip

(4) In a more qualitative cut-and-paste fashion, $\BB(\R_n,0)$ is obtained by
gluing along the boundary the closure of $\R^n \setminus D^n$ with $\BB(D^n,0)$
and this last can be identified with the {\it mapping cylinder} of the natural
degree two covering map $$c: S^{n-1}\to \PP^{n-1}(\R)$$ $S^{n-1}=\partial D^n$.

\smallskip

Consider now $\R^k\subset \R^{k+n}=\R^k\times \R^n$ (defined as usual by the equation $x_i=0$, $i>k$).
$\R^{k+n}=\R^k\times \R^n$ can be considered as the total space of the product vector bundle  over the manifold 
$X=\R^k$,
with fibre $\R^n$. Then define the {\it blowing up of $\R^{k+n}$ with centre $X=\R^k$} by 
$$ \BB(\R^{k+n},\R^k):= \R^k\times \BB(\R^n,0)$$
endowed with the restriction of the natural projection
$$\rho=\rho_{n,k}: \BB(\R^{k+n},\R^k)\to \R^{k+n} \ . $$
The above properties extend directly; set 
$$D_{n,k}=\rho^{-1}(\R^k)$$
then:

(1) The restriction
$$\rho: \BB(\R^{k+n},\R^k) \setminus D_{n,k} \to \R^{k+n}\setminus \R^k$$
is a diffeomorphism; 

(2) $D_{n,k}=\R^k\times \PP^{n-1}(\R)$
and it is the total space
of the projectivization of the above trivial vector bundle;  

(3) $\BB(D^n,\R^k)$ is the mapping
cylinder of the natural degree two covering map 
$$c: \R^k \times S^{n-1} \to \R^k \times \PP^{n-1}(\R) $$ 
and  $\BB(\R^{k+n},\R^k)$ can be obtained
by gluing  $\BB(D^n,\R^k)$ to the closure of $\R^{k+n}\setminus (\R^k \times D^n)$, along the boundary.

\noindent Moreover:

(4) If $\R^k\subset \R^{k+h} \subset \R^{k+n}$, $h<n$, then  
the closure in $\BB(\R^{k+n},\R^k)$ of $\rho^{-1} (\R^{k+h} \setminus \R^k)$
is equal to $\BB(\R^{k+h},\R^k)$, $\rho_{h,k}$ is the restriction of $\rho_{n,k}$,
$\BB(\R^{k+h},\R^k)$ intersects transversely $D_{n,k}$ at $D_{h,k}$.

(5) Given $\R^k\times \R^{s}\times \R^h$, then $\BB(\R^{k+s},\R^k)$ and
$\BB(\R^{k+h},\R^k)$ are disjoint submanifolds of $\BB(\R^{k+s+h}, \R^k)$.
Note that $\R^{k+s}\cap \R^{k+h}=\R^k \subset \R^{k+s+h}$, hence
$\R^{k+s}\cup \R^{k+h}$ is `singular' along $\R^k$. Blowing up with centre
$\R^k$ is a way to `desingularize' it.   
\medskip
 
Let $M$ be a compact boundaryless smooth $(k+n)$-manifold and $X\subset M$
a proper $k$-submanifold. We define the {\it blowing up of $M$ with centre $X$}
$$\rho=\rho_{M,X}: \BB(M,X)\to M$$
as follows: recall that a tubular neighbourhood 
$$\pi: U\to X$$ 
of $X$ in $M$
is by construction isomorphic to a neighbourhood fibred by $n$-disks of
the $0$-section (identified with $X$) of a rank $k$ vector sub-bundle
$$p: E\to X$$ 
of the restriction of $T(M)$ to $X$, such that 
$$\partial \pi: \partial U\to X$$
is isomorphic to 
to the unitary bundle 
$$up: UE\to X$$
with fibre $S^{n-1}$. There is a natural degree $2$ covering map 
$$c: UE\to \PP(E)$$
such that $up= \pp \circ c$. Then $\BB(M,X)$ is obtained by gluing
the mapping cylinder of this map $c$ 
 to the closure of $M\setminus U$, along its boundary.  The above 
 $(\BB(\R^{k+n},\R^k),\rho_{n,k})$
provides the local model for $(\BB(M,X), \rho_{M,X})$, so that
\smallskip

$\bullet$  $\BB(M,X)$ is a smooth compact $(k+n)$-manifold as well as $\rho$ is a smooth map;
\smallskip

$\bullet$ Denote by $D(M,X)=\rho^{-1}(X)$. Then   
the restriction $$\rho: \BB(M,X)\setminus D(M,X) \to M\setminus X$$ is a diffeomorphism;
\smallskip

$\bullet$   The restriction $\rho: D(M,X)\to X$ is isomorphic to the projectivized bundle $\pp: \PP(E)\to X$.

\begin{remark}{\rm
If $X$ is a hypersurface of $M$ ($\dim X= \dim M -1$), then $\rho: \BB(M,X)\to M$ is a global
diffeomorphism.}
\end{remark}

\bigskip

If $Y$ is a subset of $M$, the {\it strict transform} $\tilde Y$ of $Y$ in $\BB(M,X)$ is by definition
the closure of
$\rho^{-1}(Y\setminus X)$. Then we have: 

\smallskip

$\bullet$  Let $M$ be as above, $N\subset M$ a proper submanifold of $M$ and $X\subset N$
a proper submanifold of $N$ (whence of $M$). Then the strict trasform $\tilde N$ in $\BB(M,X)$
is a proper submanifold diffeomorphic to $\BB(N,X)$, moreover $\tilde N$ intersects transversely $D(M,X)$ at $D(N,X)$.

\smallskip

$\bullet$ If $N$ and $N'$ are proper submanifolds of $M$ which intersect transversely at $X=N\cap N' \neq \emptyset$,
then the strict transforms $\tilde N$ and $\tilde N'$ are disjoint in $\BB(M,X)$. Note that $N\cup N'$ is not
a submanifold of $M$ because it is `singular' along $X$. So by blowing up the singularity and taking the strict transforms
we can `desingularize' it.

\medskip

When $X=\{x_0\}\subset M$ is reduced to one point, blowing up $X$  is related to the connected sum.
We have
\begin{proposition}\label{BU-point} 
(1) If $\dim M=m$ is even, then $\BB(M,x_0)\sim M \# \PP^m(\R)$ (recall that $\PP^m(\R)$ is not orientable).

(2) If $M$ is oriented and $\dim M=m$ is odd, then:

(a) $\BB(M,x_0)$ is oriented
in such a way that the restriction $$\rho: \BB(M,x_0)\setminus D(M,x_0) \to M\setminus \{x_0\}$$
preserves the orientation;

(b) Let us stipulate that $S^{m}$ is oriented as the boundary of $D^{m+1}$ oriented by the standard orientation of $\R^n$,
and that
$\PP^{m} (\R)$ is oriented in such a way that the standard covering map $S^m \to \PP^m(\R)$
preserves the orientation; then
$$\BB(M,x_0)\sim M \# - \PP^m(\R)$$ where `$-$' indicates the opposite orietation and we are dealing
with the oriented connected sum.
\end{proposition}
\Dim Forget for a while orientation questions. By taking a chart $\sim \R^m$ of $M$ at $x_0\sim 0$,
we can assume that $D^m$ is a tubular neighbourhood of $x_0$. Recall that 
$\BB(D^m,0)\subset \R^m\times \PP^{m-1}$, this last endowed  with `mixed' coordinates $(x,t)$.
Let $z=(z_1,\dots, z_{m+1})$ be homogeneous coordinates on $\PP^m(\R)$, take the standard
affine chart $U = \{t_{m+1}\neq 0\}$; $U\sim \R^m$, with coordinate $y_1=(z_1/z_{m+1},\dots, y_m=z_m/z_{m+1})$.
 Then it is enough to prove that there is a diffeomorphism
$$\phi: \BB(D^m,0)\to \overline {\PP^m(\R)\setminus D^m}$$
which is the identity on $\partial D^m$. The diffeomorphism $\phi$ can be defined explicitely
as follows:
$$\phi(x_1,\dots, x_m, t_1,\dots , t_m)=(t_1,\dots, t_m, t(\sum_{j=1}^m x_j^2))\in \PP^m(\R)$$
where
$t= t_i/z_i$ if $z_i\neq 0$, $i=1,\dots, m$. The verifications that $\phi$ is well defined, its image is 
 $\overline {\PP^m(\R)\setminus D^m}$ and that it is a diffeomorphism are left as an exercise.
 Coming back to the orientation question: if $m$ is even then $\PP^m(\R)$ is non orientable,
 hence the connected sum with it is well defined. In the oriented case we easily check
 that $\BB(D^m,0)$ and $\overline {\PP^m(\R)\setminus D^m}$ induce opposite orientations
 on the common boundary $\partial D^m$. Hence the diffeomorphism $\phi$ reverses
 the orientation and (b) follows.

\cvd

\subsection{On complex blowing up}\label{complex-BU} The (complex) blowing up $\BB_\C(M,X)$ 
can be performed in the category of complex manifolds as well. At least the basic  $\BB_\C(\C^n,0)$
is defined by the very same formulas of $\BB(\R^n,0)$, in terms of complex coordinates.
Hence we can define the blowing up $\BB_\C(M,x_0)$ of a complex manifold at a point $x_0$. 
More generally,  il $M$ is an {\it oriented} $2n$-smooth manifold, 
$x_0\in M$ we can define
$\BB_\C(M,x_0)$ (up to oriented diffeomorphism)
by taking an oriented chart $\R^{2n}\sim \C^n$ at $x_0\sim 0$, perform $\BB_\C(D^{2n},0)$ and glue
it  to $\overline {M\setminus D^{2n}}$. We have
\begin{proposition} \label{C-BU-point} Let $M$ be a compact oriented $2n$-manifold, $x_0\in M$.
Then  $\BB_\C(M,x_0)\sim M \# - \PP^n(\C)$.
\end{proposition}
\Dim As in the proof of Proposition \ref{BU-point}, the key point is to construct a suitable
diffeomorphism 
$$\phi_\C: \BB_\C(D^{2n},0)\to \overline {\PP^n(\C)\setminus D^{2n}} \ . $$
In fact the formula that defines $\phi$ above works as well, provided that it is considered in terms of 
the complex coordinates and we replace each $x_j^2$ with $|x_j|^2$.

\cvd

\smallskip

\begin{remark}\label{Alg-BU} {\rm Blowing up works in the category of (real or complex) 
regular algebraic varieties. In fact algebraic geometry is the first source of this construction
and we have just developed a smooth version. Note that in the algebraic setting,
$\BB(M,X)\setminus D(M,X)$ is a Zariski open set of the regular algebraic variety $\BB(M,X)$
as well as $M\setminus X$ is a Zariski open set of the regular algebraic variety $M$;
the restriction of $\rho$  is an algebraic isomorphism between these Zariski open sets,
hence (essentially by definition) $M$ and $\BB(M,X)$ are {\it birationally equivalent }. 
$M$ is said to be {\it rational} if it is birationally equivalent to the projective space
of the same dimension. Blowing up a projective space along regular centres is a basic way
to construct rational varieties.}
\end{remark}

\chapter{Transversality}\label{TD-TRANSVERSE}

We have already employed some instances of transversality and related concepts. 
Here we will treat this topic more systematically. 
First we point out so called ``basic transversality theorems'' which to a large extent will suffice to
our aims. Then we will develop some complements.

\section{Basic transversality}\label{basic-transverse}
We consider the following setting.

$\bullet$ $M$ is a smooth $m$-manifold with (possibly empty) boundary $\partial M$;

$\bullet$ $N$ is a smooth boundaryless $n$-manifold and $Z\subset N$ is a proper $r$-submanifold of $N$, hence $Z$ is both boundaryless and
a closed subset of $N$;

$\bullet$ $f:M \to N$ is a smooth map. If the boundary is non empty, then $\partial f$ denotes the restriction of
$f$ to $\partial M$.

\begin{definition}\label{trans-defi}{\rm We say that $f$ is {\it tranverse to $Z$} (and we write $f \pitchfork Z$) if 
\begin{enumerate}
\item  For every $x\in M$
such that $y=f(x)\in Z$, we have 
$$T_yN= T_yZ + d_xf(T_xM) \ . $$

\item For every $x\in \partial M$ such that $\partial f(x)\in Z$, we have
$$T_yN= T_yZ + d_x\partial f(T_x \partial M) \  , $$
in other words, $\partial f \pitchfork Z$  by itself. Obviously, if
$\partial M = \emptyset$, then this second requirement is empty.
\end{enumerate}

We denote by $\pitchfork(M,N;Z)$ the subspace of $\Ee(M,N)$
formed by the maps transverse to $Z$. If $A$ is a subset of $M$
we denote by  $\pitchfork_A(M,N;Z)$ the space of maps which verify the transversality
conditions for every $x\in A$ or $\in A\cap \partial M$, so that  
$\pitchfork(M,N;Z)= \pitchfork_M(M,N;Z)$.}
\end{definition}

\medskip

{\it Some particular cases:}
\smallskip

- If $f(M)\cap Z = \emptyset$, then $f\pitchfork Z$ ;
\smallskip

- If $Z=\{y_0\}$ a single point then $f\pitchfork Z$ if and only if $y_0$ is a regular value of both $f$
and $\partial f$.
\smallskip

- If also $M$  is a boundaryless submanifold of $N$ and $f$ is the inclusion, then $f\pitchfork Z$
(and we write also $M\pitchfork Z$)
if and only if for every $x\in M \cap Z$, $T_xN = T_xM + T_xZ$; if $\dim M + \dim Z = \dim N$, then
$T_xN = T_xM \oplus T_xZ$.

- The basic local models for $M\pitchfork Z$, and in fact for the whole transversality stuff,
is given by the possible mutual position of two affine subspaces, say $A$ and $B$,
in some $\R^n$. If $\dim A + \dim B < n$, then $A\pitchfork B$ if and only if $A\cap B = \emptyset$.
If $A\cap B \neq \emptyset$, up to translation we can assume that they are linear subspaces which are transverse
if and only if $\R^n = A+B$. Note that $A\cap B$ is also a linear subspace and, by elementary
linear algebra
$\dim A\cap B = \dim A + \dim B - n \geq 0$.  

\medskip

There are two kinds of basic transversality theorems; roughly speaking, 
they respectively claim that transversality
implies nice geometric features of the map $f$,  and that (at least when $M$ is compact) it is  a {\it generic and stable}
property: up to arbitrarily small perturbation every map $f$ becomes transverse, and trasversality
cannot be destroyed by small perturbations. 

In the given setting we have: 
\begin{theorem}\label{firstT} {\rm (First basic transversality theorem)} 
(1) If $f:M \to N$ is transverse to $Z$, then $(Y,\partial Y):=(f^{-1}(Z),\partial f^{-1}(Z))$ 
is a proper submanifold of $(M,\partial M)$; moreover $\dim M - \dim Y = \dim N -\dim Z$.
\smallskip

(2) If $(M,\partial M)$, $N$ and $Z$ are oriented then $Y$ and $\partial Y$
are orientable and we can fix an orientation procedure in such a way that
$\partial Y$ is the oriented boundary of $Y$.  
\end{theorem}

\Dim When $Z=\{y_0\}$ consists of one points, then the theorem is equivalent
to  Proposition \ref {b-summ}. Let us reduce the general  to this special
case. As $Z$ is a closed subset then also $f^{-1}(Z)$ and $\partial f^{-1}(Z)$
are closed sets. Being a proper submanifold is a local property. For every $z\in Z$
there is a chart of $N$, $\phi: W\to U\times U' \subset  \R^r\times \R^{n-r}$,  such that 
$\phi(z)= 0\in U\times U'$, and $\phi(W\cap Z)= U\times \{0\}$. Let 
$p:U\times U' \to U'$ be the projection. Then it is easy to see that the restriction
of $f$ to $f^{-1}(W)$ is transverse to $Z$ if and only if $p\circ \phi \circ f$ is transverse to
$\{0\}$. This is enough to achieve point (1). As for the orientation, let us orient
$\R^{n-r}$ is such a way that the given orientation of $\R^n$ (i.e. of $N$) is the direct sum
of the given orientation of $\R^r$ (i.e. of $Z$) followed the selected one on $\R^{n-r}$.
Then we can apply to $p\circ \phi \circ f$ the orientation rule of point (2) of Proposition
\ref{b-summ} to orient the intersection of $(Y,\partial Y)$ with $W$; by construction
these local orientations are globally coherent. 

\cvd  

\medskip

\begin{remark}\label{orient-int}{\rm It is useful to make explicit the
orientation rule in the case of transverse intersection $M\pitchfork Z$ of submanifolds
of $N$. For every $x\in M\cap Z$, $T_xN = T_xM + T_xZ$, and by assumption the
linear spaces $T_xN$, $T_xM$ and $T_xZ$ are oriented (in a globally coherent way)
and the last two intersect tranversely in the first. We have to orient $T_xM\cap T_xZ$.
So we have reduced the problem to the basic situation of two transverse oriented 
linear subspaces $(A,\omega_A)$ and $(B,\omega_B)$ in $\R^n$ 
(endowed say with the standard orientation $\omega_n$). 
Given any orientation  $\omega_{A\cap B}$ on the intersection, it can be extented
in an unique way to $A$ and $B$ in such a way that 
$\omega_A= \omega_{A\cap B} \oplus \omega'$
and $\omega_B= \omega_{A\cap B} \oplus \omega"$. Then 
$\omega_{A\cap B}\oplus \omega' \oplus \omega"$ determines an orientation
on the whole $\R^n$. Finally we select the orientation $\omega_{A\pitchfork B}$
such that the orientation of $\R^n$ obtained so far coincides with the given $\omega_n$.
Note that in the non oriented setting $M\pitchfork Z = Z\pitchfork M$, but the orientation
depends on the order; this can be checked straighforwardly in the linear local model; 
we get
$$ M\pitchfork Z = (-1)^{(\dim N - \dim M)(\dim N - \dim Z)} Z\pitchfork M \ . $$
}
\end{remark} 
\medskip

A very important consequence of Theorem \ref{firstT} is the following
{\it parametric transversality theorem}. In a sense it represents the bridge
between the two kinds of transversality theorems. 
Keeping the above setting, consider 
furthermore a  boundaryless ``parameter'' 
smooth manifold $S$, so that $M\times S$ has boundary equal to $\partial M \times S$.
We have   
\begin{theorem}\label{parametricT} Let $F: M\times S \to N$ be transverse to $Z$.
For every $s\in S$, set $f_s:M\to N$ the restriction of $F$ to $M\sim M\times \{s\}$. Then
the set of parameters $s\in S$ such that $f_s$ is not transverse to $Z$ is negligible in $S$.
\end{theorem}

\Dim  Let $(Y,\partial Y)=(F^{-1}(Z), \partial F^{-1}(Z))$ be the proper submanifold
of $(M\times S, \partial M \times S)$ accordingly with Theorem \ref{firstT}.
Set $\pi: Y \to S$ the restriction to $Y$ of the projection $p:M\times S \to S$.
We claim
that for every regular value $s$ of both $\pi$ and $\partial \pi$ (i.e. such that $\pi \pitchfork \{s\}$), then 
$f_s$ is transverse to $Z$. The thesis will follow from the Sard-Brown theorem. Let us justify the claim.
Let $x\in M$ be such that $f_s(x)=F(x,s)=z\in Z$. As $F\pitchfork Z$, for every
$w\in T_zN$, there are $(u,v)\in T_xM \times T_sS$ and $t\in T_zZ$ such that
$$ w= d_{(x,s)}F(u,v)+t \ . $$
The differential 
$$ d_{(x,s)} p : T_xM \times T_sS \to T_sS$$
is just the projection onto the second factor, and $d_{(x,s)}\pi$ is obtained by restriction. 
As $s$ is a regular value of $\pi$,
then there exists a vector of the form $(u',v)\in T_{(x,s)}Y$. By definition of $Y$,
 $t':= d_{(x,s)}F(u',v) \in T_zZ$. Finally we readily verify that
 $$w=  d_{(x,s)}F(u-u',0)+  d_{(x,s)}F(u',v) + t =  d_xf_s(u-u') + (t' - t) \ . $$
 This proves that $T_zN = d_xf_s(T_xM)+T_zZ$. 
 By using that $s$ is also a regular value of $\partial \pi$,
 the very same argument shows that $\partial f_s \pitchfork Z$. This achieves the proof.
 
 \cvd

\medskip
 
To state the second transversality theorem, we refine the setting. That is we
assume furthermore that
\begin{enumerate}
\item $M$ is compact;
\item $N$ can be embedded in some  $\R^h$ being also a 
closed subset.
\end{enumerate}

In many application also $N$ and $Z$ will be compact. In any case these assumptions
allows to apply to $N$ the results of Section \ref{exhaustive}. In this refined
setting we have:

\begin{theorem}\label{secondT} {\rm (Second basic transversality theorem)} 
(1) The set $\pitchfork (M,N;Z)$ of smooth maps transverse to $Z$ is open
and dense in $\Ee(M,N)$.
\smallskip

(2) Let $f\in  \Ee(M,N)$ be such that $\partial f: \partial M \to N$ is transverse
to $Z$. Denote by $\Ee(M,N,\partial f)$ (resp. $\pitchfork (M,N,\partial f;Z)$) the space of smooth maps 
that coincide with $\partial f$ on $\partial M$ (and are transverse to $Z$). Then   $\pitchfork (M,N,\partial f; Z)$ 
is open dense in $\Ee(M,N,\partial f)$.
\smallskip

(3) For every $h\in \Ee(M,N)$ (resp. $h\in \Ee(M,N,\partial f)$) there is 
$g\in \ \pitchfork (M,N; Z)$ ($g\in \ \pitchfork (M,N,\partial f; Z)$)
smoothly homotopic to $h$. 
\end{theorem}

\Dim  Let us consider first the openess in both items (1) and (2). 
As $M$ is compact, in early chapters we have already achieved it in 
the case of summersions; this easily implies the Theorem when $Z=\{y_0\}$ 
consists of one point.
By using the local reduction argument to this case as in the proof of Theorem
\ref{firstT}, for every $f \in \ \pitchfork (M,N; Z)$, we can  
find a finite covering of $M$ by compact sets $K$ such that $f$ reduces 
to the special case on a neighbourhood of each $K$ in $M$. Then, for every
$K$,  there is a open neighbourhhood $\Uu_K$ of $f$ in $\Ee(M,N)$  
formed by maps which verify the transversality conditions at every $x\in K$.
Then the intersection of these finite family of open sets $\Uu_K$ is a open neighbourhhood of
$f$ in $\Ee(M,N)$ contained in $\pitchfork (M,N; Z)$; hence it is open. The same argument
applies to $\pitchfork (M,N,\partial f; Z)$.
\smallskip

Let us come now to the density stated in (1). We consider first the special case 
when $N=\R^n=\R^r \times \R^{n-r}$ and $Z= \R^r=\R^r\times \{0\}$. Let
$f\in \Ee(M,\R^n)$. Then clearly the map
$$F: M\times \R^n \to \R^n, \ F(x,s)= f(x)+s $$
is transverse to $\R^r$ (in fact it is a summersion onto the whole $\R^n$) 
and we can apply to it the parametric transversality
Theorem \ref{parametricT}. Then for every $\epsilon >0$ there is $s\in \R^n$
such that $||s||<\epsilon$ and $f_s\pitchfork Z$. As $M$ is compact, by taking  $\epsilon$
small enough, then $f_s=f+s$ can be arbitrarily close to $f$ in the $C^\infty$-topology.

We are going to apply the same argument in the general case, by means of
a more elabotare construction. Let $f\in \Ee(M,N)$. 
For the moment assume for simplicity that $N\subset \R^h$ is compact
and take a tubular neighbourhood $\pi_N:U_N \to N$ of $N$ in $\R^h$ constructed by means of
the standard riemannian metric $g_0$ on $\R^h$ and some $\epsilon_0 >0$. Consider the restriction
of the map defined above
$$ F: M \times B^h(0,\epsilon) \to \R^h, \ F(x,s)= f(x)+s \ . $$
The parameter space is now restricted to the open ball of ray $\epsilon$; as $M$ and $N$ are compact,
then if $\epsilon$ is
smalls enough, the image of $F$
is contained in $ U_N$ and we can define
$$\hat F: M\times B^h(0,\epsilon) \to N, \ \hat F(x,s)=\pi_N(F(x,s)) \ . $$
As both $F$ and $\pi_N$ are summersions, also $\hat F$ is a summersion, hence
$\hat F \pitchfork Z$, and we can apply again
Theorem \ref{parametricT}. For $s$ generic and small enough, $\hat f_s \pitchfork Z$
and is arbitrarily close to $f$. If $N$ is not compact, by using the considerations of 
Section  \ref{exhaustive}, there is a compact submanifold with boundary $N' \subset N$
such that $f(M)\subset {\rm Int} (N')$ and we can repeat the above argument by using 
a  tubular ``neighbourhood'' $\pi_{N'}: U_{N'}\to N'$. Alternatively, we can use (instead of
$\pi_N$) the
projection $\pi: N_\epsilon \to N$ defined on the $\epsilon$-neighbourhood 
of $N$ determined by a suitable smooth positive function $\epsilon: N \to \R$,
and the modified maps 
$$\hat F: M \times B^h(0,1)\to N, \ \hat F(x,s)=\pi(f(x)+\epsilon(x)s) \ . $$
\smallskip

Let us face now the density stated in (2). We follow the same scheme,
by suitably modifying the map $\hat F$. Let $f\in \Ee(M,N)$ be such that
$\partial f \pitchfork Z$. By using the same consideration developed
to prove the openess, it is easy to verify that $f\pitchfork Z$ provided
that it is restricted to a small collar $C$ of $\partial M$. By slightly
modifying the construction of a collar bump function, we can construct
a smooth function $\gamma: M\to [0,1]$ such that $\gamma$ is constantly
equal to $0$ on a smaller closed collar $C'\subset C$, $\gamma$ is positive on the complement of $C'$
and constantly equal to $1$ outside $C$. Again assume for simplicity that $N\subset \R^h$
is compact  and let $\pi_N: U_N\to N$ as above.
Then define
$$\hat F: M\times B^h(0,\epsilon) \to N, \ \hat F(x,s) = \pi_N(f(x) + \gamma^2(x)s) \ . $$
We claim that $\hat F\pitchfork Z$, then for generic $s$ small enough, $\hat f_s = \pi_N\circ (f+\gamma^2s)$
belongs to $\pitchfork (M,N,\partial f; Z)$ and is arbitrarily close to $f$. We can complete the discussion
to $N$ non compact as above. It remains to justify that $\hat F \pitchfork Z$.
The restriction of $\hat F$ to $\{x; \ \gamma^2(x)\neq 0\}\times B^h(0,\epsilon)$
is a summersion because  for for every fixed $x$,  $s \to \gamma^2(x)s$ is a diffeomorphism
onto its image, the map $F(x,t)=\pi_N(f(x)+t)$ is a summersion, and $\hat F$ is obtained by composition.
It follows that if $\hat F(x,s)=z\in Z$ and $\gamma^2(x)\neq 0$, then the transversality conditions
are verified at $(x,s)$. Assume now that $\hat F(x,s)=z\in Z$ and $\gamma^2(x)= 0$, that is $x\in C'$.
We note that $d_x \gamma^2 = 2\gamma(x)d_x\gamma$, hence  it vanishes on $C'$.
By using this fact it is not hard to verifies that for every $(u,v)\in T_x(M)\times T_sB^h(0,\epsilon)$,
$$ d_{(x,s)}\hat F(u,v) = d_xf(u)$$
hence these differentials have the same image in $T_zN$. As $f$ restricted to $C'$ is tranverse to $Z$,
then also the restriction of $\hat F$ to $C'$ is transverse to $Z$.
\smallskip

Concerning point (3), referring for instance to $f\in \Ee(M,N)$ and to the above proof of the density,
we note that $f= \hat f_0$, and  it is homotopic to $\hat f_s \pitchfork Z$ via
the path $\hat f_{\sigma (t)}$, $\sigma (t)=  (1-t)s$, $t\in [0,1]$. On the other hand, we know
in general that if $g$ is close enough to $f$, then they are homotopic (recall Lemma \ref{near-homotopic}).

The proof is now complete.

\cvd

\smallskip

\begin{remark}\label{non-emb-density} {\rm The proof of the openess does not use that
$N$ is embedded. We sketch here an ``abstract'' (similar) proof of the density of (1) and (2) 
in Theorem \ref{secondT} too.
For semplicity we consider statement (1) and assume that $M$ is boundaryless 
(we left to the reader the task to
adapt the discussion to the other situations). Let $f\in \Ee(M,N)$.
By compacteness of $M$ there is a nice atlas $\Nn$ of $M$
$$\{ \phi_j: W_j \to B^m(0,1)\}_{j=1,\dots, s}$$ 
and a family $\Ff$ of charts of $(N,Z)$ of the form
$$\{\alpha_j: (V_j, Z\cap V_j)\to (\R^a\times \R^{n-r}, \R^a )\}$$ 
such that for every $j$, $f(W_j)\subset V_j$ so that we have the
the family 
$$\{ f_j:  U_j \to \R^r\times \R^{n-r} \}  $$
of associated representations of $f$ in local coordinates
supported by $(\Nn,\Ff)$. Recall that every $K_j=\overline B_j \subset W_j$ is compact and this provides a finite compact
covering of $M$. The  subset $A_{\Nn,\Ff}$ of $\Ee(M,N)$ formed by the maps admitting local representations
supported by  $(\Nn, \Ff)$ is open and non empty as it contains $f$. By applying to every $\Ee(W_j,V_j)$
the special case of the density considered in the proof of Theorem \ref{secondT} (1) and by using the bump function $\gamma_j$ 
in order to extend locally defined maps to maps in $\Ee(M,N)$, we realize that for every $j$,
$\pitchfork_{K_j}(M,N;Z)\cap A_{\Nn,\Ff}$ is dense In $A_{\Nn,\Ff}$. We know that it is also open. 
Then the intersection of these finite family of open and dense sets is open and dense in $A_{\Nn,\Ff}$,
and is contained in $\pitchfork (M,N;Z)$ because the $K_j$ cover the whole of $M$ .

\cvd}   
\end{remark}

\begin{remark}\label{generic}{\rm The meaning of the transversality theorems has been precised.
We have already recalled  that for example {\it any} compact subset $K\subset B^m(0,1)$
can be realized as $K=f^{-1}(0)$ for some smooth function $f: \overline B^m(0,1) \to \R$; compared with the
tame behaviour  of $K$ when $f\pitchfork \{0\}$, this shows that non transverse situations can be really
weird. On the other hand, by Theorem \ref{secondT} remarkably any weird non transversal situation can be made
stably tame up to arbitrarily small perturbations (at least when $M$ is compact).}
\end{remark} 

\section{Miscellaneous transversalities}\label{miscellaneaT}
Transversality is a profound, potent and pervasive paradigm beyond the basic results stated
in the previous section. Without any pretention of completeness we collect here a few instances
of further applications.

\subsection{Jet trasversality}\label{jet}
First we perform some constructions within the smooth category of open sets considered in Chapter
\ref{TD-LOCAL}. In particular, we refer to the Taylor polynomials defined in Section
\ref{diff-calc}.  Recall that a {\it homogeneus polynomial maps of degree $k\geq 1$} 
$$\pG: \R^m \to \R^n$$
is of the form $\pG(x)=\phi(x,\dots,x)$, where $\phi: (\R^m)^k\to \R^n$ is a (necessarily unique)
{\it symmetric} $k$-linear map. The set $\Pp_k(m,n)$ of these homogeneus polynomial maps 
has a natural structure of finite dimensional  real vector space endowed with a standard basis
so that it is identified with $\R^{\dim \Pp_k(m,n)}$.
A {\it polynomial map of degree} $\leq r$, $p:\R^m \to \R^n$, is of the form 
$$p= p_0+p_1+\dots + p_r$$
where $p_0\in \R^n$ and for $k\geq 1$, $p_k$ is homogeneous polynomial map of degree $k$.
Denote by $J^r(m,n)$ the set of these polynomial maps.
We can use the natural identification
$$J^r(m,n) = \prod_{k=0}^r \Pp^k(m,n)$$
to give it a finite dimensional real vector space structure
and $J^r(m,n)$ is identified with $\R^{\dim J^r(m,n)}$.

\begin{remark}{\rm With some effort
one can compute the dimension: 
$$\dim J^r(m,n) = n\binom{r+m}{n} \ . $$
}
\end{remark}

\smallskip

Let $U\subset \R^m$, $V\subset \R^n$ be non empty open sets. 
Then we can define the open set of $\R^m \times J^r(m,n)$ by 
$$J^r(U,V):= \{(x,p)\in U\times J^r(m,n); \ p_0\in V\} \ . $$
Given a smooth map $f:U\to V$, define the smooth map
$$j^rf: U \to J^r(U,V), \ j^rf(x)= \Tt_rf(x) $$
sending every point of $U$ to the Taylor polynomial of $f$ at $x$ of degree
$\leq r$.
\smallskip

{\it (Composition rule)} Let $U\subset \R^m$, $V\subset \R^n$ and
$W\subset \R^h$ be non empty open sets.  Set
$$J^r(U,V,W)=\{((y,q), (x,p))\in J^r(V,W)\times J^r(U,V); \ p_0=y\} \ . $$
Let $f:U\to V$, $g: V\to W$ be smooth maps. By a suitable extension to
higher order derivatives of the chain rule, one can find an unique polynomial
map (the explicit expression is called {\it Faa di Bruno formula})
$$\PG^r: J^r(U,V,W)\to J^r(U,W)$$ 
such that
$$j^r(g\circ f)(x) = \PG^r(j^rg(y),j^rf(x)) \ . $$
\smallskip

As a particular application of the composition rule we have:
\smallskip

{\it (Change of coordinates)} Let $U, U' \subset \R^m$, $V,V' \subset
\R^n$ be non empty open sets; $\phi: U\to U'$, $\psi: V\to V'$ be
diffeomorphisms. Then for every $r$, there is a unique smooth
diffeomorphism
$$j^r_{\psi,\phi}:J^r(U,V)\to J^r(U',V')$$
such that
$$j^r_{\psi,\phi}(j^rf(x))=j^rf'(x')$$
where
$$x'=\phi(x), \ f'=\psi\circ f \circ \phi^{-1} \ . $$

Now we can globalize the above local considerations, extending what we have done for
the (co)-tangent map.

Let $M$, $N$ be smooth manifolds of dimension $m$ and $n$
respectively. Define on $M\times \Cc^\infty(M,N)$ the following relation:
\smallskip

$(x,f)\sim_r (x',f')$ if $x=x'$, $f(x)=f'(x)$ and there are {\it compatible
representations in local coordinates} of $f$ and $f'$ at $x=x'$, $y=f(x)$, that
is defined on the same charts of $M$ and $N$ respectively:
$$ f_{U,V}, f'_{U,V}: U\to V, \ U\subset \R^m, \ V\subset \R^n$$
such that
$$j^rf_{U,V}(x)=j^rf'_{U,V}(x) \ . $$

\smallskip

By using the change of coordinates rule, it is easy to check that this
defines an equivalence relation and that if $(x,f)\sim_r (x',f')$,
then the above defining property holds {\it for every} pair of
compatible representations in local coordinates. We denote the
equivalence class of $(x,f)$ by $j^rf(x)$ and it is called the {\it
  $r$-jet} of $f$ at $x$, $J^r(M,N)$ is the {\it space of $r$-jets}
from $M$ to $N$. For every smooth map $f:M \to N$, the map
$$ j^rf: M \to J^r(M,N)$$ is called the {\it $r$-jet extension} of
$f$.  Clearly $J^0(M,N)=M\times N$. For every $r\geq 1$, $J^r(M,N)$
has a natural structure of smooth manifold of dimension
$$\dim J^r(M,N)=\dim M + \dim J^r(m,n) \ . $$
Local coordinates $U$ and $V$ for
$M$ and $N$ carry local coordinates $J^r(U,V)$ for $J^r(M,N)$.
This provides a smooth atlas
of $J^r(M,N)$ and we have already settled the change of coordinates
rules.  We see above also the local representations of an extension
$j^rf$, which is a smooth map indeed.  There is a natural smooth
projection
 $$\sigma_r: J^r(M,N)\to M $$ and a sequence of smooth``forgetting"
maps which factorize $\sigma$:
 $$ M \leftarrow J^1(M,N)\leftarrow \dots \leftarrow J^r(M,N) \ . $$
The map $\sigma_r$ is a smooth fibration with fibre {\it
  diffeomorphic} to $J^r(m,n)$; note that in spite of the fact that
$J^r(m,n)$ is a vector space with a preferred basis, for $r>1$
$\sigma_r$ is {\it not} a vector bundle.  The atlas of $J^r(M,N)$ is
fibred but the changes of coordinates do not preserve the linear
structure of the fibre. Every jet extension $j^rf:M \to J^r(M,N)$ is a
section of such smooth fibre bundle. Also every map $J^s(M,N)\to
J^{s-1}(M,N)$ is a smooth fibration with fibre $\Pp^s(m,n)$.
 \smallskip
 
 We are ready to state a version of the so called {\it jet
   transversality theorem}.  Let $M$ and $N$ be smooth boundaryless
 manifolds and $Z$ be a submanifold of $J^r(M,N)$. Denote
 by $$\pitchfork j^r(M,N,Z)$$ the set of smooth map $f\in \Ee(M,N)$
 such that $j^rf\pitchfork Z$. We have:

 \begin{theorem}\label{jet-transv}
   Let $M$ be a compact smooth boundaryless manifold and $N$ be a
   boundaryless proper smooth submanifold of some $\R^h$.  Let $Z$ be
   a proper submanifold of $J^r(M,N)$. Then $\pitchfork j^r(M,N,Z)$ is
   open and dense in $\Ee(M,N)$.
\end{theorem} 

\Dim We limit to an outline. Note that also $J^r(M,N)$
can be embedded as a proper submanifold of some $\R^k$.  When $r=0$,
Theorem \ref{jet-transv} incorporates (1) of Theorem \ref{secondT} (at
least when $M$ is boundaryless). Openness is not hard.  As for the density, Theorem
\ref{secondT} ensures that every $j^rf$ can be approximated by a
smooth map $g:M\to J^r(M,N)$ transverse to $Z$, but the statement of
theorem \ref{jet-transv} requires furthemore that $g$ is the $r$-jet
extension of some map $\tilde f: M\to N$. So jet-transversality is not
an immediate consequence of standard transversality.  Nevertheless the
structure of the proofs is basically the same.  A first, fundamental
case to deal with is when $N=\R^n$.  In the proof of Theorem
\ref{secondT} the key point was the application of parametric
transversality to the deformations of a given map $f:M\to \R^n$ of the
form $f+s$, $s\in \R^n$. In the present situation the main difference
consists in using polynomial deformations of the form $f+p_0+p_1+\dots
+ p_r$, where $p=p_0+\dots + p_r$ varies among the polynomial maps $p:
\R^n \to \R^n$ of degree $\leq r$.  Provided this new ingredient, the
proof theorem \ref{secondT} can be repeated with minor changes.

\cvd

\smallskip

\subsection{Transversality to stratifications}\label{stratification}
In several situations it is convenient to extend the notion of
``general position" (i.e.  of transversality) with respect to suitable
``stratification'' either of $N$ for the standard transversality or of
$J^r(M,N)$ for jet-transversality.  We do not intend to present here a
consistent treatment of stratification theory. We limit to a few
suggestion. At a first sight a stratification of a smooth manifold
$X$ is a partition $\Ss=\{S_j\}$ by means of boundaryless, connected
{\it not necessarily proper} smooth submanifolds of $X$, called the
{\it strata} of the stratification.  In fact one usually requires
more; reasonable requirements are:
\begin{itemize}
\item The stratification is locally finite;

\item {\it (Frontier condition)} The frontier $\bar S_j \setminus S_j$
of every stratum $S_j$ is union of strata of strictly lower dimension;
  
\item For every $0\leq s \leq \dim X$, denote by $X^s$ the $s$-{\it skeleton} of the stratification
  that is the union of strata of dimension less or equal to $s$. Then $X^s$ is a closed subset of $X$.

\end{itemize}

\smallskip

For example if $S\subset X$ is a boundaryless proper submanifold, then
$\{X\setminus S, S\}$ is a stratification of $X$; the open simplices
of a smooth triangulation of $X$ as it is described in Section
\ref{comb-chi} form a stratification; in this case  every stratum of dimension
geater or equal to $1$ is not a proper submanifold.
\smallskip

Given a stratification $\Ss$ of $N$, denote by $\pitchfork(M,N,\Ss)$
the subspace of $\Ee(M,N)$ formed by the map $f:M\to N$ wich are
transverse to every stratum of $\Ss$ (we write $f\pitchfork \Ss$).
Similarly, for every $r\geq 1$, given a stratification $\Ss$ of
$J^r(M,N)$ we define $\pitchfork j^r(M,N,\Ss)$.
\smallskip

{\it (Nice stratifications)} We define this notion in a quite implicit
way. Assume $N$ satifies the hypotheses of Theorem
\ref{jet-transv}. We say that a stratification $\Ss$ as above is
    {\it nice} if for every compact boundaryless smooth manifold $M$,
    $\pitchfork(M,N,\Ss)$ (resp. $\pitchfork j^r(M,N,\Ss)$) is {\it
      open and dense} in $\Ee(M,N)$ and, moreover, for every such a map
    $f$ transverse to $\Ss$, $f^{-1}(\Ss)$ (resp. $j^rf^{-1}(\Ss)$) is
    a nice stratification of $M$.
    \smallskip
    
    A key question is to determine further explicit  (as mild as possible)
    conditions in order that a stratification $\Ss$ is nice.
    Roughly speaking such conditions should imply that the transversality
    to any stratum $S_j$ forces at least locally at $S_j$ the transversality
    to every stratum $S_i$ such that $S_j$ is in the frontier of $S_i$. We will
    not face this rather deep question (see also \cite{Wall2}). We
    limit to state some results where the stratifications are
    nice, without justifying this fact.

\subsection{A classification of map singularities}\label{sing-class}
An important field of application of jet-transversality (in the
stratified extension) is the study of {\it singularities of smooth
  maps} (see \cite {A2}). The idea is that, under suitable hypotheses,
for a ``generic'' map $f: M\to N$, the source manifold $M$ carries a
nice stratification such that the increasing codimension of the strata
corresponds to more and more `deep'  classes of singular points of $f$
determined by a certain specific {\it lack of transversality}.
Moreover, the occurrence of such singular points cannot be eliminated
by means of small perturbations of the map.

\smallskip

{\it (Classification by the differential rank)} A first coarse
classification is in term of the rank of differentials.  Let
$M$ and $N$ be boundaryless manifolds.  Let $f:M\to N$ be a smooth
map. A point $x\in M$ is said of {\it class $\Sigma^i$ (with respect
  to $f$)} if $\dim \ker d_xf = i$. For every $i$, denote by
$\Sigma^i(f)$ the subset of $M$ of points of class $\Sigma^i$.  They
form a partition of $M$. If $f$ is arbitrary this partition might be
weird. However we have:

\begin{proposition}\label{partition-rank}
  Let $M$ and $N$ verify the hypotheses of Theorem \ref{jet-transv}.
  Then there is an open dense set $\Rr$ in $\Ee(M,N)$ such that for
  every $f\in \Rr$, the connected components of the $\Sigma^i(f)$'s
  form a nice stratification of $M$.  Moreover, every $\Sigma^i(f)$ is
  a submanifold of $M$ of dimension given by
$$ \dim M - \dim \Sigma^i(f) = (\dim M-r)(\dim N -r), \ r= \dim M - i \ . $$  
\end{proposition}

\smallskip

In fact one defines a suitable nice stratification $\Ss_\Sigma$ of
$J^1(M,N)$ and for every generic $f$ we consider $\Rr=
(j^1f)^{-1}(\Ss_\Sigma)$. In local coordinates $J^1(U,V)$, $\Ss_\Sigma$
corresponds to the stratification of the matrix space $M(n,m,\R)$ by
the matrix rank.

\cvd

\smallskip

\begin{example}\label{Whitney ex}{\rm
(1) If $N=\R$, then $J^1(M,\R)$ is naturally identified with the
    cotangent bundle and $df=j^1f$; $f$ is a Morse function if and
    only if $j^1f$ is transverse to the zero section of the
    bundle. Hence the result about Morse function (at least when $M$
    is compact boundaryless) of Chapter \ref{TD-COMP-EMB} can be
    reobtained as a special case of jet-trasversality.
    \smallskip

(2) {\it (Whitney fold)} Consider the map $f:\R^2 \to \R^2$
    defined by $$f(x_1,x_2)=(x_1^3+x_1x_2, x_2)\ . $$ The set of singular
    points is the parabole $S:=\{3x_1^2+x_2=0\}$. The nice
    stratification of the source $\R^2$ is given by $\Sigma^0(f)=\R^2
    \setminus S$, $\Sigma^1(f)=S$.
    \smallskip

    (3) {\it (Whitney umbrella)} See also Section
    \ref{n-in-2n-1} . Consider the map $f:\R^2 \to \R^3$ defined by
    $$f(x_1,x_2)=(x_1x_2, x_2, x_1^2) \ . $$ The point $0\in \R^2$ is the
    only one at which $f$ is not an immersion; hence the nice
    stratification of $\R^2$ is given by $\Sigma^0(f)=\R^2 \setminus
    \{0\}$, $\Sigma^1(f)=\{0\}$.  }
\end{example}

\smallskip

The above examples show that the stratification by the differential
rank is in general too coarse.  In the Whitney fold, $0\in
\Sigma^1(f)=S$ is clearly special: $\ker d_0f=T_0S$ while for other
$x\in S$, $\R^2 = \ker d_xf + T_xS$. In the Whitney umbrella a
refinement of the startification can be rather obtained by noticing
that the line $\{x_2=0\}$ is the locus where the map is not injective.
\smallskip

If $f\in \Rr$ as in Proposition \ref{partition-rank}, a tentative
refinement of the stratification $\{\Sigma^i(f)\}$ would be defined by
recurrence as follows: assume that for every multi-index of length
$k$, $I=(i_1,\dots, i_k)$ is defined $\Sigma^I(f)\subset M$, then for
every multi-index of length $k+1$, $\tilde I = (i_1,\dots, i_k,
i_{k+1})$, set $\Sigma^{\tilde I}(f):=
\Sigma^{i_{k+1}}(f|\Sigma^I(f))$. It is not evident that this
eventually produces a nice (sub) stratification. The correct way to do
(see \cite {Bo}) is to extend the above stratification $\Ss_\Sigma$ 
of $J^1(M,N)$, to get a nice stratification $\tilde \Ss_{\Sigma}$ and
extend Proposition \ref{partition-rank}.

\subsection{Multi-transversality} Assume for example that $f:M\to N$ is an immersion,
and for simplicity $M$ is connected.  Then the nice stratification of $M$
consists of one stratum $\Sigma^0(f)=M$. This does not give any
information about the image of $f$. Clearly this last might be ``non
generic''. We say that an immersion $f$ is {\it in general position} if for every
$k\geq 2$, whenever $y=f(x_1)=f(x_2)= \dots = f(x_k)$ and the points
$x_1,\dots, x_k$ are distinct, then
$$T_yN = df_{x_k}(T_{x_k}M)+ \cap_{j=1}^{k-1} df_{x_j}(T_{x_j}M) \ . $$
For example if $\dim N = 2\dim M$, then the multiple points $y$ are isolated
and are image of exactly two points of $M$.
The following is a basic example of multi-transversality result.

\begin{proposition}\label{generic-imm} Let $M$, $N$ verify the hypotheses of Theorem
  \ref{jet-transv}. Assume that the open set Im$(M,N)$ of immersions of $M$ in $N$
  is non empty. Then the set of immersions in general position is open and dense in
  Im$(M,N)$.
\end{proposition}

The general concept of multi-jet-transversality was introduced in
\cite{Ma2}.  One considers the products $J^r(M,N)^k$, $k\geq 1$. Then
for every $f\in \Ee(M,N)$ we have the product map $(j^rf)^k: M^k \to
J^r(M,N)^k$. So for every submanifold $V$ of $J^r(M,N)^k$ we can
consider $f$ such that $(j^rf)^k\pitchfork V$.  The submanifolds $V$
of most interest are as follows:
\smallskip

- Given submanifolds $V_i$ of $J^r(M,N)$, $1\leq i \leq k$, consider
the product $\prod V_i \subset J^r(M,N)^k$;
\smallskip

- There is a natural projection $\tau_k: J^r(M,N)^k \to N^k$. Then
take $$V= \tau_k^{-1}(\Delta_k N)\cap \prod V_i$$ where
$\Delta_k(N)=\{(y,\dots, y)\}\subset N^k$ is the (multi) diagonal of
$N^k$.
\smallskip

Multi-transversality to such a manifold $V$ means that the following
conditions are satisfied:
\smallskip

- If $f(x_i)=y$ for every $i=1,\dots k$, then $f$ is transverse to
$V_i$ at $x_i$ with pre-image say $X_i = f^{-1}(V_i)$;
\smallskip

- The images say $B_i$ of $T_{x_i}X_i$ in $T_y N$ satisfy
$$ (\oplus_i B_i)\oplus T_{(y,\dots, y)}\Delta_k(N) = (T_yN)^k \ . $$ 

Finally, in the same hypotheses, one gets \cite{Ma2} a multi-transverse
version of Theorem \ref{jet-transv}.

\chapter{Morse functions and handle decompositions}\label{TD-HANDLE}
Let us call  {\it smooth triad} a triple $(M,V_0,V_1)$ where
$M$ is a compact smooth $m$-manifold with (possibly empty) boundary, $V_0$ and 
$V_1$ are union of connected components of $\partial M$, so that
the boundary is the {\it disjoint union} $$\partial M = V_0 \amalg V_1 \ . $$
A boundaryless $M$ corresponds to the triad $(M,\emptyset,\emptyset)$.
We stress that different ordered bipartitions of the components of $\partial M$
give rise to different triads. For example if $\partial M \neq \emptyset$, then  $(M,\partial M, \emptyset)$ 
and $(M,\emptyset, \partial M)$ are different triads. 
We know from Proposition \ref{partial-morse-dense} that generic Morse
functions form a dense open set in $\Ee(M,V_0,V_1)$, the space of functions
$f:M\to [0,1]$ such that $V_j=f^{-1}(j)$, $j=0,1$ and without critical
points on a neighbourhood of $\partial M$.
Let $f:M\to [0,1]$ be such a generic Morse function on the triad $(M,V_0,V_1)$.
We have a finite set of non degenerate critical points $p_0,\dots, p_s$
of indices $q_0,\dots, q_s$, and critical values $c_r=f(p_r)$, such that
$0<c_r < c_{r+1}<1$, $r=0,\dots, s-1$. For every $X\subset [0,1]$,
denote $V_{X}:=f^{-1}(X)$. For every regular value $a$ of $f$,
$V_a$ is a compact boundaryless submanifold of $M$ of dimension $m-1$.
If $0\leq a <b \leq 1$ are regular values, then we have the subtriad $(V_{[a,b]},V_a,V_b)$. 

The following lemma ultimately is an instance of a fibration theorem.
We give a ``non embedded'' proof by assuming a few results of analysis about the existence, 
the uniqueness and the regular dependence on the data for ordinary differential equations. 

\begin{lemma}\label{cylinder}  {\rm (Cylinder Lemma)}
Assume that $[a,b]\subset [0,1]$ does not contain any critical value of $f$.
Then there is a diffeomorphism $\psi: V_a \times [a,b]\to V_{[a,b]}$ such that 
$f\circ \psi(y,t)=t$  for every $y\in V_a$.
\end{lemma}
\Dim Fix an auxiliary riemannian metric $g$ on $M$ and let $\nabla_gf$ the associated
gradient field of $f$, which is non zero everywhere on $V_{[a,b]}$. We can normalize it
by taking for every $p\in V_{[a,b]}$, 
$$\nu(p)= \nabla_g f(p)/||\nabla_g f(p)||_{g(p)} \ . $$
Every integral curve $\alpha$ of $\nu$ verifies $f(\alpha(s))=s + c$, $c$ being
a constant. Possibly by means of the change of parameter 
$\beta(t)=\alpha(t-c)$, we can assume that $f(\alpha(t))=t$. Since $V_{[a,b]}$ is compact
every maximal integral curve is defined on the whole $[a,b]$. Then
for every $y\in V_a$ there is a unique maximal integral curve
of $\nu$ 
$$\alpha_y: [a,b]\to V_{[a,b]}$$ 
such that $\alpha(a)=y$, and  $f(\alpha(t))=t$
for every $t\in [a,b]$. The required diffeomorphism is defined 
by  $\psi(y,t)=\alpha_y(t)$, with inverse $\psi^{-1}(x)= (\alpha_x(a), f(x))$,
where $\alpha_x$ is the unique maximal integral curve passing through $x\in V_{[a,b]}$.

\cvd

\begin{remark}\label{abs-collar}{\rm Via the existence of embeddings of compact manifolds
in some $\R^n$, we  are currently exploiting the results obtained for compact embedded
manifolds. However, in several situations we could provide also an ``abstract'' treatment.
For example, the {\it existence} of collars of $\partial M$ in $M$ is an immediate
consequence of  Lemma \ref{cylinder}, provided that one knows that
$\Ee(M,\partial M, \emptyset)$ is non empty. This last fact can be obtained
as follows:  fix a nice atlas of $M$. Define local functions as follows: if $(W,\phi)$ is an internal chart,
then $f_j$ is the constant function equal to $1/2$. If 
$$\phi_j: (W_j,W_j\cap \partial M)\to (B^m\cap \HH^m, B^m\cap \partial \HH^m)$$ is a chart at the boundary,
then $f_j$ is the restriction of the projection of $B^m$ to the $x_m$-axis. By using
the partition of unity subordinate to the atlas, define $$f=\sum_j \lambda_j f_j \ . $$
One can check directly that $f$ has the desired properties.

Strictly speaking in Chapter \ref{TD-COMP-EMB} we have proved
the collars {\it uniqueness up to isotopy} only for the ones realized by means of that 
(embedded) construction.
In fact it holds  in full generality (see \cite{Mu}). However we do not really need this
fact, so we omit the somewhat technical proof. Also the density of Morse functions
can be obtained in an abstract way; the result about the generic linear projections to lines
gives us a ``local'' density for  representations in local coordinates, then one uses
nice atlas and partitions of unity to get the global result.}
\end{remark}

\section{Dissections carried by generic Morse functions}\label{triad-dissection}
First we fix a nice atlas with collars on the triad $(M,V_0,V_1)$ {\it adapted to
the given Morse function $f:M\to [0,1]$}. This means the following facts:
\begin{itemize}
\item The collars are of the form
$V_{[0,\epsilon_0]}$, $V_{[1-\epsilon_0,1]}$, for some $\epsilon_0>0$,
$\epsilon_0 < c_0=f(p_0)$, $c_s=f(p_s)< 1-\epsilon_0$;
\item every critical point 
$p_r$ of $f$ is contained in a unique internal normal chart $(W_r,\phi_r)$, in such a way that
$B_r\cap B_{r'} = \emptyset$ if $r\neq r'$ (recall that $B_r = \phi_r^{-1}(B^m(0,1/3)$);
\item  Every $(W_r,\phi_r)$ is such that $(f\circ \psi_r-c_r): B^m(0,1/3))\to \R$ is in normal form
according to Morse's Lemma of Section \ref{morse} at $0=\phi_r(p_r)$.
\end{itemize}
\smallskip

Certainly such an adapted atlas exists. Then we take $\epsilon>0$ such that
\begin{itemize}
\item  $\epsilon_0 < c_0-\epsilon$, $c_0+\epsilon < c_1-\epsilon$, $\dots$, $c_{s-1}+\epsilon < c_s-\epsilon$, $c_s+\epsilon < 1-\epsilon_0$;
\item for every $r=0,\dots, s$, $V_{c_r-\epsilon}\cap B_r\neq \emptyset$ and $V_{c_r+\epsilon}\cap B_r\neq \emptyset$, so that
$V{[c_r-\epsilon,c_r+\epsilon]}$ is the union  of  $V{[c_r-\epsilon,c_r+\epsilon]} \cap B_r$ and its complement.
\end{itemize}
\smallskip

So we have the {\it dissection} of the triad $(M,V_0,V_1)$ associated to the Morse function $f$:
$$ V_{[0,c_0-\epsilon]} \cup V_{[c_0-\epsilon,c_0+\epsilon]}\cup V_{[c_0+\epsilon, c_1-\epsilon]}
\cup V_{[c_1-\epsilon,c_1+\epsilon]}\cup \dots \cup V_{[c_s-\epsilon,c_s+\epsilon]}\cup V_{[c_s+\epsilon,1]}  \ . $$
By applying the cylinder and Thom's lemmas, we have that
\begin{itemize}
\item $V_{[0,c_0-\epsilon]} \sim V_0\times [0,c_0-\epsilon]$, $V_{[c_s+\epsilon,1]} \sim [c_s+\epsilon, 1]\times V_1$; 
\item for every $r=0,\dots, s-1$, 
$V_{[c_r+\epsilon, c_{r+1}-\epsilon]} \sim V_{c_r+\epsilon}\times [c_r+\epsilon , c_{r+1} - \epsilon]$;
\item $V_{[0,c_r+\epsilon]} \sim V_{[0, c_{r+1}-\epsilon]}$.
\end{itemize}
\smallskip

For every $r=0,\dots, s-1$, $(V_{[c_r-\epsilon, c_r+\epsilon]}, V_{c_r-\epsilon}, V_{c_r+\epsilon})$ is an 
{\it elementary triad} in the sense that it carries a Morse function
(the restriction of $f$) with {\it only one} critical point ($p_r$ of a given index $q_r$). 
\smallskip

{\bf Adapted gradient fields.} By using the above adapted nice atlas of $(M,V_0,V_1)$ with respect to $f$, we can construct an adapted
riemannian metric $g$ on $M$, so that for every $r=0,\dots ,s$, the gradient field $\nabla f:=\nabla_gf$ has the normalized expression in the local
coordinates over $B_r$:
$$2(-x_1,-x_2,\dots, -x_{q_r}, x_{q_r+1},\dots, x_m )$$
while the collars of $V_0$ and $V_1$ are obtained by integrating such a (normalized) field as in the proof of Lemma \ref{cylinder}. 
\smallskip

So the key point will be to understand what happens up to diffeomorphism by passing from $V_{[0,c_r-\epsilon]}$ to 
$V_{[0,c_r+\epsilon]}$ (equivalently to $V_{[0,c_{r+1}-\epsilon]}$)
through such an elementary triad. It is evident that  the choice of the 
parameters $\epsilon_0$ and $\epsilon$ is immaterial. An answer is given by the following Proposition. 
We refer to notions introduced in Chapter \ref{TD-CUT-PASTE}.
The  proof is extracted from \cite{Pa}. 
\begin{proposition}\label{attaching-q-handle}
  Let $f:M\to [0,1]$ be a generic Morse function on the triad
  $(M,V_0,V_1)$ and consider an associated dissection. Let $p$ be a
  critical point of $f$ of index $q$, and $c'$ be the next critical
  value of $f$ after $c=f(p)$. Then
\smallskip

(1) $V_{[0,c+\epsilon]}$ is diffeomorphic to $V_{[0,c'-\epsilon]}$. 
\smallskip

(2) Up to diffeomorphism, $V_{[0,c+\epsilon]}$ is obtained by
attaching a $q$-handle (of dimension $m$) to $V_{[0,c-\epsilon]}$
along $V_{c-\epsilon}$.
\end{proposition}
\Dim As already remarked, (1) follows from the Cylinder and Thom's lemmas.

As for (2), take a
nice atlas associated to the given Morse dissection of $(M,V_0,V_1)$.
Take the Morse chart $(\psi(B),\phi)$ at $p$, so that in that local
coordinates, $\phi(p)=0$, and $$\hat f= f\circ \psi : B\to \R$$ has
the normal form
$$ \hat f(x_1,\dots x_m)= -(x_1^2+\dots + x_q^2)+(x_{q+1}^2+\dots +
x_m^2)+c \ . $$
According to our usual conventions, $B$ should be
$B^m(0,1/3)$, but up to reparametrization we can normalize the picture
as follows. First we simplify the notations by setting
$$(x_1,\dots,x_q,x_{q+1},\dots x_m)= (X,Y)\in \R^q \times \R^{m-q} \ . $$
Then we can assume that: 

- $f:M\to [a_0, a_1]$ for suitable $ a_0 <-1$,  $1 < a_1$;

- $B= B^m(0,2)$,  $\hat f(0)=c=0$, $\epsilon = 1$;

- $B\cap \phi (W\cap V_{[a_0, -1]})= \{(X,Y)\in B; \ -||X^2||+||Y||^2\leq -1\}$;

- $B\cap \phi (W\cap V_{[a_0, 1]})= \{(X,Y)\in B; \ -||X^2||+||Y||^2\leq 1\}$.
\smallskip

The standard handle $H^q=D^q\times D^{m-q}$ is contained into 
$$B\cap \phi(W \cap V_{[-1, 1]})= \{(X,Y)\in B; \  -1 \leq  -||X^2||+||Y||^2\leq 1\}$$
and $H^q$ intersects  $\{ -||X^2||+||Y||^2=\pm1\}$ along the union of its $a$ and $b$-spheres.
Moreover, if $H'= (\R^q \times D^{m-q})\cap  \{ -1 \leq  -||X^2||+||Y||^2\leq 1\}$,
then $V_{[a_0,-1]}\cup \psi(H')$ is 
a submanifold with corners of $V_{[a_0,1]}$ obtained by attaching the $q$-handle to $V_{[0,-1]}$ 
along $V_{-1}$. The idea is to modify the inclusion of $H'$ to an embedding $j$ of $H^q$
(actually an embedded corner smoothing) in such a way that: 
\begin{enumerate}
\item $\Hh:=j(H^q) \subset \{(X,Y)\in B; \  -1 \leq  -||X^2||+||Y||^2 < 1\}$.
\item $\Hh \cap \{-||X^2||+||Y||^2=-1\} = j(\Tt_a)$, the image of the $a$-tube.
\item The embedding $j$ is still equal to the identity at the core of the handle. 
\item $\tilde M:=V_{[a_0,-1]}\cup \psi(\Hh)$ is 
a smooth submanifold of $V_{[a_0,1]}$ obtained by attaching the $q$-handle to $V_{[0,-1]}$ 
along $V_{-1}$, having the restriction of $j$ to $\Tt_a$ as attaching map.
 \item $V_{[a_0,1]}\setminus \tilde M$ is a collar of $V_{1}$ in $V_{[a_0,1]}$.
 \end{enumerate}

\noindent Take the $1$-dimensional bump function $\gamma= \gamma_{1/2,1}$; then define
$$\hat g: B\to \R; \ \hat g(X,Y)= -||X||^2+||Y||^2 -
\frac{3}{2}\gamma(||Y||^2) \ . $$
Clearly
$$ \{\hat g \leq -1\}= \{\hat f \leq -1\}\cup (\{\hat f \geq -1\}\cap \{\hat g \leq -1\}) :=
\{\hat f \leq -1\} \cup \Hh$$
and $\Hh$ intersects $\{\hat f \leq -1\}$ at $\{\hat f = -1\}$; $\{\hat g\leq -1\}$ is contained in the interior of $\{f\leq 1\}$, and 
  $\{\hat f\leq 1\} = \{\hat g\leq 1\}$.

{\bf Claim:} {\it $\Hh$ is  $q$-handle
  attached to $\{\hat f\leq -1\}$ along $\{\hat f=-1\}$, via a characteristic
  map $H:D^q\times D^{m-q}\to \Hh$ which is the identity on the core $D^q \times \{0\}$.}

We are going to write down the explicit formulas establishing the claim. Several verifications
are understood; for all details (in a more general setting) we refer to \cite{Pa}. 
The smooth function $\sigma: [0,1]\to \R$  is uniquely defined by the equation
$$ \frac{\gamma(\sigma(s))}{1+\sigma(s)}= \frac{2}{3}(1-s) \ . $$
The function $\sigma$ is strictly increasing, $\sigma(0)= \frac{1}{2}$, $\sigma(1)=1$
and moreover we have that for every $(X,Y)\in \Hh$,
$$ ||Y||^2< \sigma(\frac{||X||^2}{1+||Y||^2})\ . $$
By using $\sigma$ and its properties, we can give the explicit characteristic
map
$$H: D^q\times D^{m-q}\to \Hh$$

$$ H(X,Y)= (\sqrt{\sigma(||X||^2)||Y||^2 + 1}\ X, \sqrt{\sigma(||X||^2)}\ Y)$$
which restricts to the attaching map

$$h: S^{q-1}\times D^{m-q}\to \partial \Hh\subset \{\hat f =-1\} $$

$$h(X,Y)=(\sqrt{||Y||^2+1} \ X, Y) \ . $$

Let us consider now
$$ M':= [\{f\geq -1\}\cap (M\setminus \psi(B))] \cup \psi(\{(X,Y)\in B; \ \hat g\geq -1\}) \ . $$
By construction, the functions $f$ and $\hat g \circ \phi$ match on $M'$, giving us a global 
function $g: M' \to \R$, such that
$$ \{f \leq 1\} = \{f\leq -1\}\cup \psi(\Hh)\cup \{p\in M'; \ -1 \leq g \leq 1\} \ . $$
The final remark is that $[-1,1]$ does not contain critical values of $g$. It is enough to verify
it for $\hat g$ on $B$. In fact
$$ \nabla \hat g(X,Y)= 2(-X,Y) - 2(0, \gamma'(||Y||^2)Y)$$
which vanishes only at $0$ because $\gamma' \leq 0$ on $(0,+\infty)$.
Summarizing, as $$\{f\leq -1\}\cup \psi(\Hh)$$ is obtained by attaching a $q$-handle to
$\{f\leq -1\}$ along  $\{f= -1\}$, by applying the Cylinder Lemma to $g$ over $[-1,1]$
we conclude that also $\{f\leq 1\}$ is obtained by attaching a $q$-handle to $\{f\leq -1\}$
along  $\{f= -1\}$. Ultimately, by restoring the usual notations, $V_{[0,c+\epsilon]}$
is obtained by attaching a $q$-handle to $V_{[0,c-\epsilon]}$ along $V_{c-\epsilon}$.

\cvd

\begin{remark}\label{far-attaching}
  {\rm With the notations of (the proof of) Proposition
    \ref{attaching-q-handle}, we realize that the core $D^q\times \{0\}$
    of the
    $q$-handle $\Hh$ is formed by the integral lines of the
    adapted gradient field $\nabla f$ which start at a point of
    $V_{c-\epsilon}$ and end in the critical point $p$.  If
    $c-\epsilon > \delta>0$ is any value such that
    $[\delta,c-\epsilon]$ does not contain any critical value of $f$,
    then again by the Cylinder Lemma, $V_{[0,c+\epsilon]}$ is also obtained
    by attaching a $q$-handle say $\Hh'$ to $V_{[0,\delta]}$ along
    $V_{\delta}$. As well as the core of $\Hh$ and the relative
    attaching map $h$ look ``simple and local'', the core and the
    relative attaching map $h'$ of $\Hh'$ can be ``far from
    $V_{c-\epsilon}$ and complicated''.  In fact $h'$ is obtained by
    composing $h$ with the diffeomorphism between $V_{c'-\epsilon}$
    and $V_\delta$ provided by the Cylinder Lemma; again the core of
    $\Hh'$ is formed by the integral lines of the adapted gradient
    $\nabla f$ (used in the Cylinder Lemma) which start at a point of
    $V_\delta$ and end in $p$.}
 \end{remark}
 
 \section{Handle decompositions}\label{handle-decomp}
Let $(M,V_0,V_1)$ be a triad a before. By definition, a {\it handle decomposition} of the triad is
a sequence of nested triads of the form
$$ (M_0, V_0, V_{1,0})\subset (M_1,V_0,V_{1,1}) \subset (M_2, V_0, V_{1,2})\subset \dots \subset (M_k,V_0,V_{1,k})$$
such that 
\begin{itemize}
\item $V_{1,k}=V_1$, and $(M_k,V_0,V_1)$ is diffeomorphic to $(M,V_0,V_1)$  via a diffeomorphism which is the identity
in a neighbourhood of $V_0\amalg V_1$;
\item For every $r=0,\dots, k-1$, $(M_{r+1},V_0,V_{1,r+1})$ is obtained  
by attaching a $q$-handle (of dimension $m$) to $(M_{r},V_0,V_{1,r})$ along $V_{1,r}$ (for some $q$).
\end{itemize}
\smallskip

Two handle decompositions are diffeomorphic if they are related by a diffeomorphism which is the identity
near the boundary and respects the sequences of nested triads. 
We can also {\it normalize} the form of a given handle decomposition by stipulating that it starts with a ``right'' collar $C_0$
of $V_0$ and ends with a ``left'' collar $C_1$ of $V_1$.
\smallskip

As an immediate Corollary of Proposition \ref {attaching-q-handle} we have the {\it existence} of handle decompositions
for every triad.

\begin{corollary}\label{existence-handle} Every triad $(M,V_0,V_1)$ admits handle decompositions.
\end{corollary}
\Dim Take a dissection carried by any generic Morse function on the triad. The sequence of nested submanifolds
$$V_{[0, c_0-\epsilon]}\subset V_{[0,c_1-\epsilon]} \subset V_{[0, c_2-\epsilon]} \subset \dots \subset V_{[0,1]}$$
leads to a desired handle decomposition.

\cvd

Sometimes a handle decomposition of $(M,V_0,V_1)$ (in normalized form) is formally indicated as
$$C_0\cup H_1^{q_1}\cup H_2^{q_2}\cup \dots \cup H_k^{q_k} \cup C_1$$ 
where $C_0$ and $C_1$ are the respective collars of $V_0$ and $V_1$, and
for $r=0,\dots, k-1$, $M_r= C_0\cup H_1^{q_1}\cup H_2^{q_2}\cup \dots \cup H_r^{q_r}$,
$M_{r+1}$ is obtained by attaching the $q_{r+1}$-handle $H_{r+1}^{q_{r+1}}$ to $M_r$, along
$V_{1,r}$. Sometimes we will omit to indicate the index $q_r$.
\medskip

{\bf The dual decompositions.} Given a triad $(M,V_0,V_1)$, the {\it dual triad} is by definition
$(M,V_1,V_2)$.  Given a decomposition $\Hh$ of the triad $(M,V_0,V_1)$ formally indicated as
$$C_0\cup H_1^{q_1}\cup H_2^{q_2}\cup \dots \cup H_k^{q_k} \cup C_1$$
we can consider the dual decomposition $\Hh^*$ of $(M,V_1,V_0)$ obtained
by going from $C_1$ to $C_0$ in the opposite direction. Every $q$-handle $H^q$ of $\Hh$
is converted into a ``dual'' $(m-q)$-handle $(H^*)^{m-q}$ of $\Hh^*$ where the
core and the cocore exchange their roles. If $\Hh$ is associated to a Morse function $f$,
then $\Hh^*$ is associated to the function $f^*= 1-f$.
\medskip

Once we have obtained the existence of handle decompositions,
we will develop our discussion in terms of these last, somehow forgetting the Morse functions.
To this respect Morse functions have been rather a tool in order to produce handle decompositions. 
On another hand, one can prove that 
\smallskip

{\it For every handle decomposition of a triad, there is a Morse
function that recovers it, so that every $q$-handle corresponds to a critical point of index $q$.}
\medskip

So handle decompositions and Morse functions (with the associated dissections) {\it basically are equivalent
stuff}. This means that any manipulation in terms of handle decompositions should have a counterpart
in the realm of Morse functions. One can find such a purely Morse function approach in
\cite{M3}. However, dealing directly with handle decompositions is often easier and topologically transparent
with respect to its Morse function counterpart which can be demanding. Moreover, handle technology works 
as well even for other categories of manifolds (like the {\it piecewise-linear} (PL) one, see \cite{RS}) where there is 
{\it not} a Morse function counterpart. 
For these reasons we will not pursue the equivalence between Morse function and handle approaches,
preferring the latter. 

 \section{Moves on handle decompositions}\label{handle-moves} 
There are two basic ways to modify a given handle decomposition of a triad $(M,V_0,V_1)$ (up to diffeomorphism
equal to the identity on a neighbourhood of $V_0\amalg V_1$).

\medskip

{\it Handle sliding.} This is a synonymous of modifying the attacching map 
of a handle, say $H_r$,  in the decomposition staying in the same isotopy class.  
We have already noticed in Chapter \ref{TD-CUT-PASTE} that up to diffeomorphism
this does not modify $M_r$, then we can continue the decomposition by composing
the subsequent attaching maps with such a diffeomorphism, finally obtaining
a decomposition diffeomorphic to the given one (possibly by attaching a final collar of $V_1$
in order to normalize the form).

\smallskip

Before stating the other modification, let us give a definition. 

\begin{definition}\label{comp-handles} {\rm Let 
$$ \dots \cup H_r^{q_r}\cup H_{r+1}^{q_{r+1}}\cup \dots$$
be a fragment of a handle decomposition of a triad $(M,V_0,V_1)$. Assume that $q_{r}=q$, $q_{r+1}=q+1$. 
Both the embedded $b$-sphere $S_b$ of $H_r^{q}$ (which is diffeomorphic to $S^{m-q-1}$) 
and the embedded $a$-sphere $S_a$ of $H_{r+1}^{q+1}$ (which is diffeomorphic to $ S^{q}$ ) 
are submanifolds of  the $(m-1)$-manifold $V_{1,r}$,  and $\dim S_b +\dim S_a = m-1$. 
Then the adjacent handles  $H_{r}\cup H_{r+1}$ form {\it a pair of 
complementary handles} provided that $S_b$ and $S_a$ intersect transversely in $V_{1,r}$
at {\it exactly one point}. Note that under the above dimensional assumptions, 
by transversality and up to handle sliding we can assume
anyway  that $S_b$ and $S_a$ intersect transversely at a finite number of points. }
\end{definition}

\smallskip

{\it Cancelling/inserting pairs of  complementary handles.}  We can state the basic handle cancellation result. 

\begin{proposition}\label{cancellation} If 
$$\ \dots \cup H_r^{q}\cup H_{r+1}^{q+1}\cup \dots \ $$
is a pair of  complementary handles  in a handle decomposition
of $(M,V_0,V_1)$, then $(M_{r-1},V_0,V_{1,r-1})$ is diffeomorphic to $(M_{r+1},V_0,V_{1,r+1})$.
Hence we can cancel the pair and get a handle decomposition of the form
$$C_0\cup H_1^{q_1}\cup \dots \cup H_{r-1}^{q_{r-1}}\cup H_{r+2}^{q_{r+2}} \dots \cup H_k^{q_k} \cup C_1 \ . $$ 
Reciprocally, we can freely insert a pair of  complementary handles between
any two adjacent handles  into a given decomposition.
\end{proposition} 

\medskip

We postpone the proof below. 
\smallskip

A key problem is to study the handle decompositions of a given triad up to the {\it move-equivalence}
relation generated by such basic moves.  In fact by using Cerf's theory \cite{Ce2} (see \cite{Kirby}),
one can prove the following fact. 

\begin{theorem}\label{completehandlemove} Any two handle decompositions of a triads $(M,V_0,V_1)$
are move-equivalent to each other.
\end{theorem}

We will not prove  nor use such a rather demanding result. We limit to some remarks and simple applications.

\medskip

$\bullet$ For every handle decomposition $\Hh$ of $(M,V_0,V_1)$
set
$$\chi(\Hh)= \sum_q (-1)^q |\Hh_q| $$
where $|\Hh^q|$ denotes the number of $q$-handle of $\Hh$.
Obviously this {\it characteristic} of $\Hh$ is  move-equivalence {\it invariant}. We will see later that
$\chi(\Hh)$ has in fact an intrinsic topological meaning (see Remark \ref{handle-chi}).

\medskip

$\bullet$ The following is a first important application of sliding handle in order to specialize the handle decompositions. 
Let us give first a definition
\begin{definition}\label{orderedH}{\rm A handle decomposition of $(M,V_0,V_1)$ is said
{\it ordered} if
\begin{itemize}
\item For every $q=0,\dots, m-1$, the $q+1$ handles are attached after the $q$-handles;
\item  For every $q=0,\dots, m$, the $q$-handles are attached simultaneously.
Precisely, if $\Hh^q$ denotes the pattern of $q$-handle, $M_{q-1}=C_0\cup \Hh_0 \cup \dots \cup \Hh_{q-1}$
then the attaching maps of the handles in $\Hh_{q}$ have disjoint images in $V_{1,q-1}$.
\end{itemize}
}
\end{definition}

\begin{proposition}\label{reorder} {\rm (Reordering)} By handle sliding, every handle decomposition of $(M,V_0,V_1)$
can be transformed into an ordered decomposition. 
\end{proposition}
\Dim Let 
$$ \dots \cup H_r^{q_r}\cup H_{r+1}^{q_{r+1}}\cup \dots$$
be a fragment of a given handle decomposition $\Hh$. Set $q_{r}=p$ , $q_{r+1}=q$, and
assume that $p\geq q$. Then the embedded $b$-sphere $S_b$ of $H_r^{p}$ is diffeomorphic to $S^{m-p-1}$ 
while the embedded $a$-sphere $S_a$ of $H_{r+1}^{q}$ is diffeomorphic to $ S^{q-1}$. Then
$\dim S_b +\dim S_a \leq m-2 <m-1$. Up to handle sliding, we can assume that $S_b$ and $S_a$ are transverse submanifolds of 
the $(m-1)$-manifold $V_{1,r}$, so that $S_b \cap S_a = \emptyset$. There is a tubular neighbourhood $U$ of $S_b$ contained in the $b$-tube
$T_b$ around $S_b$, such that $S_a\cap U=\emptyset$; $T_b$ itself is a tubular neighbourhood of $S_b$. Hence by the uniqueness of the tubular neighbourhood
up to isotopy and the extension of isotopy to diffeotopy, there is a diffeotopy of $V_{1,r}$ which keeps $S_b$ fixed and pushs  the complement of $U$ in $T_b$
(hence $S_a$) outside $T_b$. It follows that up to handle sliding the two handles have now disjoint attaching tubes so that we can attach them in the inverse order or even
simultaneously. The proposition follows by several applications of this argument.

\cvd

\smallskip 

\begin{remark}{\rm In terms of Morse functions, the last proposition corresponds to the existence of Morse functions such that critical points of the same
index share the same critical value, and the critical values strictly increase together with the corresponding indices.}
\end{remark} 

\medskip

{\it Proof of  Proposition \ref {cancellation}.} Let us consider first the simplest case $q=0$. Attaching a $0$-handle means
``create'' a new disjoint $m$-ball component 
$$H_r^0=D^m= \{0\}\times D^m \ . $$ 
The whole boundary $S^{m-1}$ forms the
$b$-sphere. If the $1$-handle $H_{r+1}^1$ is complementary to $H_r^0$, then its attaching map embedds one component of
$$\partial D^1 \times D^{m-1}= \{-1,1\}\times D^{m-1}$$ into $S^{m-1}$, while the other component is embedded into
$V_{1,r-1} = V_{1,r}\setminus S^{m-1}$. The partial attachment of $D^{1}\times D^{m-1}$ to $D^m$ is a shelling (refer to
Section \ref {connected-sum}) 
of $D^m$ producing another diffeomorphic copy of $D^m$. Then the remaining component of the attaching map   
finally produces a shelling of $M_{r-1}$ hence a diffeomorphic copy of it.  The same facts hold in the general case
by a more elaborate argument. Assume first that the complementary handles have normalized attaching
maps as follows. Let us decompose the $b$-sphere $S_b$ of $H_r^{q}$ as $S_b=D^+_b \cup D^-_b$, where both
$D^\pm_b$ are diffeomorphic to $D^{m-q-1}$ and intersect along an equatorial $(m-q-2)$-sphere.
Then the $b$-tube around $S_b$ is given as $T_b=D^q\times (D^+_b\cup D^-_b)$. Similarly for the $a$-sphere and the $a$-tube
of $H^{q+1}_{r+1}$, let $S_a= D^+_a\cup D^-_a$, $D^\pm_a \sim D^q$, $D^+_a\cup D^{-}_a \sim S^{q-1}$, $T_a= (D^+_a \cup D^-_a)\times D^{m-q-1}$.
Assume that the intersection, say $A$, between the image of the attaching map of $H^{q+1}_{r+1}$ and $T_b$ is equal to $D^q\times D^+_b$,
and that the inverse image of $A$, say $\hat A$, is equal to $D^+_a\times D^{m-q-1}$, so that $\hat A \sim A$ and $\hat A \cap S_a=D^+_a$ is mapped
onto $D^q\times \{x_0\}$, $x_0$ being the `centre' of $D^+_b$. In such a normalized situation, we can factorize the attachment of the pattern
made by the two complementary handles as follows:
\begin{enumerate}
\item First glue $H^{q+1}_{r+1}$ to $H^q_{r}$ by using as attaching map the restriction of the whole attaching map to $\hat A$.
This is a shelling of a disk, so it results a smooth $m$-disk with a residual attaching zone contained in the boundary and diffeomeorphic
to a $(m-1)$-disk.

\item Perform the residual attachment; actually this is a further shelling over $M_{r-1}$.
\end{enumerate}
\smallskip

This achieves the result in the normalized situation. In our hypothesis, a priori we have such a normalized situation
provided we replace the whole $b$-tube $T_b$ with a smaller tubular neighbourhood $U$ of $S_b$ contained in $T_b$.
Now, similarly to the proof of Proposition   \ref{reorder}, by the uniqueness of the tubular neighbourhood
up to isotopy and the extension of isotopy to diffeotopy, there is a diffeotopy which keeps $S_b$ fixed and transforms $U\cup H^{q+1}_{r+1}$
to a pair of complementary handles in normal situation. This completes the proof.

\cvd

\smallskip

$\bullet$ A measure of the complication of a given handle decomposition is the {\it total number of handles}. For example if it is equal to
$0$, then $(M,V_0,V_1)$ is diffeomorphic to the {\it product triad} $(V_0\times [0,1],V_0,V_0)$, in particular $V_0$ and $V_1$
are diffeomorphic; if a boundaryless $M$ has a decomposition with only one $0$-handle and one $m$-handle, then $M$ is a twisted sphere. 
A natural task is to try to reduce such a complication by applying to a given decomposition some instances
of the basic moves. The following is a first simple but useful step in this direction.
\begin{proposition} \label{0-m-elimination} {\rm (Cancellation of $0$- and $m$-handles)}  Assume that $M$ is connected.
Then:

(1) For every triad of the form $(M,\emptyset,\emptyset)$ (i.e. $M$ is boundaryless),  every handle decomposition $\Hh$ is move-equivalent to
 an ordered decomposition $\Hh'$ with only one $0$-handle and only one $m$-handle.

(2) For every triad of the form $(M,\emptyset,\partial M)$, $\partial M \neq \emptyset$,  every handle decomposition $\Hh$ is move-equivalent to
an ordered decomposition $\Hh'$ with only one $0$-handle and without $m$-handles.

(3) For every triad of the form $(M,\partial M,\emptyset)$, $\partial M \neq \emptyset$,  every handle decomposition $\Hh$ 
is move-equivalent to an ordered decomposition $\Hh'$ with only one
$m$-handle and without $0$-handles.

(4) For every triad of the form $(M,V_0, V_1)$, both $V_0$ and $V_1$ being non empty,  every handle decomposition $\Hh$ 
is move-equivalent to an ordered decomposition $\Hh'$ without both $0$- and $m$-handles. 
\end{proposition}
\Dim Let us prove (1) and (2) simultaneously. By handle sliding we can assume that the decompostion is ordered.
Assume that  we have attached a certain number of $0$-handles, that is we have created a set of disjoint components
diffeomorphic to $D^m$. The only way to restore the fact that $M$ is connected is by means of the $1$-handles.
By successive application of elimination of complementary $H^0\cup H^1$ or reordering we eventually rich two possible
situations: either we remain with only one $0$-handle and this happens when $V_0=\emptyset$ (if there are no longer 
complementary $H^0\cup H^1$ to eliminate then $M$ would be not connected), or we remain with no $0$-handles and this 
happens when $V_0\neq \emptyset$ and the $1$-handles
connect to each other all the components of $C_0$. To deal with the the $m$-handles is enough to apply the same
argument to the dual decompostion. 

\cvd

\begin{remark}{\rm In terms of Morse functions, for example the first case of the above proposition
corresponds to the existence of functions with only one minimum and one
maximum. Similarly for the other cases.}
\end{remark} 

\subsection{The CW complex associated to an ordered decomposition}\label{CW}
Let $M$ be boundaryless. Let 
$$H^0\cup \{H^1\} \cup \{H^2\} \dots \cup \{H^{m-1}\} \cup H^m$$ 
be an ordered handle decomposition of the triad
$(M, \emptyset, \emptyset)$ with one $0$-handle and one $m$-handle; $\{H^j\}$ means a (possibly empty) pattern of $i_j$ $j$-handles attached simultaneously.
Every handle $H$ has a natural retraction 
$$r:H\to {\rm core}(H)\cup a-{\rm tube}(H)$$ 
which realizes a homotopy equivalence.
By using the notations fixed above, we are going to construct inductively homotopy equivalence $$l_j: W_j\to K_j$$ 
where $K_0$ consists of one point and $K_j$ will be obtained by attaching $i_j$ $j$-cells to $K_{j-1}$; we eventually get a homotopy equivalence 
$$l: M \to K, \ K=K^m$$ where (by the very definition of this term) $K$
is a {\it finite CW-complex of dimension $m$}. So, let $K_0$ be the core of $H^0$; then $l_0: M_0\to K_0$ is an instance of retraction $r$ as above.
Assume we have defined $l_{j-1}: M_{j-1}\to K_{j-1}$. Then $$M_j=M_{j-1}\cup_{\{h_j\}} \{H^j\}$$ is homotopy equivalent (via the retraction $l_j:=l_{j-1}\circ \{r_j\}$)
to $$K_j= K_{j-1}\cup_{\{g_j\}} \{D^j\}$$ where  $\{g_j\}$ is the restriction of $l_{j-1} \circ \{h_j\}$.

Assume now  that $\partial M$ is not empty and consider the triad $(M,\partial M, \emptyset)$. In such a case the ordered handle decomposition
has no $m$-handles, hence there is an homotopy equivalence $l: M\to K$ where {\it $K$ is a finite CW-complex of dimension $d\leq m-1$}.    

\section{Compact $1$-manifolds}\label{1-man}
We use the handle technology developed so far in order to classify compact $1$-manifolds up to diffeomorphism. This is simple and intuitive; 
nevertheless it is a fundamental result with many applications (see Chapter \ref{TD-COBORDISM}). It is not restrictive to assume that these manifolds are connected.

\begin{proposition}\label{1-class} 
(1) A compact connected boundaryless $1$-manifold is diffeomorphic to $S^1$.

(2) A compact connect $1$-manifold with non empty boundary is diffeomorphic to the interval $D^1$.
\end{proposition}
\Dim In both cases apply Proposition \ref{0-m-elimination}. In the second case there is a handle decomposition of $(M,\emptyset,\partial M)$
formed by one   $0$-handle (of dimension $1$).  Hence $(M,\emptyset,\partial M)$ is diffeomorphic to $(D^1,\emptyset, \{\pm 1\})$.
In the second case there is a handle decomposition of $(M,\emptyset,\emptyset)$
formed by one   $0$-handle and one $1$-handle (of dimension $1$). Hence $M$ is a twisted $1$-sphere and we know from
Chapter \ref{TD-CUT-PASTE} that it is diffeomorphic to $S^1$.

\cvd

\chapter{Bordism}\label{TD-BORDISM}
For every $m\geq 0$, denote by $\Ss_m$ the class of smooth compact (not necessarily connected) boundaryless
$m$-manifolds. A natural question would be to classify the elements of $\Ss_m$ up to {\it diffeomorphism}. 
We can also specialize the question to the class $\Oo_m$ of oriented manifolds up to oriented diffeomorphism. 
Sometimes we will use $\Mm_m$ to indicate indifferently either $\Ss_m$ or $\Oo_m$.
It turns out that beyond
$m\leq 2$ these are  very demanding, even hopeless questions. Then it is natural to relax the diffeomorphism
to  a suitable equivalence up to (possibly oriented) {\it bordism}. 

On another hand, homotopy groups $\pi_m(X,x_0)$ of any pointed
topological space $(X,x_0)$ provide the basic examples of
topological/algebraic functors and are constructed by implementing the
following idea: to get information about a complicated ``unknow''
space $X$, continuously map into it ``tame'' spaces (the $m$-sphere)
and study the behaviour of these {\it singular tame} objects in $X$ up
to homotopy which is a sort of basic prototype of bordism between
maps. Note that the singular ``tame'' objects are in general not so
simple in spite of the tame source spaces because the maps and their
images in $X$ can be complicated. The same idea can be implemented by
considering singular smooth $m$-manifolds in $X$, that is continuous
maps $f:M\to X$ where $M\in \Mm_m$, up to suitable {\it bordism of
  maps} (naturally extending the bordism of manifolds mentioned
above).  This leads in a simple way to further topological/algebraic
functors; once also the relative theory for topological pairs $(X,A)$
has been developed, then one easily checks that these functors verify
the so called {\it Eilenberg-Steenrod axioms} which characterize {\it
  generalized homology theories}.  Of course all this specializes to
the case when $X$ itself belongs to $\Mm_k$, for some $k$. We will
develop this topological/differential specialization mainly in Chapter
\ref{TD-COBORDISM}.

\section{The bordism modules of a topological space} 
Let $X$ be a topological space. For every $m\geq 0$ a {\it singular
  smooth $m$-manifold in $X$} is a continuous map $f: M \to X$ where
$M\in \Ss_m$ Denote by $$\Ss_m(X)$$ the set of such singular manifolds
to which {\it we formally add the empty set}.
\begin{definition}\label{sing-boundary}{\rm $(M,f)\in \Ss_m(X)$ is a {\it singular boundary}
    if there are a compact smooth {(m+1)}-manifold with boundary $(W,\partial W)$,
    a diffeomorphism $\rho:M\to \partial W$,
a continuous map $F: W\to X$ such that $F\circ \rho = f$.}
\end{definition}
\smallskip

Let us put on $\Ss_m(X)$ the following relation:
\medskip

We say that $(M_0,f_0)$ is {\it bordant with $(M_1,f_1)$} and we write
$(M_0,f_0)\sim_b (M_1,f_1)$ if the disjoint union $(M_0,f_0)\amalg (M_1,f_1)$ is a singular boundary.
It is consistent to state that $(M,f)\sim_b \emptyset$ 
if and only if $(M,f)$  is a singular boundary.

\medskip

We claim that this is an {\it equivalence relation}:
\smallskip

$\bullet$  The cylinder $(M\times [0,1], F)$ , $F(x,t)=f(x)$ for every $t\in [0,1]$, establishes that
$(M,f)\sim_b (M,f)$, $\rho: M\amalg M \to (M\times\{0\})\amalg (M\times \{1\})$ being the natural inclusion.
\smallskip

$\bullet$ As the disjoint union is symmetric, then also $\sim_b$ is obviously symmetric.
\smallskip

$\bullet$ {\it Transitivity follows by gluing smooth manifolds along boundary components}. 
Precisely, assume that 
$(W_0, F_0)$, $\rho_0: M_0\amalg M_1 \to \partial W_0$   
realize $(M_0,f_0)\sim_b (M_1,f_1)$,
while   $(W_1, F_1)$, $\rho_1 : M_1\amalg M_2 \to \partial W_1$ realize $(M_1,f_1)\sim_b (M_2,f_2)$.
Then $F_0$ and $F_1$ match to define a smooth map $F_2$ on  
$W_2:= W_0 \amalg _\psi W_1$, where $\psi$ is the composition of the restriction of $\rho_0^{-1}$ to
$\rho_0(M_1)$ with the restriction of $\rho_1$ to $M_1$. Finally $(W_2,F_2)$ together with the disjoint
union of $\rho_0$ restricted to $M_0$ and $\rho_1$ restricted to $M_2$ realize $(M_0,f_0)\sim_b (M_2,f_2)$.

\medskip

We denote by $\eta_m(X)$ 
the quotient set $\Ss_m(X)/\sim_b$,  by 
$[M,f]$ 
the equivalence class of $(M,f)$. 

The disjoint union is an operation on $\Ss_m(X)$. It is immediate that it descends to the quotient, that is
$ [M,f]+[N,g]:= [M\amalg N, f\amalg g]$
is a well defined  operation on $\eta_m(X)$. We have
\begin{proposition}\label{eta+} $(\eta_m(X), +)$ is an abelian group.
\end{proposition}
\Dim The operation $+$ is associative and commutative because the disjoint union is associative and commutative. 
$[\emptyset]$ that is the class of the singular boundaries, is the zero element. For every $\alpha=[M,f]$, $-\alpha=\alpha$,
in fact by using the cylinder as above we see that $[M,f]+[M,f]=0$.
\cvd

\medskip

Since for every $\alpha$, $\alpha=-\alpha$,  then $(\eta_m(X), +)$ can be enhanced to be 
a $\Z/2\Z$-module, that is a $\Z/2\Z$-{\it vector space} $(\eta_m(X),+,\cdot)$; 
we call it the {\it unoriented $m$-bordism space} of $X$.

\medskip

\subsection{The oriented bordism $\Z$-modules} 
We follow the same sheme by using oriented manifolds.

We denote by $\Oo_m(X)$ the set of {\it oriented} singular $m$-manifolds $f: M\to X$,
that is  $M\in \Oo_m$. 

$(M,f)$ is a {\it singular oriented boundary} if $(W,F)$, 
$\rho:M \to \partial W$ are as above and we require furthermore that $(W,\partial W)$
is oriented and $\rho$ preserves the orientation. 

The relation $(M_0, f_0)\sim_{ob} (M_1,f_1)$ on $\Oo_m(X)$
is defined by requiring that $(M_1,f_1)\amalg (-M_2,f_2)$ is a singular oriented boundary. 
The verification that it is an equivalence relation is similar: 
\smallskip

- the cylinder can be naturally oriented in such a way that its
oriented boundary is $M\amalg -M$. 
\smallskip

- To get the symmetry it is enough to replace $W$ with $-W$.
\smallskip

- As for the transitivity, we glue again $W_0$ and $W_1$ by taking into account that the gluing
diffeomorphism $\psi$ reverses necessarily the orientation: in $\partial W_0$ there is a copy of $-M_1$
while in $\partial W_1$ there is a copy of $M_1$. Hence the gluing can be performed in the oriented category.
\medskip

We denote by $\Omega_m(X)$ the quotient set.
Again the operation $+$ on $\Omega_m(X)$ is induced by the disjoint union on $\Oo_m(X)$. 
It results a commutative group (i.e. a $\Z$-{\it module}) $(\Omega_m(X),+)$. 
Again $0=[\emptyset]$, that is the class of the singular oriented boundaries. By means
of the oriented cylinder we see that  $-[M,f]= [-M,f] $. 
This is the $m$-{\it oriented bordism module} of the topological space $X$. 

There is a natural group homomorphism $$\sigma_m: \Omega_m(X)\to \eta_m(X)$$ which
maps the class of $(M,f)$ in $\Omega_m(X)$ to its class in $\eta_m(X)$, just by ``forgetting the orientation''.

As many considerations run formally in the same way for both bordism versions, 
sometimes we will indifferently indicate by $\Mm_m(X)$ either $\Ss_m(X)$ or $\Oo_m(X)$, and by  $\Bb_m(X)=\Bb_m(X;R)$ 
either the quotient $R$-module $\eta_m(X)$ or $\Omega_m(X)$, $R=\Z/2\Z, \Z$.

\begin{lemma}\label{b-vs-diffeo} Let $\phi: N\to M$ be a diffeomorphism (preserving the orientation in the oriented
  setting); $f:M\to X$, $m=\dim M$. Then $[M,f]=[N,f\circ \phi] \in
  \Bb_m(X)$
  \end{lemma}
  \Dim The cylinder $(M\times [0,1], f\circ \pi)$ ($\pi:
  M\times [0,1]\to M$ being the projection), and $\rho: M \amalg N
  \to (M\times \{0\})\amalg (M\times \{1\})$, $\rho= {\rm id}_M
  \amalg \phi$, realize $(M,f)\sim_\Bb (N, f\circ \phi)$.

  \cvd

  \smallskip

  \begin{remark}\label{rho=id}{\rm Let $(M,f)$ be a singular boundary in $X$.
      Let $((W,\partial W),F)$ and $\rho:M\to \partial W$ realize $(M,f)\sim_\Bb \emptyset$.
      By applying Lemma \ref{b-vs-diffeo} we have 
      $$(M,f)\sim_\Bb (\partial W, \partial F)$$ 
      and this is realized by a cylinder; obviously 
      $((W,\partial W),F)$ and ${\rm id}_{\partial W}$ realize 
      $$(\partial W,\partial F)\sim_\Bb \emptyset \ . $$ 
      By applying to this situation the gluing
      argument employed to show the transitivity, we can conclude that
      it is not restrictive to require that $M=\partial W$ and $\rho= {\rm id}_M$
    }
    \end{remark}

{\bf An important special case.} When $X=\{x_0\}$ consists of one point, 
then the maps are immaterial and, by definition, $\Bb_m:=\Bb_m(\{x_0\})$ is the quotient of 
$\Mm_m$ up to {\it bordism of manifolds}. It follows from Lemma \ref{b-vs-diffeo} that
the bordism extends the diffeomorphism equivalence in the category. 

\section{Bordism covariant functors}\label{b-functor}  We have the following Proposition.
All verifications are straighforward consequence of the very definitions.

\begin{proposition}\label{functor-b} For every $m\geq 0$, 
$$ X \  \Rightarrow \ \Bb_m(X)$$
$$ g:X\to Y \ \Rightarrow \ g_* :\Bb_m(X)\to \Bb_m(Y), \ g_*([M,f])= [M, g\circ f] $$
is a covariant functor from the category of topological spaces and continuous maps to the category of
$R$-modules and $R$-linear maps. That is
$$ (g\circ h)_* = g_*\circ h_*$$
$$({\rm id}_X)_* = {\rm id}_{\Bb_m(X)} \ . $$
\end{proposition}

\cvd

\smallskip

In particular if $g:X \to Y$ is a homeomorphism, then $g_*$ is a $R$-linear isomorphism
with inverse $(g^{-1})_*$. Considered up to linear isomorphism, $\Bb_m(X)$ {\it is an invariant of
  the topological type of $X$}. The family introduced above of ``forgetting'' linear maps 
  $$\{\sigma_m: \Bb_m(X;\Z)\to \Bb_m(X;\Z/2\Z)\}$$
is {\it functorial}, that is they form commutative squares together with the respective 
families of $g_*$'s;   in form
of a slogan: ``$g_* \circ \sigma = \sigma \circ g_*$''.

\section{Relative bordism of topological pairs}\label{relative-b}
We consider topological pairs $(X,A)$ where $A$ is a subspace of $X$
and the class $\Mm_m^\partial$ of compact smooth $m$-manifolds with boundary
$(M,\partial M)$. This incorporates the ``absolute situations'' by identifying
$X$ with the pair $(X,\emptyset)$ and a boundaryless manifold $M\in \Mm_m$ with
$(M,\emptyset)$. 

A {\it relative singular $m$-manifold in $(X,A)$} is a continuous map
of pairs $$f:(M,\partial M)\to (X,A)$$ where by definition $f(\partial M)\subset A$ and
$(M,\partial M)\in \Mm_m^\partial$. We set $\Mm_m(X,A)$ the collection
of such relative singular $m$-manifolds.

\begin{definition}\label{relative-boundary} {\rm $f:(M,\partial M)\to (X,A)$ is a {\it
relative singular boundary} if there are continuous pair maps $F:(W,V)\to (X,A)$,
$\rho: (M,\partial M)\to (Z,\partial Z)$ such that:
\begin{enumerate}
\item $(W,\partial W) \in \Mm^\partial _{m+1}$;
\item $(V,\partial V)$ and $(Z,\partial Z)$ are smooth $m$-submanifolds of $\partial W$
such that $$\partial W= V\cup Z, \ V\cap Z = \partial V = \partial Z \ ; $$
\item $\rho: (M,\partial M)\to (Z,\partial Z)$ is a  smooth diffeomorphism (preserving
the orientation in the oriented case). In particular if $\partial M$ is empty, 
then $V$ and $Z$ are also boundaryless, $\partial W= V \amalg Z$
and $F(V)\subset A$.   
\end{enumerate}
}
\end{definition}

\smallskip

We put on  $\Mm_m(X,A)$ the equivalence relation $(M_0,\partial M_0,f_0)\sim_\Bb (M_1,\partial M_1,f_1)$
if and only if $(M_0,\partial M_0,f_0)\amalg (-M_1,\partial M_1,f_1)$ is a relative singular boundary
(in the unoriented case the sign ``$-$'' is immaterial). The verification that it is an equivalence relation (in particular
the transitivity) incorporates some instances of corner smoothing, accordingly with Remark \ref {gen-glue}.

The disjoint union on $\Mm_m(X,A)$ descends
to a operation $+$ on the quotient set that eventually makes it 
a $R$-module $\Bb_m(X,A)=\Bb_m(X,A;R)$, called the {\it realtive $m$-bordism $R$-module} of the topological
pair $(X,A)$.   

Proposition \ref{functor-b} extends directly:

\begin{proposition}\label{relative-b-functor} For every $m\geq 0$,
$$ (X,A) \  \Rightarrow \ \Bb_m(X,A)$$
$$ g:(X,A)\to (Y,B) \ \Rightarrow \ g_*: \Bb_m(X,A)\to \Bb_m(Y,B), \ g_*([M\partial M,f])= [M,\partial M, g\circ f] $$
is a covariant functor from the category of pairs of topological spaces and continuous pair maps to 
the category of $R$-modules and $R$-linear maps.
\end{proposition}

\section{On Eilenberg-Steenrood axioms}\label{ES-axiom}
The singular homology (sometimes called ``Betti homology'')
with coefficients in the ring $R$ is a determined family of functors
(indexed by $m\geq 0$) of
the same kind of
Propositions \ref{functor-b}, \ref{relative-b-functor}.  The
(E-S)-axioms are abstractions of some properties verified by the singular
homology functors and which deserve the name because all models (no
matter how they have been produced) that fulfill such axioms are
isomorphic to each other, at least if one restricts to pairs of
compact CW-complexes (see \cite{Hatch}). It turns out that the most
critical one is the so called {\it dimension axiom}; every model which
verifies the other axioms (with the possible exception of
``dimension'') is called a {\it generalized homology theory}. We are
going to check that this is the case of bordism.  The verifications
are of geometric/topological nature and often immediate consequences
of the definitions.
\medskip

{\bf The homotopy axiom.} {\it If $g_0, g_1: (X,A)\to (Y,B)$ are homotopic through pair maps, then $g_{0,*} =g _{1,*}$.}
\medskip

We have to show that for every $[M,\partial M,f]\in \Bb_m(X,A)$, 
$$[M,\partial M, g_0\circ f]=[M,\partial M, g_1\circ f]  \ {\rm  in} \  \Bb_m(Y,B) \ . $$
Given a homotopy $$G: (X\times [0,1],A\times [0,1])\to (Y,B)$$ between $g_0$ and $g_1$,
then $$F: (M\times [0,1], \partial M \times [0,1])\to (Y,B), \  f_t =g_t\circ f$$ together 
with the natural inclusion of $(M,\partial M)\amalg (M,\partial M)$ in $\partial (M\times [0,1])$
realize that $(M,\partial M, g_0\circ f) \sim_\Bb (M,\partial M, g_1\circ f)$.

\cvd

\smallskip

This implies that if $g:(X,A)\to (Y,B)$ is a relative homotopy equivalence, then $g_*$ is a $R$-linear
isomorphism. Up to isomorphism, {\it the bordism modules are invariants of the homotopy type rather than
the topology type}.

\medskip

{\bf Direct sum over path connected components.} For every topological space $X$,
$\Bb_m(X)$ is isomorphic to the direct sum of the modules $\Bb_m(X_c)$, where
$X_c$ varies among the path connected components of $X$.  This follows from
the fact that continuous maps send every path connected component of a manifold $M$
into one path connected component of $X$. A similar fact holds in the relative version.

\medskip

{\bf Long exact sequence.} For every $m\geq 1$ there is the natural well defined $R$-linear map
$$ \partial: \Bb_m(X,A)\to \Bb_{m-1}(A), \ \partial ([M,\partial M,f])=[\partial M, \partial f] \ . $$
Denote by $i_*: \Bb_m(A)\to \Bb_m(X)$, $j_*: \Bb_m(X,\emptyset)\to \Bb_m(X,A)$
the $R$-linear maps induced by the inclusions. Then we have a {\it bordism long sequence} of linear maps
$$\cdots \rightarrow \Bb_m(A) \xrightarrow {i_*} \Bb_m(X)\xrightarrow {j_*} \Bb_m(X,A)\xrightarrow {\partial} \Bb_{m-1}(A)\rightarrow \cdots$$
which ends on the right with the $0$ $R$-module.

Recall that a sequence of linear maps
$$A\xrightarrow{\alpha} B \xrightarrow {\beta} C$$
is {\it exact in $B$} if $\ker(\beta)= \alpha(A)$. Then we have:
\smallskip

\noindent {\it
(1)The long sequences are {\it functorial}: if $g:(X,A)\to (Y,B)$ then the respective long sequences
  together with the family of linear maps $\{g_*\}$ form commutative squares.

\noindent (2)  Every bordism long sequence is exact everywhere.}
\smallskip

\noindent Fuctoriality is immediate consequence of the definitions.
The verifications of exactness are simple and useful exercises. Let us show for example that the above long sequence
is exact in $\Bb_m(X,A)$. If $[N,g]\in \Bb_m(X)$ then
$N$ is boundaryless, so it is clear that $\partial \circ j_*([N,g])=0 \in \Bb_{m-1}(A)$. On the other hand, 
If $(M,\partial M, f)$ is in the kernel of $\partial$ and  $(W,\partial W, F)$
realizes that $(\partial M, \partial f)$ is a boundary, then by gluing
$W$ and $M$ along $\partial M$, we get $\tilde f: \tilde M \to X$, $\tilde M$ being boundaryless, $\tilde f$ obtained 
by matching $f$ and $F$, such that $j_*([\tilde M, \tilde f])= [M,\partial M, f] \in \Bb_m(X,A)$.

\cvd 

\medskip

{\bf Excision.}  Let $Z\subset A \subset X$ be a triad of topological space. Assume that the closure $\bar Z$ of 
$Z$ in $X$ is contained in the interior $\mathring {A}$ of $A$.
Then we have
\smallskip 

{\it For every $m\geq 0$, the linear map induced by the inclusion
$$ i_*: \Bb_m(X\setminus Z, A\setminus Z)\to \Bb_m(X,A) $$
is an isomorphism. We say that $Z$ {\rm is excisable}.}

\medskip

Let us prove first that it is surjective. Let $[M,\partial M,f]\in \Bb_m(X,A)$.
The manifold $M$ can be endowed with a distance $d$ compatible with its topology 
so that $(M,d)$ is a compact metric space; for example embedd $M$ in some
$\R^n$ an take the distance induced by the euclidean distance.
$K:=f^{-1}(\bar Z)$ is a compact set contained in the open set $\tilde A:= f^{-1}(\mathring{A})$.
The distance function from $K$ 
$$\delta: M \to \R$$ 
is  non negative, continuous and  $K=\{\delta = 0 \}$.
Then there is a smooth approximation say $g: M \to \R$
and a regular value $\epsilon>0$ of both $g$ and $\partial g$,
sufficiently close to $0$, such that $\tilde M:= \{g\geq \epsilon\}$ is a compact $m$-submanifold
with corners such that $\partial \tilde M = \{g=\epsilon\}$ is contained in $\tilde A$. Up to smoothing the corners,
if  $\tilde f$ is the restriction
of $f$ to $\tilde M$,  we finally have 
that $[\tilde M, \partial \tilde M,\tilde f]\in \Bb_m(X\setminus Z, A\setminus Z)$
and $i_*([\tilde M, \partial \tilde M, \tilde f])=[M, \partial M, f]\in \Bb_m(X,A)$. 
To prove the injectivity we apply the same argument
to $(W,\partial W, F)$ which realizes that a $(M, \partial M, f)\in \Mm_m(X\setminus Z, A\setminus Z)$ 
is a relative singular boundary in $(X,A)$.

\cvd  

\medskip

{\bf About the dimension axiom.} This axiom for the  singular homology (with coefficients
in $R$) determines the homoloy modules of a singleton. Precisely, the $0$-module is
isomorphic to $R$, while the others are all trivial.

For every $X$,  $\Bb_0(X)$ has a clear topological meaning. In fact, by using the classification of compact
$1$-manifolds (Proposition \ref {1-class} ), it is easy to check that 
it is isomorphic to the direct
sum  $\oplus_{\pi_0(X)} R$, where $\pi_0(X)$ is the set
of path connected components of $X$. In particular  $\Bb_0 = R$.

On the other hand,  we do not know for the moment if the modules $\Bb_m$, $m>1$ are all trivial.
In fact we will see in Section \ref {EP-bord}  that {\it they are not}.
\medskip

The (E-S)-axioms establish in more or less
explicit way relations between the modules $\Hh_*(X)$ in any
(generalized) homology theory
of a given space and the ones of the presumably simpler pieces of some 
suitable decomposition of $X$. If 
also ``dimension'' holds, then in many cases they allow to compute 
(up to linear isomorphism) the modules of $X$. Without ``dimension''
things are more complicated. The first interesting cases to face
are $X=S^n$  or the pair $(X,A)=(D^n, S^{n-1})$. 
These are the building blocks of CW-complexes.
\smallskip

$\bullet$ As the $n$-disk is contractible for every $n\geq 0$, by
``homotopy'' $\Hh_m(D^n)\sim \Hh_m$ for every $m\geq 0$.

$\bullet$ For every $n\geq 1$,  we can decompose $S^n$ as the union of the closed northern and southern hemispheres (both diffeomorphic to $D^n$)
$$S^n= D^+ \cup D^-, \ D^+\cap D^- = S^{n-1} \ .$$
We claim that the inclusion induces isomorpfisms
$$i_*: \Hh_m(D^+,S^{n-1})\to \Hh_m(S^n, D^-) \ . $$ We cannot apply
directly ``excision'' to $Z= \mathring {D}^-$. We can do it by using
instead $\tilde Z\subset D^-$ equal to the complement of a small
collar of $S^{n-1}$ in $D^-$. Finally we use ``homotopy'' and the fact
that $(S^n\setminus \tilde Z, D^-\setminus \tilde Z)$ retracts to
$(D^+, S^{n-1})$ to achieve the required isomorphisms.

$\bullet$ Again for $n\geq 1$, we have the exact long sequence of the pair $(D^n,S^{n-1})$
$$\cdots \rightarrow \Hh_m(S^{n-1}) \xrightarrow {i_*} \Hh_m (D^n) \xrightarrow {j_*} \Hh_m(D^n,S^{n-1})\xrightarrow {\partial} \Hh_{m-1}(S^{n-1})\rightarrow \cdots$$
and the one of the pair $(S^n,D^-)$
$$\cdots \rightarrow \Hh_m(D^-) \xrightarrow {i_*} \Hh_m(S^n) \xrightarrow {j_*} \Hh_m(S^n,D^-)\xrightarrow {\partial} \Hh_{m-1}(D^-) \rightarrow \cdots$$

$\bullet$ If the theory $\Hh$ verifies also ``dimension'',
by simple algebraic considerations one realizes  that for $n\geq 1$,
\begin{itemize}
\item $\partial: \Hh_m(D^n,S^{n-1})\to \Hh_{m-1}(S^{n-1})$ is an isomorphism for $m\geq 2$;
\item $\j_*: \Hh_m(S^n)\to \Hh_m(S^n,D^-)$ is an isomorphism for $m\geq 2$;
\item for every $m\geq 1$, $\Hh_m(S^n)\sim \Hh_{m-1}(S^{n-1})$
(immediately for $m\geq 2$, with a little extrawork for $m=1$).
\end{itemize}

\noindent Then by a simple induction we  can finally  achieve the computation: 
\smallskip

{\it For every $n\geq 1$, $m=0, n$,
$$ \Hh_m(S^n)\sim \Hh_m(D^n,S^{n-1}) \sim R \  . $$
For   every $n\geq 1$, $m\geq 1$,  $m\neq n$,
$$\Hh_m(S^n)\sim \Hh_m(D^n,S^{n-1})=0  \ . $$}
\smallskip

If the theory (like the bordism) does not verify ``dimension'' the
considerations based on the other axioms hold as well but are not
immediately conclusive.

\section{Bordism non triviality}\label{b-nontriv} 
By combining the axioms with the specific way the bordism has been
defined, we will provide a few evidences that it is not trivial.

$\bullet$ Assume that $X$ is path connected. Consider the long exact
sequence of a pair $(X,x_0)$ for some base point in $X$,
$$\cdots \rightarrow \Bb_m \xrightarrow {i_*} \Bb_m (X) \xrightarrow
{j_*} \Bb_m(X,x_0)\xrightarrow {\partial} \Bb_{m-1} \rightarrow
\cdots$$ it is immediate by the bordism definition that $\partial = 0$
(hence $j_*$ is onto) and that $i_*$ is injective.  Hence every
$\Bb_m(X)$ contains a submodule isomorphic to $\Bb_m$ which in general
is not trivial. Note that since $X$ is path connected, by ``homotopy''
this submodule does not depend on the choice of the base point $x_0$.
When $R=\Z/2\Z$ (algebra is simpler in the case of vector spaces) we
have $\eta_m(X) \sim \eta_m \oplus \eta_m(X,x_0)$.

$\bullet$ Assume that $X$ is a compact connected boundaryless
(possibly oriented) smooth $m$-manifolds. Then by the approximation
theorems of continuous maps by smooth maps, {\it it is not restrictive
  to assume that all maps entering the bordism treatment are
  smooth}. We have

\begin{proposition}\label{fund-class}
  $[X,{\rm id}_X]\in \Bb_m(X)$ is non trivial and does not belong to
  $\Bb_m \subset \Bb_m(X)$.  In particular $\dim \eta_m(X) \geq 1+
  \dim \eta_m$.
\end{proposition}
\Dim Assume that it is trivial; then there is a smooth map $F: W\to
X$, such that $\partial W= X$ and $F_{|X}={\rm id}_X$.  Let $p\in
X$. Clearly it is a regular value for $\partial F$. Apply to $F$ the
transversality theorems relatively to $\partial F$.  Then we can
assume that $F\pitchfork \{p\}$, $Y=F^{-1}(p)$ is a proper
$1$-submanifold of $(W,X)$ and $p\in Y$.  By the classification of
compact smooth $1$-manifolds, $p$ is contained in an interval
component $I\subset Y$, hence there is another $p' \in \partial I
\subset X$ such that $p'\neq p$ and $\partial F(p)=p=\partial
F(p')=p'$. This is absurd.  This proves that $[X,{\rm id}_X]\neq
0$. Let $c:N \to \{ p \}$ be a constant map representing some element
of $\Bb_m \subset \Bb_m(X)$. Let $q \neq p$ so that it is a regular
value for both ${\rm id}_X$ and $c$.  If $(W,F)$ would realize a
bordism between $(X,{\rm id}_X)$ and $(N,c)$, by applying again the
relative first transversality theorem to $(W,F)$ we should deduce that
$\partial F^{-1}(q)=\{ q \}$ is a boundary; again by the
classification of compact $1$-manifolds this is absurd.

\cvd

\smallskip

By a similar argument, we have the following generalization. 

\begin{proposition}\label{fund-class2}  In the setting of Proposition \ref{fund-class}
Let $[N]\in \Bb_k$ be non trivial, and consider $(N\times X, {\rm id}_X\circ \pi)$, $\pi$ being the projection
to $X$. Then $[N\times X,  {\rm id}_X\circ \pi] \in \Bb_{m+k}(X)$ is non trivial.
\end{proposition}

\cvd

\smallskip

The class $[X, {\rm id}_X]\in \Bb_m(X)$ is called the bordism {\it fundamental class} of the (possibly oriented)
manifold $X$. If $X$ has non empty boundary similar facts hold for $[X,\partial X, {\rm id}_X]\in \Bb_m(X,\partial X)$. 
\smallskip

$\bullet$  {\it (On the bordism modules of spheres)} For every $n\geq 1$, consider  $X=S^n$ or $(D^n,S^{n-1})$ as above.  
If $m<n$,
by transversality we can assume that every class $\alpha$ in $\Bb_m(S^n)$ is represented by a smooth and non surjective
map $f: M\to S^n$; say that $\infty \notin f(M)$. Then $f$ factorizes through
$\R^n \subset \R^n\cup \infty = S^n$, hence it is homotopic to a constant map.  By ``homotopy''  
$\alpha$ belongs to $\Bb_m \subset \Bb_m(S^n)$, hence if $m<n$, $\Bb_m(S^n)=\Bb_m$.

Referring to the long exact sequence for the pair $(S^n,D^-)$, using that $D^-$ is contractible and ``homotopy'',
we have that $\partial=0$ so that $j_*$ is onto; and $i_*$ is injective. In particular we have
$$\eta_m(S^n)\sim \eta_m \oplus \eta_m(S^n,D^-)\sim \eta_m \oplus \eta_m(D^n,S^{n-1})$$
where for the last isomorphism we have applied ``excision'' and ``homotopy" as above.

Referring to the long exact sequence for the pair $(D^n,S^{n-1})$, we see that $i_*$ is onto, hence $j_*=0$, $\partial$ is
injective. Hence we have in particular that
$$\eta_{m-1}(S^{n-1}) \sim \eta_{m-1}\oplus \eta_m(D^n,S^{n-1}) \ ; $$ 
hence 
$$\eta_{m-1}(S^{n-1})\oplus \eta_m \sim \eta_m(S^n)\oplus \eta_{m-1}\ . $$
\smallskip

By a similar inductive argument already used to compute $\Hh_*(S^n)$ when the theory $\Hh$ verifies also
``dimension'', we can eventually achieve the determination of $\eta_*(S^n)$.
\begin{proposition} (1) For every $m\geq 0$, $\eta_m(S^0)= \eta_m \oplus \eta_m$.
 
(2) For every $n\geq 1$, for every $0\leq m <n$,
$\eta_m (S^n)=\eta_m$.

(3) For every $n\geq 1$, $k\geq 1$, 
$$\eta_{n+k}(S^n)= \eta_k \oplus \eta_{n+k} \ . $$ 
Precisely every class in $\eta_{n+k}(S^n)$ either belongs to $\eta_{n+k}$ 
or is of the form $[N\times S^n,  {\rm id}_{S^n} \circ \pi]$
as in Proposition \ref{fund-class2}

\end{proposition}

\cvd

\smallskip

It is already clear from these few remaks that the determination of $\Bb_m$, for every $m\geq 0$, that is of the actual
failure of ``dimension'' is a key point of this story. 

\section{Relation between bordism and homotopy group functors}\label{b-vs-homot} Here we assume some familiarity with the homotopy group $\pi_m(X,x_0)$, $m\geq 1$,
of the {\it pointed} topological space $(X,x_0)$ (see for instance \cite{Hatch}).
When $m=1$ it is called the {\it fundamental group}. Let us recall anyway a few facts.

$\bullet$ As a set $\pi_m(X,x_0)$ is formed by the classes $<f>$ of pointed continuous maps $f:(S^m,p)\to (X,x_0)$
considered up to pointed homotopy. It is endowed with a natural group operation
``$\cdot$'' well defined on any given representatives. The $1$ element is the class of the constant
pointed map. They are abelian for $m\geq 2$ while the fundamental group is not in general.
If $X$ is path connected, up to group isomorphism they do not depend on the choice of the base
point. 
\smallskip

$\bullet$ Similarly to the bordism, we have for every $m\geq 1$ a covariant functor
$$ (X,x_0) \ \Rightarrow \ \pi_m(X,x_0)$$
$$ g:(X,x_0)\to (Y,y_0) \ \Rightarrow \ g_*: \pi_m(X,x_0)\to \pi_m(Y,y_0), \ g_*(<f>)=<g\circ f> $$
from the category of pointed topological spaces and pointed continuous maps to the category of
groups (abelian for $m\geq 2$) and group homomorphisms.

\smallskip

$\bullet$ There is a relative version for pointed pairs $(X,A,x_0)$ ($x_0\in A$) of topological spaces.
Then the elements of $\pi_m(X,A,x_0)$ are relative homotopy classes $<f>$ of maps $f:(D^m,S^{m-1},p)\to (X,A,x_0)$
As usual the ``absolute'' theory is incorporated by identifying $(X,x_0)$ with $(X,x_0,x_0)$. If $A\neq \{x_0\}$,
then $\pi_m(X,A,x_0)$ is abelian for $m\geq 3$. Similarly to the bordism, for every $m\geq 2$ there is a natural
homomorphism
$$\partial: \pi_m(X,A,x_0)\to \pi_{m-1}(A,x_0), \ \partial (<f>)=<\partial f> \ . $$
Together with the homomorphisms 
$$i_*: \pi_m(A,x_0)\to \pi_m(X,x_0), \ \  j_*: \pi_m(X,x_0)\to \pi_m(X,A,x_0)$$
induced by the inclusions,
they give rise to the homotopy long exact sequence of the pointed pair $(X,A,x_0)$
$$\cdots \rightarrow \pi_m(A,x_0) \xrightarrow {i_*} \pi_m (X,x_0) \xrightarrow {j_*} \pi_m(X,A,x_0)\xrightarrow
{\partial} \pi_{m-1}(A,x_0)\rightarrow \cdots$$

\medskip

For every $m\geq 1$ it is well defined the map (in the oriented case we stipulate that $D^m$ inherits the standard orientation of $\R^m$)
$$ h_m: \pi_m(X,A,x_0)\to \Bb_m(X,A), \ h_m(<f>)= [D^m,S^{m-1},f]$$
obtained by ``forgetting the base points''. It is well defined because homotopy  is a special case of bordism
where only the cylinders are permitted. In fact

\begin{proposition}\label{h-omom} (1) For every $m\geq 1$, $h_m$ is a group homomorphism.

  (2) The family  of homomorphisms $\{h_m\}$ is {\rm functorial}, in a slogan: ``$g_* \circ h = h \circ g_*$'',
  and commutes with the respective long exact sequences.
\end{proposition}
\Dim Both the respective morphisms $g_*$ and long exact sequences have the very same definition on representatives.
Then (2) follows because the $h$'s are well defined. As for (1), for simplicity we consider the absolute case $m=1$, 
but the argument generalizes without difficulty. Realize an elementary bordism $W$ between $S^1 \coprod S^1$
and $S^1$ obtained by attaching a $1$-handle to $(S^1 \coprod S^1)\times [0,1]$ along
$(S^1 \coprod S^1)\times \{1\}$. There is a properly embedded arc $D\sim D^1$ (essentially the core of the handle)
which intersects $(S^1 \coprod S^1)\times \{0\}$ at two points belonging to different components
  and a properly embedded arc $D'$ dual to $D$ (essentially the co-core of the handle) which intersects
  the other component of $\partial W$ in two points. $W\setminus (D\cup D')$ is diffeomorphic
  to the cylinder $C:=((S^1\setminus \{p\}) \coprod (S^1\setminus \{p\})) \times [0,1]$. Let
  $f_0,f_1:(S^1,p)\to (X,x_0)$.  Up to the natural
  identification, this induces a map $F: C\to X$, $F(x,t):= f_0\coprod f_1(x)$ which extends to a 
  continuos map $F: W\to X$, by setting it constantly equal to $x_0$ on $D\cup D'$.
  This establishes a bordism between $(S^1,f_0)\coprod (S^1,f_1)$ and a determined map $g:S^1 \to X$.
  Recalling the definition of the operation on $\pi_1(X,x_0)$ (see \cite{Hatch}) it is immediate that 
  $$[S^1,g]= h_1(<f_0>\cdot <f_1>)$$ 
  hence 
  $$h_1(<f_0>\cdot <f_1>)= h_1(<f_0>)+h_1(<f_1>)$$ 
  as desired.

  \cvd

  \smallskip

  In general the study of both $\ker(h_m)$ and its image is a difficult question, even if $X$ is a compact
  smooth manifold. We can say something more for $m=1$.
\smallskip

  {\bf On the $1$-bordism.} It is evident that the homorphism
  $$\sigma_1: \Omega_1(X) \to \eta_1(X)$$ 
  is surjective: given $[M,f]$ in $\eta_1(X)$ it is enough to arbitrarily orient the components of $M$
(each diffeomorphic to $S^1$) to get $[\tilde M,f]$ in $\Omega_1(X)$ such that
$\sigma_1([\tilde M,f])=[M,f]$.  We have
\begin{proposition}\label{tau-onto} Assume that $X$ is path connected. Then the homomorphism $h_1:\pi_1(X,x_0)\to \Omega_1(X)$
    is surjective, hence the oriented bordism $\Omega_1(X)$ is a abelian quotient group of $\pi_1(X,x_0)$. By composition
with the surjective homomorphism $\sigma_1$, the same fact holds for $\eta_1(X)$.
\end{proposition}
\Dim Let $[S^1,f]\in \Omega_1(X)$. Let $p\in S^1$ the base point, $q=f(p)$.
Up to isotopy, hence up to bordism, we can assume that $f$ is constantly
equal to $q$ on a closed interval $p\in J \subset S^1$. Let $J=J_1\cup J_2$, $J_1\cap J_2=\{p\}$.
Let $\gamma_i: J_i\to X$ be a continuous path joining $q$ and the base point $x_0$ and such that $\gamma_i(p)=x_0$. Then
define $f':(S^1,p)\to (X,x_0)$ to be equal to $\gamma_i$  on $J_i$  and equal to
$f$ outside $J$. Clearly $[S^1,f']$ belongs to the image of $h_1$. We claim that $[S^1,f]=[S^1,f']$.  In fact it is not hard
to prove that they are homotopic. For a  general $[M,f]$ we can assune that  $M$ is union of a finite number of copies
$S^1_j$, of $S^1$. Consider the corresponding pointed copies $(S^1_j,p_j)$. Let $q_j=f(p_j)$.
By applying the above construction for every $j$, we can assume that $[M,f]$ is the sum of classes each one being the
image via $h_1$ of some $\alpha_j\in \pi_1(X,x_0)$. Finally by applying inductively on the number of components
the same argument used above to show that $h_1$ is a homomorphism,
we conclude that $[M,f]$ is the image of the product of such $\alpha_j$'s.

\cvd

\medskip   

We will complete the analysis of $\Omega_1(X)$ as a quotient of the fundamental group in Chapter \ref{TD-SURFACE}, 
Proposition \ref {omega1-ab}.
  
\section{Bordism categories}\label{b-cat} There is another important way to organize bordism matter.
As usual $\Mm_m$ either denotes $\Ss_m$ or $\Oo_m$, $\Bb$ either denotes $\eta$ or $\Omega$. 
For every $m\geq 0$, we define the {\it bordism category} ${\bf CAT}_\Bb(m+1)$.
\medskip

$\bullet$ $\Mm_m$ is the class of {\it objects} (recall that also $\emptyset$ is an object).
\smallskip

$\bullet$ For every couple of objects $M,N\in \Mm$, a {\it morphism} (``arrow'') 
$M \mapsto N$ is of the form
$$([\rho_0], [\rho_1], [W,V_0,V_1])$$
where $(W,V_0,V_1)$ is a triad of compact smooth manifolds (recall that $V_0$ and $V_1$
are union of components of $\partial W$, and $\partial W = V_0 \amalg V_1$) considered
up to diffeomorphisms which are isotopic to the identity on a neighbourhood of the boundary;
$\rho_0:M\to V_0$ and $\rho_1:N\to V_1$ are diffeomorphisms (preserving the orientation in the oriented setting)
considered up to isotopy.

\smallskip

$\bullet$ Two arrows $f:M\mapsto N$, $g: M'\mapsto N'$ can be composed if $N=M'$.
In such a case if $f=([\rho_0],[\rho_1], [W,V_0,V_1])$, $g=([\rho'_0],[\rho'_1],[W',V'_0,V'_1])$
then 
$$g\circ f=([\rho_0],[\rho'_1], [\tilde W, V_0,V'_1])$$ 
where
$$\tilde W = W\amalg_\psi W', \ \psi=\rho'_0 \circ \rho_1^{-1}: V_1\to V'_0 \ . $$ 
It is consistent because $\tilde W$ obtained by gluing is defined up to diffeomorphism
relatively to the boundary and only depends on the isotopy class of the gluing diffeomorphism. 
Note again that gluing can be performed in the oriented setting.
\smallskip

$\bullet$ For every object $M\in \Mm_m$, $M \neq \emptyset$, the {\it unit arrow} is 
$$1_M = ([{\rm id}_M], [{\rm id}_M], [M\times [0,1], M\times \{0\}, M\times \{1\}]) \ . $$  
\medskip

The discussion made in Chapter \ref{TD-HANDLE} about Morse functions on triads, dissections and  
handle-decompositions can be rephrased within the bordism category: every arrow is composition
of {\it elementary arrows} that is supported by triads admitting a handle decomposition with only one handle (of some index).

\section{A glance to TQFT}\label{quick-TQFT} 
A $(m+1)$ {\it topological quantum field theory} (TQFT) is a kind of non trivial
representation of ${\bf CAT}_\Bb(m+1)$ in the category of vector spaces
on some scalar field $K$. In last decades this has emerged as
a potent paradigm, the source of a plenty of  so called ``quantum invariants'' 
for $3$-dimensional manifolds and the right conceptual framework
for deep $4$-dimensional invariants.
The actual categorial definition involves many subtleties and is technically 
quite demanding (see for instance \cite{Tur}). Here we limit to a rough outline of the
main features. 

First we note that the objects $\Mm_m$ of a bordism category are endowed 
with the disjoint union operation ``$\amalg$''.

Let $K$ be a field and denote by $\Vv_K$ the category having as objects the
class $V_K$  of {\it finite dimensional} $K$-vector spaces and as morphisms the 
$K$-linear maps. Also $V_K$  is endowed with an operation ``$\otimes$"
given by the tensor product.

A  $(m+1)$ TQFT is a morphism of categories
$${\bf CAT}_\Bb(m+1) \  \Rightarrow  \Vv_K$$
which verifies certain conditions:
\begin{itemize}
\item To every object $M\in \Mm_m$ is associated
an object $Z(M)\in V_K$.
\item To every arrow $f:M\mapsto N$ in ${\bf CAT}_\Bb(m+1)$
is associated a linear map $Z(f): Z(M)\to Z(N)$, in such a way
that the composition is respected:
 $$ Z(g \circ f) = Z(g)\circ Z(f) \ . $$
 \item The correspondence $M\ \Rightarrow Z(M)$
 respects the operations on the objects:
 $$Z(M\amalg N)=Z(M)\otimes Z(N) \ . $$
 \end{itemize}
 
 \noindent  Moreover there are the following `non triviality requirements':
 \begin{itemize}
 \item $Z(\emptyset)=K$ (the space of ``states'' of the ``quantum'' empty set is non trivial).
 \item $Z(1_M) = {\rm id}_{Z(M)}$.
 \item $Z(M)$ is not constantly equal to $K$ and $Z(f)$ is not constantly
 equal to ${\rm id}_K$.
 \end{itemize}
 
 \noindent In the oriented setting, on $\Oo_m$ there is the involution $M\to -M$.
On $V_K$ there is the  duality ``involution'' $Z\to Z^*$ (where $Z$ is canonically 
identified with its bidual space $(Z^*)^*$).
Then here we require also 
\begin{itemize}
\item $Z(-M)=Z(M)^*$.
\end{itemize}       
          
\smallskip

One realizes quickly that the existence of such TQFT is not evident at all.
A possible attack could be to associate to all connected $M\in \Ss_m$ 
(possibly equipped with one fixed orientation) a same
vector space $Z(M)=V$ (so that $Z(-M)=V^*$ in the oriented setting). 
As every $M$ is the disjoint union of its connected
components,  $Z(M)$ is the tensor product of some copies of $V$ or $V^*$.
Then one could try to define first the elementary $Z(e)$ associated
to the elementary arrows in ${\bf CAT}_\Bb(m+1)$, perhaps in such a way that
they depend only on the handle index. A generic $Z(f)$ should be necessarily
a composition of such elementary morphisms. The key hard point is that
the decomposition by elementary arrows in   ${\bf CAT}_\Bb(m+1) $ is far
to be unique (as well as any triad supports a plenty of Morse functions)
but the resulting composite $Z(f)$ should not depend on the choice of
the decomposition. This means that our elementary $Z(e)$'s must verifies a huge 
collection of (a priori implicit) relations. 
For instace if we take $V=K^n$ for some $n$, $V^*=M(1,n,K)$, 
the unknown $Z(e)$'s in matrix form, we should find non trivial solutions of a huge system
of matrix equations.  It is not evident that such a solution exists (even if we take
$V=K$). 

\medskip

Every TQFT (if any) associates to every $M\in \Mm_{m+1}$, a scalar $\mu(M)$
which is {\it an invariant up to (possibly oriented) diffeomorphism}. In fact
as $M$ is compact and boundaryless, $(\emptyset,\emptyset, [M,\emptyset,\emptyset])$
is an arrow $f_{[M]}:\emptyset \mapsto \emptyset$, then $Z(f_{[M]}): K\to K$ and $\mu([M]):= Z(f_{[M]})(1)$.
\medskip

We will point out a ``baby'' (non trivial) TQFT in Chapter \ref{TD-EP}.

\chapter{Smooth cobordism}\label{TD-COBORDISM}
We specialize the bordism modules $\Bb_m(X,R)$ introduced in Chapter
\ref{TD-BORDISM} to $X$ which varies among the boundaryless compact
smooth manifolds. More precisely if $X$ is {\it not oriented} (even
{\it non orientable}), then we consider $\eta_m(X)=\Bb_m(X;\Z/2\Z)$,
if $X$ is {\it oriented}, we consider $\Omega_m(X)=\Bb_m(X;\Z)$.  A
first important fact, already used in Section \ref {b-nontriv}, is
that by means of the approximation theorems of continuous by smooth
maps, we can assume that all maps entering the definition of the
bordism modules are smooth; moreover, in dealing with functoriality we
can also assume that the maps $g:X\to Y$ are smooth. So all discussion
will have a differential/topolological character. The main issue of
this chapter is that by means of tranversality these ``smooth''
bordism modules (renamed ``cobordism'' modules up to a suitable
reindexing) can be embodied into {\it contravariant functors} and
their direct sum can be endowed with a functorial {\it graded ring}
structure.  This multiplicative structure is a substantial enrichement
of the theory and will lead to several important applications.

\section{Map transversality}\label{map-transv} We consider the following variant of 
the basic transversality setting  (Section \ref{basic-transverse}):
\smallskip

$\bullet$ All involved smooth manifolds admit an embedding in some $\R^n$ being furthermore
a closed subsets. This is certainly the case if a manifold is compact.

$\bullet$ All  involved smooth maps are {\it proper} (i.e. the inverse image of a compact set is compact).
Of course this is the case if the source manifold is compact. General topology tells us that proper maps 
between manifolds are {\it closed} (i.e. the image of a closed set is closed).
\smallskip

$\bullet$ $N$ and $Z$ are boundaryless smooth manifolds, $M$ is a compact smooth manifold with (possibly
empty) boundary $\partial M$.

$\bullet$ $f:M\to N$, $g:Z\to N$ are smooth maps.
\smallskip

\noindent In such a situation, we can define the product map
$$ (f\times g):M\times Z \to N\times N, \ (f\times g)(x,z)=(f(x),g(z))$$
and denote by $$\Delta_N=\{(y,y)\in N\times N\}$$ the {\it diagonal submanifold} of $N\times N$,
which is obviously diffeomorphic to $N$ by the canonical diffeomorphism 
$$N \to \Delta_N, \ y\to (y,y) \ . $$ 
Recall that $\partial (M\times Z)=\partial M \times Z$.

\begin{definition}\label{trans-map-defi}{\rm We say that $f$ is {\it tranverse to $g$} (and we write $f \pitchfork g$) if 
$(f\times g)\pitchfork \Delta_N$. This incorpotares that $\partial f\pitchfork g$.}
\end{definition}
\smallskip

By using that $T_{(y,y)} \Delta_N= \Delta_{T_yN} \subset T_yN\oplus T_yN$ one readily checks that:

\begin{lemma}\label{transv-lemma} $f\pitchfork g$ if and only if for every $(x,z)\in M\times Z$ such that $f(x)=g(z)=y$, then
$T_yN= d_xf(T_xM)+d_zg(T_zZ)$, and for every $(x,z)\in \partial M \times Z$ such $\partial f(x)=g(z)=y$,
then  $T_yN= d_x\partial f(T_x\partial M)+d_zg(T_zZ)$.
\end{lemma}

\cvd

We have the following version of the first transversality theorem:
\begin{theorem}\label{firstT2} In the given setting: 
\smallskip

(1) If $f\pitchfork g$ then 
$$(Y,\partial Y)=((f\times g)^{-1}(\Delta_N), (\partial f \times g)^{-1}(\Delta_N))$$
is a compact proper submanifold of $(M\times Z,\partial M \times Z)$. Moreover,
$$\dim (M\times Z)- \dim (Y)= \dim (N\times N)- \dim (N)= \dim (N) \ . $$

(2) If all involved manifolds are oriented, then $Y$ and $\partial Y$ are orientable
and we can fix an orientation procedure such that  $\partial Y$ becomes the oriented boundary of $Y$.

\end{theorem}
\Dim With the exception of the compactness of $Y$, all statements in (1) are direct consequence of
Theorem \ref{firstT} (and  they hold also without assuming that $g$ is proper). On the other hand, the compacteness
of $Y$ follows from the compactness of $M$ and the properness of $g$. Point (2) is a direct consequence 
of point (2) of Theorem \ref{firstT}, once $N\times N$ is endowed with the product orientation of two copies
of the given orientation on $N$, $\Delta_N$ is oriented in such a way that the canonical diffeomorphism
 is orientation preserving. 

\cvd

\begin{remark}{\rm  If $Z\subset N$ is a submanifold and $g$ is the inclusion,
then $$Y=\{(x,z)\in M\times Z; \ f(x)=z\}$$ that is the graph of the restriction of $f$ to $f^{-1}(Z)$.
If $Z$ is also a closed subset of $N$, then we are in the setting fixed above, and the projection
of $Y$ in $M$ is equal to $f^{-1}(Z)$ and is a proper submanifold of $(M,\partial M)$ recovering
the conclusion of Theorem \ref{firstT}.}
\end{remark}

\smallskip

We denote by $\pitchfork(M,N;g)$ the subspace of $\Ee(M,N)$
formed by the maps transverse to $g$. If $\partial f \pitchfork g$,
then we denote by  $\Ee(M,N,\partial f)$ (resp. $\pitchfork(M,N,\partial f;g)$)
the subspace of $\Ee(M,N)$ ($\pitchfork(M,N;g)$) formed by the maps that coincide with
$\partial f$ on $\partial M$. We have the following version of Theorem \ref{secondT}.

\begin{theorem}\label{secondT2}  In the given setting:

(1) $\pitchfork(M,N;g)$ is open dense in $\Ee(M,N)$.
\smallskip

(2) $\pitchfork(M,N,\partial f;g)$ is open dense in $\Ee(M,N,\partial f)$.
\smallskip

(3) For every $h \in \Ee(M,N)$ (resp. $h \in \Ee(M,N,\partial f)$) there is $\tilde h \in \ \pitchfork (M,N;g)$ ($\tilde h \in \ \pitchfork(M,N,\partial f;g)$) 
smoothly homotopic to $h$.
\end{theorem}
\Dim The proof is not a direct consequence of the {\it statement} of Theorem \ref{secondT} but it is a consequence of its proof
which can be adapted with minor changes.

\cvd

\smallskip

\section{Cobordism contravariant functors}\label{cb-functors}
Let $X$ be a compact boundaryless smooth manifold. Let $[M,f]\in \Bb_m(X;R)$ (either $R=\Z/2\Z$ or $R=\Z$ according
to the convention fixed at the beginning of the Chapter). Then we say that $[M,f]$ is of {\it codimension $k$  in $X$}
if
$$k= {\rm codim}_X[M,f]:= \dim (X) -m \ . $$
We can consider
the modules $\Bb_m(X;R)$ indexed by $\Z$ by stipulating that $\Bb_m(X;R)=0$ if $m<0$.  If $k$ is the codimension, set 
$$\Bb^k(X;R):= \Bb_m(X;R)$$
so we have a formal reidexing by $\Z$ of the family of bordism modules of $X$ in terms of the codimension, so that  
$\Bb^k(X;R)=0$ if $k> \dim X$.
To stress it we say that $\Bb^k(X;R)$ is the $k$-{cobordism module} of $X$ (over $R$). 
Formally for every $k\in \Z$, there are  tautological reindexing isomorphisms 
$$d: \Bb_{\dim (X)-k}(X;R)\to \Bb^k(X;R), \ D: \Bb^k(X;R)\to \Bb_{\dim (X)-k}(X;R)$$ 
$ d(\alpha)=D(\alpha)= \alpha $.

For every $k\in \Z$ we want  to enhance the object correspondence
$$ X \ \Rightarrow \  \Bb^k(X;R)$$
with a correspondence
$$ g:X\to Y \ \Rightarrow \ g^*: \Bb^k(Y;R)\to \Bb^k(X;R)$$
to build a contravariant functor from the category of compact boundaryless (possibly oriented) smooth manifolds 
and smooth maps to the category of $R$-modules and $R$-linear maps. Hence we want that
$$ (g\circ h)^*= h^* \circ g^*$$ whenever the composition makes sense, and
$$ {\rm id}_X^* = {\rm id}_{\Bb^k(X;R)} \ . $$
We have to define the induced linear maps $g^*$. We implement the following procedure,
basically it is the same ``pull-back'' construction that we have used for  vector bundles.

\smallskip

$\bullet$  If $k> \dim (Y)$, then $g^*: \{0\} \to \Bb^k(X;R)$ is uniquely determined. 
\smallskip

$\bullet$ Assume that $k\leq \dim (Y)$ and let  $\alpha \in \Bb^k(Y;R)$. Fix a representative
$$\alpha = [M,f] \ . $$ Hence $M$ is 
compact boundaryless (possibly oriented) of dimension $m=\dim (Y) - k$. By the transversality theorems,
up to homotopy hence up to bordism, we can assume that $f\pitchfork g$. Then $$V= (f\times g)^{-1}(\Delta_Y)$$
is a compact boundaryless (possibly oriented) submanifold of $M\times X$ such that 
$\dim (M\times X)-\dim(V)= \dim (Y)$, that is 
$$\dim (X)-\dim(V) = \dim (Y) - \dim (M) = k  \ . $$
Hence $[V,p_X]\in \Bb^k(X;R)$, where $p_X$ is the restriction of the projection $M\times X \to X$.

We have
\begin{proposition}\label{g*} Let $g:X\to Y$ be a smooth map between compact boudaryless (possibly
oriented) smooth manifolds. Let $\alpha  \in \Bb^k(Y;R)$. Let $[V,p_X]\in \Bb^k(X;R)$
obtained by means of any implementation of the above ``pull-back'' procedure starting from a representative $\alpha=[M,f]$. 
Then 
 
 (1) The map
 $$g^*:\Bb^k(Y;R)\to \Bb^k(X;R), \ g^*(\alpha)=[V,p_X]$$ 
 is well defined   (it does not depend on the arbitrary choices of a given
 implementation).   
 
 (2) $g^*$ is $R$-linear.
 
 (3) For every $X$ 
 $$ {\rm id}_X^* = {\rm id}_{\Bb^k(X;R)} \ . $$

(4) Whenever the composition makes sense $$ (g\circ h)^*= h^* \circ g^* \ . $$ 

(5) Let $n=\dim X$; if $[X,g_0]=[X,g_1]\in \Bb_n(Y;R)$, then $g^*_0=g^*_1$.
In particular this holds if $g_0$ and $g_1$ are homotopic.
\end{proposition}
\Dim Assume that $g^*$ is well defined and prove items (2)-(4). The procedure
distributes on the addends of a dijoint union so (2) follows easily.

As for (3) Every $[M,f]$ is tranverse to ${\rm id}_X$, hence $V$ is the graph of $f$
and clearly $[V,p_X]=[M,f]$.

Concerning (4), If $g^*([M,f])=[M',f']$, $h^*([M',f'])=[M",f"]$ the representatives being
obtained by iterated application of the pull-back procedure, then $f"\pitchfork (g\circ h)$
and $[M",f"]$ results from an implementation of the procedure applied to
$[M,f]$ and $g\circ h$.

Let us show now (1), that is $g^*$ is well defined. Let $(V,p_X)$ and $(V',p'_X)$ be obtained
by implementing the procedure starting from representatives $(M,f)$ and $(M',f')$,
$f\pitchfork g$, $f'\pitchfork g$; let $(W,F)$ realizes a bordism of $(M,f)$ with $(M',f')$.
By applying the transversality theorems we can assume that $F\pitchfork g$. 
Then $((F\times g)^{-1}(\Delta_Y), P_X)$ realizes a bordism of $(V,p_X)$ with
$(V',p'_X)$.

Finally (5) follows by the very similar argument used for (1): 
if $(W,F)$ realizes a bordism of $(X,g_0)$ with $(X,g_1)$, then we can assume that $F$
verifies suitable transversality conditions, so that $(f\times F)^{-1}(\Delta_Y)$ 
leads to a bordism of $((f\times g_0)^{-1}(\Delta_Y), p_X)$ with $((f\times g_1)^{-1}(\Delta_Y),p_X)$.

\cvd

\subsection{Reduction mod$(2)$}\label{red-mod-2}
When $X$ is oriented, we already known the natural ``forgetting'' homomorphisms 
$$\sigma: \Bb^k(X;\Z)\to \Bb^k(X;\Z/2\Z) \ . $$
These are functorial, that is 
\begin{proposition} For every  smooth map $g:X\to Y$ between oriented
compact boundaryless manifolds, for every $\alpha \in \Bb^k(Y;\Z)$
then $g^*(\sigma(\alpha))= \sigma(g^*(\alpha))$, where the first $g^*$ refers
to the $\Z/2\Z$-cobordism, the second to the $\Z$-cobordism.
\end{proposition}
\Dim The construction of $g^*(\sigma (\alpha))$ is obtained by the one
of $g^*(\alpha)$ just by forgetting the orientation.

\cvd

\section{The cobordism cup product}\label{cb-prod} Let $X$ be as above.
For every $r,s \in \Z$, we are going to define a bilinear map
$$\sqcup: \Bb^r(X;R)\times \Bb^s(X;R) \to \Bb^{r+s}(X;R) \ . $$
Let us describe the procedure that defines this ``cup'' product.
\smallskip

$\bullet$ If at least one among $r$ and $s$ is bigger than $\dim (X)$, then $\alpha \sqcup \beta = 0$.
\smallskip

$\bullet$ Let $(\alpha, \beta) \in \Bb^r(X;R)\times \Bb^s(X;R)$ and assume that both $r$ and $s$ are
$\leq \dim (X)$. Fix representatives $\alpha=[M,f]$ and $\beta=[N,h]$. We claim that 
$$[M\times N, f\times h]\in \Bb^{r+s}(X\times X; R) \ . $$
In fact
$$2\dim (X) - (\dim(M) + \dim (N))= 2\dim (X) - ( \dim (X) -r + \dim (X) -s ) = r+s \ . $$
\smallskip

$\bullet$ Let $\delta_X: X\to X\times X, \ \delta_X(x)=(x,x)$ be the canonical diffeomorphism onto the diagonal
$\Delta_X$. Finally take 
$$\delta_X^*[M\times N, f\times h]\in \Bb^{r+s}(X;R) \ . $$
We stress that we are actually using the contravariant nature of the cobordism functors.
 
\begin{remark}\label{cup-rep}{\rm If $f\pitchfork h$ we can explicitly describe representatives
of  $\delta_X^*[M\times N, f\times h]$. In fact in such a case $(f\times h)\pitchfork \delta_X$.
Then $\delta_X^*[M\times N, f\times h]= [\tilde V,p_X]$ where
$$\tilde V=\{(x,p,q)\in X\times M \times N;\ f(p)=h(q)=x\} \ . $$
Let 
$$V = (f\times h)^{-1}(\Delta_X)=\{(p,q)\in M\times N; \ f(p)=h(q)\} \ . $$
Then $\tilde V$ is the graph of $u:=f_{|V} = h_{|V}$, $V$ and $\tilde V$ are canonically diffeomorphic,
and $$[\tilde V,p_X]=[V, u]\in \Bb^{r+s}(X;R) \ . $$ In particular if $f$ and $h$ are the inclusions of
two transverse submanifolds $M$ and $N$ of $X$ and $j$ is the inclusion of $M\pitchfork N$,
then 
$$\delta_X^*[M\times N, f\times h]= [M\pitchfork N,j]\ . $$
}
\end{remark}

We have

\begin{proposition}\label{cup} Let $X$ be a compact boundaryless (possibly oriented) smooth manifolds.
Let $(\alpha,\beta)\in \Bb^r(X;R)\times \Bb^s(X;R)$, $\delta_X^*[M\times N, f\times h]\ \in \Bb^{r+s}(X;R)$ 
be obtained by any implementation of the above procedure applied to arbitrary representatives $\alpha=[M,f]$, $\beta=[N,h]$.
Then:
\smallskip

(1) The class $\alpha \times \beta :=[M\times N, f\times h]$, whence the map
$\alpha \sqcup \beta := \delta_X^*[M\times N, f\times h]\ $ 
are well defined  
(they do not depend on the arbitrary choices of a given
 implementation).
 
 (2) $\sqcup$ is bilinear.
 
 (3) For every $(\alpha,\beta)\in \Bb^r(X;R)\times \Bb^s(X;R)$, 
 $$ \alpha \sqcup \beta = (-1)^{rs}\beta \sqcup \alpha \ . $$ 
 
 (4) $\sqcup$ is {\rm functorial}, that is for every $g:X\to Y$, for every $(\alpha, \beta)\in  \Bb^r(Y;R)\times \Bb^s(Y;R)$,
 $$g^*(\alpha)\sqcup g^*(\beta)= g^*(\alpha \sqcup \beta) \ . $$
 
 \end{proposition}
 \Dim Again assume that $\sqcup$ is well defined and prove the other items. By the transversality theorems
 the assumption allows us to use representatives which verify all suitable transversality conditions.   
 The disjoint union distributes to the product
 of manifods; (2) follows easily. Item (3) is a local verification and reduces to Remark \ref{orient-int}. Let $(M,f)$, $(N,h)$
 be representatives of $\alpha$ and $\beta$ such that $f\pitchfork g$, $h\pitchfork g$ and $f\pitchfork h$. It follows
 that $(g\times g)\circ \delta_X \pitchfork (f\times h)$. By combining the two procedures that define $g^*$ and $\sqcup$
 starting from such representatives in general position we obtain representatives for both terms of the equality of (4) that are evidently
 bordant to each other (in the same spirit of Remark \ref{cup-rep}). It remains to prove that $\sqcup$ is well defined. 
 As $\delta^*_X$ is well defined, it is enough
 to show that 
 $$\alpha \times \beta:=[M\times N,f\times h]\in \Bb^{r+s}(X\times X;R)$$ only depends on the class $\alpha$ and $\beta$.
 By symmetry it is enough to show that it does not depend on the choice of a representative of $\alpha$.
 If $(W,F)$ realizes a bordism of $(M,f)$ with $(M',f')$ then $(W\times N, F\times h)$ realizes a bordism 
 of $(M\times N, f\times h)$ with $(M'\times N, f'\times h)$.

\cvd

\subsection{Reduction mod$(2)$} Similarly to Proposition \ref {red-mod-2} we have
\begin{proposition}\label {red-mod-2-2}  For every  compact oriented boundaryless manifold $X$,
for every $(\alpha,\beta)\in \Bb^r(X;\Z)\times \Bb^s(X;\Z)$, 
$\sigma(\alpha)\sqcup \sigma(\beta)= \sigma(\alpha \sqcup \beta)$,
where the first $\sqcup$ refers
to the $\Z/2\Z$-cobordism, the second to the $\Z$-cobordism.
\end{proposition}

\Dim The construction of $\sigma(\alpha)\sqcup \sigma(\beta)$ is obtained by the one
of $\alpha \sqcup \beta$ just by forgetting the orientation.

\cvd

\subsection {The cobordism ring}\label{cb-ring}
The collection of the above cup products gives a globally defined  product
$$\sqcup: \Bb^\bullet(X;R)\times \Bb^\bullet(X;R)\to \Bb^\bullet(X;R)$$
on the direct sum $R$-module
$$ \Bb^\bullet (X;R):=\oplus_{k\in \Z} \Bb^k(X;R) \ . $$
$(\Bb^\bullet(X;R),+,\sqcup)$ is called the {\it graded $R$-cobordism ring of $X$}
(it is a graded algebra when $R=\Z/2Z$). 
\smallskip

Similarly the collection of above $g^*$'s  defines a global graded ring homomorphism 
$$ g^*: \Bb^\bullet(Y;R)\to \Bb^\bullet(X;R) \ . $$ 
We can summarize the above achievements as follows:
$$ X \ \Rightarrow \Bb^\bullet(X;R)$$
$$ g:X\to Y \  \Rightarrow g^*:  \Bb^\bullet(Y;R)\to \Bb^\bullet(X;R) $$
{\it define a contravariant functor from the category of compact boundaryless (possibly oriented) smooth
manifolds and smooth maps to the category of graded rings and graded ring homomorphisms.}
\smallskip

\begin{remark}{\rm A graded ring verifying the {\it non} commutative
rule  (3) in Proposition \ref{cup} is sometimes called a ``commutative'' graded ring.}
\end{remark}

\begin{remark}\label{ring-point}{\rm A particular case of the above constructions is when
$X$ is reduced to one point. In this case the product
$$ \Bb^r(R) \times \Bb^s(R) \to \Bb^{r+s}(R)$$
for every couple of indices $r,s\leq 0$ is just defined by the product of representatives
$$ [M]\sqcup [N] = [M\times N] \ . $$}
\end{remark}

 \begin{remark} \label{non-compact} (Non compact $X$) {\rm Referring to the setting of the tranversality 
 theorems of  Section \ref{map-transv}, we can extend the range of
 cobordism functors and product to the category of boundaryless possibly non compact 
 manifolds $X$  but which can be embedded anyway in some $\R^k$ being also a closed subset,
 and smooth {\it proper} maps between these manifolds.}
 \end{remark}
 
\section{Duality, intersection forms}\label{intersection-form} Assume that $X$ is connected (possibly oriented), 
 $\dim (X)=n$.
 Then $$\Bb^n(X;R)\sim \Bb_0(X;R) \sim R \ . $$ If $R=\Z/2\Z$, we have a generator 
 $\beta_X$ of $\Bb^n(X;\Z/2\Z)$ represented as
 $\beta_X=[x,i]$ where $x \in X$ and $i$ is the inclusion (it does not depend on the choice
 of $x$ because $X$ is path connected). If $R=\Z$ we have two generators of
 the form $[\pm x, i]$. As usual we encode the point sign by associating to $+x$  the 
 orientation on $T_xX$ carried by the global orientation of $X$ and this selects again one
 generator $\beta_X$. By this choice
 of generators we have fixed in both cases an identification of $\Bb^n(X;R)$ with $R$.
 
 For every $r,s$, set $p=n-r$, $q=n-s$.
 Let $r,s$ be such that $r+s=n$ (hence also $p+q= n$, $p=s$, $q=r$). Then
 $$\sqcup : \Bb^r(X;R)\times \Bb^s(X;R) \to R \ . $$
 Note in particular that
 $$  d(\alpha_X)\sqcup \beta_X=1$$
 where $\alpha_X=[X,{\rm id}_X]\in \Bb_n(X;R)$ is the bordism fundamental class of $X$
 and $d:\Bb_n(X;R)\to \Bb^0(X;R)$ is the tautological isomorphism.
 
 By using the tautological isomorphisms, all this can be lifted to a bilinear map
 $$\bullet: \Bb_{p}(X;R)\times \Bb_{q}(X;R)\to R$$
 or to a bilinear pairing
 $$ \sqcap: \Bb^r(X;R)\times \Bb_{q}(X;R)\to R \ . $$
 This last induces a linear map ($q=r$)
 $$ \phi^r: \Bb^r(X;R)\to {\rm Hom}(\Bb_{r}(X;R),R), \ \gamma\to \phi_\gamma, \ \phi_\gamma(\sigma)=\gamma \sqcap \sigma \ . $$
 Recall that by applying the Hom functors we have a basic way to convert the covariant bordism funtors into cotravariant ones
 $$X \ \Rightarrow {\rm Hom}(\Bb_m(X;R), R)$$
 $$ g:X\to Y \ \Rightarrow g^t_*: {\rm Hom}(\Bb_m(Y;R),R)\to {\rm Hom}(\Bb_m(X;R),R)$$
 $g^t_*(\gamma)= \gamma \circ g_*$.
 The homomorphisms $\phi^r$, $g^t_*$ and $g^*$ are compatible; in a slogan:
 ``$ \phi^r \circ g^* = g^t_* \circ \phi^r$''
 
 The map $\phi^r$ is in general not injective nor surjective. A reason is the possible existence of non trivial submodules
 of $\Bb_*(X;R)$ isomorphic to $\Bb_*=\Bb_*(\{x_0\};R)$. The image  via the tautological
 isomorphism of such  submodule in $\Bb^r(X;R)$ is contained in the kernel of $\phi^r$. 
 If $R=\Z/2\Z$, so that $\Bb_{r}$ can be realized as a direct addend
 of $\Bb_{r}(X;\Z/2\Z)$, then any functional $\gamma$ which holds $1$ on $\Bb_{r}$ and such that 
 $\Bb_{r}(X;\Z/2\Z)= \Bb_{r}\oplus \ker \gamma$ does not belong to the image of $\phi^r$.
 If $R=\Z$, then the {\it torsion submodule} of $\Bb^r(X;\Z)$ is contained in the kernel of $\phi^r$.
 For every $r$ we set
 $$\Hh^r(X;R):=   \Bb^r(X;R)/\ker (\phi^r) $$
 and extending the usual reindexing set
 $$ \Hh_{n-r}(X;R):= \Hh^r(X;R)$$
 where in this last equality only the  $R$-module structure is considered,
 forgetting the multiplicative structure. Then the above map $\phi^r$ induces
 an injective $R$-linear map
 $$ \hat \phi^r: \Hh^r(X;R)\to {\rm Hom}(\Hh_{r}(X;R),R) \ . $$
   
 If $X$ is connected (possibly oriented), then
 $$\Hh^0(X;R) \sim R$$
 and is generated by the  fundamental class. 
 The map $\sqcap$ can be formally generalized by composing $\sqcup$ with the tautological isomorphisms
 $$ \sqcap : \Bb^r(X;R)\times \Bb_q(X;R) \to \Bb_{2n-(r+q)}(X;R) \ . $$
 In particular
 $$ \sqcap: \Bb^r(X;R)\times \Bb_n(X;R)\to \Bb_{n-r}(X;R)$$
 and it is a consequence of the definitions that for every $\sigma \in \Bb^r(X;R)$
 $$\sigma \sqcap \alpha_X = D(\sigma) \ . $$
 
 If $\dim (X)=2m$ then we can consider
 $$ \sqcup: \Bb^m(X;R)\times \Bb^m(X;R)\to R$$
 or equivalently
 $$\bullet : \Bb_m(X;R)\times \Bb_m(X;R)\to R $$
 this second is also called the {\it $R$-bordism intersection form of $X$}.
 Note that these forms are symmetric on $\Z/2\Z$, while on $\Z$ they are
 symmetric (resp. antisymmetric) if $m$ is even  ($m$ is odd).  The kernel
 of $\phi^r$ coincides in this case with the {\it radical} of the form, hence
 the induced form (also called ``intersection form")
 $$ \sqcup: \Hh^m(X;R) \times \Hh^m(X;R)  \to R $$
 determines an inclusion of  
 $\Hh^m(X;R)$ as a submodule of its
 dual module
 $$\hat \phi^m: \Hh^m(X;R)\to {\rm Hom}(\Hh_m(X;R),R) \ . $$

 \section{Cobordism theory for compact manifolds with boundary}\label{relative-cb}
 First let us strengthen the notion of map between pairs of spaces 
 $$h:(X,A)\to (Y,B)  \ ; $$
 it is a {\it strict} pair map if (as usual) $h(A)\subset B$ and furthermore $h(X\setminus A)\subset Y\setminus B$.
 
 We consider the category of compact smooth (possibly oriented) manifolds with (possibly empty)
 boundary $(X,\partial X)$ and smooth  strict pair maps $$h: (X,\partial X)\to (Y, \partial Y) \ . $$
  For example the inclusion of a {\it proper} submanifold $(M,\partial M)$ in $(X,\partial X)$ is a 
 typical example of strict  map. We stress that a strict map
 $f:(M,\emptyset)\to (X,\partial X)$ sends the boundaryless $M$ in the interior Int$(X)$ of $X$.
 
 The non compact manifold Int$(X)$ verifies the conditions of Remark \ref{non-compact};
 for example if $X\subset \R^k$ for some $k$ (this is possible because $X$ is compact) and 
 $h:\R^k\to \R$ is a non negative smooth function such that $\partial X = h^{-1}(0)$,
 then the restriction to $X\setminus \partial X$ of  $\R^k\setminus \partial X \to \R^{k+1}, \ x\to (x,1/h(x))$
 is an embedding of Int$(X)$ onto a closed subset of $\R^{k+1}$. 
 
 The usual definitions of the absolute or relative bordism modules $\Bb_m(X;R)$ or $\Bb_m(X,\partial X;R)$ can be enhanced 
 by stipulating that all involved pair maps are smooth and strict. By using the approximation
 theorem of continuous maps by smooth maps and the boundary collars to push into the interior what is
 necessary in order to make strict any given ``singular'' smooth manifold in $X$ or in $(X,\partial X)$,  it is not hard check
 that: 
 \smallskip
 
 {\it These enhanced modules are actually isomorphic to the original ones and moreover,
 $\Bb_m(X;R)$ is naturally isomorphic to $\Bb_m({\rm Int}(X);R)$.}
 \smallskip
 
  The reindexing $\Bb^k(X;R) = \Bb_m(X;R)$ or  
 $\Bb^k(X,\partial X;R)=\Bb_m(X,\partial X;R)$, $k= \dim (X)-m$
 is made as usual with respect to the codimension in $X$.
 \smallskip
 
Let $g: (X,\partial X)\to (Y,\partial Y)$ be a smooth strict map in our category. We want
to extend the definition of the induced linear morphism 
$$g^*:  \Bb^k(Y,\partial Y;R)\to \Bb^k(X,\partial X;R) \ . $$
For every strict pair map $h: (N,\partial N)\to (Y,\partial Y)$ we denote as usual $$\partial h: \partial N \to Y$$
the restriction of $h$ to the boundary; then set 
$$\partial \partial h: \partial N \to \partial Y$$
such that $\partial h = j\circ \partial \partial h$, where $j$ is the inclusion of $\partial Y$ in $Y$.
Close to Lemma \ref{transv-lemma}, 
we say that $(f,\partial f)\pitchfork (g,\partial g)$ if and only if
\begin{enumerate}
\item  for every $(p,x)\in M\times X$ such that $f(p)=g(x)=y$, $T_yY= d_pf(T_pM)+d_xg(T_xX)$;
for every $(p,x)\in \partial M\times X$ such that $\partial f(p)=g(x)=y$, $T_yY= d_p\partial f(T_p\partial M)+d_xg(T_xX)$;
coherently with the notations of Section \ref{map-transv}, we summarize this item by ``$f\pitchfork g$'';
\item $g\pitchfork f$;
\item $\partial \partial f \pitchfork \partial \partial g$ (in the usual sense).
\end{enumerate}

Let $(M,\partial M,f)$ be a smooth and strict representative of a given $\alpha \in \Bb^k(Y,\partial Y;R)$.
By suitably and straighforwardly adapting the transversality theorems, we can assume that
$(f,\partial f)\pitchfork (g,\partial g)$. Set $V=\{(p,x)\in M\times X; \ f(p)=g(x)\}$. Then 
$$g^*(\alpha):= [V,\partial V,p_X]$$
well defines our desired linear map $g^*$.
\smallskip

Now, by formally using the very same definition given when $X$ is boundaryless, 
we (partially) extend the cup product as follows:
$$\sqcup: \Bb^r(X,\partial X;R)\times \Bb^s(X;R) \to \Bb^{r+s}(X;R)$$
$$\sqcup:\Bb^r(X;R)\times \Bb^s(X;R)\to \Bb^{r+s}(X;R) \ . $$
Then we have a linear map
$$ \phi^r: \Bb^r(X,\partial X;R)\to {\rm Hom}(\Bb_{r}(X;R),R)$$
which restricts to (we keep the same name)
$$ \phi^r: \Bb^r(X;R)\to {\rm Hom}(\Bb_{r}(X;R),R) \ . $$
Finally we have the induced injective map
$$ \hat \phi^r: \Hh^r(X,\partial X;R) \to {\rm Hom}(\Hh_{r}(X;R),R) \ . $$
 
 \chapter{Applications of cobordism rings}\label{TD-CB-APPL}
In this chapter we will see several, sometimes very classical, applications of the cobordism theory,
especially of its multiplicative structure.

\section{Fundamental class revised, Brouwer's fixed point Theorem}\label{fund-class-2}
Here we recover Proposition \ref{fund-class} in terms of cobordism. Let $X$ be a boundaryless connected (possibly oriented)
smooth $n$-manifold. Let $[X,{\rm id}_X]\in \Bb^0(X;R)$ (often we will simply write $[X]$). Let $\beta_X\in \Bb^n(X;R)$ the generator given in
Section \ref{intersection-form} in order to fix an identification $\Bb^n(X;R)=R$. We have already remarked that
$$ [X] \sqcup \beta_X = 1\in R$$
hence in particular $[X]\neq 0$.  On the other hand, if $\gamma$ belongs to the image via the tautological
isomorphism $d: \Bb_n(X;R)\to \Bb^0(X;R)$ of the natural submodule isomorphic to $\Bb_n$, then
$$\gamma \sqcup \beta_X=0$$
hence $[X] \neq \gamma$. If $X$ has non empty boundary $\partial X$, we can consider 
$$[X,\partial X] \in \Bb^0(X,\partial X;R)$$
and we have again
$$ [X,\partial X] \cup \beta_X = 1 \in R \ . $$
The following is a very classical topological application of such a fundamental class.
\begin{theorem}\label{brouwer} {\rm (Brouwer fixed point theorem)} For every continuous map $$f:D^n \to D^n$$
there is $x\in D^n$ such that $f(x)=x$.
\end{theorem}
\Dim  The case $n=0$ is trivial. For $n>0$, assume that there is such an $f$ without any fixed point. Define 
$F: D^n \to S^{n-1}$ by setting
$F(x)$ equal to the unique point of intersection between $S^{n-1}=\partial D^n$ and the ray emanating from
$f(x)$ and passing through $x$. As $f$ is continuous, it is easy to verify that also $F$ is continuous and that 
$\partial F = {\rm id}_{S^{n-1}}$. Hence $[S^{n-1}]$ should be trivial in $\Bb_{n-1}(S^{n-1})$ against 
Proposition \ref{fund-class}.

\cvd

\section{A separation theorem}\label{sep-teo}
It is evident that an equatorial $S^{n-1}\subset S^n$ divides this last into two connected components.
If $n\geq 2$, every connected hypersurface in $S^n$ shares the same behaviour.
  
\begin{proposition} (1) Let $M\subset S^n$ be a compact boundaryless connected submanifold, $\dim (M)=n-1$,
$n\geq 2$.  Then $S^n\setminus M$ has exactly two connected components $W$, $W'$
and the closures are compact submanifolds with boundary such that
$\partial \bar W = \partial \bar W'= M$.

(2) Let $M\subset \R^n$ be a compact boundaryless connected submanifold, $\dim (M)=n-1$, $n\geq 2$.
Then $\R^n\setminus M$ has two connected components, one say $W$ has compact closure
and $\partial \bar W = M$.
\end{proposition}
\Dim The item (2) follows from (1) by considering $\R^n\subset \R^n\cup \infty = S^n$, such that
$\infty$ does not belong to $M$. As for (1), we know by Section \ref{b-nontriv} that 
$[M]:=[M,i_M]\in \Bb_{n-1}(S^n;\Z/2\Z)\sim \Bb^1(S^n;\Z/2\Z)$ ($i_M$ being the inclusion) belongs 
to the submodule isomorphic to $\Bb_{n-1}$. Hence we know that $[M]$ belongs to the kernel of the map
$\phi: \Bb^1(S^n;\Z/2\Z)\to {\rm Hom}(\Bb_1(S^n;\Z/2\Z),\Z/2\Z)$. Assume that $S^n\setminus M$ is
connected. Take a small simple arc $\gamma$ intersecting transversely $M$ at one point.
The endpoints of $\gamma$ belong to $S^n\setminus M$, hence $\gamma$ can be completed to
a smooth simple curve $\hat \gamma$ in $S^n$ that intersects transversely $M$ at one point. It follows
that $\phi_{[M]}([\hat \gamma, i_{\hat \gamma}])=1$ and this is a contradiction. Hence
$S^n\setminus M$ is not connected.  A tubular neighbourhood $U$ of $M$ in $S^n$ is diffeomorphic
to $M\times (-1,1)$, in fact  $M\times [0,1)$ can be identified with a collar of $M$ in $\bar W$, where
$W$ is a component of $S^n\setminus M$. Since $ U\setminus M$ has evidently two connected components,
then $S^n\setminus M$ has at most two components and this achieve the proof.

\cvd

\section{Intersection numbers}\label{inter-numb}
Let $X$ be a compact connected (possibly oriented) boundaryless smooth $n$-manifold.
Let $M$ and $N$ be compact boundaryless (possibly oriented) submanifolds of $X$, $\dim M= p$, $\dim N= q$.
Assume that $p+q=n$. Then 
$$ [M]\bullet [N] \in R$$
is the $R$-{\it intersection number} of the two submanifolds. Obviously it is invariant up to isotopy of $M$ or $N$
in $X$ (isotopy is a  particular instance of bordism). Hence if $[M]\bullet [N] \neq 0$, then there is no any isotopy
that makes $M$ and $N$ apart. In particular, if $M=N$ (hence $n=2m$), then $M\bullet M$ is called the {\it self-intersection number}
of the submanifold $M$.

\subsection{Lefschetz's number and fixed point theorem}\label{Lefscetz} 
Let $X$ be as above a connected compact boundaryless $n$-manifold.
Let $f:X\to X$ be a smooth map.
Consider the submanifolds $\Delta_X$ and $G(f)$ of $X\times X$, $G(f)$ being the graph of $f$. If $n= \dim (X)$,
then 
$$L_2(f):=[\Delta_X] \sqcup [G(f)] \in \Bb^{2n}(X\times X;\Z/2\Z)=\Z/2Z $$
is called {\it the Lefschetz number of $f$} mod$(2)$.
This is invariant (in particular) if $f$ is considered up to homotopy. As usual this allows us to define this number
also when $f$ is merely a continuous map. It is clear that if $\Delta_X \cap G(f)= \emptyset$ (that is if $f$ has no
fixed points), then $L_2(f)=0$. Viceversa we have the following ``fixed point theorem":
\smallskip

{\it If $L_2(f)\neq 0$, then $f$ has a fixed point.}

\smallskip

\noindent If $M$ is oriented, we can define the Lefschetz number 
$$L(f)\in \Bb^{2n}(X\times X;\Z)=\Z, \ L(f)=L_2(f) \ {\rm mod}(2)$$
in the oriented setting, and repeat verbatim the above considerations.

\section{Linking numbers}\label{linking}
Let $X$ be as in Section \ref{inter-numb}, $n\geq 3$. Let $(M,\partial M)$ be a $(n-k)$-compact submanifold (possibly oriented) of $X$
with non empty boundary, $n-k\geq 1$. $M$ is called a $R$-{\it Seifert surface of $T=\partial M$ in $X$}. 
Let $U$ be a ``small'' tubular neigbourhood of $T$ in $X$, such that $\partial U \pitchfork M$.
The closure $(Y,\partial Y)$ ($\partial Y = \partial U$) of $X\setminus U$ is a compact $n$-manifold with non empty boundary; the closure
$(N,\partial N)$ of $M\setminus U$ is a proper $(n-k)$-submanifold of $(Y,\partial Y)$. Then $[N,\partial N]\in \Bb^k(Y,\partial Y;R)$
(we omit to indicate the inclusion map).
Let $Z$ be a compact boundaryless (possibly oriented) proper $k$-submanifold of $(Y,\partial Y)$. Hence $[Z]\in \Bb^{n-k}(Y;R)$. 
Then
$$ lk_{M}(T,Z):= [N,\partial N]\sqcup [Z] \in R$$
is called the $R$-{\it linking number of $Z$ with $T$ with respect to the Seifert surface $M$}.
By the uniqueness  of tubular neighbourhoods up to isotopy, it is well defined. Moreover, it is invariant
up to isotopy $Z$ in Int$(Y)$. In some case the linking number does not depend on the choice of the Seifert
surface. For example we have
\begin{proposition}\label{link-in-sphere} In the above setting, assume that $X=S^n$. Then
$$lk(T,Z):= lk_M(T,Z)\in R$$
is well defined, that is it does not depend on the choice of a Seifert surface of $T$ in $S^n$.
\end{proposition}
\Dim Let $T=\partial M= \partial M'$. By (abstractly) gluing $M$ and $M'$ along $T$ and taking the union of the inclusions,
we get say $[W,f]\in \Bb^k(S^n;R)$. Let us consider $[Z]\in \Bb^{n-k}(S^n)$.
We have already noticed that
$$ [W,f]\sqcup [Z]=0\in R \ . $$
On the other hand, it follows from the very geometric definition of the cobordism cup product that 
$$ [W,f]\cup[Z] = lk_W(T,Z)-lk_{W'}(T,Z)$$
and the Proposition follows.

\cvd

\begin{remark}\label{3D-link}{\rm A classical example of linking number is the case $X=S^3$ and $T$, $Z$ (possibly oriented)
disjoint {\it knots} in $S^3$ (that is disjoint submanifolds diffeomorphic to $S^1$). It is a classical well-known fact
(see \cite{Rolf}) that a knot in $S^3$ admits a Seifert surface. Hence we eventually have
$$lk(T,Z) = lk(Z,T)\in R \ . $$
Another classical situation is when $X=S^n$, $T\sim S^p$, $Z\sim S^q$ and these last are {\it unknotted spheres}
in $S^n$, that is they are the boundary of embedded $(p+1)$ or $(q+1)$ smooth disks respectively.} 
\end{remark}

\section{Degree}\label{degree}
Let $X$ and $Y$ be compact connected boundaryless (possibly oriented) smooth $n$-manifolds,
$g:X\to Y$ be a continuous map. Let us fix generators $\beta_X$ of $\Bb^n(X;R)=R$ and $\beta_Y$
of $\Bb^n(Y;R)=R$ as in Section \ref{intersection-form}. Consider
$$ g^*: \Bb^n(Y;R)\to \Bb^n(X;R)$$
then define the $R$-{\it degree} of $g$ by:
$$ \deg_R(g):=g^*(\beta_Y) \in R \ . $$
Although we have already given an operative definition of $g^*$ in full generality, it is convenient to spell it again in the present
situation: fix $y_0\in Y$; up to homotopy make $g$ smooth and transverse to $y_0$ (equivalently move a little $y_0$ to make it
a regular value of $g$); then $g^{-1}(y_0)=\{x_1,\dots x_r\}$ is a finite set of points; in the oriented setting they are oriented,
that is endowed with signs $\epsilon_j$, $j=1,\dots, r$; on $R=\Z/2Z$ the degree is equal to $r$ mod$(2)$; on $\Z$
the degree is the sum of the signs $\epsilon_j$.
\smallskip

Now we list a few properties of the degree.

\smallskip

\noindent $\bullet$ If $g$ is not surjective, then $\deg_R g = 0$.
\smallskip

\noindent $\bullet$ If $g:X \to Y$ is a diffeomorphism, then $\deg_R(g)=\pm 1$.
\smallskip

\noindent $\bullet$ If $h\circ g$ and the the degrees of all involved maps make sense, then
$$ \deg_R(h\circ g)= \deg_R(h)\deg_R(g)$$
that is the degree is multiplicative under composition. This follows immediately from
functoriality.
\smallskip

\noindent $\bullet$ If $g$ and $h$ are homotopic, then 
$$\deg_R(g)=\deg_R(h)$$ this follows from (5) of Proposition \ref{g*}. 
\smallskip

\noindent $\bullet$ To define the degree of a map $f:X\to Y$ it is not strictly necessary
that $X$ is connected. In fact we can define 
$$\deg_R (f)= \sum_{X_c} \deg_R(f_{|X_c})$$
where $X_c$ varies among the connected components of $X$. By extending the above
homotopy invariance, we have: If $[X_0,f_0]=[X_1,f_1]\in \Bb_n(Y;R)$ then
$$\deg_R(f_0)=\deg_R(f_1) \ . $$
\smallskip

\noindent $\bullet$ {\it For every oriented $X$ as above, $n\geq 1$, for every $r\in \Z$ there is $g:X\to S^n$ 
such that $\deg_\Z(g)=r$.} 

First we prove it when $X=S^n$, by induction on $n\geq 1$. Consider $S^1$
as the unitary circle of $\C$. The restriction of $z\to \bar z$ to $S^1$ has $\Z$-degree equal to $-1$.
For every $r\geq 1$, the restriction of $z\to z^r$ has $\Z$-degree equal to $r$. As the degree is multiplicative
under composition  this achieves the result for $n=1$. For a given $r\in \Z$, let $g:S^n\to S^n$ be of degree
equal to $r$; we have to construct $\hat g: S^{n+1}\to S^{n+1}$ having the same degree.
Take $\hat g$ which fixes  the northern and southern poles and holds $\hat g (x)= tg(x/t)$ on $S^{n+1}\cap \{x_{n+2}=t\}$,
for every $t\in (-1,1)$. One checks that it has $\Z$-degree equal to $r$ as well.
To finish it is enough to construct $g:X \to S^n$ of $\Z$-degree equal to $\pm 1$. Fix  a smooth $D^n$ contained in a chart
of $X$. By using a tubular neighbourhood $U$ of $\partial D^n$ in $X$, it is not hard to construct a smooth map
$g:X\to S^n$ such that the restriction of $g$ to $D^n$ is a diffeomorphism onto $D^-=\{x\in S^n| \ x_{n+1}\leq 0\}$,
and holds constantly  the northern pole of $S^n$ on the complement of $D^n\cup U$ in $X$. Such a $g$ does the job. 

\begin{remark}{\rm For arbitrary oriented $X$ and $Y$ as above, it is in general a hard question to determine
the set of $r\in Z$ which can be realized as the $\Z$-degree of some $g:X\to Y$.}
\end{remark}
\smallskip

\noindent $\bullet$ Again in the case $X=Y=S^n$, $n\geq 1$. If $\rho:S^n\to S^n$ is the restriction of a reflection
of $\R^{n+1}$ along a linear hyperplane, then $\deg_\Z(\rho)=-1$. Denote by $a_n: S^n \to S^n, \ a_n(x)=-x$
the {\it antipodal map}; $a_n$ is the composition of the restriction of $n+1$ reflections (e.g. the reflections
along the hyperplanes $\{x_j=0\}$, $j=1,\dots, n+1$). Then we have
$$ \deg_\Z(a_n)= (-1)^{n+1} \ . $$
\smallskip

\noindent $\bullet$ In the setting of Remark \ref{3D-link}, let $S^n= \R^n\cup \infty$, $\infty \in S^n\setminus (T\cup Z)$.
$$L: T\times Z \to S^{n-1}, \ L(t,z)= \frac{t-z}{||t-z||}$$
then one can prove that
$$\deg_Z (L) = \pm lk(T,Z)$$
we left it as a (non trivial) exercise.

\subsection{A proof of the fundamental theorem of algebra}\label{fund-teo-alg}
The fundamental theorem of algebra states that every non constant complex polynomial $p(Z)\in \C[Z]$
has a complex root $a$, $p(a)=0$. There are several proofs; here is a topological/differential one based on the degree.

Let $p(Z)$ of degree $m\geq 1$. It is not restrictive to assume that $$p(Z)=Z^m + \sum_{j=1}^m a_jZ^{m-j} $$ is monic. 
Define the homotopy through polynomial maps:

$$ p_t(z)= tp(z) +(1-t)z^m= z^m+t( \sum_{j=1}^m a_jz^{m-j}), \ t\in [0,1] \ . $$
By the compactness of $[0,1]$, the ratios $p_t(z)/z^m$ tend uniformly to $1$ when $|z|\to +\infty$.
Hence there is $R$ bigh enough such that for every $t\in [0,1]$, the roots of $p_t(Z)$ are in the 
open ball $B_R=\{|z|<R\}$, with boundary $S_R \sim S^1$. Hence 
$$p_t/|p_t|: S_R\to S^1$$
is a well defined smooth map  for every $t$, so that $p_1/|p_1|(z)= p(z)/|p(z)|$ and $p_0/|p_0|(z)= z^m/R^m$ are homotopic
to each other. It is immediate that $$\deg_\Z (p_0/|p_0|)=m$$ hence also $\deg_\Z(p/|p|)=m$.
On the other hand, if $p(Z)$ has no roots, then $p/|p|$ can be extended to the whole closed ball $\bar B_R$,
it would be homotopically trivial, hence $\deg_\Z(p/|p|)=0$, a contradiction.

\cvd  

\section{The Euler class of a vector bundle}\label{char-VB}  Let $$\xi:= \pi: E\to X$$ be a vector bundle of {\it rank} $k$
(that is $k$ is the dimension of the fibre) over a compact boundaryless smooth $n$-manifold $X$. $X$ is considered
as a submanifold of $E$ via the canonical {\it zero section} $s_0:X\to E$. Then $$[X]\in \Bb^k(E;\Z/2\Z)$$ and set
$$w^k(\xi):=s_0^*([X])\in \Bb^k(X;\Z/2\Z) \ . $$ 
This is called the {\it Euler class of the vector bundle $\xi$}.
Let us spell how to get nice representatives of this last cobordism class.
\begin{lemma}\label{transv- section}  (1) The subset $\pitchfork\Gamma (\xi,X)$  made by the sections
$s: X\to E$ of $\xi$ such that $s\pitchfork X$ is open and dense  in $\Gamma(\xi)$. 

(2) Two sections transverse to $X$ are homotopic to each other through sections of $\xi$.
\end{lemma}
\Dim As $X$ is compact, the openess is now a routine fact. Let us show the density. Let $s:X\to E$ be any section. 
By transversality theorems, there is a map $z:X\to E$ close to $s$, $z\pitchfork X$, $z$ not necessarily
a section. If $z$ is close enough to $s$, then $h= \pi \circ z$ is a diffeomorphism onto $X\subset E$. Then 
$z\circ h^{-1}: X\to E$ is a section close to $s$ and transverse to $X$. Every section is homotopic to $s_0$
via a natural fibrewise radial homotopy.

\cvd

Let $s:X\to E$ be any section of $\xi$ transverse to $X$. Then its zero set 
$$Z_s=\{x\in X| \ s(x)=0\}$$ 
is a proper submanifold
of $X$ of dimension $n-k$. It follows from the very definition of $s_0^*$ that

\begin{lemma}\label{w-rep} For any section $s:X\to E$, $s\pitchfork X$, we have
$$w^k(\xi) = [Z_s]\in \Bb^k(X;\Z/2\Z) \ . $$
\end{lemma}

\cvd

\smallskip

\begin{proposition} For every couple $\xi$, $\rho$ of vector bundles on $X$ of rank $r$ and $s$
respectively, then
$$ w^{r+s}(\xi\oplus \rho)= w^r(\xi)\sqcup w^s(\rho) \ . $$
\end{proposition}
\Dim By using sections $s$ and $s'$ of $\xi$ and $\rho$  transverse to $X$ in $E(\xi)$ and
$E(\rho)$ respectively and such that
$s\oplus s'$ is transverse to $X$ in $E(\xi \oplus \rho)$, then
$$ Z_{s\oplus s'}=Z_s \pitchfork Z_{s'}$$ 
we conclude by means of Lemma \ref{w-rep}.

\cvd

\smallskip

It is evident that if there exists $s$ such that $Z_s=\emptyset$, then $w^k(\xi)=0$. Then:

{\it The non vanishing of the Euler class $w^k(\xi)$ is a basic
obstruction to the existence of a nowhere vanishing section of the vector bundle $\xi$.}

\smallskip

If $k>n=\dim X$, then  for every $s$ as above $Z_s=\emptyset$ and this fits with $\Bb^k(X;\Z/2\Z)=0$.
It follows that  $\xi$ of rank $k>n$ is strictly isomorphic to $\eta\oplus \epsilon^{n-k}$, $\eta$ being of rank $n$; in other words

{\it Every vector bundle over $X$ is stably equivalent to a vector bundle of rank $\leq \dim (X)$}.

\smallskip

\begin{proposition}\label{g*g*} Let $g:X\to Y$ be a smooth map between compact boundaryless smooth manifolds. 
Let $\xi$ be a rank $k$ vector bundle  over $Y$. Then
$$ w^k(g^*(\xi))=g^*(w^k(\xi))\in \Bb^k(X;\Z/2\Z) \ . $$
\end{proposition}
\Dim We stress that the first $g^*$ refers to the vector bundle pull-back while the second refers to the
cobordism pull-back. The two pull-back procedures are formally very similar and the equality is a
direct consequence.

\cvd

\smallskip  

{\bf Manifolds with boundary.} If the compact manifold $(X,\partial X)$ has non empty boundary then, for every rank $k$ vector bundle $\xi$ on $X$, 
the same procedure defines the relative Euler class
$$ w^k(\xi) \in \Bb^k(X,\partial X;\Z/2\Z) \ , $$
if $i:\partial X \to X$ is the inclusion then (as a particular case of the above proposition)
$$i^*(w^k(\xi))= w^k(i^*(\xi))\in \Bb^k(\partial X;\Z/2\Z) \ . $$
 \smallskip

{\bf Universal basic cobordism classes.} If $g:X\to \GG_{h,k}$ is any classifying map of $\xi$ so that $\xi$ is strictly equivalent to $g^*(\tau_{h,k})$
then $w^k(\xi)= g^*(w^k(\tau_{h,k}))$, $w^k(\tau_{h,k})\in \Bb^k(\GG_{h,k};\Z/2\Z)$. So these last  can be considered as
the {\it universal Euler classes of vector bundles}.
\smallskip

{\bf The total cobordism characteristic classes of projective spaces.} 
Consider the particular case of the {\it real projective space} $\PP^n(\R)=\GG_{n+1,1}$ with the tautological
line bundle $\tau_{n+1,1}$. Then 
$$\gamma^1:= w^1(\tau_{n+1,1})=[Z^1]\in \Bb^1(\PP^n;\Z/2\Z)$$
where $Z^1\sim \PP^{n-1}(\R)$ is any projective hyperplane in $\PP^n(\R)$.
For every $s\geq 1$,
$$ \gamma^s:=\sqcup_{j=1}^s \gamma^1 = [Z^s]$$
where $Z^s\sim \PP^{n-s}(\R)$ is any codimension $s$ projective subspace of $\PP^n(\R)$.
Set $\gamma^0:=[Z^0]=[\PP^n(\R)]$ the $\Z/2\Z$-fundamental class. Clearly if $s\leq n$,
$$\gamma^s\sqcup \gamma^{n-s}=1$$
hence they do not belong to $\ker (\phi^s)$ and $\ker (\phi^{n-s})$ respectively.
If $s>n$, $\gamma^s=0$. By definition
$$\sum_{s=0}^n \gamma^s \in \Bb^\bullet(\PP^n(\R);\Z/2\Z)$$
is the {\it total $\Z/2\Z$-cobordism characteristic class of $\PP^n(\R)$}.
If necessary we write $\gamma^s=\gamma^s_n$ in order to stress that it refers to
$\PP^n(\R)$. Then if we consider any linear inclusion $j:\PP^{k}(\R)\to \PP^n(\R)$, $k\leq n$,
$\PP^k(\R)= Z^{n-k}$ as above, then for every $m\geq 0$,
$$ \gamma_{k}^m= j^*(\gamma_n^m)\ . $$
\smallskip

\subsection{Oriented vector bundles}\label{oriented-bundles} A rank $r$ vector bundle $\xi$ over $X$ 
is oriented if it is defined by a maximal fibred atlas with  GL$^+(k,\R)$ cocycle. If the base manifold
is also oriented, then the total space manifold is naturally oriented itself. If $X$ is compact boundaryless, 
then we can repeat the above constructions in the oriented setting. This define the {\it oriented  Euler class}
$$e^r(\xi):= j^*([X])\in \Bb^r(X;\Z) \ . $$
$\omega^r(\xi)$ is the image of $e^r(\xi)$ via the natural forgetting map $\Bb(X;\Z)\to \Bb(X;\Z/2\Z)$.
For every pair of oriented bundles over $X$ of rank $r$ and $s$ respectively
$$ e^{r+s}(\xi \oplus \rho)= e^r(\xi)\sqcup e^s(\rho)\in \Bb^{r+s}(X;\Z) \ ; $$
for every $f:X\to Y$ smooth maps between oriented compact boundaryless manifolds, for every oriented rank $r$ vector bundle
$\xi$ bundle over $Y$,
$$ g^*(e^r(\xi))= e^r(g^*(\xi)) \in \Bb^r(X;\Z) \ . $$
Similarly we have relative oriented Euler classes $e^k(\xi)\in \Bb^k(X,\partial X;\Z)$ when $X$ has non empty boundary

\cvd 

A case of main interest  is the tangent bundle of $X$; then 
$$w^n(X):= w^n(T(X))\in \Bb^n(X;\Z/2\Z)=\Z/2\Z $$
{\it provides a basic obstruction to the existence of nowhere vanishing tangent vector fields on $X$}.

If we consider the rank $1$ determinant bundle of $X$, then 
$$w^1(X):= w^1(\det T(X))\in \Bb^1(X;\Z/2\Z)$$
{\it provides a basic obstruction in order that $X$ is orientable}. We will see in Corollary \ref{w1}
that it is a {\it complete} obstruction.
\smallskip
 
We will develop the case of (real and complex) rank $1$ bundles (also called {\it line bundles}) in Chapter
\ref{TD-LINE-BUND}. 
We will  develop the study of the Euler class of the tangent bundle of $X$  in Chapter \ref{TD-EP}.

\section{Borsuk-Ulam theorem}\label{Borsuk}
By definition a map $f:S^n \to S^m$ is {\it antipodal preserving} if for every $x\in S^n$,
$$ f(-x)=-f(x) \ . $$ 
\begin{proposition}\label{no-antipodal} For every $n\geq 1$, there does not exist any continuous 
antipodal preserving map $f:S^n\to S^{n-1}$.
\end{proposition}
\smallskip

The following corollary is known as the {\it Borsuk-Ulam theorem}.

\begin{corollary}\label{B-U} For every $n\geq 1$, for every continuous map $f: S^n \to \R^n$, there exists $x\in S^n$
such that $f(x)=f(-x)$.
\end{corollary}

\smallskip

For example, assuming that the surface of the  earth is a round sphere and that temperature and pressure vary
continuously on it in space and time, then at every instant there is a couple of antipodal points at which we have
the same couple of temperature and pressure values.
\medskip

{\it Proof of BUT.} By contradiction, if a given $f$ does not verifies the consclusion of the Corollary,
then 
$$g: S^n\to \R^n, \ g(x)= f(x)-f(-x)$$
is continuous, nowhere vanishing, and for every $x\in S^n$,
$$ g(-x)=f(-x)-f(x)=-g(x) \ . $$
Then
$$\hat g: S^n \to S^{n-1}, \ \hat g(x)= g(x)/||g(x||$$
is continuous and would be antipodal preserving, against Proposition \ref{no-antipodal}

\cvd

\smallskip

{\it Proof of Proposition \ref{no-antipodal}.} 
To lighten the notations, in this proof we will use $\eta_k(*)$ instead $\Bb_k(*; \Z/2\Z)$,
and write $\PP^m$ instead of $\PP^m(\R)$.
\smallskip

The case $n=1$ is evident because $S^1$ is connected
while $S^0=\{\pm1 \}$ is not. 
\smallskip

For $n=2$ we use some basic facts about the fundamental group of
a manifold and its action on a universal covering space. Assume that there is such a continuous antipodal preserving 
map $f:S^2\to S^1$. It induces a map $\hat f: \PP^2\to \PP^1\sim S^1$ such that the following
diagram commutes, the vertical maps being the natural degree $2$ covering maps:
$$ \begin{array}[c]{ccc}
S^2 &\stackrel{ f }{\rightarrow}&  S^1\\
\downarrow\scriptstyle{p_2}&&\downarrow\scriptstyle { p_1}\\
\PP^2&\stackrel{\hat f}{\rightarrow}& \PP^1 \end{array} \ . $$

We know that $\pi_1(\PP^2,x_0)\sim \Z/2Z$, generated by 
the class of a projective line passing through the base point,
while $\pi_1(\PP^1,\hat f(x_0)) \sim \Z$, generated by the the identity loop. Hence the induced
homomorphism 
$\hat f_*:\pi_1(\PP^2,x_0)\to \pi_1(\PP^1, \hat f(x_0))$ is necessarily trivial.
On the other hand, take the two antipodal points $x, -x \in S^2$ over $x_0$ and an arc $\sigma$ in $S^2$
that joins them.  Then $p_2(\sigma)$ represents a non trivial element of $\pi_1(\PP^2,x_0)$,
because it acts non trivially on $S^2$ which is the universal covering of the projective plane. 
The class $\hat f_*(<p_2(\sigma)>)$ is represented by $p_1\circ f\circ \sigma$
and again it is non trivial because it acts non trivially on the universal covering space of $\PP^1$ that dominates the
covering  $p_1$. This is agaist the fact that $\hat f_*=0$.
\smallskip

If $n>2$ we have a similar commutative diagram
$$ \begin{array}[c]{ccc}
S^n &\stackrel{ f }{\rightarrow}&  S^{n-1}\\
\downarrow\scriptstyle{p_n}&&\downarrow\scriptstyle { p_{n-1}}\\
\PP^n&\stackrel{\hat f}{\rightarrow}& \PP^{n-1} \end{array} $$
where both vertical maps are now universal covering maps.
Both fundamental groups are isomorphic to $\Z/2\Z$ and the very
same argument used above shows that 
$$\hat f_*: \pi_1(\PP^n,x_0)\to \pi_1(\PP^{n-1},\hat f(x_0))$$
is an isomorphism. 
Any surjective  homomorphism $g:\Z/2\Z \to G$ either is
an isomorphism or $G=0$ and $g$ is trivial. For every $m>1$, the surjective homorphism
$$ \hat h:=\sigma_1 \circ h_1: \pi_1(\PP^m ,x_0)\to \eta_1(\PP^m)$$
is non trivial (the class of a projective line $Z^{m-1}$ passing through the base point 
is sent by $\hat h$ to the non trivial class $[Z^{m-1}]\in \eta_1(\PP^m)$, for
via the tautological isomorphism $[Z^{m-1}]=\gamma^{m-1}_m \in \eta^{m-1}(\PP^m)$,
and we know that $\gamma^{m-1}_m \sqcup \gamma^1=1$). Hence $\hat h$ is an isomorphism
anf $\hat f$ induces an isomorphism (we keep the notation)
$$ \hat f_*: \eta_1(\PP^n)\to \eta_1(\PP^{n-1}) \ . $$
For every $m>1$, Hom$(\eta_1(\PP^m),\Z/2\Z)\sim \Z/2\Z$. Then in our situation
$$f^t_*: {\rm Hom}(\eta_1(\PP^{n-1}),\Z/2\Z)\to {\rm Hom}(\eta_1(\PP^n),\Z/2\Z)$$ 
is also an isomorphism. For every $m>1$,
$$\hat \phi: \eta^1(\PP^m)/\ker(\phi)\to {\rm Hom} (\eta_1(\PP^m),\Z/2\Z)$$
is an isomorphism and $\eta^1(\PP^m)/\ker(\phi)$ is generated by $\gamma_m^1$.
Then on one hand we would have 
$$ \hat f^*(\gamma^1_{m-1}) = \gamma_{m}^1, \ \hat f_* (\hat \phi (\gamma^1_{m-1}))=\hat \phi(\gamma^1_m)  $$
on another hand 
$$0= \hat f^*(0)=\hat f^*(\sqcup_{s=1}^{m} \gamma^1_{m-1})= \sqcup_{s=1}^m \gamma^1_m=1$$ 
and this is a contradiction.

\cvd

\chapter{Line bundles, hypersurfaces and cobordism}\label{TD-LINE-BUND}
In this chapter $X$ will denote a compact boundaryless smooth manifold and we also assume that
$X$ is connected (in general we can apply the next arguments to every
connected component). We will use indifferently the notations $\eta_j(X)$ or 
$\Bb_j(X;\Z/2\Z)$ (resp. $\Omega_j(X)$ or $\Bb_j(X;\Z)$)
and so on. Recall also 
$$\Hh^r(X;R):= \Bb^r(X;R)/\ker (\phi^r)$$ 
defined
in Section \ref {intersection-form}. 
By means of the Euler classes of line bundles over $X$ one can
achieve a good understanding of $\eta^1(X)=\Bb^1(X;\Z/2\Z)$. If $X$ is oriented, we will get
information about $\Omega^1(X)$ and by using {\it complex} line bundles also about $\Omega^2(X)$. 

\section{Real line bundles and hypersurfaces}\label{real-LB}
 Let $X$ be as above. Denote by 
 $$\Vv_1(X)$$ 
 the set of rank $1$ real vector bundles on $X$ (also called (real) {\it line bundles}) 
 considered up to strict equivalence. We know from Chapter \ref{TD-EMB-VB}
that 
$$\Vv_1(X) \sim [X, \PP^\infty(\R)]$$
where this last is  the space of homotopy classes of {\it classifying maps} $f\in \Ee(X,\PP^\infty(\R))$,
and the bijective correspondence is given via the pull back of the tautological line bundle:
$$ [X,\PP^\infty(\R)] \to \Vv_1(X), \ [f] \to [f^*(\tau_{\infty,1})] \ . $$ 
Moreover, by Section \ref{truncated-VB} we know that we can ``truncate''  the classifying maps so that eventually
$$ \Vv_1(X)\sim  [X,\PP^{m(n)}(\R)]$$
where $m=m(n)$ is big enough only depending on $n$.
Often we will confuse a class with a given representative 
(say we write $f$ instead of $[f]$, $\xi$ instead of $[\xi]$, and so on).
Recall that the {\it tensor product} defines an operation
$$ \otimes: \Vv_1(X)\times \Vv_1(X) \to \Vv_1(X), \ (\xi,\beta)\to \xi \otimes \beta \ . $$ 
In Section \ref {char-VB}, we have defined a map
$$w^1: \Vv_1(X)\to \eta^1(X), \ \xi \to w^1(\xi)$$
which associates to every line bundle its Euler class.
Precisely $w^1(\xi)$ can be represetended as
$$w^1(\xi)=[Z]$$
where $Z$ is a smooth compact {\it hypersurface} in $X$ given as the zero set
$Z=Z_s$ of any section $s\in \Gamma(\xi)$ transverse to $X$ in $E(\xi)$,
where $X$ is canonically embedded in the total space of $\xi$ by the zero section $s_0$.
Moreover, if $Z_0$ and $Z_1$ are two such zero sets, then we can realize the equality
of their bordism classes 
$[Z_0]=[Z_1]\in \eta^1(X)$ by means of {\it embedded bordisms}: 
\smallskip

{\it There exists a proper hypersurface $(Y,\partial Y)$ of $(X\times [0,1], (X\times \{0\}) \amalg (X\times \{1\}))$
such that $\partial Y = Z_0 \amalg Z_1$, $Z_i \subset X\times \{i\}$. The map which interpolates the 
two inclusions $j_i: Z_i \to X$
is the projection onto $X$.}
\smallskip

So we denote by 
$$\eta^1_{{\rm Emb}}(X)$$ 
the {\it set} of proper smooth hypersurfaces of $X$ considered up to embedded bordism. 
There is a natural projection 
$$\pG: \eta^1_{{\rm Emb}}(X)\to \eta^1(X)$$
so that  the above map $ w^1$ factorizes as
$$ w^1= \pG \circ \hat w^1$$
through  a well defined map
$$\hat w^1: \Vv_1(X)\to \eta^1_{{\rm Emb}}(X) \ . $$

We have
\begin{proposition}\label{LB=Hyper} 
 (1) The map $\hat w^1:\Vv_1(X)\to \eta^1_{{\rm Emb}}(X)$ is bijective.
 \smallskip
 
 (2)  For every couple $(\xi,\beta)\in \Vv^2_1$, 
 $$w^1(\xi \otimes \beta)= w^1(\xi)+w^1(\beta) \ . $$

(3) The projection  $\pG$ maps $\eta^1_{{\rm Emb}}(X)$ onto a $\Z/2\Z$-submodule, say $\HH^1(X;R)$,
of $\Bb^1(X;\ Z/2\Z)$,  
the one made by the (unoriented) cobordism classes that can be represented
by embedded hypersurfaces).
\end{proposition}

\Dim Let us describe the inverse map of $\hat w^1$. For every proper hypersurface
$Z$ of $X$ we have to construct a line bundle $\xi_Z$ on $X$ such that $Z=Z_s$
for some $s\in \Gamma(\xi_Z)$, $s\pitchfork X$. We can find a finite nice atlas
of $(X,Z)$, $\{(W_j,\phi_j)\}$ such that for every $j$, there is a summersion 
$f_j: W_j \to \R$, such that $W_j\cap Z = \{ f_j=0\}$. On $W_i\cap W_j$,
by Remark \ref{division2} (2) the ratio $f_i/f_j$ defined a priori outside
the zero set of $f_j$, extends to a well defined, smooth and 
nowhere vanishing function 
$$g_{i,j}: W_i\cap W_j \to \R, \ g_{i,j}(x)= f_i(x)/f_j(x) \ . $$
Hence
$$ \{ g_{i,j} :W_i\cap W_j \to \R^*\}$$ 
actually defines a cocycle of a line bundle $\xi_Z$ on $X$ which
has the desired properties by construction.

As for (2), we can assume that $\xi$ and $\beta$ are defined by means
of cocycles $\{\mu_{i,j}\}$ and $\{\nu_{i,j}\}$ respectively over a same nice
atlas of $X$. Then $\{\mu_{i,j}\nu_{i,j}\}$ is a cocycle for $\xi\otimes \beta$.
Then if $\{s_i\}$ and $\{s'_{i}\}$ are representations in local coordinates
of sections $s$ and $s'$ of $\xi$ and $\beta$ respectively ,
such that $s\pitchfork X$, $s'\pitchfork X$, $s\pitchfork s'$, then 
$\{s_is'_i\}$ detemines a section say $ss'$ of $\xi\otimes \beta$
such that $[Z_s]= w^1(\xi)$, $[Z_{s'}]= w^1(\beta)$; by perturbing
$ss'$ to get $s"\pitchfork X$, eventually $Z_{s"}$ represents $w^1(\xi\otimes \beta)$
and $[Z_{s"}]= [(Z_s,i) \amalg (Z_{s'},i')]$. In fact $Z_{s"}$ can be considered as
an embedded desingularization in $X$ of $Z_s\cup Z_{s'}$, which
is singular along the codimension $2$ submanifold $Y=Z_s\pitchfork Z_{s'}$.
\smallskip
  
Item (3) is a consequence of (1) and (2).

\cvd 

\section{Real line bundles and Rep$(\pi_1, \Z/2\Z)$}\label{pi-1-rep}
Recall that we are assuming that $X$ is connected.
We denote by Rep$(\pi_1(X), \Z/2\Z)$ the set of group homomorphisms
(the base point of $X$ is understood). Recall the linear map
$$ \phi: \eta^1(X)\to {\rm Hom}(\eta_1,\Z/2\Z), \ \phi_\gamma(\sigma)= \gamma \sqcap \sigma \ . $$
Recall the surjective homomorphism
$$\hat h: \pi_1(X)\to \eta_1(X) \ . $$
Then we define the map
$$ \kappa: \Vv_1(X)\to {\rm Rep}(\pi_1(X),\Z/2\Z), \  \kappa(\xi)= \phi_{w^1(\xi)} \circ \hat h \ . $$
Here is a concrete way to describe $\kappa(\xi)$.
As $\pi_1(\PP^\infty(\R))=\Z/2\Z$, then $\Vv_1(S^1)$ consists of two line bundles: the trivial 
and the non trivial one which has the total space diffeomorphic to an open M\"obius band.
If $\sigma=<f:S^1\to X>\in \pi_1(X)$, then $\kappa(\xi)(\sigma)=1$ if and only if $f^*\xi$
is non trivial.
We have
\begin{proposition}\label{v1-rep} The map $\kappa: \Vv_1(X)\to {\rm Hom}(\pi_1(X),\Z/2\Z)$
is bijective.
\end{proposition}
\Dim We have already remarked in Example \ref{P-infinito} that 
$\PP^\infty(\R)$ is a $K(\Z/2\Z,1)$ space. It is a fundamental property
of such a space that for every 
$$\sigma \in {\rm Rep}(\pi_1(X),\Z/2\Z)$$ 
there is a unique  
$$f\in [X, \PP^\infty(\R)]$$ such that 
$$\sigma= f_*:\pi_1(X)\to \pi_1(\PP^\infty(\R)) \ . $$
Then 
$$\sigma \to  \xi_\sigma := f^*(\tau_{\infty,1})$$ defines the inverse map of $\kappa$.
Equivalently, we can describe $\kappa^{-1}$ in terms of degree $2$ covering maps.
It is known that there is a bijection between the degree $2$ covering maps over $X$
(up to strict equivalence) and Rep$(\pi_1(X),\Z/2\Z)$. For every line bundle $\xi$,
$\kappa(\xi)$ corresponds to the double covering of $X$ given by the unitary bundle with fibre
$S^0$ associated to $\xi$. Viceversa every degree $2$ covering of $X$ can be considered as a
fibre bundle defined by a cocycle over a finite open covering of $X$ with values in the multiplicative
subgroup $\{\pm 1\}$ of $\R^*$. So it can be considered as the unitary bundle associated to the
line bundle determined by the same cocycle.

\cvd

\smallskip

Referring to Proposition \ref{LB=Hyper}, we have the following immediate corollaries.
 
\begin{corollary}\label{eta1} (1) The map 
$\pG: \eta^1_{{\rm Emb}}(X)\to \HH^1(X;\Z/2\Z)\subset \Bb^1(X;\Z/2\Z)$
is bijective.
\smallskip

(2) $\HH^1(X;\Z/2\Z) \sim \Hh^1(X;\Z/2\Z) \sim {\rm Hom}(\eta_1(X), \Z/2\Z)$.
\smallskip

(3) $\Vv_1(X)\sim \Hh^1(X;\Z/2\Z)$.
\end{corollary}

\cvd

\smallskip

Another consequence of the above discussion is that 
\medskip

{\it$\Hh^1(X;\Z/2\Z)$ is finite dimensional}. 
\smallskip

\noindent For
as $X$ is compact, then $\pi_1(X)$ is finitely generated, hence $\eta_1(X)=\hat h (\pi_1(X))$ is a finite
dimensional $\Z/2\Z$-vector space as well as $\Hh^1(X;\Z/2\Z)$.

\medskip

By applying the above results to the determinant line bundle of $X$ we have

\begin{corollary}\label{w1} A compact connected boundaryless smooth manifold $X$
is orientable if and only if $w^1(X)=0\in \Hh^1(X;\Z/2\Z)$.
\end{corollary}

\cvd 

\section{Oriented hypersurfaces and $\Omega^1$}\label{omega1}
Assume that $X$ is oriented. Then we have the $\Z$-linear map
$$\phi: \Omega^1(X)\to {\rm Hom}(\Omega_1(X),\Z)$$
and via the homomorphism
$$h: \pi_1(X)\to \Omega_1(X)$$
we define a map
$$ \kappa: \Omega^1(X)/\ker (\phi) \to {\rm Rep}(\pi_1(X);\Z) \ . $$
As usual 
$$[X,S^1]$$
is the set of homotopy classes in $\Ee(X,S^1)$. Denote by $\beta_{S^1}$
the usual generator of $\Omega^1(S^1)$ which fixes the identification $\Omega^1(S^1)=\Z$.
We have the $\Z$-linear map
$$\wG: [X,S^1]\to \Omega^1(X), \ f \to f^*(\beta_{S^1}) \ . $$
In fact $$f^*(\beta_{S^1})=[Z]$$
where $Z$ is an oriented proper hypersurface of $X$ of the form
$$Z= f^{-1}(s_0)$$
$s_0$ being any regular value of $f$. 
We denote by
$$ \Omega^1_{{\rm Emb}}(X)$$
the set of of proper oriented hypersurfaces of $X$ considered up to {\it oriented
embedded bordism} (this notion is the natural enhancement of the unoriented one
given above). Then we have the projection
$$\pG:   \Omega^1_{{\rm Emb}}(X) \to \Omega^1(X)$$
such that $\wG$ factorizes as $\hat \wG \circ \pG$
for a well defined map
$$ \hat \wG: [X,S^1]\to  \Omega^1_{{\rm Emb}}(X) \ . $$
Finally we have the map
$$\rG: [X,S^1]\to {\rm Rep}(\pi_1(X);\Z), \  f\to f_*: \pi_1(X)\to \pi_1(S^1)=\Z \ . $$ 

We have
\begin{proposition}\label {hyper-omega1}
(1) The map $ \hat \wG: [X,S^1]\to  \Omega^1_{{\rm Emb}}(X)  $
is bijective.
\smallskip

(2) The map  $\rG:  [X,S^1]\to {\rm Rep}(\pi_1(X);\Z)$ is bijective.
\smallskip

(3) The map $ \kappa: \Hh^1(X;\Z) \to {\rm Rep}(\pi_1(X);\Z) $
is bijective.
\smallskip

(4) The projection $\pG:  \Omega^1_{{\rm Emb}}(X) \to \Omega^1(X)$
is injective onto a $\Z$-submodule say
$$\HH^1(X;\Z) \subset \Bb^1(X;\Z) \ . $$

(5) $\HH^1(X;\Z) \sim  \Hh^1(X;\Z) \sim {\rm Hom}(\Omega_1,\Z)$.
\smallskip

(6) $\Hh^1(X;\Z)$ is finitely generated.

\end{proposition}
\Dim Let us define the inverse map of $\hat \wG$. This is a first sample
of a general construction that we will study with all details in Chapter \ref{TD-PT}.
So we limit here to indicate the main points.
Let $Z$ be a proper oriented hypersurface of $X$. As both $X$ and $Z$ are
oriented, we can fix a global trivialization $t: Z\times (-1,1) \to U$ of a tubular neighbourhood
of $Z$ in $X$. Let $s_-$ the southern pole of $S^1$, $s_+:=\infty$ the northern one.
Let $D\sim (-1,1)$ be an open interval  in $S^1$ centred at $s_-$. Then the composition
of $t^{-1}$ with the projection onto $(-1,1)$ define a local summersion $f: U \to D\subset S^1$.
By using  a suitable partition of unity as usual, we can globally define $f_Z: X \to S^1$
such that $f_Z$ is constantly equal to $\infty$ on the complement of $U$,
equals $f$ on $t((-1/2,1/2))$ and $f^{-1}(s_-) = Z$.  One verifies that the homotopy
class of such a map $f_Z$ is invariant up to oriented embedded bordism of hypersurfaces,
so $[Z]\to [f_Z]$ eventually defines the inverse map of $\hat \wG$. This achieves (1).

As for (2), it is well known that $S^1$ is a $K(\Z,1)$ space. Hence for every
$\sigma \in {\rm Rep}(\pi_1(X);\Z)$, there is $f^\sigma: X \to S^1$, uniquely defined up to homotopy,
such that $\sigma = f^\sigma_* :\pi_1(X)\to \pi_1(S^1)=\Z$. This defines the
inverse map of $\kappa$. 

The item (3) follows from (1) and (2) by readly noticing that if $[Z]=\hat \wG([f])$, then 
$f_*=\phi([Z])$. Items (4) and (5) are basically a rephrasing of the previous ones;
(6) follows again from the fact that $X$ is compact, hence $\pi_1(X)$ is finitely
generated and the homomorphism $h$ is surjective. 

\cvd

\section{Complex line bundles and $\Omega^2$}\label{complexLB}
Assume again that $X$ is oriented. Denote by 
$$\Vv_1(X,\C)$$
the set of {\it complex} line bundles over $X$ considered up to strict equivalence.
Similarly to the real case,
$$\Vv_1(X,\C) \sim [X, \PP^\infty(\C)]$$
where this last is  the space of homotopy classes of {\it classifying maps} $f\in \Ee(X,\PP^\infty(\C))$,
and the bijective correspondence is given via the pull back of the tautological complex line bundle:
$$ [X,\PP^\infty(\C)] \to \Vv_1(X,\C), \ [f] \to [f^*(\tau^\C_{\infty,1})] \ . $$ 
Moreover, we can ``truncate''  the classifying maps so that eventually
$$ \Vv_1(X,\C)\sim  [X,\PP^{m(n)}(\C)]$$
where $m=m(n)$ is big enough only depending on $n=\dim (X)$.
Every complex line bundle $\xi$ underlies a rank $2$ {\it oriented} real bundle
$\xi_\R$. Viceversa, every rank $2$ oriented real bundle can be endowed with
a structure a complex line bundle by reducing the structural group to $SO(1)$
and by identifying the rotation by $\pi/2$ to the product by $\sqrt{-1}$. 
Then we can define
$$e^2: \Vv_1(X;\C)\to \Omega^2(X), \ \xi \to e^2(\xi_\R)$$
which associates to every $\xi$ the oriented Euler class of its ``realification''. 
Precisely $e^2(\xi)$ can be represetended as
$$e^2(\xi)=[Z]$$
where $Z$ is a proper codimension $2$ oriented smooth submanifold of $X$
given as the oriented zero set
$Z=Z_s$ of any section $s\in \Gamma(\xi_\R)$ transverse to $X$ in $E(\xi_\R)$.
If $Z_0$ and $Z_1$ are two such zero sets, then we can realize the equality
of their bordism classes 
$[Z_0]=[Z_1]\in \Omega^2(X)$ by means of {\it oriented embedded bordisms}
via  proper oriented codimension $2$ submanifold $(Y,\partial Y)$ of 
$(X\times [0,1], (X\times \{0\}) \amalg (X\times \{1\}))$.
Similarly as above denote by $\Omega^2_{{\rm Emb}}(X)$
the set of codimension $2$ oriented proper submanifolds of $X$
considered up to embedded oriented bordism, and
$$\pG: \Omega^2_{{\rm Emb}}(X)\to \Omega^2(X)$$
the natural projection. The map $e^2$ factorizes as
$\pG \circ \hat e^2$ where
$$ \hat e^2: \Vv_1(\C)\to \Omega^2_{{\rm Emb}}(X) \ . $$
Recall the $\Z$-linear map
$$\phi^2: \Omega^2 \to {\rm Hom}(\Omega_2(X),\Z)$$
which composed with $e^2$ and the homomorphism
$$h: \pi_2(X)\to \Omega_2(X)$$
leads to the map
$$ \kappa: \Vv_1(\C) \to  {\rm Rep}(\pi_2(X),\Z) \ . $$
Finally, analogously  to the real case, $\PP^\infty(\C)$ is a $K(\Z,2)$-space (\cite {Hatch}), hence
it is defined and is  {\it bijective} the map
$$\rG: [X,\PP^\infty(\C)]\to {\rm Rep}(\pi_2(X),\Z), \ f\to f_*:\pi_2(X)\to \pi_2(\PP^\infty(\C))=\Z \ . $$
By combining these facts similarly to the real case we have
\begin{proposition}\label{v1C-omega2}
(1) The map $\hat e^2: \Vv_1(\C)\to  \Omega^2_{{\rm Emb}}(X)$ is bijective.
\smallskip

(2) For every $(\xi, \beta)\in \Vv_1^2(\C)$, $e^2(\xi\otimes_\C \beta)= e^2(\xi)+ e^2(\beta)$.
\smallskip

(3) The map $\kappa: \Vv_1(\C) \to {\rm Rep}(\pi_2(X),\Z) $ is bijective.
\smallskip

(4) The projection $\pG$ is injective and maps $\Omega^2_{{\rm Emb}}(X)$ onto a $\Z$-submodule, say
$\HH^2(X;\Z)$ of $\Bb^2(X;\Z)$.

(5) $\HH^2(X;\Z)\sim \Hh^2(X;\Z) \sim {\rm Hom}(\Omega_2(X)/\phi^2, \Z)$. 

\end{proposition}

\cvd

\subsection{Relative case}\label {relative} If $(X,\partial X)$ is compact with non empty boundary,
possibly oriented, this is part of of the setting of Section  \ref{relative-cb}. So one can elaborate
a relative version of the previous results. We limit to state the existence of isomorphisms

$$ \Hh^1(X, \partial X;\Z/2\Z)\to {\rm Hom}(\Hh_1(X;\Z/2\Z),\Z/2\Z)$$

$$ \Hh^1(X,\partial X;\Z)\to {\rm Hom}(\Hh_1(X;\Z),\Z)$$

$$ \Hh^2(X,\partial X;\Z) \to {\rm Hom}(\Hh_2(X;\Z),\Z) \ . $$

\section{Seifert's surfaces}\label{seifert} 
Let $X$ be a compact oriented boundaryless manifold. By applying similar arguments about complex line bundles or
rank $2$ oriented real bundles, we want to prove the 
following proposition.

\begin{proposition}\label{seifert-surf} Let $Y\subset X$ be a proper oriented
codimension $2$ submanifold of $X$. Assume that $[Y]\in \ker(\phi)$, that is
$[Y]=0\in \Hh^2(X;\Z)$. Let $\pi: U\to Y$ be a tubular neighbourhood of $Y$ in $X$.
Let $W= X\setminus {\rm Int}(U)$ with boundary $\partial W = \partial U$.
Then there exists a compact  oriented hypersurface with boundary $\tilde Z$
of $X$ such that $\partial \tilde Z = Y$. Such a $\tilde Z$ is called a {\rm Seifert
surface} of $Y$. Precisely, $\tilde Z$ is transverse to $\partial W$, 
$(Z,\partial Z):= (\tilde Z \cap W, \tilde Z \cap \partial W)$ is a proper oriented
hypersurface in $(W,\partial W)$, $U\cap \tilde Z$ is a collar of $Y$ in $\tilde Z$.
\end{proposition}
\Dim Let $i: Y\to X$ be the inclusion. Any tubular neighbourhood $p:U\to Y$ of $Y$ in $X$
can be associated to a direct sum decomposition of the form
$$ i^*(T(X))=T(Y)\oplus \xi_\R$$
where $\xi_\R$ is the ``realification'' of a complex line bundle on $Y$.
As $[Y]\in \ker(\phi)$, then $e^2(\xi)=0$, hence $\xi$ is trivial so that 
$U$ admits  global
trivializations  which induce trivializations of $\partial W$. Let us fix one 
$h_0: \partial W \to Y\times S^1$. Fix one oriented fibre $D\sim D^2$ of $\pi$
with oriented boundary $S\sim S^1$. We claim that $[S]$ is of infinite order in $\Omega_1(W)$. By contradiction,
let us assume that say $p\neq 0$ parallel copies of $S$ are the boundary of a singular manifold $g: (V,\partial V)\to (W,\partial W)$.
Then by gluing $V$ and $p$ parallel copies of $D$ along the boundary, we would get an ``absolute''  singular $2$-manifold $(\tilde V,\tilde g)$
in $X$ such that $[Y]\sqcup [\tilde V,\tilde g]=p$, against the fact that $[Y]\in \ker(\phi)$.  As $[S]$ is indivisible in $\Omega_1(W)$,
there exists $\psi \in {\rm Hom}(\Omega_1(W),\Z)$ such that $\psi([S])=1$. We know that $\psi$ is realized by a map $f_\psi: (W,\partial W)\to S^1$
transverse to a given point $q\in S^1$. Denote by $j: \partial W \to W$ and $r: S\to \partial W$ the two inclusions. Then $\gamma:=j^t(\psi)$
is realized by the restriction $f_\gamma$ of $f_\psi$ to $\partial W$, while the restriction of $f_\gamma$ to $S$ realizes $(j\circ r)^t(\phi)$
and is homotopic to the identity. Up to modify the given trivialization $h_0$ by a suitable one say $h$, $f_\gamma$ factorizes as $p\circ h$,
where $h: \partial W \to Y\times S^1$ and $p: Y\times S^1\to S^1$ is the projection onto the second factor. Then $(Z,\partial Z)=(f_\phi^{-1}(q), f_\gamma^{-1}(q))$
and $\tilde Z$ obtained by gluing along $\partial Z$ the mapping cylinder of the restriction of $\pi$ to it
achieve the proof.

\cvd

From the last step of the above proof we have the following corollary.

\begin{corollary}\label{seifert-cor} Let $X$ be an oriented compact $n$-manifold with boundary $\partial X$.
Let $Z$ be e proper oriented submanifold of dimension $n-2$ of $\partial X$. Assume that
$[Z]=0$ in $\Hh^2(X;\Z)$. Then there is a proper oriented hypersurface $(W,\partial W)$ such that
$Z=\partial W$.
\end{corollary}

We have also the following version of  Corollary \ref{seifert-cor} 
when $Z$ is of codimension $2$ in $\partial X$.
  
  \begin{proposition}\label{seifert-surf2}  Let $X$ be an oriented compact $n$-manifold with boundary $\partial X$.
Let $Z$ be e proper submanifold of dimension $n-3$ of $\partial X$. Assume that
$[Z]=0$ in $\Hh^3(X;\Z)$. Then there is a proper codimension-2 oriented submanifold $(W,\partial W)$ 
of $(X,\partial X)$ such that $Z=\partial W$.
\end{proposition} 
\Dim  The hypotheses put us in a situation analogous to the last step in the proof of Proposition \ref{seifert-surf}, that is to
Corollary \ref{seifert-cor}. Here $S^1$ is replaced by $\PP^n(\C)$ ($n$ big enough) in the sense that both carry special instances of the 
{\it Pontryagin-Thom's construction} which will be considered in Chapther \ref{TD-PT} in full generality. Let $f_0: Z\to \PP^{n-1}(\C)$ be a classifying
map of the oriented normal rank-$2$ bundle of $Z$ in $\partial X$. Note that $\PP^{n}(\C)\setminus \{x_0\}$, $x_0\in \PP^n(\C) \setminus \PP^{n-1}(\C)$,
is diffeomorphic to the total space of the tautological vector bundle on $\PP^{n-1}(\C)$. Hence $f_0$ extends to a map
$f: \partial X \to \PP^{n}(\C)$ such that $f \pitchfork \PP^{n-1}(\C)$ and $Z= f^{-1}(\PP^{n-1}(\C)$. As $[Z]=0$ in $\Hh^3(X;\Z)$, if
$n$ is big enough then $f$ can be extended to a map $F: X\to \PP^n(\C)$ which we can assume transverse to $\PP^{n-1}(\C)$.
Finally $W=F^{-1}(\PP^{n-1}(\C)$ does the job.

\cvd

As a corollary we have a weak version of Proposition \ref{seifert-surf} when $Y$ has codimension $3$.

\begin{corollary}\label{seifert-cor3}  Let $Y\subset X$ be a proper oriented
codimension $3$ submanifold of $X$. Assume that the normal bundle of $Y$ in
$X$ has a non vanishing section $s$ and let $Y'=s(Y)$ be a copy of $Y$ in the boundary
$\partial U$ of a tubular neighbourhood of $Y$ in $X$. Assume that  $[Y']=0$
in $\Hh^3(X\setminus {\rm Int}(U); \Z)$. Then there is a proper oriented codimension-$2$ 
submanifold $(W,\partial W)$ of $(X\setminus {\rm Int}(U),\partial U)$ such that $\partial W = Y'$.
\end{corollary}

\begin{remark}\label{unor-seifert} (Non orientable Seifert surfaces)
  {\rm In the statement of Proposition
    \ref{seifert-surf} do not assume that $X$ and $Y$ are orientable
    and use $\Hh^2(X;\Z/2\Z)$ instead. It is natural to inquire about
    the existence of possibly non orientable Seifert surface. We see
    an immediate obstruction: if a Seifert surface exists and $
    i^*(T(X))=T(Y)\oplus \xi$ is as above (where $\xi$ is now not
    necesseraly trivial nor orientable), then $\xi$ has a nowhere
    vanishing section.  The above proof can be adapted to show that
    this is really the only obstruction.}
  \end{remark}

\chapter{Euler-Poincar\'e characteristic}\label{TD-EP}
$X$ will denote a compact connected oriented boundaryless smooth $n$-manifold. Then also the tangent
bundle $\pi: T(X)\to X$ is tautologically an {\it oriented} rank $n$ vector bundle on $X$: the orientation of $X$
determines in a coherent way an orientation on every fibre $T_pX$ of $T(X)$. Then we can consider the oriented
Euler class
$$ e^n(X)\in \Omega^n(X)=\Bb^n(X;\Z)=\Z \ . $$
By a traditional change of notation 
$$\chi(X):= e^n(X)\in \Z$$
is called the {\it Euler-Poincar\'e characteristic of $X$}. If $X$ is not connected, $\chi(X)$ is defined as
the sum of the characteristics of its connected components.

Recall that $\chi(X)$ is computed by means of any section $s$ of $T(X)$ transverse to $X$. In other words,
$\chi(X)$ is the {\it self-intersection number}  of $X$ in $T(X)$. Such a section $s\pitchfork X$ is 
a tangent vector fields on $X$ with only {\it non-degenerate zeros}: $s$ can be expressed in local coordinates at every
such a zero $p\sim 0$ in the form $$s(x)=(x, f_p(x))$$ where $f_p: (\R^n,0)\to (\R^n,0)$ is a diffeomorphism.
The sign $\epsilon (p)=\pm 1$ of the zero $p$, so that $$\chi(X)= \sum_{p; s(p)=0} \epsilon (p)$$ is readily computed as
$$\epsilon (p) = {\rm sign}(\det d_0f_p) \ . $$ 

\section{E-P characteristic via Morse functions}\label{E-P-morse}
Let $f:X\to \R$ be a Morse function with critical points $p_1,\dots,p_r$ of index $q_1,\dots , q_r$.
Let $\nabla_g f$ be an adapted gradient field of $f$ as in Section \ref {triad-dissection}.
Then $p_1,\dots, p_r$ are also the zeros of this field. It is easy to check by using the Morse local
coordinates that they are non degenerate zeros and their sign is given by
$$\epsilon (p_j) = (-1)^{q_j} \ . $$
Hence we have
$$\chi(X) = \sum_{j=1}^r (-1)^{q_j} \ . $$
This has the following interesting corollary.
\begin{corollary}\label{odd} If $\dim (X)=n$ is odd, then $\chi(X)=0$.
\end{corollary}
\Dim Consider the Morse function $1-f$; $f$ and $1-f$ have the same critical points $p_1,\dots, p_r$,
of index $q_j$ and $n-q_j$, $j=1,\dots,r$, respectively. Then
$$\chi(X)=\sum_{j=1}^r (-1)^{q_j} = \sum_{j=1}^r (-1)^{n-q_j}$$
as $n$ is odd, this implies that $\chi(X)=-\chi(X)$.

\cvd

\begin{remark}\label{handle-chi}{\rm If we consider the handle decomposition of $X$, say $\Hh$,  associated
to a Morse function $f$, the above expression of $\chi(X)$ can be rephrased in terms of
handle indices, that is $\chi(X)=\chi(\Hh)$ (see Section \ref {handle-moves}) . 
The characteristic $\chi(\Hh)$ is defined for every handle decomposition, not necessarily
associated (a priori) to any Morse function.  We know that it is invariant for the (handle) move-equivalence. 
}
\end{remark}  

\section{The index of an isolated zero of a tangent vector field}\label{zero-index}
We are going to reformulate the sign $\epsilon (p)$ of a non degenerate zero of a tangent vector field on $X$
in a way which will make sense also for any {\it isolated} zero (not necessarily non degenerate). 
Let $p$ be an isolated zero of a vector field $s$. Let us implement the following procedure:
\begin{enumerate}
\item Take local coordinates of $X$ at $p\sim 0$, so that $s$ is of the form
$$ s(x)= (x,f_p(x))$$
where
$$f_p: (\R^n,0)\to (\R^n,0)$$
is a smooth map such that $f^{-1}_p(0)=\{0\}$.
\item Then it is well defined the smooth map
$$ f_p/||f_p||: S^{n-1}\to S^{n-1} \ . $$
\item We can assume that the standard orientation of $\R^n$
associated to the standard basis is coherent with the global orientation
of $X$, so that $S^{n-1}$ is oriented as the boundary of the oriented disk
$D^n \subset \R^n$.
\item Finally set
$$i_p = \deg (f_p/||f_p||)\in \Z \ . $$
A priori this mights depend on the particular choices made in the implementation.
\end{enumerate}

We have
\begin{lemma}\label{well-def-ind} (1) $i_p(s):= i_p = \deg (f_p/||f_p||)\in \Z $ is well defined
(i.e. it does not depend on the specific implementation of the procedure) and is
called {\rm the index of the isolated zero $p$} of the tangent vector field $s$.
\smallskip

(2) If $p$ is a non degenerate zero of $s$, then (with the notations fixed above)
$$i_p(s) = \epsilon (p) = {\rm sign}(\det d_0f_p) \ . $$
\end{lemma}
\Dim  Let $\phi: (\R^n,0)\to (\R^n,0)$ be a change of coordinates relating
two different implementations. Then $\Dd^n:=\phi^{-1}(D^n)$ is a smooth oriented $n$-disk
around $0$, with oriented boundary $\Sigma$ diffeomorphic to $S^{n-1}$.
Let $s(x)=(x,f_p(x))$ be the expression of $s$ in the source local coordinates.
Set $$g:=f_p/||f_p||: \R^n\setminus \{0\} \to S^{n-1} \ . $$
It is clear that $i_p$ computed with respect to the target local coordinates
is equal to the degree of the restriction of $g$ to $\Sigma$. 
So we have to prove that this degree equals
$i_p$ computed with respect to source local coordinates. There is $1>\epsilon>0$
small enough such that the closed $n$-disk $\epsilon D^n$ (with boundary
$\epsilon S^{n-1}$) is contained in the interior of  $\Dd^n$.
Then the restriction of $g$ to $D^n \setminus {\rm Int}(\epsilon D^n)$
establishes an oriented bordism of $g_{|S^{n-1}}$ with $g_{|\epsilon S^{n-1}}$ ; 
similarly the restriction of $g$ to $\Dd^n \setminus {\rm Int}(\epsilon D^n)$
establishes a bordism of $g_{|\Sigma}$ with  $g_{|\epsilon S^{n-1}}$.
Then we can conclude by applying twice the invariance of the degree up to
bordism. This achieves (1).

As for (2), assume that $f_p$ is a diffeomorphism. The result is immediate
if $f_p$ is a linear isomorphism. Then we can conclude by means
of the results of Section \ref{linearization} and the invariance properties of
the degree again.

\cvd

\section{Index theorem}\label{Hopf-index}
Let $s$ be a tangent vector field on $X$ with only isolated zeros, say $p_1,\dots, p_r$
(there is a finite number because $X$ is compact). Then we can set
$$\chi(X,s)=\sum_{j=1}^r i_{p_j}(s) \ ; $$
if $s\pitchfork X$, that is all zeros are non degenerate, then we know that
$$\chi(X,s)=\chi(X)$$
has an intrinsic meaning, not depending on the field $s$.
Next theorem extends this fact to an arbitrary field as above.
\begin{theorem} For every tangent vector field $s$ on $X$ with only
isolated zeros, we have
$$\chi(X,s)=\chi(X) \ . $$
\end{theorem}
\Dim For every zero $p_j$ of $s$ fix an implementation of the procedure
that computes $i_{p_j}(s)$. Hence $i_{p_j}(s)= \deg (g_j:S^{n-1}_j\to S^{n-1})$. 
We can also assume that these charts are pairwise disjoint. Let $\tilde s$ be
a section of $T(X)$, $\tilde s \pitchfork X$, very close to $s$. Then the non degenerate
zeros of $\tilde s$ distribute in bunches $z_{j,1},\dots , z_{j,r_j}$, contained
in the interior of the $n$-disk $D^n_j$, $j=1,\dots, r$. Fix one of these zeros $p=p_j$
and consider the  corresponding $z_1,\dots, z_{r_j}\in D^n=D^n_j$.
We can take a system of pairwise disjoint small $n$-disks $D^n_{i}$ centred at $z_{i}$, contained in the interior of
$D^n$. As $s$ and $\tilde s$ are homotopic along $S^{n-1}=\partial D^n$, then we can use
$\tilde s$ instead of $s$ in order to compute $i_{p_j}(s)$ via the degree. On the other hand,
we can use the restriction of $\tilde s$ to $\partial D^n_i$ in order to compute the index
of the non degenerate zero $z_i$ of $\tilde s$. The normalized field is defined on
$D^n \setminus (\cup_i {\rm Int} (D^n_i))$ and this establishes a bordism between
the restriction on the boundary components. By the invariance of the degree up to bordism,
we realize that
$$ i_p(s)= \sum_i i_{z_i}(\tilde s) \ . $$
By taking the sum over all zeros of $s$ we eventually get
$$\chi(X,s)= \sum_{j,i} i_{z_{j,i}}(\tilde s) = \chi(X) \ . $$

\cvd 

\section{E-P characteristic for non oriented manifolds}\label{E-P-non-or}
Let us fix first the behaviour of $\chi$ with respect to the change of orientation.
So let $X$ be as above and $-X$ denotes it endowed with the opposite orientation. 
\begin{lemma}\label{EP_-X} $\chi(X)=\chi(-X)$.
\end{lemma}
\Dim Use a same given tangent field $s$ with isolated zeros to compute both characteristic
numbers. As $T(X)$ is tautologically oriented in agreement with the orientation of $X$,
it is immediate that  the index  of every zero of $s$ does not depend on the choice of
this orientation.

\cvd

\smallskip

In fact the computation of the index of an isolated zero $p$ of $s$ is a purely local stuff:
\smallskip

{\it One does not really need a {\rm global} orientation of $X$ to compute it; a local
orientation of $X$ at $p$ suffices and the same argument of the above lemma
shows that it does not depend on the choice of such local orientation.}

\smallskip

This suggests that the procedure to compute $\chi(X)$ can be extended to every
$X$ not oriented and even non orientable; it is enough to replace in the computation
of the indices a global orientation of $X$ (if any) with an arbitrary system of local
orientations at the zeros of a given tangent field $s$ with isolated zeros.
Then we have defined in general $\chi(X,s)$, which a priori  depends on the choice
of $s$. In fact it does not. If $X$ is orientable we have already achieved this
result. Assume that $X$ is connected and non orientable. Let $p: \tilde X \to X$
be the degree $2$ orientation covering of $X$, where $\tilde X$ is the connected orientable
total space of the unitary determinant bundle of $X$. Every field $s$ on $X$ as above
lifts to a field $\tilde s$ on $\tilde X$ so that every isolated zero $p$ of $s$ lifts
to a couple $p_{\pm}$ of isolated zeros of $\tilde s$. It follows from the very
definition that 
$$ i_p(s)= i_{p_\pm}(\tilde s)$$
so eventually
$$ \chi(X,s) = \frac{1}{2}\chi(\tilde X,\tilde s)= \frac{1}{2} \chi(\tilde X) \ . $$
Recall that if $X$ is orientable then $\tilde X$ consists of two copies of $X$
so that also in this case
$$\chi(X)=\frac{1}{2}\chi(\tilde X) \ . $$
Summing up
$$ \chi(X):= \frac{1}{2} \chi(\tilde X)\in \Z $$
is always a well defined characteristic number of $X$, and in every case
($X$ being orientable or not) can be computed as the sum of indices
of any tangent vector field $s$ on $X$ with isolated zeros. 

Recall that we have also the non oriented cobordism Euler class
$$ w^n(X)\in \eta^n(X)=\Bb^n(X;\Z/2\Z)= \Z/2\Z \ . $$
Clearly
$$ w^n(X)= \chi(X) \ {\rm mod} (2) $$
and sometimes one writes
$$ \chi_{(2)}(X):=w^n(X) \ . $$

\section{Some examples and  properties of $\chi$}\label{EP-examples}
\medskip

$\bullet$ The unit sphere $S^n$ admits a Morse function with just one minimum
and one maximum, then
$$\chi(S^n)= 1 + (-1)^n$$
and it is zero when $n$ is odd (as it must be), while $\chi(S^n)=2$ if $n$ is even.
This implies that an even dimensional sphere does not admit any nowhere
vanishing tangent vector field. In fact we have
\smallskip

{\it $S^n$ admits a nowhere vanishing tangent vector field if and only if $n$ is odd.}
\smallskip

We have to hexibit such a tangent vector field on $S^n$ when $n$ is odd.
For $n=1$, let $S^1\subset \R^2$ the unit circle. For every $p=(x,y)\in S^1$,
set $s(p)=(-y,x)$,this  does the job. In general for every
$p=(x_1,y_1,\dots, x_{n+1},y_{n+1})\in S^n \subset \R^{n+1}$, set
$s(p)= (-y_1,x_1,\dots , -y_{n+1},x_{n+1})$.
\smallskip

$\bullet$ If $\pi: \tilde X \to X$ is a degree $d$ covering map, then
$$ \chi(\tilde X)= d\chi(X) \ . $$
In fact we can argue as made above for the degree $2$ covering
maps, by lifting to $\tilde X$ any tangent vector field $s$ with isolated zeros on $X$;
every  zero $p$ of $s$ lifts to $d$ isolated zeros of $\tilde s$ sharing the same index of $i_p(s)$.
In particular by considering the natural degree $2$ covering map $\pi: S^n \to \PP^n(\R)$, we have
$\chi (\PP^n(\R))=0$ if $n$ is odd, while $\chi(\PP^n(\R))= 1$ if $n$ is even.
\smallskip

$\bullet$ Consider the complex projective space $\PP^n(\C)$ as the quotient
space of the unitary sphere $S^{2n+1}\subset \C^{n+1}$. One verifies (do it by exercise by using the
standard atlas of $\PP^{n}(\C)$ with $n+1$ complex affine charts)
that 
$$f([z_0,z_1,\dots,z_n])= \sum_{j=0}^n (j+1)|z_j|^2$$ 
defines a Morse function on $\PP^n(\C)$ with exactly $n+1$ critical points
$$p_0=[1,0,\dots, 0],\dots, p_n= [0,\dots, 0, 1]$$ 
and every even index between $0$
and $2n$ occurs exactly once. Hence
$$\chi(\PP^n(\C))=n+1 \ . $$
\smallskip

$\bullet$ The characteristic $\chi$ is {\it multiplicative with respect to the product of manifolds}.
That is, if $X$ and $X'$ are compact boundaryless manifold as above, then
$$ \chi(X\times X')= \chi(X)\chi(X') \ . $$
In fact if $s$ ($s'$) is a tangent field on $X$ (on $X'$) with non degenerate zeros
$p_1,\dots, p_r$ ($p'_1,\dots, p'_{r'}$), then $s\times s'$ defines a field on 
$X\times X'$ with $rr'$ non degenerate zeros $(p_j,p'_i)$, $j=1,\dots , r$, $i=1,\dots , r'$,
each one having index 
$$i_{(p_j,p'_i)}(s\times s')= i_{p_j}(s)i_{p'_i}(s') \ . $$
For example 
$$\chi (X\times S^1)=0$$ 
for every $X$ (in fact we can explicitly define
a nowhere vanishing tangent vector field on $X\times S^1$ which restricts to the avove
standard field on every fibre $\sim S^1$).

Whenever both $n$ and $m$ are even, then 
$$\chi(\PP^n(\R) \times \PP^m(\R))=1 \ . $$

\section{The relative E-P characteristic of a triad, $\chi$-additivity}\label{EP-triad}
Here we adopt the setting of Chapter \ref{TD-HANDLE}. By definition
a {\it relative tangent vector field on a triad} $(W,V_0,V_1)$ of compact
smooth manifolds, at the boundary $\partial W = V_0 \amalg V_1$  looks like a gradient
of a smooth function $f: W \to [0,1]$ such that $V_j=f^{-1}(j)$, $j=0,1$, and has no critical
points on a neighbourhood of $\partial W$. 
Hence it is ingoing $W$ along $V_0$ and outgoing along $V_1$.
An adapted gradient of any Morse function
on the triad is a typical example of such a field. By using these fields we can develop
with minor changes a notion of {\it relative of E-P characteristic for triads}. 
Assuming first that $W$ is oriented (with oriented boundary),
by using relative fields
with only non degenerate zeros
we can define the self-intersection number 
$$\chi(W,V_0)\in \Z$$
of $W$ in $T(W)$ relatively to $V_0$; it is well defined as
does not depend on the choice of the non degenerate field.
Then we can extend the Hops index theorem which allows us to compute 
$\chi(W,V_0)$ by means of any relative field with isolated zeros; finally we can
extend the definition of $\chi(W,V_0)\in \Z$ to non oriented and even non orientable
triads. Of course every $W$ with non empty boundary gives rise to
several triads $(W,V_0,V_1)$; among these: $(W,\emptyset, \partial W)$
and $(W,\partial W,\emptyset)$. The notation
$$\chi(W):= \chi(W,\emptyset,\partial W)$$
is compatible with 
$$\chi(W)=\chi(W,\emptyset,\emptyset)$$ 
when $W$ is boundaryless.

If $f: W\to [0,1]$ is a Morse function on the triad $(W,V_0,V_1)$, then
$\hat f = 1-f$ is a Morse function on $(W,V_1,V_0)$. By using  respective
adapted gradient fields to compute the relative charcteristics we get
\begin{lemma}\label{vuoto-bordo}
$$ \chi(W,V_0) = (-1)^{\dim(W)}\chi(W,V_1) \ . $$
\end{lemma}

\cvd
\smallskip

Note that $\chi(D^n)=1$: use a Morse function on $(D^n,\emptyset, S^{n-1})$
with just one minimum. 
\smallskip

If $X$ is boundaryless and $Y$ is with boundary, then the very
same argument used when also $Y$ is boundaryless
allows to extend the multiplicative property.
\begin{lemma}\label{mult}
$$ \chi(X\times Y)= \chi(X)\chi(Y)\ . $$
In particular
$$\chi(X\times D^n)=\chi(X) \ . $$
\end{lemma}

\cvd
\smallskip

\subsection{Additive property of $\chi$}\label{additiveEP}
If $(W,V_0,V_1)$, $(W',V_0',V'_1)$ are triads and
$\phi:V_1\to V'_0$ is a diffeomorphism, we get a new
composite triad $(W",V_0,V'_1)$, where
$W"= W \amalg_\phi W'$. Any couple of relative fields $v$ and $v'$
with isolated zeros on the given two triads respectively, can be glued together to
produce a relative field $v"$ having as zeros the union of the zeros of $v$ and $v'$
each one keeping its index. Then we have
$$\chi(W",V_0,V'_1)= \chi(W,V_0,V_1)+\chi(W',V'_0,V'_1) \ . $$
This additive property of $\chi$ has remarkable consequences.

\subsection{A baby TQFT}
In Section \ref{quick-TQFT} we have roughly outlined the axioms of a so called TQFT
and posed the question about the existence of any ``non trivial'' one.  Here we use $\chi$
to provide a baby but non trivial example. Consider {\bf CAT}$_\eta(n+1)$.
Associate to every object $M$ the vector space $Z(M)=\C$. To every arrow $f$
carried by any triad $(W,M_0,M_1)$, associate the unitary $\C$-linear map 
$$Z(f): Z(M_0) \to Z(M_1), \ z \to e^{i\chi(W,M_0)}z \ . $$
By using the additive property of $\chi$ it is easy to check that all axioms are verified.
This shows at least that there are not logical contradictions within the given pattern
of axioms.

\section{E-P characteristic of tubular neighbourhoods and the Gauss map}\label{EP-tube}
The above equality $\chi(X\times D^n)=\chi(X)$ is a special case of the following
\begin{proposition}\label{tube} Let $p:U\to X$ be a closed tubular neighbourhood
of a submanifold $X$ of some $Y$. Then $\chi(U)=\chi(X)$.
\end{proposition} 
\Dim It is enough to show the equality for an $\epsilon$-neighbourhood
$N_\epsilon(X)$ of the zero section $X$ of a vector bundle bundle $\pi:E\to X$
endowed with a field of positive definite scalar products on every fibre.
Let $v$ be a tangent field on $X$ with non degenerate zeros. Define the field
on $N_\epsilon(X)$
$$w(z)= (z-p(z))+ v(p(z)) \ . $$
One checks that $w$ is a field on the triad $(N_\epsilon(X),\emptyset,\partial N_\epsilon(X))$,
the zeros of $w$ coincide with the zeros of $v$, are non degenerate and keep the sign.
The Proposition follows.
 
 \cvd
\smallskip

In the special case  $X\subset \R^k$, assume that $U$ has been constructed
by means of the standard metric on $\R^k$.
By removing from the interior of $U$
a system of pairwise disjoint small open disks $\Dd_p$ around every zero of
$w$, we get a manifold $W$ with boundary $(\amalg_p S^{k-1}_p)\amalg \partial U$
on which the normalized field $\wG:=w/||w||$ is well defined, as well
as the map $\wG: W \to S^{k-1}$. The restriction of $\wG$ to $\partial U$
is a field of unitary vectors pointing out from $U$ (in fact normal to the boundary). 
This restriction, say $g_{\partial U}$
is called the {\it Gauss map of the hypersurface $\partial U$}.
By computing the characteristic as sum of zero indices and by means of the
bordism invariance of the degree, finally we have
\begin{corollary} Let $p:U\to X$ be a tubular neighbourhood of $X$ in $\R^k$.
Then
$$\chi(X) = \deg (g_{\partial U})$$
where $g_{\partial U}$ is the Gauss map of the hypersurface $\partial U$.
\end{corollary}

\cvd
\smallskip

If $X$ itself is an oriented hypersurface in $\R^k$, we can define
its Gauss map $g_X: X\to S^{k-1}$ as the unitary field of
normal vectors along $X$ such that followed by the orientation of $X$
produce the given (standard) orientation of $\R^k$ along $X$.
In this case $\partial U$ consists of two parallel copies of $X$ with opposite
orientations. The last corollary specializes to
\begin{corollary}\label{gauss-map} Let $X$ be an oriented hypersurface
of $\R^k$ then
$$\chi(X)= 2\deg(g_X)$$
where $g_X$ is the Gauss map of $X$.
\end{corollary}

\cvd

\begin{remark}\label{proj-again}{\rm We can compute inductively
the characteristic of  real and complex projective spaces by
decomposing  $\PP^n(K)$ as the union of a tubular neighbourhood
of $\PP^{n-1}(K) \subset \PP^n(K)$ and its complement, and
applying Proposition \ref{tube} together with the additivity
of $\chi$.}
\end{remark} 
 
\section {Non triviality of $\eta_\bullet$ and $\Omega_\bullet$}\label{EP-bord}
The integer E-P characteristic is {\it not}  invariant up to bordism.
For example $[S^2]=[S^1\times S^1]=0\in \eta^2$, but
$\chi(S^2)=2 \neq 0 = \chi(S^1\times S^1)$. On the other hand
the E-P characteristic mod$(2)$ is bordism invariant.
\begin{proposition}\label{EP-2-invariant} Let $[X]=0\in \eta^{n}$.
Then $\chi_{(2)}(X)=0\in \Z/2\Z$.
\end{proposition}
\Dim If $n$ is odd, we know in general that $\chi(X)=0$.
Assume that $n$ is even.  Let $X= \partial W$, $W$ being a compact manifold with
boundary. Take the double $D(W)$. $D(W)$ can be presented as the composition
of the triad $$(X\times D^1, \emptyset, (X\times \{-1\}\amalg (X\times \{1\}))$$
followed by two copies of the triad $$(W, X,\emptyset)$$ 
glued to $X\times D^1$
along $X\times \{\pm 1\}$ respectively.  By the additive property
$$\chi(D(W))= \chi(X\times D^1) + 2\chi(W,X,\emptyset) \ . $$
By Lemmas \ref{mult}, \ref{vuoto-bordo} and the facts that $n+1$ is odd
and the double is boundaryless we have
$$ \chi(X)= \chi(D(W)) - 2\chi(W,X,\emptyset)= 2\chi(W)\in \Z$$
so that
$$ \chi_{(2)}(X)=0 \in \Z/2\Z \ . $$

\cvd
\smallskip

As an immediate corollary we have the {\it non triviality of $\eta_{2n}$ and $\Omega_{4n}$}.

\begin{corollary} For every even $n\geq 1$, $\eta_{2n} \neq 0$ and $\Omega_{4n}\neq 0$
\end{corollary}
\Dim We know that $\chi(\PP^{2n}(\R))=1$, hence $[\PP^{2n}(\R)]\neq 0 \in \eta_{2n}$.
Similarly $\chi_{(2)}(\PP^{2n}(\C))=1$, hence $[\PP^{2n}(\C)]\neq 0 \in \Omega_{4n}$.
\cvd
\medskip

By using the multiplicative property of $\chi$ and the obvious fact that it is
additive under disjoint union we also have
\begin{corollary} $\chi_{(2)}: \eta^\bullet \to \Z/2\Z$ is a well defined non trivial
{\rm ring} homomorphism.
\end{corollary}

\cvd
\smallskip

In fact every $\PP^a(\R)\times \PP^b(\R)$, $a$, $b$ even, is non trivial in $\eta_{a+b}$.
For example in $\eta_4$ we have the non trivial $[\PP^4(\R)]$, $[\PP^2(\R)\times \PP^2(\R)]$, $[\PP^2(\C)]$.
At present we are not able to decide if they are equal or not. 
Similarly we have that $\PP^a(\C)\times \PP^b(\C)$, $a$, $b$ even, is non trivial in $\Omega_{2(a+b)}$.

\section{Combinatorial E-P characteristic}\label{comb-chi}
We have treated the E-P characteristic of  smooth manifolds
in purely differential/topological terms. However, 
the reader is probably aware that the name E-P characteristic
is used in other different settings. Probably she/he has at least encountered 
a combinatorial formula producing the value $2=\chi(S^2)$
for every polyhedral realization of the sphere as the boundary of
a convex polytope in $\R^3$. In this very  sketchy Section we would outline
a few bridges between such different ways to recover the E-P characteristic.

\subsection{Piecewise smooth triangulations and the combinatorial characteristic}\label{simp} 
Recall that a $m$-{simplex} $\sigma$ in some euclidean space $\R^h$, $h\geq m$
is the convex hull of $m+1$ affinely independent points (that is they span
an $m$-dimensional affine subspace of $\R^h$). These are called the {\it vertices}
of  $\sigma$. By removing one vertex, say $p$, we detemine a $(m-1)$ simplex
$\sigma_p$ which is the $(m-1)$ face of $\sigma$ opposite to the vertex $p$.
By iterating the face operation we get the iterated $k$-faces of $\sigma$, $0\leq k \leq m$,
where the vertices are the $0$-faces and $\sigma$ itself is the unique $m$-face.
By definition a {\it finite simplicial complex} is a finite family $\Kk$ of simplexes 
in some $\R^h$ such that
\begin{itemize}
\item $\Kk$ is closed with respect to the iterated faces.
\item Two simplexes of $\Kk$ may intersect each other only at a common
iterated face.
\end{itemize}

The union $|\Kk|$ of the simplexes of $\Kk$ is a subspace of $\R^h$
called the {\it geometric support} of the complex $\Kk$.
\medskip

Let $X$ be a compact boundaryless smooth manifold.
A {\it piecewise smooth triangulation of $X$} is given by
a homeomorphism
$$\tau: |\Kk|\to X$$
where $\Kk$ is a finite simplicial complex in some $\R^h$ and the
restriction of $\tau$ to every $n$-symplex of $\Kk$ is a smooth embedding
in $X$. If $\partial X\neq \emptyset$, we require furthermore that $\tau_{|\tau^{-1}(\partial X)}$
is a triangulation of $\partial X$.

One can prove
\begin{proposition}\label{whitney-field} Let $\tau: |\Kk|\to X$ be a piecewise smooth triangulation
of the compact boundaryless smooth $n$-manifold $X$. Then there is a tangent, so called {\rm Whitney}
vector field $v_\tau$ on $X$ whose zero set coincides with the set of images
of the barycenters $\hat \sigma$ of the simplexes $\sigma$ of $\Kk$  and every zero has index equal to
$(-1)^{\dim (\sigma)}$.
\end{proposition}

\medskip

As a Corollary of the above Proposition and the Index Theorem, we have

\begin{corollary} For every piecewise triangulation $\tau: |\Kk| \to X$ as above,
$$\chi(X) = \sum _{j=0}^n (-1)^j c_j:= \chi(\Kk)$$
where $c_j$ is the number of $j$ dimensional simplexes of $\Kk$.
In particular the {\rm combinatorial characteristic} $\chi(\Kk)$ does not depend
on the choice of the triangulation of $X$.
\end{corollary}

\cvd
\smallskip

A few comments about the proof. The Whitney field $v_\tau$  can be explicitly given by means of barycentric
coordinates on the simplexes of $\Kk$, see for instance \cite{HT}; every barycenter of a $n$-simplex of 
$\Kk$ corresponds to a source of $v_\tau$, every vertex of $\Kk$ corresponds to a pit, in general 
every barycenter of a $j$-simplex corresponds to a saddle point with a $j$ dimensional space of 
ingoing directions tangent to the simplex and a $n-j$-dimensional space of outgoing directions
tranverse to the simplex.
\medskip

For the {\it existence} (and a suitable form of ``uniqueness up to subdivision'') 
of piecewise smooth triangulations see \cite{Mu}. 

\subsection{Homological characteristic}\label{hom-chi}
Here we want to recover the combinatorial characteristic in an algebraic/topological
setting. 

Fix any field $F$  (for example $F=\Z/2\Z, \Q, \R, \C$). 

Given a triangulation of $X$ as above, we can define the 
{\it simplicial homology of $\Kk$ with coefficients in $F$} as follows:
\smallskip

$\bullet$ Give every simplex $\sigma$ of $\Kk$ an orientation, induced by the
choice of an orientation on the affine subspace of $\R^h$ spanned by the vertices
of $\sigma$.
\smallskip

$\bullet$ For every $0\leq j\leq n$, set $C_j(\Kk;F)$ the finite dimensional $F$-vector space
having as a basis the oriented $j$-simplexes of $\Kk$ (note that $-\sigma$ considered
as the simplex endowed with the opposite orientation is confused with
$-1\sigma$ i.e. the product of $\sigma$ with the scalar $-1\in F$). Hence 
$$\dim C_j(\Kk;F)=c_j \ . $$
\smallskip

$\bullet$ Every $(j-1)$-face $\sigma'$ of the oriented $j$-simplex $\sigma$ of $\Kk$ 
inherits a boundary orientation accordingly with our usual convention.
Hence $\sigma'$ has two orientations, the one fixed above and the boundary
orientation. Give it the sign $\epsilon(\sigma',\sigma)=  1$ if these 
orientations agree to each other, the sign $-1$ otherwise. Then define the unique 
$F$-linear map
$$\partial_j: C_j(\Kk;F)\to C_{j-1}(\Kk,F)$$ 
which on every oriented $j$-simplex $\sigma$ holds:
$$\partial_j (\sigma)= \sum_{\sigma'} \epsilon(\sigma',\sigma)\sigma'$$
where $\sigma'$ varies among the $(j-1)$ faces of $\sigma$.
It is not hard to verify that
$$\delta_{j-1}\circ \delta_j = 0$$
basically because two $(j-1)$ faces of the oriented $j$-simplex 
$\sigma$ both endowed with the boundary orientation, induce
opposite boundary orientations on their common $(j-2)$ face of $\sigma$.
Hence we can define the quotient $F$-vector spaces
$$H_j(\Kk;F)= \ker (\delta_j)/{\rm Im}(\delta_{j+1}) $$
and these are the desired simplicial $F$-homology spaces of the complex
$\Kk$. By using the elementary dimension formula for any finite dimensional
linear map $f:V\to W$:
$$ \dim (V)= \dim (\ker(f)) + \dim ({\rm Im}(f))$$ it is not hard to check
that the $F$-{\it homological characteristic}
$$ \chi(H_\bullet(\Kk;F)):=\sum_{j=0}^n (-1)^j \dim H_j(\Kk;F) = \sum_{j=0}^n (-1)^j\dim C_j(\Kk;F)$$
hence it equals the combinatorial characteristic so that
$$ \chi(H_\bullet(\Kk;F)) = \chi (X) \ . $$
Remarkably it does not depend on the choice of the triangulation of $X$ and not even of the field $F$. 
It is a fundamental and basic result of algebraic topology (see \cite{Hatch}, \cite{Mu2}) that even the single
dimensions (also called the $F$-{\it Betti numbers} of $X$) 
$$\dim H_j(X;F):=\dim H_j(\Kk;F)$$
do not depend on the choice of the triangulation, although they depend on $F$.

\chapter{Surfaces}\label{TD-SURFACE}
We are going to apply  several tools developed in the previous Chapters in order to classify
the compact surfaces (i.e. smooth $2$-manifolds) and also to determine both bordisms
$\eta_2$ and $\Omega_2$.

Let $M$ be a compact connected boundaryless surface.

$\bullet$ We know from Chapter \ref{TD-HANDLE} that 
\smallskip

{\it $M$ admits a `reduced' ordered handle decomposition with one $0$-handle, followed by say $\kappa$ disjoint 
$1$-handles
and one final $2$-handle, where  $\kappa := \kappa(M)$ is intrinsically determined by
$$\kappa(M) = 2 -\chi(M)  \ . $$}
\smallskip

In fact recall that for any handle decomposition $\Hh$ of $M$, its characteristic 
$$\chi(\Hh):= \sum_{j=0}^2 (-1)^{j}b(j)$$
$b(j)$ being the number of index $j$ handles, is preserved by the basic moves on 
handle decompositions; if $\Hh$ is associated to a Morse function
on $M$, then $\chi(\Hh)=\chi(M)$; finally we can get a reduced ordered decomposition of $M$
by performing some basic moves on any given decomposition.

\begin{remark}{\rm For any ordered handle decomposition of $M$ with one $0$-handle, one
$2$-handle and a few disjoint $1$-handles, it is not hard to triangulate $M$ in the following
way: take a vertex internal to every handles; take as further vertices on the boundary of the $0$-handle
the union of the boundaries  of the attaching $1$-disks of the $1$-handles; they also provide
a triangulation of the boundary of every $1$-handle; triangulate both the one $0$-handle and
every $1$-handle by the cones on the boundary with centre at the respective internal vertex;
these triangulations match and give a triangulation of the union of 
the $0$-handle with the $1$-handles; the resulting surface has as boundary a triangulated circle; 
finally complete it to a triangulation of the whole of $M$
by means again of the cones with centre at the internal vertex of the $2$-handle. 
By using the combinatorial computation of $\chi(M)$  applied to such a
triangulation one can easily check that the number of $1$-handles is always equal to $\kappa(M)$.}
\end{remark}
\smallskip 

$\bullet$ Recall from Section \ref{twisted-sphere} that 
\smallskip

{\it In dimension $2$ connected sum and weak connected sum are equivalent to each other; 
moreover every twisted $2$-sphere
is {\rm diffeomorphic} to the standard $S^2$.}

\smallskip

So let $\gamma \subset M$ be any dividing  
connected simple curve $\gamma$,
that is 
$$M\setminus \gamma = N_1\amalg N_2$$
where $N_j$ is a non empty connected open set of $M$
and the closure $\bar N_j$ is a compact smooth manifold
with boundary $\partial \bar N_j = \gamma$; let $M_j$ be the
boundaryless surface obtained from $\bar N_j$ by filling $\partial \bar N_j$
with a $2$-disk  glued along the boundary; then (up to diffeomorphism)
$$M\sim M_1 \# M_2 \ . $$
To be more precise, the result is uniquely determined if
at least one among $M_1$ and $M_2$ is non orientable, or both are orientable and admit orientation 
reversing diffeomorphisms. In general let us say that
$M$ is ``a" connected sum of $M_1$ and $M_2$ 
(this precision will be eventually immaterial). By the additive property of $\chi$,
we have
$$ \kappa(M)=\kappa(M_1)+\kappa(M_2) \ . $$ 

\smallskip

$\bullet$  Let us consider $\eta_1(M)=\Bb_1(M;\Z/2\Z)$. 
\begin{lemma}\label{finite-dim}  $\eta_1(M)$ is $\Z/2\Z$-vector space of finite dimension $\leq \kappa(M)$.
\end{lemma}
\Dim By Section \ref{b-vs-homot} there is a {\it surjective} 
homomorphism (a base point being understood)
$$ \pi_1(M)\to \eta_1(M) \ . $$
 By using a reduced ordered handle decomposition of $M$ as above and applying (an elementary version of) Van Kampen theorem
we see that $\pi_1(M)$ has a presentation with $\kappa$ generators and one relation; for the union of the $0$-handle with the $\kappa$ $1$-handles
has the homotopy type of a wedge of $\kappa$ copies of $S^1$ whose fundamental group is a free group with $\kappa$ generators; the defining relation
between them is given by the attaching map of the $2$-handle. The Lemma follows.

\cvd
 
\smallskip

\begin{lemma}\label{connected-rep}
Every $\alpha \in \eta_1(M)$ can be represented by a {\rm connected} simple smooth curve 
$C$ traced on $M$. 
\end{lemma}
\Dim We already know from general results in Chapter \ref{TD-LINE-BUND}
that a codimension $1$ class can be represented by hypersurfaces. 
In the present $2$-dimensional 
situation we can get an elementary direct proof of this fact as follows. 
Certainly $\alpha = [f: \tilde C \to M]$ where $\tilde C$ is a finite union of copies of $S^1$.
By a standard `general position' argument (see Section \ref {miscellaneaT}) we can assume that up to homotopy, 
hence up to bordism, $f:\tilde C \to M$ is a generic immersion possibly having only simple double points in its image $f(\tilde C)\subset M$.
In local coordinates every crossing of $f(\tilde C)$ is of the form $\{xy=0\} $ and has two local `simplifications' of the form $\{xy\pm\phi(x,y)\epsilon=0\}$
where $\epsilon>0$ is small enough and $\phi$ is a suitable bump function with support in a small disk centred at $0$.  
By locally simplifying every crossing of $f(\tilde C)$ (choose arbitrarily one way) we get a $1$-submanifold $C'$ of $M$. 
It is not hard to verify that $\alpha = [f:\tilde C\to M]=[C'] \in \eta_1(M)$, 
this is left as an exercise. In general $C'$ is not connected.  In order to modify $C'$ to get a connected representative $C$ of $\alpha$,
first we can remove all dividing components of $C'$ (keeping the name); if $C'$ is not connected then apply the following argument 
that decreases the number of components by $1$. We can find  two components $C_1$ and $C_2$ of $C'$ 
which can be connected by a smooth arc $I$ whose internal part is embedded into $M\setminus C'$, one endpoint $x_j$ is on $C_j$, $j=1,2$,
and is tranverse to $C_1\cup C_2$. $I$ can be thikened to an embedded $1$-handle $H\sim I\times [-1,1]$ which intersects
$C_j$ at $\{x_j\} \times [-1,1]$ and is contained in $M\setminus C'$ elsewhere. Then consider 
$$ C":= (C'\setminus (C_1\cup C_2))\cup C^*$$
where
$$C^* = ((C_1\cup C_2)\setminus  H)\cup (I\times \{\pm 1\})$$
up to corner smoothing. Hence $C_1\cup C_2$ has been replaced by the connected
curve $C^*$.  Again it is not hard to show that $[C']=[C"] \in \eta_1(M)$.
By iterating the procedure we eventually get a required connected representative
$C$ of $\alpha$.

\cvd

$\bullet$ Consider now the symmetric intersection form  (Section \ref{intersection-form})
$$ \bullet = \bullet_M: \eta_1(M)\times \eta_1(M)\to \Z/2\Z \ . $$
We have
\begin{lemma}\label{2Dnon-degenerate} The intersection form on $\eta_1(M)$ is non degenerate.
\end{lemma}
\Dim We have to show that if  $\alpha\neq 0$ in $\eta_1(M)$, then there is $\beta \in \eta_1(M)$ such that $\alpha \bullet \beta =1\in \Z/2\Z$.
Let $C\subset M$  be a connected smooth representative of $\alpha$ as in Lemma \ref{connected-rep}.
As $\alpha \neq 0$, then $M\setminus C$ is connected (otherwise $C$ would be the boundary of the closure of a component
of $M\setminus C$, so that $[C]=0$). Take a fibre $I$, necessarily tranverse to $C$ at one point, of a tubular neighbourd of $C$ in $M$.
Also $M\setminus (C\cup I)$ is connected, so that 
the endpoints of the interval $I$ can be connected by a smooth simple arc $\gamma$ whose internal part is contained in 
$M\setminus (C\cup I)$. Then (possibly by corner smoothing) $C':=I \cup \gamma$ is a smooth boundaryless curve in $M$
which intersects $C$ transversely at one point, hence $[C]\bullet [C']=1$.

\cvd

The next Lemma follows  from Chapter \ref{TD-LINE-BUND} .

\begin{lemma}\label{tube-alternative} Let $C\subset M$ be a connected smooth boundaryless curve.
Then there are two possibilities: either $[C] \bullet [C] =1$ and this happens if and only
if $C$ has  tubular neighbourhood in $M$ diffeomorphic to a M\"obius band, or $[C] \bullet [C] =0$   and this happens if and only
if $C$ has  a product tubular neighbourhood in $M$.
\end{lemma}

\cvd

The following Lemma  is obvious

\begin{lemma}\label{int-form-invariance}
 If $f: M\to M'$ is a surface diffeomorphism, then $$f_*: (\eta_1(M),\bullet_M)\to (\eta_1(M'),\bullet_{M'})$$
 is an isometry, that
 is $f_*$ is a $\Z/2\Z$-linear isomorphism and for every $\alpha,\beta \in \eta_1(M)$, 
 $$\alpha \bullet_M \beta = f_*(\alpha)\bullet_{M'} f_*(\beta) \ . $$
 \end{lemma}
 
 \cvd
 
 Hence {\it the isometry class of the non degenerate symmetric intersection form on $\eta_1(*)$ is an invariant up to diffeomorphism}.
\smallskip
 
 In what follows we will make the abuse of confusing a form with its isometry class. If $(V,\rho)$ and $(V',\rho')$
 are finite dimensional $\Z/2\Z$-vector spaces endowed with non degenerate symmetric forms, we can define the {\it orthogonal direct sum}
 $(V,\rho)\perp (V', \rho')$ which denotes the non degenerate symmetric form $\rho \perp \rho'$ on $V\oplus V'$ 
 that restricts  to $\rho$ (resp. $\rho'$) on $V$ ($V'$) and  such that $V$ and $V'$ are orthogonal to each
 other. We have
 \begin{lemma}\label{connected-direct-sum}   If the surface $M$ is a connected sum $$M\sim M_1  \# M_2$$ then
 (up to isometry) $$(\eta_1(M),\bullet_M) = (\eta_1(M_1), \bullet_{M_1}) \perp (\eta_1(M_2), \bullet_{M_2})  \ . $$ 
 \end{lemma}
 \Dim We can assume that the connected  sum has been realized from a connected dividing curve $\gamma$ in $M$ as at the beginning
 of the section (we adopt those notations).  It is easy to see that the linear map $i_*: \eta_1(N_j) \to \eta_1(M_j)$ induced by the inclusion is an isomorphism,
 $j=1,2$.  Denote by $V_j$ the image of $\eta_1(N_j)$ in $\eta_1(M)$ by the inclusion. It is evident that $V_1$ and $V_2$ are
 orthogonal to each other with respect to $\bullet_M$. It is enough to show that $\eta_1(M)=V_1 + V_2$,  whence $\eta_1(M)=V_1\perp V_2$
 because $\bullet_M$ is non degenerate, and that  $V_j$ is actually isomorphic to $\eta_1(N_j)$, $j=1,2$. 
Let $\alpha \in \eta_1(M)$ and $C\subset M$ be a smooth
 representative as above. By transversality we can assume that $C\pitchfork \gamma$. As $[\gamma]=0$ in $\eta_1(M)$,
 then $C\cap \gamma$ consists of an even number  of points $\{ p_1,\dots, p_{2d}\} $. We can assume that they are the endpoints
 of a family $\{I_1,\dots, I_d\}$ of pairwise disjoint intervals embedded into $\gamma$.  
 Take a `small' tubular neighbourhood $U\sim \gamma \times [-1,1]$
 of $\gamma$ in $M$. Then $M\setminus U$ consists of two connected components $W_1$ and $W_2$ such that $W_j\subset N_j$.
 The boundary of $W_j$ is a parallel copy $\gamma_j$ of $\gamma$. Denote by $I_{i,j}$, $j=1,2, \ i=1,\dots d$, the parallel copy in $\gamma_j$
 of the interval $I_i$. Finally for $j=1,2$, set
 $$ C_j = (C\cap W_j)\cup (\bigcup_{i=1}^d I_{i,j}) \ . $$
 Up to corner smoothing, $C_j$ is a smooth curve (not necessarily connected) in $N_j$ and it is easy to see that
 $$ [C_1 \amalg C_2]=[C]\in \eta_1(M)$$
 this shows that $\eta_1(M)=V_1+V_2$. Finally let  $\alpha \in \eta_1(N_1)\sim \eta_1(M_1)$
 and denote by $\alpha'$ its image in $\eta_1(M)$ by the inclusion. If $\alpha$ is not zero,
 as $\bullet_{M_1}$ is non degenerate, then there is $\beta \in \eta_1(N_1)$ such that $\alpha \bullet_{M_1} \beta = 1$;
 due to the geometric way one computes the intersection forms, it follows that also $\alpha' \bullet_M \beta'=1$, whence
 $\alpha'$ is non zero.
 
 \cvd 
 
 \smallskip
 
 We are going to see that the isometry class of the intersection form contains all relevant information about the diffeomorphism
 class. 
 
 \section{ Classification of symmetric bilinear forms on $\Z/2\Z$}\label{algebraic-form}
 Here we classify up to isometry non degenerate symmetric bilinear forms on finite dimensional $\Z/2\Z$-vector spaces.
 We denote by $\UU$ the unique $1$ dimensional isometry class; by $\HH$ the isometry class
 of {\it hyperbolic planes}, i.e.  $2$-dimensional spaces endowed with a non degenerate symmetric form admitting a basis made by isotropic
 vectors (recall that a vector $v$ is {\it isotropic} for a form $\beta$ if $\beta(v,v)=0$). Note that although $\HH$ is non degenerate
 it is {\it totally isotropic} (every vector is so), this depends on the fact that the characteristic of the field $\Z/2\Z$ is equal to
 $2$, in characteristic $\neq 2$ the zero form is the only  totally isotropic one by the so called `polarization formula' . 
 For every $n\geq 1$, denote by $n\UU$ (resp. $n\HH$) the orthogonal direct sum of $n$ copies of $\UU$ (resp. of $\HH$).
 We have
 \begin{proposition}\label{classification} Let $(V,\beta)$ be a finite dimensional $\Z/2\Z$-vector space endowed with a non degenerate
 symmetric bilinear form, $\dim V >0$. Then we have one of the following exclusive occurences:
 \smallskip
 
 (1) $(V,\beta)$ admits an orthogonal basis so that it is isometric to $n\UU$, $n=\dim V$, and this happens if and only if it is not totally isotropic.
 \smallskip
 
 (2) $\dim V=2n$, $(V,\beta)$ is isometric to $n\HH$, and this happens if and only if it is totally isotropic.
 \end{proposition}
 \Dim Assume first that $(V,\beta)$ is totally isotropic. Let $\Bb=\{v_1,\dots v_k\}$ be a basis of $V$, $\Bb^*=\{v^*_1,\dots, v^*_k\}$
 its dual basis, $w_1$ the vector which represents the functional $v^*_1$ by means of the non degenerate form $\beta$.
 Then the subspaces spanned by $\{v_1,w_1\}$ endowed with the restriction of $\beta$ is a hyperbolic plane $\HH$. As this last is non
 degenerate, then (up to isometry)
 $$  (V,\beta)= \HH \perp \HH^\perp$$
 all spaces being endowed with the restriction of $\beta$. Clearly also the restriction to $\HH^\perp$ is non degenerate and
 totally isotropic, $\dim \HH^\perp = \dim V -2$. So we can achieve the item (2) by induction on the dimension.
 Assume now that $v\in (V,\beta)$ is not isotropic. Then the subspace spanned by $v$ endowed with the restriction of $\beta$
 represents $\UU$ and (up to isometry)
 $$ (V,\beta)=\UU \perp \UU^\perp \ . $$
 By iterating the argument we get that either
 $$(V,\beta)= n\UU, \ n=\dim V$$
 and we have done, or
 $$(V,\beta)= k\UU \perp T$$
 for some $k\geq 1$ where $T$ is totally isotropic, $\dim T>0$. 
 We apply (2) to $T$, and get
 $$ (V,\beta)= k\UU \perp h \HH $$
 for some $k, h\geq 1$. Finally item (1) is achieved by means of the following Lemma.
 Note by the way that it also shows that $\perp$ does not verify the `cancellation properties'.
 
 \begin{lemma}\label{basic-id} Up to isometry $\UU\perp \HH = 3\UU$.
 \end{lemma}
 \Dim  Let $\Dd = \{u,w,t\}$ be a basis for $\UU \perp \HH$ adapted to the decomposition so that $\{w,t\}$
 is a basis of the hyperbolic plane.  Let $N$ be the subspace spanned by  $\{ u + w, u + t \}$.
 One readily verifies that this last is a orthogonal basis of $N$ so that $N=2\UU$. 
 Then $\UU \perp \HH = N \perp N^\perp$ and the last space is $1$ dimensional and non degenerate,
 so eventually $\UU \perp \HH = 3\UU$.
 
 \cvd
 \smallskip
 
 Also the proof of Proposition \ref{classification} is now complete.
 
 \cvd
 
 \smallskip
 
 \section{Classification of compact surfaces}\label{surface-class}
 We are going to prove the following topological classification theorem.
 \begin{theorem}\label{classification2}
 (0) Let $M$ be a compact connected boundaryless surface. Then
 the following facts are equivalent to each other.
 \begin{itemize}
 \item $M$ is diffeomorphic to $S^2$;
 \item $\kappa(M)=0$;
 \item $\dim \eta_1(M)=0$.
 \end{itemize}
 \smallskip
 
 (1) For every $n\geq 1$, the isometry class $n\UU$ is realized by
 the intersection form of $\eta_1(n\PP^2(\R))$ where $n\PP^2(\R)$ denotes the connected
 sum of $n$ copies of the real projective plane.
 
 \smallskip
 
 (2) For every $n\geq 1$, the isometry class $n\HH$ is realized by
 the intersection form of $\eta_1(n(S^1\times S^1))$, where $n(S^1\times S^1)$ denotes the connected
 sum of $n$ copies of the torus.
 
 \smallskip
 
 (3) Two compact connected boundaryless surfaces $M$ and $M'$ are diffeomorphic
 if and only if the intersection forms on $\eta_1(M)$ and $\eta_1(M')$ respectively
 are isometric to each other.
 \end{theorem}
 
 \smallskip
 
 This theorem has several interesting corollaries.
 
 \begin{corollary}\label{class-cor}
 In the hypotheses of Theorem \ref{surface-class}:
 \smallskip
 
 (1) $\dim \eta_1(M)= \kappa(M)= 2-\chi(M)$.
 If $M$ is orientable then $\kappa(M)=2g(M)$ is even
 ($g(M)$ is called the {\rm genus} of $M$).
 \smallskip
 
(2) Two surfaces $M$ and $M'$ are diffeomorphic if and only if 
 $\chi(M)=\chi(M')$ and either they are both orientable or
 non orientable. 
 \smallskip
 
 (3) Every orientable surface $M$ admits orientation reversing diffeomorphisms.
 Hence the connected sum of two surfaces $M=M_1  \# M_2$ is always
 uniquely defined up to diffeomorphism.
 \smallskip
 
 (4) Every $M$ can be embedded into $\R^4$. If $M$ is orientable then it can be
 embedded into $\R^3$. 
  
 \end{corollary}
 
 \smallskip
 
 {\it Proofs.} First item $(0)$ of Theorem \ref{classification2}, that is
 the characterization of the $2$-sphere up to diffeomorphism. If $\kappa(M)=0$, then $M$ 
 has a handle decomposition with only one $0$-handle and one $2$-handle. 
 So it is a twisted $2$-sphere, whence
 it is diffeomorphic to $S^2$. Then $M$ is simply connected, hence
 $\dim \eta_1(M)=0$. Let us show now that if $\kappa(M)>0$ then $\dim \eta_1(M)>0$.
Take a reduced ordered handle decomposition with $\kappa(M)$
 $1$-handles. The core of every $1$-handle can be completed with a simple arc embedded into
 the $0$-handle to get a connected simple smooth curve $C$ in $M$. There are two possibilities:
 either for one such a curve $[C]\bullet_M [C]=1$ or there are two such curves $C$ and $C'$ such that
 $[C]\bullet [C']=1$ (here we use that the boundary of the union of the $0$-handle with the disjoint
 $1$-handles must be connected). In any case $\dim \eta_1(M)>0$.  
 The other implications of item $(0)$ are evident.  
  
Let us show now that  that $\UU$ and $\HH$ can be realized.
 $\PP^2(\R)$ can be obtained by gluing a $2$-disk along the boundary of a M\"obius
 band. By Van Kampen theorem we realize that $\pi_1(\PP^2(\R))\sim \Z/2\Z$ and is generated
 by the core $C$ of the M\"obius band. Another way to check this fact is by means of the orientation
 covering $S^2\to \PP^2(\R)$. Then also $\eta_1(\PP^2(\R))\sim \Z/2\Z$, generated by $[C]$
 and $[C]\bullet [C]=1$. The above M\"obius band can be realized by attaching one $1$-handle
 to an initial $0$-handle, and we get $\PP^2(\R)$ by adding one final $2$-handle; this provides a 
 reduced ordered handle decomposition with $\kappa (\PP^2(\R)) = 1$ $1$-handle.
 By the way we realize also that if $\kappa(M)=1$ then $M$ is diffeomorphic to $\PP^2(\R)$.
 
 The fundamental group $\pi_1(S^1\times S^1)\sim \Z\oplus \Z$
 generated by the simple loops $C_1=S^1\times \{b_0\}$, $C_2=\{a_0\}\times S^1$ with base point
 $(a_0,b_0)$. It is immediate that $[C_1]\bullet [C_2]=1$ in $\eta_1(S^1\times S^1)$, 
 while $[C_j]\bullet [C_j]=0$, $j=1,2$. Hence $[C_1]$ and $[C_2]$ are non zero and linearly independent,
 $\dim \eta_1(S^1\times S^1)=2$ and the intersection form realizes $\HH$.
 The union $B$ of a tubular neighbourhood $U_1$ of $C_1$ with a tubular neighbourhood $U_2$ 
 of $C_2$ can be realized
 by attaching two disjoint $1$-handles to one initial $0$-handle, and we get  $S^1\times S^1$
 by adding one final $2$-handle; this provides a  reduced ordered handle decomposition with
 $\kappa (S^1\times S^1)=2$ $1$-handles.
 
 Now items (1) and (2) of Theorem \ref{classification2} follow from  Lemma \ref{connected-direct-sum}.
 Note that every $n\PP^2(\R)$ is not orientable (because it contains a connected curve $C$ such that $[C]\bullet [C]=1$)
 while every $n(S^1\times S^1)$ is orientable, and that all items of Corollary \ref{class-cor} hold at least if we limit to consider
 surfaces $M$, $M'$ belonging to the families of $n\PP^2(\R)$'s or $n(S^1\times S^1)$'s.  
 
 It remains to prove item (3) of  Theorem \ref{classification2}. This is the main point. 
 Thanks to the above characterization of the $2$-sphere, we can assume that $\dim \eta_1(M)>0$.
 We will follow
 the proof of the algebraic classification Theorem \ref{classification}, pointing out step by step a topological
 counterpart. We have already obtained the counterpart of $n\UU$ and $n\HH$.
 Assume first that $(\eta_1(M),\bullet_M)$
 is totally isotropic. Then every connected smooth simple curve $C\subset M$ has a product tubular neighbourhood,
 that is equivalently $[C]\bullet_M [C]=0$. Take such a curve $C$ such that $[C]\neq 0$. By the proof of Lemma     
\ref{connected-rep}, there is another connected curve $C' \subset M$ which intersects transversely $C$ at
one point (so that $[C]\bullet [C']=1$, also $[C']\neq 0$ while $[C']\bullet_M [C']=0$).
We check straightforwardly  that the union $\tilde B$ of a tubular neighbourhood $U$ of $C$ with a tubular neighbourhoos
$U'$ of $C'$  is diffeomorphic to the union $B$ of tubular neighbourhoods of the geometric generators of
$\pi_1(S^1\times S^1)$ considered above. Hence the boundary of $\tilde B$ is a connected dividing curve in $M$
and this gives rise to a connected sum decomposition 
 $$M \sim (S^1\times S^1) \# M'$$
 and we know that 
 $$\kappa (M')=\kappa(M)-2 \ . $$
 Again by Lemma   \ref{connected-direct-sum}, $(\eta_1(M'),\bullet_{M'})$ is also totally
 isotropic. Then we can conclude by induction on the dimension that in the totally isotropic case
 $$M\sim n(S^1\times S^1), \  2n=\kappa(M)=2-\chi(M) \ . $$
 
 Assume now that there is $\alpha \in \eta_1(M)$ such that $\alpha \bullet_M \alpha=1$.
 Let $C\subset M$ be a connected simple smooth representative of $\alpha$.
 Then a tubular neighbourhood $U$ of $C$ is a M\"obius band, its boundary is a dividing
 curve, we have a connected sum decomposition 
 $$M \sim \PP^2(\R)  \# M'$$
 and we know that  
 $$\kappa (M')=\kappa(M)-1 \ . $$
 By iterating the argument either we get
 $$M \sim \kappa(M) \PP^2(\R)$$
 and we have done, or
 $$M\sim k\PP^2(\R)  \# M'$$
 for some $k\geq 1$, where $\dim \eta_1(M')>0$ and $\bullet_{M'}$ is totally isotropic.
 By applying the above result in this case we eventually get
 $$M\sim k\PP^2(\R) \# h(S^1\times S^1)$$
 $$ \kappa(M)=k+2h $$
 for some $k, h\geq 1$. We conclude by applying the following final Lemma.
 Note by the way that it shows  also  that $\#$ does not verify the `cancellation property'.
 
  \begin{lemma}\label{basic-rel2} $  \PP^2(\R) \# (S^1\times S^1) \sim 3\PP^2(\R)$.
 \end{lemma}
 \Dim  First we outline a bare hands proof. After we will outline onother (but actually
 equivalent) based on a  transparent geometric construction by using the blowing up of Section \ref{BU}.
 
 {\it First proof.} Consider $S^1\times S^1$ with the geometric generators 
 $C_1$ and $C_2$ of $\pi_1(S^1\times S^1)$ transveserly intersecting 
 at the base point $(a_0,b_0)$ as above. Remove a open $2$-disk $D$ centred at $(a_0,b_0)$ and glue a M\"obius band $\Mm$
 along the boundary to get $(S^1\times S^1) \# \PP^2(\R)$. Then $(C_1\cup C_2)\setminus D$ can be completed
 by means of two fibres of the natural fibration of $\Mm$ over its core and get two disjoint simple curves $\tilde C_1$
 and $\tilde C_2$ in $(S^1\times S^1) \# \PP^2(\R)$ which intersect the core of $\Mm$ transversely at one point respectively. 
 One checks that these curves have disjoint M\"obius band tubular neighbourhoods $U_1$ and $U_2$ respectively which can be filled to give
 two copies of $\PP^2(\R)$; moreover,
$(S^1\times S^1) \# \PP^2(\R)\setminus (U_1\cup U_2)$ is connected. By filling each boundary component with a $2$-disk
we get a connected boundaryless surface $Z$ such that 
$$(S^1\times S^1) \# \PP^2(\R)\sim \PP^2(\R) \# Z \# \PP^2(\R)$$
and
$\kappa(Z)=1$ so that eventually $Z\sim \PP^2(\R)$.

\smallskip

{\it Second proof.}  Consider the product $\PP^1(\R)\times \PP^1(\R) \sim S^1\times S^1$, endowed with 
 a couple of homogeneous coordinates $(t,s)=((t_1,t_2),(s_1,s_2))$. Let $\PP^3(\R)$
 with homogeneous coordinates $x=(x_1,x_2,x_3,x_4)$. Define
 $$\psi: \PP^1(\R)\times \PP^1(\R)\to \PP^3(\R)$$
 $$\psi(t,s)=(t_1s_1, t_1s_2, t_2s_1, t_2s_2) \ . $$
 One  verifies that $\psi$ is a well defined smooth embedding onto
 the quadric $Q\subset \PP^3(\R)$ defined by the homogeneous equation $x_1x_4=x_2x_3$.
 Let $p_0=(1,0,0,0)\in Q$ and consider the ``stereographic projection'' 
 $$\phi: V\setminus \{p_0\} \to P$$
 where $P\sim \PP^2(\R)$ is the projective plane $P\subset \PP^3(\R)$ defined
 by the equation $x_1=0$. Denote by $T$ the plane tangent to $Q$ at $p_0$.
 It is defined by the equations $x_4=0$. The intersection $T\cap Q$
 consists of the union of the two lines passing through $p_0$, $l_1=\{x_4=x_2=0\}$
 and $l_2=\{x_4=x_3=0\}$. $T\cap P$ is the line $l_0=\{x_1=x_4=0\}$
 One verifies that the restriction of $\phi$ is a diffeomorphism
 $$\phi: Q\setminus (l_1\cup l_2)\to P\setminus l_0 \ . $$
 Let us blow up $\PP^3(\R)$ at the point $p_0$ and take the strict transform
 $\tilde Q$. We know from the results of Section \ref{BU} that 
 $\tilde Q \sim (S^1\times S^1) \# \PP^2(\R)$. Blow up $\PP^3(\R)$ 
 at the two points $p_1=l_1\cap P = (0,1,0,0)$ and $p_2=l_2\cap P=(0,0,1,0)$.
 Take the strict transform $\tilde P\sim 3\PP^2(\R)$. Finally one verifies
 that $\phi$ extends to a diffeomorphism
 $$\tilde \phi: \tilde Q \to \tilde P \ . $$  
 
\cvd

\smallskip

Also the proof of Theorem \ref{classification2} and of Corollary \ref{class-cor} is now complete.

\cvd

\smallskip

The above classification extends to {\it compact connected surfaces with boundary}. We limit to a few
indications. Details are left to the reader.
\smallskip

$\bullet$  Let $M$ be a compact connected smooth surface with $r\geq 1$ boundary components.
Denote by $\hat M$ the boundaryless surface obtained by filling every boundary component
with a $2$-disk. {\it Viceversa} $M$ is obtained from $\hat M$ by removing the interior of $r$ disjoint
closed $2$-disks. By the uniqueness of the disks up to isotopy, $M$ is determined up to diffeomorphism
by $r$ and the diffeomorphism type of $\hat M$.

\smallskip

$\bullet$  The radical Rad$(\bullet_M)\subset \eta_1(M)$ of the intersection form $\bullet_M$
is of dimension $r-1$ and is generated by the boundary components of $M$. The non degenerate
form $\hat \bullet _M$ uniquely induced up to isometry by $\bullet_M$ on $\eta_1(M)/{\rm Rad}(\bullet_M)$
is isometric to $\bullet _{\hat M}$. Hence $M$ is determined up to diffeomorphism by the isometry class
of the  intersection form $\bullet_M$, that is by $\dim {\rm Rad}(\bullet_M)$ and the isometry class
of $\bullet_{\hat M}$.
\smallskip

 $\bullet$ Two compact connected smooth surfaces with boundary $M$ and $M'$ are diffeomorphic
 if and only if they have the same number of boundary components, $\chi(M)=\chi(M')$, and
 either they are both orientable or non orientable.
 
 \section{$\Omega_1(X)$ as the abelianization of the fundamental group}\label{omega1-ab}
 Recall that in Proposition  \ref{tau-onto} we have established a natural epimorphism
  $$h_1:\pi_1(X,x_0)\to \Omega_1(X)$$ $X$ being a path-connected topological space.
  Now we are able to determine the  kernel of this epimorphism.
  \begin{proposition}\label{ab} The kernel $\ker h_1$ coincides with the commutator subgroup
  of $\pi_1(X,x_0)$, hence $\Omega_1(X)$ is the abelianization of the fundamental group.
  \end{proposition}
  \Dim Let $\gamma: (S^1,p)\to (X,x_0)$ be a homotopically non trivial loop which represents $0\in \Omega_1(X)$.
  Then $\gamma$ can be extended to a map $h: \Sigma \to X$ where $\Sigma$ is a compact orientable
  surface with boundary $\partial \Sigma = S^1$ such that by attaching a $2$-disk along $\partial \Sigma$, 
  we get  a boundaryless compact orientable surface $\tilde \Sigma$ of genus say $g\geq 1$.
  By using the concrete models for such a surface provided by the classification theorem, we see that there is embedded in
  $\tilde \Sigma$ a wedge of $2g$-simple loops based at $p$, not intersecting $D^2\setminus \{p\}$, 
  such that by cutting the surface along these loops we get a $4g$-gone and $\gamma$ retracts onto that wedge within $\Sigma$.
  Finally one realizes that these loops can be distribute in two family  say $a_1,\dots a_g$, $b_1,\dots , b_g$,
  in such a way that the above retraction realizes a homotopy between $\gamma$ and the composite
  loop $$a_1b_1a_1^{-1}b_1^{-1}a_2b_2a_2^{-1}b_2^{-1}\cdots a_gb_ga_g^{-1}b_g^{-1} \ . $$ The proposition follows.
  
  \cvd 
  
  \smallskip
  
  The above proposition means that every  homomorphism $\phi: \pi_1(X,x_0)\to G$ where $G$ is abelian factorizes
  as $\phi= \hat \phi \circ h_1$, $\hat \phi: \Omega_1(X)\to G$. 
 
  \section{$\Omega_2$ and $\eta_2$}\label{2-bordism}
 Here we consider again boundaryless compact surfaces.
 As a corollary of the classification we have
 
 \begin{theorem}\label{omega-eta-2}
 (1) $\Omega_2=0$;
 \smallskip
 
 (2) $\eta_2 \sim \Z/2\Z$ and is generated by $[\PP^2(\R)]$.
 \smallskip
 
 (3) $\psi: \eta_2\to \Z/2\Z, \ \phi([M]):= \chi_{(2)}(M)$
 is a well defined isomorphism.
 
 \end{theorem}
 
 \Dim Recall that $$[M_1\# M_2]=[M_1]+[M_2] \in \eta_2$$ 
 (resp. $\in \Omega_2$ in the oriented setting). 
 It follows from the classification that every compact connected oriented surface
 is the boundary of an oriented $3$-manifold (in fact $n(S^1\times S^1)$ can be embedded in $S^3=\R^3\cup \infty$ 
 and divides it). Hence $\Omega_2=0$. 
 
 On the other hand  for every compact connected surface $M$,
 $$[M \# 2\PP^2(\R)]= [M] \in \eta_2$$ 
 and 
 $$M \# 2\PP^2(\R)\sim (\kappa(M)+2)\PP^2(\R) $$
 by the classification.
 Hence
 $$ [M]= \chi_{(2)}(M) [\PP^2(\R)]\in \eta_2 \ . $$
 As $[\PP^2(\R)]\neq 0$ then items (2) and (3) follow.
 
 \cvd
 
 \smallskip
 
 \subsection{$\eta_2$ as a Witt group}\label{Witt}
 Apparently Theorem \ref{omega-eta-2} is exhaustive.
 However the topological classification of surfaces runs
 parallel to the algebraic classification on $\Z/2\Z$-symmetric
 bilinear forms up to isometry. We would like to recast also
 the content of  Theorem \ref{omega-eta-2} within this vein. 
 
Denote by $I(\Z/2\Z)$ the set of isometry classes of non degenerate
symmetric bilinear forms defined on $\Z/2\Z$-vector spaces of arbitrary
finite dimension. $I(\Z/2\Z)$ is a semigroup provided it is endowed
with the operation $\perp$. $S\in I(\Z/2\Z)$ is said {\it neutral} 
if $\dim S = 2m$ is even and there is a subspace $Z\subset S$, $\dim Z=m$
such that $Z=Z^\perp$. It follows from Theorem \ref{classification} that $S$ is neutral
if and only if either $S= 2m\UU$ or $S=m\HH$, for some $m$. Put on $I(\Z/2\Z)$ the equivalence
relation $X\sim X'$ if and only if there are neutral spaces $S,S'$ such that 
$$X\perp S = X' \perp S' \ . $$
Denote by $W(\Z/2\Z)$ the quotient set. For every $X\in I(\Z/2\Z)$, $X\perp X$ is neutral,
hence $\perp$ descends to $W(\Z/2Z)$ and makes it an abelian group called the
{\it Witt group} of the field $\Z/2\Z$; $0\in W(\Z/2\Z)$ is the class of neutral spaces, and  for every 
$[X]\in W(\Z/2\Z)$, $-[X]=[X]$. It follows from Theorem \ref{classification} that
$$ r_{(2)}: W(\Z/2\Z)\to \Z/2\Z, \  r_{(2)}([X]):=  \dim X \ {\rm mod} (2)$$
is a well defined isomorphism of groups.
Finally the content of Theorem \ref{omega-eta-2} can be rephrased as follows
\begin{theorem}\label{eta-2} $$ \wG: \eta_2 \to W(\Z/2\Z), \ \wG([M])= [ \bullet_M]$$
is a well defined isomorphism; moreover
$$ r_{(2)}\circ \wG = \chi_{(2)} \ . $$
\end{theorem}

\cvd

\subsection{A direct derivation of $\Omega_2$ and $\eta_2$}\label{direct-bord}
Theorem \ref{omega-eta-2} has been derived as a corollary of
the classification. Here we outline a direct derivation;
the mechanism is interesting: starting from a $2$-dimensional handle decomposition
of $M$ it produces a $3$-dimensional handle decomposition in such a way that the surface is the boundary; 
somehow  {\it $M$ builds its 'simplest bulk'}. 
  
For every compact connected $M$ as usual, take a reduced ordered handle decomposition with $\kappa=\kappa(M)$ $1$-handles.
Hence starting from $(M_0,\partial M_0)=(D^2,S^1)$, we have a sequence $(M_i,\partial M_i)$, $i=1,\dots, \kappa$, obtained by
attaching one $1$-handle to $(M_{i-1},\partial M_{i-1})$; finally $M$ is obtained by attaching a $2$-handle to $(M_\kappa, \partial M_\kappa)$.
Consider the product $3$-manifold $W=M\times [0,1]$. On the copy $M'=M\times \{1\}$ of $M$, consider the family of pairwise disjoint 
non necessarily connected curves $\partial M_i$, $i=0,\dots , \kappa$. There is a system of pairwise disjoint tubular neighbourhoods 
$U_i\sim \partial M_i \times [-1,1]$ of these curves in  $M'$. Let us attach to $W$ along $M'$ a family of $\kappa$ disjoint
three dimensional $2$-handles, each one attached along $U_i$, $i=0,\dots, \kappa$. In this way we get a $3$-manifold
$W'$ such that 
$$\partial W' = (M\times \{0\})\amalg M"$$
where $M"$ has $\kappa+2$ connected components, each one associated to  one of the handles of the original decomposition of $M$.
It is not hard to see that a component of $M"$ corresponding either to the $0$-handle or the $2$-handle of $M$ is diffeomorphic
to $S^2$. For a component associated to a $1$-handle there are two possibilities:
\begin{enumerate}
\item Starting from an annulus $A\sim S^1\times [0,1]$ we attach the $1$-handle along $S^1\times \{1\}$ in such a way that the
resulting surface is orientable; then this surface is a `pant' $P$ and the corresponding component of $M"$ is obtained by filling
every component of $\partial P$ with a $2$-disk, so that it is diffeomorphic to $S^2$. 

\item Starting from an annulus $A\sim S^1\times [0,1]$ we attach the $1$-handle along $S^1\times \{1\}$ in such a way that the
resulting surface is non orientable; then this surface is a M\"obius band $\Mm$ and the corresponding component of $M"$ is obtained by filling
$\partial \Mm$ with a $2$-disk, so that it is diffeomorphic to $\PP^2(\R)$.
\end{enumerate}

It follows that $M$ is bordant with the disjoint union of $k$ copies of $\PP^2(\R)$ for some $k\geq 0$. This is enough to conclude
that $\eta_2\sim \Z/2\Z$ and is generated by $[\PP^2(\R)]$. 

Assume now that $M$ is orientable.  Hence $W$ is orientable, and also $W'$ is orientable
because attaching a $2$-handle does not destroy the orientability. Also $\partial W'$ is orientable
so that $M"$ is a disjoint union of $2$-spheres. This is enough to conclude that $\Omega_2=0$. 
But we can say more. Let $W"$ be obtained from $W'$ by filling every component of $M"$ with a $3$-disk. By construction, $W"$ is obtained from $W$
by attaching a few disjoint $2$-handles followed by a few $3$-handles. By considering the dual decomposition, we see that $W$ is obtained
starting from a few $0$-handles followed by a few disjoint $1$-handles. By cancellation of $0$-handles we can assume that there is
only one $0$-handle. By sliding handles, we realize that up to diffeomorphism $W":= \Hh_h$  is uniquely determined by the number
$h$ of $1$-handles, it is called a {\it handlebody} of {\it genus} $h$, and $M=\partial \Hh_h$. By some consideration about the Euler-Poincar\'e
characteristic, one finally realizes that $\kappa(M)=2h$; in this way we have re-obtained a classification up to diffeomorphism, at least
in the orientable case.
 
\section{Stable equivalence - Rational  models ($2D$ Nash's conjecture)}\label{stable-2-Nash}
The classification of surfaces up to diffeomorphism  contains a coarse classification up to {\it stabilization}:
let us say that two (compact connected boundaryless, as usual) surfaces $M$ and $M'$ are {\it stably equivalent} if there are
$n,m\in \N$ such that 
$$M \# n\PP^2(\R)\sim M' \# m\PP^2(\R)\ . $$
Then we have as an immediate corollary of the full classification that {\it every surface is stably equivalent to each other}.
In the orientable setting we have a similar result up to stabilization by some $n(S^1\times S^1)$.

This coarse classification deserves to be pointed out because it is a sort of toy model of phenomena occurring
for example in dimension $4$ (in spite of the fact that a full classification is not known in such a case), and also because it
has a different flavour once we interpret $\# \PP^2(\R)$ as the {\it blowing up at a point}, accordingly to Section \ref{BU}.
Then a stable equivalence between $M$ and $M'$ is realized by a $\tilde M$ which dominates both being
obtained by blowing up some points of each respectively; equivalently we can say that $M'$ is obtained from $M$
by firstly blowing up some points of $M$ and then performing a certain {\it blowing down} to $M'$.

Recall (Remark \ref{Alg-BU}) that a compact  real algebraic set $X$ is {\it rational} if it {\it birationally equivalent} to the projective space
of the same dimension, say $\PP^n(\R)$; that  is $X$  contains a non empty Zariski open set which is algebraically isomorphic
to a Zariski open set in the projective space of the same dimension. If $X=B(\PP^n(\R),Y)$ is obtained from $\PP^n(\R)$ by some blowing up
along a regular algebraic centre (in particular a finite set of points), then $X$ is a rational regular algebraic set. 
A so called ``Nash's conjecture'' stated in \cite{Na}
asked if every compact smooth manifold admits up to diffeomorphism any  rational regular real algebraic model. 
We have a rather complete answer in the case of surfaces: 

\smallskip

$\bullet$  Every non orientable surface $M\sim \PP^2(\R)\# n\PP^2(\R)\sim B(\PP^2(\R),Y)$, $\kappa(M)=n+1$,
$Y$ consisting of $n$ points, has a rational model;
\smallskip

$\bullet$ If $M$ is orientable $M  \# \PP^2(\R)$ admits a rational model $B(\PP^2(\R),Y)$ where $Y$ consists of $2n=\kappa(M)$ points.
One can ask if  $Y$ as above can be chosen in such a way that a blowing down that returns $M$ can be done in the algebraic
setting, providing a rational model for $M$ itself. For example in the second proof of  Lemma \ref{basic-rel2}, we see such a mechanism
which produces $\PP^1(\R)\times \PP^1(\R)\sim S^1\times S^1$ by blowing down $B(\PP^2(\R), \{p_1,p_2\})$, collapsing to a point $p_0$
the strict trasform of the line of $\PP^2(\R)$ passing through the points $p_1$ and $p_2$. One can prove in general
that if $Y=\{p_1,\dots,p_{2n}\}$ is contained in a projective line $l\subset \PP^2(\R)$, then by blowing down to a point $p_0$ the strict transform $\tilde l$
of $l$ in $B(\PP^2(\R),Y)$ we get a rational algebraic set $X$, which is {\it homeomorphic} to $M$ via a algebraic homeomorphism which restricts to an 
algebraic isomorphism between regular Zariski open sets
$$B(\PP^2(\R),Y)\setminus \tilde l \to X\setminus \{p_0\} \ . $$ However, {\it if $n>1$ $X$  is not regular} 
as it has one isolated singularity at $p_0$.  These rational models with one isolated singularity are the best we can do
because it is known since Comessati \cite{COM} that $S^1\times S^1$ is the only {\it orientable} surface admitting a {\it regular} rational model. 

\section{Quadratic enhancement of surface intersection forms}\label{quadratic}
Let $(\eta_1(M), \bullet_M)$ be as above, where $M$ is a compact connected boundaryless surface. 
In several situations one is interested to the embeddings or immersions of $M$ in a given higher dimensional manifold,
considered up to suitable equivalence relations which often enhance the abstract surface bordism. In such situations so called 
{\it quadratic enhancements} of the intersection form naturally    
arise. In this section we will develop a few aspects  of the abstract theory of such
structures. Many proofs are  simple exercises and we will omit them.
Later in the text we will see concrete applications (see Sections \ref {small-k}, \ref {Pin-}, \ref {tear}, \ref{16}).

\medskip

Let $(V,\beta)$ be a finite dimensional $\Z/2\Z$-vector space endowed with a non degenerate symmetric bilinear form $\beta$.
\smallskip

{\it (Totally isotropic case)} Assume first that $\beta$ is totally isotropic, so that $(V,\beta)$ is isometric to $g\HH$, $\dim V= 2g$.

\begin{definition}\label{quadratic-form} {\rm a map $q: V \to \Z/2\Z$ is a {\it quadratic enhancement} of $(V,\beta)$ (sometimes we
simply say ``of $\beta$")  if for every
$x,y \in V$,
$$ q(x+y)= q(x)+q(y) + \beta(x, y)  \ . $$}
\end{definition}

\smallskip 
 
We can enhance the equivalence relation ``up to isometry" to the set of such triples: 
$$f: (V_1,\beta_1,q_1)\to (V_2,\beta_2,q_2)$$ 
is an {\it isometry} if and only if 
$$f:(V_1,\beta_1)\to (V_2,\beta_2)$$ 
is an isometry in the usual sense
and moreover, for every $x\in V_1$, $q_1(x)= q_2(f(x))$. We denote by $$I^\HH_q(\Z/2\Z)$$ the set of isometry classes of these triples. The operation ``$\perp$"
gives it a semigroup structure.
\smallskip

It is rather easy to enhance the results of Section \ref {algebraic-form} (in the totally isotropic case); as usual sometimes we
will  confuse representatives with their isometry classes:
\medskip

-  Up to isometry there are exactly two quadratic enhancement of $\HH$ (endowed with the standard hyperbolic basis say $\{e_0,e_1\}$):
\begin{itemize}
\item  $q_0(e_0)=q_0(e_1)=0, \ q_0(e_0+e_1)=1$; denote by $\HH^{0,0}$ the corresponding equipped space;
\item  $q_1(e_0)=q_1(e_1)=q_1(e_0+e_1)=1$; denote by $\HH^{1,1}$ the corresponding equipped space.

\end{itemize}

Then every triple $(V,\beta,q)$ 
is isometric to
$$ m\HH^{0,0}\perp n\HH^{1,1}$$
for some $m,n \in \N$ such that $2(m+n)=2g = \dim V$. Such integers $m$ and $n$ are not unique; in fact we have
\begin{lemma}\label{2H00=2H11} $\HH^{0,0}\perp \HH^{0,0} =  \HH^{1,1}\perp \HH^{1,1}$.
\end{lemma}

\cvd

\begin{proposition}\label{ARF}  (1) $${\rm Arf}:  ( I^\HH_q(\Z/2Z), \perp) \to (\Z/2Z,+),  \ {\rm Arf}([V,\beta,q])= n \ {\rm mod}(2)$$
provided that $[V,\beta,q] = m\HH^{0,0}\perp n\HH^{1,1}$ for some $(m,n)\in \N^2$, is a well defined surjective semigroup
homomorphism.
\smallskip

(2) ${\rm Arf}([V,\beta,q])=1$ if and only if $|q^{-1}(1)|>|q^{-1}(0)|$;  ${\rm Arf}([V,\beta,q])=0$ if and only if $|q^{-1}(1)|<|q^{-1}(0)|$.
\smallskip
 
(3) If $[V,\beta]=g\HH$ and the $j$-copy of $\HH$ is endowed with its standard hyperbolic basis $\{e^j_0,e^j_1\}$, $j=1,\dots , g$,
then 
$$ {\rm Arf}([V,\beta,q])= \sum_j q(e^j_0)q(e^j_1) \ . $$
\end{proposition}

\cvd

\smallskip

Arf is called the {\it Arf invariant}.

We can define the  {\it Witt group} associated to the semigroup $$(I^\HH_q(\Z/2\Z),\perp) \ . $$
$[V,\beta,q]\in I^\HH_q(\Z/2\Z)$, $\dim V=2g$, is said {\it neutral} 
if there is a subspace $Z\subset V$, such that $\dim Z=g$,
$Z=Z^\perp$ and $q$ vanishes on $Z$. 
Put on $I^\HH_q(\Z/2\Z)$ the equivalence
relation $X\sim X'$ if and only if there are neutral spaces $S,S'$ such that 
$$X\perp S = X' \perp S' \ . $$
Denote by $W^\HH_q(\Z/2\Z)$ the quotient set. 
The operation $\perp$ descends to $W^\HH_q(\Z/2Z)$ and makes it an abelian group. 
We have:
\begin{proposition} The Arf homomorphism passes to the quotient
$${\rm Arf}: W^\HH_q(\Z/2\Z)\to \Z/2\Z$$
and is in fact a group isomorphism. The Witt group is generated by $\HH^{1,1}$.
\end{proposition}

\cvd
\smallskip

 We know that  $(I^\HH(\Z/2\Z),\perp)$  is isomorphic to the semigroup of {\it orientable} compact connected boundaryless
surfaces (considered up to diffeomorphism) endowed with the ``$\#$" operation.  The isomorphism is given by associating
to every surface $M$ the class of $(\eta_1(M),\bullet_M)$. So the above algebraic discussion can be rephrased in such a more
topological setting. In particular the bases evoked in item (3) of Proposition \ref{ARF} can be realized geometrically:
if $M$ is a surface of genus $g$ then we can find two families of $g$ smooth circles 
$\{A_1,\dots , A_g\}$ and $\{B_1,\dots, B_g\}$ such that

- $A_i \cap A_j = B_i\cap B_j=\emptyset$ if $i\neq j$, 

- $A_i$ and $B_j$ intersect transversely at one point if and only if $i=j$, otherwise $A_i \cap B_j = \emptyset$.
\smallskip

Then these $2g$ circles form a basis of $\eta_1(M)$; if $q$ is a quadratic enhancement of $\bullet_M$, then
$${\rm Arf}(q)= \sum_j q([A_j])q([B_j]) \ . $$

$W^\HH_q(\Z/2\Z)$ can be considered as a formal non trivial refinement of $\Omega_2=0$.

\medskip

{\it (General case)} Now we consider arbitrary non degenerate
spaces $(V,\beta)$.
In this generality the notion of quadratic enhancement is subtler, due to the presence of non isotropic elements.
The key point is to consider  $\Z/4\Z$ instead of $\Z/2\Z$-valued  forms $q$.

\begin{definition}\label{q4} {\rm A map
$$ q: V \to \Z/4Z$$
is a quadratic enhancement of $\beta$ if for every $x,y \in V$,
$$q(x+y)=q(x)+q(y)+2\beta(x,y)$$
where $a\to 2a$ is the unique non trivial homomorphism $\Z/2\Z \to \Z/4\Z$.}
\end{definition} 

\medskip

\begin{remark}\label{2-4}{\rm Assume that $(V,\beta)$ is totally isotropic. If $\bar q: V \to \Z/2\Z$ is a quadratic enhancement of $\beta$
in the early sense, then $q=2\bar q$ is a quadratic enhancement in the new sense. On the other hand,
if $q: V\to \Z/4\Z$ is as in Definition \ref {q4}, then it takes only even values and there is a unique $\bar q: V\to \Z/2\Z$
such that $q=2\bar q$. So if we restrict to totally isotropic spaces we recover the previous setting.}
\end{remark}

\smallskip

The set of quadratic enhancement of $(V,\beta)$ has a structure of {\it affine space over $V$}.  That is we have

\begin{lemma}\label{V-action} There are $2^{\dim V}$ {\rm mod} (4) quadratic enhancements of $(V,\beta)$; if $q$ is one the others are of the
form
$$q'(x)=q(x)+2\beta(u,x)$$
for a unique $u\in V$.
\end{lemma}
\Dim $l(x):=2^{-1}(q'(x)-q(x))$ is linear hence represented by a unique $u\in V$ by means of the non degenerate form $\beta$.

\cvd 

The notion of isometry of triples extends {\it verbatim} and we denote by 
$$(I_q(\Z/2\Z,\Z/4\Z), \perp)$$ 
the semigroup of isometry classes. We have
\smallskip

- Up to isometry on $\HH$ there are two $\Z/4\Z$-valued quadratic enhancenents, that is 
$\qG_j = 2q_j$, $j=0,1$, where $q_j: \HH \to \Z/2\Z$ have been already defined above. We keep the
notations $\HH^{j,j}$ for the associated equipped spaces.

\smallskip

- Up to isometry, on $\UU$ there are two quadratic enhancement $q^{\pm}: \UU \to \Z/4\Z$, $q^{\pm}(1)=\pm 1$.
Denote by $\UU^\pm$ the corresponding equipped spaces.

\smallskip

Hence for every $(V,\beta)$ totally isotropic we still have 

$$ [V,\beta,q]= m\HH^{0,0} \perp n\HH^{1,1}, \ 2(m+n)=\dim V \ ; $$

If $(V,\beta)$ is not totally isotropic, then

$$[V,\beta,q] = a\UU^- \perp b\UU^+$$
for some $(a,b)\in \N^2$,  $a+b = \dim V $. As above we are not claiming that $(a,b)$ is unique.

\smallskip

In any case we say that $[V,\beta,q]$ is {\it neutral} if  there exists a subspace $Z\subset V$ such that $Z = Z^\perp$ (so that $\dim V=2h$ is even
and $\dim Z = h$) and $q$ vanishes on $Z$. As above we can define the  {\it Witt group} denoted by 
$$W_q(\Z/2\Z,\Z/4\Z)$$ 
as  a quotient of the
semigroup $(I_q(\Z/2\Z, \Z/4\Z), \perp)$.

For every $[V,\beta,q]\in I_q(\Z/2\Z, \Z/4\Z)$, for every $x\in V$, define
$$\psi_{[V,\beta,q]}(x):= \exp(\frac{i\pi}{2}q(x))=i^{q(x)} \ . $$
Finally set
$$\gamma([V,\beta,q]):= (\frac{1}{\sqrt{2}})^{\dim V} \sum_{x\in V} \psi_{[V,\beta,q]}(x)\ . $$
This is called the multiplicative {\it Brown invariant} of $[V,\beta,q]$.

For every $k\geq 2$, denote by $U_k$ the multiplicative subgroup of $U(1)$ formed by the $k$th-roots of $1$.
Denote by 
$$\alpha_k: (\Z/k\Z, +)\to U_k$$
the natural isomorphism of groups.

\begin{lemma}\label{brown=arf} If $(V,\beta)$ is totally isotropic so that $q=2\bar q$ for a unique $$\bar q: V\to \Z/2\Z$$ then 
$$\gamma([V,\beta,q])= \alpha_2({\rm Arf}([V,\beta, \bar q]) \ . $$
\end{lemma}

\cvd
\smallskip

Hence the Brown invariant extends the Arf one.  

For every $X=[V,\beta,q]$, set $-X:= [V,\beta,-q]$. We have

\begin{lemma}\label{B-property} Let $X,Y \in I_q(\Z/2\Z, \Z/4\Z)$. Then:
\smallskip

(1) $\gamma(X \perp Y)=\gamma(X)\gamma(Y)$;
\smallskip

(2) If $X$ is neutral, then $\gamma(X)=1$;
\smallskip

(3) $4X=4(-X)$.
\end{lemma}

\Dim $(1)$  follows from the very definition. 

As for $(2)$, let $X=[V,\beta,q]$, $Z\subset V$,
$\dim V =2n$, $\dim Z = n$, $Z=Z^\perp$, $q$ vanishing on $Z$. For simplicity we omit
the index $X$  in denoting $\psi$.
Let $V= Z \oplus L$ be any direct sum decomposition. Then
$$\gamma (q)= (\frac{1}{\sqrt{2}})^{2n} \sum_{z\in Z, l\in L} \psi(z+l) = (\frac{1}{\sqrt{2}})^{2n} \sum_{z\in Z, l\in L} \psi(l)(-1)^{\beta(l,z)}=$$
$$  (\frac{1}{\sqrt{2}})^{2n} [\sum_{l\in L\setminus \{0\}}(\sum_{z\in Z} (-1)^{\beta(l,z)}\psi(l))+ |Z|)] = $$
$$    (\frac{1}{\sqrt{2}})^{2n} |Z|=1$$
where the fourth equality depends on the fact that for every $l \neq 0$, $z\to \beta(l,z)$ defines a linear form $\phi$ on $Z$, and 
$\dim \ker(\phi)=\dim Z -1$ as $\beta$ is non degenerate.

As for $(3)$, it is enough to show that $4\UU^+=4\UU^-$. Let $\{e_1, e_2 , e_3, e_4\}$ be the standard basis of $4\C\sim \C^4$.
Let $\rho_j: \C \to \C^4 $ be the the linear embedding such that $\rho_j(1)=e_j$. Then one verifies that the linear isomorphism
$$\rho=(\rho_1,\dots,\rho_4):\C^4 \to \C^4$$ induces a required isomorphism
$$ \rho: 4\UU^+ \to 4\UU^- \ . $$ 

\cvd

\smallskip

Finally we can state the main result of this matter.

\begin{theorem}\label{brown-iso}  The Brown semigroup morphism $\gamma$ passes to the quotient 
and in fact it determines a group isomorphism
$$\tilde \gamma:=\alpha_8^{-1} \circ \gamma: W_q(\Z/2\Z,\Z/4\Z)\to \Z/8\Z \ . $$ 
The Witt group is generated by $\UU^+$.
\end{theorem}
\Dim $\UU^+ \perp \UU^-$ is neutral, then the Witt group is cyclic generated by $\UU^+$.
By the previous lemma, $8\UU^+$ is neutral, hence the order of $\UU^+$ divides $8$.
Finally by direct computation $\gamma(\UU^+)= \exp (\frac{i\pi}{4})$ that is it is a primitive 
fourth root of $1$.

\cvd

\smallskip

The following Corollary is easy.

\begin{corollary} The Brown invariant of $q$, the dimension of $V$ and the fact that $\beta$
is or not totally isotropic form a complete set of invariants which classifies $[V,\beta,q]\in I_q(\Z/2\Z, \Z/4\Z)$.
\end{corollary}

\cvd

\smallskip

By rephrasing everything in the topological $2$-dimensional setting, we can say that the Witt group $W_q(\Z/2\Z, \Z/4\Z)\sim \Z/8\Z$ is a formal enhancement
of the Witt group $W(\Z/2\Z) \sim \eta_2\sim \Z/2\Z$. 

We conclude this section by outlining a constructive way to build quadratic enhancements
of $(M,\bullet_M)$ for a given compact boundaryless surface $M$ (see \cite {KT}, Lemma 3.4). It is enough to define a function $q$ which associates an element in $\Z/4\Z$ to every 
disjoint union of smooth circles on $M$ (considered up to ambient isotopy) provided that the the following conditions are satisfied:
\begin{enumerate}
\item The function $q$ is {\it additive on disjoint unions}: if $L_1\amalg L_2$ is again a disjoint union of circles, then $q(L_1\amalg L_2)=q(L_1)+q(L_2)$;
\item If $K_1$ and $K_2$ are two circles that cross transversely at $r$ points, then by resolving (in one of the two possible ways) each crossing we get a disjoint union $L$
of embedded circles. Then $q(L)=q(K_1)+q(K_2)+2r$  mod$(4)$.
\item If $K$ is is a smooth circle that bounds a $2$-disk in $M$, then $q(K)=0$.
\end{enumerate}

\smallskip

In such a situation, a  quadratic enhancement  of  $(\eta_1(M),\bullet_M)$   is defined by setting $q(\alpha)=q(C)$ where $C$ is any smooth circle representing $\alpha$.

\chapter{Bordism characteristic numbers}\label{TD-ETA-CHAR}
Let us give a definition of $\eta$-characteristic number modeled on the Euler-Poincar\'e characteristic mod$(2)$, 
$\chi_{(2)}$.
As usual denote by $\Ss_n$ the class of compact boundaryless smooth $n$-manifolds.
For every $X\in \Ss_n$, let 
$$t_X: X \to \GG_{m,n}$$
be a ``truncated'' classifying map of the tangent bundle $T(X)$,
where $m=m(n)$ big enough only depends on $n$. 
An $\eta$-{\it characteristic number} is a function
$$\cG: \Ss_n \to \Z/2\Z$$
such that
\begin {enumerate}
\item It is of the form
$$\cG(X)= \cG_\alpha(X):=\sum _j  t_X^*(\alpha)\sqcap [X_j]$$
for some $\alpha \in \eta^n(\GG_{m,n})$,
where $X_j$ varies among the connected components of $X$. 
Clearly such a $\cG(X)$ is a diffeomorphism
invariant.

\item If $[X]=0\in \eta_n$, then $\cG(X)=0$. It follows that
$\cG$ induces a $\Z/2\Z$-linear map
$$ \cG: \eta_n \to \Z/2\Z \ . $$
\end{enumerate}
\medskip

Here is another characteristic $\eta$-number besides $\chi_{(2)}$. 
For every $X$, consider the $n$th-power (with respect to the $\sqcup$ product) 
$$ w^1(X)^n $$
of the Euler class of the determinant line bundle of $X$.

\begin{proposition} $\cG_{w^1(X)^n}$ is a $\eta$-characteristic number,
different from $\chi_{(2)}$.
\end{proposition}
\Dim To see that it is characteristic, it is enough to show that if $X=\partial W$ is a  boundary, then
$\cG_{w^1(X)^n}(X)=0$. Note that 
$$w^1(X) = j^*w^1(W)\in \eta^1(W,\partial W)$$
where $j:\partial W \to W$ is the inclusion. Then $w^1(X)^n= (j^*(w^1(W)))^n$,
and $w^1(X)^n$ is represented by the boundary of the proper $1$-dimensional
submanifold of $(W,\partial W)$ which represents $w^1(W)^n\in \eta^{n}(W,\partial W)$,
hence it consists of an even number of points.
To see that it is different from $\chi_{(2)}$, consider for example
$w^1(\PP^4(\R))^4=1$ while we can show (do it by exercise)
that $w^1(\PP^2(\R)\times \PP^2(\R))^4=0$.
We know that both characteristic mod$(2)$ are equal to $1$.
Hence $[\PP^4(\R)]$ and $[\PP^2(\R)\times \PP^2(\R)]$ are non trivial independent elements
of $\eta_4$. Similarly $w^1(*)^n$ distinguishes $[\PP^4(\R)]$ from $[\PP^2(\C)]$. 
The same argument extends to any couple $\PP^{a+b}(\R)$, $\PP^a(\R)\times \PP^b(\R)$
(hence to $\eta_{a+b}$)
where both $a$ and $b$ are even.

\cvd
\smallskip
\section{Stable $\eta$-numbers}\label{stable}
It is not so easy  to check directly if a function of the form $\cG_\alpha$ as
above is a characteristic number or not (that is if it vanishes on boundaries). On the other hand, this becomes almost
immediate if we consider so called ``stable classes'' in the grassmannian cobordism.
Consider the ``stabilized tautological bundle'' 
$$\tau_{m,n}\oplus \epsilon^1 \  ; $$
this corresponds to an evident classifying map 
$$s_{n}: \GG_{m,n}\to \GG_{m+1,n+1} \ . $$
Then $\alpha \in \eta^k(\GG_{m,n})$
(not necessarily $k=n$)
is by definition a {\it stable class} if
 $$ \alpha = s_{n}^*(\tilde \alpha) $$
for some  $\tilde \alpha \in \eta^k(\GG_{m+1,n+1})$.
The sum and the product of stable classes are stable.
A class of the form $\alpha=(s_{n+j}\circ \dots \circ s_n)^*(\tilde \alpha)$
is stable for every $j\geq 0$.

For every $X\in \Ss_n$, the classifying map of the {\it stable tangent bundle}
$$T(X)\oplus \epsilon^1$$  is the composite
map 
$$s_X:=s_{n} \circ t_X  \  . $$
We have

\begin{proposition}\label{stable-number}  For every $n\geq 0$, if $\alpha \in \eta^n(\GG_{m,n})$ is a stable class
then $\cG_\alpha$ is a ({\rm stable} by definition) $\eta$-characteristic number defined on $\Ss_n$.
\end{proposition}
\Dim  Assume that $X=\partial W$; then 
$$j^*(T(W))=T(X)\oplus \epsilon^1$$
so that $s_X = s_W \circ j$, where $j$is the inclusion.
It follows that 
$$t_X^*(\alpha)= j^*(t_W^*(\alpha))$$ 
hence $t_X^*(\alpha)$  is represented by the boundary of a singular proper $1$-submanifold of $(W,\partial W)$
which represents $t_W^*( \alpha)$.

\cvd

\smallskip
{\bf A construction of  stable characteristic classes.}  
In general, if $\alpha\in \eta^k(\GG_{m,n})$
is a stable class, then $t^*_X(\alpha)$ is called a {\it stable characteristic class} of $X$.
This can be extended by dealing with the classifying map of arbitrary vector bundles $\xi$ on $X$ 
and leads to the notion of stable characteristic classes of $\xi$. For simplicity we will assume that
$X$ is connected
\smallskip

For every {\it line bundle} $\xi$ on a $X$,  define the {\it total basic cobordism class}  
$$w(\xi)=  \sum_{j=0}^n  w^1(\xi)^j \in \eta^\bullet(X) $$
where we stipulate that $w(\xi)^0:= [X]$.
If we have the direct sum $\xi=\xi_1 \oplus \xi_2$ of two line bundles
set its total cobordism class 
$$ w(\xi_1\oplus \xi_2):= w(\xi_1)\sqcup w(\xi_2) \in \eta^\bullet(X)$$
and define $w^j(\xi_1\oplus \xi_2)\in \eta^j(X)$, $j=0,\dots, n$, the $j$th-homogeneous
term of $w(\xi_1\oplus \xi_2)$. This can be inductively extended to
every direct sum of line bundles on $X$, $\xi=\xi_1\oplus \dots \oplus \xi_r$, $r\leq \dim X$.
As $w(\epsilon^h)= [X]$, we see that all classes defined so far are stable classes of $\xi$.

\begin{remark}{\rm  The stable classes defined sofar might depend {\it a priori} on the given splitting of $\xi$ as 
direct sum of line bundles. It is a non trivial fact that they do not. This is part of the construction
of the so called {\it Stiefel Whitney classes} of vector bundles (see \cite{MS}) which we will not develop.}
\end{remark}

\smallskip 

For every $\PP^{a}(\R)$, for every $n>0$,
denote by $\beta$ the bundle of rank $n$ on  
$\prod_{j=1}^n \PP^a(\R)$ given by the product of $n$
copies of the tautological line bundle over $\PP^a(\R)$.
Then $\beta$ is a direct sum of $n$
line bundles. Assume that $m$ is big enough so that
we have a truncated classifying map of $\beta$
$$h_\beta: \prod_{j=1}^n \PP^a(\R) \to \GG_{m,n} \ ;$$
then for every $w^j(\beta)$ defined as above with respect to the given splitting, 
every $\alpha \in \eta^j(\GG_{m,n})$ such that
$$w^j(\beta) = h_\beta^*(\alpha)$$
is a stable class. For every direct sum of line bundles on some $Y$,
of the form $g^*(\beta)$, then $g^*(w^j(\beta))\in \eta^j(Y)$ is stable.
If $\xi$ is a vector bundle on $X$ and,
referring to Proposition \ref{flag} and adopting those notations,
$f_\xi:F(\xi)\to X$ such that $f_\xi^*(\xi)=g^*(\beta)$ ($Y=F(\xi)$) splits as a sum of line bundles,
then every class $\alpha \in \eta^j(X)$ (if any) such that $f_\xi^*(\alpha)=g^*(w^j(\beta))$
is stable.

\begin{remark}{\rm It is not evident that the construction outlined above leads to
non trivial stable classes. Actual non triviality again is part of the construction
of Stiefel-Whitney classes that we will not develop here.}
\end{remark}
 
\section{Completeness of stable $\eta$-numbers}\label{stable-complete}
This ``completeness'' refers to the fact that the necessary condition
to be a boundary stated in Proposition \ref{stable-number}
is also sufficient. 
This is an important theorem due to R. Thom \cite{T}.
The original proof is an application of the {\it Pontryagin-Thom construction} 
that allows to rephrase  the study of the cobordism ring $\eta^\bullet$
in terms of the homotopy theory of certain so called {\it Thom's spaces} (see Chapter \ref{TD-PT}).
Here we propose an elementary proof extracted from \cite{BH} which ultimately
uses only transversality. Let us state this theorem.

\begin{theorem}\label{stable-number3} $[X]=0 \in \eta_n$ if and only if 
every stable $\eta$-characteristic number vanishes on $X$.
\end{theorem}

\medskip

 It is enough to show the ``if''  implication. This will be an immediate consequence
of  the next two lemmas.

By the classification of compact $1$-manifolds, if $n=0$ then $X$ is a boundary if and only if it consists of an {\it even} number of points, 
thus it is easy to check that Theorem \ref{stable-number3} holds true for $n=0$.  If $\dim X >0$,
there is a {\it special case} such that
the stable characteristic numbers clearly vanish, that is when  $(X,s_X)$ is 
bordant with a costant map $(N,c)$; in other words $[X,s_X]$
belongs to a copy of $\eta_n$ embedded in $\eta_n(\GG_{m+1,n+1})$.
First let us prove that $X$ is a boundary under such a stroger  hypothesis.

\begin{lemma}\label{special-case} Let $\dim X >0$ and $F: Q \to \GG_{m+1,n+1}$ 
realize a bordism of $(X,s_X)$ with a constant map $c:N\to \GG_{m+1,n+1}$.
Then $N$ (hence $X$) is a boundary.
\end{lemma}

\Dim
The map $F$ pulls back the tautological bundle over the grassmannian to a rank  
$(n + 1)$ vector bundle $\xi$ on $Q$ which restricts to
$\tau_X:=T(X)\oplus \epsilon^1$ on $X$ and to a trivial bundle $\epsilon^{n+1}$ on $N$.  
Denote by $D(\xi)$, $S(\xi)=\partial D(\xi)$, the total spaces of the unitary $(n+1)$-disk and $n$-sphere 
bundles of $\xi$ respectively.
Similarly denote the restrictions $D(\tau_X)$, $S(\tau_X)$ and $D(\epsilon^{n+1})$, $S(\epsilon^{n+1})$.
Let $\iota$ be the fibrewise antipodal involution on $\xi$. 
Then $S(\xi)$ is a compact $(2n+1)$-manifold with boundary 
$$\partial S(\xi)=S(\tau_X) \amalg S(\epsilon^{n+1})$$ 
equipped with the involution
$\iota_S$ (the restriction of $\iota$). Consider the $(2n+1)$-manifold with boundary
$$Y= X\times X \times [-1,1]$$
equipped with the involution
$$\sigma (x,y,t)=(y,x,-t) $$
so that $\partial Y$ is an invariant set of $\sigma$.
The fixed point set of $\sigma$ is given by 
$$\Xx= \Delta_X \times \{0\} =\{(x,x,0)\} \subset Y$$ 
which can be naturally identified with $X$ itself.
We can find a tubular neighbourhood $U$ of $\Xx$ in $Y$ such that 
by removing the interior of $U$ from $Y$ we get a compact $(2n+1)$-manifold say $Z$,
with boundary 
$$\partial Z = \partial U \amalg \partial Y$$  
such that $(Z,\partial Z)$ is invariant for $\sigma$
and the restriction of $\sigma$ to $\partial U$ can be identified with the restriction of $\iota_S$ to $S(\tau_X)$.
Then we can glue $Z$ and $S(\xi)$ along $\partial U \sim S(\tau_X)$ and get a compact $(2n+1)$-manifold $W$ with boundary
$$\partial W =  \partial Y \amalg S(\epsilon^{n+1})$$
equipped with a smooth fixed point free  involution say $\sigma_W$, which coincides with  $\sigma \amalg \iota_S$ on
$\partial W$. Then  the quotient space $\Ww:=W/\sigma_W$ is a smooth manifold with boundary such that the quotient map
$$ q: W \to \Ww$$ 
is a degree $2$ smooth covering map. We note that the restriction of $q$ to $\partial Y$ is a trivial
covering, while 
$$S(\epsilon^{n+1})/\sigma_W \sim N\times \PP^n(\R)$$ and the restriction of $q$ to $S(\epsilon^{n+1})\sim N\times S^n$
may be identified with the map $ {\rm Id}_N\times s$, $s:S^n \to \PP^n(\R)$ being the standard double covering.
The associated real line bundle on $\Ww$ (see Chapter \ref {TD-LINE-BUND}) is the pull back by a classifying map
$$\phi: \Ww \to \PP^a(\R)$$ for some $a$ big enough, considered up to homotopy.
By the above remark about the restriction of the covering to $\partial \Ww$, we can assume that $\phi_{|\partial Y}$ is a 
constant map, while $\phi_{|N\times \PP^n(\R)}$
is the composition of the projection $N\times \PP^n(\R) \to \PP^n(\R)$ followed by the inclusion $\PP^n(\R) \subset \PP^a(\R)$.
Let $\PP^{a-n}(\R)$ be a projective subspace of $\PP^a(\R)$ which intersects $\PP^n(\R)$ transversely at one point $x_0$. 
We can also assume that $\phi(\partial Y)\cap \PP^{a-n}(\R)=\emptyset$, so that $\phi_{|\partial \Ww}\pitchfork \PP^{a-n}(\R)$
and 
$$\phi_{|\partial Y}^{-1}(\PP^{a-n}(\R))= N\times \{x_0\}\sim N \ . $$
By using usual transversality theorems, finally we can also assume that
the whole map $\phi$ is transverse to $\PP^{a-n}(\R)$ so that the proper $(n+1)$-submanifold
$(R,\partial R)$ of $(\Ww,\partial \Ww)$ given by $R=\phi^{-1}(\PP^{a-n}(\R))$ 
is such that $N\times \{x_0\} = \partial R$. This achieves the proof of  Lemma \ref{special-case}.

\cvd

\medskip

As Theorem \ref{stable-number3} holds true for $n=0$, 
we will argue by induction on the dimension $n\geq 0$.
The inductive step is provided by the following lemma combined with Lemma \ref{special-case}.
 
 \begin{lemma}\label{Induction} Let  $\dim X=n>0$.
 Assume that all stable $\eta$-characteristic numbers of $X$
 vanish, and that Theorem \ref{stable-number3}
 holds true for all dimensions $m$ smaller than $ n$. Then $(X,s_X)$ is bordant
 with a constant map $c: N\to \GG_{m+1,n+1}$.
 \end{lemma}
 \Dim  This proof will be somewhat scketchy an definitely not self-contained
 within the content of this text. Let us given a triangulation $\Kk$
 of $\GG_{m+1,n+1}$ made by smoothly embedded simplices, whose existence has been
 evoked in Section \ref {simp} (without a proof). The interior of every such a $h$-simplex is a submanifold of
 $\GG_{m+1,n+1}$ diffeomorphic to $\R^h$ and called a (open) $h$-cell of $\Kk$.
 Alternatively one can use the open cells of the natural cellular decomposition of the Grassmannian
 depicted in Section   \ref{cell-decomp}. For every $0\leq h \leq \dim \GG_{m+1,n+1}$, the union of the 
 cells of dimension less or equal to $h$
 is  called the $h$-{\it skeleton} $\Kk_h$ of $\Kk$.
 Fix a base point $x_e$ in every open cells $e$, call it the ``centre'' of the cell. 
 For every $h$ as above, by removing the centre from every $h$-cell, we get 
 a subspace $\tilde \Kk_h$ of  $\Kk_h$
 which retracts to $\Kk_{h-1}$. By basic transversality, we may assume that the smooth map $s_X$ 
 misses the centre of every cell of dimension greater than $n= \dim X$; hence, up to
 (continuous) homotopy we can assume that
 $$s_X: X \to  \GG_{m+1,n+1}$$ 
 {\it is continuos with values in 
 the $n$-skeleton $\Kk_n$, it is smooth on $s_X^{-1}(\Kk_n \setminus \Kk_{n-1})$
 and is transverse to the centre $x_e$ of every $n$-cell $e$.}
 \smallskip

 We claim  that for every $n$-cell $e$, the $0$-submanifold $Y:=s_X^{-1}(x_e)$ of $X$ consists of an {\it even}
 number of points, that is it is a $0$-dimensional boundary. 
 In fact, by collapsing $\Kk_n\setminus \{e\}$ to one point,
 we get a projection 
 $$p_e: \Kk_n \to S^n$$ 
 which restricts to a smooth embedding of the $n$-cell $e$ onto $\R^n \subset \R^n \cup \infty = S^n$,
 so that we will confuse $x_e$ with $p_e(x_e)$. 
 Then
 $$Y= (p_e\circ s_X)^{-1}(x_e)$$
 and one easily realizes that 
 $$[Y]= s_X^*(p_e^*([x_e])\in \eta^n(X)$$ 
 which vanishes as it is a stable $\eta$-characteristic number of $X$.
 Fix a small $n$-ball $D$ around $x_e$ in $e$. Then
 $$s_X^{-1}(D)= (\tilde D_1\cup \tilde D_2)\cup \dots \cup (\tilde D_s \cup \tilde D_{s+1})$$
 and the restriction of $s_X$ to every $\tilde D_j$ is a diffeomorphism onto $D$.
 Remove from $X$ the interior of every $\tilde D_j$ and pairwise glue together
 the boundary components $\partial \tilde D_j$ and $\partial \tilde D_{j+1}$, $j=1,\dots, s$
 by means of the above identifications with $\partial D$.
 Do it simultaneously at the centre of every $n$-cell. Then we get
 a boundaryless $n$-manifold $N_1$ such that the map $s_X$       
 descends to a stable classifying map 
 $$s_1=s_{N_1}:N_1 \to \Kk_n\subset \GG_{m+1,n+1}$$ which misses the centres 
 of every $n$-cell, hence up to homotopy we may assume that
 {\it $s_1$ takes values in $\Kk_{n-1}$, it is smooth on $s_{1}^{-1}(\Kk_{n-1}\setminus \Kk_{n-2})$
 and is transverse to the centre of every $(n-1)$-cell.
 Moreover, it is not hard to check that by construction $(X,s_X)$ is bordant with $(N_1,s_{1})$, so that
 also all stable $\eta$-characteristic numbers of $N_1$ vanish.}
 
 \smallskip
 
 Now we would proceed by induction on the codimension  of the skeleton
 to eventually reach  $(N_n,s_{n})$ which takes values in $\Kk_0$ and is bordant
 with $(N_{n-1}, s_{n-1})$ (hence with the initial $(X,s_X)$). As the grassmannian is connected, $(N_n,s_{n})$ 
 will be homotopic the a required constant map $(N,c)$, $N=N_n$ (this last step is not necessary
 if we use the natural cellular decomposition which has only one $0$-cell).
 
 So let us assume inductively that for some $h \geq 1$ we have obtained
 $$s_h=s_{N_h}: N_h\to \Kk_{n-h}\subset \GG_{m+1,n+1}$$ bordant with $(X,s_X)$, 
 which is smooth on $s_{h}^{-1}(\Kk_{n-h}\setminus \Kk_{n-h-1})$ and transverse
 to the centre $x_e$ of every $(n-h)$-cell $e$. The stable $\eta$-characteristic numbers of $N_h$
 vanish. By a similar  augument as above, 
for every such a cell $e$,
there is a collapsing projection
$$p_e: \Kk_{n-h}\to S^{n-h}$$
which restricts to a smooth embedding of the cell $e$ onto $\R^{n-h} \subset \R^{n-h} \cup \infty = S^{n-h}$;
by confusing $x_e $ with $p_e(x_e)$, set $Y=(p_e \circ s_h)^{-1}(x_e)$.   
We claim that
this $h$-submanifold  $Y$ of $N_h$ is a boundary. By the inductive assumption 
of Lemma \ref{Induction}, it is enough to show that every stable $\eta$-characteristic
number  of $Y$ vanishes. We note that, by using the terminology defined in Chapter \ref {TD-PT}, $Y$ is {\it framed} 
that is it has a trivialized tubular neighbourhood
$U\sim Y\times D^{n-h}$ in $N_h$ such that the restriction of $s_h$ to $U$ can be identified
with the projection $Y\times D^{n-h} \to D^{n-h}$, where $D^{n-h}$ is a small disk in $e$
around $x_e$. This implies that a stable classifying map $s_{Y}$ for $Y$
is given by $s_h\circ j$, where $j:Y \to N_h$ is the inclusion.
Then it is enough to show that  for every $\alpha  \in \eta^h(\GG_{m+1,n+1})$,
$s_Y^*(\alpha)\sqcap [Y] \in \Z/2\Z$ vanishes. By the geometric definition of the 
cobordism products, we realize
that as an element of $\Z/2\Z$, $s_Y^*(\alpha)\sqcap [Y]$ equals  $s^*_h(p_e^*[x_e]\sqcup \alpha)\sqcap [N_h]$
which vanishes being a stable $\eta$-characteristic number of $N_h$.
Then $Y$ is a boundary of a manifold $W$.
We make a surgery on $N_h$ by replacing the above product neighbourhood $U\sim Y \times D^{n-h}$ 
with $W\times \partial D^{n-h}$; do it simultaneously at every $(n-h)$ cell. we get a manifold
$N_{h+1}$; the map $s_{h}$ descends to $s_{h+1}: N_{h+1}\to \GG_{m+1,n+1}$ which can be identified
with the projection   $W\times \partial D^{n-h} \to \partial D^{n-h}$ at every $(n-h)$ cell.
By construction $(N_{h+1},s_{h+1})$ is bordant with $(N_h,s_h)$ and this 
eventually achieves the inductive step.

The proofs of Lemma \ref{Induction} and of Theorem \ref{stable-number3} are now complete.

\cvd

\section{On $\Omega$-characteristic numbers}\label{omega-char}
Recall that a compact boundaryless $n$-manifold $X$ is parallelizable
if the tangent bundle admits a global trivialization so that its total space is
diffeomorphic to $X\times \R^n$; in such a case $X$ is orientable.
If $X$ is parallelizable then any classifying map $t_X: X \to \GG_{m,n}$
of $T(X)$  is homotopic to a constant map as well as any stable
classifying map $s_X: X \to \GG_{m+1,n+1}$. Then if $X$ is parallelizable and $\dim X =n>0$
certainly it verifies the hypothesis of Lemma \ref{special-case}, hence 
$[X]=0 \in \eta_n$. We can strenghten this result.

\begin{proposition}\label{parallelizable} Let $X$ be a parallelizable and oriented
compact boundaryless $n$-manifold, $n>0$. Then $[X]=0\in \Omega_n$.
\end{proposition}
\Dim  It is enough to prove the statement when $X$ is connected.
We will use and refine the proof of Lemma \ref{special-case}.
If $\dim X =n$ is {\it even}, we can apply such a proof starting from
a homotopy $F: X\times [0,1]\to \GG_{m+1,n+1}$ between 
$s_X$ and a constant map. Clearly $X\times [0,1]$ is orientable.
At the end of the proof we may assume that both $\PP^a(\R)$
and $\PP^{a-n}(\R)$ are odd dimensional, hence they are both orientable.
Then we conclude by means of the oriented version of the transversality
theorems.  

If $\dim X =n$ is {\it odd} we modify the construction as follows:
we consider $$Y=X\times X$$ endowed with the involution
$\sigma(x,y)=(y,x)$.  The fixed point set consists of the diagonal
$\Delta_X$ which is naturally identified with $X$ itself.
A tubular neighbourhood $U$ of $\Delta_X$ can be identified
with the unitary disk bundle of $T(X)$, hence with the product
$X\times D^n$. By removing the interior of $U$ from $Y$, we get
a compact $2n$-manifold $W$ with boundary $\partial W = X\times S^{n-1}$;
$\sigma$ restricts to a fixed point free involution on $W$, and can be identified 
with the  fibrewise antipodal map on $\partial W$, that is the trivial unitary
sphere bundle of $T(X)$. Then the proof runs similarly to the one of
Lemma \ref{special-case}. At the end we can assume that both
$\PP^a(\R)$ and $\PP^{a-n+1}(\R)$ are orientable and conclude
again by oriented transversality.

\cvd

\medskip

Every $\eta$-characateric number lifts to an $\Omega$-characteristic number
(with the obvious meaning of the term)
via the forgetting projection $$\Omega_\bullet \to \eta_\bullet \ . $$
If the manifold $X$ is oriented we can consider also the 
{\it complexification} $T_\C(X)$  of the tangent bundle: every real vector
bundle $\xi$ can be complexified to $\xi_\C$ via the inclusion $\R\subset \C$
so that  every
real cocycle defining $\xi$ can be considered as a cocycle defining $\xi_\C$.
Then $T_\C(X)$ corresponds to a classyfying map
$$t_{X,\C}: X\to \GG_{m,n}(\C) \ . $$
We can apply almost verbatim the above discussion about (stable)
characteristic numbers in the present complexified setting
(by replacing in particular real with complex line bundles).
This gives rise to further $\Omega$-characteristic numbers with values in 
$\Z$ instead of $\Z/2\Z$. We call generically {\it stable $\Omega$-characteristic
number} one belonging to the union of such two families.

\begin{remark}\label{sw}{\rm The classical treatment of stable characteristic
numbers (classes) takes places in the {\it singular cohomology ring} of
real or complex grassmannians with $\Z/2\Z$ or $\Z$ coefficients (see \cite{MS}, \cite{BT}); 
they are called {\it Stiefel-Whitney} and {\it Pontryagin} numbers
respectively. As we do not assume any familiarity with singular
cohomology, above we have just `lifted' some facts of such a
theory in terms of the cobordism rings that we have
introduced in a self-contained way. 
In the case of $\eta$, ``to lift'' is quite appropriate
because one can prove (it is non trivial) that for every 
compact boundaryless $n$-manifold $X$, the $\Z/2\Z$-cohomology
$H^j(X;\Z/2\Z)$ can be different from $0$ only if $0\leq j \leq n$,
is finite dimensional and coincides with the quotient space $\Hh^j(X;\Z/2\Z):=\eta^j(X)/\ker (\phi^j)$.
Hence stable $\eta$-characteristic numbers and Stiefel-Whitney numbers
are basically the same. In the oriented case, $\Omega$-characteristic numbers
are not exhaustive.  Our presentation of the matter is necessarily incomplete.
}
\end{remark}  

By using  Stiefel-Whitney and Pontryagin numbers,
we have the following oriented version of Theorem \ref{stable-number3}.
The  proof \cite{Wall} is more complicated. 
Parallelizable manifolds as in Proposition \ref{parallelizable}   represent the basic instance
for this theorem.

\begin{theorem}\label{omega-stable} Let $X$ be a compact oriented boundaryless $n$-manifold.
Then $[X]=0 \in \Omega_n$ if and only if 
all Stiefel-Whitney and Pontryagin numbers  of $X$ vanish.
\end{theorem}

\cvd

\chapter{The Pontryagin-Thom construction}\label{TD-PT}
The original Pontryagin construction was inventend to rephrase the
study of the homotopy groups of spheres in terms of a certain more
geometric (hence presumably more accessible at that time, about 1938)
bordism theory.  Viceversa, later Thom's extension of Pontryagin
construction was mainly intended as a way to rephrase the study of
$\eta^\bullet$ (or $\Omega^\bullet$) in terms of the homotopy groups
(becomed more accessible at that time, about 1954, after Serre's
Thesis) of certain so called Thom's spaces which in a sense generalize
the spheres.  So the P-T construction is a powerfull bridge between
two different ways to approach a same ``mathematical reality''.
\medskip

Let us start by describing the Pontryagin construction
(introduced in 1938; see the later 
exposition in \cite{Pont}, and also \cite{M1}) . 
We are primarily interested here in the determination of the homotopy groups 
$$\pi_m(S^n, p)$$
for $m,n \geq 1$. By suitable approximation theorems, 
we know that we can manage with them in purely
differential/topological way. We know that 
$$\pi_1(S^1,p)\sim \Z,  \ \pi_m(S^1,p)=1 \ {\rm for} \  m>1 , \ $$
$$\pi_m(S^n,p)=1 \ {\rm for}\  n\geq 2, \ 1\leq m <n \ . $$
Hence we will assume that $m\geq n > 1$.
In such a case $\pi_m(S^n,p)$ is abelian, the base point is immaterial
and the group can be identified with $[S^m,S^n]$,  the set of smooth homotopy
classes of maps $f:S^m\to S^n$. Moreover, it is convenient to extend
the discussion to $[M,S^n]$ where  $M$ is any {\it compact, connected 
boundaryless} smooth $m$-manifold, $m\geq n \geq 1$.

\section {Embedded and framed bordism}\label{emb-fram-bord}
We have already encountered instances of embedded bordism within a
given manifold in Chapter \ref {TD-LINE-BUND}. Let us state it in
general.

\begin{definition}\label{emb-bord}
{\rm Let $M$ be a compact connected boundaryless $m$-manifold.  Let
  $0\leq k < m$. Given compact boundaryless smooth $k$-submanifolds
  $V_0$, $V_1$ of $M$, we say that {\it $V_0$ is bordant with $V_1$
    within $M$} (and we write $V_0\sim_{b,M} V_1$) if there is a
  smooth triad $(W,V_0,V_1)$, properly embedded into $M\times
  [a_0,a_1]$, for some $a_0<a_1$, such that for $j=0,1$,
$$\partial W \cap (M\times \{a_j\})= V_j \ . $$
}
\end{definition}
\smallskip

The relation ``$\sim_{b,M}$" is an equivalence relation on the set of
compact boundaryless $k$-submanifolds of $M$: every such a $V$ is in
relation with itself because the cylinder $V\times [a_0,a_1]$ properly
embeds into $M\times [a_0,a_1]$; the relation is obviously symmetric;
as for the transitivity, up to isotopy we can normalize the proper
embeddings of the triads $(W,V_0,V_1)$ in such a way that they are
locally cylinder-like as above near the boundary. Given properly
embedded triads $(W,V_0,V_1)$ in $ M\times [a_0,a_1]$,
$(W',V'_0,V'_1)$ in $M\times [a'_0, a'_1]$ respectively, such that
$V_1=V'_0$, then we can construct $(W'',V_0, V'_1)$ in $M\times [a_0,
  a_1+a'_1-a'_0]$ just by stacking $M\times [a'_0,a'_1]$ over $M\times
[a_0,a_1]$.
\smallskip

We denote by 
$$\eta^{\rm emb}_k(M)$$ 
the quotient set.

\smallskip

By restriction to {\it oriented} $k$-submanifolds of $M$, we can get
an oriented version of the above definition leading to a quotient set
$$ \Omega^{\rm emb}_k(M) \ . $$
Stress that we are not assuming that $M$ is oriented.

\smallskip

Let $M$ be as above.
\begin{definition}\label{framing}{\rm
A compact boundaryless $k$-submanifold $V\subset M$ is {\it framed} if
it is endowed with a {\it framing}. This last is of the form
$$\fG = (s_1,\dots,s_{m-k})$$
where
\begin{enumerate}
\item Every $s_j$ is a nowhere vanishing section of the bundle
  $i_V^*T(M)$, $$i_V: V \to M$$ being the inclusion;
\item For every $x\in V$, the vector $s_1(x),\dots, s_{m-k}(x)$ are linearly independent in $T_xM$;
\item For every $x\in V$, $T_xM = T_xV \oplus F_x$, where $F_x:= {\rm Span}\{s_1(x),\dots, s_{m-k}(x)\}$.
\end{enumerate}
}
\end{definition}
\smallskip

Hence $x\to F_x$ defines a smooth field of transverse $(m-k)$-planes
along $V$ tangent to $M$. The framing provides a global trivialization
of the bundle $i_V^*T(M)$, hence of every tubular neighbourhood of $V$
in $M$ constructed by means of such a field.  This means in particular
that a necessary (and sufficient) condition in order that $V$ admits a
framing is that it has globally trivializable tubular neigbourhoods in
$M$.
\smallskip

We are going to specialize and enhance the embedded bordism to framed
submanifolds.  First let us extend the definition of framing to
properly embedded triads. Let $(W,V_0,V_1)$ be a properly embedded
$(k+1)$-triad in $M\times [a_0,a_1]$; from now on we will assume by
default that the embedding is normalized, i.e. cylinder-like near the
boundary as above.  A framing of the triad in $M\times [a_0,a_1]$ is
of the form
$$\fG = (s_1,\dots,s_{m-k})$$ where these are pointwise linearly
independent sections of the bundle $$i^*_W T(M\times [a_0,a_1]) \ , $$
induce a smooth field of transverse $(m-k)$-planes along $W$ tangent
to $M\times [a_0,a_1]$, and we require furthermore that the
restriction of $\fG$ to the boundary defines a framing of $V_j$ in
$M$, $j=0,1$.

\begin{definition}\label{framed-eta} {\rm
    Let $(V_0,\fG_0)$ and $(V_1,\fG_1)$ framed $k$-submanifolds of
    $M$.  We say that {\it $(V_0,\fG_0)$ is framed bordant with
      $(V_1,\fG_1)$ within $M$ } and we write $$(V_0,\fG_0)\sim_{fb}
    (V_1,\fG_1) \ , $$ if there is a properly embedded framed triad
    $((W,V_0,V_1), \fG_W)$ in some $M\times [a_0,a_1]$ such that the
    restriction of the framing $\fG_W$ to the boundary coincides with
    the union of the framings $\fG_0$ and $\fG_1$.  }
\end{definition} 
\smallskip

Similarly as above, one cheks that this defines an equivalence
relation on the set of framed $k$-submanifolds of $M$, and we denote
by
$$ \eta^\Ff_k(M)$$
the quotient set.

We develop now also an {\it oriented} version; we stress that we do it
{\it provided that $M$ itself is oriented}.  So not only we require
that for every framed $k$-submanifolds $(V,\fG)$, $V$ is also
oriented; furthermore we impose that for every $x\in V$, the given
orientation of $T_xV$ followed by the tranverse orientation of $F_x$
determined by $\fG(x)$ coincides with the given orientation of
$T_xM$. Hence it is enough to require that $V$ is {\it orientable}
because the framing together with the orientation of $M$ select the
preferred orientation of $V$. Note that this way recovers the
orientation procedure stated in Theorem \ref{firstT}.  To define the
pertinent relation ``$\sim_{fob}$" we use {\it oriented} and framed
triads $((W,V_0,V_1),\fG_W)$ properly embedded in some $M\times
[a_0,a_1]$, so that the oriented boundary $\partial W = V_0\amalg
-V_1$. This leads to the quotient set
$$ \Omega^\Ff_k(M) \ . $$      

\section{The Pontryagin map}\label{pont-map}
Let us keep the above setting. We establish the following
procedure.
\smallskip

$\bullet$ Fix $x_0 \in S^n$. For every $\alpha \in [M,S^n]$, thanks to
transversality take $f: M \to S^n$ belonging to $\alpha$ and such that
$f\pitchfork \{x_0\}$.

\smallskip

$\bullet$ $V:=f^{-1}(x_0)$ is submanifold of $M$ of dimension $\dim V=
k:=m-n$. Fix a positive basis $\Bb$ of $T_{x_0}S^n$ (as usual the
unitary sphere is the oriented boundary of the unit disk $D^{n+1}$ of
$\R^{n+1}$ endowed with the standard orientation).  For every $x\in
V$, set
$$\fG(x)= (d_xf)^{-1}(\Bb) \ ; $$ by the very definition of
transversality, this defines a framing $\fG$ of $V$ in $M$. Hence we
have constructed a framed $k$-submanifold $(V,\fG)$ of $M$. We denote
by $[V,\fG]$ its class in $\eta^\Ff_k(M)$.

\smallskip

$\bullet$ If $M$ is oriented, then $V$ is also orientable and we
select a preferred orientation via the usual rule stated in Theorem
\ref{firstT}. This eventually leads to
$$[V,\fG]\in \Omega^\Ff_k(M) \ . $$

\smallskip

We have:

\begin{proposition}\label{pont-maps} 
(1) If $M$ is oriented, let us associate to every $\alpha \in [M,S^n]$ a class
$\pG_\Omega(\alpha)=[V,\fG] \in \Omega^\Ff_k(M)$ by means of an arbitrary
implementation of the procedure stated above. Then this actually well defines
the {\rm Pontryagin map} 
$$\pG_\Omega: [M,S^n]\to \Omega^\Ff_k(M) \ . $$

(2) If $M$ is non orientable, let us associate to every $\alpha \in [M,S^n]$ a class
$\pG_\eta(\alpha)=[V,\fG] \in \eta^\Ff_k(M)$ by means of an arbitrary
implementation of the procedure stated above. Then this actually well defines
the {\rm Pontryagin map} 
$$\pG_\eta: [M,S^n]\to \eta^\Ff_k(M) \ . $$
\end{proposition}

\Dim Every implementation of the procedure involves a few arbitrary
choices. We have to check that they are immaterial with respect to the
framed bordism class of the resulting framed (possibly oriented)
manifold $(V,\fG)$. Given $\alpha \in [M,S^n]$, let us assume first
that two implementions just differ by the choice of the maps $f_0$ and
$f_1$ in $\alpha$ and transverse to $x_0\in S^n$.  By the basic
transversality theorems, we can assume that a homotopy $F:M\times
[0,1]\to S^n$ which connects $f_0$ to $f_1$ is also tranverse to
$x_0\in S^n$; hence $W=F^{-1}(x_0)$ endowed with the framing $x\to
(d_xF)^{-1}(\Bb)$ gives rise to a framed cobordism between
$(V_0,\fG_0)$ and $(V_1,\fG_1)$ constructed by means of $f_0$ and
$f_1$ respectively. Assume now that the two implementations just
differ by the choice of the positive bases $\Bb_0$ and $\Bb_1$ of
$T_{x_0}S^n$. Then the resulting framed manifolds $(V,\fG_0)$ and
$(V,\fG_1)$ just differ by the framing.  As GL$(n,\R)$ is connected,
there is a smooth path $\Bb_t$, $t\in [0,1]$, of such bases connecting
$\Bb_0$ and $\Bb_1$.  Clearly this gives rise to a $1$-family of
framed manifolds of the form $(V,\fG_t)$, and eventually to a framing
of $V\times [0,1]$ properly embedded into $M\times [0,1]$ which
realizes a framed bordism between $(V,\fG_0)$ and $(V,\fG_1)$.
Finally, let us assume that we deal with two different points $x_0,
x_1 \in S^n$. By the homogeneity of $S^n$, there is a diffeotopy
$h_t$, $t\in [0,1]$, of $S^n$ such that $h_0={\rm Id}_{S^n}$,
$h_1(x_0)=x_1$. Given $f_0\in \alpha$, $f_0\pitchfork \{x_0\}$,
clearly also $f_1:= h_1\circ f_0$ belongs to $\alpha$ and
$f_1\pitchfork \{x_1\}$. Thanks to the above results, it is enough to
show that the framed manifold $(V_0,\fG_0)$ constructed by using $x_0,
\Bb, f_0$ and the framed manifold $(V_1,\fG_1)$ constructed by means
of $x_1, \Bb_1:=d_{x_0}h_1(\Bb), f_1$ belong to the same framed
cobordism class. This is easy to achieve by using the $1$-parameter
family of framed manifolds $(V_t,\fG_t)$ constructed by means of
$x_t:=h_t(x_0), \Bb_t:=d_{x_0}h_t(\Bb), f_t:=f_0\circ h_t$.  We have
understood that all these considerations work as well in the oriented
setting, as it is easy to see. The proposition is proved.

\cvd

\smallskip

We can state the main result of this Pontryagin construction.

\begin{theorem}\label{main-pont} Let $M$ be a compact, connected and boundaryless
  smooth $m$-manifold, $m\geq n \geq 1$, $k= m-n$. Then:

1) If $M$ is oriented, then the Pontryagin map
 $$\pG_\Omega: [M,S^n]\to \Omega^\Ff_k(M) $$
 is bijective.
 \smallskip
 
 2) If $M$ is non orientable, then the Pontryagin map
 $$\pG_\eta: [M,S^n]\to \eta^\Ff_k(M) $$
 is bijective.
 \end{theorem}
 
 Before giving a proof, let us state immediately an interesting
 corollary, early due to Hopf.
 
 \begin{corollary}\label{HOPF} Assume that $\dim M = \dim S^n \geq 1$. Then:
 
 1) If $M$ is oriented, then $f_0, f_1: M \to S^n$ are homotopic to each other if and only if
 $$ \deg_\Z (f_0)=\deg_\Z (f_1) \ . $$
 
 2) If $M$ is non orientable, then $f_0, f_1: M \to S^n$ are homotopic to each other if and only if
 $$ \deg_{\Z/2\Z} (f_0)=\deg_{\Z/2\Z} (f_1) \ . $$
 \end{corollary}
 \Dim It is enough to show that if the two maps have the same degree,
 then they are homotopic.  As $M$ and the sphere have the same
 dimension, the respective framed manifolds $(V_0,\fG_0)$ and
 $(V_1,\fG_1)$ constructed by means of $f_0$ or $f_1$ consist of a
 finite number of (possibly oriented) points.  Then it follows from
 the very definition of $\deg_R$, $R=\Z,\Z/2\Z$, that they are framed
 bordant (possibly in the oriented setting) if and only if the two
 maps have the same degree. The result follows by Theorem
 \ref{main-pont}.
 
 \cvd
 
 \smallskip  
 
 {\it Proof of Theorem \ref{main-pont}:} We will deal simultaneously
 with both Pontryagin's maps, understanding the necessary refinement
 in the oriented setting. Let us show first that the Pontryagin maps
 are onto. Let $(V,\fG)$ be a framed $k$-submanifold of $M$. It is
 enough to prove that there is a map $f:M \to S^n$ such that
 $[(V,\fG)]$ is produced by some implementation of the procedure used
 to define the Pontryagin maps, starting from the map $f$. As usual
 let us decompose the sphere as $S^n=D^+ \cup D^-$ such that $D^+\cap
 D^- = S^{n-1}$. By the stereographic projection from the northern
 pole, we can identify $D^-$ with the unit disk $D^n$; take $x_0=0\in
 D^n\subset S^n$.  By using the framing $\fG$, we can define a global
 trivialization
$$\tau: V\times D^n \to U$$
of a tubular neighbourhood of $V$ in $M$, such that the restriction of $\tau$ to $V\times \{0\}$
is the identity. Then we can define the map
$$\tilde f: U\to D^n, \ \tilde f (u):= \pi \circ \tau^{-1}$$
$\pi$ being the projection $V\times D^n \to D^n$. By construction:
\begin{itemize}
\item $\tilde f \pitchfork \{0\}$.
\item $\tilde f^{-1}(0)=V$.
\item Up to framed bordism (use again that GL$(n,\R)$ is connected), 
the framing $\fG$ can be recovered by the usual construction applied to $0$, $\tilde f$
and a basis $\Bb$ of $T_0D^n$.
\end{itemize}
\smallskip

By using a collar of $\partial U$ in $M$ and a collar bump function, it is not hard to extend
$\tilde f$ to a smooth map
$$f: M\to S^n$$
such that
\begin{itemize}
\item $f=\tilde f$ on $U$;
\item The map$f$ sends the complement of $U$ in $D^+$ and is constantly equal to the northern pole 
of $S^n$, say $\infty$, on the complement
of a slightly bigger tubular neighbourhood of $V$ in $M$;
\item $f^{-1}(0)=\tilde f^{-1}(0)=V$.
\end{itemize}

\smallskip

By construction such a map $f$ has the desidered property. So we have proved that
the Pontryagin maps are onto.
\smallskip

Let us prove now that they are injective. Let us say that a map $f: M \to S^n$
is {\it in standard form} if it has the qualitative properties of the map $f$
constructed above in order to prove the surjectivity. Let us prove first
the result for the restriction to the homotopy classes that admit representatives
in normal form. 

\begin{lemma} Assume that $f_0,f_1: M\to S^n$ are in standard form,
let $\alpha_0$ and $\alpha_1$ be the respective homotopy classes, and
assume that $\pG_*(\alpha_0)=\pG_*(\alpha_1)$. Then $\alpha_0= \alpha_1$.
\end{lemma}
\Dim Let $(V_0,\fG_0)$ and $(V_1,\fG_1)$ be framed manifolds obtained
by implementing the procedure with respect to $0, \ \Bb$ and $f_0$ or
$f_1$.  By hypothesis there is a properly embedded framed triad
$((W,V_0,V_1), \fG_W)$ in $M\times [0,1]$ which realizes a framed
bordism between them. Let us apply to the triad the construction used
above to define $\tilde f$. This produces a suitable map
$$ \tilde F: U_W \to D^n$$ where $U_W$ is properly embedded relative
tubular neighbourhood of $W$ in $M\times [0,1]$ which restricts to
tubular neigbourhoods $U_j$ of $V_j$ in $M$, $j=0,1$. As well as we
have extended above $\tilde f$ to $f:M\to S^n$ (in normal form), we
can extend $\tilde F$ to
$$ F: M\times [0,1] \to S^n $$ in relative normal form with respect to
$U_W$. As $f_0$ and $f_1$ are themselves in normal form by hypothesis,
up to diffeotopy we can assume that the restriction of $F$ to the
boundary recovers the given maps $f_0$ and $f_1$. Them $F$ establishes
a required homotopy between them.

\cvd

\smallskip

To achieve the proof of the main theorem, it is enough now to prove
that the assumptions in the above lemma are not restrictive.  Let $g:
M\to S^n$, it is not restrictive to assume that $g\pitchfork 0$, and
let $(V,\fG)$ obtained by implementing the usual procedure with
respect to $0, \ \Bb$ and $g$. Let $f: M \to S^n$ be a map in normal
form obtained as in the proof of surjectivity from $(V,\fG)$. Up to
diffeotopy we can assume that the tubular neighbourhood $U$ of $V$
which supports $\tilde f$ coincides with $g^{-1}(D^-)$ and that
eventually $g$ and $f$ coincide on $U$, both $f$ and $g$ send the
complement of $U$ in $D^+$ which retracts to $\infty$. Using this
facts it is an exercise to show that $f$ and $g$ are homotopic. This
completes the proof of the main Theorem \ref{main-pont}.

\cvd

\smallskip

\section{Characterization of combable manifolds}\label{combing}
Recall that a manifold is {\it combable} if it carries a nowhere vanishing
tangent vector field. We are now able to characterize this property.

\begin{theorem}\label{comb} Let $M$ be a compact connected boundaryless
smooth manifold. Then $M$ is combable if and only if $\chi(M)=0$.
In particular if $m=\dim M$ is odd, then $M$ is combable.
\end{theorem}
\Dim We already know that $\chi(M)=0$ is a necessary condition.  Let
us prove the other implication. Let $\vG$ any tangent vector field on
$M$ with isolated zeros. By using the homogeneity of $M$, up to a
diffeotopy we can assume that there is a chart $\phi: W \to \R^m$ such
that the zeros $x_1,\dots, x_k$ of $\vG$ are contained in $W$ and
their images are contained in the unitary disk $D^m \subset \R^m$.
For simplicity, let us keep the name $\vG$ for its expression in such
local coordinates, and $x_j$ for the images of the zero sets in $D^m$.
We can fix an auxiliary riemannian metric $g$ on $M$ which looks as
the standard euclidean metric $g_0$ on a neighbourhood of $D^m$. Fix a
system of small pairwise disjoint disks $D_j\subset D^m$, centred at
the $x_j$, $j=1,\dots, k$. The field $\hat \vG:= \vG/||\vG||_g$ is
well defined on $M\setminus \cup_j {\rm Int}(D_j)$ and homotopic to
the restriction of $\vG$.  The restriction of $\hat \vG$ to $D^m
\setminus \cup_j {\rm Int}(D_j)$ defines a map
$$\rho:  D^m \setminus  \cup_j {\rm Int}(D_j) \to S^{m-1} \ . $$
Assume at first that $M$ is oriented. By the bordism invariance of the degree
we have
$$ \deg_\Z (\rho_{|\partial D^m})= \sum_j \deg_\Z (\rho_{|\partial
  D_j})$$ and the second term is equal to $\chi(M)=0$. By Corollary
\ref{HOPF}, $\rho_{|\partial D^m}$ is homotopically trivial, hence can
be extended to a map $\hat \rho: D^m \to S^{m-1}$. By matching this
last map with the restriction of $\hat \vG$ to $M \setminus {\rm
  Int}(D^m)$, we eventually get a nowhere vanishing vector field on
$M$. If $M$ is not orientable, arguing similarly as in the proof of
Proposition \ref{unique-disk2} we can assume that the local picture at
$D^m$ agrees with the one in the oriented case, so we can conclude as
well.  \cvd

\smallskip

The above result extends to triads with a very similar proof.

\begin{proposition}\label{comb-triad} A smooth triad $(W,V_0,V_1)$ carries
a nowhere vanishing triad tangent vector field if and only if 
the relative characteristic $\chi(W,V_0)=0$.
\end{proposition}

\cvd

\smallskip

\section{On (stable) homotopy groups of spheres}\label{hg-sphere}
Accordingly with the basic motivation of the Pontryagin construction,
let us manage with
$$ \pi_m(S^n) \sim [S^m,S^n] \sim \Omega _{m-n}^\Ff(S^m)$$ for $m\geq
n >1$, in terms of framed bordism.  The first step is to transport on
$\Omega_{m-n}^\Ff(S^m)$ the group operation of $\pi_m(S^n)$. Recall
that the operation of the ordinary bordism modules is induced by the
disjoint union of representatives; moreover disjoint union and
connected sum belong to the same bordism class; this implies that
every ordinary bordism class can be represented by {\it connected}
manifolds.  The operation of the framed bordisms of the spheres is in
fact an embedded version of the disjoint union, again with the help of
connected sum.  Let $(V_1,\fG_1)$ and $(V_2,\fG_2)$ oriented framed
$(m-n)$-submanifolds of $S^m$, then the operation on
$\Omega_{m-n}^\Ff(S^m)$ is defined by
$$[V_1,\fG_1]+[V_2,\fG_2] = [(V_1,\fG_1) \amalg (V_2,\fG_2)]$$ where
we assume at first that the given framed manifolds are embedded into
two disjoint copies of $S^m$, and the disjoint union $(V_1,\fG_1)
\amalg (V_2,\fG_2)$ means the framed submanifold of $$S^m=S^m \# S^m$$
understanding that the connected sum is performed at disks which are
respectively disjoint from the two given framed submanifolds. It is
not hard to verify that this operation is well defined and recovers
(via the Pontryagin map) the usual operation of the homotopy group
$\pi_m(S^n)$.  By forgetting the embedding, we have immediately a
homomorphism of $\Z$-modules
$$ \phi_{k}: \Omega_{k}^\Ff(S^m)\to \Omega_{k}, \ k=m-n \ . $$ 

\begin{remark}\label{2connected-rep}{\rm
    In the ordinary setting we have noticed that every class has
    connected representatives. By means of embedded connected sums
    performed by attaching embedded $1$-handles, we can obtain that
    also every class in $\Omega^\Ff_k(S^m)$ has representative
    $[V,\fG]$ with connected $V$. This is easy if we forget the
    framing, a bit more demanding taking it into account. We left the
    details by exercise.}
\end{remark}

\smallskip 

As an immediate corollary of Corollary \ref{HOPF}, we have
\begin{proposition}\label{pi-m}
  For every $m\geq 2$, $\deg: \pi_m(S^m)\to \Z$ is an isomorphism of
  $\Z$-modules, and $[S^m,{\rm id}_{S^m}]$ is a generator of
  $\pi_m(S^m)$.
\end{proposition}

\cvd

\medskip

The same result was already known for $m=1$.

\subsection {The $J$-homomorphism}\label{J}
For every $m,n \geq 1$, there is an important homomorphism early defined by Whitehead
$$J: \pi_{m}(SO(n)) \to \pi_{m+n}(S^n)$$
which can be naturally expressed in terms of
$$J: \pi_m(SO(n)) \to \Omega_m^{\\Ff}(S^{m+n}) \ . $$
In fact by taking the usual equatorial embedding $S^m\subset S^{m+n}$,
every $\alpha \in \pi_m(SO(n))$ can be considered as a framing $\fG_\alpha$ of $S^m$ in $S^{m+n}$;
hence $J(\alpha)= [S^m, \fG_\alpha]$.

\subsection{Freudenthal's homomorphism and stable homotopy groups}\label{stable}
Let $S^m\subset S^{m+1}$ be the usual equatorial embedding. Set
$m=k+n$, so that $m+1= k+(n+1)$. If $(V,\fG)$ is an oriented framed
$k$-submanifold of $S^m$, then we can consider the framed
$k$-submanifold of $S^{m+1}$, say $(V, \sG\fG)$, where the framing
$\sG\fG$ is obtained by completing $\fG$ with the unitary normal
vectors along $S^m$ which point toward the northern pole of $S^{m+1}$.
It is easy from the definition of the operation that this induces a
$\Z$-modules homomorphism
$$\sG : \Omega_{k}^\Ff(S^m)\to \Omega_k^\Ff(S^{m+1})$$
whence, via the Pontryagin map,
$$ \sG: \pi_{n+k}(S^n) \to \pi_{n+1+k}(S^{n+1}) $$ called {\it
  Freudenthal suspension homomorphism}.  By using the same ``general
position argument'' used for the weak Whitney embedding theorem
(Corollary \ref {Weak -Whitney}) we have:

\begin {proposition}\label{freud-iso} For every $k\geq 1$, 

1) If $n\geq k+1$ then 
$$ \sG: \pi_{n+k}(S^n) \to \pi_{n+1+k}(S^{n+1})  $$
is surjective;
\smallskip

2)  If $n\geq k+2$ then 
$$ \sG: \pi_{n+k}(S^n) \to \pi_{n+1+k}(S^{n+1})  $$
is an isomorphism.
\end{proposition}

\cvd

\smallskip

One says that for every $k\geq 0$, the homotopy groups
$\pi_{n+k}(S^n)$ {\it stabilize} for $n\geq k+2$, being all isomorphic
to the (by definition) {\it stable homotopy group} denoted by
$\pi_k^\infty$.

By keeping the above notations, it is convenient to organize the
groups$$ \pi_{n+k}(S^n) \sim \Omega_k^\Ff(S^{n+k})$$ as being indexed
by the couples of integers $(k,n)$, $k\geq 0$, $n\geq 2$, endowed with
the lexicographic order. So for every $k$, by increasing $n$ we
encounter a few groups in the ``unstable regime", until we reach
$$\pi_k^\infty \sim \pi_{2+2k}(S^{k+2})\sim \Omega_k^\Ff(S^{2+2k}) \ . $$

\subsection{Homotopy groups of spheres for small $k$}\label{small-k}
``Small'' will mean $k\leq 3$.
Pontryagin himself succeeded to compute by geometric means the cases
$k\leq 2$ via his own construction. We will limit to a few indications
about these cases, the reader would fill all details by exercise or
refer to the exposition \cite{Pont} which contains detailed proofs.

\medskip

($k=0$) In agreement with Proposition \ref{pi-m}, the situation
stabilizes immediately:
$$ \pi_0^\infty \sim \pi_2(S^2)\sim \Z \ . $$

\medskip

($k=1$) The group in the unstable regime is
$$\pi_3(S^2) \sim \Omega_1^\Ff(S^3)$$ while $$\pi^\infty_1\sim \Omega^\Ff_1(S^4)\sim
\pi_4(S^3)\ . $$ Let us analyse the first one.  Every finite family of
embedded say $r$ smooth circles in $S^3$ can be transformed into the
boundary of $r$ pairwise disjoint embedded smooth $2$-disks by means
of a generic homotopy which is an embedding for every $t\in [0,1]$
with the exception of a finite number of instants at which two
branches of two circles (possibly the same one) cross each other with
distinct tangents.  Such a generic homotopy induces an embedded
framed bordism. So $\Omega_1^\Ff(S^3)$ is generated by classes of the form
$[S^1,\fG]$, where $S^1$ is the standard $S^1\subset S^2 \subset S^3$
via equatorial embeddings, hence such representatives only differ by the
framings. We can take as reference framing $\fG_0$ the one having as
first component a transverse field along a collar of $S^1$ in the
standard $2$-disk $D^+ \subset S^2$.
In fact $[S^1,\fG_0]$ corresponds
to $1\in \pi_3(S^2)$.  In this way every framing is of the form
$\fG=h_\fG\fG_0$ for a map $$h_{\fG}: S^1\to SO(2) \ . $$ As $SO(2)\sim S^1$,
the class $\alpha_\fG$  of $h_\fG$ belongs to  $\Z \sim \pi_1(SO(2))$.
We claim that $[S^1,\fG_1]=[S^1,\fG_2]\in \Omega^\Ff_1(S^3)$ if and only if $\alpha_{\fG_1}=\alpha_{\fG_2}$.
In fact if $f:S^3 \to S^2$  corresponds to $(S^1,\fG)$ via the Pontryagin construction,
then one realizes that $\alpha_{\fG}$ coincides with the linking number of two generic fibres of
$f$ over two distinct regular values (this is called the {\it Hopf number} of $f$).
Two maps with different Hopf number are not homotopic to each other.  
Then enventually we have that
$$\pi_3(S^2)\sim \Omega^\Ff_1(S^3)\sim \Z \ . $$ We can also exhibit a
geometric very interesting generator. This is the so called {\it
  Hopf map}: let $S^3$ be realized as the unitary sphere in $\C^2$ and
recall that $$\PP^1(\C)\sim S^2$$ the so called {\it Riemann
  sphere}. Then the mentioned map is
$$\hG: S^3 \to S^2$$
given by the restriction of the natural projection $\C^2\setminus \{0\} \to \PP^1(\C)$.
One can see that $\hG$ is a fibre bundle map with fibre $S^1$; the union of two
distinct fibres is the so called  (oriented) {\it Hopf link}
formed by two simply linked unknotted knots in $S^2$ with linking number equal to $1$.

With similar and easier considerations (now every embedding of $S^1$
is ``standard" by dimensional reasons), we see that
$\Omega^\Ff_1(S^4)$ is generated by classes of the form $[S^1,\fG]$,
and every framing induces a classifying map $\alpha_\fG \in
\pi_1(S0(3))$; we know that $SO(3)\sim \PP^3(\R)$ (see Example
\ref{so(3)}), so that $\pi_1(SO(3))\sim \Z/2\Z$, and eventually
 $$\pi^\infty_1\sim \Omega^\Ff_1(S^4)\sim \pi_4(S^3)\sim \Z/2\Z \ . $$
 Again we can exhibit geometric generators.
 We have  $$ \sG^{n-2}: \pi_3(S^2) \to \pi_{n+1}(S^n)$$
 then 
 $$\sG^{n-2}([\hG])=[\hG_n]$$  for a suitable ``suspended Hopf map"
 $$\hG_n: S^{n+1}\to S^n$$
 eventually generates $\pi_{n+1}(S^n)$ for $n\geq 3$.

\medskip

($k=2$) We have $\pi_4(S^2)$ and $\pi_5(S^3)$ in the unstable range, while
$\pi^\infty_2\sim \pi_6(S^4)$. It turns out that they are all isomorphic to $\Z/2\Z$.
Again we can exhibit geometric generators. In fact the class of the map
$$ \gG:=\hG\circ \hG_3: S^4 \to S^2$$
generates $\pi_4(S^2)$, while
$$ \sG^{n-2}([\gG]):=[\gG_n]$$
generates $\pi_{n+2}(S^n)$ for $n\geq 2$.
\medskip 

This is subtler to establish than the previous cases. It follows by the following
steps.

(a) The map 
$$\pi_ 4(S^3)\to \pi_4(S^2), \  [\alpha: S^4\to S^3]\to [\hG\circ \alpha]$$
is an isomorphism. Assuming it, $\pi_2(S^4)\sim \pi_4(S^3)\sim \Z/2\Z$
by the case $k=1$.
\medskip

(b) One constructs an explicit isomorphism
$$\delta: \pi_6(S^4)\to \Z/2\Z \ . $$
\medskip

(c) One proves that
$$\sG: \pi_4(S^2)\to \pi_5(S^3)$$
is onto. Assuming (a) (b), (c) and recalling that $\sG: \pi_5(S^3)\to \pi_6(S^4)$
is onto by Proposition \ref {freud-iso}, it follows that also $\pi_5(S^3)\sim \Z/2\Z$.
\medskip

Let us outline now a proof of these steps.

\medskip

(a): A basic fundamental tool in homotopy theory is the so called {\it
  homotopy long exact sequence of a fibre bundle} (see for instance \cite{Hu}, \cite{Hatch}). 
  We apply
  it to the Hopf fibration $\hG: S^3 \to S^2$ with fibre $S^1$;
  extract from the exact sequence the strings
  $$\dots \to \pi_n(S^1)\to \pi_n(S^3)\to \pi_n(S^2)\to \pi_{n-1}(S^1) \to \dots $$
  where the middle homomorphism is $\hG_*$ induced by $\hG$.
  As $\pi_m(S^1)=1$ for $m\geq 2$, we get that for $n\geq 3$,
  $$ \pi_n(S^3)\sim \pi_n(S^2)$$
  in particular
  $$\pi_4(S^3)\sim \pi_4(S^2) $$
  as desired. Note that this also proves again
  that $\pi_2(S^3)=\pi_2(S^2)\sim \Z$.

  \medskip

  (b) This is the most interesting step. To construct the isomorphism
  $\delta$ we will use several facts about surfaces discussed in
  Chapter \ref{TD-SURFACE}. Let $(V,\fG)$ be a framed surface in
  $S^6$, representing a class in $\Omega^\Ff_2(S^6)$. Assume that $V$
  is connected, then it is orientable of a certain genus $g\geq 0$. By
  dimensional reasons, up to diffeotopy $V$ is embedded in a standard
  way in $S^3 \subset S^6$. So only the framing contribution is
  relevant. Let $C$ be a compact oriented smooth circle on $V$.  The
  restriction of the framing $\fG=(s_1,\dots, s_{4})$ to $C$ can be
  completed by adding $s_5$ that is a normal field along $C$ tangent
  to $V$ which together with an oriented field tangent to $C$ gives
  the orientation of $T_xV$ at every $x\in C$. In this way we have
  constructed a framed circle $(C,\fG_C)$ representing an element of
  $\Omega^\Ff_1(S^6)\sim \Z/2\Z$.  Hence we can associate to
  $(C,\fG_C)$ the corresponding value $q(C):=q( [C,\fG_C])\in \Z/2Z$.
  Actually such a value does not depend on the orientation of $C$. If
  $L=\amalg_j C_j$ is a disjoint union of smooth circles on $V$, set
  $$q(L):= \sum_j q(C_j)\in \Z/2\Z \ . $$
  It is an istructive exercise to check that  the function $q$ defined so far  
  verifies the  conditions stated at the end of Chapter \ref{TD-SURFACE}; hence 
  
  \begin{lemma}\label{quad} The map
 $$q_{(V,\fG)}:  \eta_1(V)\to \Z/2\Z, \ q_{(V,\fG)}(\alpha)=q(C)$$
    provided that $C$ is any smooth circle on $V$ which represents $\alpha$,
    is a well defined quadratic
  enhancement of $(\eta_1(V),\bullet)$
  \end{lemma}
  
  \cvd
  
   \smallskip
  
  Then we can associate to $(V,\fG)$, the Arf invariant
  Arf$(q_{(V,\fG)}),\in \Z/2Z$. With more work one eventually realizes
  (recall also Remark \ref {2connected-rep}) that
  
  \begin{proposition}\label{delta}
  $$ \delta: \Omega^\Ff_2(S^6) \to \Z/2Z, \ \delta(\alpha) = {\rm Arf}(q_{(V,\fG)}) $$
  provided that $(V,\fG)$ represents $\alpha$ and $V$ is connected, is a well defined
  isomorphism.
  \end{proposition}
  
  Thus $\Omega^\Ff_2(S^6)$ is isomorphic to the Witt group $W^\HH_q(\Z/2\Z)$
  and realizes in a geometric way the formal non trivial enhancement of $\Omega_2=0$
  mentioned in Section \ref{quadratic}. $\Omega^\Ff_2(S^6)$ is generated by
  a framed torus $S^1\times S^1$ embedded in the standard way into $S^3\subset S^6$,
  such that the framing realizes $\HH^{1,1}$. Let us outline now the key step in 
  the proof of Proposition \ref{delta}. Let $(V,\fG)$ be as above, let $C$ be smooth circle
  traced on $V$, and assume that $q([C,\fG_C])=0$. Abstractly we can attach a $2$-handle
  to $V\times [0,1]$ at $V\times \{1\}$ in such a way that the embedded attaching tube
  is a tubular neighbourhood of $C$ in $V$. In this way we have constructed a triad
  $(W,V, V')$ such that  $g(V')=g(V)-1$. By easy dimensional reasons, we can extend
  the embedding $V\subset S^6$ to a proper embedding of the triad $(W,V,V')$ into
  $S^6\times [0,1]$. Then one realizes that the condition $q([C,\fG_C])=0$ is sufficient
  (and necessary) in order that this can be enhaced to a framed bordism between
  $(V,\fG)$ and $(V',\fG')$ for some framing $\fG'$. Moreover,   
  ${\rm Arf}(q_{(V,\fG)})=    {\rm Arf}(q_{(V',\fG')})$. By applying several times this
  argument, starting with an arbitrary $(V,\fG)$ we eventually reach either a framed
  sphere which represents the null class or a generating framed torus.

  \medskip
  
  (c) Here we will be very very sketchy. Given $f: S^5 \to S^3$, let
  $p, q\in S^3$ regular values such that both inverse images $V_p$ and
  $V_q$ are contained in $\R^5\subset S^5$.  As $\dim S^3$ is odd, the
  map
  $$ V_p\times V_q \to S^4, \ (x,y)\to \frac{y-x}{||y-x||}$$ has
  vanishing $\Z$-degree. Given $[V,\fG]\in \Omega^\Ff_2(S^5)$,
  $V\subset \R^5$, there is a generic projection of $V$ in $\R^4$; we
  can simplify the crossing points in the image of $V$ in $\R^4$ and
  eventually get $(V',\fG')$ framed bordant with $(V,\fG)$, such that
  $V'\subset S^4 \subset S^5$.  Let $f:S^5\to S^3$ be associated to
  $(V',\fG')$ via the Pontryagin construction. Assuming that $V'=V_p$,
  the vanishing of the degree of the map constructed as above with
  respect to $f$ eventually allows us to construct a framing
  $(V',\fG")$ representing an element in $\Omega^\Ff_2(S^4)$ whose
  suspension equals $[V',\fG']$.
   
 \medskip

($k=3$) This remarkably more complicated case was settled (by using
 the Pontryagin construction) by Rohlin in a series of four papers in
 1951-52 of great historical importance, mostly for the relation with
 the theory of $4$-manifolds.  We refer to \cite{GM} for the
 translation (in french) of these papers and wide deep commentaries.
 Here we limit to state the final results. We will come back on it in
 Chapter \ref{TD-4}, Section \ref{16}.

qThere is a quaternionic version of the Hopf map (recall Example \ref{so(3)})
$$ \hG^\HH: S^7 \to S^4$$ obtained in the following way. Let us
identify $\R^4$ with $\HH^2$, with quaternionic coordinates
$(q_0,q_1)$.  The unitary sphere $S^7$ is defined by the equation
$|q_0|^2+|q_1|^2=1$.  The group of unitary quaternion ($|q|=1$)
$SU(2)$ acts on $S^7$ by left multiplication. The quotient space is
diffeomorphic to $S^4$ and $\hG^\HH$ is just the quotient
projection. It is a fibre bundle map with fibre $S^3$. Then we have:

\smallskip

- $\pi_6(S^3)\sim \Z/12\Z$; 

\smallskip

- $\pi_7(S^4)\sim \Z \times \Z/12\Z $ where the first free factor is
generated by $[\hG^\HH]$, the finite factor is generated by the
suspension of a generator of $\pi_6(S^3)$;
\smallskip

- For every $n\geq 5$, $\pi_{n+3}(S^n)\sim \Z/24\Z$ and is generated
by $\sG^{n-2}([\hG^\HH])$.

\medskip

This geometric way of determining the homotopy groups of spheres
has been worked out only for $k\leq 3$ as we have outlined
above. Presumably the difficulty would increase too much with $k$. On
the other hand, the main interest (especially from the view point of {\it
low dimensional} differential topology) of such a direct method
consists in the method itself. Since Serre's thesis (\cite{Se}) powerful
tools (including the use of so called {\it spectral sequences}) have
been developed in homotopy theory; being just interested to the final result,
the above cases $k\leq 3$ become first ``trivial''
applications of these potent methods.  Moreover, one gets general
structural information; for example we have the following Serre's
result:

\begin{proposition}\label{Serre}  For every $k\geq 0$ and $n> 1$,
  the homotopy group $\pi_{n+k}(S^n)$ is finite with the following exceptions:
\smallskip

- $k=0$, as $\pi_n(S^n)\sim \Z$;
\smallskip

-  $k=2h-1$, $n=2h$, $h>0$, where $\pi_{n+k}(S^n)\sim \Z \oplus F$, $F$ being a finite group.   
\end{proposition}

\cvd

\smallskip

A great amount of researches concerns the determination of the
$p$-components of these homotopy groups for all primes $p\geq 2$.

Nevertheless, in spite of such powerful tools (see \cite{To}), the
full determination of the groups $\pi_{n+k}(S^n)$ has been not
achieved (not even of the stable groups $\pi_k^\infty$); in fact their
behaviour for increasing $k$ is apparently quite irregular, does not
present any kind of `stabilization'.
 
\section{Thom's spaces}\label{Thom-space}
Here the purpose is to rephrase for every $k> 0$ the determination of
the bordism $\Z/2Z$-vector spaces $\eta_k$ in terms of the homotopy
groups of certain so called {\it Thom's spaces}, say
$\TT^\eta_k$. Having as ideal model the Pontryagin construction, $S^n$
would be the ``Thom space'' for the framed bordism
$\Omega^\Ff_k(S^{n+k})$.
\medskip

To rich a setting closer to the Pontryagin construction, let us
recover first the ``absolute" bordism in terms of embedded one into
spheres.  For every sphere $S^m$, $m>k$, consider the sets
$\eta^{emb}_k(S^m)$ defined in Section \ref {emb-fram-bord}.  By means
of the embedded disjoint union already used above to define the
operation on $\Omega_k^\Ff(S^m)$, we can endow $\eta^{emb}_k(S^m)$
with a $\Z/2\Z$-vector space structure, so that the natural map
obtained by forgetting the embedding is a $\Z/2\Z$-linear map:
$$ \phi_{k,m}: \eta^{emb}_k(S^m)\to \eta_k \ . $$

Via the usual equatorial embedding $S^m \subset S^{m+1}$, we get
linear maps
$$ \sG_{k,m}: \eta^{emb}_k(S^m)\to \eta^{emb}_k(S^{m+1}) \ . $$ By
means of general position considerations as in the weak Whitney
embedding theorem, and dealing also with proper embeddings into $S^m
\times [0,1]$, we easily have:

\begin{lemma}\label{emb-eta}

1) If $m\geq 2k+1$, then $\phi_{k,m}$
is onto;

2) If  $m\geq 2k+2$, then $\phi_{k,m}$ is a isomorphism; moreover
$\phi_{k,m}= \phi_{k,m+1}\circ \sG_{k,m}$.
\end{lemma}

\cvd

\smallskip

$\bullet$ From now on we stipulate that for every $k>0$ we will take $m\geq 2k+2$, and set $h=m-k$.

\smallskip

Let $M$ be a $(r+h)$-manifold which is the interior of  a (possibly boundaryless)
compact smooth manifold with boundary; let $Y\subset M$ a
boundaryless compact $r$-submanifold. The following facts are now wellknown:   
\smallskip

{\it If $f: S^m \to M$ is transverse to $Y$, then $V_f=f^{-1}(Y)$ is a compact boundaryless 
$k$-submanifold of $S^m$; if $f_0$ and $f_1$ are homotopic and both transverse to the zero section, then
$[V_{f_0}]=[V_{f_1}] \in \eta^{emb}_k(S^m)$. Then by applying the transversality theorems we well define
the map
$$ [S^m, M]\to \eta^{emb}_k(S^m), \  \alpha=[f:S^m \to M] \to [f^{-1}(Y)]$$
provided that $f$ is any representative of $\alpha$ transverse to $Y$. Recall that in our situation
$$[S^m,M] \sim \pi_m(M) \ . $$}

\medskip 

This would suggest to look for such a pair $(M,Y)$ (if any) such that the above map is bijective.
Recall that  the pair $(S^n,\{x_0\})$ has played this role with respect to the framed 
bordism $\Omega^\Ff_k(S^{n+k})$.

With this perspective in mind, let us recall a construction already employed in Section \ref {smooth-Cr}. 
For every $(k,m)$ as above, $h= m-k$, take the tautological vector bundle
$$\tau: \Vv(\GG_{m,h})\to \GG_{m,h} \ ,$$
the grassmannian $\GG_{m,h}$ being identified with the zero section of this bundle.
As usual present the sphere as $S^m= \R^m \cup \infty$; up to diffeotopy every compact boundaryless $k$-submanifold
$V$ of $S^m$ misses $\infty$, that is $V\subset \R^m \subset S^m$.  Let
$$\nu: V\to \GG_{m,h}, \ \nu(x)= (T_xV)^\perp $$
be the orthogonal distribution of $h$-planes along $V$ with respect to a riemannian metric on $\R^m$, for instance the standard
one $g_0$. We can use $\nu$ to build a tubular neighbourhood $p:U\to V$ of $V$ in $\R^m$ and this can be incorporated
into a commutative diagram of maps
$$ \begin{array}[c]{ccc}
U&\stackrel{ \tilde f }{\rightarrow}& \Vv(\GG_{m,h})\\
\downarrow\scriptstyle{p}&&\downarrow\scriptstyle { \tau}\\
V&\stackrel{\nu}{\rightarrow}&\GG_{m,h} \end{array}$$
where the image of $\tilde f$ is a tubular neigbourhood of the zero section
in $\Vv(\GG_{m,h})$, $\tilde f$ is a fibred map onto its image, hence transverse to 
$\GG_{m,h}$, and $\tilde f^{-1}(\GG_{m,h})=V$.
Although it would be tempting to take $(M,Y)=(\Vv(\GG_{m,h}), \GG_{m,h})$,
one immediately realizes that there are not reasons that $\tilde f$ can be extended to the whole
of $S^m$. The situation is very similar to the step in the proof of the surjectivity of the
Pontryagin map when we have constructed the map also called $\tilde f: U \to \R^n$,
where $(\R^n,\{0\})$ played the role of $(\Vv(\GG_{m,h}), \GG_{m,h})$. The key fact
that allowed us to extend that $\tilde f$ to a map $f: S^m \to S^n$, was that the complement
of the image of $\tilde f$ retracts to the northern pole of $S^n$; note that $S^n= \R^n \cup \infty$ can be considered as the
one-point compactification of $\R^n$. This suggests a very simple way to compactify $\Vv(\GG_{m,h})$ in order to
make the extension of the map $\tilde f$ possible. Set
$$ \TT^\eta_{m,h} := \Vv(\GG_{m,h})\cup \infty$$
that is the one-point compactification. This space has some remarkable features
\begin{itemize}
\item It is no longer a manifold; however the only non manifold point
is just the added point at infinity; 

\item This point $\infty$ has a fundamental system of conical neighbourhhoods
centred at it and with base diffeomorphic to  the total
space of the unitary bundle  of the tautological bundle $\tau$;

\item The one-point compactification (which is isomorphic to the sphere $S^h$)
of every fibre of $\tau$ is embedded into $\TT_{m,h}$ which can be considered
as the wedge of such infinite family of $h$-spheres, based at $\infty$; 

\item $\TT^\eta_{m,h} \setminus \GG_{m,h}$ retracts to $\infty$.
\end{itemize}

So although it is not a manifold, $\TT^\eta_{m,n}$ is a ``honest'' rather tame path connected compact space
(in particular it has a structure of finite CW complex) whose homotopy groups are suited to be treated
by the powerful tools mentioned above.

Then arguing similarly to the Pontryagin construction, we can extend the above map
$$\tilde f: U\to \Vv(\GG_{m,h})$$
to a map
$$ f: S^m \to \TT^\eta_{m,h}$$
such that the complement of $U$ is mapped into the complement of the image of $U$ in $\TT^\eta_{m,h}$,
$f$ is constantly equal to $\infty$ on the complement of a sligthly bigger tubular neighbourhood say $U'$
of $V$ in $S^m$, $f$ is smooth on $U'$. Let us say that a map sharing these properties of $f$ is {\it in
standard form}. Similarly to the end of the proof of Theorem \ref {main-pont}, we have
\begin{lemma} \label{standard} Every $\alpha \in [S^m, \TT^\eta_{m,h}]$ has representatives in standard form.
\end{lemma}
\Dim Let $\alpha=[g: S^m\to \TT^\eta_{m,h}]$. Up to a first homotopy
we can assume that $g$ is smooth on $D^-\subset S^m$ (as usual
$D^-\sim D^m$), $g^{-1}(\infty) \cap D^- = \emptyset$ and $g_{|D^-}$
is transverse to $\GG_{m,h}$. Then we can construct $f: S^m \to
\TT^\eta_{m,h}$ in normal form which coincides with $g$ on the tubular
neighbourhood $U$ of $V=g^{-1}(\GG_{m,h})$ involved in the
construction of $\tilde f$, whence of $f$ itself. Set $A=
g(U)=f(U)$. As $\TT^\eta_{m,h} \setminus A$ is contractible to
$\infty$, we can conclude that $g$ and $f$ are homotopic.

\cvd

\smallskip

We summarize the above discussion in the following main result of the
present section. Thanks to Lemma \ref{standard} the proof runs
parallel to the one of Theorem \ref{main-pont}, details are omitted.
\begin{theorem}\label{main-thom} For every $k>0$, $m\geq 2k+2$, $h=m-k$, the map
$$ \tG_{m,h}: [S^m,\TT^\eta_{m,h}]\to \eta^{emb}_k(S^m), \ \tG_{m,h}(\alpha)= [f^{-1}(\GG_{m,h})]$$
provided that $f: S^m \to \TT^\eta_{m,h}$ is any representative in normal form of $\alpha$, is well defined
and eventually establishes group isomorphisms
$$ \pi_m(\TT^\eta_{m,h}) \sim \eta^{emb}_k(S^m) \sim \eta_k \ . $$
\end{theorem}

\cvd

\smallskip

Every such a $\TT^\eta_{m,h}$ is called a Thom spaces for $\eta_k$. Sometimes one prefers to write them as
$\TT^\eta_{ k+h,h}$; the homotopy groups $\pi_{k+h}(\TT^\eta_{k+h,h})$ stabilize when $h\geq k+2$.

\subsection{On Thom's spaces for $\Omega_k$}\label{T-Omega}
First we identify $\Omega_k$ with $\Omega^{emb}_k(S^m)$, $m\geq 2k+2$. Then we replace
the tautological bundle $\tau$ with the tautological bundle of the grassmannian of {\it oriented}
$h$-planes in $\R^m$ (see Chapter \ref{TD-SMOOTH-MAN})
$$\tilde \tau: \Vv(\tilde \GG_{m,h}) \to \tilde \GG_{m,h} \ . $$
Note that the fibres of this bundle are tautologically oriented. 
Set $\TT^\Omega_{m,h}$ the one-point compactification of $\Vv(\tilde \GG_{m,h})$.
For every $[V]\in \Omega_k(S^m)$,
in a very similar way as above, we can construct 
$$\tilde f: U\to \Vv(\tilde \GG_{m,h})$$
which extends to a map in normal form 
$$ f: S^m \to \TT^\Omega_{m,h}$$
in such a way that the given orientation of $V$ coincides with the one obtained by
the usual rule already employed in the Pontryagin construction by means of the orientation
of $S^m$ and the transverse orientation to $V$ induced, in that case, by the framing. Arguing similarly
to the $\eta$-case we eventually get:
\begin{theorem}\label{main-thom-omega} For every $k>0$, $m\geq 2k+2$, $h=m-k$, the map
$$ \tilde \tG_{m,h}: [S^m,\TT^\Omega_{m,h}]\to \Omega^{emb}_k(S^m), \ \tilde \tG_{m,h}(\alpha)= [f^{-1}(\tilde \GG_{m,h})]$$
provided that $f: S^m \to \TT^\Omega_{m,h}$ is any representative in normal form of $\alpha$, is well defined
and eventually establishes group isomorphisms
$$ \pi_m(\TT^\Omega_{m,h}) \sim \Omega^{emb}_k(S^m) \sim \Omega_k \ . $$
\end{theorem}
 
 \cvd
 
 \smallskip  

Every such a $\TT^\Omega_{m,h}= \TT^\Omega_{k+h,h}$ is called a Thom spaces for $\Omega_k$; again 
the homotopy groups $\pi_{k+h}(\TT^\Omega_{k+h,h})$ stabilize when $h\geq k+2$.

\subsection{Determination of $\eta_\bullet$}
The homotopy groups $\pi_m(\TT^\eta_{m,h})$ look qualitatively simpler
than in the case of spheres as we already know for example that they
are finite dimensional $\Z/2\Z$-vector spaces. In fact they can be
computed by advanced homotopy theory methods (\cite {Se}), providing
the full determination of
$$\eta_\bullet = \oplus_k \eta_k \ . $$
Recall that $\eta_\bullet$ has furthermore a $\Z/2\Z$- graded algebra structure where the product
is induced by the cartesian product of manifolds:
$$ [V]\cdot [W] = [V\times W] \ . $$
This has been noticed in Remark \ref {ring-point}; here we omit the cobordism reindexing
$\eta_k=\eta_k(x_0)\sim \eta^{-k}(x_0)=\eta^{-k}$. In \cite{T} one eventually determines these algebra.
Here we limit to the statement:
\begin{theorem}\label{eta-alg}  
The $\Z/2\Z$-graded algebra $\eta_\bullet$ is isomorphic to the polynomial 
algebra 
$$\Z/2\Z[X_i; \ i\in J]$$
where 
$$J= \N \setminus \{2^j-1; \ j\in \N\} \ . $$ 
\end{theorem}

\cvd

\smallskip

We can also give explicit geometric generators (see \cite {M5}).  For
every $m\leq n$, let $H_{m,n}$ denote the regular real algebraic
hypersurface in the product of projective spaces $\PP^m(\R)\times
\PP^n(\R)$ defined in terms of the respective homogeneous coordinates
$(w_0, \dots, w_m)$ and $(z_0, \dots, , z_n)$ as the locus
$$H_{m,n} = \{ w_0z_0 + w_1z_1 + \dots + w_mz_m = 0\} \ . $$

Set
$$ \{X_{2j}:= [\PP^{2j}(\R)], \ j>1\}$$
$$ \{X_{2^{k+1}+1}:=[H_{2^k,2+2^k}], \ k>1\} \ . $$ To show that their
union is a family of independent generators it is enough to show that
for every $i\in J$ there exists a unique representative $X_i=[V_i]$ in
the family and that for every finite product of such $V_i$, there is a
non vanishing (stable) $\eta$-characteristic number (recall Section
\ref{stable-complete}).  This last task is easy for even indices $2j$
and the E-P characteristic mod$(2)$ suffices.  In general it is easier
if one would dispose of the cohomological formulation in terms od
Stiefel-Whitney numbers (see Remark \ref {sw}).

\smallskip
As a remarkable qualitative consequence we have
\begin{corollary}\label{alg-bord-gen}
  For every $k\geq 0$, every $\alpha \in \eta_k$ can be represented by
  regular real algebraic sets (projective indeed).
\end{corollary}

\smallskip

The determination of $\Omega_\bullet$ can be performed in the same
vein, however the proof, even the statement are more complicated (see
\cite{Wall}).

\subsection{On Nash-Tognoli theorem}\label{nash-tognoli}
We have discussed in Chapter \ref{TD-COMP-EMB} how every compact
boundaryless $m$-submanifold $M$ of $\R^n$ can be approximated by a
Nash manifold $M'$ (normal if the embedding dimension is big enough).
As already said, in his
paper \cite {Na}, Nash stated also a few conjectures/questions towards
potential improvements of this result (see also Sections
\ref{stable-2-Nash}, \ref{tear}). The most natural conjecture was that
$M$ can be approximated by a regular real algebraic set (not only by some
``analytic sheet'' of it).
A first step was accomplished In \cite{Wa2} by proving the conjecture under
the restrictive hypotheses that the embedding dimension is big enough (as for normality),
and
$[M]=0\in \eta_m$ i.e. it is a boundary $M=\partial W$.  Roughly, one
realizes the double $D(W)\subset \R^n$ in such a way that $M$ is the
transverse intersection of $D(W)$ with a hyperplane $P$.  Then one show that
$D(W)$ can be approximated by a normal Nash manifold $N$ made by regular
components of a real algebraic set  $X$ such that $X\setminus N$
is far from the hyperplane.  Finally
$M'=P\pitchfork X$ is a required regular real algebraic approximation
of $M$. Corollary \ref{alg-bord-gen} can be rephrased by saying that the
conjecture hols {\it up to bordism}. By
using this fact, the actual conjecture has been proved in general \cite{Tog},
again assuming that the embedding dimension is big enough.
By that Corollary, there is a regular $m$-dimensional
real algebraic set $\Sigma$ such that $M\amalg \Sigma = \partial W$. A
suitable {\it relative approximation theorem} allows us to refine the
above construction in such a way that
$$P\pitchfork X = M' \amalg \Sigma \ ; $$ as both $M'\amalg \Sigma$
and $\Sigma$ are regular algebraic sets, it is not hard to conclude
that also $M'$ is regular algebraic so that it is a required
approximation of $M$.  In \cite{Ki}, one refines the Nash-Tognoli
theorem in the projective setting, and proves that $M\subset
\PP^n(\R)$ can be approximated by regular algebraic subsets of the
projective space.  For more details about this matter see \cite{BCR}.

\chapter{High dimensional manifolds}\label{TD-HIGH}
``High" means of dimension greater or equal to $6$.
The reason of this specific opposition
``low dimensions less or equal to $5$" {\it vs} ``high dimensions greater or equal to $6$" 
mainly depends on the fact 
that in high dimension Smale's \cite{S2} $h$-{\it cobordism
theorem} holds and, moreover, we have a ``stable'' differential/topological  proof,
 in the sense that it works in the same way for every high dimension. Such a 
 proof definitely does not work for low dimensions. In  dimension $5$ the $h$-cobordism theorem
 fails and this reflects specific phenomena of persistent {\it geometric} intersection between surfaces embedded 
 in boundaryless compact simply connected $4$-manifolds, although they have 
 vanishing {\it algebraic} intersection number. 
 In dimension $4$ the proof does not apply  because
 of specific geometric linking phenomena between knots in $S^3$ with vanishing (algebraic) 
 linking number; the validity of the $4$-dimensional $h$-cobordism theorem still is an {\it open question}.
 The $3$-dimensional $h$-cobordism theorem is equivalent to the celebrated
 {\it Poincar\'e conjecture}; this last has been proved rather recently by means of deep 
 $3$-dimensional methods of {\it geometric analysis}.
 In a sense dimension $5$ is really in the border between the two regimes; as already said, 
it is infuenced by  the behaviours of four dimensional manifolds; on the other hand,  with some 
specific additional care, shares some remarkable behaviours with higher dimensions.  

In this Chapter we will not provide a proof
of the $h$-cobordism theorem (see \cite{M3} for a  proof in terms of
Morse functions, see \cite{RS} for a proof in terms of handle decompositions
which actually works also for PL manifolds); rather we will focus a key point
 where the high dimensional assumption is crucial.

Together with Chapter  \ref{TD-SURFACE}, Chapters \ref{TD-3} and \ref{TD-4} will be devoted
to some aspects of low dimensional theory.

\section{On the $h$-cobordism theorem} \label{h-cob}
Let us start with a definition.

\begin{definition}\label{h-cob-defi}{\rm Let $(W,V_0,V_1)$ be a smooth $m$-dimensional triad ($m=\dim W$).
 It is a $h$-{\it cobordism} if both inclusions $j_i: V_i\to W$, $i=0,1$, are homotopy equivalences (i.e. they
 have an inverse up to smooth homotopy $r_i: W\to V_i$ such that (by definition) $r_i \circ j_i$ is homotopic to ${\rm id}_{V_i}$,
 $j_i \circ r_i$ is homotopic to ${\rm id}_W$).}
 \end{definition}

\smallskip

The basic example of $h$-cobordism is a cylinder $(V\times [0,1],V, V)$. The general vague question is
under which minimal hypothesis the cylinders are the unique instance of $h$-cobordism up to diffeomorphism of triads.
We can formulate the following more specific question:

\begin{question}\label {h-cob-ques}(Simply connected $m$-dimensional $h$-cobordism question) {\rm Let $(W,V_0, V_1)$ be a $h$-cobordism, $\dim W=m$; assume that $W$ 
(whence both $V_0$ and $V_1$)
is {\it simply connected}. Is it true that the triad is diffeomorphic to the cylinder $(V_0\times [0,1], V_0,V_1)$, so that, in particular, $V_0$ is diffeomorphic to $V_1$?}
\end{question}

\medskip

Note that the question is empty for $m=2$. Assume the positive answer, let us derive some important consequences.

\begin{proposition}\label{h-app}
Assume that $m$-dimensional simply connected $h$-cobordisms are diffeomorphic to cylinders. Then we have:

\smallskip

(1) {\rm (Characterization of the $m$-disk)} Every contractible compact $m$-manifold $M$ with simply connected boundary is diffeomorphic to the closed disk $D^m$.

\smallskip

(2) {\rm (Generalized Poincar\'e conjecture)} If $\Sigma$ is a compact $m$-manifold which is homotopically equivalent to $S^{m}$ (i.e. it is a {\rm homotopy sphere}), 
then it is {\it homeomorphic} to $S^{m}$.
\smallskip

(3) {\rm (Smooth Schoenfliess property)} If $\Sigma$ is a smooth embedded $(m-1)$-sphere in $S^{m}$, then there is a diffeotopy of $S^{m}$ that sends
$\Sigma$ onto the standard equator $S^{m-1}\subset S^{m}$.
\end{proposition}

\smallskip

{\it Sketch of proof.} Some of the facts claimed below are not so evident; to prove them one would dispose of more advanced algebraic/topological
tools; we limit to an outline.  
\smallskip

(1) Remove from $M$ a $m$-disk $D$ standarly embedded into a chart of $M$. Set $W= M \setminus {\rm Int}(D)$. The triad $(W, \partial D, \partial M)$
is a simply connected $h$-cobordism, hence it is diffeomorphic to the cylinder $(S^{m-1} \times [0,1], S^{m-1}, S^{m-1})$ and $M$ is 
diffeomorphic to the manifold obtained by gluing
$D$ to this cylinder by a diffeomorphism $\phi: \partial D \to S^{m-1}\times \{1\}$; it is not hard to conclude that $M$ is diffeomorphic to $D^m$.
\smallskip

(2) Remove from  $\Sigma$ a standard $m$-disk $D$ in a chart as above. $M= \Sigma \setminus {\rm Int}(D)$ verifies the hypothesis of item (1),
then it is diffeomorphic to a disk, $\Sigma$ is eventually a twisted sphere (see Section \ref {twisted-sphere} )  and we know that it is homeomeorphic
(not necessarily diffeomorphic)  to $S^m$.

\smallskip

(3)  By the separation theorem of Section \ref{sep-teo}, $S^m \setminus \Sigma$ has two connected components, the closure of each one of these components
verifies the hypothesis of item (1), hence it is an embedded smooth $m$-disk in $S^m$ and we conclude by means of  the uniqueness of disks up to diffeotopy.  

\cvd

\smallskip

\begin{remark}{\rm The above proposition shows that the $h$-cobordism question is strictly related to (in fact motivated by)  basic fundamental questions about
the topology of smooth manifolds. For example for $m=3$,  if $(W,V_0,V_1)$ is a simply connected $h$-cobordism, then $V_0\sim V_1\sim S^2$ by the classification
of surfaces. As a $3$-dimensional twisted sphere is a true sphere, it follows  that a positive answer  to question \ref{h-cob-ques} for $m=3$ is equivalent to the validity of the original celebrated {\it Poincar\'e conjecture},
with furthermore the refinement that for $m=3$ every smooth homotopy sphere $\Sigma$ is {\it diffeomorphic} to $S^3$.  Probably the reader is aware that this has been
proved by G. Perelmann at the beginning of the new century, by achieving the  program based on the Ricci flows of riemannian metrics on $3$-manifolds,
early introduced by R. Hamilton. We stress that this peculiarly $3$-dimensional geometric/analytic approach is very far from the differential/topological methods discussed in this text.
As the $3$-dimensional Poincar\'e conjecture is true, then if $(W,V_0,V_1)$ is a simply connected $4$-dimensional $h$-cobordism, then
$V_0\sim V_1 \sim S^3$. Thus, as a twisted $4$-sphere is a true sphere, a positive answer to question \ref{h-cob-ques} for $m=4$ is equivalent to the fact that 
every smooth $4$-dimensional homotopy sphere is actually {\it diffeomorphic} to $S^4$. This still is an {\it open question}, as well as the validity of 
the $4$-dimensional smooth  Schoenfliess property. On the other hand we recall  that the {\it purely topological} $4$-dimensional Poincar\'e conjecture
 (even dealing with topological not necessarily smooth $4$-manifolds)  has been proved in 1982 by M.H. Freedman \cite{ Fr}.

}
\end{remark}

Now we can state the {\it high dimensional simply connected $h$-cobordism theorem}.

\begin{theorem}\label{h6} Let $(W,V_0,V_1)$ be a simply connected $h$-cobordism,
$\dim W \geq 6$. Then it is diffeomorphic to the cylinder $(V_0\times [0,1], V_0, V_0)$.
\end{theorem}

\smallskip

Hence all consequences stated in Proposition \ref{h-app} hold for $m\geq 6$. 
We have mentioned before that although  the $h$-cobordism theorem fails for $m=5$,
nevertheless this dimension shares some behaviour with higher dimensions. 
Referring to the statement of Proposition \ref{h-app}, we recall for example (without proof) that:

\medskip

(1) The characterization of the $5$-disk holds under the {\it stronger} hypothesis that the boundary of the contractible $5$-manifold $M$
is {\it diffeomorphic} to $S^4$;
\smallskip

(2) The $5$-dimensional generalized Poincar\'e conjecture holds true;
\smallskip

(3) The $5$-dimensional smooth Schoenfliess property holds true. 

\medskip

\subsection{On the proof of the high dimensional $h$-cobordism theorem}\label{on-h-proof}
The strategy to prove the $h$-cobordism theorem is based on handle decompositions (refer to Chapter
\ref{TD-HANDLE}). Given a simply connected $h$-cobordism $(W,V_0,V_1)$, $\dim W = m$, one can start with an ordered
handle decomposition 
$$C_0\cup H_1^{q_1} \cup \dots H_k^{q_k} \cup C_1$$
without  $0$- and $m$-handles (Proposition \ref{0-m-elimination}). If necessary we can also assume that
the handles of a given index $q<m$ are attached simultaneously at pairwise disjoint attaching tubes.
Note also that in the hypothesis of the theorem, all involved manifolds ($W$ and all submanifolds $W_r$ of $W$ 
obtained by attaching till the  $r$th-handle) are orientable.
We dispose of  two basic handle moves in order to try to make it simpler and simpler. 
If we succeed to eventually reach a decomposition without handles of any index, then the theorem will be proved.
A priori the only way we dispose to reduce the number of handles is the cancellation of pairs of complementary handles.
The core of the proof is a much more flexible {\it cancellation theorem} which applies in the setting of the theorem.
Consider a fragment of a given handle decomposition of the form
$$ \dots \cup  H_r^q \cup H_{r+1}^{q+1} \cup \dots \ $$
Then both the (embedded) $b$-sphere $S_b$ of $H^q_r$ and the $a$-sphere $S_a$ of $H_{r+1}^{q+1}$ are submanifolds
of $\partial W_r$ and $\dim S_b + \dim S_a = \dim \partial W_r = m-1$. So fixing auxiliary orientations, we can
compute  their intersection number in $\partial W_r$, $[S_b] \bullet [S_a]\in \Z$.

\begin{definition}\label{alg-comp}{\rm In the situation depicted above, we say that   $H_r^q \cup H_{r+1}^{q+1}$ is a pair of 
{\it algebraically} complementary handles if $[S_b] \bullet [S_a]= \pm 1$.}
\end{definition}

\smallskip

Obviously this extends the notion of complementary handles. Now we can state a {\it stronger cancellation theorem}.

\begin{theorem}\label{strong-canc} Let $(U, Z_0,Z_1)$ be a smooth triad of dimension $m$ 
which admits a handle decomposition 
$$C_0 \cup H^q \cup H^{q+1}\cup C_1$$
made by two algebraically complementary handles.
Assume that both $Z_0$ and $Z_1$ are 
simply connected, and that 
$$m\geq 6, \ q\geq 2, \ m-q \geq 4 \ . $$ Then the given triad is diffeomorphic
to the cylinder $(Z_0 \times [0,1], Z_0,Z_0)$.
\end{theorem} 

\medskip

The idea in order to prove the stronger cancellation theorem is clear. By transversality and handle sliding, we can assume 
that $S_b\pitchfork S_a$ in $\partial M$, $M:= C_0\cup H^q$ and that the intersection  consists of an odd number of signed points, 
such that the sum of the signs is equal to $\pm 1$.
So by means of handle sliding, one would progressively cancel pairs of intersection points of {\it opposite sign}, so that at the end one reaches a 
decomposition made by two genuine complementary handles which can be cancelled. In the discussion on the strong Whitney 
embedding theorem (Section \ref{whitney-trick}) of compact $n$-manifolds into $\R^{2n}$, for $n\geq 3$,  
we have already mentioned the so called ``{\it Whitney trick}'' as a tool in order to cancel pairs of crossing points.
The hypotheses of the stronger cancellation theorem allow to apply it.
This will be discussed with some care in the next section. 

\section{Whitney trick and unlinking spheres into a sphere}\label{WT-U}
First we state a lemma under  the hypotheses of the stronger cancellation theorem.

\begin{lemma}\label{pre-WT} In the hypotheses of  Theorem \ref{strong-canc}, denote by $\partial (C_0\cup H^q)= Z_0 \amalg M$, so that the $b$-sphere
$S_b$ of $H^q$ and the $a$-sphere $S_a$ of $H^{q+1}$ are transverse submanifolds of $M$. Then $M\setminus (S_b\cup S_a)$ is simply connected.
\end{lemma}
 \Dim Set $m=n+1$. Denote by $S'_a$ the $a$-sphere of $H^q$.
 Its codimension $\dim Z_0 - \dim S'_a = n- (q-1) \geq 4$.
 Then by transversality also $Z_0\setminus S'_a$ is simply connected; as both $Z_0\setminus S'_a$ and $M\setminus S_b$ 
retract onto $Z_0\setminus {\rm Int}(T'_a)$, it follows that also $M\setminus S_b$ is simply connected.
The codimension of $S_a$ is $\dim M - q= n-q\geq 3$. So by the same transversality argument we have
that $(M \setminus S_b)\setminus S_a = M\setminus (S_b \cup S_a)$ is simply connected.

\cvd

\smallskip

Referring to the last lemma, we can abstractly formalize  some features of the situation occurring on the manifold $M$.
\smallskip

By {\it a situation $(M,R,S,\pm x)$ of type $(n,r)\in \N^2$} we mean:
\smallskip

- $M$ is a connected oriented boundaryless smooth manifold of dimension  
$ n $;

\smallskip

- $R$ and $S$ are boundaryless compact connected oriented submanifolds of $M$ such that
$\dim R = r$, $\dim S=s$, $n> s\geq r>0$, $r+s=n$, $R\pitchfork S$.
\smallskip

- $M\setminus (S\cup R)$ is simply connected;
\smallskip

- $x_\pm \in R\cap S$ are intersection points of {\it opposite sign}.  

\medskip

\begin{remarks}\label{simply-c} {\rm (1) In a situation of type $(n,r)$ as above, 
if both codimensions of $S$ and $R$ are greater or equal to $3$, then
by an  usual transversality argument, $M\setminus (S \cup R)$ is simply connected if and only if $M$ is simply connected.
\smallskip

(2) In situations arising under the hypotheses of Theorem \ref{strong-canc}, 
we have furthemore that $n\geq 5$ and $r\geq 2$.
}
\end{remarks}

\smallskip

{\it (Whitney disk)} Let $(M,R,S,\pm x)$ be a situation of type $(n,r)$. By a {\it Whitney disk} $D$
for $(M,R,S,\pm x)$ we mean the realization of the following pattern (recall Section \ref{whitney-trick}) 

\smallskip

(1) There is an embedded smooth circle $\gamma$ in $R\cup S$ with two corners
at $\pm x$; these divide $\gamma$ in two arcs with closures say $\gamma_R$ and $\gamma_S$
respectively; $\gamma_R$ (resp. $\gamma_S$) is contained into an smooth open $r$-disk ($s$-disk)
$U_R\subset R$ ($U_S\subset S$); $U_R\cup U_S$ is a neighbourhood of $\gamma$ in $R\cup S$;
$U_R \pitchfork U_S = \{\pm x\}$ and $U_R\cup U_S$ does not contain other points of $R\cap S$; 
\smallskip

(2) There are:

\smallskip 

- a $2$-disk $\Dd$ in $\R^2$ with boundary $\partial \Dd$ with two corners
$a_1$, $a_2$ which is contained in the union of two smooth arcs say $\lambda_R$, $\lambda_S$ 
in $\R^2$ which intersect transversely at $\{a_1,a_2\}$; 
\smallskip

- an embedding $\psi: U \to M$
where $U$ is a closed $2$-disk in $\R^2$ containing  $\Dd\cup (\lambda_R \cup \lambda_S)$,
such that 
\begin{itemize}
\item $\psi(\lambda_*)\subset U_*, \  *=R, \ S$;
\item  $\psi(\partial \Dd,\{a_1,a_2\})=(\gamma,\{q_1,q_2\})$;
\item for every $x\in \lambda_*$, $d_x\psi(T_xU)\cap T_{\psi(x)}U_* = d_x\psi(T_x\lambda_*)$;  
\item $\psi({\rm Int}(\Dd)) \subset M\setminus (R\cup S) $.
\end{itemize}
We summarize (1) and (2) by saying that the smooth $2$-disk with corners $D:= \psi(\Dd)$ is 
{\it properly embedded} into $(M,R\cup S)$ and {\it connects the crossing points} $\pm x$.
Moreover, we require:  
\smallskip

(3) We can extend the embedding $\psi$ to a  parametrization of a neigbourhood
of $D$ in $M$ by a {\it standard model}, that is
to an embedding
$$\Psi: U\times \R^{r-1}\times \R^{s-1}\to M$$
such that $\Psi(\lambda_R\times \R^{r-1}\times \{0\})=U_R$ and $\Psi(\lambda_S\times\{0\} \times \R^{s-1})=U_S$.
\medskip

\begin{remark}{\rm We stress that the existence of a Whitney disk (in particular condition (3)) for  a situation 
$(M,R,S,x_0,x_1)$ implies that the two points are necessarily of opposite sign.} 
\end{remark}
 
{\it (Whitney trick)}  The Whitney trick applies to $(M,R,S,\pm x)$ at a Whitney disk connecting $\pm x$:
thanks to the standard model,
such a Whitney disk can be easily used as a guide to construct  an isotopy of $R$ in $M$  with support not intersecting the other points of $R\cap S$ and 
carrying $R$ to $R'\pitchfork S$ such that $R'\cap S = R\cap S \setminus \{\pm x\}$ (recall Figure \ref{WT} of Chapter \ref {TD-CUT-PASTE},
by renaming $R=P$, $S=Q$). 

\smallskip

\begin{definition}\label{WT(n,r)} {\rm For every type $(n,r)$ as above, we say that  $\WW \TT(n,r)$
holds if every situation $(M,R,S,\pm x)$ of type $(n,r)$ admits a Whitney disk.}
\end{definition}

\smallskip

We are going to relate the validity of  $\WW \TT(n,r)$ with a certain {\it unlinking property} of {\it unknotted spheres} into a sphere.

A smooth $p$-sphere $\Sigma \subset S^k$, $k>p\geq 1$, is {\it unknotted} if it is the boundary of a smooth $(p+1)$-disk embedded into
$S^k$.  The following lemma is easy, by using the unicity of disks up to diffeotopy.

\begin{lemma}\label{unk} Let $\Sigma \subset S^k$ be unknotted. Let $D$ be a smooth $k$-disk in $S^k$ disjoint
from $\Sigma$. Then $\Sigma$ is the boundary of  a smooth $(p+1)$-disk embedded into $S^k \setminus D$.
\end{lemma}

\cvd

\smallskip

A {\it link of unknotted spheres $(S^k,\Sigma,\Sigma')$ of type $(k,p)\in \N^2$}  consists of 
two {\it disjoint} unknotted smooth spheres  $\Sigma, \Sigma'  \subset S^k$ such that
 $$p=\dim \Sigma, \ q= \dim \Sigma', \ p\leq q, \  k = p+q +1 \ . $$
 
Such a link  $(S^k,\Sigma,\Sigma')$ is  (geometrically) {\it unlinked} if the two spheres are the boundary
of disjoint $(p+1)$- and $(q+1)$-disks respectively. By using again the unicity of disks up to diffeotopy, we have

\begin{lemma}\label{unk2} Up to diffeotopy there is a unique unlinked link of type $(k,p)$.
\end{lemma}

\cvd 

For every link $(S^k,\Sigma, \Sigma')$, give the spheres auxiliary orientations; 
then we can define the {\it linking number} (recall Section \ref{linking}
and Remarks \ref{3D-link} )
$$ lk(\Sigma, \Sigma') \in \Z \ . $$
A link is {\it algebraically unlinked} if 
$$lk(\Sigma, \Sigma')=0 \ . $$ 
We know (see the end of Section \ref{degree}) that the choice of auxiliary orientations is immaterial
with respect to the vanishing of the linking number; moreover this property is symmetric:
$lk(\Sigma, \Sigma')= 0$ if and only if $lk(\Sigma', \Sigma)=0$. 
Obviously, geometrically unlinked links are algebraically unlinked. 

\begin{definition}\label{U}{\rm For every link type $(k,p)\in \N^2$, we say that  the {\it unlinking property} $\UU(k,p)$ holds, if
every link (of unknotted spheres) $(S^k,\Sigma, \Sigma')$ of type $(k,p)$ which is algebraically unlinked is in fact geometrically unlinked.}
\end{definition}

\smallskip

It follows from the above discussion that Theorem \ref {strong-canc} will be a corollary of  item (1) in the next proposition. 

\begin{proposition}\label{Conc-inductive} (1) For every  type $(n,r)$ such that $n\geq 5$ and $r\geq 2$, $\WW\TT(n,r)$ holds.
\smallskip

(2) For every link type $(k,p)$ such that $k\geq 4$, $\UU(k,p)$ holds.
\end{proposition}
\Dim First let us prove that   $\UU(k,1)$ holds for every $k\geq 4$;
for  consider an algebraically unlinked link $(S^k, \Sigma, \Sigma')$, $\dim \Sigma =1$, $\dim \Sigma' = q= k-2\geq 2$.
Then $S^k \setminus \Sigma'$ is homotopically equivalent to the standard $S^1\subset S^k$ and the embedding of
$\Sigma$ in $S^k \setminus \Sigma'$ is homotopically trivial; as $k> 2\dim \Sigma + 1 =3$, then $\Sigma$ is isotopic
in $S^k\setminus \Sigma'$ onto a geometrically unlinked circle.
\smallskip

Next we prove the following claim.

\smallskip

{\bf Claim 1}  {\it For every $n\geq 5$, If  $\ \WW\TT(n,r)$ holds, then $\UU(n,\min(r, q))$, $q=n-r-1$, holds. }
\smallskip

{\it Proof of the claim:}    Consider an algebraically unlinked link $(S^n, \Sigma, \Sigma')$, $\dim \Sigma =r$, $\dim \Sigma' = q$.
Assume for simplicity that $r\leq q$. Let $D\subset S^n$ be a $(q+1)$-disk such that $\partial D = \Sigma'$.
Then the intersection number $[\Sigma]\bullet [D]$ in $S^n\setminus \Sigma'$ is equal to $0\in \Z$.
Then as $\WW\TT(n,r)$ holds, $\Sigma$ is isotopic to say $\Sigma"$ such that $\Sigma"\cap D= \emptyset$.
We can assume that $\Sigma"$  is embedded into $S^n \setminus B$ where $B\sim D^n$ is a $n$-disk of $S^n$
which thickens $D$. Then we conclude by means of Lemma \ref{unk}.

\smallskip

Next we propose two ways to conclude. The first way consists in a direct proof of item (1);
then item (2) will follow as a corollary of Claim 1 and the case $\UU(k,1)$ already achieved.
By the second way  both statements will be proved simultaneously by implementing the {\it concatenated inductive scheme}
obtained by combining Claim 1 with the following Claim 2 (the case $\UU(k,1)$  being the initial step of this induction):
\smallskip

{\bf Claim 2} {\it For every $k\geq 4$, if $\UU(k,p)$ holds, then $\WW\TT(k+1,p+1)$
holds.}

\smallskip

The second way makes fully manifest the strict relationship
between  $\WW\TT$ and $\UU$.  The presentation of this second way is very close to Chapter 5 of \cite{RS}.

\smallskip

{\it Proof of item (1):}  As $n\geq 5$, by general position we can assume that points (1) and (2) in the definition of a Whitney
disk  for $(M,R,S,\pm x)$ are fulfilled. It remains to achieve point (3). This is rephrased in terms of a suitable configuration
of subbundles of $T(M)$ over $(D,\partial D)$. We can assume that an auxiliary
Riemannian metric $g$ on $M$ is fixed in such a way that $R$ and $S$ are orthogonal at their intersection points, the normal bundles
and the associated tubular neighbourhoods are constructed by means of $g$. 
We use the notation $\nu_XY$ to mean the normal bundle of $Y$ in $X$.
The tangent bundle $T(R)$ splits over $\gamma_R$
as 
$$T(R)|\gamma_R = T(\gamma_R)\oplus E_R$$ 
where $E_R$ is a rank-$(r-1)$ subbundle of $(\nu_M D)|\gamma_R$.
Thus $E_R$ is tangent to $R$ and normal to $D$. 
The normal bundle $\nu_M S$ splits over $\gamma_S$ as 
$$(\nu_M S)|\gamma_S = \nu_D \gamma_S \oplus E_S$$
where $E_S$ is a rank-$(r-1)$ subbundle of $(\nu_M D)|\gamma_S$.
Thus $E_S$ is normal to both $S$ and $D$. 
$E_R$ and $E_S$ match at the intersection points
$\pm x$, so that we have a rank-$(r-1)$ bundle $E$ defined over the whole
$\partial D$. By construction $E$ is tangent to $R$ and normal to $S$.
We claim that $E$ can be extended to a subbundle of the whole $\nu_M D$.
By means of a trivialization of $\nu_M D$ we can encode $E$ as a map
$E: \partial D \to \GG_{n-2,r-1}$. Then $E$ extends if and only if it is homotopically
trivial. It is known that under our dimensional hypotheses (see for instance \cite {Steen})
$$\pi_1(\GG_{n-2,r-1}) = \Z/2\Z$$ 
and that $E$ as above is homotopically trivial if and only
if the corresponding rank-$(r-1)$ bundle is {\it orientable}. This is actually the case because
the intersection points have opposite signs. At this point it is not hard to build compatible
trivializations of the bundles considered so far and achieve point (3) in the definition
of a Whitney disk.

\smallskip

{\it A sketch of proof of Claim 2:} Let $(M,R,S,\pm x)$ of type $(k+1,p+1)$, $k\geq 4$. 
Argue as in the above proof of item (1), so that we can assume that points (1) and (2) in the definition of a Whitney
disk   for $(M,R,S,\pm x)$ are fulfilled. Again it remains to achieve point (3). Assume that it holds. 
We analyze the standard model and then we transport the conclusions in $M$ around the disk $D$ by means of the embedding $\Psi$.
Up to corner smoothing, $B:=U\times D^{p}\times D^{k-p-1}$ is a $k+1$-disk, and 
$$(\partial B, \partial (\lambda_R\times D^{p}\times \{0\}), \partial  (\lambda_S\times \{0\} \times D^{k-p-1}))$$
is diffeomorphic to an unlinked link of type $(k,p)$. Moreover, the whole $B$ can be recostructed from such an unlinked link.
Assume now that a priori only (1) and (2) are verified. We can nevertheless find a smooth  
$(k+1)$-disk $B$ in $M$ around $D$, which retracts to $D$, such that 
$$\partial B \pitchfork U_R:= \Sigma_R, \ \partial B \pitchfork U_S:= \Sigma_S$$  
are smooth spheres in the sphere $\partial B\sim S^k$ forming a link of type $(k,p)$.
In order to incorporate it in a standard model, by Lemma \ref{unk2} it is enough to prove that it is unlinked. 
As $\pm x$ have opposite signs, it follows that the link is algebraically unlinked, and we can conclude
because $\UU(k,p)$ holds  by the hypothesis of Claim 2.

\smallskip

The Proposition is proved.

\cvd

\smallskip

\begin{remark}\label{no-low}{\rm As for low dimensions we note that:  
\smallskip

 \begin{figure}[ht]
\begin{center}
 \includegraphics[width=5cm]{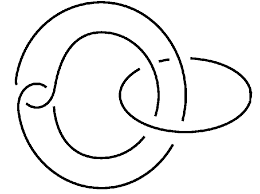}
\caption{\label{Wlink} Whitehead's link.} 
\end{center}
\end{figure}

- $\UU(3,1)$ fails. The simplest counterexample is the so called {\it Whitehead link}; 
several classical knot invariants show that it is geometrically linked in spite of the fact
that it is algebraically unlinked (see \cite{Rolf}).

\smallskip

- Trying to perform the construction in order to approach $\WW\TT(4,2)$,  it is not hard to realize item (1)
in the definition of Whitney disk; however (2) and even more (3) are very problematic - in fact we will see
in Chapter \ref{TD-4} that there are actual obstructions.}
\end{remark}  

\chapter{On $3$-manifolds}\label{TD-3}
In this chapter we will apply several results estasblished so far to
compact $3$-manifolds. We stress that we will develop a few themes
based on classical differential/topological tools, mainly on
transversality.  In no way we will touch Thurston's {\it
  geometrization} approach that has dominated the study of
$3$-manifolds in last decades. We will not even touch fundamental
results in $3$-dimensional geometric topology such as the
decomposition in prime manifolds or the so called $JSJ$-decomposition.
We will provide elementary and selfcontained proofs of the primary
fact that compact orientable boundaryless $3$-manifolds are parallelizable.
An important amount of the chapter will be devoted to several proofs
of ``$\Omega_3=0$" and of the equivalent Lickorish-Wallace theorem
about $3$-manifolds up to surgery equivalence respectively. Every
proof will illuminate different facets of the matter.  We will study
the behaviour of surfaces immersed or embedded in
a given $3$-manifold $M$, including  the determination of the bordism
group $\Ii_2(M)$ of immersed surfaces. 
An emerging theme will be the quadratic enhancement of the intersection
forms of surfaces immersed in  $3$-manifolds. This will occur
also in the classification of $3$-manifolds up to 
equivalence relations defined in terms of blowing up along smooth
centres.
 
 \section{Heegaard splitting}\label{HEE}
Let $M$ be a connected, {\it orientable}, boundaryless compact
$3$-manifold. We know that there is an ordered handle decomposition
$\Hh$ of $M$ with only one $0$-handle, only one $3$-handle, and such
that both $1$- and $2$-handles respectively are attached
simultaneously at disjoint attaching tubes. Denote by $M_1$ the
submanifolds with boundary of $M$ obtained by attaching the
$1$-handles at the boundary of the unique $0$-handle.  As $M$ is
orientable, then also $M_1$ is orientable; by the uniqueness of disks
up to diffeotopy applied to the attaching tubes of $1$-handles and
handle sliding, $M_1$ only depends up to diffeomorphism on the number
say $g\geq 0$ of $1$-handles and is called a {\it handlebody of genus
  $g$}, denoted by $\HG_g$.  Its boundary $\Sigma= \partial M_1$ is a
surface of genus $g$, that is diffeomorphic to the connected sum of
$g$ copies of the torus $S^1\times S^1$. If $g=1$, $\HG_1= D^2\times
S^1$ is also called a {\it solid torus}.  Consider the dual handle
decomposition $\tilde \Hh$, so that the $2$-handles of $\Hh$ become
the $1$-handles of $\tilde \Hh$. Apply the above discussion to $\tilde
M_1$. Then $\partial M_1 = \partial \tilde M_1=\Sigma$ and also
$\tilde M_1$ is a handlebody of genus $g$. Then $$M= M_1 \cup \tilde
M_1$$ is called a {\it Heegaard splitting} of $M$ of genus $g$ and the
separating surface $\Sigma$ is the corresponding {\it Heegaard
  surface}.

So every such an $M$ admits a Heegaard splitting of some genus and we
can define the {\it Heegaard genus} $g_H(M)$ of $M$ as the minimum $g$
such that $M$ has a splitting of genus $g$.  As it often happens such
an invariant is easy to define but in general hard to compute or even
to estimate.

Up to diffeomorphism, a Heegaard splitting of $M$ can be described
equivalently as follows: fix a standard model $\HG_g$ of genus $g$
handlebody (for instance embedded in $\R^3$ and endowed with the
standard induced orientation); let $\Sigma_g = \partial \HG_g$ with
the boundary orientation. Fix an auxiliary smooth automorphism
$\gamma$ of $\Sigma_g$ which  reverses the orientation.
Then there is an orientation {\it preserving}
(say ``positive")
smooth automorphism $\phi \in {\rm Diff}^+(\Sigma_g)$ such that
$$ M=M_1\cup\tilde M_1 \sim \HG_g \amalg_{\gamma \circ \phi} \HG_g \ . $$ 
Moreover, we know that up to diffeomorphism the last term only depends on the
isotopy class of $\phi$; in other words, define
$${\rm Mod}(\Sigma_g):= {\rm Diff}^+(\Sigma_g)/{\rm
  Diff}^0(\Sigma_g)$$ that is the quotient group mod the normal
subgroup of automorphisms isotopic to the identity. This is called the
{\it mapping class group} of $\Sigma_g$ (also called its {\it modular
  group}) and is an object of main importance and interest. Then every
splitting is of the form
$$M \sim \HG_g \amalg_{[\phi]}\HG_g \sim \HG\amalg_{\gamma \circ \phi} \HG_g,  \  
[\phi]\in {\rm Mod}(\Sigma_g) \ . $$

\begin{example}{\rm (1)
    If $g_H(M)=0$ then $M$ is a twisted hence a true smooth
    $3$-sphere.

(2) The $3$-manifolds such that $g_H(M)=1$ are classified and called
    {\it lens spaces} \cite{Brod}.  Let us recall the main
    facts. Realize the torus as the quotient manifold $\R^2/\Z^2 $.
    The matrix group $SL(2,\Z)$ acts linearly on $\R^2$ by preserving
    the lattice $\Z^2$. Then the action descends to the quotient.  In
    fact one can prove that
$$ {\rm Mod}(\Sigma_1)\sim SL(2,\Z) \ . $$ Fix an identification of
    the torus as the boundary $\Sigma_1$ of $\HG_1$ in such a way that
    the circle image in $\R^2/\Z^2$ of the $x$-axis of $\R^2$ becomes
    a {\it meridian} $m$ that is it bounds a $2$-disk properly
    embedded into $(\HG_1,\Sigma_1)$, while the image of the $y$-axis
    is a longitude $l$ which intersects transversely $m$ at one point;
    $m$, $l$ form a basis of $\Omega_1(\Sigma_1)\sim \Z^2$. Let $A\in
    SL(2,\Z)$, so that $A(m)=pm+ql$, $gcd(p,q)=1$. Denote by $L_{p,q}$
    the resulting lens space obtained by using $A$ as gluing map. It
    is not hard to check via Van Kampen theorem that
    $\pi_1(L_{p,q})\sim \Z/p\Z$.  Then $L(p,q)$ is diffeomorphic to
    $L(p,q')$ if and only if
$$\pm q'= q^{\pm1} \ {\rm mod} (p) \ . $$
For higher genus  the situation is much more complicated.}
\end{example} 

\subsection{Heegaard diagrams and a diagramatic ``calculus"}\label{H-diag}
Heegaard splittings can be encoded by means of suitable Heegaard
diagrams.
\begin{definition}\label{H-diag2}
  {\rm A {\it genus $g$ Heegaard diagram} consits of a triple
    $(\Sigma, C^-, C^+)$ where
\begin{enumerate}
\item $\Sigma$ is a surface of genus $g$; 
\item
  $C^{\pm}= \{c^\pm_1,\dots, c^\pm_g\}$ is a family of $g$ disjoint
  simple smooth circles on $\Sigma$ whose union does not divide $\Sigma$,
  that is by removing from
  $\Sigma$ the interiors of small pairwise disjoint annular
  neighbourhoods of these circles we get a $2$-sphere with $2g$ holes;
\item
  $C^-\pitchfork C^+$, that is the union of the $c^-_j$'s is
  transverse to the union of the $c^+_j$'s.
\end{enumerate}
}
\end{definition}

Given a Heegaard diagram we can construct a $3$-manifold $M$ endowed
with an Heegaard splitting as follows.  Take the product $\Sigma
\times [-1,1]$ and stipulate that the circle of $C^\pm$ are traced
on $$\Sigma \times \{\pm 1\}:= \Sigma^\pm\ . $$ $\Sigma$ is identified
with the separating surface $\Sigma \times \{0\}$.  Then take a system
of pairwise disjoint annular neighbourhoods say $T_j^\pm$ for $C^\pm$
on $\Sigma^\pm$.  Consider the $T^+_j$ as a system of attaching tubes
of disjoint $2$-handle attached to $\Sigma\times [0,1]$ at
$\Sigma^+$. Thanks to the properties of the circles in $C^+$ this
produces a $3$-manifold with boundary diffeomorphic to $\Sigma \amalg
S^2$.  By filling the spherical component with a $3$-handle we get the
piece $\tilde M_1$ of the desired handle decomposition of $M$. Doing
similarly on the other side $\Sigma \times [-1,0]$ we get the piece
$M_1$ and eventually the splitting
$$M\sim M_1 \cup \tilde M_1$$ with Heegaard surface $\Sigma$. 

\begin{remark}\label{use-Smale}{\rm The fact that the resulting 
$3$-manifold is {\it unique up to diffeomorphim} follows from Smale
theorem recalled in Proposition \ref{on-twisted}, (1), $m=3$.}
\end{remark}

On the
other hand, every Heegaard splitting with Heegaard surface $\Sigma$
gives rise to an encoding Heegaard diagram, possibly by handle sliding
in order to reach the transversality requirement of the definition.

\smallskip

{\it (Heegaard diagram moves)} The elementary handle moves induce
elementary moves on Heegaard diagrams which keep the resulting
manifold $M$ fixed up to diffeomorphism.
\begin{itemize}
\item
  Handle sliding produces the following diagram moves (called {\it
    $H$-diagram sliding}):

1) of course we can modify $C^\pm$ up to ambient isotopy (keeping that
$C^+ \pitchfork C^-$);

2) more substantially we have: let $T^\pm_j$ and $T^\pm_i$ be disjoint
annular neighbourhoods of two circles of $C^\pm$ as above.  Connect
these annuli by attaching an embedded $1$-handle $H$ at $\partial
(T^\pm_j \amalg T^\pm_i)$ in such a way that apart the attaching
segments, $H$ is contained in $\Sigma \setminus \cup_{s=1}^g
T^\pm_s$. The boundary of $T^\pm_j\cup T^\pm_i \cup H$ contains a
component say $c'_j$ which is the embedded connected sum of a parallel
copy of $c^\pm_j$ with a parallel copy of $c^\pm_i$.  Then get a new
$C^\pm$ just by replacing $c^\pm_j$ with $c'_j$.

\item
  Cancellation/introduction of a pair of complementary handles
  produces the following diagram move.  Consider the
  diagram $$(S^1\times S^1, c^-=S^1\times \{y_0\}, c^+=\{x_0\} \times
  S^1) \ . $$ Given any diagram $(\Sigma,C^-,C^+)$ of genus $g$,
  replace $\Sigma$ with $\Sigma \# (S^1\times S^1)$ provided that the
  sum is performed at $2$-disks disjoint from $C^-\cup C^+$ and
  $c^-\cup c^+$ respectively; then add to $C^\pm$ the circle $c^\pm$
  to get the new diagram of genus $g+1$. In terms of the resulting
  $3$-manifolds we replace $M$ with $M \# S^3 \sim M$. This move is
  called {\it elementary stabilization}.

\end{itemize}

The stabilization shows by the way that for every $g\geq g_H(M)$, $M$
admits Heegaard splitting of genus $g$. In particular
$S^3$ has splittings of every genus. One can prove (see \cite{Sing}):

\begin{theorem}
  Two Heegaard diagrams encode Heegaard splittings of a same
  $3$-manifold $M$ (considered up to diffeomorphism) if and only if
  they become equal up to finite sequences of $H$-diagram sliding or
  stabilizations.
\end{theorem}

\cvd

\smallskip
\begin{remark}{\rm Once the existence of Heegaard splitting has been easily
established, several non trivial questions naturally arise such as:

- For a given $M$, estimate  in effective terms its genus $g_H(M)$;

- For every $g\geq g_H(M)$, study the Heegaard splittings of $M$ of
genus $g$ up to ambient isotopy.

Concerning the second question a complete answer is known for the $3$-sphere and lens spaces
defined above, that is for manifolds such that  $g_H\leq 1$; we have: 

\smallskip

{\it For every $g\geq 1$, $S^3$ and every lens space have up to diffeotopy a unique Heegaard
splitting of genus $g$.}

\smallskip

On the other hand, for $g\geq 2$, there are manifolds with non isotopic splittings of genus $g$.

We refer to the body and the references of \cite{BO} for more information about this question. }
\end{remark} 

\subsection{From Heegaard diagrams to spines and $\Delta$-complexes}\label{spine}
The aim of this section, mainly of technical nature, is to show other ways
to present $3$-manifolds derived from Heegaard splittings.
We refer to \cite{BP} for a wide treatment of the topic touched in this
section.  We will use some of this facts in Section \ref {3-parallel}.

Let $(\Sigma, C^-, C^+)$ be a Heegaard diagram of $M$ as above.  Up to
$H$-sliding, we can assume that not only $C^-\pitchfork C^+$, also that
every component (called a {\it region}) of $\Sigma \setminus (C^-\cup
C^+)$ is a open $2$-disk. By following the reconstruction of the
Heegaard splitting
$$M\sim M_1 \cup \tilde M_1$$ of $M$ encoded by the diagram, we see
that the core of every $2$-handle attached to a circle $c$ of
$C^+\times \{1\}$ can be extented by means of the annulus $c\times
[0,1]$ and we get an embedded $2$-disk in $\tilde M_1$ which
intersects tranversely $\Sigma = \Sigma \times \{0\}$ at $c$. Do it
for every $c$ in $\C^+$ and similarly for every $c$ in $\C^-$ getting
a disk in $M_1$.  Denote by $\PP$ the union of $\Sigma$ with all such
disks. $\PP$ is a kind of singular surface embedded into $M$ with the
following properties:
\begin{enumerate}
\item $S(\PP):= (C^-\cup C^+) \subset \Sigma$ is the singular locus of
  $\PP$;
\item $V(\PP):= C^-\cap C^+$ is the singular locus of $S(P)$, its
  points are the {\it vertices} of $\PP$. The components, each
  diffeomorphic to the open $1$-disk $(-1,1)$, of $S(\PP)\setminus
  V(\PP)$ are the {\it edges} of $\PP$; at every vertex there are four
  edge germs.
\item The components, each diffeomorphic to an open $2$-disk, of
  $\PP\setminus S(\PP)$ are the {\it regions} of $\PP$. Along every
  edge there are three region germs. At every vertex there are six
  region germs.
\item If $B^+$ and $B^-$ are the $0$ and $3$-handles of the splitting,
  then $\PP$ is a retract by deformation of
$$\hat M := M\setminus ({\rm Int}(B^-)\cup {\rm Int}(B^+) \ . $$ In
  fact there is a {\it normal retraction} $r: \hat M \to \PP$ such
  that: the fibre over a region point is diffeomorphic to $[-1,1]$;
  the fibre over an edge point is a {\it tripode} that is the wedge of
  three segments $[0, 1]$ with common endpoint $0$; the fibre over a
  vertex is a wedge of four such segments $[0,1]$; $\hat M$ can be reconstructed
  as being the mapping cilynder of such normal retraction.
\end{enumerate}

We summarize all this by saying that $\PP$ is a {\it standard spine}
of $\hat M$.  By using the language of $CW$-complexes, $\PP$ is the
$2$-skeleton of such a complex over $M$ which is obtained by attaching
two $3$-cells to it. 

Now we give $\PP$ an additional structure called a {\it branching}.
Give $\Sigma$, hence every region of $\PP$ contained in $\Sigma$, an
orientation; give every circle $c$ in $C^-\cup C^+$ an orientation,
hence give the region of $\PP$ bounded by $c$ the orientation with the
prescribed boundary orientation. In this way $S(\PP)$ is union of
oriented circles crossing transversely on $\Sigma$ at some vertices;
every region of $\PP$ is oriented in such a way there is a {\it
  prevailing orientation} induced on every edge of $\PP$ and this
agrees with the one of the circle $c$ in $S(\PP)$ which contains the
edge. Notice that at every vertex the four configurations at the edge
germes automatically match.  We call this system of region
orientations a {\it branching} $\bb$ of $\PP$ and we summarize by
saying that the standard spine $\PP$ has be enhanced to be a {\it
  branched standard spine} $(\PP,\bb)$.  The terminology is justified
because the branching encodes a way to convert $\PP$ to be a
(oriented) {\it branched surface} embedded in $M$. This means that
$\PP$ can be moved in $M$ is such a way that, although being singular,
nevertheless it is well defined everywhere on $\PP$ a smooth field of
oriented tangent $2$-planes. In our specific situation, we can keep
$\Sigma$ fixed and isotopically move every other region $R$ bounding
some circle $c$ to becomes tangent to $\Sigma$ along $c$ over the side
of $c$ in $\Sigma$ which carries together with $R$ the prevailing
boundary orientation.

The branched spine $(\PP,\bb)$ can be considered as the $2$-skeleton
of the {\it dual cell decomposition} to a $\Delta$-{\it complex structure}
over $M$ in the sense of \cite{Hatch}. This is a kind of {\it
  triangulation} of $M$ obtained as follows. Select one base point in
each edge, region of $\PP$ and in the interior of the $3$-balls
$B^\pm$.

Recall that the standard $j$-simplex $\Delta^j$ is contained in the
affine hyperplane $\{x_0+x_1+\dots + x_j=1\}$ of $\R^{j+1}$ and is the
convex hull with ordered vertices of the vectors $e_0, e_1, e_2,
\dots, e_j$ of the standard basis.  Every $h$-face of $\Delta^j$,
$h=0,1,2,\dots, j$, is the $h$-simplex with $h+1$ vertices obtained by
omitting $j+1-(h+1)$ vertices of $\Delta^j$. For every such a $h$-face
$F$, there is a {\it canonical affine parametrization}
$$\phi_F: \Delta^h \to F$$ defined on the standard $h$-simplex and
preserving the vertex ordering.  A {\it singular $j$-simplex in $M$}
is a continuous map $\sigma: \Delta^j \to M$. For every $h$-face $F$
of $\Delta^j$,
$$(\Delta^h, \sigma \circ \phi_F)$$ 
is the corresponding {\it singular face} of the singular simplex.

Then we can associate to every $x\in V(\PP)$ a ``dual'' singular
$3$-simplex $(\Delta^3,\sigma_x)$ in such a way that the following
properties are verified
  
\begin {enumerate}
\item For every $x$, the restriction of $\sigma_x$ to the {\it
  interior} of every $h$-face of $\Delta$, $h=0,1,2,3$, is a smooth
  embedding into $M$.
\item For every vertex $x$ of $\PP$, the image by $\sigma_x$ of every
  vertex of $\Delta^3$ is one of the base points of $B^\pm$; $x$
  belongs to the image of the interior of $\Delta^3$; the image of
  every open edge of $\Delta^3$ is transverse to one dual region of
  $\PP$ which has $x$ in its closure, exactly at the region base
  point; the image of every open $2$-face of $\Delta^3$ is trasverse
  to one edge of $\PP$ which has $x$ in its closure, exactly at the
  edge base point.
\item Giving every image of an open edge of $\Delta^3$ the orientation
  dual to the $\bb$-orientation of the dual region, and the edge
  itself the orientation determined by the vertex order of $\Delta^3$,
  then the embedding of the open edge by $\sigma_x$ is orientation
  preserving.
\item If the image of two singular open $h$-faces by some $\sigma_x$,
  $\sigma_{x'}$ (possibly $x=x'$) share the same dual $(3-h)$-cell of
  $\PP$, then the whole singular faces coincide.
\item Varying $x$ in $V(\PP)$, the images of the several open
  $h$-faces form a {\it partition} of $M$.
\item Up to a piecewise smooth homeomorphism, $M$ is obtained by
  gluing the abstract $3$-simplices associated to the vertices of
  $\PP$ at common singular faces.
 \end{enumerate}
 
 \smallskip
 
 We can modify a branched standard spine $(\PP,\bb)$ of $\hat M$,
 associated as above to a Heegaard diagram of $M$, to become a
 branched standard spine $(\PP_0, \bb_0)$ of $M_0$, where $M_0$ is of
 the form
 $$M_0 = M \setminus {\rm Int} (B)$$ where $B$ is some smooth $3$-ball
 in $M$. So in particular $\PP_0$ will be the $2$-skeleton of a
 $CW$-complex over $M$ with a unique $3$-cell. Do it as follows. Take
 a point $p$ on an edge of $\PP$ and locally insert an embedded
 triangle $T$, whose interior is contained in $M\setminus \PP$, $p$ is
 a vertex of $T$, $T$ intersects transversely $\PP$ at two edges with
 $p$ as common endpoint, contained respectively into germs of regions
 both inducing the prevailing orientation on the edge. Then $T$ has a
 ``free'' edge $l$. Attach an embedded $1$-handle with core parallel
 to $l$, intersecting transversely $T$ along its $b$-tube at $l$, with
 attaching tube on $\PP$.  Then $\PP_0$ results from $\PP$ by such a
 surgery. It is easy to see that the handle has fused the two
 components of $\partial \hat M$ into one spherical boundary component
 of a $M_0$ of the desired form. By construction $\PP_0$ is a standard
 spine of $M_0$ and it carries a branching $\bb_0$ which agree with
 $\bb$ on the regions that have not be effected by the surgery. The
 above considerations apply to $(\PP_0,\bb_0)$ as well, so that we
 have a dual $\Delta$-complex structure on $M$ with only one singular
 $0$-simplex.
 
 We know that $M$ is combable. Here we construct a nowhere
 vanishing tangent vector field by means of $(\PP_0,\bb_0)$.
 The tangent oriented $2$-planes distribution along the branched
 surface $\PP_0$, has an orthogonal distribution of unitary tangent
 vector (with respect to an auxiliary riemann metric on $M$). This can
 be extended to a {\it generic traversing} unitary tangent vector
 field $v_0$ on $M_0$. This means that the following properties hold:
 \begin{enumerate}
 \item Every integral line of $v_0$ is a segment with endpoints on
   $\partial M_0$.
 \item $v_0$ is simply tangent to $\partial M_0$ at the disjoint union
   $\Ss$ of some smooth circles. For every $y\in \Ss$, the integral
   line passing through $y$ is tangent to $\partial M_0$ and trasverse
   to $\Ss$.
 \item Generic integral lines are not tangent to $\partial M_0$;
   generic tangent integral lines are tangent to $\partial M_0$ at one
   point; a finite number of exceptional integral lines is tangent at
   exactly two points.
 \end{enumerate}
 
 We can assume that the image of every singular edge of the
 $\Delta$-complex structure dual to $(\PP_0,\bb_0)$ intersects $M_0$
 at the integral line of $v_0$ though the base point of dual region
 and that this line is not tangent to $\partial M_0$.
\smallskip
 
 We have
 \begin{proposition}\label{v0-extend}
   $v_0$ extends to a unitary tangent vector field $v$ 
defined on the whole of $M$.
\end{proposition}
\Dim We can assume that $B$ is in a chart of $M$ and that the
auxiliary metric looks standard in that coordinates. So the
restriction of $v_0$ to $\partial M_0\sim S^2$ is encoded by a map $h:
S^2\to S^2$ and can be extended over $B$ if and only if its degree
vanishes.  Assume that $M_0$ is endowed with a framing (we will see later that
this is always true), then the whole $v_0$ can be encoded by a map
 $$H: M_0 \to S^2$$ which extends $h$. Usual invariance of the degree
up to bordism shows that the degree of $h$ vanishes indeed. 

\cvd

 \subsection{Non orientable Heegaard splitting}\label{non_or_H}
If $M$ is compact connected boundaryless and {\it non} orientable,
then by using a nice handle decomposition as above we see that
$$M\sim M_1\cup \tilde M_1$$ where $M_1$ is non orientable and is
obtained by attaching say $h+1$ disjoint $2$-handles to the unique
$0$-handle at the boundary $\partial D^3 = S^2$ (and similarly for
$\tilde M_1$ with respect to the dual decomposition).  Up to handle
sliding, we can assume that only one of these $2$-handles destroyed
the orientability and that $M_1\sim \tilde M_1$ only depend (up to
diffeomorphism) to the nunber $h+1$. Let us call it a {\it non
  orientable handlebody of genus $h$}.  The separating (non
orientable) Heegaard surface is now diffeomorphic to
$$ \tilde \Sigma_h:= (\PP^2(\R) \# \PP^2(\R)) \# h(S^1\times S^1)
\ . $$ The readear would imagine how to develop a non orientable
version of Heegaard diagrams and diagram moves supported by such
surfaces. Stabilization extends verbatim; a bit of care is necessary
for the sliding diagram moves.

\section{Surgery equivalence}\label{surgery}
We define a ``surgery'' equivalence relation on compact connected
boundaryless $3$-manifolds in terms of certain special $4$-dimensional
triads; the main application will be a characterization of $3$-dimensional
orientable
boundary as the manifolds which are surgery equivalent to the sphere
$S^3$.
   
\begin{definition}\label{surg-def}{\rm 
Let $M_0$ and $M_1$ be compact connected boundaryless non empty
$3$-manifolds. We say that $M_1$ can be obtained by {\it
  (longitudinal) surgery (along a framed link) of $M_0$} (and we write
$M_1 \sim_\sigma M_0$) if there exists a $4$-dimensional triads
$(W,M_0,M_1)$ which admits a handle decomposition $\Hh$ consisting
only of $2$-handles attached simultaneously at disjoint attaching
tubes.}
\end{definition}
\smallskip

To justify the terminology let us analyze the situation of the above
definition. The decomposition is of the form
$$C_0 \cup (\cup_{j=1}^d H^2_j)\cup C_1$$ where $C_0=M_0\times [0,1]$
and $C_1=[-1,0]$ are respective collars of $M_0$ and $M_1$ in $W$.
The union of the embedded attaching spheres of the $2$-handles
$$L= \cup_{j=1}^d K_s$$ is a so called {\it link} in $M_0\sim
M_0\times \{1\}$. Every component $K_s$ is a {\it knot} in $M_0$.
Moreover, we have a family of disjoint attaching
tubes $T_s$ each one equipped with a trivialization (also called a
``framing'') by $S^1\times D^2$, so that $K_s \sim S^1\times
\{0\}$. $M_1$ is obtained from $M_0$ by removing the interior of these
attaching tubes and attaching back a copy of $D^2\times S^1$ to every
boundary component $\partial T_s$, in such a way that a meridian $S^1
\times \{x_0\}$ of $D^2\times S^1$ is mapped onto a
longitude $\l_s \sim S^1\times \{y_0\}$, $y_0\in \partial D^2$ of
$K_s$ determined by the framing (such a longitude is unique up to isotopy). 

This defines an equivalence relation; in particular $M_1\sim_\sigma
M_0$ implies $M_0 \sim_\sigma M_1$ because the dual decomposition of
such an $\Hh$ also consists of $2$-handles only.  If $M_0 \sim_\sigma
M_1$, then $M_0$ is orientable if and only if $M_1$ is orientable and
in such a case any special triad connecting them is necessarily orientable.

Let us restrict for a while to
the orientable case. We have (see \cite{ Wa})

\begin{proposition}\label{wallace0}
  Let $M_0$, $M_1$ be compact connected orientable boundaryless
  $3$-manifold.  Then $M_1\sim_\sigma M_0$ if and only if there is an orientable
  $4$-dimensional triad $(W, M_0,M_1)$; that is for suitable orientations,
  $[M_0]=[M_1]\in \Omega_3$.
\end{proposition}

\begin{corollary}\label{wallace}
  $M\sim_\sigma S^3$ if and only if for every orientation of
  $M$, $[M]=0 \in \Omega_3$.
\end{corollary}

{\it Proofs:} Let us prove the corollary, assuming the proposition. If
$M_1\sim_\sigma S^3$, then by completing with one $4$-handle attached
at $S^3$ the dual $\Hh^*$ of a special decomposition $\Hh$ of a given
triad $(W,S^3, M)$, we get a triad $(V,M,\emptyset)$ so that
$M=\partial V$.  On the other way round, assume that $M=\partial V$
for some orientable connected $4$-manifold $V$. Then the triad
$(V,\emptyset, M)$ admits an ordered handle decomposition with one
$0$-handle, and no $4$-handles. By removing the $0$-handle we get an
orientable triad $(W,S^3,M)$ and we conclude by applying to it the
proposition.

Let us prove now the proposition. One implication is trivial. On the
other hand, let us start with any orientable triad $(W,M_0,M_1)$. It
has an ordered handle decomposition without both $0$ and
$4$-handles. Moreover, we can assume that all handles of a given index
are attached simultaneously at disjoint attaching tubes.  The idea is
to {\it trade} first every $1$-handle for a $2$-handle in such a way
that the $4$-manifold $W$ possibly changes but its boundary is kept
fixed. Every $1$-handle does not destroy the orientability.  Moreover,
by the uniqueness of disks up to diffeotopy we can assume that all
attaching tubes of the say $d$ $1$-handles are contained in a smooth
$3$-disk $D$ in $M_0\sim M_0\times \{1\}$; then after having attached
the $1$-handles to $C_0=M_0\times [0,1]$ at $M_0\sim M_0\times \{1\}$,
we get a $4$-manifold $W_1$ such that $\partial W_1$ is the connected
sum of $M_0$ with $d$ copies of $S^2\times S^1$.  A $4$-manifold $V_1$
with the same boundary can be obtained by surgery along a link $L$ in
$M_0$ formed by $d$ unknotted and unlinked components contained in the
above disk $D$, such that each component $K_s$ is endowed with the
framing associated to the distinguished longitude carried by a collar
in a $2$-disk $D_s$ in $D$ such that $\partial D_s = K_s$.  The rest
of the handle decomposition is unchanged and we get a $4$-dimensional
triad $(W',M_0,M_1)$ having an ordered handle decomposition $\Hh'$
without $0$, $1$ and $4$-handles.  In order to trade also the
$3$-handles for some $2$-handles, we manage similarly by using the
dual decomposition of $\Hh'$. Similarly as above we eventually get a
triad $(W",M_1,M_0)$ with a handle decomposition $\Hh"$ consisting
only of $2$-handles. The proposition is proved.

\cvd

Now we state two main theorems of this chapter.

\begin{theorem}\label{LW} {\rm  (Lickorish-Wallace)} Every orientable
  connected compact boundaryless $3$-manifold $M$ is surgery
  equivalent to $S^3$ ($M\sim_\sigma S^3$).
\end{theorem}

\begin{theorem}\label{omega_3=0} $\Omega_3=0$.
\end{theorem}

\smallskip

By Corollary \ref{wallace} they can be considered as a corollary of
each other.  This actually happened. For example Lickorish proved
Theorem \ref{LW} as an application of his main results about the
generators of the mapping class groups of surfaces, and by the way he
got a (new) proof that $\Omega_3 = 0$. On the contrary, Wallace
obtained the result via the above Corollary \ref{wallace}, as it was
already known (by several different proofs) that $\Omega_3=0$. We will
develop diffusely this theme.

\subsection{Non orientable surgery}\label{non-or-surg}
There is a non orientable version of Corollary \ref{wallace}.  Denote
by $\MG$ the non orientable $3$-manifold which is the boundary of the
non orientable $4$-manifold $\VG$ (unique up to diffeomorphism) with a
handle decomposition consisting of one $0$-handle and one $1$-handle.
In fact $\MG$ is the non orientable total space of a fibration over
$S^1$ with fibre $S^2$. Then we have (the proof is similar to the orientable case):
 
\begin{proposition}\label{wallace2} Let $M$ be a compact connected
  boundaryless non orientable $3$-manifold.
 Then  $M\sim_\sigma \MG$ if and only if $[M]=0\in \eta_3$.
 \end{proposition}

 \cvd  

\section{Proofs of  $\Omega_3 = 0$}\label{Omega3}
In this section we discuss a few ``direct'' proofs of Theorem \ref {omega_3=0}, 
so that Theorem \ref {LW} will result as a corollary.
\smallskip

$\bullet$ {\it (Via immersions in $\R^5$ and Seifert's surfaces)} This
is the first proof of $\Omega_3=0$ (Rohlin 1950, see his papers
translated in \cite{GM}).  If a compact connected orientable
boundaryless $3$-manifold $\hat M$ is {\it embedded} in $\R^5$, then
by Proposition \ref{seifert-surf} it admits an orientable Seifert's
surface $W$ so that in particular $\hat M=\partial W$.

\begin{remark}\label{on_seifert}{\rm
    Rohlin used a different argument to show the existence of
    Seifert's surfaces based on the estension of a combinatorial
    method due to Kneser to desingularize embedded simplicial cycles
    in triangulated manifolds to the codimension $2$ oriented and
    relative case (see \cite{GM} for an exhaustive discussion of this
    point).}
\end{remark}

In order to prove the theorem it is enough to show that for every
orientable $M$ there is an orientable triad $(V,M,\hat M)$ such that
$\hat M$ is embedded into $\R^5$. It was known since \cite{Whit3}
(1944) (recall Section \ref {n-in-2n-1}) that for every such an $M$
there is a generic immersion $f: M\to \R^5$; this also follows from
Smale-Hirsch immersion theory because we will see in Section
\ref{3-parallel} that $M$ is parallelizable.  We can conclude by
applying the ``embedding up to surgery" of Section \ref
{n-in-2n-1-emb}.

\smallskip

$\bullet$ {\it (Via vanishing of characteristic numbers)} In a sense
the most ``modern'' proof (being a special case of a general
determination of bordism groups based on Thom's spaces and
characteristic numbers) is the one obtained by applying Proposition
\ref{parallelizable}, as we will see in Section  \ref{3-parallel}
that orientable
$3$-manifolds are parallelizable.

\section{Proofs of Lickorish-Wallace theorem}\label{LW-proof}
In this section we discuss a few ``direct'' proofs of Theorem \ref
{LW}, so that Theorem \ref {omega_3=0} will result as a corollary.

These proofs are based on Heegaard splittings.

\smallskip

{(\it Via Dehn twists)} This is original Lickorish's proof
\cite{Lick}. A main Lickorish result establishes a distinguished set
of generators of the mapping class group Mod$(\Sigma_g)$.  Let $C$ be
a smooth circle on the surface $\Sigma_g$. Assume that $C$ is {\it
  essential} that is it is not the boundary of a smooth disk embedded
into $\Sigma_g$. Fix an auxiliary trivialization $$\psi: S^1\times
[-1,2] \to U$$ of a tubular neighbourhood of $C$.  Give $S^1\times
[-1,2]$ the coordinates $(e^{i\theta},t)$, $\theta \in [0,2\pi]$. Let
$\rho:[-1,2]\to [0,1]$ be a smooth non decreasing function such that
the restriction to $[0,1]$ is a diffeomorphism onto the image, it is
constantly equal to $0$ on $[-1,0]$, constantly equal to $1$ on
$[1,2]$. Then define the diffeomorphism
$$\tau_C: \Sigma_g \to \Sigma_g$$ which is the identity ourside $U$,
and is defined on $U$ as $\psi \circ h \circ \psi^{-1}$, where
$$h(e^{i\theta},t)= (e^{i(\theta+2\pi\rho(t))},t) \ . $$ $\tau_C$ and
$\tau_C^{-1}$ are called Dehn's twists along $C$.  Their classes in
Mod$(\Sigma_g)$ do not depend on the arbitrary choices we made,
including the fact that $C$ is considered up to ambient isotopy. Let
us call Dehn's twists also these classes.  Then we have:
\begin{theorem}
  {\rm Mod}$(\Sigma_g)$ is generated by the Dehn twists along essential
  smooth circles.
\end{theorem} 

\cvd

In fact the result is more precise because it shows that a determined
finite set of twists suffices. Anyway, we assume this theorem and we
show how to deduce that $M\sim_\sigma S^3$.

\begin{lemma}\label{tunnel}
  Let $[\psi]=[\tau_k]\circ \cdots \circ [\tau_1]$ be an element of
  {\rm Mod}$(\Sigma_g)$ expressed as composition of $k$ Dehn's twists. Then
  there exist two systems of $k$ disjoint solid tori $V_1, \dots, V_k$
  and $V'_1, \dots, V'_k$ in the interior of the handlebody $\HG_g$
  such that $\psi$ extends to a diffeomorphism
  $$ \bar \psi: \HG_g \setminus \cup_j {\rm Int}(V_j) \to
  \HG_g \setminus \cup_j {\rm Int}(V'_j) \ . $$
\end{lemma}
\Dim If $k=0$, then $\psi$ is isotopic to
the identity and the statement is trivially verified.  Assume that
$k=1$, $\psi= \tau=\tau_C^{\pm 1}$. Consider a collar $C(\Sigma_g)\sim
\Sigma_g \times [0,1]$ of $\Sigma_g=\partial \HG_g$ in $\HG_g$. Set
$V\sim  U(C)\times [1/2,1]\subset C(\Sigma)$ (up to corner smoothing)
where $U(C)$ is a annular neighbourhood of $C$ in $\Sigma_g$.  Set
$V'=V$. Then an extension of $\tau$ is obtained by taking a parallel
copy of $\tau$ on every leaf $U(C)\times \{s\}$, $0\leq s \leq 1/2$,
and setting $\bar \tau$ equal to the identity on the remaining part
of $\HG_1\setminus {\rm Int}(V)$.  If $k=2$ we can extend $\tau_2$
along $C_2$ by the same method, provided that the ``tunnel" $V_2$ is
more deeply in the interior of $\HG_g$ so that $V_1 \cap
V_2=\emptyset$ and $\bar \tau_1={\rm id}$ along $V_2$. Then set
$V'_2=V_2$, $V'_1=\bar \tau_2(V_1)$, so that $\bar \tau_2 \circ \bar
\tau_1$ is a desired estension of $\psi$. By iterating the same
method, by induction we get the resul for every $k\geq 0$.

\cvd

\smallskip

Consider any genus $g$ Heegaard splitting presented as above in the
form
$$ M \sim \HG_g \amalg_{[\phi]} \HG_g, \ \ [\phi]\in {\rm
  Mod}(\Sigma_g) \ . $$
  We know that also $S^3$ admits a genus $g$ splitting say  
$$S^3= \HG_g \amalg _{[\phi']} \HG_g \ . $$
Set $\psi= \phi^{-1}\circ \phi' = (\phi^{-1}\circ \gamma^{-1})\circ (\gamma \circ \phi')$. 
Apply the above lemma
to $\psi$. Then we get an extension
$$\bar \psi: \HG_g \setminus \cup_j{\rm Int}(V_j)\to 
\HG_g \setminus \cup_j{\rm Int}(V'_j)$$
which by construction extends to a diffeomorphism
$$\bar \psi: S^3 \setminus \cup_j{\rm Int}(V_j)\to 
M \setminus \cup_j{\rm Int}(V'_j) $$
and this readily shows that $M\sim_\sigma S^3$.

\cvd

\smallskip

{\it (By induction on a Heegaard diagram complexity)} Last but not least,
we present the clever proof of \cite{Rourke}. 
Let us fix an orientation of $M$; it is understood that
all manifolds produced by the following construction are oriented and
that the orientations are compatible.  Actually we are going to
realize that $S^3 \sim_\sigma M$.

\begin{lemma}\label{con_sum}
  If $M = M_1 \# M_2$ and $S^3\sim_\sigma M_j$, $j=1,2$, then
  $S^3\sim_\sigma M$.
\end{lemma}
\Dim As $S^3=S^3 \# S^3$, the lemma follows immediately.

\cvd

\smallskip

We write
$$M = M(x,y)$$ to mean that $M$ is encoded by a genus $g$ Heegaard
diagram
$(\Sigma,x,y)$ where $x=\{x_1,\dots, x_g\}$, $y=\{y_1,\dots, y_g\}$
are the two non dividing families of simple smooth circles on the
surface $\Sigma$ early denoted by $C^-$ and $C^+$ respectively.
Recall that $x\pitchfork y$.

Let $z=\{z_1,\dots,z_g\}$ be another family of $g$ smooth circles on
$\Sigma$ which does not divide the surface. Assume that $z\pitchfork
x$ and $z\pitchfork y$. Recalling the reconstruction of $M=M(x,y)$
from the diagram, we can assume that $z$ is traced on the Heegaard
surface $\Sigma \sim \Sigma \times \{0\}$. Give an orientation every
$z_j$, fix a system of disjoint tubular neigbourhoods $U_j$ of every
$z_j$ in $M$ such that $\partial U_j \pitchfork \Sigma$ along a pair
of curves parallel to $z_j$, and select the longitude $l_j \subset
\partial U_j$ given by the component of $\partial U_j \cap \Sigma$
whose orientation is parallel to the one of $z_j$.
For every $j$, up to isotopy there is a unique framing
$\rho_j: S^1\times D^2\to U_j$ so that the longitude $l_j$ is carried
by $\rho_j$; thus we have determined a framed link $L:= \cup_j (z_j,l_j)$
in $M=M(x,y)$. These trivializations are used as attaching maps of disjoint
$2$-handles so that we have constructed a special triad
$$(W,M,\tilde M), \ \tilde M \sim_\sigma M \ . $$
The following simple lemma, which is in fact the core of the proof, 
establishes a key relationship between surgery
equivalence and Heegaard splitting. In the situation depicted so far we have

\begin{lemma}\label{surg_split} $\tilde M \sim M(x,z) \# M(z,y)$.
\end{lemma} 
\Dim  Denote by $M_0(x,z)$ the manifold with spherical boundary
obtained by removing from $M(x,z)$ the interior of a smooth embedded $3$-disk.
Similarly for $M_0(z,y)$. It follows straightforwardly by comparing the reconstruction of
$M(x,z)$ and $M(z,y)$ from the diagrams and the construction of $\tilde M$
by surgery on $M$ along the framed link  $L:= \cup_j (z_j,l_j)$ that, up
to diffeomorphism, $\tilde M$
is obtained by gluing $M_0(x,z)$ and $M_0(z,y)$ by a diffeomorphism
between the boundaries. With the terminology of Section \ref {twisted-sphere},
$\tilde M$ is a weak connected sum of $M(x,z)$ and $M(z,y)$. 
Then by Smale theorem (Proposition \ref{on-twisted}, (1), $m=3$)
it is a true connected sum.

\cvd

\smallskip

The last ingredient is a suitable {\it measure of the complexity of
  the Heegaard diagrams}. Let $(\Sigma,x,y)$ be such a diagram of
genus $g$. Recall that every $x_i\cap y_j$ is a finite set and denote
by $|x_i\cap y_j|$ the number of elements (we stress that it is the
``geometric'' number, no algebraic intersection numbers are involved).
Then set
$$c(\Sigma M, x, y):= (g,r:= \min_{i,j}|x_i\cap y_j|)\in \N^2 $$
where $\N^2$ is endowed with the lexicographic order. We will achieve
the result by (double) induction on the complexity $c$ of a given
splitting of $M$.

The initial step is when  $g=0$;  in such a case by the very definition
$M$ is a twisted $3$-sphere, so it is a true smooth sphere again by
Smale theorem (Proposition \ref{on-twisted}, (2), $m=3$);  the empty surgery does the job.

Let $M=M(x,y)$ of complexity $c=(g,r)$ and assume that $S^3\sim_\sigma M'$
for every $M'$ admitting an encoding diagram of complexity
$c'=(g',r') < c=(g,r)$.

If $c=(g,1)$, then the given diagram is a stabilization of a diagramm
of genus $g-1$, hence $S^3\sim_\sigma M$ by the inductive hypothesis.

If $c=(g,0)$, it is not restrictive to assume that $x_1\cap y_1=\emptyset$.
\smallskip

{\bf Caim 1.} {\it There exists a non separating circle $z_1$ on $\Sigma$
  which intersects each of $x_1$ and $y_1$ transvesely at a single point.}

\smallskip

Assuming this fact, extend $z_1$ to a non dividing family $z$ of $g$
circles on $\Sigma$, $z\pitchfork x$ and $z\pitchfork y$.
Then both $M(x,z)$ and $M(z,y)$ have encoding diagrams with $r=1$
and we conclude by applying the previous case and Lemmas
\ref{con_sum}, \ref {surg_split}.
\smallskip

Assume that $r>1$. It is not restrictive to assume that $r=|x_1\cap y_1|$.
\smallskip

{\bf Claim 2.}  {\it There exists a non separating circle $z_1$ on $\Sigma$
which intersects each of $x_1$ and $y_1$ transvesely at a number of points
stricly less than $r$.}
\smallskip

Assuming this fact, extend $z_1$ to a non dividing family $z$ of $g$
circles on $\Sigma$, $z\pitchfork x$ and $z\pitchfork y$.
Then both $M(x,z)$ and $M(z,y)$ have encoding diagrams of the same genus $g$
but with strictly smaller complexity anyway. Then by the inductive
hypothesis $S^3$ is surgery equivalent to both and again we can conclude
by applying Lemmas
\ref{con_sum} and \ref {surg_split}.

\smallskip  

It remains to prove the two claims. As for Claim 1, there are two
possibilities, either $\Sigma':=\Sigma \setminus (x_1\cup y_1)$ is
connected or non connected. Take a small segment $\gamma$ in $\Sigma$
tranvese to $x_1$ at one point, with endpoints $p_0,p_1$; similarly
let $\gamma'$ be transverse to $y_1$ at one point, with endpoints
$p'_0,p'_1$. If $\Sigma'$ is not connected, up to reordering, we can
assume that the couples of endpoints $p_0, p_0'$ and $p_1,p'_1$ belong
to different connected components.  Then in both cases a smooth circle
$z_1$ in $\Sigma$ with the required properties can be obtained of the
form
$$z_1= \gamma \cup \alpha \cup \gamma' \cup \alpha'$$ 
where $\alpha$ is a smooth arc which connects $p_0$ and $p_0'$, while $\alpha'$
is such an arc connecting $p_1$ and $p_1'$. 

As for Claim 2, let $A$ and $B$ two points of $x_1\cap y_1$ which are
adjacent in $x_1$. Then there is an arc $\alpha$ in $x_1$ which
intersects $y_1$ only at its endpoints $A$ and $B$.  These points also
divide $y_1$ in two arcs $\beta$ and $\gamma$.  As $y_1$ does not
separate $\Sigma$, there is at least one of these arcs, say $\beta$,
such that $\alpha \cup \beta$ does not separates $\Sigma$.  Then we
can construct $z_1$ made by a parallel copy $\alpha'$ of $\alpha$
which near $A$ is in the direction of $\beta$, completed by a segment
$\beta'$ close to $\beta$. One realizes that $z_1$ intersects $x_1$ in
at most $r-1$ points and intersects $y_1$ in at most one point. So
$z_1$ has the desired properties. This proof of Theorem \ref{LW} is
now complete.

\cvd

\subsection{On Kirby's calculus}\label{kirby-cal}
We have proved that for every orientable compact, connected,
boundaryless $3$-manifold $M$ there is a special triad
$(W,S^3,M)$ which realizes the surgery equivalence $S^3\sim_\sigma M$,
so that $W$ admits an ordered handle decompositions consisting only
of $2$-handles. Every such a handle decomposition with say $k$ handles
is encoded by a framed link $L$ in $S^3$ with $k$ constituent knot $K_j$, $j=1,\dots , k$. 
For every $K_j$, its framing is encoded by a parellel longitude $l_j$; fixing
an auxiliary parallel orientation of both $K_j$ and $l_j$, 
this last is encoded by the linking number $L(K_j,l_j)$, that is, equivalently, by
the intersection number of $l_j$ with any oriented Seifert surface of $K_j$ in $S^3$.
The natural question is how two such  framed links are related two each other.
Certainly a given handle deconposition can be modified by handle sliding and this
can be translated in terms of the corresponding framed links. Moreover we must
consider the possibility of modifying the special triad without changing its
boundary. A distinguished way to do it  consists in attaching  a $2$-handle
with attaching circle contained and unknotted in a $3$-ball disjoint from the other
link components, and with framing equal to $\pm 1$. One realizes that
this does not modify the boundary  while we pass from $W$ to $W\cs \pm\PP^2(\C)$.
This is called an {\it elementary blow-up move}. We can consider also
the inverse (negative) move of removing such a handle.
An important Kirby's result \cite {Kirby2} can be formulated, somewhat qualitatively, as follows.

\begin{theorem}\label{Kirby-C} Two framed link $L_1$ and $L_2$ in $S^3$
encode a realization of $S^3\sim_\sigma M$ if and only if they are related
to each other by a finite sequence of modifications which either
translate $2$-handle sliding or are  positive/negative elementary blow-up moves.
\end{theorem}

The proof is rather demanding and is based on Cerf's theory \cite{Ce2}.
After such a qualitative statemet, successive efforts have been devoted
to convert it into an efficient diagrammatic calculus on framed links in $S^3$.
Kirby himself found a generator (called ``band move") for the handle sliding;
this is {\it not a `local' move}, and resembles a move described above on Heegaard
diagrams. Later in \cite{FR} one points out an {\it infinite} family of {\it local} moves
generating the whole calculus. Finally in \cite{Mart2} one has determined
a generating {\it finite} family of {\it local} moves.  

\section{On $\eta_3=0$.}\label{eta3=0}
Referring to Section \ref {non-or-surg}, the following two theorems
can be obtained as a corollary of each other.

\begin{theorem}\label{LW2}
  Every non orientable compact connected boundaryless $3$-manifold $M$
  is surgery equivalent to $\MG$ ($M \sim_\sigma \MG$).
\end{theorem}

\begin{theorem}\label{eta_3=0} $\eta_3=0$.
\end{theorem}

\subsection{On some proofs of $\eta_3=0$} In the spirit of
the above discussion about $\Omega_3=0$, we give here a few indication
about ``direct'' proofs of Theorem \ref{eta_3=0}. Certainly it is
contained in the general statement of Thom's Theorem \ref{eta-alg} and
in a sense this is the first proof of this result. However, Rohlin
claimed, without further explaination (see \cite{GM}), that the method
he had used to prove $\Omega_3=0$ allows to prove the same in the non
orientable case. This is not so immediate. Starting from a general
immersion of $M$ (non orientable) in $\R^5$, the ``embedding up to
bordism'' works as well and we can assume that $M$ is actually
embedded into $\R^5$.  However, (recall Remark \ref{unor-seifert}), if
a tubular neighbourhood $U$ of $M$ in $\R^5$ is associated to a
splitting $T(M)\oplus \xi$ of the restriction of $T(\R^5)$ to $M$, we
cannot assume in general that $\xi$ has a nonwhere vanishing section
and hence we cannot assume that there is a possibly non orientable
Seifert surface.  To conclude it would be enough to find $M'$ embedded
in some $5$-manifold $X$ such that $[M]=[M']\in \eta_3$, $[M']=0\in
\Hh^2(X,\Z/2\Z)$, and there is a splitting $T(M)\oplus \xi'$ of the
restriction of $T(X)$ to $M'$ such that $\xi'$ has a nowhere vanishing
section. This can be achieved as follows (see also the suggestion at
pag. 91 of \cite{GM}).  Let $M$ embedded in $\R^5$ be as above.
Consider the Euler class of $\xi$ belonging to $\eta_1(M)$. This is
represented by smooth circle $C$ on $M$. Take the blow up say $X$ of
$\R^5$ along $C$ (see Section \ref {blow-up}); let $M'$ be the blow up
of $M$ along $C$ which is embedded into $X$ as the strict transform of
$M$. One can check that $M'\subset X$ verify the required
properties. In particular $[M']=[M]+[S^1 \times \PP^2(\R)]=[M]\in
\eta_3$.

\subsection{On some proofs that $M \sim_\sigma \MG$}
Lickorish extended in \cite{Lick2} his main result on the generators
of the mapping class groups to non orientable surfaces.  This allows
him to extend also the proof about the surgery equivalence to the non
orientable case.

In \cite{AG} the simpler clever proof of \cite{Rourke} has been
extended  to the non orientable
case.

 \section {Combing and framing}\label{3-parallel}
 A main result of this section will be that every compact connected
 orientable boundaryless $3$-manifold $M$ is {\it parallelizable}.
 Current modern proofs of
 this primary result in $3$-dimensional differential topology
 (originally attributed to Stiefel \cite{Sti}) use either a mixture 
 of {\it spin structures} and of {\it Stiefel Whitney classes} theory
 (see for instance \cite{Ge}, Section 4.2), or a refinement due to Kaplan 
 \cite{Ka} of Lickorish-Wallace theorem  by means of the
 so called {\it Kirby calculus} (see also \cite{FM}, Section 9.4.). 
 We do not dispose of this artillery.
 Nevertheless, by following \cite{BL} we will provide two selfcontained
 elementary proofs, revealing by the way different aspects
 of the question. The first proof uses some ideas of the last mentioned approach, however
 it avoids the use of both Lickorish-Wallace Theorem and Kirby calculus. 
 The second proof will result
 from a parallel discussion about  combing and framing $3$-manifolds.
We will also provide a classification
of combings with respect to a given auxiliary reference framing. 

From now on $M$ will denote a compact connected orientable
boundaryless $3$-manifold.  Alike every odd dimensional manifold, $M$
is combable, then it carries nowhere vanishing tangent vector fields
$v$. These are considered up to smooth homotopy through such fields
and called {\it combings} of $M$. We will systematically confuse a
homotopy class with  suitable representatives. As we know, a {\it framing}
$\Ff$ of $T(M)$ is a triple $(v,w,z)$ of pointwise linearly
independent tangent vector fields. Also framings are considered up to
homotopy; the three components of $\Ff$ determine a same combing of
$M$.  Fixing any auxiliary riemannian metric $g$ on $M$, we can assume
that a given combing is (represented by) an unitary field with respect
to $g$, and every framing is represented by pointwise orthonormal
fields. A framing, if any, determines also an orientation of $M$ (so
that orientability of $M$ is a necessary condition). If $M$ is {\it
oriented} and parallelizable, then there are framings which induce
the given orientation. From now on we will assume that $M$ is
{\it oriented}, by fixing an auxiliary orientation.
 
 \subsection{Framing via even surgery}\label{even-framing}
The first remark is that it is enough to prove that $M$ is 
 almost-parallelizable. A {\it quasi-framing} of $M$ is a framing of $T(M)$
over a submanifold of the form
$$ M_0:= M\setminus {\rm Int}(B)$$
where $B$ is a smooth $3$-disk in $M$. We say that $M$
is {\it almost-parallelizable} if admits a quasi-framing.
In such a case, by the uniqueness of the disk up to ambient isotopy,
we see that the choice of the disk $B$ is immaterial. We have

\begin{lemma}\label{Q-F} $M$ is parallelizable if and only if it is
almost-parallelizable.
\end{lemma}

\Dim An implication is trivial. As for the other implication, we can
assume that $B$ is contained in a chart of $M$ and looks standard
therein as well as the auxiliary metric. Then the restriction of a
quasi-framing $\Ff'$ to $\partial B=S^2$ is encoded by a map
 $$ \rho: S^2 \to SO(3) \ . $$ We know that $SO(3) \sim \PP^3(\R)$
 (Example \ref {so(3)}), with $S^3$ as universal covering space, hence $\pi_2(SO(3))\sim
\pi_2(S^3)=0$.  It follows that $\rho$ extends to $\hat \beta: B\to
SO(3)$, and that $\Ff'$ extends to a framing $\Ff$ of the whole
$T(M)$.
 
 \cvd
 
 \smallskip
 
 Let $M$ be obtained by longitudinal
 surgery along a framed link $L$ in $S^3$; we write
 $$ M = \chi(S^3,L) \ . $$
 $M$ is the final boundary of a triad $(W,\emptyset, \chi(S^3,L))$
  where  $W$ is obtained by attaching disjoint $2$-handles to $D^4$ at $S^3=\partial D^4$.
  Every $2$-handle $D^2\times D^2$ determines a constituent knot $K$ of $L$, so that
  $\partial D^2 \times D^2 \sim N(K)$, $N(K)$ being a tubular neighbourhood of $K$
  in $S^3$, $\partial D^2 \times \{0\}$ being identified with a longitude $l_K$ on $\partial N(K)$
  olong $K$.  The framing of every component $K$ of $L$ is encoded by the linking number 
 $n_K\in \Z$ between $K$ and the longitude $l_K$, where $K$ and $l_K$
 are co-oriented in such a way that the projection of $L_K$ onto $K$ is of degree $1$.
 We say that the surgery is {\it even} if for every constituent knot $K$ of $L$,
 $n_K \in 2\Z$. We have 
 
 \begin{proposition}\label{even-framed} Let $(W,\emptyset, M)$
 be the triad associated to an even surgery $M=\chi(S^3,L)$.
 Then $W$ is parallelizable.
 \end{proposition}
 \Dim To simplify the notation, we give the proof for a one-component link but this generalizes
straightforwardly. So let $L=(K,n)$, $n\in 2\Z$. Both $D^4$ and $D^2\times D^2$ are parallelizable,
so we have to show that they carry some framings which match on $N(K)$.
Fix a reference framing $\Ff_0$ on $D^4$; the restriction to $N(K)$ of any framing $\Ff$ on the $2$-handle
is encoded by a map $\rho: N(K)\to SO(4)$.  
Viewing $S^3$ as the group of unit quaternions one can construct a 
$2$-fold covering map $S^3 \times S^3 \to SO(4)$ showing that $\pi_1(SO(4)) = \Z/2\Z$ (see Example \ref {so(3)}).
As the solid torus $N(K)$ retracts to $K\sim S^1$, 
$\rho$ determines an element of $\Z/2Z$, and the two framings coincide on $N(K)$ if and only if this is equal to $0$. 
It can be readily seen that this element is equal to the number $n$ mod $(2)$. 

\cvd

\smallskip

\begin{corollary}\label{stable-par} Let $M=\chi(S^3,L)$ be an even surgery.
Then $M$ is stably-parallelizable (i.e. $T(M)\oplus \epsilon^1$ is a product
bundle).
\end{corollary}
\Dim Let $(W,\emptyset, \chi(S^3,L))$ be as above. Then
$T(W)_{M}=T(M)\oplus \nu$ where $\nu$ is a trivial normal line bundle
of $M=\partial W$ in $W$. We know by the proposition that $T(W)$
is a product bundle.

\cvd

\smallskip

\begin{lemma}\label{stable-almost} If $M$ is stably parallelizable then it is 
almost-parallelizable.
\end{lemma}
\Dim As $T(M)\oplus \epsilon^1 = M\times \R^4$, every $T_xM$ is an
oriented $3$-plane in $\R^4$. So we have a smooth classifying map
$\rho: M \to S^3$ where the sphere is considered as the space of
oriented $3$-planes in $\R^4$,  and $T(M)$ is the pull back of the
corresponding tautological bundle (see Chapter \ref{TD-SMOOTH-MAN}).
Now we know that $M_0$ retracts onto a $2$-dimensional spine 
$\PP_0$ as in Section \ref{spine}.  Hence the restriction of $\rho$
to $\PP_0$ is not surjective, then it is homotopic to a constant map, 
the restriction of $TM$ to $\PP_0$ whence to $M_0$ is a product bundle.
 
\cvd

\smallskip

\begin{remark}{\rm Lemma \ref{stable-almost} holds in every dimension $n$;
the key point is that $M\setminus {\rm Int}(B^n)$ has the homotopy type of
a CW-complex of dimension less or equal $n-1$ (see Section \ref{CW}).
}
\end{remark}

Recall the notion of weak connected sum
given in Section \ref{twisted-sphere}. We know by Smale theorem that
$3$-dimensional weak connected sums are veritable connected sums,
but we do not need this fact in the present discussion. The following lemma is trivial. 

\begin{lemma}\label{trivial}  If there
exists $M'$ such that a weak connected sum of $M$ and $M'$ is parallelizzable, then 
 $M$ is almost parallelizable.
\end{lemma}

\cvd

\smallskip

A Heegaard splitting (of some genus $g$) of $M$ can be encoded by a non dividing family say $L$
of $g$ smooth circles on the boundary $\partial \HG_g$ of an handlebody $\HG_g$.
We can assume that $\HG_g$ is embedded in a standard way in $S^3$ so that
$\HG'_g:=\overline { S^3 \setminus \HG_g}$ is also a handlebody of genus $g$, and
we have a Heegaard splitting of $S^3$. Give every component $K$ of $L$
the framing carried by a tubular neighbourhood of $K$ in $\partial \HG_g$.
Then we have a framed link $L$ in $S^3$. By applying the proof of Lemma \ref {surg_split} 
we readily have
\begin{lemma}\label{M'}
$\chi(S^3,L)$ is a weak connected sum of $M$ and $M'$, for some $M'$.
\end{lemma}

\cvd

\smallskip

So, by combining the above lemmas, 
in order to show that $M$ is almost parallelizable (hence parallelizable)
it is enough to show that we can implement the above construction
in such a way that the surgery $\chi(S^3,L)$ is even.
Fix any embedding $L\subset \partial \HG_g \subset \HG_g \subset S^3$
as above. Fix a system $\mu=\{m_1,\dots, m_g\}$ of $g$ meridians on $\partial \HG_g$
(which bound $2$-disks properly emebedded in $\HG_g$)
and a dual system of $g$ meridians $\lambda=\{l_1,\dots, l_g\}$ for the complementary
handlebody $\HG'_g$. A Dehn twist on $\partial \HG_g$ along a curve $m_i$  extends to a diffeomorphism
of the whole $\HG_g$. Hence we can modify the family $L$ by applying any finite sequence of such
Dehn twists, keeping the fact that $\chi(S^3,L)$ is a weak connected sum of $M$ and $M'$, for some $M'$. 
We are reduced to prove the following lemma.

\begin{lemma}\label{even-by-twist} Up to a suitable finite sequence of Dehn
twists along the meridians in $\mu$, $\chi(S^3,L)$ is an even surgery. 
\end{lemma}
\Dim The question can be reduced to $\Z/2\Z$-linear algebra on $\eta_1(\partial \HG_g)$.
Start with any surgery $\chi(S^3,L)=M\cs M'$ as above. 
The union of curves in the families $\mu$ and $\lambda$ form a symplectic basis of
$\eta_1(\partial \HG_g)$ with respect to the intersection form. So, by confusing classes 
mod $(2)$ and representatives and
setting $L=\{K_1,\dots, K_g\}$, we have the $\Z/2\Z$-linear combinations:
$$ K_j = \sum_{i=1}^g (a^j_i m_i + b^j_il_i) \ . $$
The framing mod $(2)$ of $K_j$ is given by
$$n_j = \sum_i a^j_ib^j_i \in \Z/2\Z \ . $$
A Dehn twist $T_j$ along $m_i$ acts on $\eta_1(\partial \HG_g)$ so that
$$T_i(l_i)= l_i + m_i$$
while it is the identity on the other $2g-1$ elements of the given basis.
All intersection numbers mod $(2)$ of the curves of $L$ vanish, that is 
$$ K_r\bullet K_s = 0, \ r, s =0,\dots g \ . $$
This means that the coefficients of the above lnear combinations verify
the system of conditions: 

\begin{equation}\label{C}
\sum_{i=1}^g (a_i^rb_i^s + a_i^sb_i^r) = 0,  \  r, s =0,\dots g \ . 
\end{equation}

We allow ourselves to apply twist combinations of the form
$T^{x_1}_1\dots T^{x_g}_g$.
Then we want to show that the $\Z/2\Z$-linear non homogeneous system
\begin{equation}\label{B}
\sum_{i=1}^g (x_i + b_i^{r})a_i^r = 0,  \  r =1,\dots g \ . 
\end{equation}
admits a solution in $(\Z/2\Z)^g$. Note that we tacitly use several times that
$z=z^2$ for every $z\in \Z/2\Z$.
If for every $r$ all $a^r_i=0$, then every $(x_1,\dots,x_g)$ is a solution.
Otherwise we can assume that $a^1_1=1$. Then the solution of the equation
$$\sum_{i=1}^g (x_i + b_i^{1})a_i^1 = 0$$ 
are of the form $x_1= \sum_{j=2}^g c_j x_j$.
By replacing in the other equations and using the relations \ref{C},
we are reduced to solve a system in $x_2,\dots, x_g$ of the same form
$$\sum_{i=2}^g (x_i + \tilde b_i^r)\tilde a_i^r=0, \ r=2,\dots, g$$
with 
$$ \tilde a^r_i = a_1^ra_i^1 + a_i^r, \ \tilde b^r_i = a_1^rb^1_i + b^r_1 \ . $$
One ferifies directly that these new coefficients formally satify the 
corresponding conditions \ref{C}. So we can conclude by recurrence.

 \cvd

\smallskip

\begin{remark}{\rm It is proved in \cite{Ka}, see also \cite{FM}, that for every $M$
as above there is an even surgery $M=\chi(S^3,L)$. Starting from any surgery presentation of $M$
with associated triad $(W,\emptyset, M)$ (which exists by Lickorish-Wallace Theorem), 
the proof consists in an algorithm which
modifies the triad to some $(W',\emptyset, M)$ associated to an even surgery.  
More precisely, by using some notions that we will define in Chapter \ref {TD-4},
one proves firts that every $L$ contains a so called {\it characteristic sub-link} and that
the surgery is even if a characteristic sub-link is empty. Then the algorithm
reduces progressively the number of components of a characteristic sub-link 
by means of certain moves on the handle decompositions 
(organized in an efficient so called `Kirby calculus') which 
may change the $4$-manifold $W$ by keeping the triad boundary fixed.
Note that this proof does {\it not} use the harder fact that 
Kirby calculus connects any two surgery presentations of $M$ \cite{Kirby2}.
}
\end{remark}
\smallskip

Our first proof that $M$ is parallelizable is now complete.

\cvd

\smallskip

Next we will elaborate on the second proof.

\subsection{On the cobordism ring of an orientable $3$-manifold}\label{relevantcobord}
 We specialize the results of Chapters \ref{TD-LINE-BUND}.
 In the present situation the relevant co-bordism modules are
 $$ \Hh^j(M;\Z/2\Z), \ \Hh^j(M;\Z), \ j=0,1,2, 3 \ . $$
We summarize here some properties which we will use.

 \smallskip
 
 - $\Hh^3(M;\Z/2\Z)\sim \Hh^0(M;\Z/2\Z)\sim \Z/2\Z$ by the isomorphism
 which associates the usual generator of $\Hh^3(M;\Z/2\Z)$ to the
 fundamental class mod $(2)$ $[M]$; similarly over $\Z$.
 \smallskip
 
- $\Hh^2(M;\Z/2\Z))=\eta^2(M)=\eta_1(M)$
 \smallskip
 
- $\Hh^2(M;\Z/2\Z)\sim \Hh^1(M;\Z/2\Z)$ in a natural way: if
$\alpha=[F]\in \Hh^1(M;\Z/2\Z)$ we can assume that the embedded
surface $F\subset M$ is connected and does not divide $M$ if
$\alpha\neq 0$. If $\gamma$ is a smooth simple arc in $M$ trasverse to
$F$ at one point, it can be completed to a smooth circle $c$ by means
of an arc $\gamma'$ contained in $M\setminus F$ so that $[F]\sqcup
[c]=1$. Viceversa, if $[c]\neq 0 \in \Hh^2(M;\Z/2\Z)$, then it is part
of a basis $\Bb$ of $\Hh^2(M;\Z/2\Z)$ which is finite dimensional. The
functional $[c]^*$ belonging to the dual basis composed with the
natural homomorphism $\pi_1(M)\to \eta_1(M)$ defines a $\Z/2\Z$-valued
representation of the fundamental group that can be realized by a
connected hypersurface $F$, so that in particular $[F]\sqcup
[c]=1$. Moreover we can assume that $F$ intersects transversely $c$ at
one point: if $F$ intersects $c$ at an odd number of points, we can
reduce them to one by attaching suitable embedded $1$-handles along
$c$ and performing surgeries of $F$.
 \smallskip
 
- If $c$ is a connected oriented smooth circle in $M$ such that $[c]=0
\in \Hh^2(M;\Z)$ then there is an oriented Seifert surface for $c$ in
$M$; if $[c]=0 \in \Hh^2(M;\Z/2\Z)$ then there is a possibly non
orientable Seifert surface for $c$ in $M$;
 \smallskip
 
- Consider the natural forgetting morphism $\Hh^2(M;\Z)\to
\Hh^2(M;\Z/2\Z)$. We have

\begin{lemma}\label{evenclass}
A class $\alpha\in \Hh^2(M;\Z)$ belongs to the kernel of the
forgetting morphism $\Hh^2(M;\Z)\to \Hh^2(M;\Z/2\Z)$ if and only if
$\alpha$ is an {\rm even class} that is there is $\beta \in
\Hh^2(M;\Z)$ such that $\alpha = 2\beta$.
\end{lemma}
\Dim We can assume that $\alpha$ is represented by a connected
oriented smooth circle $c$. By hypothesis $c$ is the boundary of a
possibly non orientable connected compact surface $F$ embedded in
$M$. If $F$ is orientable, then $\alpha = 0$ and we have done.  If $F$
is not orientable, it follows from the classification of surfaces that
there is a smooth $1$-submanifold $C$ on ${\rm Int}(F)$ such that a
tubular neighbourhood $U(C)$ of $C$ in $F$ is union of M\"obius
strips, and $F\setminus C$ is orientable. Then orient $F\setminus C$
in such a way the oriented $c$ inherits the boundary orientation, and
orient consequently $C':= \partial U(C) \subset F\setminus {\rm
  Int}(U(C))$.  Then $[c]=[C'] \in \Hh^2(M;Z)$ and $[C']=2[C"]$ where
$C"$ is the union of the cores of $U(C)$ oriented in such a way that
the restriction of the projection of $C'$  onto every core is of positive degree.

 \cvd
 
 \smallskip

\subsection{Combings and orthogonal plane distributions}\label{plane-field}
Let $v$ be a combing of $M$. Fix an auxiliary metric as above.
We have the distribution of orthogonal tangent $2$-planes
$$\{P_x:={\rm span}(v(x))^\perp\}_{x\in M} \ . $$ These planes
$P_x$ are oriented by the unique orientation which added to
 $v(x)$ agrees with the given orientation on $T_xM$.
 This actually defines an oriented rank-$2$ vector bundle $\xi_v$ on $M$
 whose strict equivalence class does not depend on the choice of the combing
 representative nor of the auxiliary metric. We consider the  oriented Euler class
 $$ e^2(\xi_v)\in \Omega^2(M)=\Omega_1(M) \ . $$ In fact
 $e^2(\xi_v)\in \Hh^2(M;\Z)$.  If $\xi_v$ has a non vanishing unitary
 section $w$ orhogonal to $v$, then $(v,w)$ extends to the unique
 orthonormal framing $\Ff=(v,w,z)$ of $T(M)$ such that the
 orientations are compatible. So $\xi_v$ is trivial if and only if it
 admits a nowhere vanishing section $w$ as above. We know from section
 \ref{complexLB} that
 \begin{lemma}\label{vanish-e}
   The bundle $\xi_v$ has a non vanishing section, if and only if the
   Euler class $e^2(\xi_v)$ vanishes.
 \end{lemma}
 
 \cvd
 
\smallskip

 As usual, $\omega^2(\xi_v)\in \Hh^2(M;\Z/2\Z)$ is the image of $e^2(\xi_v)$ 
 via the natural forgetting map.
 \smallskip

 {\bf Combing comparison class.}  We can associate to an ordered pair
 of unitary combings $(v, v')$ of $M$ a smooth section $v\x v'$ of
 $\xi_v$ as follows.  At a point $x\in M$ where $v(x) \neq\pm v'(x)$,
 $v\x v'(x) \in P_v(x)\subset T_x M$ is the ``vector product'' of
 $v(x)$ and $v'(x)$, i.e.~the only tangent vector such that
\begin{itemize}
\item 
$\|v\x v'(x)\|_{g(x)}^2 = 1 - g(v,v')^2$;
\item
$v\x v'(x)$ is $g(x)$-orthogonal to $v(x)$ and $v'(x)$; 
\item 
$(v(x), v'(x), v\x v'(x))$ is an oriented basis of $T_ x M$.
\end{itemize}
At a point $x\in M$ where $v(x) = \pm v'(x)$, we set $v\x v'(x) = 0$. 

If the two unitary combings $v$ and $v'$ are generic,
the section $v\x v'$ of $F_v$ is transverse to the zero section and
the zero locus 
\[
C := \{x\in M\ |\ v\x v'(x) = 0\} \subset M
\]
is a disjoint collection of simple closed curves. Moreover, $C =
C_+\cup C_-$, where
\[
C_+ = \{x\in M\ |\ v(x) = v'(x)\}\quad\text{and}\quad
C_- = \{x\in M\ |\ v(x) = -v'(x)\}.
\]
By the very definition of $e^2(\xi_v)$, $C$ can be oriented to
represent the Euler class of $\xi_v$.  Indeed, let $E(\xi_v)$ denote the
total space of $\xi_v$, $M_0\subset E(\xi_v)$ the zero-section and $M_1 =
v\x v'(M)\subset E(\xi_v)$. Under the natural identification of $M$ with
$M_0$ the submanifold $C$ is identified with $M_0\cap M_1$. By
transversality, for each $x\in M_0\cap M_1$ the natural projection
$p_x : T_x E(\xi_v)\to P_v(x)$ maps isomorphically the image under $(v\x
v)'_*$ of the fiber $N_x(C)$ of the normal bundle of $TC\subset TM|_C$
onto $P_v(x)$. Therefore, the given orientation on $\xi_v(x)$ can be
pulled-back to $N_x(C)$ and, together with the orientation of $T_x M$,
it induces an orientation on $T_x C$ in a standard way.

\begin{definition} \label{compare}{\rm
An ordered pair of unitary combings $(v,v')$ of $M$ such that $v\x v'$
is a section of $\xi_v$ transverse to the zero section will be called a
{\em generic pair} of unitary combings.  We define the {\em comparison
  class} $\alpha(v,v')\in \Omega^2(M)$ of a generic pair of unitary
combings as the class $[C_-]$ carried
by the collection of curves $C_-$ oriented as part of the oriented
zero locus of $v\x v' : M\to \xi_v$ representing $e^2(\xi_v)$.}
\end{definition}

\begin{lemma}\label{l:comparison-properties}
Let $(v,v')$ be a generic pair of unitary combings of $M$. Then, 
\[
\alpha(v,v') = -\alpha(v',v)\quad\text{and}\quad \alpha(v,-v') = \alpha(v', -v).
\]
\end{lemma} 

\Dim
For each $x\in C$ the equality $\xi_v(x)=\xi_{v'}(x)$ holds, with the
orientations of $\xi_v(x)$ and $\xi_{v'}(x)$ being the same or different
according to, respectively, whether $x\in C_+$ or $x\in C_-$. We may
choose a tubular neighborhood $U = U(C)$ such that the restrictions of
the tangent plane fields $P_v |_U$ and $P_{v'} |_U$ are so close that
there is a vector bundle isomorphism $\varphi : \xi_v
|_U\stackrel{\cong}{\to} \xi_{v'} |_U$
which is the identity map on the
intersections $P_v (x)\cap P_{v'} (x)$, $x\in U$, is
orientation-preserving near $C_+ = \{x\in M\ |\ v(x) = v'(x)\}$ and
orientation-reversing near $C_- = \{x\in M\ |\ v(x) = -v'(x)\}$. Since
$\varphi\circ (v\x v') = v\x v' = - v'\x v$ and $- v'\x v$ is obtained
by composing the section $v'\x v$ with the orientation-preserving
automorphism of $F_{v'}$ given by minus the identity on each fiber,
the orientation on $C_-$ as part of the zero locus of $v\x v' : M\to
\xi_v$ is the opposite of its orientation as part of the zero locus of
$v'\x v = -v\x v': M\to \xi_{v'}$.  This implies $\alpha(v,v') =
-\alpha(v',v)$.  Similarly, the orientation on $C_+$ as part of the zero
locus of $v\x (-v') : M\to \xi_v$ coincides with its orientation as part
of the zero locus of $(-v')\x v = v'\x (-v): M\to \xi_{v'}$, which
implies $\alpha(v,-v') = \alpha(v',-v)$.

\cvd

\smallskip

\begin{lemma}\label{l:compare} 
Let $(v,v')$ be a generic pair of unitary combings of $M$. Then, 
\[
e^2(\xi_v)-e^2(\xi_{v'})= 2\alpha(v,v').
\]
\end{lemma}

\Dim
According to the definitions we have 
\[
e^2(\xi_v) = \alpha(v,v') + \alpha(v,-v')\quad\text{and}\quad
e^2(\xi_{v'}) = \alpha(v',v) + \alpha(v',-v).
\]
The statement follows applying Lemma~\ref{l:comparison-properties}
after taking the difference of the two equations.

\cvd

\smallskip

{\bf Combing Pontryagin surgery.}  Let $v$ be a unitary combing of $M$
and $C\subset M$ an oriented, simple closed curve such that the
positive, unit tangent field along $C$ is equal to $v|_C$ and there is
a trivialization
\[ 
j:D^2\times S^1 \stackrel{\cong}{\to} U(C)
\]
of a tubular neighborhood of $C$ in $M$ such that
\[
v \circ j= j_*(\partial /\partial \phi),
\]
where $\phi$ is a periodic coordinate on the $S^1$-factor of $D^2\x
S^1$. Let $(\rho, \theta)$ be polar coordinates on the
$D^2$-factor. Following terminology from~\cite{BP}, we say that a
unitary combing $v'$ is obtained from $v$ by {\it Pontryagin surgery}
along $C$ if, up to homotopy, $v'$ coincides with $v$ on $M\setminus
U(C)$ and
\[  
v' \circ j =
j_*\left
(-\cos(\pi\rho)\frac{\partial}{\partial \phi} -\sin(\pi\rho)\frac{\partial}{\partial \rho}\right)
\] 
on $U(C)$. 

\begin{remark}{\rm
A basic fact not used in this paper is that any two combings of $M$
are obtained from each other, up to homotopy, by Pontryagin
surgery ~ \cite{BP}.}
\end{remark}

\begin{lemma}\label{pont-comb-surg} Let $v$ be a unitary combing of $M$ and
$\beta \in \ H^2(M;\ZZ)$. Then, possibly after a homotopy of $v$, there is a unitary 
combing $v'$ such that $(v,v')$ is a generic pair of unitary combings and 
\[
\alpha(v, v')= \beta. 
\]
\end{lemma}

\Dim Let $C\subset M$ be an oriented simple closed curve representing
the Poincar\'e dual of $\beta$ and let $j : D^2\x S^1 \to U(C)$ be a
trivialization of a neighborhood of $C$.  Without loss of generality
we may assume that the pull-back $j^*(g)$ of the auxiliary metric $g$
on $M$ is the standard product metric on $D^2\x S^1$.  After a
suitable homotopy of $v$ the assumptions to perform Pontryagin surgery
on $v$ along $C$ are satisfied.  Consider a normal disc $D_{\phi_0} =
j(D^2\x\{\phi_0\})$ and let $p = D_{\phi_0}\cap C$. Then, $T_p
D_{\phi_0}$ coincides, as an oriented $2$-plane, with $P_v(p)$ as well
as with the $g(p)$-orthogonal subspace of $T_p C$ inside $T_p M$. Let
$v'$ be a unitary combing obtained from $v$ by first performing a
Pontryagin surgery on $U(C)$ and then applying a small generic
perturbation supported on a small neighborhood of $M\setminus U(C)$.
Then, $(v,v')$ is a generic pair of unitary combings and $C = \{x\in
M\ |\ v(x) = -v'(x)\}$.  By the definition of $\alpha(v,v')$, to prove
the statement it suffices to show that the given orientation of $C$
coincides with its orientation as part of the zero set of $v\x v' :
M\to \xi_v$.  Near $C$ we have
\[
(v\x v')\circ j = j_*\left(-\sin(\pi\rho) \frac{\partial}{\partial \theta}\right) = 
j_*\left(\frac{\sin(\pi\rho)}{\rho}\left(y\frac{\partial}{\partial x}
- x\frac{\partial}{\partial y}\right)\right),
\]
where $x = \rho\cos\theta$ and $y=\rho\sin\theta$ are rectangular
coordinates on the $D^2$-factor. Observe that $j_*$ sends the pair
$(\partial/\partial x, \partial/\partial y)$ to an oriented framing of $\xi_v$.  Using
the resulting trivialization of $\xi_v$ we can write locally the
restriction of $v\x v'$ to to the disc $D_{\phi_0}$ followed by
projection onto $\xi_v$ as follows:
\[
v\x v'|_{D_{\phi_0}} : (x,y)\mapsto  \frac{\sin(\pi\rho)}{\rho} (y, -x) = \pi (y, -x)\
+\ \text{higher order terms}.
\]
It is easy to compute that $(v\x v')_* \circ j_*$ sends $\partial/\partial x$
to $-\pi\partial/\partial y$ and $\partial/\partial y$ to $\pi\partial/\partial x$, and since
the matrix
$\left(\begin{smallmatrix}0&\pi\\-\pi&0\end{smallmatrix}\right)$ has
  determinant $\pi^2>0$ this shows that the restriction of $(v\x
  v')_*$ to the normal bundle to $C$ composed with the projection onto
  $\xi_v$ is orientation-preserving along $C$, concluding the proof.

  \cvd

\smallskip

We are ready to state a main theorem of this section.
 \begin{theorem}\label{3-par} Let $M$ be a compact connected
 oriented boundaryless $3$-manifold. The the following facts
 are equivalent to each other and all hold true.
 \smallskip
 
 (1) $M$ is parallelizable.
 \smallskip
 
 (2) There exists a combing $v$ of $M$ such that $e^2(\xi_v)=0$.
 \smallskip
 
 (3) There exists a combing $v$ of $M$ such that $e^2(\xi_v)$
 is an {\rm even class} 
 \noindent that is of the form $e^2(\xi_v)=2\beta$
 for some $\beta \in \Hh^2(M;\Z)$.
 \smallskip
 
 (4) For every combing $v$ of $M$, $e^2(\xi_v)$ is an even class.
 \smallskip
 
 (5) For every combing $v$ of $M$, $\omega^2(\xi_v)=0 \in \Hh^2(M;\Z/2\Z)$.
 \smallskip
 
  \end{theorem} 
 
 \Dim First we prove the equivalence between the five statements. 
 We will prove $ (j) \Leftrightarrow (j+1)$ for $j=1, \dots, 4$. 

 \smallskip
 
 $(1)\Rightarrow (2)$: If $\Ff=(v,w,z)$ is a framing of $M$, then $e^2(\xi_v)=0$.
 
 $(1)\Leftarrow (2)$: we have already remarked above that if $e^2(\xi_v)=0$ then
 $v$ can be extended to a global framing $\Ff=(v,w,z)$.

 $(2)\Rightarrow (3)$:  this is trivial.
 
 $(2)\Leftarrow (3)$: If $e^2(\xi_v)= 2\beta$, then by
 applying the Pontryagin surgery to $v$ and
 the class $-\beta$, we get $v'$ such that
 $$ e^2(v')= -2\beta + e^2(v)= 0 \ . $$
 
 $(3)\Rightarrow (4)$: If $e^2(\xi_v)=2\beta$ and $v'$ is another
 combing, then by Lemma \ref {l:compare} 
 $$e^2(v')= 2(\alpha(v,v') - \beta) \ . $$
 
 $(3) \Leftarrow (4)$: this is trivial.
 
 $(4)\Rightarrow (5)$: this is trivial.
 
 $(4)\Leftarrow (5)$: this follows from Lemma \ref {evenclass}. 
 
 \smallskip
 
 The equivalence between the five statements is achieved. Now it is
 enough to show that at least one among them holds true. We are going
 to prove that statement (5) holds true:
 
 \begin{proposition} \label{(5)}
   For every combing $v$ of $M$, $\omega^2(\xi_v)=0 \in
   \Hh^2(M;\Z/2\Z)$.
\end{proposition}
 
 \smallskip
 
Equivalently, we have to show that for every compact closed surface
$F$ embedded in $M$, possibly $F$ non orientable, then
$$ \omega^2(\xi_v)\sqcup [F]=0 \in \Z/2\Z $$
that is 
$$ \omega^2(i^*\xi_v)\sqcup [F] = 0$$ where $i: F\to M$ is the
inclusion, and it is not restrictive to assume that $F$ is connected.
\smallskip

Consider the restriction $i^*T(M)$ of the tangent bundle of
$M$ to $F$. Similarly consider $i^*\xi_v$. Then we have the following
two splittings as direct sum:
$$ i^*T(M)= i^*\xi_v \oplus \epsilon^1 = T(F)\oplus \nu $$ where $\nu$
denotes the orthogonal line bundle along $F$, and $\epsilon^1$ is
the restriction to $F$ of the trivial line bundle which has $v$ as
nowhere vanishing section. Here is the key lemma:

\begin{lemma}\label{key}
For every combing $v$ of $M$ and every compact closed embedded
  surface $F$ we have
  $$ \omega^2(i^*\xi_v)\sqcup [F] = \omega^2(T(F))\sqcup [F] + 
  (\omega^1(\det T(F))\cup \omega^1(\nu))\sqcup [F]  \ . $$
\end{lemma}
\smallskip

{\bf Claim:} {\it Lemma \ref{key} $\Rightarrow$ Proposition \ref{(5)}:}
\smallskip

{\it Proof of the Claim:}
If $F$ is {\it orientable}, then the identity of Lemma \ref{key} reduces to
$$ \omega^2(i^*\xi_v)\sqcup [F] = \omega^2(T(F))\sqcup [F] = \chi_2 (F)$$
and we conclude because 
$ \chi(F)$ is even. If $F$ is {\it non
  orientable}, then $F\sim h\PP^2(\R)$, that is the connected sum of
$h$ copies of the projective plane.  As $M$ is orientable, then $\nu$
is isomorphic to the determinant line bundle $\det T(F)$, hence also
in this case
$$\omega^2(i^*\xi_v)\sqcup [F]  =
\chi_2(F) + (\omega^1(F)\cup \omega^1(F)) \sqcup [F] = 2-h + h = 0
\ {\rm mod} (2) \ . $$

\cvd

\medskip

\noindent {\it Proof of Lemma \ref{key}:}
Consider again the two splittings
$$i^*T(M)=i^*\xi_v \oplus \epsilon^1 =  T(F)\oplus \nu  $$
realized geometrically by a field of splittings
$$T_xM = P_x\oplus l(x) = T_xF\oplus \nu(x), \ x\in F$$
where $l(x)$ is the (oriented) line spanned by $v(x)$, while $\nu(x)$
is the (unoriented) line orthogonal to $T_xF$.
Let $s$ be a generic section of $i^*(\xi_v)$, that is a field
of vectors $s=\{s(x) \in P_x\}_{x\in F}$. For every $x\in F$, the direct sum
$T_xF \oplus \nu(x)$ induces the decompositions
$$ s(x) = s_F(x)+ s_\nu(x), \ v(x)= v_F(x) + v_\nu(x) \ . $$
By transversality we can assume that:
\begin{enumerate}
\item $\{s=0\}$ is a finite number of points representing
  $\omega^2(i^*\xi_v)$.
  \item $s_\nu=\{s_\nu(x)\}$ and $v_\nu=\{v_\nu(x)\}$ are generic
    sections of $\nu$, so that both are smooth curves on $F$
    representing $\omega^1(\nu)$ and moreover are transverse to each
    other in $F$, so that their intersection represents
    $\omega^1(\eta)\cup \omega^1(\eta) = \omega^1(\det T(F))\cup
    \omega^1(\nu)$.
\item $\{s=0\}\cap \{v_\nu=0\} = \emptyset$.
\item $s_F=\{s_f(x)\}$ is a generic section of $T(F)$ so that
  $\{s_F=0\}$ is a finite number of points representing
  $\omega^2(T(F))$.
\end{enumerate}
  
For every finite set $X$, let $\#X$ denote the number of its elements mod $(2)$.
Then we have
$$ \omega^2(i^*\xi_v)\sqcup [F] = \#\{s=0\}, \ \omega^2(T(F))\sqcup [F] =\#\{s_F=0\}$$ 
$$(\omega^1(\det T(F))\cup \omega^1(\nu))\sqcup [F] =\#(\{v_\nu=0\}\cap \{s_\nu=0\}) \ . $$
So we have to prove that
$$ \#\{s=0\}= \#\{s_F=0\} + \#(\{v_\nu=0\}\cap \{s_\nu=0\}) \ . $$
On the other hand, obviously
$$ \{s_F=0\} = (\{v_\nu=0\}\cap \{s_F=0\})\amalg (\{v_\nu \neq 0\}\cap \{s_F=0\})\ . $$
We claim that 
$$ \{v_\nu \neq 0\}\cap \{s_F=0\} = \{s=0\}$$
in fact, by item (3) above
$$\{s=0\}= \{v_\nu \neq 0\}\cap \{s=0\} \ ; $$
clearly
$$ \{v_\nu \neq 0\}\cap \{s=0\}\subset  \{v_\nu \neq 0\}\cap \{s_F=0\} \ ;$$
on the other hand if $s(x)\neq 0$, then $s_F(x)\neq 0$, because the projection
$P_x \to T_xF$ is an isomorphism being $v_\nu(x)\neq 0$.
It remains to check that
$$\#(\{v_\nu=0\}\cap \{s_F=0\}) = \#(\{v_\nu=0\}\cap \{s_\nu=0\})
\ . $$ Set $C=\{v_\nu=0\}$ and $j: C \to F$ the inclusion of this
smooth curve; for every $x\in C$, the line $\nu(x)$ is contained in
$P_x$ an we have the splitting as direct sum
$$ P_x = (P_x\cap T_xF)\oplus \nu(x) \ . $$
Hence we have a splitting as direct sum
of line bundles
$$ j^*\xi_v = \lambda \oplus j^*\nu$$ These two lines bundle are
isomorphic to each other; in fact along every component of $C$,
$j^*\xi_v$ is trivial because it is oriented, then the two line
bundles are both trivial or both non trivial; eventually
$$\omega^1(\lambda)\sqcup [C] = \omega^1(j^*\nu)\sqcup [C] \ . $$ We
conclude by noticing that the restriction of $s_F$ and $s_\nu$ are
respectively generic sections of these line bundles.

\cvd

\medskip

The proof of Proposition \ref{(5)}, hence of the main Theorem \ref{3-par}
is now complete.
 
\cvd

\smallskip

\begin{remark}\label{more-compare}{\rm Lemma \ref{l:compare}
    shows in particular that the class $2\alpha(v,v')$ does not
    depend on the choice of the generic pair of combing representatives $v$ and $v'$.
If $\Ff=(v,w,z)$ is a framing of $T(M)$, and $v'$ is any other combing,
then $e^2(v')= 2\alpha(v',v)$.
Thanks to the framing, $v'$ is encoded by a map
$s: M\to S^2$
and it is not hard to verify (do it by excercise) that 
$\alpha (v',v)= s^*(u) \in \Omega^2(M)$,
where $u$ is the usual standard generator of $\Omega^2(S^2)\sim \Z$.
More generally, if $\tilde v$ is another combing encoded by the map say
$\tilde s: M \to S^2$,
then $\alpha(\tilde v, v')= \tilde s^*(u)-s^*(u)$ which by the way shows
that the comparison class itself only depends on the combings as homotopy
classes.}
  \end{remark}

\subsection{Classification of framings}\label{Fram-class}
We provide a classification of the framings on $M$ with respect to a
given reference framing $\Ff_0$. Then any other framing $\Ff$ is
encoded by a map
$$\rho_\Ff: M \to SO(3)$$ considered up to homotopy. The set
$[M,SO(3)]$ can be endowed with a group structure by pointwise
multiplication. As $SO(3)\sim \PP^3(\R)$ there is a natural
homomorphism (see Section \ref {real-LB})
$$\psi: [M,SO(3)]\to \Hh^1(M;\Z/2\Z), \ [h] \to h^*([\PP^2(\R)]) \ . $$
Denote by $p: S^3 \to SO(3)\sim \PP^3(\R)$ the universal covering.
Recall that by Corollary \ref{HOPF} 
$$[M,S^3]\sim \Omega^\Ff_0(M) \sim \Z$$ 
every homotopy class being classified by the common
$\Z$-degree of its representative maps. There is  a natural homomorphism
$$\phi: [M,S^3]\to [M,\PP^3(\R)], \ [f] \to [p\circ f] \ . $$
Finally we can state
\begin{proposition}\label{Fram-exact} The homomorphism sequence
$$ 0\to \Z \xrightarrow{\phi} [M,\PP^3(\R)] \xrightarrow{\psi} \Hh^1(M;\Z/2\Z) \to 0 $$
is exact.
\end{proposition}
\Dim If $p\circ f$ is homotopic to a costant map, then the homotopy
can be lifted to $S^3$, hence $f$ is homotopically trivial and $\phi$
is injective.

Given $g: M \to SO(3)$, $\psi([g])=0$ if and only if $g$ lifts to
$S^3$, hence the kernel of $\psi$ is the image of $\phi$.

We are left to prove that $\psi$ is surjective. We use a spine $\PP_0$
of $M_0$ constructed in Section \ref{spine}.
First one proves that every homomorphism $\alpha:
\pi_1(\PP_0)\to \Z/2\Z$ is induced by a map $j: \PP_0 \to
\PP^2(\R)$. Let $a: (S^1,e) \to (\PP^1(\R), x_0)$, $\PP^1(\R) \subset
\PP^2(\R)$, be a loop which generates $\pi_1(\PP^2(\R)\sim \Z/2\Z$.
We choose a maximal tree $T$ in the singular set of the spine $\PP_0$,
and define $j:{\rm Sing}(\PP_0)\to \PP^2(\R)$ by setting it constantly
equal to $x_0$ on $T$, while on every other edge of the singular set
it is either equal to the constant map or to $a$ according to the
value of $\alpha$ on the loop determined by such an edge. On the
boundary of every region of $\PP_0$ there is an even number of edges
at which $j$ is not constant, hence the map $j$ extends to the whole
of $\PP_0$.  Now we consider $\PP^2(\R)\subset \PP^3(\R)\sim SO(3)$.
The map $j$ extends to $M_0$, and finally to the whole of $M$ because
$\pi_2(SO(3))=0$.

 \cvd
 
 \smallskip

 \subsection{Classification of combings}\label{combings}
 Fix a reference framing $\Ff_0$ of $M$ as above. The set of combings
 of $M$ can be identified with $[M, S^2]\sim \Omega^\Ff_1(M)$ by the
 Pontryagin construction of Chapter \ref{TD-PT}.  We want to make it
 explicit.  There is a natural forgetting projection
$$\pi: \Omega^\Ff_1(M) \to \Omega_1(M) \ . $$ In fact $\pi(v)=
 v^*(u)\in \Omega^2(M)=\Omega_1(M)$, where $u=[y_0]$ is a standard
 generator of $\Omega^2(S^2)$.  We have already remarked that
$$ e^2(\xi_v)= 2 \pi(v) \ . $$ The projection $\pi$ is onto. So we
 have to understand the fibre $\pi^{-1}(x)$ of every $x\in
 \Omega_1(M)$.  If we consider the comparison class $\alpha(v,v')$ as
 the first obstruction in order that the combings coincide, to
 distinguish the combings in a same fibre we have to point out a
 secondary comparison invariant.  Given an oriented framed knot
 $(K,\fG)$ in $M$ which projects to $x$, we can modify the framing to
 $(K,n\fG)$ by adding $n$ twists to the given framing.  This gives a
 transitive action of $\Z$ on such a fibre. We have to understand when
 $(K, \fG)$ and $(K,n\fG)$ represent the same element of
 $\Omega^\Ff_1(M)$.  Assume this is the case, realized by a framed
 surface $S$ in $M\times I$. By taking the double of $M\times I$,
 diffeomorphic to $M\times S^1$, the double $\Sigma$ of $S$ embedded
 therein is an oriented boundaryless surface in $M\times S^1$ such
 that $[\Sigma]\bullet [\Sigma]=n \in \Omega_0(M\times S^1)$.  We
 have $$([\Sigma]- \lambda)\bullet [M\times \{1\}]=0$$ where $\lambda
 = [ K\times S^1]$. Then
$$ ([\Sigma] - \lambda])\bullet ([\Sigma -\lambda])= [\lambda]\bullet [\lambda]=0$$
$$n= 2([\Sigma]-\lambda)\bullet \lambda = [\Sigma - \lambda]\bullet e^2(\xi_v) \ . $$
Then there are two cases:
\smallskip

- $\pi(v)$ is a torsion element, then also $e^2(\xi_v)$ is so, and then $n=0$.

- $e^2(\xi_v)$ is not a torsion element; if $d$ is the biggest integer
such that $\pi(v) = d\beta$ for some $\beta$, then
$$n = 0 \  {\rm mod} \ (2d) \ . $$

Summarizing, we have 

\begin{proposition}\label{comb-class}  (1) Every framing $\Ff_0$ on $M$ determines a
  surjective map
$$ \pi: \Omega_1^\Ff(M)\to \Omega^2(M)$$ such that for every combing
  $v\in \Omega_1^\Ff(M)$, $2\pi(v)= e^2(\xi_v)$.
\smallskip

(2) If $e^2(\xi_v)$ is a torsion element, set $d=0$; then for every
$v, v_0 \in \pi^{-1}(\pi(v_0))$ it is defined a {\it secondary
  comparison invariant} $h(v,v_0)\in \Z/2d\Z = \Z$ such that $v=v_0$
iff and only if $h(v,v_0)=0$.
\smallskip

(3) If $e^2(\xi_v)=2\pi(v)$ is not a torsion element, let $d$ be the
maximum integer such that $\pi(v)=d\beta$ for some $\beta$, then it is
defined a {\it secondary comparison invariant} $h(v,v_0)\in \Z/2d\Z$
such that $v=v_0$ iff and only if $h(v,v_0)=0$.

\end{proposition}

\cvd

\smallskip

\begin{remark}\label{2-torsion}{\rm If
    $\Omega_1(M)$ has no non trivial elements of order $2$, then the
    map $\pi$ does not depend on the choice of the framing $\Ff_0$. On
    the other hand, let $M= \PP^3(\R)$. Fix a trivialization $b:
    U\PP^3(\R)\to \PP^3(\R)\times S^2$ of its unitary tangent bundle
    (associated to a framing $\Ff_0$).  Identify $\PP^3(\R)$ with
    $SO(3)$. Consider a new trivialization $c$ defined by
    $c(b^{-1}(p,y))=(p, py)$ Let $v$ be a combing encoded by a
    constant map with respect to $b$.  Then $\pi_b(v)=0$. On the other
    hand $\pi_c(v)$ is represented by the loop in $SO(3)$ given by the
    rotation in a certain plane, hence it is not trivial.  }
\end{remark}

\smallskip 

Finally we want to outline that the Pontryagin surgery acts transitively.

\begin{proposition}\label{P-surg-trans} Let $v$, $v_0$ be combings of $M$.
Then they are connected by a finite sequence of combing Pontryagin
surgeries.
\end{proposition}
\Dim  Up to Pontryagin
surgery we can assume that the first comparison obstruction vanishes:
$\alpha(v,v_0)=0$.  Fix a reference framing $\Ff_0$ as above. Then
combings are encoded by $[M,S^2] \sim \Omega_1^\Ff(M)$, and we can
assume that $v, v_0$ belong to a same fibre of $\pi:
\Omega^\Ff_1(M)\to \Omega^2(M)$.  It remains to prove that up to
further Pontryagin surgeries say on $v_0$ which stay in the given
fibre, also the second comparison invariant $h(v,v_0)$ vanishes.  As
$\alpha(v,v_0)=0$, we can assume that $v$ and $v_0$
coincide on $M_0 = M \setminus {\rm Int}(B)$ where $B$ is a standard
$3$-disk in a chart of $M$ diffeomorphic to $\R^3$ and moreover they
are constantly equal to a base point $s_0\in S^2$ on $\partial B \sim
S^2$.  
As $B/\partial B \sim S^3$ and is endowed with the base
point $p_0=[\partial B]$, then $v$ and $v_0$ determines two elements
$\bar v, \bar v_0 \in \pi_3(S^2)$.  We know that this last is
isomorphic to $\Z$ and is generated by the Hopf map $\hG$; then $\bar
v=n\hG$, $\bar v_0 = n_0 \hG$.  It is not hard to verify that (with
the notations of Proposition \ref{comb-class})
$$ h(v,v_0)= n-n_0 \ {\rm mod}\ ( 2d)$$ where $d$ only depends on the
given fibre of $\pi$. Then we are essentially reduced to prove that
starting from the map $c_0: S^3 \to S^2, \ c_0(x)=s_0$, for every
$n\in \Z$, we can realize a map $f: S^3 \to S^2$ such that $[f]=[n\hG$
  by means of a finite sequence of Pontryagin surgeries. Assume that
  $B\subset \R^3$ is a suitably big radius; consider the following
  loops in $\R^3$:
$$\gamma_\pm: [0,2\pi]\ni \phi \to 3(0,\cos (\phi), \pm \sin(\phi))\in
  \R^3 \ . $$ Parametrize a tubular neighbourhood of $\gamma_\pm$ as:
$$j_\pm: [0,2]\times[0,2\pi]\times [0,2\pi]\ni (\rho,\theta,\phi) \to$$
$$ \to (3+\rho\cos(\theta))(0,\cos(\phi),\pm
  \sin(\phi))+(\rho\sin(\theta),0,0)\in \R^3 \ . $$ Now, by taking
  convex combinations in $S^2$ on the region $1\leq \rho \leq 2$, we
  can construct a homotopy between the constant field $s_0$ and the
  field
$$e^{(0)}_\pm(j_\pm(\rho,\theta,\phi))=(0,-\sin(\phi),\pm
  \cos(\phi))=\dot \gamma_\pm(\phi)/3 \ . $$ Up to rescaling the
  field, we can apply the Pontryagin surgery along the tube $\{\rho
  \leq 1\}$.  This produces another field $e^{(1)}$ which coincides
  with $e^{(0)}$ outside the tube and is given there by:
$$e^{(1)}_\pm(j_\pm(\rho,\theta,\phi))=$$
$$ = -\cos(\pi\rho)(0,-\sin(\phi),\pm \cos(\phi))-
  \sin(\pi\rho)(\sin(\theta),\cos(\theta)\cos(\phi),\pm
  \cos(\theta)\sin(\phi)) \ . $$ The value $(-1,0,0)$ is regular and
  the inverse image is the curve
  $$\delta_\pm:[0,2\pi]\ni \phi \to j_\pm(1/2,\pi/2,\phi)=(1/2,
  3\cos(\phi), \pm 3\sin(\phi)) \ . $$ By direct computation one
  checks that the framing on $\delta_\pm$ is given by the normal field
  $$\nu_\pm(\phi)= -\frac{\sin(\phi)}{\pi}(1,0,0)-
  \frac{\cos(\phi)}{2}(0,\cos(\phi),\pm \sin(\phi))$$
so that one finally checks that
$$lk(\delta_\pm, \delta_\pm + \nu_\pm) = \mp 1 \ . $$ We can therefore
conclude that starting from the constant field $c_0$, the element of
$\pi_3(S^2)$ which corresponds to the integer $n$ can be realized by
$|n|$ Pontryagin surgeries.

\cvd

\section{What is the simplest proof that $\Omega_3=0$?}\label{simplest-omega3}
We have discussed several proofs that $\Omega_3=0$ and of the
equivalent Lickorish-Wallace theorem on surgery equivalence. By travelling through again
these proofs we can ask about the ``simplest one'' that is, more precisely, the
one with minimal mathematical background. Rohlin's first proof certainly uses
non trivial fact about immersions of $3$-manifolds in $\R^5$.
Lickorish's proof arises as a corollary of an important  result on the surface
mapping class group which nevertheless is rather expensive if one is just interested
about the corollary. The proof in \cite{Rourke} is certaily very simple and
self-contained, provided one assumes Smale theorem.
Then the most basic proof would be obtained by combining one with minimal background
of parallelizability of $3$-manifolds (as in Section \ref{3-parallel}) and the
specialization to the $3$-dimensional case of Proposition \ref{parallelizable}. 

\section{The bordism group of immersed surfaces into a $3$-manifold}\label{2_imm_3}  
Let $S$ be a compact boundaryless surface and $M$ be a connected
boundaryless $3$-manifold.  As usual $[S,M]$ denotes the set of
homotopy classes of maps $f:S\to M$. By using Section \ref {n-in-2n-1}
(see in particular Remark \ref{all-W}) we know that every class
$\alpha \in [S,M]$ contains generic immersions whose local models are
the same as for immersions in $\R^3$ described therein. Generic
immersions in a given homotopy class can be considered up to the finer
relation of {\it regular homotopy}.  This is a particular case of
Smale-Hirsch theory, but the resulting classification is a bit
implicit; several efforts have been made to make it more transparent.
Closer to the themes of the present text, we can consider generic
immersions of compact boundaryless surfaces into a given $3$-manifold
up to a notion of bordism which extends the one of embedded bordism.
In this section we mainly refer to \cite{HH}, \cite{Pi}, \cite{BS}. We
will refer to these papers for details of some proofs. Nevertheless,
we hope to eventually provide a substantial report.

Let us recall first the notion of {\it regular homotopy}.

\begin{definition}\label{regular-homot}
  {\rm Let $\alpha \in [S,M]$; we say that two generic immersions
    $f_0,f_1 : S\to M$ belonging to $\alpha$ are {\it regularly
      homotopic} if there are connected by a homotopy $f_t$, $t\in
    [0,1]$, such that $f_t$ is an immersion for every $t$. We denote
    by $\Rr[S,M]_\alpha$ the set of regular homotopy classes in
    $\alpha$, and by $[f]_r$ the class of a generic immersion
    belonging to $\alpha$.}
\end{definition}

Let us define now the $\iG$-bordism.
  
\begin{definition}\label{bord-imm}
  {\rm Let $f_j:S_j\to M$, $j=0,1$, be generic immersions of surfaces
    into the $3$-manifold $M$. Then $f_0$ is $\iG$-{\it bordant} with
    $f_1$ if there is a $3$-dimensional triad $(W,S_0,S_1)$ and an
    immersion $F: W \to M\times [0,1]$ such that $F\pitchfork M\times
    \{0,1\}$ and $f_j\times \{j\} = F_{|S_j}$, $j=0,1$}.
\end{definition}
\smallskip

Some first remarks:
\smallskip

$\bullet$ As usual, $\iG$-bordism is an equivalence relation.  Denote
by $[f]_\iG$ the equivalence class of a generic immersion $f$.

$\bullet$ If $\phi:S\to S$ is a smooth diffeomorphism, then for every
generic immersion $f: S\to M$, $f$ is $\iG$-bordant with $f\circ
\phi$: the bordism relation incorporates reparametrizations of
surfaces, so that for every immersion $f$, the intrinsic object of
interest is rather its image $f(S)\subset M$ which is a kind of
singular surface in $M$.

$\bullet$ If $f_0, f_1: S\to M$ are connected by a regular homotopy
$F: S\times [0,1] \to M$, then
$$F \times {\rm id}: S\times [0,1] \to M\times [0,1]$$ realizes a
$\iG$-bordism of $f_0$ with $f_1$. Hence in a sense $\iG$-bordism
embodies regular homotopy, but we stress that reparametrization is not
included in the definition of regular homotopy.

$\bullet$ Denote by $\Ii_2(M)$ the set of $\iG$-bordism classes. The
disjoint union defines an abelian {\it semigroup} structure
$(\Ii_2(M),+)$ with $0$ the class of the empty immersion:
$$ [S_1,f_1]_\iG+[S_2,f_2]_\iG = [S_1 \amalg S_2, f_1\amalg f_2]_\iG \ . $$
{\it A priori it is not evident that it is a group},
that is it is not clear how to define the inverses $-[f]_\iG$.

$\bullet$ By using $1$-handles embedded in $M$ we can define a {\it
  connected sum} between immersions $f_1\#f_2: S_1 \# S_2 \to M$ such
that
$$[ S_1 \# S_2, f_1\#f_2]_\iG= [S_1,f_1]_\iG+[S_2,f_2]_\iG \in
\Ii_2(M) \ ; $$ it follows that every class in $\Ii_2(M)$ can be
represented as $[S,f]_\iG$ where $S$ is connected, and the operation
$+$ is induced by $\#$ as well.
\smallskip

We will be mainly concerned with compact $3$-manifolds $M$ and we
distinguish two cases depending on $M$ being orientable or non
orientable.  When $M$ is orientable, a main ingredient of the
discussion will be a certain quadratic enhancement of the intersection
form of surfaces associated to every such an immersion. We will
discuss diffusely the orientable case following \cite{HH}, \cite{Pi},
\cite{BS}.  Later we will give a few indications about the non
orientable one.

An important special case is $M=S^3$ \cite{Pi}. In this case, for
every surface $S$ there is only one homotopy class of maps $f:S \to
S^3$, and via the usual inclusion $\R^3 \subset \R^3\cup \infty =
S^3$, we easily see by transversality that $\Rr[S,S^3]=\Rr[S,\R^3]$
and $\Ii_2(S^3)= \Ii_2(\R^3)$.

\subsection{From immersions in orientable $3$-manifolds 
to quadratic enhancements of surface intersection forms}\label{Pin-}
Let us recall the current setting:

$\bullet$ $M$ is an {\it orientable} connected compact boundaryless $3$-manifold; 

$\bullet$ $S$ is a compact and boundaryless surface, not necessarily
orientable. For a while we will assume also that $S$ is connected.

$\bullet$  $f: S\to M$ is a generic immersion. 

\smallskip

We know that $M$ is parallelizable, 
so let us fix an auxiliary {\it framing} $\Ff$
of $M$, that is a trivialization of the tangent bundle $T(M)$,
considered up to homotopy of framings.  This includes also the choice
of an orientation of $M$.  The framing $\Ff$ can be equivalently
identified with an ordered triple $\Ff=(v,w,z)$ of pointwise linearly
independent tangent vector fields on $M$. By taking an auxiliary
riemannian metric $g$ on $M$, we can also assume that these fields are
pointwise orthonormal.

Let $K$ be a smooth knot in $M$ ($K\sim S^1$). Give $K$ an auxiliary
orientation. The restriction of $v$ along $K$ can be considered as a
map $v: K\to S^2$, then up to homotopy of framings we can assume that
$v$ coincides along $K$ with the positive unitary tangent field on
$K$; thus $\nG_\Ff:=(w,z)$ is along $K$ an ordered couple of pointwise
orthonormal vectors normal to $K$, i.e. it is a {\it normal framing};
it determines a tubular neighbourhood $N(K)$ of $K$ in $M$ equipped
with a trivialization. If $\nG=(w_1,z_1)$ is any other normal framing
along $K$, then by using $\nG_\Ff$ as a reference, we encode $\nG$ by
a map $\rho: K \to SO(2)\sim S^1$ and we associate to $\nG$ the degree
$\phi(\nG):=\deg_\Z(\rho) \in \Z$, so that obviously
$\phi(\nG_\Ff)=0$. This number can be equivalently obtained as follows.
The framing $\nG_\Ff$, that is its first component $w$, determines a longitude 
$l_\Ff$ on $\partial N(K)$
oriented in such a way that the projection onto $K$ is of degree $1$.
Another framing $\nG$ also determines a longitude $l_\nG$.
Then
$$ \phi(\nG)= [l_\nG]\bullet [l_\Ff]\in \Omega_0(\partial N(K))\sim \Z$$
where $\partial N(K)$ is endowed with the boundary orientation. 

We say that $\nG$ differs from $\nG_\Ff$ by $\phi(\nG)$ 
positive or negative twists along $K$. Clearly we can modify
$\nG$ by adding an arbitrary number of twists.
We stipulate that $\nG_\Ff$ is the basic {\it odd} normal
framing of $K$ determined by $\Ff$ and that a normal framing is odd
if it differs from $\nG_\Ff$ by an even number of twists.
Otherwise a framing is {\it even}.   So we have distributed the normal framings to $K$
into two classes; we note that these classes of {\it odd/even framings
do not depend on the choice of the auxiliary orientation on $K$}.
If we apply
this construction to $S^1=\partial D^2 \subset \R^2 \subset \R^3$ with
respect to the standard constant framing of $\R^3$, we realize that
even (resp. odd) normal framings along $S^1$ are characterized by the
property that they cannot (they can) be extended to a framing of the
restriction of $T(\R^3)$ to the spanning $2$-disk $D^2$. The typical
even framing along $S^1$ has as field $w$ the ingoing normals to
$S^1$, tangent to $D^2$; the associated longitude is determined by a collar
of $S^1$ in $D^2$.  

Consider now a smooth circle $C$ on the surface $S$. By trasversality we can
assume that the restriction $f_{|C}$ of the
immersion is an embedding of $C$ onto a knot $K \subset f(S)\subset M$
which extends to an embedding of a tubular neighbourhood $U(C)$ of $C$
in $S$ onto a band $B(K)$ in $f(S)$, with core $K$. We can assume that $B(K)$
is the transverse intersection with $f(S)$ of a neighbourhood $N(K)$ of $K$ in $M$ as above.
We can apply to this knot $K$ the above considerations.  
Give $C$, hence $K$ an auxiliary orientation. Let us orient
$\partial B(K)$ in such a way that  the natural projection onto its core $K$
is a degree-$2$ covering. Fix an {\it even} normal framing $\Ff_e$
along $K$, with associated longitude $l_{\Ff_e}$. 
For every normal framing $\nG$ define as above $\phi_e(\nG)\in \Z$
with respect to $\Ff_e$.  We can consider the integer 
$$[\partial B(K)] \bullet [l_{\Ff_e}]\in \Omega_0(\partial N(K))\sim \Z \ . $$
Then   set
$$q_f(C): = [\partial B(K)] \bullet [l_{\Ff_e}] \ \ {\rm mod}(4) \ . $$ 
If $U(C)$ is annular, then a normal framing, say $u$, of
$C$ in $S$ gives rise to a normal framing $\nG_f=(w,z)$ of $K$ in $M$,
provided that $w$ is the immage of $u$ by the differential of $f$, and
$(v,w,z)$ agrees with the given orientation of $T_xM$ along $K$, where
$v$ is tangent to $K$ as above. Then 
$$ [\partial B(K)] \bullet [l_{\Ff_e}] =2\phi_e(\nG_f) \ . $$
We can say that $q_f(C)$ counts the number mod$(4)$
of  {\it half-twists} the band $B(K)$ makes along its core $K$.
The same interpretation makes sense also when $U(C)$
is a M\"obius strip. In this case  $[\partial B(K)] \bullet [l_{\Ff_e}]$ is odd.
 
\begin{remark}\label{in R3}{\rm If $M=\R^3$, $q_f(C)$ is the linking number mod$(4)$
between $\partial B(K)$ and the core $K$ of the band (co-oriented as before).}
\end{remark}  
 
\smallskip

If $L=\amalg_j C_j$ is the finite disjoint union of smooth circles on $S$, set
$$q_f(L)=\sum_j q_f(C_j) \ . $$  

We have

\begin{lemma}\label{q_f}
(1) The procedure described above well defines a function $q_f$ which
  associates to every finite disjoint union of smooth circles on the
  surface $S$ considered up to ambient isotopy, an element $q_f(C) \in
  \Z/4\Z$.

(2) The function $q_f$ verifies the conditions stated at the end of
  Chapter \ref{TD-SURFACE}; hence by setting for every $\alpha \in
  \eta_1(S)$, $q_f(\alpha):= q_f(C)$, where $C$ is any smooth circle
  on $S$ representing $\alpha$, we well define a quadratic enhancement
  of $(\eta_1(S),\bullet_S)$.
\end{lemma}

As for item (1), it is a bit complicated to show that $q_f(C)$ in
invariant up to ambient isotopy. In fact a generic isotopy between two
copies of $C$ which embed in $M$ by the restriction of $f$, might pass
though non injective immersions and we have to check that this
accidents are immaterial with respect to the value of $q_f$.  As for
item (2), basically one is reduced to a local analysis at a single
crossing point (by the way also the choice of the simplification of
the crossing turns to be immaterial out); this is not very hard. We
left the details as an exercise.

\cvd

\begin{remarks}\label{Spin-Pin}{\rm
1) The choice of the framing $\Ff$ is not immaterial, in the sense
that the quadratic form $q_f$ mights depend on such a choice.
However, it will be immaterial with respect to the statement of main
Theorems \ref{RH} and \ref{i-class}. In the case of $S^3=\R^3\cup
\infty$ we will deal with the unique framing (up to homotopy) of
$\R^3$.

2) The above construction would be placed in a more conceptual
framework in terms of {\it spin structures} on $M$ and induced {\it
  pin$^-$} on $S$. In fact (addressed to a reader who knows this matter), 
 given $f: S\to M$ as above, as $M$ is oriented,
$f^*T(M)= T(S)\oplus \Lambda(S)$ where this last is the determinat
bundle of $S$.  For every spin structure $\Theta$ on $M$, we have the
pull-back spin structure $f^*(\Theta)$ on $f^*T(M)$, and there is a
natural bijection between the spin structures on $T(S)\oplus
\Lambda(S)$ and the pin$^-$ structures on $S$; moreover these last are
in natural bijection with the quadratic enhancements of the
intersection form of $S$.  Rather than the framing $\Ff$ itself, above
we have used the spin structure carried by it. In this framework the
statement of last lemma becomes conceptually clear and even simpler to
prove. However, to our present aims {\it we have preferred the above direct
operative presentation, without introducing the general theory}. A
reader interested to it is mainly addressed to \cite{KT}.

3) The constructions of the present section work as well if $M$ is any 
framed $3$-manifold, not necessarily compact.} 

\end{remarks}
\begin{figure}[ht]
\begin{center}
 \includegraphics[width=12cm]{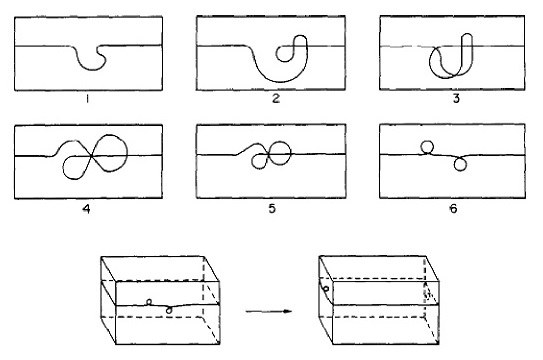}
\caption{\label{kinkbox}  A kink box.} 
\end{center}
\end{figure}

\subsection{Adding kinks}\label{add-kink}
Let $f:S\to M$ be a generic immersion, $S$ connected. Let $C$ be a
smooth circle on $S$ such that $f$ restricts to an embedding of a
small tubular neighbourhood $U(C)$ of $C$ in $S$.  We are going to
modify the immersion $f$ by {\it adding a kink along $C$}. This nice
and crucial construction has been introduced in \cite{HH}.  Denote by
$K=f(C)$, $B(K) = f(U(C))$. $U(C)$ either is an annulus or a
M\"obius strip.  As $M$ is orientable, then any tubular neighbourhood
$N(K)$ of $K$ in $M$ is diffeomorphic to the product $S^1\times
D^2$. As usual we can assume that $\partial N(K)$ is transverse to $f(S)$ and
that $B(K)=N(K)\cap f(S)$. We have two possible models for the pair
$(N(K),B(K))$, depending on $U(C)$ being orientable or not. Consider
$(D^2,X)$ where $X= \{(x_1,x_2)\in D^2; x_1x_2=0\}$. $X=X_1\cup X_2$,
$X_1=\{x_2=0\}$, $X_2=\{x_1=0\}$.

\smallskip

$\bullet$ If $U(C)$ is an annulus then the model for $(N(K), B(K))$
is the mapping cylinder of ${\rm id}: (D^2,X_1)\to (D^2,X_1)$.
\smallskip

$\bullet$ If $U(C)$ is a M\"obius strip then the model for $(N(K),
B(K))$ is the mapping cylinder of ${\rm -id}: (D^2,X_1)\to (D^2,X_1)$.
\smallskip

Accordingly there are two models for adding a kink along $C$. Let
$\tilde X_1$ be the image of an immersion $\alpha: [-1,1]\to D^2$ such
that $\tilde X_1$ is contained in $x_2\geq 0$, is symmetric with
respect to the $x_2$-axis, has one double point, and coincides with
the inclusion of $X_1$ near the end-points. Denote by $-\tilde X_1$
its image by $-{\rm id}$.

If $U(C)$ is an annulus the kink model is very simple: take the
mapping cylinder of ${\rm id}: (D^2,\tilde X_1)\to (D^2,\tilde X_1)$.

If $U(C)$ is a M\"obius strip, then the kink model is more
complicated (see  \cite{HH} pages 104-105); one constructs a so called
``kink box" that is a determined immersion of  the $2$-disk in $D^3$
with one triple point.  A way to visualize this immersion is given in Figure \ref{kinkbox}. 
 First we consider the immersion of $D^2$ into $D^3=D^2\times D^1$ described by the movie
in the first two rows; it results the bottom left-hand picture; then we apply an isotopy 
to it and reach the eventual kink box of the bottom right-hand picture.
We can consider it as an immersion  $X_1\times [-1,1]$ in $D^2\times [-1,1]$ such that
for some $\epsilon >0$:
\begin{enumerate}
\item The image of $X_1\times [-1,-1+\epsilon ]$ coincides with the embedding of
  $\tilde X_1\times [-1, -1+\epsilon]$;
\item The image of $X_1\times [1-\epsilon ,1]$ coincides with the embedding of
  $-\tilde X_1\times [1-\epsilon, 1]$;
\item The image along the boundary of $D^2\times [-1,1]$ coincides with
  the inclusion of $X_1\times [-1,1]$;
\item  There is one triple point in the middle.
 \end{enumerate}
 
Denote by $Z$ the image of this immersion.

Then the kink model is obtained by taking
$$(D^2\times [0,1], Z)/(x_1,x_2,0)\sim (-x_1,-x_2,1) \ . $$ $Z$
projects to a new immersion of $U(C)$ which agrees with $B(K)$ along
the boundary.
\smallskip

By using these models we can modify the given immersion $f:S\to M$
just along $U(C)$ and get $f_C:S\to M$.  It is clear by the
construction that $f_C$ is homotopic to $f$.

\subsection{Determination of $\Rr[S,M]_\alpha$}\label{det-RH}
We give here a first remarkable application of adding kinks.  Let
$f:S\to M$ be a generic immersion as above and $q_f$ the associated
quadratic enhancement of $(\eta_1(S),\bullet_S)$. We know by Lemma
\ref{V-action} that every other enhancement is {\it abstractly} of the
form
$$q'(x)=q_f(x) +2 x\bullet u$$ for a unique $u\in \eta_1(S)$.  Adding
kinks is a natural way to realize it geometrically, by keeping the
homotopy class $\alpha$ of $f$ fixed.  Assume that $u=[C]$, $C$ being
a smooth circle on $S$ to which we can apply the kink construction.
If $C'$ is another smooth circle on $S$ which intersects transversely
$C$ at one point. Denote as above $U(C')$ a small tubulat
neighbourhood of $C'$ in $S$. Then it is immediate that $f(U(C'))$ and
$f_C(U(C'))$ differ by one full twist. Recalling the geometric
definition of $q_f$ in terms of counting half twists mod $(4)$, one
easily realizes that
$$q_{f_C}([C'])=q_f([C])+ 2[C']\bullet[C] \ {\rm mod} (4)$$ 
as desired. 

This result is the key to prove

\begin{theorem}\label{RH}
  Let $S$ be a compact connected boundaryless surface, $\alpha \in
  [S,M]$. Denote by $Q(S)$ the set of quadratic enhancements of
  $(\eta_1(S),\bullet_S)$. Then the map
$$\qG: \Rr[S,M]_\alpha \to Q(S), \  \qG([f]_r) = q_f$$
is well defined and bijective.
\end{theorem}
\Dim An outline: it is not hard to check that it is well defined. We
already know that the map $\qG$ is onto. The proof that it is injective
is non trivial and consists in rephrasing Smale-Hirsch immersion
theory in terms of the quadratic enhancement. This theory
provides a simply transitive action of $\eta_1(S)$ on
$\Rr[S,M]_\alpha$; a main result of \cite {HH} is that this action can
be realized by adding kinks as well as the one on $Q(S)$; so
eventually $\qG$ is an equivariant bijection.

\cvd 

\begin{figure}[ht]
\begin{center}
 \includegraphics[width=12cm]{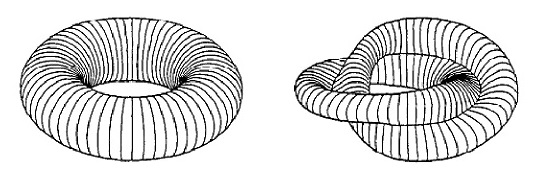}
\caption{\label{tori}  Immersed tori.} 
\end{center}
\end{figure}

\begin{figure}[ht]
\begin{center}
 \includegraphics[width=8cm]{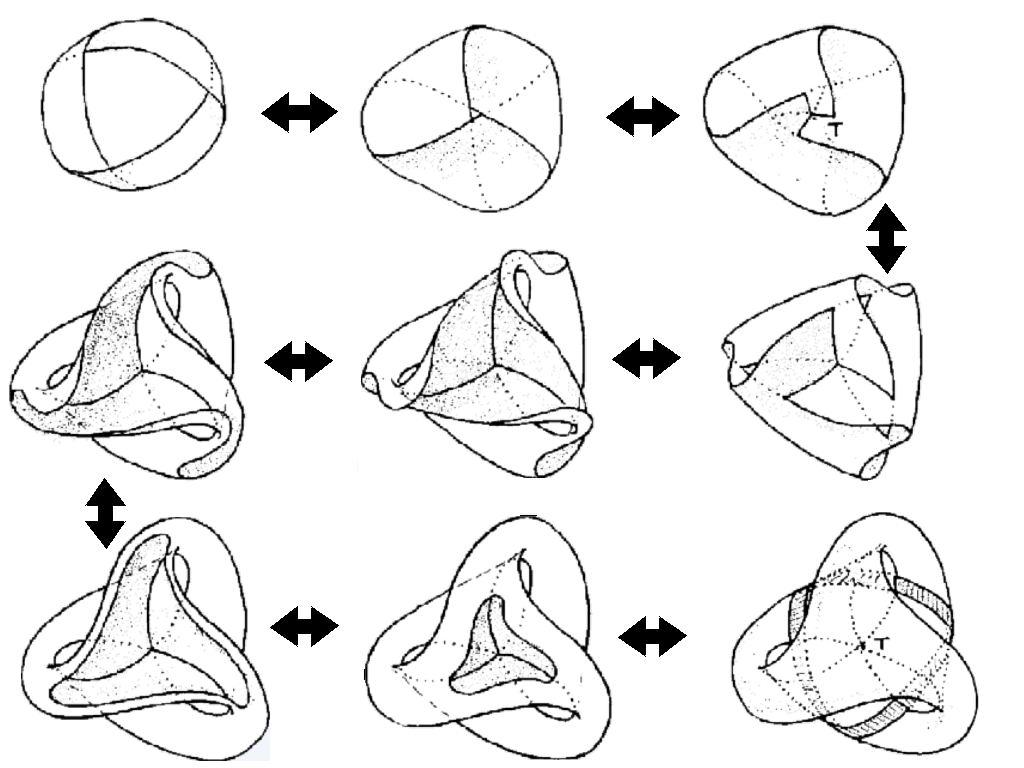}
\caption{\label{Boy}  Boy's surface.} 
\end{center}
\end{figure}

\begin{remarks}\label{basic-in-r3} (Basic immersed surfaces in $\R^3$) \rm{ We refer to \cite{Pi}. 

1) By Theorem \ref{RH}, $\Rr[S^2,\R^3]$ is trivial i.e. it is reduced to one point.
A regular homotopy connecting the standard inclusion $i$ of $S^2$ in $\R^3$
with $-i$ is called  a {\it sphere
eversion} whose surprising existence was discovered by S. Smale \cite{S0}.

2) The elementary surface bricks, besides the sphere,  are the torus $S^1\times S^1$ and the
    projective plane $\PP^2(\R)$.  We denote by $T$ the standard
    embedding of the torus in $\R^3$ bounding a solid torus. We denote
    by $\tilde T$ the immersion obtained by adding a kink along a
    meridian of $T$ and then along the priviliged longitude of $T$
    which bounds a $2$-disk in the complement of the solid
    torus. These realize the two quadratic enhancements of
    $(\eta_1(S^1\times S^1),\bullet)$ (up to isometry) - $T$ and
    $\tilde T$ are illustrated in Figure \ref{tori}.
    
There is a famous immersion of the projective plane with one triple point
called {\it Boy's surface} - see for instance the body and the references of \cite{Ap}).
Figure \ref{Boy} suggests how to construct it. Such an immersion denoted by 
$B$ and $\bar B$ the mirror of $B$, that is $B$
composed with a reflection at a hyperplane of $\R^3$, realize
the two quadratic enhancements of $(\eta_1(\PP^2(\R),\bullet)$. }
   
\end{remarks}

\subsection {Determination of $(\Ii_2(M),+)$}\label{I_det}
First we will point out a few invariants up to $\iG$-bordism.
\smallskip

{\bf The Arf-Brown invariant.} Let $f:S\to M$ be a generic immersion,
$S$ connected, with the associated $q_f$.  Accordingly with section
\ref{quadratic}, we can consider the Arf-Brown multiplicative
invariant
$$ \gamma(f):= \gamma(q_f) \in U_8 \sim \Z/8\Z $$ where for simplicity
we have written $\gamma(q_f)$ instead of $\gamma(S,\bullet_S,q_f)$.
If $f: S \to M$, where $S=\amalg_j S_j$ is union of several connected
components, then set
$$\gamma(f):= \prod_j \gamma(f_j)$$
where $f_j= f_{|S_j}$.
We have
\begin{lemma}\label{A-B_inv} Let $f_j: S_j\to M$ be generic immersions, $j=0,1$. If 
$\ [f_0]_\iG = [f_1]_\iG$, then $\gamma(q_{f_0})=\gamma(q_{f_1})$.
\end{lemma}
\Dim Let $(W,S_0,S_1)$, $F:W \to M\times [0,1]$ be as in Definition
\ref {bord-imm}, and let $t: M\times [0,1]\to [0,1]$ be the
projection. Without loss of generality we can assume that $t\circ F$
is a Morse function on the triad. Then consider the possible accidents
when passing though a critical point of $t\circ F$. Modifications
occur locally in a chart of $M$ at the critical point. We use the
notations of Remark \ref{basic-in-r3}. At local minima/maxima a new
spherical component appears/disappears. For the other kinds of
critical point, there are three possibilities: 

- one performs the
immersed connected sum of two components of the surface; 

- one performs the connected sum with either a standard torus $T$ or a Klein bottle
immersion $B \# \bar B$. 
\smallskip

\noindent In every case the value of $\gamma$ does not
change (for all details one can see \cite{Pi}, pp. 432-433).

\cvd

\smallskip

So we have detected a first main $U_8$-valued invariant
$\gamma([f]_\iG)$ defined on $\Ii_2(M)$.
\smallskip

{\it From now on we will use the standard isomorphism $U_8 \sim (\Z/8\Z,+)$ and hence adopt the additive notation.}
\smallskip

{\bf Other invariants.}  Let $f: S\to M$ be a generic immersion ($S$
not necessarily connected).  It is obvious just by forgetting part of
the structure of $[f]_\iG$, that $[f]=[S,f]\in \eta_2(M)$ is invariant
under $\iG$-bordism.  Recall the quotient module $\Hh^1(M,\Z/2\Z)$ of
$\eta^1(M)= \eta_2(M)$ defined in Corollary \ref {eta1}; recall also
that the cup product $\sqcup$ descends to this quotient with values in
$\eta^2(M)=\eta_1(M)$. Keep the notation $[f]$ for its image in
$\Hh^1(M,\Z/2\Z)$.

Denote by $\Sigma\subset S$ the non-injectivity locus of $f$.  We
claim that the image $\Sigma_f:=f(\Sigma)$ determines an element
$[\Sigma_f] \in \eta_1(M)$.  In fact the components of
$f^{-1}(\Sigma_f)$ are of two kinds:

1) they are member of a couple $\tilde C = C \amalg C'$ such that
$f(C)=f(C')$ and $f$ is generically $1-1$ on such a $C$.  In such a
case select one $C$ in each couple;

2) Components $\tilde C$ such that $\tilde C= f^{-1}(f(\tilde C))$ and in such a case $f$ is generically 2-1 on $\tilde C$.
\smallskip

Then select one component $C$ in every couple $\tilde C$ of the first kind; for the second kind
one finds a quotient $C$ of $\tilde C$ such that $f$ induces a map (we keep the name)
$f: C\to M$, such that $f(C)=f(\tilde C)$ and $f$ is generically 1-1.
Then set
$$[\Sigma_f]:= \sum_{\tilde C}  [C,f] \in \eta_1(M) \ . $$

The triple points of $f(S)$ determines a class $t_f \in \eta_0(M)\sim
\Z/2\Z$.  We have

\begin{lemma}\label{hom-inv} If  $[f_0]_\iG = [f_1]_\iG$, then $[\Sigma_{f_0}]=[\Sigma_{f_1}] \in \Hh^1(M,\Z)$
and $t_{f_0}=t_{f_1}\in \eta_0(M)$.
\end{lemma}
\Dim Let $(W,S_0,S_1)$, $F:W \to M\times [0,1]$ be as in Definition
\ref {bord-imm}. We can assume that also $F$ is generic. Then
$F(\Sigma_F)$ is a kind of singular surface properly embedded into
$M\times [0,1]$ such that $F(\Sigma_ F) \cap (M\times \{0,1\})=
f_0(\Sigma_0) \amalg f_1(\Sigma_1)$; by using the regular surface
$F^{-1}(\Sigma_F))$ we can explicitly define a triad which connects
the sum of the components that form $[\Sigma_{f_0}]$ and
$[\Sigma_{f_1}]$ respectively.  Similarly for the triple points.

\cvd

\smallskip

Consider the product
$$\Gamma(M)= \eta_1(M) \times \Hh^1(M,\Z/2\Z) \times \Z/8\Z$$
endowed with the {\it twisted} group structure defined by the operation:
$$(\delta, h, a)+(\delta',h',a'):= (\delta+\delta'+ h \sqcup h', h+h', a+a') \ . $$
We can state now the main result of this section.
\begin{theorem}\label{i-class} The map $\psi: \Ii_2(M) \to \Gamma(M)$ well defined by
$$ [f]_\iG \to ([\Sigma_f], [f], \gamma(f) )$$
is a semigroup isomorphism. In particular the semigroup $(\Ii_2(M),+)$ is a group.
Moreover, the invariant $t_{[f]_\iG}$ is determined by the others.
\end{theorem}

\smallskip

The rest of this section is occupied by the proof of Theorem \ref{i-class}.
It is immediate that $\phi$ is a semigroup homomorphism.

\medskip

{\bf The $3$-sphere.} If $M=S^3$, Theorem \ref{i-class} specializes to

\begin{theorem}\label{i-sfera} The map $\phi: \Ii_2(S^3)\to \Z/8\Z$, $\phi([f]_\iG) = \gamma(f)$
is a group isomorphism.
\end{theorem}
\Dim This is a main result of \cite{Pi} to which we refer for all
details. We can use $\R^3$ instead of $S^3$.  Note that we know a
priori that $\Ii_2(\R^3)$ is a group : inverses are obtained by mirror image
along a hyperplane.  By using connected
sums (or disjoint unions) of the basic immersed surfaces of Remark
\ref{basic-in-r3} it is easy to prove that $\phi$ is onto. By
Proposition \ref{RH} (and the classification of surfaces) one realizes
that every generic immersion $f:S\to \R^3$ is regularly homotopic to a
connected sum of several copies of the standard embedding $T$ and one
among the following eight surfaces
$$B, \ \bar B, \ K_0, \ K_+, \ K_-, \ K_+ \# B, \ K_0 \# \tilde T,
\ K_- \# \bar B $$ where $K_0= B \# \bar B$, $K_+= B \# B$, $K_-= \bar
B \# \bar B$.  Up to $\iG$-bordism the $T$-components are immaterial
and one eventually gets that eight explicit generators suffice and
this achieves the desired bijection onto $\Z/8\Z$.

\cvd

\smallskip

{\bf The map $\psi$ is onto.} Let us prove now in general that the map $\psi$ is surjective.

\begin{lemma}\label{psi-onto}  The map $\psi: \Ii_2(M) \to \Gamma(M)$ is onto.
\end{lemma}
\Dim As $M = M\# S^3$, we see that $\Ii_2(M)$ contains the subgroup
$\Ii_2(S^3)$; it consists of the classes with a representative
contained in a $3$-disk of $M$.

It contains also the {\it subset} $E(M)$ given by the classes which are represented by
{\it embedded} surfaces. By the description of $\Hh^1(M,\Z/2\Z)$ as the embedded surfaces
in $M$ up to embedded bordism, we see that $E(M)$ is in fact the image of $\Hh^1(M,\Z/2\Z)$
in $\Ii_2(M)$ by a natural quotient map. 

Let $(\delta, h, a)\in \Gamma(M)$. Represent $h$ by an embedding $e: S\to M$.
Represent $\delta$ by a knot $K$ in $M$. Consider the boundary $\Tt\sim S^1\times S^1$
of a tubular neighbourhood of $K$ in $M$. Add a kink along a longitude $K'$ of $K$  on $\Tt$
and get a generic immersion $j: \Tt\to M$. By construction $\delta = [\Sigma_j]$,
while $[j]=0\in \Hh^1(M,\Z/2\Z)$, hence $[j]\sqcup [e]=0$.  By the elementary fact that $\gamma$ 
is onto in the case of $S^3$,
there is $s:S\to S^3$ such that $\gamma(s)= a- \gamma(e) + \gamma(j)$.
Clearly $[s]=0$ and $[\Sigma_s]=0$. 
Finally 
$$\psi([j]_\iG + [e]_\iG + [s]_\iG)=(\delta, h, a) \ . $$  

\cvd

\smallskip

{\bf A normal decomposition of $\iG$-bordism classes.}  Now the idea
is that every $[f]_\iG$ admits  a certain {\it normal decomposition}
modelled on the classes used to prove that $\psi$ is
onto. Precisely we have the following key proposition.

\begin{proposition}\label{normal-dec} Every $[f]_\iG$ can be represented by a sum
$$ [f]_\iG = [j]_\iG + [e]_\iG + [s]_\iG$$ where $j: \Tt\to M$ is
  obtained by adding a kink along a longitude $K'$ on the boundary
  $\Tt$ of a tubular neighbourhood of a knot $K$ in $M$, $[e]_\iG \in
  E(M)$, $[s]_\iG \in \Ii_2(\R^3)$ where $\R^3$ is a chart of
  $M$. Moreover, we can choose the decomposition in such a way that
  $q_j(K')=0$ hence so that $\gamma(j)=0$.
\end{proposition}
\Dim We will proceed in several steps. We adopt the notations of
Remark \ref{basic-in-r3}, in particular $B$, $\bar B$ are the two
versions of Boy's surface.
\smallskip

{\bf Step 1.} {\it $[f]_\iG = [f']_\iG + [s]_\iG$, where $f'$ has no triple
points and $[s]_\iG \in \Ii_2(\R^3)$.}  

Notice that $K_0=B \# \bar B$
is regularly homotopic to the usual immersion of the Klein bottle in
$\R^3$ without triple points (and a plane of symmetry) and recall that
$[K_0]_\iG =0$.  Similarly if $x_0$ is a triple point of $f$, either
$f \# B$ or $f \# \bar B$ is regularly homotopic to $\tilde f$ with
one triple point less than $f$, and either $[f]_\iG = [\tilde f]_\iG +
[\bar B]_\iG$ or $[f]_\iG = [\tilde f]_\iG + [B]_\iG$.  So the step is
achieved by induction on the number of triple points.
\smallskip

Then the double line locus $\Sigma_{f'}$ consists of the disjoint
union of a finite number of embedded circles in $M$. If $K$ is such a
circle, then it has a neighbourhood in $f(S)$ which is a bundle over
$K$, sub-bundle of a tubular neighbourhood of $K$ in $M$, with fibre
isomorphic to $X=\{ (x_1,x_2)\in D^2; \ x_1x_2=0\}$.  We can count the
number mod$(4)$ of {\it quarter turns} this configuration does when
moving along $K$ .  Denote it by $l(K)\in \Z/4\Z$; it characterizes
the bundle. The cases $l(K)=0, 2$ correspond to the situation where
$f': (f')^{-1}(K)\to K$ is a trivial double covering; if $l(K)=0$ then
the two components of this inverse image have annular tubular
neighbourhoods in $S'$; if $l(K)=2$, both have M\"obius strip
neighbourhoods. The cases $l(K)=1,3$ correspond to a non trivial
double covering.  
\smallskip

{\bf Step 2.}  {\it $[f]_\iG = [f']_\iG + [s]_\iG$ as in Step 1 and
moreover we can require that $\Sigma_{f'}$ is connected.} 

If $\Sigma_{f'}$ is not connected, there are two components $K$ and $K'$
and points $p\in K$, $p'\in K'$ belonging to the closure of a same
connected component of $M\setminus {\rm Im}(f')$.  So there is a
smooth simple arc $\sigma$ in $M$, connecting $p$ and $p'$ and without
any further intersection with ${\rm Im}(f')$.  Locally in chart of $M$
at $p$, the image of $f'$ looks like two transverse planes $P_1$ and
$P_2$.  Similarly at $p'$, with planes $P'_1$ and $P'_2$.  Remove from
the image of $f'$ the intersection, say $B_p$, of the interior of a
small $3$-ball centred at $p$, with transverse boundary spheres.  The
closure of $B_p$ is the union $D_1\cup D_2$ of two $2$-disks,
$D_j\subset P_j$, $j=1,2$, which intersect transversely at a segment
of $K$. Do similarly at $p'$.  Possibly up to reordering the planes,
we can attach two embedded $1$-handles $H_j$ along the arc $\sigma$,
$j=1,2$, with attaching tube $T_{a,j}=D_j\cup D'_j$, and transverse
$b$-tubes such that $T_{b,1}\pitchfork T_{b,2}$ consists of two
disjoint double arcs having as endpoints the four points of $(D_1\cap
D_2)\cup (D'_1\cap D'_2)$. Ultimately, (up to corner smoothing) we get
the immersed surface
$${\rm Im}(\tilde f):= ({\rm Im}(f') \setminus (B_p\cup B_{p'}))\cup (T_{b,1}\cup T_{b,2})$$
which by construction is $\iG$-bordant with $f'$, alike $f'$ has no triple points,
the two knots $K$ and $K'$ of $\Sigma_{f'}$ have fused into one knot $K''$ of $\Sigma_{\tilde f}$, so that
this last has one compent less. The step is achieved 
by induction on the number of components of $\Sigma_{f'}$.
We stress that by the above construction we have furthermore that
$$ l(K")= l(K)+ l(K') \ . $$
\smallskip

{\bf Step 3.}  {\it Let  $[f]_\iG = [f']_\iG + [s]_\iG$ be as in Step 2 (i.e. with $\Sigma_{f'}=K$ connected) 
and assume that $l(K)= 0, 2$. Then it is not restrictive to assume that $l(K)=0$}.

By using the results about the group $\Ii_2(\R^3)$ we see that there is an immersion $s_0$ of the Klein
bottle in a chart of $M$, without triple points and having connected $\Sigma_{s_0}=K_0$  
such that $l(K_0)=2$.  Take
$$[ f' \# s_0]_\iG + [s]_\iG - [s_0]_\iG = [f]_\iG$$
and apply Step 2 to  $f' \# s_0$. This achieves the step.
\smallskip

Let  $[f]_\iG = [f']_\iG + [s]_\iG$ be as in Step 3, so that $l(K)=0$.
Set $q_f(K):= q_f(C)$, where $C$ is a component of $(f')^{-1}(K)$.
It is well defined and either $q_f(K)=0$ or $q_f(K)=2$.

{\bf Step 4.}  {\it Let  $[f]_\iG = [f']_\iG + [s]_\iG$ be as in Step 3, so that $l(K)=0$.
Then it is not restrictive to assume that $q_f(K)=0$}

There is an immersion $s_1$ of the torus in a chart of $M$, without triple points
and with connected $\Sigma_{s_1}=K_1$ such that $l(K_1)=0$ and $q_{s_0}(K_1)=2$.
If $q_f(K)=2$, take
$$[ f' \# s_0]_\iG + [s]_\iG - [s_0]_\iG = [f]_\iG$$
and apply Step 2 to  $f' \# s_0$. This achieves the step.
\smallskip

{\bf Step 5.} {\it   Let  $[f]_\iG = [f']_\iG + [s]_\iG$ be as in Step 2 (i.e. with $\Sigma_{f'}=K$ connected) 
and assume that $l(K)= 0, 2$. Then Proposition \ref{normal-dec}  holds in this case.}

By Steps 3 and 4 we can assume that $l(K)=0$ and $q_f(K)=0$. Perform a {\it Rohlin surgery} along
$K$ (recall Section \ref{n-in-2n-1-emb}). This splits $f'$ in two disjoint immersed surfaces:
an embedding $e$ and the immersion $j$ of torus having as immage a product sub-bundle of (the interior of) a tubular
neighbourhood  $N(K)$ of $K$ in $M$ with fibre a lemniscate; the germ of $j$ along $K$ equals the germ of $f'$. 
It is easy to see that $j$ is obtained by adding a kink along a longitute $C$ on the boundary 
$\Tt$ of a smaller tubular neighbourhood $N'(K)\subset N(K)$ and that $q_\Tt(C)=q_f(K)=0$.
By construction $[f]_\iG = [j]_\iG + [e]_\iG + [s]_\iG$.
The Proposition is proved under such restrictive hypotheses.
 
\smallskip

\begin{figure}[ht]
\begin{center}
 \includegraphics[width=7cm]{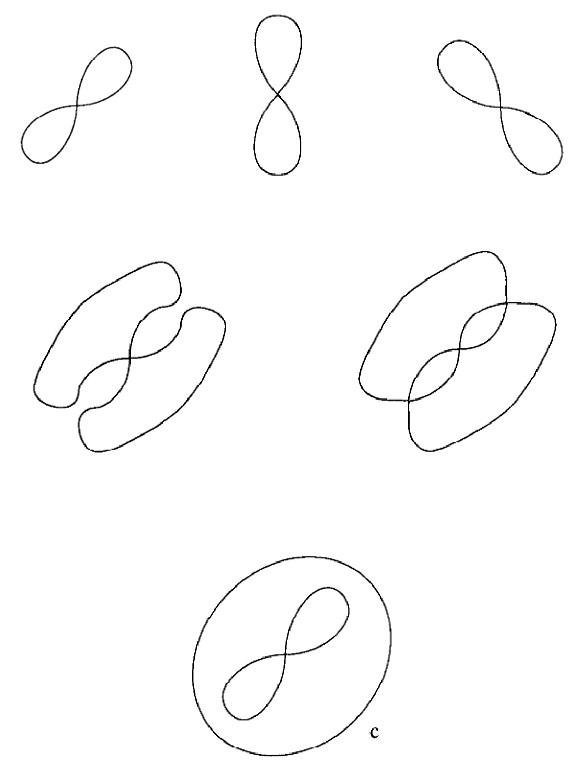}
\caption{\label{lastsurface}  An auxiliary immersed surface.} 
\end{center}
\end{figure}

To proceed we need the following lemma.

\begin{lemma} \label{s2} There is an immersion $s_2$ in $\R^3$ of a surface $F$ of Euler-Poincar\'e characteristic
$\chi(F)=-1$  such that:

1) $s_2$ has one triple point;

2) $\Sigma_{s_2}$ consists of the union of a smooth circle $K_2$ endowed with an  $X$-bundle neighbouhood  in the image of $s_2$ 
such that $l(K_2)=1$,
and a lemniscate in a $2$-disk $D$ contained in the image of $s_2$, intersecting $K$ at the triple point; $D$ is transverse
to $K$ and  the germ of the lemniscate at the triple point is a fibre of the $X$-bundle along $K$.
\end{lemma}
\smallskip

\Dim First we construct an immersion of a surface $G$ with boundary in $D^2\times D^1$. This is given by the movie
of Figure \ref{lastsurface}. Note that at the initial time $t=-1$ and at the final time $t=1$ of the movie we see two copies of a same
lemniscate $L$; in the final configuration $L$ is encircled by a smooth circle $c$. 
Finally we complete $G$ by filling the curve
$c$ by a $2$-disk, and identifying by the identity of $\R^2$ the two copies of $L$ over $-1$ and $1$ respectively. 
One readily check that this is the image of an immersion of a surface $F$ with the required properties.

\cvd

\smallskip

We denote by $\bar s_2$ the mirror image of the immersion $s_2$ as above.

\smallskip 

{\bf Step 6.} {\it Proposition \ref{normal-dec} holds in full generality.}

It remains to prove it when  $[f]_\iG = [f']_\iG + [s]_\iG$ is again as in Step 2, but we assume now that 
$l(K)= 1, 3$. Let $l(K)=1$. By realizing $s_2$ in a chart of $M$, 
take
$$[ f' \# \bar s_2]_\iG + [s]_\iG + [s_2]_\iG = [f]_\iG$$
and apply Step 2 to  $f' \# \bar s_2$. In this way we reach a 
decomposition $[f]_\iG = [f"]_\iG + [s']_\iG$, where $\Sigma_{f"}$ is qualitatively similar to 
the one of $s_2$, that is it consists of the union of a smooth circle $K"$ endowed with an  $X$-bundle neighbouhood  in the image of $f"$ 
and a lemniscate in a $2$-disk $D$ contained in the image of $f"$, intersecting $K"$ at one triple point; $D$ is transverse
to $K"$ and  the germ of the lemniscate at the triple point is a fibre of the $X$-bundle along $K"$. Moreover, $l(K")=0$.
By applying Step 4, we can also assume that $q_{f"}(K")=0$.
Now, although there is a triple point, we can apply Steps 5 along $K"$. This produces a decomposition
of the form $[f]_\iG = [j]_\iG + [g]_\iG + [s']_\iG$ where $[j]_\iG$ has the required final properties, while
$\Sigma_g$ is contained in $D$ and consists of the union of a lemniscate fibre of $j$ an two further simple double circles.
We can eliminate such circle by applying again Steps 4, 5; eventually we get a required decomposition 
$$ [f]_\iG = [j]_\iG + [e]_\iG + [s"]_\iG \ . $$
If at the beginning $l(K)=3$, we manage similarly by exchanging the roles of $s_2$ and $\bar s_2$ respectively.
This achieves Step 6.

\begin{remark}{\it We stress that when $l(K)=0,2$, the images of $j$ and $e$ in the normal decomposition obtained above
are disjoint. When $l(K)=1,3$, they intersect producing one triple point. In the first case $ [\Sigma_f]\bullet [f]=0 \in \eta_0(M)\sim \Z/2\Z$,
in the second    $ [\Sigma_f]\bullet [f]=1$.}
\end{remark}

\smallskip

The proof of Proposition \ref{normal-dec} is now complete.

\cvd

\smallskip

{\bf The map $\psi$ is injective.}
We have
\begin{lemma}\label{psi-inj} The map $\psi: \Ii_2(M)\to \Gamma(M)$ is injective.
\end{lemma}
\Dim  We can use normal decompositions of $\iG$-bordism classes.
Assume that 
$$\psi([j]_\iG + [e]_\iG + [s]_\iG)= \psi([j']_\iG + [e']_\iG + [s']_\iG) \ . $$
As $[e]=[e']\in \Hh^1(M,\Z/2\Z)$ then they are bordant by means of an embedded bordism, hence
$\gamma(e)=\gamma(e')$. As $\gamma(j)=\gamma(j')=0$, then $\gamma(s)=\gamma(s')$ and by 
Theorem \ref{i-sfera}, we have $[s]_\iG = [s']_\iG$. It remains to prove that  $[j]_\iG = [j']_\iG$. 
Now $[j]_\iG + [j']_\iG = [j \# j']_\iG$
and this last can be obtained from the embedding $\Tt  \# \Tt'$ by adding kinks along two disjoint 
circles $K'$, $K"$ at which the quadratic enhancement vanishes. Let $C$ be a smooth circle
on $\Tt  \# \Tt'$ such that $[K']+[K"]=[C]\in \eta_1(\Tt  \# \Tt')$. Then  up to regular homotopy $j \# j'$ 
can be obtained by adding a kink to $\Tt  \# \Tt'$ along $C$. It follows from the hypotheses that $[C]=0 \in \eta_1(M)$
and that the quadratic enhancement of $\Tt  \#  \Tt'$ vanishes on $C$. We claim that 
in such a situation $[j \# j']_\iG=0$. As the same considerations hold for $ [j \# j]_\iG$,
we will eventually concude $[j]_\iG=-[j]_\iG$ and hence  that $[j]_\iG = [j']_\iG$ as desired.

We need the following lemma.

\begin{lemma}\label{even-framing} Let $F$ be a compact  surface with connected boundary
embedded into a framed $3$-manifold $N$ ($F$ might be non orientable and $N$ non compact).
Then the normal framing of $C=\partial F$ determined by a collar in $F$ is even with respect to the ambient
framing.
\end{lemma}
\Dim  We can extend the embedding of $F$ to a generic immersion of the double $D(F)$ of $F$ into $N$.
If $F$ is orientable, up to corner smoothig, we can take the boundary of a tubular neighbourhood 
of $F$ in $N$; if $F$ is not orientable, we can take an immersion which looks like in the orientable case
along the boundary and have double lines in the interior of $F$. We use the ambient framing
to define a quadratic enhancement  $q_{D(F)}$ of the intersection form of the double. As 
$[C]=0 \in \eta_1(D(F))$, then $q_{D(F)}(C)=0$. This means exactly that the collar normal framing
is even.

\cvd

\smallskip  

To simplify the notations, denote by $f:S\to M$ the embedding $\Tt  \#  \Tt'$, so that $q_f(C)=0$.
As $[C]=0\in \eta_1(M)$, then there is a (possibly non orientable) embedded Seifert surface $F\subset M$ such that
$\partial F= C$. Apply Lemma \ref{even-framing} to $F$. As also $q_f(C)=0$, then both the normal framings of
$C$ determined by a tubular neighbourhood in $S$ and by a collar in $F$ respectively differ to each othet by an even number 
of twists. It follows that we can  ``roll up''  $F$ in a tubular neighbourhhod $U$ of $C$
in $M$, in such a way that $F$ is transverse to $S$ along $C$, and intersects transversely $S$ outside $U$.

Assume first that $F=D$ is a $2$-disk.  Let $\tau$ be a Dehn twist on $S$ along $C$.
For every $\alpha \in \eta_1(S)$, 
$$\tau_*(\alpha)=\alpha + ([C]\bullet \alpha)[C] \ . $$
As $q_f(C)=0$, by recalling the geometric definition of $q_f$, we readily see that
$$ q_{f_C}=q_{f\circ \tau} \ . $$
We claim that $f_C$ and $f\circ \tau$ are homotopic (equivalently $f$ and $f\circ \tau$ are homotopic).
To prove the last statement, let $U$ denote now a tubular neighbourhhod of $C$ in $S$; there is a natural
map $h: U \to D$ which realizes a homotopy to a point of $f_{|C}$. Then $f$ and $f\circ \tau$ are homotopic
to maps $f'$ and $f"$ such that:

- they coincide outside $U$;

- $f'_{U}$ and $f"_{U}$ factor though $h$.

Since $D$ is contractible they are homotopic relatively to $S\setminus U$. 
By Theorem \ref{RH}, $[f_C]_r = [f\circ \tau]_r$, hence $[f_C]_i = [f]_i$.
\smallskip

It remains to reduce to such a special case $F=D$. To this aim, consider a generic Morse function 
$$r: F \to [0,1]$$
such that $r^{-1}(0)= C$ and $r$ has no minima and only one maximum. Then we can find a non critical value $\lambda \in [0,1)$
such that $D=r^{-1}([\lambda, 1])$ is a $2$-disk embedded in $M$ with boundary denoted by $\hat C$.
By following the level lines of $r$ between $0$ and $\lambda$ we can modify $(S,f,C)$  into a $(\hat S, \hat f, \hat C)$ 
such that $[f]_\iG = [\hat f]_\iG$. Between two consecutive critical values we can extend the isotopy between level lines 
to a diffeotopy of $M$. At a critical point the analysis is local in a chart of $M$: the critical level of $r$ containing a crossing
point $x_0$ is contained in a ``critical" surface $S'$ with one isolated singular point at $x_0$ isomorphic to a cone centred at $x_0$ 
and bases at two disjoint circles; $F$ and $S'$ intersect along such a critical level, transversely outside $x_0$. 
By such a local analysis one realizes that $q_{\hat f}(\hat C)=0$ and that $[f_C]\iG = [\hat f_{\hat C}]_\iG$.
So we have reduced to the special case $F=D$ and the Lemma is proved.

\cvd

\smallskip

The proof of the main Theorem \ref{i-class} is now complete.

\cvd

\smallskip 

\subsection{More quasi-framing}\label{more-QF}
Now we give a further proof of  the existence of a quasi-framing
on $M$
based on some constructions established in Section  \ref{Pin-}.
 \smallskip
 
 By contradiction, assume that there is $v$ such that 
 $$\beta:= \omega^2(\xi_v)\neq 0 \  .$$
 Let $K$ be an oriented knot in $M$ which represents $e^2(\xi_v)$. By forgetting the orientation,
 $K$ represents $\omega^2(\xi_v)$. Then it follows from the hypotheses that (see Section \ref {relevantcobord}):
 \begin{enumerate}
 \item There is a framing $\Ff'$ of $T(M)$ over $M\setminus K$.
 \item There is  a (possibly non orientable) compact boundaryless surface $F$ embedded into $M$
 such that $F\pitchfork K$ at exactly one points.
 \end{enumerate}
 \smallskip
 
 Let $N(K)\sim S^1 \times D^2$ be a
 tubular neighbourhood of
 $K$ in $M$ transverse to $F$. By removing the interior of $N(K)$
 from $F$, we can assume to get a surface $F_0$ properly embedded in
 $$M':=M\setminus {\rm Int}(N(K))$$
 such that $C:=\partial F_0$ is a meridian of
 $\partial N(K)$ bounding a fibre $D$ of $N(K)$. 
 As in Section \ref{Pin-}, we can use the framing $\Ff'$
 to construct a quadratic enhancement of the intersection form of
 every surface immersed into $M'$. By Lemma \ref{even-framing}, we see that
 the normal framing of $C$ determined by a collar of $C$ in $F_0$ - 
 equivalently by a collar of $C$ in the meridian disk $D$ -
 is {\it even} with respect to $\Ff'$, and it is also {\it even} with respect to a framing of a $3$-ball
 containing $D$. Then the normal framing determined by $\Ff'$ is odd within the $3$-ball,
consequently  $\Ff'$ can be extended over a
 neighbourhood $U\sim D \times [-1,1]$ of $D$ in $N(K)$; as  the
 closure of $N(K)\setminus U$ is a closed $3$-ball, we have eventually
 obtained an almost-framing of $M$. By Lemma \ref {Q-F}   and $(1) \Rightarrow (5)$ of Theorem \ref{3-par}, 
 we get that $\omega^2(\xi_v)=0$ againts the assumption that $\omega^2(\xi_v)\neq 0$. This is a contradiction.

\cvd

\smallskip

\subsection{On $\Ii_2(M)$ for a non orientable $3$-manifold}\label{i-non-or} 
 If $M$ is non orientable  the structure of $\Ii_2(M)$ is eventually simpler.
 Consider the product 
 $$ \Gamma_0(M)= \eta_1(M) \times \Hh^1(M,\Z/2\Z)\times \Z/2\Z$$
 with the twisted group structure given by the operation
 $$(\delta, h, a)+(\delta',h',a'):= (\delta+\delta'+ h \sqcup h', h+h', a+a') \ . $$
 Then we have  \cite{G}
 \begin{theorem}\label{i-class2} Let $M$ be a non orientable compact connected boundaryless $3$-manifold.
 The map 
 $$\psi_0: \Ii_2(M)\to \Gamma_0(M), \ \psi_0([S,f]_\iG) = ([\Sigma_f],[f], \chi_{(2)}(S))$$
 is a well defined semigroup isomorphism (hence $\Ii_2(M)$ is eventually a group).
 \end{theorem}
 
 To a large extent the proof is an adaptation of the above one when $M$ is orientable, 
 but one has to face  several differences (the existence of knots in $M$ with
 solid Klein bottle tubular neighbourhoods, the absence of framing of $M$ and so on).
 The basic reason for the final simpler form of $\Ii_2(M)$ is that the subgroup 
 of the immersed surfaces in a $3$-ball of $M$  is a quotient of 
 $\Ii_2(\R^3) \sim \Z/8\Z$ isomorphic to $\Z/2\Z$. For $\Ii_2(\R^3)$ is generated
 by the Boy surface $B$ and as $M$ is non orientable there is a diffeotopy of $M$
 which sends a $3$-disk of $M$ containing a copy of $B$ into itself reversing
 the orientation; hence $[B]_\iG = [\bar B]_\iG = -[B]_\iG$.   
 
 \section{Tear and smooth-rational equivalences}\label{tear}
The notion of blowing up a manifold along a smooth centre has been
defined in Section \ref {BU}. In Section \ref {stable-2-Nash} we have
interpreted the stable equivalence between surfaces in terms of
blowing up of points which are the only possible smooth centres in
such a case. If $M$ is now a compact boundaryless $3$-manifold besides
the points we have also any link of knots in $M$ as a possible smooth
centre.  In this section, referring to \cite{BM}, we widely study some
equivalence relations generated by blowing up $3$-manifolds along
smooth centres (and diffeomorphisms). We will discuss also applications
of this study to the so called $3$-{\it dimensional Nash's rationality
  conjecture}.
\medskip

\subsection{$3$-dimensional blowing-up-or-down}\label{BUD}
We denote by $\Mm_3$ the class of all compact
connected boundaryless $3$-manifolds.  Let $M$ be such a manifold. A
possible smooth centre $X$ of a blow up
$$\pi: B(M,X)\to M$$
is either a finite set of points
or a link of a finite number of pairwise disjoint knots in $M$,
$L=K_1\cup \dots \cup K_s$.
We know that $D_X:= \pi^{-1}(X)$ is a hypersurface of $B(M,X)$ 
called the {\it exceptional hypersurface}. We also say that
$M$ is obtained by {\it blowing down} $\tilde M := B(M,X)$
along the hypersurface $D_X$. 

For simplicity let us analyse connected
centres. A connected smooth centre in $M$ is either a point or a knot $K$.
We know that the effect of blowing up one point consists (up to
diffeomorphism) in performing a connected sum $M\cs \PP^3(\R)$, the
exceptional hypersurface being a one-side projective plane $\PP^2(\R)$
that is a projective plane with oriented tubular neighbourhood.
\smallskip

As $M$ is not necessarily orientable then a knot $K$
either preserves the orientation, that is it has a solid torus tubular
neighbourhood in $M$, or it reverses the orientation, that is it has a
solid Klein bottle tubular neighbourhood in $M$.  In the first case
the exceptional hypersurface $D_K$ in $B(M,K)$ is a one-side torus. In
the second it is a one-side Klein-bottle. Reciprocally we have
\begin{proposition}\label{down} Let $\tilde M$ be in $\Mm_3$
and $Y$ be a hypersurface of $\tilde M$ which is either a projective
plane with oriented tubular neighbourhood, a one-side torus or
a one-side Klein bottle. Then there exists $M$ in $\Mm_3$ and 
a smooth centre $X\subset M$ such that $\tilde M = B(M,X)$
and $Y=D_X$.
\end{proposition}
\Dim If $Y\sim \PP^2(\R)$ with orientable tubular neighbourhood
$N(K)$, then $N(K)\sim \PP^3(\R)\setminus {\rm Int}(B)$
where $B$ is a $3$-ball. Then $ \tilde M= M\cs \PP^3(\R)$
for some $M$ so that $\tilde M$ is the blow up of $M$ at a point.

The standard model of a tubular neighbourhood of a one-side
torus is obtained by taking the blow up
$$\pi: N:= B(D^2\times S^1, \{0\}\times S^1)\to D^2\times S^1 \ . $$
Denote by $p:D^2\times S^1 \to S^1$ the natural projection, $D^2_x=p^{-1}(x)$.
$N$ is diffeomorphic to $\Mm \times S^1$, $\Mm$ being a M\"obius strip,
with natural projection $\tilde p: \Mm\times S^1 \to S^1$ such that $\tilde p = p\circ \pi$;
for every $x\in S^1$, $\Mm_x=\tilde p^{-1}(x)= B(D^2_x, \{0\}\times \{x\})$.
On the torus $\partial N\sim \partial D^2\times S^1$ it is defined the involution $\tau$
which restricts to the antipodal one on every $\partial D^2_x$.
$N$ (and coherently every $\Mm_x$) can be identified with the mapping cylinder of $\tau$.
The exceptional hypersurface is the torus $D= s_0 \times S^1$, where 
$s_0= \pi^{-1}(\{0\} \times \{x_0\})$ and
$x_0$ is a base point on $S^1$. The mapping cylinder structure realizes also $N$ as being
a tubular neighbourhood of $D$, endowed with its projection $q:N\to D$. The restriction
of $q$ to $\partial N$ is a fibred double covering of $D$.

If $Y\subset \tilde M$ is a one-side torus, there are in fact
{\it several ways} to fix a parametrization
$$\phi: (N, D)\to (N(Y), Y) $$ so that the blow down $\pi: N \to
D^2\times S^1$ gives rise to a blow down $\pi: \tilde M \to M$, for
some $M$ in $\Mm_3$, where $(N(Y), Y)$ is mapped onto $(N(K), K)$, $K$
is a knot in $M$ which preserves the orientation and $N(K)$ is a
tubular neighbourhood of $K$ in $M$.  To do it assume that $N(Y)$ is
constructed by using a normal line bundle $\xi$ on $Y$ in $\tilde
M$. By hypothesis, the Euler class $\omega^1(\xi)\neq 0$.  Fix any
fibration $\Ff_s$ of $Y$ by smooth circles parallel to a circle $s$
such that $\omega^2(\xi)\sqcup [s]\neq 0$.  This means that the
restriction of the line bundle $\xi$ to $s$ is not trivial. Then there
is a diffeomorphism $\phi: (N,D)\to (N(Y),Y)$ such that the fibration
$\Ff_{s_0}$ of $D$ by the circles parallel to $s_0$ is mapped to the
fibration $\Ff_s$.  To see it we can transfer the question to the
above standard model.  The fibration $\Ff_{s_0}$ of $D$ lifts by the
projection $q$ to the fibration by meridians of $\partial N \sim
\partial D^2 \times S^1$; set $m_0= \partial D^2 \times \{x_0\}$ and
denote by $\tilde \Ff_{m_0}$ such fibration.  Fix on $D$ another
fibration $\Ff_s$ parallel to an $s$ with the properties fixed above.
This lifts by the projection $q$ to a fibration $\tilde \Ff_{\tilde
  s}$ of $\partial N$ by circles parallel to a $\tilde s$ such that
$[\tilde s]=[m_0]\in \eta_1(\partial N)$. Moreover, by construction
$\tilde \Ff_{\tilde s}$ is invariant by the involution $\tau$.  We
claim that, possibly up to isotopy of $s$, there is a diffeomorphism
$h$ of the torus $\partial N$ which sends $\tilde \Ff_{m_0}$ to
$\tilde \Ff_{\tilde s}$ and extends to a diffeomorphism of $(N, D)$
sending the fibration$\Ff_{s_0}$ of $D$ to $\Ff_s$. In such a case it
is easy to see that the topological space obtained by collapsing every
fibre of $\Ff_s$ to one point results from another blow down of
$(N,D)$ obtained by the {\it flip} $\Ff_{s_0} \to \Ff_s$ of fibrations
of the exceptional hypersurface $D$.  To justify the claim, let us
identify $\partial N$ with $\R^2/\Z^2$, endowed with ``linear''
cordinates such that the line $\{y=0\}$ is mapped onto $l_0=
\{p_0\}\times S^1$, while the line $\{x=0\}$ is mapped onto $m_0$ and
the involution can be expressed as $\tau(x,y)=(x,y+1/2)$; up to
isotopy a generic diffeomorphism in the form $h(x,y)=(ax+by, cx+dy)$,
with the coefficients belonging to a matrix in $GL(2,\Z)$.  Under our
hypotheses, $h(0,y)= (by, dy)$ where $b$ is even and $d$ is odd, so
that clearly $h\circ \tau = \tau \circ h$ and this is enough to
conclude.

The discussion for the one-side Klein bottle is similar (however, see
Remark \ref{Kbottle}).

\cvd

\subsection{Tears and Dehn surgery}\label{tear-DS}
The possibility to {\it flip the fibrations of an exceptional
  hypersurface} hence to modify the corresponding blowing down
(sometimes this modification is called a {\it flop}), suggests a way
to possibly modify the topology of $3$-manifolds.

\begin{definition}\label{tear-def}
  {\rm Let $M$ be in $\Mm_3$ and $L=K_1\cup \dots \cup K_s$ be a link
    in $M$ whose constituent knots preserve the orientation. We say
    that $M'$ in $\Mm_3$ is obtained from $M$ by a {\it tear along
      $L$}, if up to diffeomorphism there is blow down flop
$$ M \leftarrow B(M,L) \rightarrow M'$$ associated to a system of
    flips of fibrations of the exceptional hypersurfaces $D_{K_i}$ as
    in the proof of Proposition \ref{down}.  In other words $(B(M,L),
    D_L)=(B(M',L'),D_{L'})$ for some link $L'=K'_1\cup \dots \cup
    K'_s$ in $M'$ whose constituent knots preserve the orientations.
  }
\end{definition}

\begin{lemma}\label{tear-equiv}
  Tears define an equivalence relation called {\rm tear equivalence}
  and we write $M\sim_t M'$.
\end{lemma}
\Dim If we move a centre by an ambient isotopy, the result of a
blowing up does not change up to diffeomorphism preserving the
exceptional hypersurfaces. Given a tear from $M$ to $M'$ (with
associated links $L$ in $M$ and $L'_1$ in $M'$) and a tear from $M'$
to $M"$ (with associated links $L'_2$ in $M'$ and $L"$ in $M"$), by
transversality we can assume that $L'_1\cap \L'_2 =\emptyset$, hence
there is a copy of $L'_2$ in $M$, and a copy of $L'_1$ in $M"$ so that
$L\cup L'_2$ and $L"\cup L'_1$ are links in $M$ and $M"$ respectively,
supporting a tear from $M$ to $M"$. This proves that the relation is
transitive. It is trivially riflessive and symmetric.

\cvd  

\smallskip

\begin{remark}\label{Kbottle}{\rm A priori one would consider also
tear along knots which reverse the orientation. However,
for such a tear $M\leftarrow \tilde M \rightarrow M'$, it turns
out that $M\sim M'$; this happens because on a Klein bottle there is 
only one isotopy class of smooth circles with annular tubular
neighbourhood. So we consider only tears along knots preserving
the orientation.}
\end{remark}
\smallskip

It is convenient to rephrase tears in terms of more usual
modifications performed on $3$-manifolds. As above, let $M$ be in
$\Mm_3$, $L=K_1\cup \dots \cup K_s$ be a link in $M$ with constituent
knots preserving the orientation. Let $N(L)=N(K_1)\amalg \dots \amalg
N(K_s)$ be a tubular neigbourhood of $L$ in $M$. Consider the manifold
with $s$ toric boundary components
$$N := M \setminus {\rm Int} N(L) \ . $$
We say that $M'$ is obtained
by a {\it Dehn surgery} of $M$ along $L$ if, up to diffeomorphism, it is obtained 
by gluing back every $N(K_i)$ to $N$ along the torus 
$\partial N(K_i)$ by means of a
diffeomorphism $h_i: \partial N(K_i)\to \partial N(K_i)$, $i=1,\dots, s$. 
$L\subset N(L)$ determines a link $L'=K'_1\cup \dots \cup K'_s$ in $M'$ and the
identity map of $N$ extends to a diffeomorphism $\psi: M\setminus L
\to M'\setminus L'$.  If $m_i$ is a meridian of $\partial N(K_i)$,
then $h_i(m_i)=s_i$ is a smooth circle on $\partial N(K_i)$.  The
fibration of $\partial N(K_i)$ by meridians parallel to $m_i$ is
mapped by $h_i$ to a fibration by circles parallel to $s_i$. These are
meridians of a tubular neighbourhood of $L'$ in $M'$.  If every $s_i$
is a {\it longitude} of $\partial N(K_i)$ then $M'$ is obtained from
$M$ by an ordinary surgery already considered above. So Dehn
surgery generalizes the ordinary surgery associated to $4$-dimensional triads. 
The
diffeomorphism $\phi$ extends to a diffeomorphism $\phi:M \to M'$ if
and only if every $s_i$ is a meridian of $\partial N(K_i)$.

Now, up to diffeomorphism, $B(M',L')$ is obtained from $B(M,L)$
by gluing back every $B(N(K_i),K_i)$ to $N$  along the torus 
$\partial N(K_i)$ by means of the same
diffeomorphism $h_i: \partial N(K_i)\to \partial N(K_i)$, $i=1,\dots, s$, as before. 

\begin{definition}\label{D-to-T}{\rm We say that a Dehn surgery {\it lifts to a tear} 
if the diffeomorphism $\tilde \phi: B(M,L)\setminus D_L \to B(M',L')\setminus D_{L'}$
which lifts $\phi: M\setminus L \to M'\setminus L'$, extends to a diffeomorphism
$\tilde \phi: B(M,L)\to B(M',L')$, preserving the exceptional hypersurface.}
\end{definition}

We have

\begin{proposition}\label{Dehn-tear}
  A Dehn surgery from $M$ to $M'$ lifts to a tear if and only if for
  every $i=1,\dots, s$, $[s_i]=[m_i]\in \eta(\partial N(K_i))=
  \Hh^1(\partial N(K_i);\Z/2\Z)$.
\end{proposition}
\Dim The condition is necessary because the meridians  generates the kernel
of the unoriented bordism morphism induced by the inclusions $\partial N(K_i)\to N(K_i)$.
The other implication rephrases the proof of Proposition \ref{down}.

\cvd 
 
\smallskip

With respect to ordinary surgery we have the following immediate
corollary.
\begin{corollary}\label{evenframe}
  Let $M'$, $M"$ be obtained by ordinary (longitudinal) surgery on $M$
  along a same link $L=\cup_i K_i$ with different normal framings
  $\{\fG'_i\}$ and $\{\fG"_i\}$ respectively.  Let $L'\subset M'$ and
  $L"\subset M"$ be the links corresponding to $L$ respectively.  Then
  $M"$ is obtained (up to diffeomorphism) from $M'$ by a tear of the
  form
$$M'\leftarrow B(M',L')=B(M",L")\rightarrow M"$$ if and only if every
  $\fG'_i$ differs from $\fG"_i$ by an even number of twists.
\end{corollary}

\cvd

\smallskip

Hence tear equivalence can be considered as a specialization of the
equivalence relation generated by Dehn surgery. As this last extends
ordinary surgery and preserves orientability, then we already know
that {\it being or not orientable is a complete invariant for Dehn
  surgery equivalence}. We are going to see that this is no longer
true for tear equivalence. We refine the `orientable/non-orientable' partition $\Mm_3 = \Mm^+_3
\amalg \Mm^-_3$ so that we eventually have three {\it types}, completely
determined by the behaviour of $\omega^1(*)$:
\smallskip

-  $\omega^1(M)=0\in \Hh^1(M;\Z/2\Z)$, that is  it is {\it orientable};

- $\omega^1(M)\neq 0$ and $\omega^1(M)^2:= \omega^1(M)\sqcup \omega^1(M)=0$,
then we say that $M$ is {\it weakly non orientable}, that is $M\in \Mm^w_3$.

-   $\omega^1(M)\neq 0$ and $\omega^1(M)^2:= \omega^1(M)\sqcup \omega^1(M)\neq 0$,
then we say that $M$ is {\it strongly non orientable}, that is $M\in \Mm^s_3$.

\medskip

{\bf Characteristic surfaces:} If $M$ is non orientable, every
hypersurface $F$ which represents $\omega^1(M)$ is called a {\it
  characteristic surface} of $M$. We can assume that $F$ is connected
and it is necessarily orientable: the boundary $\partial N(F)$ of a
tubular neighbourhood is connected and orientable as it is the
boundary of the orientable manifold $M\setminus {\rm Int} N(F)$; the
projection of $\partial N(F)$ to $F$ is $2:1$ and every orientation on
$\partial N(F)$ descends to $F$.

\smallskip

We have
\begin{proposition}\label{type-invariant} 
  Let $M\sim_t M'$ be realized by a tear
  $$M \xleftarrow{\pi} B(M,L)=\tilde M = B(M',L') \xrightarrow{\pi'} M', \ L=\cup_{i=1}^s K_i \ . $$

1) For every $j=0, \dots, 3$, $\pi^*: \Hh^j(M;\Z/2\Z)\to \Hh^j(\tilde
M;\Z/2\Z)$ is an injective homomorphism and the similar fact holds for
$\pi'$.

\smallskip

2) $\Hh^1(\tilde M;\Z/2\Z)\sim \Hh^1(M;\Z/2\Z)\oplus (\Z/2\Z)^s$ where
the last factor is generated by the components $D_{K_i}$ of $D_L$;
$\Hh^2(\tilde M;\Z/2\Z)\sim \Hh^2(M;\Z/2\Z)\oplus (\Z/2\Z)^s$ where
the last factor is generated by the fibres $\Mm_i$ of the fibrations
$\Mm\times K_i \to K_i$ of $D_{K_i}$ ; similarly for $\pi'$.
\smallskip

3) For every $j=0, \dots, 3$, there is a natural isomorphism 
$$h_j:  \Hh^j(M;\Z/2\Z)\to \Hh^j(M';\Z/2\Z)$$
such that $(\pi')^*\circ h^j = \pi^*$. Moreover $h_1(\omega^1(M))=\omega^1(M')$
and for every
$\alpha \in \Hh^1(M;\Z/2\Z)$, $h^2(\alpha \sqcup \omega^1(M))= h^1(\alpha)\sqcup \omega^1(M')$.
 
\smallskip

4) $M$, $M'$ are of the same type.
\end{proposition}

\Dim Let us justify $(1) - (3)$. 
For every $j$, every class in  $\Hh^j(M;\Z/2\Z)$ can be represented by an embedded proper $(3-j)$-submanifold $S$
transverse to the link $L$. The corresponding class in $\Hh^j(\tilde M;\Z/2\Z)$ 
is represented by the strict transform $\tilde S$ of $S$ via the blow up. If $j=2, 3$ in fact $\tilde S$
is mapped diffeomorphically onto $S$ by $\pi$. If $j= 0$, $\tilde S = \tilde M$. If $j=1$, then
$\tilde S = B(S, S\pitchfork L)$. As for $(2)$ notice that $\Mm_i \bullet D_{K_j} = \delta_{i,j}$.
As for $(3)$ consider the diffeomorphism  
$$\phi: M\setminus L \to M' \setminus L' \ , $$
If $j=2,3$, then $h_j$ is determined by the diffeomorphism $S \sim \phi(S)$.
If $j=0$, then $h_0([M])=[M']$, and notice that  $[\tilde M,\pi]=[M]$, $[\tilde M,\pi']=[M']$. 
If $j=1$, then $S$ is a hypersurface transverse to $L$. Then 
$S\setminus {\rm Int} \ N(L)$ is sent diffeomorphically onto $\bar S'$ properly embedded into 
$M'\setminus {\rm Int} N(L')$; as $\phi$ preserves the class of meridians mod $(2)$, then $\bar S'$
can be completed to a boundaryless hypersurface $S'$ transverse to $L'$.
This geometric correspondence $S\leftrightarrow S'$ induces $h_1$.
If $S$ is a characteristic surface of $M$, as the constituent knots of $L$ preserve
the orientation, we can assume that $S\cap L = \emptyset$, so that the diffeomorphic
surface $S'=\phi(S)$ does not intersect $L'$ and is a characteristic surface of $M'$.
The last statements of $(3)$ follow. 
 Clearly $(4)$ is a corollary of the other items.

\cvd

\smallskip

In what follows we will say that $S'$ obtained from $S$ as in the above proof
is obtained by {\it darning} $S$ (with respect to the given tear). 

\begin{remark}\label{no-ring}
  {\rm One would wonder about a graded ring isomorphism in above
    statement $(3)$.  But this is not true. For example $S^1\times
    S^2$ and $\PP^3(\R)$ can be obtained by ordinary surgery olong an
    unknot $K\subset \R^3 \subset S^3$ with the standard even normal
    framing $\fG_0$ and the framing which differs from it by two
    twists, respectively. By Corollary \ref{evenframe}, they are
    connected by a tear, but their $\Z/2\Z$-cobordism rings are
    different.}
\end{remark}

\subsection{$rs$-equivalence}\label{rs}
We define now a coarser equivalence relation generated by
blowing-up-or-down.

\smallskip

\begin{definition}\label{smooth-rat}{\rm
    Let $M$, $M'$ be in $\Mm_3$.  We say that,
    {\it up to diffeomorphism, $M'$ is obtained from $M$ by a finite chain of
      blowing-up-or-down} if there is a finite chain
    of the form:
$$M\rightarrow M_0 \leftrightarrow M_1 \leftrightarrow M_2 \leftrightarrow \dots \leftrightarrow M_n \leftarrow M'$$
where:
\begin{enumerate}
\item Every $M_i$ is in $\Mm_3$;
\item the right and left arrows are diffeomorphisms; 
\item for every $i\neq n$, $M_i \leftrightarrow M_{i+1}$ either is a 
  a blow up along a smooth centre
  $$M_i \leftarrow M_{i+1}=B(M_i,C_i)$$ or a blow-up 
  $$M_i=B(M_{i+1},Z_{i+1}) \rightarrow M_{i+1}$$
  so that  $M_{i+1}$ is obtained by a {\it blow down} of $M_i$.

\end{enumerate}
}
\end{definition}
\smallskip

This defines another equivalence relation called {\it
  smooth-rational equivalence} which extends the diffeomorphism one and
  also the tear equivalence. We write $M \sim_{sr} M'$. Note that noone
  of the tear invariants pointed out in Proposition \ref {type-invariant} 
  persists for the $sr$-equivalence.

  \smallskip

Our goals are to fully determine the quotient set of $\Mm_3$ 
mod $\sim_{sr}$ or mod $\sim_t$. Tear equivalence preserves the type
so we can split the study of $\Mm_3$ mod $\sim_t$ type by type.
\smallskip

The results for  $\Mm^+_3$ mod $\sim_t$ and for $\Mm_3$ mod $\sim_{sr}$ are easy to state: 

\begin{theorem}\label{t+}
  For every $M$, $M'$ in $\Mm^+_3$, then $M\sim_t M'$ if and only if
  $\dim \Hh^1(M;\Z/2\Z)=\dim \Hh^1(M';\Z/2\Z)$. If $\dim
  \Hh^1(M;\Z/2\Z) = h$, then
$$M \sim_t S^3\cs h\PP^3(\R) \ . $$
\end{theorem}
\smallskip

\begin{proposition}\label{-sr+}
  For every $M$ in $\Mm^-_3$ there exists $M'\in \Mm^+_3$ such that $M
  \sim_{sr} M'$.
\end{proposition}

As a corollary we have

\begin{theorem}\label{sr} For every $M$ in $\Mm_3$, $M\sim_{sr} S^3$.
\end{theorem}
\Dim  By assuming Theorem \ref{t+} and Proposition \ref{-sr+}. If $M$ is in $\Mm^+_3$, then
the result follows immediately from Theorem \ref{t+} as $S^3\cs h\PP^3(\R)$
is obtained by blowing-up $S^3$ at $h$ points. If $M\in \Mm^-_3$, Proposition
\ref{-sr+} reduces it to the orientable case.

\cvd 

\smallskip

The structure of $\Mm^-_3$ mod $\sim_t$ is intrinsecally more
complicated, we will face it later.

\subsection{Disorientated surfaces and weakly trivial knots}\label{disorientated}
Let $N$ be a compact $3$-manifold with possibly non empty boundary $\partial N$.
A connected properly embedded surface $F$ in $N$ is said {\it disorientated}
if it is non orientable and has an orientable neighbourhood in $N$.
\smallskip

Let $M$ be in $\Mm_3$, and $K\subset M$ be a knot which preserves
the orientation with a tubular neighbourhood $N(K)$. 
Then $K$ is said {\it weakly trivial} if there exists a longitude $l$ on 
$\partial N(K)$ which bounds a disorientated surface $F$ properly
embedded into $M\setminus {\rm Int} N(K)$.
\smallskip

The notion of tear makes sense also for a manifold with boundary $N$,
provided that the supporting link is contained in the interior of $N$.
The following proposition shows tear's power to simplify 
disorientated hypersurfaces and eventually the topology of $3$-manifolds. 

\begin{proposition}\label{tear-unknotting} Let $S\subset N$ be a disorientated
hypersurface. Assume that $S$ has at most two boundary components. 
Then there are: a link $L\subset {\rm Int} (S) \subset {\rm Int}(N)$
with constituent knots preserving the orientation, a tear
$$ N\leftarrow B(N,L)=\tilde N = B(N',L')\rightarrow N'$$
and a surface $S'\subset N'$ obtaining by {\rm darning}
$S$ (with respect to the tear) such that:
\begin{enumerate}
\item If $S$ is boundaryless then  $S'$ is a disorientated projective plane.
\item If $\partial S$ is connected then $S'$ is a disk properly
embedded in $N'$ 
\item If $\partial S$ has two components, then $S'$ is a two-sides annulus 
properly embedded in $N'$.
\end{enumerate}
\end{proposition}
\Dim $S$ is diffeomorphic to the connected sum of $s$ copies of
$\PP^2(\R)$, $s\geq 1$, from which we have removed $k$ disjoint open
$2$-disks, either $k=0,1,2$. Let $L= K_1\cup \dots \cup K_s$ be formed by
the cores of $s$ pairwise disjoint M\"obius strips $\Mm_i$ embedded in
$S$.  Each $K_j$ reverses the orientation of $S$ and preserves the
orientation of $N$ (because $S$ has an orientable neighbourhood). Then
$[\partial \Mm_i]$ is a meridian of $\partial N(K_i)$ mod $(2)$ and we can
consider the corresponding tear $ N\leftarrow B(N,L)=B(N',L')\to
N'$. Then every $K_i$ collapses to one point in a dearning surface
$S'$ properly embedded in $N'$ with orientable neigbourhood. If $k=0$ then
$S'$ is a $2$-sphere; in order to get a disorientated $\PP^2(\R)$ it
is enough to remove from $L$ one constituent knot. In the other two
cases we get either a disk or an annulus.

\cvd

\smallskip

\begin{corollary}\label{chain} For every $M\in \Mm_3$ there is a chain of the form
$$M\rightarrow M_0 \leftrightarrow M_1 \leftrightarrow \dots \leftrightarrow M_n \leftarrow M'$$
such that:
\begin{enumerate}
\item Every $M_i$ is in $\Mm_3$, the right and left arrows being diffeomorphisms;
\item $\Hh^1(M';\Z/2\Z)$ is generated by $\omega^1(M')$;
\item For every $i\neq n$, $M_i \leftrightarrow M_{i+1}$ either is:
\smallskip

- a tear;
\smallskip

- a blow up $M_i=B(M_{i+1},x_0) \to M_{i+1}$ at a point of $M_{i+1}$;
\smallskip

- a blow up $M_i=B(M_{i+1}, K) \to M_{i+1}$ along a smooth knot of $M_{i+1}$ which preserves
the orientation.
\smallskip

\end{enumerate}
\end{corollary}
\Dim If $M$ already verifies $(2)$, then take $M'=M$. Otherwise there
is a hypersurface $S$, such that $[S]\neq 0 \in \Hh^1(M;\Z/2\Z)$ and
is not a characteristic surface of $M$. We can assume that $S$ is
connected and that there is a characteristic surface $F$ such that
either
\smallskip

- $S\cap F = \emptyset$, that is $S\subset M\setminus N(F) $ for a
small tubular neighbourhood of $F$;
\smallskip

- $S\pitchfork F$ along a knot $K\subset S$ which does not divide it.
\smallskip

- In both cases $S\setminus {\rm Int} N(F)$ is properly embedded into
$M\setminus {\rm Int} N(F)$, has oriented neighbourhood therein, and
there is a smooth circle $C\subset M\setminus {\rm Int} N(F)$ with non
trivial intersection number mod $(2)$ with $S\setminus {\rm Int}
N(F)$.  By adding an embedded $1$-handle along a suitable arc of $C$,
we can also assume that $S\setminus {\rm Int} N(K)$ is
disorientated. Now, if $S$ is disjoint from $F$, by Proposition
\ref{tear-unknotting} there is a tear which converts $S$ into a
disorientated projective plane; this can be considered as the
exceptional hypersurface of a blow up of a point. In the other case
there is a tear converting $S\setminus {\rm Int} N(F)$ into an
annulus; together with $S\cap N(F)$ they form a one-side torus which
can be considered as the exceptional hypersurface of a blow up along a
knot.

\cvd

\smallskip

\begin{corollary}\label{+reduction} If $M$ is orientable and $\dim \Hh^1(M;\Z/2\Z)=h$,
then $$M\sim_t  \tilde M$$ where
$$ \tilde M = h\PP^3(\R)\cs M'$$
and $\Hh^1(M';\Z/2\Z)=0$.
\end{corollary}
\Dim As $M$ is orientable $\omega^1(M)=0$; hence the statement and the proof
of Corollary \ref {chain} tell us that $\Hh^1(M';\Z/2\Z)=0$ and that only blow up
of points does occur. As up to isotopy a point misses any possible already
present exceptional hypersurface, tears and blowing up of points commute
and the corollary follows.

\cvd

\smallskip

\begin{corollary}\label{chirurgia} Let $M$ and $M'$ in $\Mm_3$ be such that
$$\Hh^1(M;\Z/2\Z)=\Hh^1(M';\Z/2\Z)=0 \ . $$ Assume that $M'$ is obtained from 
$M$ by an ordinary (longitudinal) surgery of $M$ along a weakly trivial knot
$K\subset M$. Then $M\sim_t M'$.
\end{corollary}
\Dim By Proposition  \ref{tear-unknotting} there is a tear from $M$ to $M_1$
converting $K$ to a genuine trivial knot $K_1 \subset M_1$.
So up to tear equivalence, we can assume that $M'$ is obtained from  $M$
by an ordinary surgery along a trivial knot $K$. As they have both vanishing
$\Hh^1$ the normal framing $\fG$ of this surgery must be odd with respect
to the framing $\fG_0$ determined by a collar of $K$ in a spanning $2$-disk.
On the other hand $M$ is diffeomorphic to the manifold obtained by using
the framing $\fG_1$ which differs from $\fG_0$ by one twist. Hence
by Corollary \ref {evenframe}, there is a tear from $M$ to $M'$.

\cvd

\smallskip

As a further corollary we can prove Proposition \ref{-sr+}, which we state again 
\smallskip

{\it For every $M$ in $\Mm^-_3$ there exists $M'\in \Mm^+_3$ such that $M
  \sim_{sr} M'$.}
\smallskip

{\Dim} Assume that $M$ has a connected characteristic surface $F$ of genus $g+1>1$.
We are going to show that $M\sim_{sr} M'$ such that $M'$ either has a characteristic
surface $F'$ of genus $g$ if $g>0$ or it is orientable. Clearly this will achieve the
result by induction on $g$. First we can assume that $F$ is one-side in $M$.
In fact let $K\subset F$ be a smooth circle which does not divide $F$. Then
the strict transform $\tilde F$ of $F$ in $B(M,K)$ is a one-side characteristic surface
of the same genus. If $F$ is a one-side torus then it is the exceptional hypersurface
of a blow down to an orientable $M'$ and we have done. If $g>1$,
there is a smooth circle $C$ on $F$ which divides it by a one-side torus $T_0$
with one hole, and a bilateral surface $S_0$ of genus $g-1$ with one hole.
By adding an embedded $1$-handle as in the proof of  Corollary \ref{chain},
we can modify $S_0$ far from $C$ and make it desorientated.
Then by Proposition \ref {tear-unknotting} there is a tear from $M$ to say $M_1$
which convert $S_0$ to a $2$-disk so that $C$ becomes a trivial knot in $M_1$.
The manifold $M_2$ obtained by ordinary surgery along $C$ with normal framing
given by a tubular neighbourhood of $C$ in $F$ is tear equivalent to $M_1\cs \PP^3(\R)$,
hence it is $sr$-equivalent to $M_1$ hence to $M$. We conclude by noticing that a
characteristic surface of $M_2$ is given by the disjoint union of a surface of genus $g$
and a one-side torus which again can be considered as the exceptional hypersurface
of a blow down.

\cvd

\smallskip

\subsection{ $\Mm_3^+$ mod $\sim_t$ and $\Mm_3$ mod $\sim_{sr}$}\label{t+sr}
We are ready to prove Theorems \ref {t+} and \ref{sr}. Thanks to Corollary
\ref {chirurgia} and Proposition \ref{-sr+}, it will enough to prove the following

\begin{lemma}\label{chainH} 
 For every $M$ in $\Mm_3$ such that $\Hh^1(M;\Z/2\Z)=0$, there
 exists a sequence $S^3=M_0, M_1, \dots, M_n \sim M$, such that
 \begin{enumerate}
 \item For every $M_i$, $\Hh^1(M_i;\Z/2\Z)=0$;
 \item $M_{i+1}$ is obtained from $M_i$ by an ordinary surgery
 along a weakly trivial knot $K_{i+1}\subset M_i$.
 \end{enumerate}
 \end{lemma}
 \Dim We use some notions that we will develop in Chapter \ref{TD-4}, Section \ref{with-bound}. 
 Here
 we outline the main points. We know that $S^3\sim_\sigma M$, that
 is there is a triad $(W,S^3, M)$ with a handle decomposition made by $2$-handles
 only, so that $M$ is obtained by longitudinal surgery along a framed link
 $L=\cup_i K_i$ in $S^3$. The framing $\fG_i$ is encoded by an integer 
 which express the number of twists with respect to the framing given by
 the collar of $K_i$ in a Seifert surface. 
 The intersection form of $\Hh^2(W;\Z/2\Z)$ is represented by the linking
 matrix mod $(2)$ of this framed link $L$, so that along the diagonal
 we have the reduction mod $(2)$ of the above integers.
 As $\Hh^1(M;\Z/2\Z)=0$ then the intersection form is non degenerate.
 Possibly performing an elementary blow-up move (Section \ref {kirby-cal}),
 we can also assume that  the form is not totally
 isotropic, hence it has an orthogonal basis (see Section \ref {algebraic-form}). 
 By realizing such a change of basis
 by handle sliding, we get that every $K_i$ is the boundary of a surface $S_i$
 disjoint from the rest of the link, and the new normal framings are odd.
 So the knot $K_{i+1}$ is weakly trivial in the manifold $M_i$ obtained
 by the surgery along the partial framed link $K_1\cup \dots \cup K_i$.
 
 \cvd
 
 \smallskip 

\subsection{ $\Mm^-_3$ mod $\sim_t$}\label{-mod-t}  
This is more demanding. We will give exhaustive statements.
For detailed proofs a curious reader is addressed to \cite{BM}.

We can manage type by type. For $\Mm^s_3$ the statement is simpler;
alike the orientable case, the necessary conditions of Proposition \ref {type-invariant}
are also sufficient.
 
\begin{theorem}\label{s-t}
Let $M$, $M'$ be strongly non orientable. Then $M\sim_t M'$ if and only if
for every $j=0, \dots, 3$, there is a natural isomorphism 
$$h_j:  \Hh^j(M;\Z/2\Z)\to \Hh^j(M';\Z/2\Z)$$
such that $h_1(\omega^1(M))=\omega^1(M')$
and for every $\alpha \in \Hh^1(M;\Z/2\Z)$, $h^2(\alpha \sqcup \omega^1(M))= h^1(\alpha)\sqcup \omega^1(M')$.
\end{theorem}

\cvd

\smallskip

For weakly non orientable manifolds another tear invariant comes up. 
\smallskip

We begin with a construction that makes sense for every orientable
compact boundaryless surface $S$ embedded into any $M$ in
$\Mm_3$. Consider the subspace of $\eta_1(S)$ formed by the
$1$-boundary in $M$, that is $$\Bb(S,M) = \ker i_*$$ where $i: S\to M$
is the inclusion.  Let $\alpha \in \Bb(S,M)$. Then $\alpha = [c]$ for
some smooth circle $c$ on $S$.  By hypothesis, $c$ bounds a {\it
  membrane} $\MG\subset M$: by definition $\MG$ is a compact surface
embedded in $M$, such that $c=\partial \MG$, and moreover $\MG$ is in
``general position" with respect to $S$; this means that $S\pitchfork
{\rm Int} (\MG)$ and $S\cap \MG$ is the union $c\cup d$ where $d$ is a
smooth curve properly embedded in $S$ (i.e. $\partial d = \cap
\partial \MG$). Tubular neighbourhoods of $d$, $N(d,S)$ and $N(d,\MG)$
in $S$ and $\MG$ respectively, coincide at $\partial d$ along a
tubular neighbourhood of $\partial d=d\cap c$ in $c$. Then along the
abstract double $D(d)=d_+ \cup d_-$ of $d$ we can define a band
$N(D(d))$ equal to $N(d,S)$ on $d_+$, equal to $N(d,\MG)$ on $d_-$
glued by the indentity on $\partial d_+ = \partial d_-$.  Then we can
define by the self-intersection of $D(d)$ in $N(D(d))$
$$\rho_\MG(c) = D(d)\bullet D(d) \in \Z/2\Z \ . $$
We can pose the question under which hypotheses this construction {\it
  well defines} a homomorphism
$$ \rho_S: \Bb(S,M)\to \Z/2\Z, \ \rho (\alpha) = \rho_\MG(c), \ \alpha=
[c] \ . $$ This is widely discussed in \cite{BM}. Here we are
interested to the application of this construction to a characteristic
surface $F$ of $M$ in $\Mm^-_3$. We have

\begin{proposition}\label{rolling} Let $F$ be a characteristic surface of the
  non orientable $3$-manifold $M$.  Then $\rho_F: \Bb(F,M)\to \Z/2\Z$
  is well defined if and only if $M$ is weakly non orientable ($M\in
  \Mm^w_3$).  In such a case $\rho_F$ is a quadratic enhancement of
  the restriction, say $\beta$, to $\Bb(F,M)$ of the intersection form
  on $\eta_1(F)$.
\end{proposition}

\cvd

A first point where the vanishing of $\omega^1(M)\sqcup \omega^1(M)$ is relevant is in showing that
the value of $\rho_\MG(c)$ does not depend on the choice of the membrane
$\MG$. In fact one verifies that:
\smallskip

(i)  $ \sigma \sqcup \sigma \sqcup  \omega^1(M)+ \sigma \sqcup \omega^1(M)\sqcup \omega^1(M)=0$
for every $\sigma \in \Hh^1(M;\Z/2\Z)$  if and only if $\omega^1(M)\sqcup \omega^1(M)=0$;
\smallskip

(ii)  given two membranes $\MG$ and $\MG'$ of $c$, 
$\tau=\MG'\cup \MG$ define a  cycle mod $(2)$ in $M$ and ones verifies that
$$ \rho_{\MG'} - \rho_{\MG} = [\tau] \sqcup [\tau] \sqcup \omega^1(M)+ [\tau] \sqcup \omega^1(M)\sqcup \omega^1(M) \ . $$
This is the first step to show that $\rho(c)$ only depends on the class $[c]\in \eta_1(F)$.
\medskip

Let $M\in \Mm^w_3$, $F$, $\rho_F$, $\beta$ be as in the above proposition.
In general $\beta$ is degenerate, that is its radical $\Bb(F,M)^\perp\neq \{0\}$.
Then there are two possibilities: 
\smallskip

- $\rho_F\neq 0$ on $\Bb(F,M)^\perp$. 

\smallskip

- $\rho_F=0$ on $\Bb(F,M)^\perp$.
Set  $\hat \Bb(F,M) = \Bb(F,M)/\Bb(F,M)^\perp$.
Then $\rho_F$ descends to a homomorphism 
$$\hat \rho_F: \hat \Bb(F,M)\to \Z/2\Z$$ which is a quadratic
enhancement of the non degenerate form $\hat \beta$ induced by
$\beta$; one can define its {\it Arf invariant} (see Section \ref
{quadratic})
$$\delta_F:=\delta(\hat \rho_F)\in \Z/2\Z \ . $$

\smallskip

So we can associate to $F$ the symbol 
$$\tau_F \in \{\emptyset\} \cup \Z/2\Z$$
where $\tau_F = \emptyset$ if $\rho_F\neq 0$ on $\Bb(F,M)^\perp$, 
$\tau_F = \delta_F$ otherwise.
We have
\begin{proposition} Let $F$ be a characteristic surface of $M$ in $\Mm^w_3$.
Then
$$\tau_M := \tau_F$$ is well defined, that is it does not depend on
the choice of $F$ such that $[F]= \omega^1(M)$.
\end{proposition}

\cvd

\smallskip

Hence we have refined the type of weakly non orientable manifolds accordingly with
the value of $\tau_M$.  Finally we can complete the classification up to tear equivalence.

\begin{theorem}\label{w-t}
Let $M$, $M'$ be weakly non orientable. Then $M\sim_t M'$ if and only if
for every $j=0, \dots, 3$, there is a natural isomorphism 
$$h_j:  \Hh^j(M;\Z/2\Z)\to \Hh^j(M';\Z/2\Z)$$
such that $h_1(\omega^1(M))=\omega^1(M')$,
for every $\alpha \in \Hh^1(M;\Z/2\Z)$, 
$h^2(\alpha \sqcup \omega^1(M))= h^1(\alpha)\sqcup \omega^1(M')$
and moreover, $\tau_M = \tau_{M'}$. 
\end{theorem}

\cvd

\smallskip

Let us give more information about the eventual result.
First one finds representatives $M$ of every non orientable tear class 
endowed with a characteristic surface $F$ with minimal 
boundary space $\Bb(F,M)$. 
For every non orientable $M$, consider the pairs $(M,F)$ where
$F$ is a connected characteristic surface.
For every non orientable tear equivalence class $\alpha$, set
$$ g(\alpha) = \min \{g(F); \ (M, F), M \in \alpha\} \ . $$
We have

\begin{proposition}\label{min-g} Let $(M,F)$ be such that 
$g(F)= g(\alpha)$, $\alpha=[M]_t$. Then the {\it boundary dimension}
$$d(\alpha):= \dim \Bb(F,M)$$
is well defined (type by type) and we have:
\smallskip

(1) If $M$ is strongly non orientable,  then $d(\alpha,s)=0$;
\smallskip

(2) If $M$ is weakly orientable and $\tau_\alpha= \emptyset$, then
$d(\alpha,w,\emptyset)=1$;  
\smallskip

(3) If $M$ is weakly orientable and $\tau_\alpha= 0$, then
$d(\alpha,w,0)=0$.
\smallskip

(4) If $M$ is weakly orientable and $\tau_\alpha= 1$, then
$d(\alpha,w,1)=2$.
\end{proposition}

\cvd

\smallskip

We have given normal representatives for every 
orientable tear class $\alpha$, that is  $h\PP^3(\R)$, $h= \dim \Hh^1(M;\Z/2\Z)$,
$\alpha = [M]_t$. By elaborating on the minimizing representatives of
Proposition \ref{min-g},  we get normal representatives also for the non orientable
classes.  For every non orientable $\alpha = [M]_t$, define type by type the integer 
$$ h(\alpha,s)= \dim \Hh^1(M; \Z/2\Z)-2g(\alpha)$$
$$ h(\alpha,w, \tau_\alpha) =  \dim \Hh^1(M; \Z/2\Z)-2g(\alpha)+d(\alpha,w,\tau_\alpha) \ . $$

\begin{proposition}\label{normal-non-or} For every non orientable tear equivalence class $\alpha$       
there are explicitely given manifolds $M(\alpha,s)$ or $M(\alpha,w,\tau_\alpha)$ such that
either $$\alpha = [ h(\alpha,s)\PP^3(\R) \cs M(\alpha,s)]_t$$ or  
$$\alpha = [ h(\alpha,w,\tau_\alpha)\PP^3(\R) \cs M(\alpha,w,\tau_\alpha)]_t \ . $$
\end{proposition}

\cvd

\smallskip

We have more information about these normal representatives. Let us say that
$M$ is {\it smooth-rational elementary}  if it is obtained by means of a tower 
of blowing up along smooth centres over the standard $3$-sphere $S^3$
$$S^3  \leftarrow M_1 \leftarrow M_2 \leftarrow \dots \leftarrow M_k = M \ . $$
Then we have

\begin{proposition}\label{sr-t-rep} With the exception of the weakly non orientable class $\alpha_0$ 
such that  $\dim \Hh^1(M;\Z/2\Z)=1$, $\alpha_0 = [M]_t$, and $\tau_{\alpha_0} = 1$,
the normal representative of every tear class $\alpha$ is smooth-rational elementary.
In the exceptional case, $\alpha_0$ cannot be represented by any smooth-rational
manifold, and for the normal representative say $M_{\alpha_0}$ there is a 
smooth-rational $\tilde M_{\alpha_0}$ and a blow up 
$\tilde M_{
alpha_0} = B(M_{\alpha_0}, x_0)\to M_{\alpha_0}$, where $x_0$ is a point.
\end{proposition}

   \subsection{On $3$-dimensional Nash's rationality conjecture}\label{3DNash}
By using the classification up to tear equivalence, in \cite{BM} one gives
an answer to the so called Nash's conjecture in three dimensions.
\medskip

Let us say that a non singular $3$-dimensional real algebraic set  $X$ is 
{\it rational elementary} if it is obtained by a tower of blow up along real algebraic non singular
centres over the standard sphere $S^3$.
\smallskip

First one proves that every tear equivalence class has an explicitely given rational model which is in fact
elementary with one exception. Referring to Proposition \ref{sr-t-rep}, and using variations of Nash-Tognoli theorem (see Section \ref {nash-tognoli})
we have:

\begin{proposition}\label{rat-t} With the exception of the weakly non orientable class $\alpha_0$ 
such that  $\dim \Hh^1(M;\Z/2\Z)=1$, $\alpha_0 = [M]_t$, and $\tau_{\alpha_0} = 1$,
the normal representative of every tear class $\alpha$ can be realized to be a rational elementary
real algebraic set $Y_\alpha$. In the exceptional case,  there is 


- a  rational algebraic set $Y_0$ with one singular point $y_0$,

- a homeomorphism $h_0: Y_0 \to M_{\alpha_0}$ which is a diffeomorphism
on $Y_0 \setminus \{y_0\}$,

- an ``algebraic resolution of
singularity'' $\psi: \hat Y_0 \to Y_0$, such that $\hat Y_0$
is rational elementary and 
$\psi: \hat Y_0 \setminus \psi^{-1}(y_0) \to Y_0 \setminus \{y_0\}$
is an algebraic isomorphism.
\end{proposition}

\cvd

\smallskip

Then we have:

\begin{theorem}\label{Nash3D}
(ii) For every tear equivalence class $\alpha \neq \alpha_0$, for every
$M\in \alpha$, there is  a tear from $M$ to $Y_\alpha$ of the form
$$M\xleftarrow{\sigma} Y_M  \xleftarrow{\pG} \tilde Y_M=B(Y_\alpha, L_M) \xrightarrow{\pi} Y_\alpha$$  
where: 

- $\tilde Y_M$ is rational elementary obtained by blowing up $Y_\alpha$
along a non singular real algebraic link $L_M\subset Y_\alpha$;
\smallskip

- $Y_M$ is rational with regular $1$-dimensional singular set  
${\rm Sing} (Y_M)=\pG(D_{L_M})$
consisting of a union of non singular circles; 
\smallskip

- The surjective algebraic map $\pG$ is a `resolution of singularity', that is 
$$\pG: \tilde Y_M\setminus D_{L_M} \to Y_M \setminus {\rm Sing} (Y_M)$$
is an algebraic isomorphism between regular Zariski open sets;
\smallskip

- $\sigma$ is a homeomorphism which restricts to a diffeomorphism 
on $Y_M \setminus {\rm Sing} (Y_M)$ and on ${\rm Sing} (Y_M)$;
\smallskip

-  $\sigma \circ \pG$ is a smooth blow down.
\smallskip

(iii) As for $M\in \alpha_0$ we have a similar realization of a tear
of the form
$$M\xleftarrow{\sigma} Y_M  \xleftarrow{\pG} \tilde Y_M=B(Y_0, L_M) \xrightarrow{\pi} Y_0 \xrightarrow{h_0} M_{\alpha_0}$$
where $L_M \subset R(Y_0)$, and eventually the rational model $Y_M$ of $M$ has a further isolated
singular point and admits an algebraic resolution of singularity by means of the rational elementary $B(\hat Y_0, \hat L_M)$,
$\hat L_M = \psi^{-1}(L_M)$.

\end{theorem}

\cvd   

\smallskip

So the  theorem shows that every $M$ in $\Mm_3$ has a {\it singular} rational algebraic model
$Y_M$ with mild controlled singular set which, nevertheless, cannot be avoided by the specific
blow-up-and-dow way the model has been constructed. 
The situation is very similar to what we have done in the case
of surfaces (Section \ref  {stable-2-Nash}). In the case of surfaces Comessati tells us that  for 
genus greater than $1$,  the presence of one singular point in a rational model of an orientable surface
is not only an accident of the construction, it is intrinsecally unavoidable. 
The same question has been faced for threefolds (see  \cite{Ko});
roughly summarizing, one realizes that also in  dimension $3$, orientable manifolds admitting
a non singular rational model are very special. On the other hand, we have the following interesting fact
(see \cite{Ko2}): 
\smallskip

{\it For every $\alpha$, for every $M\in \alpha$, there are {\rm non singular} 
rational models, provided that one deals with a category of ``abstract" algebraic-like varieties
(also called {\rm Moishezon varieties})
which are only locally but not globally isomorphic to ordinary algebraic sets in some $\R^n$}.
 
 \medskip
 
 In fact in this larger setting also the singular blow down $\pG: \tilde Y_M \rightarrow Y_M$
 can be realized as a the inverse of an algebraic blow up along a non singular centre.

\chapter{On $4$-manifolds}\label{TD-4}
In this chapter we will apply several results estasblished so far to
compact $4$-manifolds. Similarly to the attitude of Chapter \ref{TD-3}
with respect to the geometrization of $3$-manifolds, we stress that we
will develop a few classical differential/topological themes, in no
way (with the exception of a final informative and discorsive section) 
we will touch the study of $4$-manifolds by means of {\it gauge
  theory} that has dominated the study of $4$-manifolds in last
decades; for a more up to date treatment of $4$-manifolds theory one
can refer for example to \cite{Sc}. In particular we will determine
$\Omega_4$, present some instances of ``classification of simply
connected $4$-manifolds up to stabilization'', and Rohlin's theorem
about the signature mod $(16)$ of $4$-manifold intersection forms.
The intersection form will be indeed the principal player.
\smallskip

We will deal with {\it oriented} $4$-manifolds.  $M$ will denote a
compact, connected, oriented, boundaryless smooth $4$-manifold.  By
using the notations and the results of Sections
\ref{intersection-form}, \ref{complexLB} and \ref{seifert} we have
that the {\it intersection form}
$$\sqcup: \Hh^2(M;\Z)\times \Hh^2(M;\Z)\to \Z$$
equivalently
$$ \bullet: \Hh_2(M;\Z)\times \Hh_2(M;\Z)\to \Z$$
is symmetric  and
induces a $\Z$-linear isomorphism
$$ \hat \phi: \Hh^2(M;\Z)\to {\rm Hom}(\Hh_2(M;\Z),\Z) \ . $$ Then the
free $\Z$-module $\Hh^2(M;\Z)=\Hh_2(M;\Z)$ is of finite rank say $n$,
and the intersection form is {\it unimodular}: for any basis of
$\Hh^2(M;\Z)$ the representing matrix $A$ belongs to $GL(n,\Z)$
i.e. $|\det A| = 1$.  Every class $\alpha \in \Hh^2(M;\Z)$ can be
represented by an oriented $2$-dimensional proper submanifold $F$;
$\alpha = [F]=0$ if and only if $F$ is the boundary of an embedded
Seifert hypersurface. Clearly the isometry class of the intersection
form is an invariant up to orietation preserving diffeomorphism. We
are in a situation formally similar to the case of compact
boundaryless surfaces $S$ with respect to the intersection form on the
$\Z/2\Z$-vector space $\eta_1(S;\Z/2\Z)=\Hh_1(S;\Z/2\Z)$. In the case
of surfaces we have seen in Chapter \ref{TD-SURFACE} that this
intersection form contains all relevant information; moreover, there
is a perfect parallelism between the abstract algebraic theory of
symmetric $\Z/2\Z$-bilinear forms and its $2$-dimensional
differential/topological realization.  We would try to pursue this
analogy as far a possible, obtaining in fact only very partial
results.

\section{Symmetric unimodular $\Z$-bilinear forms}\label{SUF}
In analogy to Section \ref{algebraic-form}, we face here the question
of the classification of finite rank, symmetric, unimodular
$\Z$-bilinear forms up to isometry.  It turns out that this abstract
classification is complete only for the class of {\it indefinite
  forms}, while the {\it definite} case is a wide largely unknown
territory. This is a main difference with respect to the
$\Z/2\Z$-case.  For more information and detailed proofs we refer the
reader to \cite{MH}.

We consider free $\Z$-modules $V$ of finite rank, endowed with a
symmetric unimodular $\Z$-bilinear form $\rho$. This means that the
$\Z$-linear map
$$V \to {\rm Hom}(V,\Z), \ v\to f_v, \ f_v:V\to \Z,
\ f_v(w)=\rho(v,w) $$ is an isomorphism. Equivalently, the symmetric
matrix $A$ representing $\rho$ with respect to any basis of $V$
belongs to $GL(n,\Z)$, $n= {\rm rank} V$, that is $|\det A| = 1$.
Isometry is defined in the usual way.  Sometimes we will make the
abuse of confusing a form with its isometry class.  Given $(V,\rho)$
and $(V',\rho')$ we can define the {\it orthogonal direct sum}
$$(V,\rho)\perp (V', \rho')$$ that is the symmetric unimodular form
$\rho \perp \rho'$ on $V\oplus V'$ that restricts to $\rho$
(resp. $\rho'$) on $V$ ($V'$) and such that $V$ and $V'$ are
orthogonal to each other.
 \smallskip
 
 \subsection{Some invariants}\label{invariant}
 We point out some isometry invariants besides the rank.
 \smallskip
 
 {\bf (Signature)} By extension of the coefficients $\Z \subset \R$,
 $V$ becomes a lattice in a $\R$-vector space $V_\R$ so that $\dim
 V_\R = {\rm rank} V = n$, and $\rho$ extends to a $\R$-bilinear non
 degenerate form $\rho_\R$. We know by Sylvester's theorem that a
 complete isometry invariant of $\rho_\R$ is given by the pair of {\it
   positivity} and {\it negativity indices} $(i_+(\rho_\R),
 i_-(\rho_\R))$, where $i_\pm(\rho_\R)$ is the maximum of dimensions
 of $\R$-linear subspaces of $V_\R$ such that the restriction of
 $\rho_\R$ to them is either positive or negative definite.  Clearly
 this pair of indices is also an isometry invariant for the
 $\Z$-bilinear form $\rho$. We set
 $$\sigma(\rho) = i_+(\rho_\R) - i_-(\rho_\R)$$ which is called the
 {\it signature} of $\rho$ (some authors call it the {\it index} of
 $\rho$).  As ; $i_+ + i_- = n$, then $\sigma \equiv n$ mod $(2)$ and
 $$ (i_+, i_-) = (\frac{n+\sigma}{2}, \frac{n-\sigma}{2}) \ . $$ 
 The signature is {\it additive with respect to orthogonal direct sum}:
 $$\sigma (\rho \perp \rho') = \sigma (\rho) + \sigma (\rho') \ . $$
   \smallskip

 We can distribute the unimodular $\Z$-forms into the 
following classes which are clearly invariant up to isometry.  
\smallskip

 {\bf (Definite/indefinite)} $(V,\rho)$ is {\it definite} either {\it
   positive} or {\it negative} if either for every $v\in V$, $v\neq
 0$, $\rho(v,v)>0$ or $\rho(v,v)<0$. Otherwise, $\rho$ is {\it
   indefinite}.
 
 \smallskip
 
 {\bf (Parity)} $(V,\rho)$ is {\it even} if for every $v\in V$,
 $\rho(v,v)\in 2\Z$ is even. If $\rho$ is not even, then it is said
 {\it odd}.  $(V,\rho)$ is even if and only if there is a basis
 $\Bb=\{v_1, \dots, v_n\}$ of $V$ such that for every $j=1,\dots, n$,
 $\rho(v_j,v_j)\in 2\Z$; in such a case this happens for every basis
 of $V$.
 
 \medskip
 
 So we have the combination sub-classes ``definite/indefinite and
 even", ``definite/indefinite and odd"; the study up to isometry can
 be made sub-class by sub-class.

 \subsection {Some basic forms}\label{basic-form}

 $ \UU_+, \ \UU_- $ are, up to isometry, the unique rank-$1$
 forms. They are both definite (of opposite sign) and odd, $\sigma
 (\UU_\pm)=\pm 1$.
 \smallskip
 
 We denote by $\HH$ the (isometry class of the) form defined on $\Z^2$
 by
 $$(x,y)\to x^tHy$$
where 

$$ H:=
\begin{pmatrix}
0&1\\
1&0
\end{pmatrix}$$ 
\smallskip

The form $\HH$ is indefinite and even; $\sigma (\HH)=0$.
\smallskip

Let us denote by $\EE_8$ the (isometry class of the) form defined on
$\Z^8$ by
$$(x,y)\to x^tEy$$
where $E=(e_{i,j})$ is the symmetric matrix $8 \times 8$ such that:
\smallskip

- For every $i$, $e_{i,i} =2$;

- For $i=1,\dots, 6$,  $e_{i,i+1} = 1$; 

-   $e_{5,8} = 1$;

-  $e_{i,j} = 0$  otherwise.

\smallskip

One verifies by direct computation that $\EE_8$ is unimodular, even,
positive definite; hence $\sigma(\EE_8)=8$.  $-\EE_8$ (that is the
isometry class of $(\Z^8, -E)$) is even, negative definite with
$\sigma (-\EE_8)=-8$.  Being even $\pm\EE_8$ is {\it not
  diagonalizable}, that is it is not isometric to $8\UU_\pm$.
\smallskip

\subsection{Full classification up to rank $4$}\label{rank4class}
We have

\begin{proposition}\label{rank4}   Isometry classes of symmetric unimodular $\Z$-bilinear forms 
of rank $n$ up to $4$ either are diagonalizable (i.e. they admit a
orthonormal basis) or are even with null signature.  The normal
representatives are respectively:
\smallskip

(1) {\bf (Diagonalizable)}  The normal representative is
$$ |\sigma|\UU_\epsilon \perp \frac{n-|\sigma|}{2}(\UU_+ \perp \UU_-)$$
where $\epsilon$ is the sign of the signature $\sigma$.
\smallskip

(2) {\bf (Even)} The normal representatives are either $\HH$ or $2\HH$.
\end{proposition}

\cvd

The key geometric fact to get this result is that for every $(V,\rho)$
such that rank$(V)\leq 4$, there is $v\neq 0$ in $V$ such that
$|\rho(v,v)|<2$; this is an application of a theorem of Minkowski on
the volume of lattice in euclidean spaces.
 
\subsection{Classification of indefinite forms}\label{indef-class}
This is summarized in the following theorem.

\begin{theorem}\label{indef-class}
(1) The triple
$$ ({\rm rank},\ {\rm signature}, \ {\rm parity})$$
is a complete invariant for the indefinite forms considered up to isometry. 
\smallskip

(2) For every indefinite isometry class we have the following distinguished representative, 
depending on the parity:
 \smallskip

  {\bf (Indefinite and odd normal representatives)} $\ $ For every
  rank $n$ and signature $\sigma$ this is
$$ |\sigma|\UU_\epsilon \perp \frac{n-|\sigma|}{2}(\UU_+ \perp
  \UU_-)$$ where $\epsilon$ is the sign of $\sigma$. Hence indefinite
  odd forms are diagonalizable, that is they admit orthonormal basis.

  \smallskip

  {\bf (Indefinite and even normal representatives)} $\ $ For every
  rank $n$ and signature $\sigma$, $\sigma \equiv 0$ mod $(8)$,
  $n-|\sigma|$ is even and non zero and the normal representative is
$$ \frac{\sigma}{8} \EE_8 \perp \frac{n-|\sigma|}{2} \HH $$
where we mean $a\EE_8 = -a(-\EE_8)$ if $a<0$.
\end{theorem}

\cvd

\smallskip

The key fact for the indefinite classification is the number-theoretic
{\it Meyer theorem} which states that for every indefinite $(V,\rho)$,
there is $v\neq 0$ in $V$ such that $\rho(v,v)=0$. If $n\leq 4$ this
follows from the above full classification.  If $n\geq 5$, via the
extension of coefficients $\Z \subset \Q$, one is reduced to prove
that, alike for $\R$-spaces, a scalar product on a $\Q$-vector space
of dimension $n \geq 5$ is definite if and only if for every non zero
vector $v$, $\rho(v,v)\neq 0$. Note that the last statement fails for
$n=4$. The proof is based on {\it Hasse-Minkowski Theorem}.  Then the
indefinite odd case follows by a rather easy inductive argument.  An
important relation to achieve the odd case is:
 $$ \HH \perp \UU_\pm = \UU_\mp \perp 2 \UU_\pm \ . $$ 
\smallskip

The classification in the indefinite and even case is more delicate,
employs the already achieved odd classification and involves in the
very statement certain congruence mod $(8)$. We limit to clarify this last point.

\subsection{Characteristic elements and congruences mod $(8)$}\label{cong-mod-8}
Let $(V,\rho)$ be as above. An element $u\in V$ is by definition {\it
  characteristic} if for every $v\in V$, $\rho(v,v)\equiv \rho(u,v)
\ {\rm mod} (2)$. We have the following so called {\it van der Blij}
lemma.

\begin{lemma} (1) For every $(V,\rho)$ there are characteristic elements.

(2) For every characteristic element $u$, $\sigma \equiv \rho(u,u) \ {\rm mod} (8)$.

(3) If $\rho$ is even then $\sigma \equiv 0 \ {\rm mod} (8)$.
\end{lemma}

\Dim (1): fix a basis of $V$, so that $V\sim \Z^n$ and let the
$n\times n$ symmetric matrix $A$ represent the form $\rho$.  By
reducing mod $(2)$, we have the $\Z/2\Z$-linear function
$(\Z/2\Z)^n\to \Z/2\Z$, $y\to y^tAy$.  As $\det A = 1 \ {\rm mod}
(2)$, there is a unique representing vector $\bar u \in (\Z/2\Z)^n$
such that for every $y$, $y^tAy=\bar u^t A y $. Every $u\in \Z^n$
whose reduction mod $(2)$ is equal to $\bar u$ is a characteristic
element of $\rho$.

As for (2), if $u$ and $u'$ are characteristic elements, so that
$u'=u+2x$ for some $x\in V$, then $\rho(u',u')=\rho(u,u)+4(\rho(u,x) +
\rho(x,x)) \equiv \rho(u,u) \ {\rm mod} (8)$. So $\rho(u,u)$ is
invariant mod $(8)$. It is additive with respect to the orthogonal
direct sum and it holds $\pm 1$ on $\UU_\pm$.  Then item (2) holds for
indefinite and odd forms thanks to the classification in this case. On
the other hand, $\rho \perp \UU_+ \perp \UU_-$ has the same signature
of $\rho$ and is indefinite and odd; so (2) holds in general.

Item (3) is an immediate corollary of (2).

\cvd

\smallskip  

\subsection{Indefinite stabilizations}\label{stab}
Given any form $\rho$ there are simple ways to transform it into an
indefinite one. The first is called {\it elementary odd
  stabilizations}:
$$ \rho \to \rho \perp \UU_\epsilon$$ for a suitable $\epsilon = \pm$,
the resulting form is indefinite and odd. The signature changes by
$\sigma \to \sigma \pm 1$.
$$ \rho \perp (\UU_+ \perp \UU_-)$$ is always indefinite odd and the signature does not change.

The {\it elementary even stabilization} is
  $$ \rho \to \rho \perp \HH$$ the resulting form is indefinite and is
even if and only if $\rho$ is even. The signature does not change.
  
  Then the classification of indefinite odd forms induces a
  classification of {\it all} forms up to such odd
  stabilizations. Similarly, the classification of indefinite even
  forms induces a classification of all {\it even} forms up to even
  stabilization. In particular we have:
  \smallskip
  
  {\it For every pair of forms $\rho$ and $\rho'$ there are $m_1,m_2,m'_1,m'_2,m \in \N$ such that
  $$ \rho \perp m_1\UU_+ \perp m_2\UU_- = \rho' \perp m'_1\UU_+ \perp m'_2\UU_- = m(\UU_+ + \UU_-) \ . $$}
  
  \smallskip

 \subsection{Neutral forms and the Witt group}\label{WittZ}
 Similarly to Section \ref{Witt}, denote by $I(\Z)$ the set of
 isometry classes of unimodular symmetric $\Z$-bilinear forms defined
 on free $\Z$-modules of arbitrary finite rank. The operation $\perp$
 makes it a semigroup. $S\in I(\Z)$ is said {\it neutral} if rank $S =
 2m$ is even and there is a submodule $Z\subset S$, rank $Z=m$ such
 that $Z=Z^\perp$. The following lemma is an immediate consequence of
 Theorem \ref {indef-class}.

 \begin{lemma}\label{indef-neutral}
   An indefinite odd class is neutral if and only if it is of the form
   $m(\UU_+ \perp \UU_-)$ for some $m\geq 1$.  An indefinite even
   class is neutral if and only if it is of the form $m\HH$ for some
   $m\geq 1$.
  \end{lemma}
  
  \cvd

\smallskip

Put on $I(\Z)$ the equivalence relation $X\sim X'$ if and only if
there are neutral spaces $S,S'$ such that
$$X\perp S = X' \perp S' \ . $$ Denote by $W(\Z)$ the quotient
set. The operation descends to $W(\Z)$ and makes it an abelian group
called the {\it Witt group} of the ring $\Z$.  All this can be
restricted to the set $I_0(\Z)$ of even classes and gives rise to the
restricted Witt group $W_0(\Z)$. Also the following proposition is an
easy consequence of Theorem \ref {indef-class}.
\begin{proposition}\label{witt-iso}
  Both following maps are well defined group isomorphisms:
$$\sigma: W(\Z) \to (\Z,+), \ \frac{\sigma}{8}: W_0(\Z) \to (\Z,+) \ . $$
Moreover, $W(\Z)$ is generated by $\UU_+$ while $W_0(\Z)$ is generated by $\EE_8$.
\end{proposition}

\cvd

\section{Some $4$-manifold counterparts}\label{partial4D}
  In analogy with the surface case, one would like to determine
  $4$-manifold couterparts of the above abstract theory, at
  least for indefinite forms where the arithmetic classification is
  complete. In particular one would wonder that every indefinite
  normal representative is realized as the intersection form
  $\bullet_M$ of some $4$-dimensional smooth manifold $M$ as
  above. Unfortunately this is too optimistic.
  \smallskip
  
  {\bf Notation:} We set $\sigma_{\bullet_M}=\sigma(M)$.
  \smallskip
    
  First we establish a topological counterpart of the operation
  $\perp$. This is analogous to surface Lemma
  \ref{connected-direct-sum}.
  
  \begin{lemma}\label{s-cs}
    Let $(M_1,\bullet_{M_1})$ and $(M_2,\bullet_{M_2})$ be
    $4$-manifolds equipped with the respective intersection forms and
    set $M= M_1 \cs M_2$. Then, up to isometry,
  $$ \bullet_{M}= \bullet_{M_1} \perp \bullet_{M_2} \ . $$
  \end{lemma}
  \Dim Let $\alpha=[F]\in \Hh_2(M;\Z)$ where $F$ is a proper oriented
  surface embedded into $M$.  Up to isotopy we can assume that
  $F\pitchfork S$, where $S$ is a smooth $3$-sphere in $M$ which
  realizes the connected sum splitting of $M$. $L= F\cap S$ is a link
  in $S\sim S^3$. Then $M$ is obtained by gluing $M'_j= M_j\setminus
  {\rm Int} D^4$, $j=1,2$, along the two boundary components of a
  tubular neighbourhood $N(S)\sim S^3 \times [-1,1]$ of $S$ in
  $M$. $F_j=F\cap \hat M'_j$ is a proper submanifold of $M'_j$ with
  boundary $L$. $F_j$ can be capped by means of a Seifert surface of
  $L$ in $S^3$. So we get boundaryless surfaces $\hat F_j$ in $M_j$
  which up to isotopy can be embedded into $M'_j$. Hence, via the
  isomorphism induced by the inclusions and a slight abuse of
  notation, we have $[F]=[\hat F_1]+[\hat F_2]$. Doing in a similar
  way for another class $\alpha'=[F']$, we get $\alpha \bullet \alpha'
  = [\hat F_1]\bullet [\hat F'_1] + [\hat F_2]\bullet [\hat F'_2]$.
  \smallskip
  
  \cvd
  
  \smallskip
  
  \begin{remark}
    {\rm We stress that we are {\bf not} claiming that every direct
      sum decomposition of an intersection form $\bullet_M$
      corresponds to a connected sum decomposition of the manifold $M$
      (see Example \ref{E8}).}
  \end{remark}
  
 It is easy to realize $\UU_\pm$ and $\HH$. In fact: 
  \smallskip
  
  $\UU_\pm$ is the intersection form of $\pm \PP^2(\C)$, where
  $\PP^2(\C)$ is endowed with the natural orientation as a complex
  manifold.  $\Hh_2(\PP^2(\C);\Z)$ is generated by $[\PP^1(\C)]$ that
  is represented by any complex line embedded into $\PP^2(\C)$. Hence
  {\it every indefinite and odd normal representative can be
    realized}.
  
  \smallskip
  
  {\bf Notation:} To simplify the notation, set $\Pp= \PP^2(\C)$ and $\Qq = - \PP^2(\C)$.
  
  \medskip
  
  $\HH$ is the intersection form of $S^2\times S^2$, where $S^2$ has
  the usual orientation and we take the product orientation.
  $\Hh_2(S^2\times S^2; \Z)$ has as basis $[S^2\times \{p\}]$ and
  $[\{p\}\times S^2]$ for any $p\in S^2$.
  
  \smallskip
  
  \begin{remark}\label{simplycon}
    {\rm Both $\PP^2(\C)$ and $S^2\times S^2$ are simply connected.
      By Van Kampen theorem, the connected sum of two simply
      connected manifolds is also simply connected. So it makes sense
      (and we will do it at some point) to restrict the discussion to
      simply connected manifolds.}
 \end{remark} 
  
  \smallskip
  
  $\HH$ and $\UU_+ \perp \UU_-$ are the basic neutral classes. As for
  their $4$-dimensional realizations we have
  
  \begin{proposition}\label{sphere-on-sphere}
    Up to isomorphism of fibre bundles, there are two distinct fibre
    bundles over $S^2$ with fibre $S^2$ and orientable total space;
    $S^2\times S^2$ and $\Pp\cs \Qq :=  S^2\tilde \times S^2$ are the
    respective total spaces.
   \end{proposition} 
    \Dim By at theorem of Smale \cite{S1} (recall also Section \ref
         {twisted-sphere}) Diff$^+(S^2)$ retracts by deformation to
         $SO(3)\sim \PP^3(\R)$. Then there are exactly two such fibre
         bundles because $\pi_1(SO(3))\sim \Z/2\Z$ (recall Section
         \ref{VB-sphere}).  $\Pp\cs \Qq$ can be obtained by the
         complex blow up of $\PP^2(\C)$ at a point.  It follows from
         the proof of Proposition \ref{C-BU-point} that it is the
         total space of a fibre bundle as in the statement of the
         proposition.  More precisely, let $\Dd$ be the unitary disk
         in an affine chart of $\Pp$ at a point $x_0\sim 0$. Then
         $\BB_{\C}(\Dd,0)$ is the oriented total space of a fible
         bundle over the Riemann sphere $S^2\sim \PP^1(\C)$ with fibre
         $D^2$; the fibres are given by the strict transform of the
         intersection with $\Dd$ of the complex lines through $0$. Set
         $\Pp_0:= \Pp \setminus {\rm Int} \Dd$.  Also $\Pp_0$ is the
         total space of a fibre bundle of the same type. Considering
         $\PP^1(\C)\subset \Pp_0$, the fibres are given by the
         intersection with $\Pp_0$ of the complex lines passing
         through $0$ and $x\in \PP^1(\C)$. The restriction of these
         fibres to $\partial \Dd$ induce the Hopf fibration $\hG: S^3
         \to S^2$. Then $\BB_{\C}(\Pp,x_0)$ is diffeomorphic to the
         double $D(\Pp_0)= \Pp_0 \amalg -\Pp_0/{\rm id}_{S^3}$ and
         hence to $\Pp \cs \Qq$.  The fibration of $\Pp\cs \Qq$ with
         fibre $S^2$ is obtained by gluing ``along the Hopf fibration"
         the two fibrations with fibre $D^2$ described so far. Finally
         $S^2\times S^2$ and $\Pp \cs \Qq$ are distinguished by the
         intersection forms.
   
   \cvd
   
   \smallskip

  Now we discuss a topological counterpart of the relation 
   $$ \HH \perp \UU_\pm = \UU_\mp \perp 2 \UU_\pm  $$
   this is analogous to surface Lemma \ref {basic-rel2}.
   
   \begin{proposition}\label{cs-relation} We have
   $$ (S^2\times S^2)\cs \Qq \sim \Pp\cs 2\Qq, \ (S^2\times S^2)\cs \Pp \sim \Qq\cs 2\Pp \ . $$
   \end{proposition}
   \Dim As $S^2\times S^2$ admits an orientation reversing
   diffeomorphism, the two relations are equivalent to each other. The
   second geometric proof of Lemma \ref{basic-rel2} applies {\it
     verbatim} to prove the first relation, provided that one replaces
   $\R$ with $\C$ everywhere.
   
   \cvd
   
  \smallskip
  
  A realization of indefinite even normal representatives, or of
  $\EE_8$ itself, possibly by means of a simply connected smooth
  $4$-manifold $M$, is much more subtle and hard question. We will
  discuss later the following fundamental Rohlin's discovery:
  \smallskip
  
  {\it If $M$ is simply connected and its intersection form is even,
    then $\sigma(M)\equiv 0$ {\rm mod} $(16)$.}
  
  \smallskip
  
 \noindent Recall that algebra tells us that the signature of an even
 form is $\equiv 0$ mod $(8)$. Then $\EE_8$ cannot be realized.  If
 $M$ is simply connected with indefinite and even intersection form,
 then this is necessarily isometric to a normal representative of the
 type
   $$ 2a\EE_8 \perp b\HH$$ for some $a\in \Z$, $b \in \N\setminus
 \{0\}$.  It is not evident (and ultimately false) that
 every such pair $(a,b)$ can be realized.  On the other hand,
 classical simply connected examples show the actual occurrence of
 $\EE_8$.
   
 \begin{example}\label{E8}
   {\rm If we relax the requirement of dealing with normal
     representatives, it is not hard to make $\EE_8$ visible. For
     example, by the indefinite and odd classification, the form of
     $M= 10\Pp \cs \Qq$ is isometric to
   $$ \EE_8 \perp \UU_+ \perp \HH \ . $$ Nevertheless, this algebraic
     decomposition does not correspond to any connected sum
     decomposition of $M$.
   
   A more substantial example, realizing a normal representative, is
   the so called {\it Kummer variety}. Let the $4$-torus $
   T^4=\R^4/\Z^4$ be realized as the product of two copies of
   $\C/(\Z\oplus i\Z)$ so that $T^4$ has a complex $2$-manifold
   structure with ``uniformizing" complex coordinates $(w_1,w_2)$. The
   involution $\tau(w_1,w_2)=(-w_1,-w_2)$ descends to $T^4$ and has 16
   fixed points. Let us perform the complex blow-up at such fixed
   points. We get a complex surface $\tilde K$, smoothly diffeomorphic
   to $T^4 \cs 16\Qq$. The exceptional complex surface over each fixed
   point is a Riemann sphere $S$ with self-intersection number in
   $\tilde K$ equal to $-1$. The involution $\tau$ lifts to an
   involution $\tilde \tau$ of $\tilde K$ which is the identity on
   each exceptional sphere. We consider the quotient
   $$K := \tilde K/\tilde \tau \ . $$ One verifies that $K$ is a
   smooth complex surface. By means of the natural projection, every
   exceptional sphere $S$ maps onto a $2$-sphere $S'$ embedded into
   $K$; the restriction of the projection on a suitable neighbourhood
   of each $S$ in $\tilde K$ is a double covering of a neighbourhood
   in $K$ of the corresponding sphere $S'$. Then the self-intersection
   number of every $S'$ in $K$ is equal to $-2$. One can verify that
   $\Hh_2(T^4;\Z)\sim \Z^6$ and is generated by six embedded $2$-tori,
   while $\Hh_2(K;\Z) \sim \Z^{22}$ generated by the image of these
   tori together with the 16 spheres $S'$. Eventually the intersection
   form of $K$ is idefinite and even with normal representative $
   -2\EE_8 \perp 3\HH$.}
   \end{example}   
   
   \subsection{On the intersection form of $4$-manifolds with boundary}\label{with-bound}
\begin{figure}[ht]
\begin{center}
 \includegraphics[width=5cm]{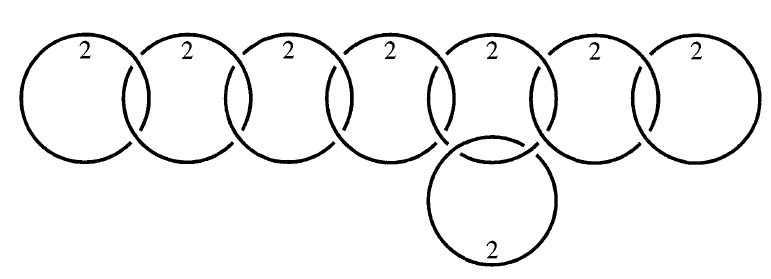}
\caption{\label{E8link} A $\EE_8$-link.} 
\end{center}
\end{figure}

   If $\partial M \neq \emptyset$, the intersection form
   $\sqcup: \Hh^2(M;\Z)\times \Hh^2(M;\Z)\to \Z$
   and the $\Z$-linear map
   $$ \hat \phi^2: \Hh^2(M;\Z) \to {\rm Hom}(\Hh_2(M;\Z),\Z)$$
   are defined as well. In general the form is not unimodular.
   If $\beta:=i_*(\alpha)\neq 0 $ in $\Hh_2(M;\Z)$ for some $\alpha \in \Hh_2(\partial M;\Z)$,
   then $\beta \sqcup \gamma = 0$ for every $\gamma$.
   On the other hand, it follows from the results of Chapter \ref {TD-LINE-BUND}
   that
   $$ \hat \phi^2: \Hh^2(M,\partial M;\Z) \to {\rm Hom}(\Hh_2(M;\Z),\Z)$$
   is an isomorphism. Hence the intersection form of $M$ is unimodular
   if and only if $j_*: \Hh_2(M;\Z)\to \Hh_2(M,\partial M;\Z)$ is an isomorphism.
   For simplicity assume that $M$ is part of a triad of the form $(M,\emptyset, V=\partial M)$
   admitting an ordered handle decomposition with one $0$-handle, some $2$-handles, say $k$,
   no $3$ and $4$-handles. In other words, by removing the $0$-handle, we realize
   a surgery equivalence $S^3 \sim_\sigma V$. Hence $V$ is connected and $M$ is simply connected.
   We claim that every symmetric $\Z$-bilinear form (not necessarily unimodular) can be
   be realized by such a $4$-manifold. Let us sketch the argument. By using Section \ref{CW}
   we see that $M$  retracts to a wedge of $k$ $2$-spheres. By using the bordism
   homotopy invariance and what we know about the bordism of $S^2$, we see that
   $\Hh_2(S^2;\Z)$ has rank $k$; a geometric basis $\alpha_1, \dots, \alpha_k$ can be obtained by completing the core of every
   $2$-handle with a Seifert surface of the corresponding attaching knot in $S^3$ (provided the handles have been ordered).
   The $k$-components framed link in $S^3$ which encodes the attaching of $2$-handles
   carries a symmetric {\it linking matrix} made by the linking numbers of pairs
   of constituent knots and, along the diagonal, by the integers encoding the framing
   of every  such a knot. With a bit of work one eventually realizes that 
   this matrix equals the matrix of the intersection form of $M$ with respect to the above
   geometric basis. In Figure \ref{E8link} we show a framed link in $S^3$ which realizes
   $\EE_8$; $\partial M$ is the {\it Poincar\'e sphere}.
   
   We have
   \begin{proposition}\label{hom-sphere} The intersection form of $M$ is unimodular if and only if
   $\Hh_1(V;\Z)=\Hh_2(V;\Z)=0$.
   \end{proposition}
   \Dim As $M$ is simply connected, $\Hh_3(M,\partial M;\Z) \sim {\rm Hom}(\Hh_1(M;\Z),\Z) =0$.
   Hence by using the bordism long exact sequence of $(M,\partial M)$, we see that
   $i_*:\Hh_2(V;\Z)\to \Hh_2(M;\Z)$ is injective; hence if the intersection form of $M$ is unimodular,
   then $\Hh_2(V;\Z)=0$. On the other hand,  if $\Hh_1(V;\Z)=0$, consider  the dual handle 
   decomposition; the cores of the $2$-handles 
   provide a basis of $\Hh_2(M,\partial M;\Z)$; by capping each of them with a Seifert
   surface in $V$ of the corresponding attaching knot, we get a further geometric basis
   of $\Hh_2(M;\Z)$ dual to the previous one.
   
   \cvd
   
   \smallskip
   
   If the intersection form of $M$ is unimodular, possibly by performing an elementary blow-up move 
   (which replaces $M$ with $M\cs \pm \PP^2(\C)$, without modifying the boundary $V$),
   we can assume that the unimodular intersection form of $M$ is diagonalizable.
   If one $2$-handle (corresponding to a costituent knot $K_i$) is slid over another, say corresponding to $K_j$,
   then the geometric basis as above,  changes by sending $\alpha_i$ to $\alpha_i + \alpha_j$, and the linking
   matrix changes by adding the $j^{th}$ row to the $i^{th}$ row, and the $j^{th}$ column to the $i^{th}$.
   It follows that we can realize a diagonalizing basis by means of handle sliding.
   
   The same discussion can be repeated (with some simplification) by replacing everywhere $\Z$ with $\Z/2\Z$.
    
\section{$\Omega_4$}\label{omega4}  
 
 We already know that $\Omega_4$ is non trivial because
 $\chi_{(2)}(\PP^2(\C))=1$. More precisely we have a surjective
 homomorphism defined by
 $$\chi_{(2)}: \Omega_4 \to \Z/2\Z, \ \chi_{(2)}([M]):= \chi_{(2)}(M) \ . $$
 
 Pontryagin remarked that there is a subtler homomorphism induced by
 the signature.  As usual
 $$ [M\cs M']= [M\amalg M']=[M]+[M'] \in \Omega_4$$
 so that every $\alpha\in \Omega_4$ can be represented by connected $4$-manifolds and
 we can replace $\amalg$ with $\cs$ to define the $\Z$-module operation on $\Omega_4$.
 Then we have
 \begin{proposition} $$\sigma: \Omega_4 \to \Z, \ \sigma (\alpha):= \sigma(M)$$
 where $M$ is any connected representative of the class $\alpha$, is a well defined and surjective
 homomorphism.
 \end{proposition}
 \Dim As the signature is additive with respect to the connected sum,
 $\sigma(M)=-\sigma (-M)$ and $\sigma(\PP^2(\C))=1$, it is enough to
 show that if $[M]=0\in \Omega_4$, then $\sigma(M)=0$.  To compute the
 signature, that is the indices $i_+, i_-$, it is enough to extend the
 coefficients $\Z\subset \Q$.  For every $\alpha \in \Hh_2(M;\Q)$
 there exists $m\in \Z$ such that $m\alpha = \alpha' \in \Hh_2(M;\Z)$,
 and $\alpha \bullet \alpha = \alpha' \bullet \alpha'/(m^2)$. If for
 every $\alpha \in \Hh_2 (M;\Q)$, $\alpha\bullet \alpha =0$, then
 $\sigma =0$. Let $M= \partial W$, $i: M \to W$ be the inclusion.  If
 $i_*(\alpha)=0$, then $\alpha'\bullet \alpha' = 0$, hence
 $\alpha\bullet \alpha =0$.  So if for every $\alpha$,
 $i_*(\alpha)=0$, then $\sigma =0$. Assume that $i_*(\alpha)\neq
 0$. Then there is $b\in \Hh_3(W,M;\Q)$ such that $\beta:=\partial b
 \in \Hh_2(M;\Q)$ and $\alpha \bullet \beta = 1$, $i_*(\beta)=0$.  Let
 $V$ be the subspace of $\Hh_2(M;\Q)$ generated by $\alpha$ and
 $\beta$. The matrix of the restriction of the intersection form on
 $V$ has $\det = -1$, hence its signature is equal to zero. As the
 restriction of the form to $V$ is non degenerate, also its
 restriction on the othogonal space $V^\perp$ is non degenerate. The
 we can iterate the construction till one finds classes such that
 $i_*(\alpha)\neq 0$.  By the additivity of the signature with respect
 to the orthogonal direct sum, we conclude that $\sigma=0$.
 
 \cvd
 
 \smallskip
 
 We are ready to state and prove the following theorem due to Rohlin
 (see \cite{GM}). We will propose his original argument.
 This is formally analogous to surface Theorems
 \ref{omega-eta-2} and \ref{eta-2}.
 
 \begin{theorem}\label{omega4=Z}
   The homomorphism induced by the signature $\sigma: \Omega_4 \to \Z$
   is an isomorphism. Hence $\Omega_4$ is generated by $[\PP^2(\C)]$
   and is naturally isomorphic to the Witt group $W(\Z)$.
 \end{theorem}
 \Dim The restriction of $\sigma$ to the submodule of $\Omega_4$
 generated by $[\PP^2(\C)]$ is an isomorphism onto $\Z$.  Hence it is
 enough to show that $\Omega_4$ is generated by $[\PP^2(\C)]$.  We
 will achieve this fact by several steps.  Let $M$ be as usual a
 compact, oriented, connected and boundaryless $4$-manifold.
 
 {\bf Step 1.} This is similar to the first step in Rohlin's proof that $\Omega_3=0$.
 That is, up to bordism, it is not
 restrictive to assume that $M \subset \R^7 \subset S^7$.
 (see also \cite{Kirby} for a somewhat different conclusion of the proof based on 
 Step 1).
 
 \smallskip
 
 {\bf Step 2.} We would like to construct along $M$ a field $v$ of
 unitary tangent vector to $S^7$ normal to $M$.  This is not possible
 in general, however we are going to see that there is $\tilde M := M
 \cs a\Pp \cs b\Qq \subset S^7$ for some $a,b \in \N$, which carries
 such a nowhere vanishing transverse field.  A first obstruction is
 given by the Euler class $e\in \Hh^3(M;\Z)$ of a normal bundle to $M$
 in $S^7$. On the other hand, $[M]=0 \in \Hh^3(S^7;\Z)$ and $e=
 i^*([M])=0$. This implies that such a field $v$ can be defined on
 $M_0 = M \setminus {\rm Int} B^4$, where $B^4$ is a smooth $4$-disk
 in $M$; in fact $M_0$ has a $3$-dimensional spine, $v$ can be always
 constructed up to the $2$-skeleton and the obstruction to extend it
 to the third skeleton belongs to $\pi_2(S^2)$ and vanishes because
 $e=0$. The restriction of $v$ to $\partial M_0$ defines an element of
 $\pi_3(S^2)$ which is in general non trivial. This is the final
 effective obstruction to extend $v$ on the whole of $M$. We know that
 $\pi_3(S^2)=\Z$ is generated by the Hopf map $\hG: S^3\to S^2$. By
 transversality we can perturb the field $v$ and assume that it is
 defined on $M'$ obtained by removing from $M$ the interior of a
 finite number of disjoint $4$-disks $B_j$ embedded into ${\rm Int} B$
 such that the restriction of $v$ to every boundary $\partial B_j$ is
 equal to $\pm \hG$.  By using the field $v$ we get an embedding of
 $M'$ into the boundary $\partial N(M)$ of a tubular neighbourhood of
 $M$ in $S^7$. By abstractly gluing to every boundary component of
 $M'$ the mapping cylinder of the corresponding map $\pm \hG$, we get
 the $4$-manifold $\tilde M := M\cs a\Pp \cs b\Qq$ for some $a,b \in
 \N$. We claim that we can assume that $\tilde M \subset \partial
 N(M)$ by extending the given embedding of $M'$.  For if $B_j \times
 D^3$ is a trivialized chart of $N(M)$ over the $4$-ball $B_j$, the
 embedding of $\partial B_j$ is for instance of the form $x\to
 (x,\hG(x))$ and $\Pp_0$ is the copy of $\Pp \setminus {\rm Int}(D^4)$
 in $\tilde M$ corresponding to $B_j$, then an embedding of $\Pp_0$ is
 given (by using suitable homogeneous coordinates $(x_0,x_1,x_2)$) by:
 $$(x_0,x_1,x_2) \to (( \frac{2x_0x_1}{\sum_{i=0}^2 |x_i|^2},
 \frac{2x_0,x_2}{\sum_{i=0}^2 |x_i|^2}),\frac{x_1}{x_2}) \in B_j
 \times \PP^1(\C) \ . $$ Clearly, the restriction $\tilde v$ to
 $\tilde M$ of a unitary normal field to the hypersurface $\partial
 N(M)$ in $S^7$ is nowhere vanishing along $\tilde M$.
 \smallskip
 
 {\bf Step 3.} The field $\tilde v$ determines an embedding of a copy
 $\hat M$ of $\tilde M$ into $\partial N(\tilde M)$ the boundary of a
 tubular neighbourhood $\pi: N(\tilde M)\to \tilde M$ of $\tilde M$ in
 $S^7$.  Set $X:= S^7 \setminus {\rm Int} N(\tilde M)$. If $[\hat M]$
 would be zero in $\Hh_4(X;\Z)$, then it should be a boundary thanks to
 Proposition \ref{seifert-surf2}, and finally $M$ bordant with $k\PP^2(\C)$
 for some $k\in \Z$. However, we cannot assume that $[\hat M]=0$.
 \smallskip

 {\bf Claim 1.}  {\it There is an oriented surface $F$ in $\tilde M$
   such that the disjoint union of inclusions $ j: \hat M \amalg \partial \pi^{-1}(F) \to
   \partial N(\tilde M)$ represents zero in $\Hh_4(X; \Z)$ (the
   $4$-manifold $S:= \partial \pi^{-1}(F)$ is oriended by the direct
   sum of the orientation of $F$ and the orientation of the normal
   bundle of $\hat M$ in $\partial N(\tilde M)$). }
\smallskip
 
 Let us prove the claim. $\Hh_4(S^7;\Z)=0$, more precisely
 $\Omega_4(S^7) \sim \Omega_4$. Hence there is an oriented triad $ (W,
 \hat M, V)$ and a map $h: W \to S^7$ where the restriction to $\hat
 M$ is the inclusion and the restriction to $V$ is a constant map. By
 transversality we can assume that the restriction of $h$ to an open
 collar of $\hat M$ in $W$ is an embedding in $X$ transverse to
 $\partial N(\tilde M)$, the image of $V$ is in the interior of $X$,
 the restriction of $h$ to the interior of $W$ is transverse to
 $(N(\tilde M)$, $\partial N(\tilde N)$ and $\tilde M$. Then $F=
 h({\rm Int}(W))\cap \tilde M$ is a surface in $\tilde M$ and $ h({\rm
   Int}(W))\cap \partial N(\tilde M)= \partial \pi^{-1}(F):=
 S$. Finally $(h^{-1}(X),h)$ realizes a bordism between $(\hat M
 \amalg \partial S, j)$ and $(V,h_|)$. The claim is proved.
 \smallskip

 With a slight abuse of notation we write $[\hat M\amalg S]$ instead of
  $[\hat M\amalg S, j]$.
 
 {\bf Step 4.}  Let $F$ be as in Claim 1. Clearly $S:= \pi^{-1}(F)$ is
 the boundary of a $2$-disk bundle.  Then it would be enough to prove
 that $[\hat M \amalg S]=0$.  We are able to do it under
 a more restrictive hypothesis. We can assume that $\hat M$ is
 transverse to $S$ in $\partial N(\tilde M)$ and that $\hat M \cap S =
 F_1$ where $F_1$ is the copy of $F$ in $\partial N(\tilde M)$
 determined by the above normal unitary field $\tilde v$ along $\tilde
 M$.
\smallskip

{\bf Claim 2.}  {\it Assume that the oriented normal bundle to $F_1$
  in $\hat M$ is isomorphic to the oriented normal bundle of $F_1$ in
  $S$. Then $[\hat M \amalg S]=0$. }
\smallskip

Let $U_1$ be a tubular neighbourhood of $F_1$ in $\hat M$, $U_2$ a
tubular neighbourhood of $F_1$ in $S$. We can construct a manifold $Y$
by gluing $\hat M \setminus {\rm Int}(U_1)$ and $S\setminus {\rm
  Int}(U_2)$ along the boundary which are isomorphic by hypothesis. In
fact $Y$ can be realized within $\partial N(\tilde M)$ in such a way
that it contains isotopic copies of the original constituent pieces.
It is not hard to check that $Y$ is bordant with $\hat M \amalg S$ and
that $[Y]=[\hat M]+[S]= 0$ in $\Hh_4(X;\Z)$. By Proposition
\ref{seifert-surf2} $Y$ is a boundary and hence also $M \amalg S$ is
so.
\smallskip

{\bf Step 5.}  In general the normal bundles of $F_1$ in $\hat M$ and
$S$ respectively are not isomorphic to each other. The oriented rank-2
normal bundle of $F_1$ in $\hat M$ is determined up to isomorphism by
the self-intersection number of $F_1$ in $\hat M$.  One realizes that
by performing a complex (anti) blow up of $\hat M$ at a point of $F_1$
we get a manifold $\hat M'$ diffeomorphic to $\hat M \# \pm \PP^2(\C)$
such that the strict transform of $F_1$ in $\hat M'$ is equal to $F_1$
and its self intersection number varies by $\pm 1$.  Moreover,
it is not restrictive to assume that $\hat M'$ is realized within $\partial
N(\tilde M)$.  By iterating this construction we eventually get $\hat
M' \sim \hat M \# p\Pp \# q \Qq = M \# k\Pp \# h\Qq$
to which Claim 2
applies. Theorem \ref{omega4=Z} is eventually achieved.

\cvd

\section{Simply connected classification up to odd stabilization}\label{odd-stab}
In this section we restrict to {\it simply connected} $4$-manifolds.
We are going to prove:

\begin{theorem}\label{scos} For every compact oriented simply connected boundaryless 
$4$-manifold $M$, there exist $(k,h), (m,n) \in \N \times \N$ such that
$ M \cs k \Pp \cs h \Qq = m\Pp \cs n\Qq$.
\end{theorem}

By using Proposition \ref{cs-relation} one can slightly refine the statement in the form:
\smallskip

{\it $\cdots$  there exists $(k,m)\in \N\times \N$ such that
  $M\cs (k+1)\Pp \cs k\Qq = (m+1)\Pp \cs m\Qq$.}
\smallskip

Theorem \ref{scos} is analogous to surface Section
\ref{stable-2-Nash}, however we have not here any {\it a priori}
information about the integers $k,h,m,n$.  By Theorem \ref{omega4=Z},
for every $M$ as above there is $l \in \Z$ such that $M\cs l\PP^2(\C)$
and this last is still simply connected; then Theorem \ref{scos} will
readily follow by combining the next proposition with Proposition
\ref{cs-relation}.

\begin{proposition}\label{sub-scos}
  Let $M$ be simply connected and a boundary. Then there are
  $(k_0,k_1), (h_0,h_1) \in \N\times \N$ such that $$M\cs k_0(S^2
  \times S^2)\cs k_1(S^2 \tilde \times S^2)\sim h_0(S^2\times S^2)\cs
  h_1(S^2 \tilde \times S^2) \ . $$
\end{proposition}
\Dim As $M$ is a boundary, there is an oriented triad $(W,M, S^4)$.
Let us take an ordered handle decomposition of $(W,M,S^4)$ without $0$-
and $5$-handles. Hence it is of the form
$$ (M\times [0,1]) \cup \{\Hh^1\} \cup \{\Hh^2)\cup \dots \cup
\{\Hh_4\} \cup ([-1,0]\times S^4)$$ where every $\Hh_j$, $j=1,\dots,
4$, denotes a pattern of $a_j$ $j$-handles attached simultaneously at
disjoint attaching tubes. We claim that we can modify the $5$-manifold
$W$ without changing the boundary $M \amalg S^4$ in such a way that it
is not restrictive to assume that $a_1=a_4=0$. To do it we apply the
``trading'' argument already used in the proof of Proposition
\ref{wallace0}.  We can assume that the attaching tube a every
$1$-handle is contained in a smooth $4$-disk of $M$. Then the new
boundary component obtained by modifying $M$ can be realized as well
by means of a $3$-handle trivially attached to $M$; thus we can trade
every $1$-handle with a $3$-handle.  By using the dual handle
decomposition we can trade every $4$-handle with a $2$-handle; so, up
to reordering, we can assume that the ordered handle decomposition of
$(W,M,S^4)$ contains only $2$ and $3$ handles. Hence $W$ can be
obtained by gluing $(M\times [0,1])\cup \{\Hh_2\}$ and $\{\Hh^3\}\cup
([-1,0]\times S^4)$ along diffeomorphic boundary components. Note that
in terms of the dual decomposition, also $\{\Hh^3\}\cup ([-1,0]\times
S^4)$ is obtained by attaching $2$-handles.  Then the following lemma
allows us to conclude.

  \begin{lemma}\label{glue2H}
    Consider the cylinder $(M\times [0,1], M_0, M_1)$, $M_j= M\times
    \{j\}$.  Let $(Y,M_0, \hat M_1)$ obtained by attaching a
    $2$-handle to $M\times [0,1]$ along $M_1$.  Assume that $M$ is
    simply connected. Then either $\hat M_1 \sim M\cs (S^2\times S^2)$
    or $\hat M_1 \sim M\cs (S^2\tilde \times S^2)$.
\end{lemma}
  \Dim As $\dim M = 4$ and $M$ is simply connected, the attaching
  $1$-sphere of the handle is isotopic to a standard $S^1$ in a chart
  of $M$. Then it is easy to check that $M_1\sim M \cs \Ff$ where
  $\Ff$ is the total space of an oriented fibre bundle over $S^2$
  with fibre $S^2$.  Then we apply Proposition \ref
  {sphere-on-sphere}. The lemma and Proposition \ref{sub-scos} are
  proved.

  \cvd

 \section{On the classification up to even stabilization}\label{even-stab}
As in the previous section we deal with simply connected
$4$-manifolds. Being very sketchy, we are going to discuss the following deeper result
\cite{Wall3}, \cite{Wall4}.

\begin{theorem}\label{even-class}
  Let $M_0$ and $M_1$ be compact oriented simply connected
  boundaryless $4$-manifolds with isometric intersection forms. Then
  there is $k\in \N$ such that $M_0 \cs k(S^2\times S^2) \sim M_1 \cs
  k(S^2\times S^2)$.
\end{theorem}

A few comments are in order:

$\bullet$ In a sense this is the strongest $4$-dimensional analogous of surface
classification in terms of the intersection form, which one has obtained by means
of classical differential/topological methods available till the ends of 70's
of the last century.

$\bullet$ Theorem \ref{even-class} implies Theorem \ref{scos}.  For up
to a suitable odd stabilization $M\cs \pm \PP^2(\C)$, this last has
the same intersection form of some $k\Pp \cs h\Qq$.  By applying to
this couple of manifolds Theorem \ref{even-class} and Proposition
\ref{cs-relation}, we get Theorem \ref{scos}. In fact a proof of
Theorem \ref{even-class} is much more demanding, it incorporates the one
of Theorem \ref{scos}, together with more advanced tools in homotopy
and homology theory beyond the limits of the present text. So we will give just
some indications. A detailed proof can be found for example in \cite{Sc}.

$\bullet$ For our main application in Section \ref{16}, the simpler classification
up to odd stabilization will suffice.

\medskip

First one proves the theorem under a stronger hypothesis. The idea is
that the $h$-cobordism theorem holds also in dimension $5$ up to even
stabilization.

\begin{proposition}\label{Wall-h}
  Let $M_0$ and $M_1$ be compact oriented simply connected
  boundaryless $4$-manifolds.  Assume that they are
  $h$-cobordant. Then there is $k\in \N$ such that $M_0 \cs
  k(S^2\times S^2) \sim M_1 \cs k(S^2\times S^2)$.
\end{proposition}

\smallskip

{\it Sketch of proof:} We know that the main difficulty to perform the
stable proof of the $h$-cobordism theorem in dimension $5$ is that we
cannot apply the Whitney trick to eliminate couples of intersection
points between the $b$-sphere $S_b$ and the $a$-sphere $S_a$ of two
{\it algebraically} complementary handles. In particular, trying to
construct a Whitney disk, we cannot avoid that such a generically
immersed $2$-disk $D$ has self-intersection points. Let $p$ such a
point. Let us make the connected sum with a copy of $S^2\times S^2$.
This contains two $2$-spheres $S_1$ and $S_2$ which intersect
transversely at one point. By means of a thin embedded $1$-handle we
connect $D$ with $S_1$ obtaining a new immersed $2$-disk $D'$ ($D'\sim
D\cs S_1$) which intersects transversely $S_2$ at one point $q$. Let
$c$ be a simple arc on $D'$ which connects $p$ and $q$ and does not
pass though other self-intersection points. By using another thin
embedded $1$-handle along $c$ we connect $D'$ with a parallel copy of
$S_2$ and get $D"$ from which both the self-intersection points $p$
and $q$ have been eliminated. Hence up to a certain number of even
stabilizations we can assume that $D$ is embedded and eventually
provides a genuine Whitney disk.

\cvd
\smallskip

The classification up to even stabilization is now a consequence of
the ``if'' implication in the the following deep Wall's theorem.

\begin{theorem}\label{if-h} Let $M_0$ and $M_1$ be compact oriented
  simply connected boundaryless $4$-manifolds.  Then they are
  $h$-cobordant if and only if they have isometric intersection forms.
\end{theorem}

Being even more sketchy: ``if'' is the hard implication; it
strenghtens a classical Whitehead theorem (based on CW complex
techniques) according to which $M_0$ and $M_1$ have the same homotopy
type.  If the intersection forms are isometric then they have in
particular the same signature, so that $M_0$ is bordant with $M_1$ by
Theorem \ref{omega4=Z} .  Arguing as in the proof of Proposition \ref
{sub-scos}, we know that there are triads $(W,M_0,M_1)$ where $W$ is
obtained by gluing some $V$ with boundary
$\partial V=M_0 \amalg (M_0\cs k(S^2\times S^2)\cs h(S^2\tilde \times
S^2))$ and some $V'$ with boundary $\partial V'=(M_1\cs k'(S^2\times
S^2)\cs h'(S^2\tilde \times S^2)) \amalg M_1$, via a diffeomorphism
$$\phi: M_0\cs k(S^2\times S^2)\cs h(S^2\tilde \times S^2)\to M_1\cs
k'(S^2\times S^2)\cs h'(S^2\tilde \times S^2) \ . $$
As $M_0$ and $M_1$ are simply connected, then also $W$ is so.
The key point is to show that, by fully exploiting the hypothesis,
amongs the triads of this kind there are such that $W$ is homologically trivial; by
standard algebraic/topological arguments this is enough to conclude
that the triad $(W,M_0,M_1)$ is a $h$-cobordism.

\section{Congruences modulo $16$}\label{16}
To introduce the theme, let us begin with a bit of history.
We have recalled in Section \ref{small-k} that by means of the hardest application
of Pontryagin method, in a series of four papers of 1951-52 (see
\cite{GM} for the translation in french and wide deep commentaries)
Rohlin eventually computed the stable homotopy group
$$\pi^\infty_3 = \pi_{n+3}(S^n)\sim \Omega_3^\Ff(S^n)  \sim \Z/24\Z,  \ n\geq 5 \ . $$
As a corollary he obtained his celebrated congruence mod$(16)$; a slightly
weaker formulation of it is as follows:

\begin{theorem}\label{rohlin} Let $M$ be a compact oriented boundaryless simply connected
$4$-manifold. Assume that its intersection form is even. Then $\sigma (M) \equiv 0$
mod $(16)$.
\end{theorem}

As $\sigma(M)$ is even, the arithmetic of unimodular forms tells us that $\sigma (M) \equiv 0$
mod$(8)$, so we can reformulate the result as
$$ \frac{\sigma(M)}{8} \equiv 0 \ {\rm mod} (2) \ . $$

This improvement by $2$ implies in particular that $\EE_8$
cannot be realized by any simply connected $4$-manifold. The derivation
of Theorem \ref{rohlin} from stably $\pi_{n+3}(S^n) \sim \Z/24\Z$
is rather demanding and uses several facts less elementary than the ones covered by the 
present text. Just to give an idea, without any pretention to be understandable,
let us sketch the argument by following \cite{MK}. 
It is shown that
$ p_1(M)= 3\sigma(M)$
where $p_1(M)$ denotes the first Pontryagin number of $T(M)$ (see Remark \ref{sw}).
This follows because both $p_1$ and $\sigma$ are bordism invariant, additive on connected sum 
and the formula holds
for the generator of $\Omega_4=\Z$. 
So it is enough to prove that $p_1(M)\equiv 0$ mod $(48)$.
One can assume that $M\subset \R^{4+n}$, $n\geq 5$.
In the hypotheses of Rohlin's theorem, one can prove that $M$ is almost parallelizable
that is the tangent bundle of $M\setminus \{x_0\}$ admit a global trivialization. 
Let $f$ be a non vanishing section of the restriction to $M\setminus \{x_0\}$ 
of the $SO(n)$ normal bundle $\nu$
of $M$ in $\R^{4+n}$. Let  $\eG$  be the obstruction to extending $f$;
it is identified with an element of  $\pi_3(SO(n))$ (which is an infinite cyclic group),
as well as the Pontryagin number $p_1(\nu)$ is identified with $\pm 2 \eG$.
Consider the $J$-homomorphism (Section \ref{J})
$J:\pi_3(SO(n))\to \pi_{3+n}(S^n)$.
One proves that $J(\eG)=0$, hence $\eG$ is divisible by $24$.
Finally one proves that  $p_1(M)=-p_1(\nu)$
because $T(M)\oplus \nu = \epsilon^{4+n}$.

An interesting feature of this history is that in the second paper of
the series, Rohlin outlined a proof of the {\it erroneous} result that
stably $\pi_{n+3}(S^n) \sim \Z/12\Z$.  Arguing as above this would
imply the non surprising congruence $\sigma(M)\equiv 0$ mod$(8)$.  In
the fourth paper, after having established the isomorphism $\sigma:
\Omega_4\to \Z$ determined by the signature (i.e. Theorem \ref
      {omega4=Z}), he firstly realized that this combined with some
      claims in his early presumed proof produced a contradiction,
      then he localized the mistake and corrected it getting the right
      group $\Z/24\Z$. In fact he pointed out that there was only one
      substantial mistake: a certain simply connected $4$-manifold $M$
      has been constructed with a characteristic element $\omega \in
      \Hh^2(M;\Z)$ of its intersection form which can be represented
      by a generic immersion $f:S^2\to M$; then by an {\it abusive}
      application of the Whitney trick in dimension $4$, he argued
      erroneously that $\omega$ was represented by an {\it embedded}
      $S^2\subset M$. This was a quite fruitful mistake: his
      correction leads to the celebrated congruence mod$(16)$ and
      provides a {\it concrete counterexample} to the applicability of
      Whitney's trick in dimension $4$. Moreover, by elaborating on
      this counterexample the authors pointed out in \cite{KM} (1961)
      an interesting extension. Recall that for every $4$-manifolds
      $M$ and for every characteristic element $\omega \in
      \Hh^2(M;\Z)$ of its intersection form
$$ \sigma(M)  -\omega \sqcup \omega \equiv 0 \ {\rm mod} (8) \ . $$
Then, assuming Theorem \ref{rohlin}, the following theorem is proved in \cite{KM}.

\begin{theorem}\label{KM} Let $M$ be a compact oriented boundaryless simply connected
$4$-manifold. Let $\omega \in \Hh^2(M;\Z)$ be a characteristic element
  of its intersection form that can be represented by an embedded
  $2$-sphere. Then
$$ \frac{\sigma(M)  - \omega \sqcup \omega }{8} \equiv 0 \ {\rm mod} (2) \ . $$
\end{theorem}

If the intersection form is even, then we can take $\omega =0$ and
recover Rohlin's theorem.  In general a characteristic element
$\omega$ as above can be represented by an oriented surface $F$
embedded in $M$ but not necessarily by a $2$-sphere. For example take
$M= \Pp \cs 8\Qq$.  If $a_0$ is the standard generator of
$\Hh^2(\Pp;\Z)$ represented by a projective complex line, and
similarly $a_j$ for the jth-copy of $\Qq$, then $\omega :=
3a_0+a_1+\dots +a_8$ is characteristic and $\omega \sqcup \omega -
\sigma(M) = 8$, hence $\omega$ cannot be represented by a $2$-sphere
by Theorem \ref{KM}. This motivates the following somewhat informal
\smallskip

{\bf Guess:} {\it (1) Let $M$ be a compact oriented boundaryless
  simply connected $4$-manifold.  Let $\omega \in \Hh^2(M;\Z)$ be a
  characteristic element of its intersection form represented by an
  embedded oriented surface $F\subset M$. Then one expects a formula
  of the type
$$ [\frac{ \sigma(M)- \omega \sqcup \omega }{8}]_{(2)} = \alpha(F) $$
  where $\alpha(F) \in \Z/2\Z$ represents an obstruction to surgery
  $F$ ``within $M$'' to get an embedded $S^2$. Moreover, having in
  mind Pontryagin's computation of $\pi^\infty_2$ depicted in Section
  \ref{small-k} (recall also the study of immersions of surfaces in
  $3$-manifolds in Section \ref{2_imm_3}), it is predictable that
  $\alpha(F)$ is the Arf invariant of some quadratic enhancement of
  $\Hh_1(F;\Z/2\Z)$ (see Section \ref {quadratic}) associated to the
  embedding of $F$ in $M$.
\smallskip

(2) Assuming the isomorphism $\sigma: \Omega_4 \to \Z$, in contrast
with the above derivation of Theorem \ref{rohlin} from the homotopic
result $\pi^\infty_3=\Z/24\Z$, the definition of $\alpha (F)$ as well
as the proof of the congruence should be geometric and possibly
elementary.}

\smallskip

 Accordingly with Freedman-Kirby \cite{FK} (1978), the realization
 therein of the above guess is derived, considerably different in
 details, from one outlined by Casson in 1974 (unpublished).
 Accordingly to the historical appendix by Kharlamov and Viro in
 \cite{GM}, Rohlin announced such a formula at the Moskow IMC 1966 but
 only in a paper of 1972 he used it to solve a conjecture by Gudkov
 concerning Hilbert's 16th problem about the configuration of ovals of
 planar even degree real algebraic curves. The study of this problem
 by means of a $4$-manifold obtained as a branched covering of
 $\PP^2(\C)$ ramified along a given non singular real
 algebraic curve in $\PP^2(\R) \subset \PP^2(\C)$ was introduced by
 Arnol'd \cite{A3} (1971). The basic congruences mod$(8)$ already
 imply non trivial prohibitions for the oval configuration; the finer
 formula as in the above guess implies stronger prohibitions.
 All this holds under weaker hypotheses relaxing the fact that $M$ is
 simply connected; for example $\Omega_1(M)=0$ suffices to define the
 quadratic enhancement by using ``membranes'' (see below) and we can
 even avoid the use of membranes by means of spin structures (see
 \cite{Kirby}).  However, we will keep $M$ to be simply connected and follow
 the treatment of Matsumoto \cite{Mat} given in a paper available in \cite{GM};
 it is the simplest one as it is readily accessible by means of the
 tools developed in the present text.

\subsection{Quadratic enhancement for characteristic surfaces}\label{q-en}
In this section $M$ will be a compact oriented connected smooth
$4$-manifold such that $\Omega_1(M)=0$ (this holds in particular if
$M$ is simply connected) and $F\subset M$ an orientable surface.  Let
$c$ be a simple connected smooth circle on $F$.  As $\Omega_1(M)=0$
and using transversality, there exists a smooth map $f:P \to M$ such
that:
\begin{itemize}
\item $P$ is an oriented compact surface with one boundary component;
\item $f(\partial P)=c$;
\item The restriction of $f$ to a collar $C$ of $\partial P$ in $P$ is an embedding;
\item $f(C\setminus \partial P) \subset M\setminus F$ and $f(C)$ is normal to $F$ along $c$;
\item $f$ is a generic immersion of $P$ in $M$;
\item $f|(P\setminus \partial P)$ is transverse to $F$.
\end{itemize}
Such a map $f$ is said a {\it membrane} along $c$. We simply write $P$
instead of $(P,f)$.  If $M$ is simply connected we can also assume
that $P$ is a $2$-disk, but this is not so important at this
point. For simplicity let us identify $c$ with $\partial P$.  The
pull-back of $T(M)$ on $P$ splits as
$$ f^*T(M)=T(P)\oplus \nu(f) $$ where $\nu(f)$ is said the {\it normal
  bundle of the membrane} and is an oriented bundle of rank $2$.  As
$P$ retracts to a wedge of a finite number of $S^1$ (to one point if
$P$ is a disk), then $\nu(f)$ is isomorphic to a product bundle. Let
us fix a global trivialization $\tau$. This induces a trivialization
of the restriction $\nu(f)|c$. Two trivializations of $\nu(f)$ differ
by a map $g:P\to SO(2)$.  The restriction $g|c$ represents $0$ in
$\Omega_1(SO(2))$, hence it is homotopically trivial (Section \ref{omega1}). 
Then the restricted trivialization $\tau_c$ does not
depend on the choice of $\tau$.  The normal bundle $\nu_c$ of $c$ in
$F$ define a rank-1 orientable sub-bundle of $\nu(f)|c$. Then denote
by $ n(P)$ the {\it number of full twists made by $\nu_c $ with
  respect to $\tau_c$, moving along $c$ in the direction given by its
  orientation as $\partial P$}.  It is not hard to check that
$[n(p)]_{(2)} \in \Z/2\Z$ does not depend on the choice of the
orientation of $P$.

Let now $a\in \Hh_1(F;\Z/2\Z)$. We know (Lemma \ref{connected-rep})
that $a=[c]$ for some simple smooth circle $c$ on $F$. Given a
membrane $P$ along $c$, set
$$ q_F(c,P)= [ n(P)]_{(2)} + [P\bullet F]_{(2)} \in \Z/2\Z \ $$
where $P\bullet F$ is in fact the intersection number between
Int$(P)$ and $F$.
We have
\begin{proposition}\label{Fquadraticform}
  Let $F\subset M$ be an oriented {\rm characteristic surface} of $M$,
  that is $\omega=[F]\in \Hh^2(M;\Z)$ is a characteristic element of
  the intersection form of $M$. Then:
\begin{enumerate}
\item For every simple smooth circle $c$ on $F$, $q_F(c):= q_F(c,P)$
  does not depend on the choice of the membrane $P$ along $c$.
\item For every $a \in \Hh_1(F;\Z/2\Z)$, for every simple smooth
  circle $c$ representing $a$ ($a=[c])$, then $q_F(a):= q_F(c)$ does
  not depend on the choice of the representative $c$.
\item The function $q_F: \Hh_1(F;\Z/2\Z) \to \Z/2\Z$ defined so far is
  a quadratic enhancement of the intersection form on
  $\Hh_1(F;\Z/2\Z)$.
\end{enumerate}
\end{proposition}
\Dim (1) Let $P$ and $P'$ be two membranes along $c$. Up to
``spinning'' $P'$ along $c$, we can assume that $P$ and $P'$ glue
along the common boundary $c$ in such a way that: (i) $\Sigma = P\cup
P'$ is a boundaryless surface generically immersed into $M$; (ii) a
tubular neighbourhood of $c$ in $\Sigma$ is an embedded annulus normal
to $F$, made by two collars $C$ and $C'$ in $P$ and $P'$ respectively,
opposite to each other.  The membranes $P$ and $P'$ determine
respective trivializations $\tau_c$ and $\tau'_c$ which induce
opposite orientations on the fibres of the bundle.  The difference
between $-\tau_c'$ and $\tau_c$ along $c$ is encoded by an element
$d\in \pi_1(SO(2))=\Z$.  One verifies that
$$ \Sigma \bullet \Sigma = d - 2P\bullet P' = d \  {\rm mod} (2)$$
$$\Sigma \bullet F = P\bullet F + P'\bullet F \  {\rm mod} (2)$$
(recall that the self-intersection of $c$ in $F$ $c\bullet c=0$ because $F$ is orientable).
As $F$ is characteristic, then 
$$\Sigma \bullet \Sigma = \Sigma \bullet F \ {\rm mod} (2) $$ 
hence
$$d= P\bullet F + P'\bullet F  \  {\rm mod} (2) \ . $$
On the other hand, 
$$ n(P') = n(P) + d   \  {\rm mod} (2) \ . $$
By combining these relations we eventually get
$$ n(P) + P\bullet F = n(P') + P'\bullet F  \  {\rm mod} (2) $$
as desired. Item (1) is proved.
\smallskip

To achieve (2) (3) we can implement the method illustrated at the end
of Section {quadratic}.  We have defined a function which associate
$q(c)\in \Z/2\Z$ to every simple smooth circle on $F$. It is clear
that $ q(c)=0$ if $c$ is the boundary of a $2$-disk embedded in $F$.
We extend additively this function to every not necessarily connected
simple curve $c=c_1 \amalg \cdots \amalg c_k$ on $F$. If $\gamma$ is
now a curve generically immersed in $F$ with a number say
$r(\gamma)\geq 0$ of normal crossings, every crossing can be
simplified in two ways. Let us call a {\it state} $s$ of $\gamma$ a
system of simplifications at every crossing.  Performing these
simplifications we get a simple curve $c_s$. Set
$$q_F (\gamma, s)= q_F(c_s) + [2r(\gamma)]_{(2)} \ . $$
Then it is enough to prove that $q_F(\gamma):= q_F(\gamma, s)$
does not depend on the choice of the state $s$. Arguing by induction of $r(\gamma)$,
we localize the question at one crossing. If $s$ and $s'$ differ just at one crossing,
then we can use membranes $P$ and $P'$ along the components of $c_s$
and $c_{s'}$ which only differ locally at the crossing. By a direct computation we can
compute $q_F(\gamma,s)$ and $q_F(\gamma,s')$ by using $P$ and $P'$
getting the desired result.

\cvd

For the definition of the Arf invariant of $q_F$ we refer to Section  \ref{quadratic}.
In the next proposition we show that the Arf invariant of $q_F$ only depends
on the characteristic element $\omega = [F] \in \Hh^2(M;\Z)$.

\begin{proposition}\label{alfa-omega} Let $F, F' \subset M$ be oriented 
characteristic surfaces of $M$ representing the same characteristic element
$\omega$ of the intersection form of $M$. Then Arf$(q_F)=$ Arf$(q_{F'})$,
so that $\alpha(\omega):= {\rm Arf}(q_F) \in \Z/2\Z$ is well defined.
\end{proposition}
\Dim We repeat an embedded bordism argument already employed in
Sections \ref {small-k}, \ref {Pin-}. We know that there is an
orientable $3$-dimensional triad $(W, F, F')$ properly embedded into
the triad $(M\times [0,1], M\times \{0\}, M\times \{1\})$ an we can
assume that the restriction to $(W,F,F')$ of the projection onto
$[0,1]$ is a Morse function. Consider the corresponding handle
decomposition of $(W,F,F')$ and the successive surgeries which produce
$F'$ from $F$.  It is immediate that either attaching a $0$-handle or
attaching a $1$-handle to different boundary connected components does
not change the value of Arf. By attaching a $1$-handle to a same
connected component, the boundary is modified by an embedded connected
sum with a copy of $T=S^1\times S^1$; we realizes that there is a
basis $l, m$ of $\Hh_1(T;\Z/2Z)$ such that the intersection form is
represented by the standard matrix $\HH$ and $m$ is the co-core of the
handle, so that $q_T(m)=0$. It follows that Arf$(q_T)=0$, so that the
total Arf does not change also in this case.  Finally we consider the
dual handle decomposition to rule out also $2$ and $3$-handles.

\cvd

\subsection{A digression in classical knot theory}\label{class-knot}
Let us recall a few facts of classical knot theory (see for instance
\cite{Kau}, \cite{Rolf}) that we will use below in the proof of the
main result.  Let $K$ be a knot in $S^3=\partial D^4$ considered up to
ambient isotopy.  Every orientend proper surface $(S,\partial
S)\subset (D^4,S^3)$ such that $\partial S = K$ is ``characteristic"
for $\Hh^2(D^4,S^3;\Z)=0$. So by a similar construction as above we
can define a quadratic form $q_S: \Hh_1(S;\Z/2\Z)\to \Z/2Z$ whose Arf
invariant $\alpha(q_S)\in \Z/2\Z$ eventually depends only on the knot
$K$ so that the {\it Arf invariant of the knot} ${\rm Arf}(K):=
\alpha(q_S)$ is well defined.  It can be computed by means of any
oriented planar diagram $\Dd$ of $K$ as follows.  We can use as $S$
the surface obtained by pushing in $D^4$ the Seifert surface of $K$ in
$S^3$ constructed by means of the Seifert algorithm via the oriented
simplification of the normal crossings of $\Dd$. If $\Dd'$ is a knot
diagram which differs from $\Dd$ just by the over/under branches at
one crossing, denote by $K'$ the corresponding knot.  Performing the
simplification at the given crossing of $\Dd$ (or of $\Dd'$, the
result is the same) we get a diagram $\Dd"$ of a link with two
oriented components $K_1$ and $K_2$. Then one realizes that the
following relation holds involving the linking number of $K_1$ and
$K_2$:
$$ {\rm Arf}(K) = {\rm Arf}(K') + [L(K_1,K_2)]_{(2)} \in \Z/2\Z \ . $$
The linking number mod $(2)$ can be easily computed by means of the diagram
$\Dd"$: the number $c$ of crossings  of $\Dd"$ whose local branches do not
belong to a same constituent knot is even and   $[L(K_1,K_2)]_{(2)}=[c/2]_{(2)}$. 
Moreover, it is well known that one gets a diagram $\Dd_0$ for
the unknot $K_0$ by switching some crossings of $\Dd$; clearly
Arf$(K_0)=0$; then the above relation allows to compute inductively
Arf$(K)$ starting from $\Dd$.

\begin{figure}[ht]
\begin{center}
 \includegraphics[width=5cm]{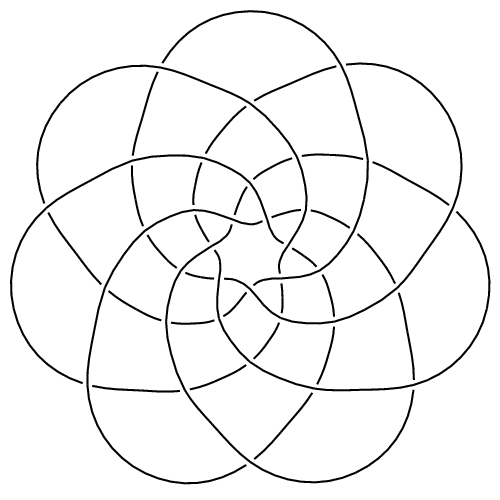}
\caption{\label{(7,6)} A standard diagram of $K(7,6)$.}
\end{center}
\end{figure}

Let $T\subset \R^3$ be the standard torus obtained by rotation of the
planar circle $\{x=0, \ (y-2)^2+ z^2 = 1\}$ around the $z$-axis .  For
every couple $(p,q)$ of coprime integers, the {\it torus knot}
$K(p,q)$ is traced on $T$ turning $p$ times in the direction of the
standard longitude of $T$, $q$ times in the direction of the meridian.
By projection onto the $(x,y)$ coordinate plane, we get a standard
diagram $\Dd(p,q)$ of $K(p,q)$.  We will be interested to the case
$K(s,s-1)$, where $s>1$ is odd (so that $(1-s^2)\equiv 0$ mod $(8)$).
It is known in knot theory (for example by applying the above method
to the diagram $\Dd(s,s-1)$) that
$$ {\rm Arf}(K(s,s-1))= [\frac{1-s^2}{8}]_{(2)} \ . $$

\subsection{The main results}\label{maincong}
We ca state now the main result of this section. 

\begin{theorem}\label{maincong16} Let $M$ be a compact oriented boundaryless simply
connected $4$-manifold. Let $\omega \in \Hh^2(M;\Z)$ be a
characteristic element of the intersection form of $M$. Then
$$ [\frac{ \sigma(M)- \omega \sqcup \omega }{8}]_{(2)} = \alpha(\omega) \ . $$
\end{theorem}

\Dim The proof is based on the classification up to odd stabilization.
First note that if $M=M_1\cs M_2$ is the connected sum of two simply
connected manifolds, then a characteristic element $\omega$ of $M$ is
the sum $\omega = \omega_1 + \omega_2$ of characteristic elements of
$M_1$ and $M_2$ respectively. So if the theorem holds for two members
of the triple $(M,\omega)$, $(M_1,\omega_1)$, $(M_2,\omega_2)$, then
it holds also for the third. By Theorem \ref{scos} we have that
$$M \cs (k \Pp \cs h \Qq) = m\Pp \cs n\Qq$$ for some $k,h,m,n \in
\N$. Then by applying inductively the above remark, it is enough to
prove the theorem for $\Pp$ and $\Qq$. If $\PP^1(\C)\subset \Pp$ is a
complex line, then every characteristic element of $\Pp$ is of the
form $\omega = s[\PP^1(\C)]$, where $s$ is an odd integer; to our aims
it is not restrictive to assume that $s\geq 1$.  The theorem clearly
holds for $s=1$, so let us assume $s>1$.  Then $\omega = [F]$ where
$F$ is any non singular complex projective curve in $\Pp$ defined as
the zero set of a homegeneous polynomial of degree $s$ in the
homogeneous complex coordinates $(z_1,z_2,z_3)$ on $\Pp$. One can
prove indeed (by using the fibration theorem \ref {fib}) that all
these curves are isotopic to each other but this is not so important
for the present discussion.  Let us consider the family of projective
complex curves
$$ F_\epsilon=\{z_1^s+z_2^{s-1}z_3 - \epsilon z_3^s=0\}$$ where
$\epsilon \in \R$, $\epsilon \geq 0$. For $\epsilon =0$, $F_0$ has one
isolated singularity at the point $x_0=(0,0,1)$ and in the affine
coordinates such that $z_3 \neq 0$, it is defined by the equation
$x^s+y^{s-1}=0$.  The best reference for the study of such isolated
singularities of complex planar curves is celebrated Milnor's book
\cite{M6}. Our case is particularly simple and the following facts are
verified.  There is a small round $4$-disk $D$ around $x_0=(0,0)$ in
such affine chart, such that:
\begin{enumerate}
\item $S^3=\partial D$ is tranverse to $F_0$ and $K:=F_0\cap S^3$ is a
  torus knot $K(s,s-1)$.
\item The pair $(D,F_0\cap D)$ is homeomorphic to the pair $(D,cK)$
  where $cK$ denotes the cone with base $K$ and centre at $x_0$.
\item $F_0 \cap (\Pp \setminus {\rm Int} (D))$ is a smooth properly
  embedded $2$-disk. Hence $F_0$ is homeomorphic to $S^2$.
 \end{enumerate}
 
 If $\epsilon >0$ is small enough, then
 
 (i) $F_\epsilon$ is non singular.
 
 (ii) $F_\epsilon \pitchfork S^3$ is an isotopic copy of $K(s,s-1)$
 and $F_\epsilon \cap D$ is properly embedded.
 
 (iii) $F_\epsilon \cap (\Pp \setminus {\rm Int} (D))$ is a smooth
 properly embedded $2$-disk.
 
 Then it is clear that
 $$\alpha(\omega) = {\rm Arf}(q_{F_\epsilon}) = {\rm Arf}(K(s,s-1)) =
 [\frac{1-s^2}{8}]_{(2)} = [\frac {\sigma(\Pp) - \omega \sqcup
     \omega}{8}]_{(2)}$$
and this achieves the case $M=\Pp$. By taking
 into account the change of orientation, the same argument holds as
 well for $M=\Qq$ and the proof is complete.
 
 \cvd 
 
 \subsection{On an extension to non orientable characteristic surfaces}\label{g-m}
We have mentioned a $4$-dimensional approach to Hilbert's 16th
problem where the congruences mod$(16)$ give non trivial
information. In this setting it is quite current to deal with {\it non
  orientable characteristic surfaces} that is representing the
reduction mod$(2)$ of any characteristic element of the intersection
form of some $4$ manifold $M$. This strongly motivates the search for
a further generalization of Theorem \ref{maincong16}.  We limit to
state it.

Let $F\subset M$ be a not necessarily orientable characteristic
surface. Assume that $\Omega_1(M)=0$. Similarly to Section
\ref{2_imm_3} and using membranes as in the above definition of $q_F$,
we can define a quadratic enhancement
$$\hat q_F: \Hh_1(F;\Z/2\Z) \to \Z/4\Z$$
of the intersection form by setting
$$\hat q_F([c])= \hat q_F(c,P)= [\hat n(P)]_{(4)} + 2\cdot([P\bullet
  F]_{(2)}+ c\bullet c) \in \Z/4\Z $$ where $\hat n(P)$ is the number
of {\it half-twists} made by $\nu_c $ with respect to $\tau_c$, moving
along $c$.
The fact that is is well defined is a bit more complicated but not so
much.

Similarly to the discussion made to define the integer
Euler-Poincar\'e characteristic also for non orientable manifolds, we
can define geometrically the self-intersection number $F\bullet F \in
\Z$ by identifying $F$ with the zero section of its normal bundle in
the oriented manifold $M$ and fixing arbitrary compatible local
orientations of $F$ and $F'$ at every point of $F\pitchfork F'$, $F'$
being a section transverse to $F$. By usual arguments this
number does not depend on the arbitrary choices made to compute
it. Recall the Arf-Brown invariant of $\hat q_F$ defined in Section
\ref{quadratic}. Here we denote it by $\hat \alpha (F)\in \Z/8\Z$.
Recall that the multiplication by $2$ determines injective
homomorphisms $\Z/2\Z \to \Z/4\Z \to \Z/8\Z \to \Z/16\Z$. Finally we
can state:

\begin{theorem}\label{g-m-cong}  Let $M$ be a compact oriented boundaryless
  simply connected
$4$-manifold. Let $F\subset M$ be a possibly non orientable surface
  which represents the reduction mod$(2)$ of any characteristic
  element $\omega$ of the intersection form of $M$.  Then
$$ [\sigma(M) - F\bullet F]_{(16)} = 2\cdot \hat \alpha(F) \ . $$
\end{theorem}
\smallskip

If $F$ is oriented we recover Theorem  \ref{maincong16},
because $F\bullet F = \omega \sqcup \omega$, $\hat q_F = 2\cdot q_F$,
$\hat \alpha(F) = 4\cdot \alpha(\omega)$. 

Theorem \ref{g-m-cong} is due to Guillou-Marin \cite{GM}.
There are
several difficulties to overcome. When $F$ is non orientable,
$F\bullet F\in \Z$ cannot be identified with the intersection number
of any bordism classes of $M$. So it is not clear how to reformulate
Proposition \ref {alfa-omega}. The reduction mod$(2)$, say
$\omega_{(2)}$, of any characteristic number $\omega$ does not depend
on the choice of $\omega$. So we should rather prove that
$[F\bullet F + 2\cdot \hat \alpha(F)]_{(16)}$ does not depend on the
choice of the (possibly non orientable) surface $F$ representing $\omega_{(2)}$.
Note also that dealing with non orientable surfaces, the embedded bordism argument
used in the proof of Proposition \ref {alfa-omega} is not immediately available
(recall Remark \ref{unor-seifert}). In the already cited paper \cite{Mat},
Matsumoto gives another proof which by an inductive argument reduces the general statement to forms 
Theorem \ref{maincong16}. In both proofs there are two other basic cases besides
$\Pp$ and $\Qq$, that is $S^4$ with suitably embedded real projective spaces
as characteristic surface.

\section{On the topological classification of smooth $4$-manifolds}\label{moreinfo}
From Rohlin's theorem (1952)  to Donaldson's work in 1982 \cite{Do}, no further
prohibitions to the realizability of unimodular forms by  boundaryless smooth $4$-manifolds
appeared. On the other hand Wall's Theorem \ref{even-class} was the strongest one about 
the extent which the intersection form determines the differential topology of a boundaryless 
$4$-manifold. At the beginning of the 80's two parallel new waves  have revolutionated
the subject. Since Donaldson's work, the introduction of new methods derived from gauge theory,
of differential-geometric/analytic nature and strongly influenced by ideas of theoretical physics, 
have produced amazing new prohibitions and powerful smooth invariants distinguishing
homeomeorphic  but non diffeomorphic smooth $4$-manifolds. Let us recall a few new prohibitions. 

\smallskip

{\bf (Donaldson 1982 \cite{Do})}  {\it If the intersection form of a simply connected, boundaryless smooth
$4$-manifold is definite then it is diagonalizable, that is of the form $k\UU_\epsilon$.}
\smallskip

Donaldson's result means that the arithmetic complication of definite forms
does not concern the intersection forms of smooth $4$-manifolds; hence the problem
of four dimensional smooth realizability is reduced to the indefinite and even case. To this respect we recall: 

{\bf (Furuta 2001 \cite{Fu})} {\it If the intersection form of a simply connected, boundaryless smooth
$4$-manifold is indefinite and even, that is of the type $2h\EE_8 \perp a \HH$, then $a\geq 2|h|+1$.}

\smallskip

The following still is an open conjecture.

\smallskip

{\bf The so called ``11/8" Conjecture:} {\it  If the intersection form of a simply connected, boundaryless smooth
$4$-manifold is indefinite and even, that is of the type $2h\EE_8 \perp a \HH$, then $a \geq 3|h|$.}

\smallskip

If the conjecture holds true, then the rank must be at least $11/8$ times $|\sigma|$. 
Furuta theorem means that the rank is at least $10/8$ times $|\sigma|$.
If the form is indefinite and even we may assume that it is of nonpositive signature 
by changing orientations if necessary, in which case $h\leq 0$. If $a \geq 3|h|$, then the form can be realized
by means of  $|h|K\cs (a-3|h|)(S^2\times S^2)$, where $K$ is the Kummer complex surface
of Example \ref{E8}. Hence a confirmation of the conjecture would achieve the realizability problem.

\smallskip

The other wave had a somewhat more conservative motivation. It was clear at least since Rohlin's `mistake', 
that there were in general actual obstructions in order to apply the Whitney trick in dimension $4$; 
nevertheless one wondered if such a `technical' difficulty could be circunvented in some way in order to prove the 
$5$-dimensional $h$-cobordism theorem. For example in Wall's theorem \ref{Wall-h} this is done by paying the price of 
performing even stabilizations. In this vein, in  73-74 A. Casson introduced so called ``flexible handles" later currently
called  ``Casson handles" (see Lecture I in the second part of \cite{GM}).  Let $M$ be a boundaryless simply connected $4$-manifold
and let $\alpha,\beta \in \Hh_2(M;\Z)$ such that $\alpha \bullet \alpha = \beta \bullet \beta = 0$, $\alpha \bullet \beta = 1$.
Then, by means of a certain `infinite construction', he produced an open set $V$ of $M$ such that
\begin{itemize}
\item $V$ has the proper homotopy type of $S^2\times S^2 \setminus \{pt\}$;
\item $\Hh_2(V;\Z)$ carries the submodute of $\Hh_2(M;\Z)$ generated by $\alpha$ and $\beta$
\end{itemize}

Moreover, he argued (Lecture III of the second part of \cite{GM}) that 
\smallskip

{\it If flexible handles $V$ are diffeomorphic to  the true $S^2\times S^2 \setminus \{pt\}$, then
we could carry out the Whitney process and cancel handles to trivialize
 five dimensional simply connected $h$-cobordisms.}

\smallskip

More information about the flexible handles (at least about its `end') would be also of main importance 
with respect to the realizability problem:
\smallskip

- If such a flexible handle $V$ would be diffeomorphic to the true $S^2\times S^2 \setminus \{pt\}$,
then we could split $M = M'\cs (S^2\times S^2)$ where $M'$ is simply connected and passing from $W$ to $W'$ 
we have surgered out a factor $\HH$ of the intersection form of $M$. 

- If $V$ is diffeomorphic to $N\setminus \{pt\}$ where $N$ is a compact boundaryless $4$-manifold, 
then $M=M' \cs N$ where $N$ has the homotopy type
of $S^2\times S^2$ and again carries $\alpha$ and $\beta$; so $M'$ has the same properties as above.

- If the end of $V$ coincides with the end of an open contractible manifold $V^*$, then by replacing $V$ with $V^*$
we get again $W'$ with $\alpha$ and $\beta$ killed. 

\smallskip

\noindent Notice that before Donaldson's result, there were not known
obstructions in order that the arithmetic splitting of an indefinite and even form  $2h\EE_8 \perp a \HH$ of some
simply connected $4$-manifold $M$ could be realized by a splitting $M' \cs a(S^2\times S^2)$.
After Donaldson we know that the above underlying hope was too optimistic, nevertheless 
the main achievement of \cite{Fr} (1982) was that 
\smallskip

{\it A flexible handle is a `true' $S^2\times S^2 \setminus \{pt\}$, provided one works in the more flexible setting of 
almost smooth $4$-manifolds.}

\smallskip

A topological manifold $N$ is {\it almost smooth} if $N \setminus \{pt\}$ has a smooth structure (which in general cannot be
extended over the whole $N$). Remarkably, more or less  at the same time it was proved in \cite{Q}:
\smallskip

{\it Every boundaryless simply connected topological $4$-manifold is almost smooth.}

\smallskip

This opens the way (via the solution of other hard technical issues) for a complete classification of topological 
simply connected $4$-manifolds, which includes  the fact that {\it every} unimodular symmetric form can be realized 
as the intersection form of a boundaryless simply connected almost smooth $4$-manifolds. Here we limit to state
a few corollaries in our favourite smooth setting.

\smallskip

{\it

(1) {\bf Topological five dimensional $h$-cobordism:}  Every smooth simply connected $5$-dimensional 
$h$-cobordism $(W, M_0, M_1)$ is {\rm homeomorphic} to the product $M_0\times [0,1]$.
In particular $M_0$ and $M_1$ are {\rm homeomorphic} to each other.

(2) {\bf A classification of smooth $4$-manifolds up to homemorphism:} Two smooth simply connected boundaryless $4$-manifolds are {\rm homeomorphic} if and only
if they have isometric intersection forms. 
}

\smallskip

The new gauge theoretical prohibitions and smooth invariants, together with the above topological classifications, lead to a dramatic failure
of the {\it smooth} five dimensional $h$-cobordism theorem and to the existence of a plenty of non diffeomorphic smooth structures on certain
topological $4$-manifolds. In particular we recall that the Kummer complex surface of Example \ref{E8} admits countably many non diffeomorphic
smooth structures \cite{FS}.
Finally we recall that the classification of topological $4$-manifolds includes  the solution of the four dimensional {\it topological}
Poincar\'e conjecture: {\it Every boundaryless topological $4$-manifold which is homotopically equivalent to $S^4$ is homeomorphic to $S^4$}.
It is not known if every {\it smooth} boundaryless $4$-manifold which is homotopically equivalent to $S^4$ is {\it diffeomorphic} to $S^4$.
This smooth four dimensional Poincar\'e conjecture presumably is the main basic open question about smooth $4$-manifolds.

 \bibliographystyle{amsalpha}

\chapter*{Appendix: baby categories}\label{2TD-CAT-APP}
Along the text we make some (very moderate indeed) use of the language of categories.
We collect in this appendix the few necessary notions.

A {\it category} $\CC$ consists of three things: 
\begin{enumerate}
\item A class of {\it objects} $X$;
\item For every ordered pair of objects $(X,Y)$, a set Hom$(X,Y)$ of {\it morphisms} (also called {\it arrows})
$f: X \mapsto Y$;
\item For every ordered triple $(X,Y,Z)$ of objects, a {\it composition function of arrows}
$$\circ : {\rm Hom}(X,Y)\times {\rm Hom}(Y,Z)\to {\rm Hom}(X,Z), \ (f,g)\to g\circ f \ . $$
\end{enumerate}

We require that the following properties are satisfied:

\begin{enumerate}
\item {\it (Associativity)} Whenever the involved compositions make sense, we have $h\circ (g\circ f) = (h\circ g)\circ f$;
\item {\it (Existence of the identity)} For every object $X$, there is a (necessarily unique) arrow $1_X\in {\rm Hom}(X,X)$
such that  $1_X\circ f= f$, $g\circ 1_X=g$, whenever the compositions make sense.
\end{enumerate}

A morphism $f\in {\rm Hom}(X,Y)$ is {\it an equivalence} in the category $\CC$ if there exists a (necessarily unique)
morphism $g\in {\rm Hom}(Y,X)$ such that $f\circ g = 1_X$ and $g\circ f = 1_Y$.
\medskip

A fundamental example  is the category of sets, denoted by {\bf SET}, which has as objects the class
of all sets, while Hom$(X,Y)$ consists of the set of all maps from $X$ to $Y$. $1_X$ is the identity map, while
the equivalences are the bijective maps. We know a lot of sub-categories of {\bf SET} obtained by specializing
both objects and arrows: the categories of groups and group homomorphisms, of vector spaces (on a given scalar field) and linear maps, 
of topological spaces and continuous maps, of smooth manifolds and smooth maps $\dots$. The equivalences
are the isomorphisms, the homemorphisms, the diffeomorphisms, \dots .

A single group $G$ can be considered as a category with just $G$
as unique object, while Hom$(G,G) \sim G$, by associating to every $h\in G$ the morphism by left multiplication by $h$,
$L_h:G\to G, \ g\to hg$. In this category all morphisms are equivalences.

Not every category is a subcategory of {\bf SET}. For example, starting from the category of topological spaces and continuous maps 
we can construct a new category with the same class of objects, and as arrows the {\it homotopy classes} of continuous maps from
$X$ to $Y$. The fact that associativity holds is left as an exercise.

If $X$ is a path connected topological space, we can consider the category whose objects are the points of $X$ and Hom$(x,y)$
consists of the homotopy classes $[\alpha]$ of paths in $X$ connecting $x$ and $y$. One can verify that every morphism in this
category is an equivalence (we say that it is a {\it groupoid}).

\medskip

Given two categories $\CC$ and $\DD$, a {\it covariant functor} $\Ff: \CC \Rightarrow \DD$ fron $\CC$ to $\DD$ assigns
to every object $X$ of $\CC$, an object $\Ff(X)$ of $\DD$, to every arrow  $f\in {\rm Hom}(X,Y)$ of $\CC$, an arrow
$\Ff(f): \Ff(X)\mapsto \Ff(Y)$ of $\DD$ in such a way that the following properties are satisfied:
\begin{enumerate}
\item For every object $X$ of $\CC$, $\Ff(1_X) = 1_{\Ff(X)}$;
\item $\Ff(g\circ f)= \Ff(g)\circ \Ff(f)$, whenever the composition is defined.
\end{enumerate}
\smallskip

A {\it contravariant functor} assigns to every $f\in {\rm Hom}(X,Y)$, an arrow $\Ff(f)\in {\rm Hom}(\Ff(Y),\Ff(X))$
in such a way that $\Ff(g\circ f)= \Ff(f)\circ \Ff(g)$. A basic example of contravariant functor if the functor from the category
of vector spaces (on a given scalar field) to itself such that for every $V$, $\Ff(V)=V^*$ the dual space, and for every
linear map $f:V\to W$, $\Ff(f)=f^t$ the transposed map of $f$, $f^t: W^* \to V^*$, $f^t(\phi) = \phi \circ f$.  

Let $\Ff$ and $\Gg$ be two say covariant functors from $\CC$ to $\DD$. 
A {\it natural transformation} $T$ from  $\Ff$ to $\Gg$ is a rule assigning to every object $X$ of $\CC$, a morphism
$T_X: \Ff(X)\mapsto \Gg(X)$ such that for every $f\in {\rm Hom}(X,Y)$ of $\CC$, $\Gg(f) \circ T_X = T_Y \circ \Ff(f)$.
If for every $X$, $T_X$ is an equivalence, then $T$ is called a {\it natural equivalence of functors}.

For example a $\Delta$-complex mentioned in the text can be abstractly defined as 
being a contravariant functor from the category $\Delta$ to the category {\bf SET}, where
$\Delta$ has as objects the ordered sets $\Delta^n=\{0,1,\dots, n-1\}$, $n\in \N$, and as arrow the
strictly increasing maps $\Delta^k \to \Delta^n$, $k\leq n$. Maps beteween $\Delta$-complexes would be defined
as natural trasnformations of the corresponding functors.

\end{document}